\documentclass{article}
\usepackage{latexsym, amscd, amsfonts, eucal, mathrsfs, amsmath, amssymb, amsthm, xypic,xr, makeidx, stmaryrd}
\def\Z{\mathbf{Z}}

\DeclareMathOperator{\acti}{ac}
\DeclareMathOperator{\CatP}{{\mathfrak P}}
\DeclareMathOperator{\Mul}{Mul}

\DeclareMathOperator{\Lin}{{\mathcal L}in}
\DeclareMathOperator{\Bimodd}{PBimod}
\DeclareMathOperator{\Bimod}{{\mathcal B}{imod}}
\DeclareMathOperator{\Spectra}{Sp}

\DeclareMathOperator{\sNerve}{N}

\DeclareMathOperator{\MMod}{\overline{{\calM}od}}
\DeclareMathOperator{\PrAlg}{ \overline{\Alg}}
\DeclareMathOperator{\tildePAlg}{ \widetilde{ \Alg}}
\DeclareMathOperator{\tildeMod}{ \widetilde{{\calM}od}}
\DeclareMathOperator{\famm}{fam}

\newcommand{\Biod}[3]{_{#1}{{\mathcal B}imod}_{#2}(#3)}

\newcommand{\mset}[1]{{(\mathcal S}et^{+}_{\Delta})_{/#1}} 

\DeclareMathOperator{\CommOp}{Comm}

\DeclareMathOperator{\LPrest}{L\mathcal{P}r}

\DeclareMathOperator{\LPress}{\widehat{\Cat}_{\infty}^{\LPrest}}

\DeclareMathOperator{\Spectralheart}{\SSet_{\infty}^{\heartsuit}}

\DeclareMathOperator{\Stabb}{\sigma}

\DeclareMathOperator{\lax}{lax}
\DeclareMathOperator{\Act}{Act}
\DeclareMathOperator{\LinSeg}{{{\mathcal L}in}_{\ast}}
\DeclareMathOperator{\FinSeg}{\Gamma}
\DeclareMathOperator{\InjSeg}{\Gamma^{si}}
\DeclareMathOperator{\PFinSeg}{\Gamma^{\ast}}

\DeclareMathOperator{\Envv}{Env}
\DeclareMathOperator{\PreOp}{{\mathcal P}Op_{\infty}}

\DeclareMathOperator{\Fintimes}{\Gamma^{\times}}

\newcommand{\Andre}{Andr\'{e}}

\DeclareMathOperator{\free}{free}
\DeclareMathOperator{\finite}{fp}

\DeclareMathOperator{\LPressStab}{\widehat{\Cat}_{\infty}^{\LPrest, \Stabb}}

\newcommand{\EInfty}{ {\mathfrak E}_{\infty}}

\newcommand{\bbox}{\boxempty}

\DeclareMathOperator{\Baar}{Bar}

\DeclareMathOperator{\disc}{disc}
\DeclareMathOperator{\mSet}{{\mathcal{S}et}_{\Delta}^{+}}
\DeclareMathOperator{\conn}{conn}

\DeclareMathOperator{\sMon}{sMon}

\newcommand{\OpE}[1]{{\mathcal E}_{#1}}

\newcommand{\seg}[1]{{\langle #1 \rangle}}
\newcommand{\Colp}[1]{\rho^{#1}}
\newcommand{\Colll}[2]{ \rho^{#1}_{#2}}
\newcommand{\nostar}[1]{{{\langle #1 \rangle}^{\circ}}}

\DeclareMathOperator{\connSpectra}{\Spectra^{conn}}

\DeclareMathOperator{\calH}{\mathcal{H}}
\DeclareMathOperator{\Mon}{Mon}

\DeclareMathOperator{\Cart}{Cart}
\DeclareMathOperator{\Alg}{Alg}
\DeclareMathOperator{\CAlg}{CAlg}

\DeclareMathOperator{\rght}{R}

\DeclareMathOperator{\Mult}{{\mathcal M}}

\newcommand{\MonLPres}{\widehat{\Cat}_{\infty}^{\pres,\otimes}}
\newcommand{\toposref}[1]{T.\ref{HTT-#1}}
\newcommand{\stableref}[1]{S.\ref{STA-#1}}
\newcommand{\monoidref}[1]{M.\ref{MON-#1}}

\newcommand{\bicatref}[1]{B.\ref{BIC-#1}}
\DeclareMathOperator{\pres}{Pr}
\newcommand{\MonStab}{\widehat{\Cat}_{\infty}^{\Stabb, \otimes}}
\newcommand{\degree}{\text{o}}
\newcommand{\bfA}{{\mathbf A}}

\DeclareMathOperator{\cDelta}{{\bf \Delta}}

\DeclareMathOperator{\Set}{\mathcal{S}et}
\DeclareMathOperator{\sSet}{\mathcal{S}et_{\Delta}}

\DeclareMathOperator{\Nerve}{N}
\DeclareMathOperator{\Cat}{\mathcal{C}at}
\newcommand{\h}[1]{\rm{h} \! #1}

\newcommand{\Adjoint}[4]{\xymatrix@1{#2 \ar@<.4ex>[r]^-{#1} & #3 \ar@<.4ex>[l]^-{#4}}}

\renewcommand{\boxtimes}{\odot}

\newcommand{\et}{\'{e}t}

\newcommand{\bigdot}{\bullet}

\DeclareMathOperator{\Stab}{Stab}

\DeclareMathOperator{\bHom}{Map}

\DeclareMathOperator{\Sym}{Sym}

 \DeclareMathOperator{\End}{End}
 
\DeclareMathOperator{\Mod}{\mathcal{M}od}

\DeclareMathOperator{\calM}{\mathcal{M}}

 \DeclareMathOperator{\G}{\mathbf{G}}

\DeclareMathOperator{\DerRing}{\mathcal{SCR}}

\DeclareMathOperator{\Q}{\mathbf{Q}}
\DeclareMathOperator{\Tor}{Tor} 
\DeclareMathOperator{\colim}{colim}

\DeclareMathOperator{\calA}{\mathcal{A}}
\DeclareMathOperator{\bd}{\partial}

\newcommand{\Sphere}{S}

\DeclareMathOperator{\calE}{\mathcal{E}}

\DeclareMathOperator{\calO}{\mathcal{O}}
\DeclareMathOperator{\calT}{\mathcal{T}}

\DeclareMathOperator{\MCat}{MCat}

\DeclareMathOperator{\calB}{\mathcal{B}}
\DeclareMathOperator{\calK}{\mathcal{K}}
\DeclareMathOperator{\Spec}{{\bf Spec}}
\DeclareMathOperator{\CatTriv}{{\bf Triv}}
\DeclareMathOperator{\Triv}{{\mathcal T}riv}
\DeclareMathOperator{\Ass}{{\mathcal A}ssoc}
\DeclareMathOperator{\CatAss}{{\bf Assoc}}
\DeclareMathOperator{\calF}{\mathcal{F}}

\DeclareMathOperator{\Hom}{Hom} 
\DeclareMathOperator{\HH}{H} 
\DeclareMathOperator{\id}{id} \DeclareMathOperator{\Fun}{Fun}
\DeclareMathOperator{\calC}{\mathcal{C}}
\DeclareMathOperator{\calI}{\mathcal{I}}
\DeclareMathOperator{\calN}{\mathcal{N}}
\DeclareMathOperator{\calJ}{\mathcal{J}}

\DeclareMathOperator{\calR}{\mathcal{R}}
\DeclareMathOperator{\SSet}{\mathcal{S}}
\DeclareMathOperator{\SCop}{{\mathcal M}Cat_{\Delta}}

\DeclareMathOperator{\calX}{\mathcal{X}}

\DeclareMathOperator{\calY}{\mathcal{Y}}

\DeclareMathOperator{\op}{op}
\DeclareMathOperator{\calD}{\mathcal{D}}
\DeclareMathOperator{\Ind}{Ind} 
\DeclareMathOperator{\calP}{\mathcal{P}} \topmargin=0in

\newtheorem{theorem}{Theorem}[subsection]
\newtheorem{lemma}[theorem]{Lemma}
\newtheorem{proposition}[theorem]{Proposition}

\newtheorem{corollary}[theorem]{Corollary}

\theoremstyle{definition}
\newtheorem{definition}[theorem]{Definition}
\newtheorem{construction}[theorem]{Construction}
\newtheorem{example}[theorem]{Example}

\newtheorem{notation}[theorem]{Notation}

\newtheorem{warning}[theorem]{Warning}
\newtheorem{remark}[theorem]{Remark}
\newtheorem{variant}[theorem]{Variant}

\makeindex

\begin{document}

\title{Derived Algebraic Geometry III: Commutative Algebra}

\maketitle
\tableofcontents

\section*{Introduction}

Let $K$ denote the functor of {\it complex $K$-theory}, which associates to every compact Hausdorff space $X$ the Grothendieck group $K(X)$ of isomorphism classes of complex vector bundles on $X$. The functor $X \mapsto K(X)$ is an example of a {\it cohomology theory}: that is, one can define more generally a sequence of abelian groups $\{ K^n(X,Y) \}_{n \in \Z}$ for every inclusion of topological spaces $Y \subseteq X$, in such a way that the Eilenberg-Steenrod axioms are satisfied (see \cite{eilenbergsteenrod}). However, the functor $K$ is endowed with even more structure: for every topological space $X$, the abelian group $K(X)$ has the structure of a commutative ring (when $X$ is compact, the multiplication on $K(X)$ is induced by the operation of tensor product of complex vector bundles). One would like that the ring structure on $K(X)$ is a reflection of the fact that $K$ itself has a ring structure, in a suitable setting.

To analyze the problem in greater detail, we observe that the functor $X \mapsto K(X)$ is {\it representable}. That is, there exists a topological space $Z = \Z \times BU$ and a universal class
$\eta \in K(Z)$, such that for every sufficiently nice topological space $X$, the pullback of $\eta$ induces a bijection
$[X,Z] \rightarrow K(X)$; here $[X,Z]$ denotes the set of homotopy classes of maps from $X$ into $Z$.
According to Yoneda's lemma, this property determines the space $Z$ up to homotopy equivalence. Moreover, since the functor $X \mapsto K(X)$ takes values in the category of commutative rings, the topological space $Z$ is automatically a {\it commutative ring object} in the homotopy category $\calH$ of topological spaces. That is, there exist addition and multiplication maps $Z \times Z \rightarrow Z$, such that all of the usual ring axioms are satisfied up to homotopy. Unfortunately, this observation is not very useful. We would like to have a robust generalization of classical algebra which includes a good theory of modules, constructions like localization and completion, and so forth. The homotopy category $\calH$ is too poorly behaved to support such a theory. 

An alternative possibility is to work with commutative ring objects in the category of topological spaces itself: that is, to require the ring axioms to hold ``on the nose'' and not just up to homotopy. Although this does lead to a reasonable generalization of classical commutative algebra, it not sufficiently general for many purposes. For example, if $Z$ is a topological commutative ring, then one can always extend the functor $X \mapsto [X, Z]$ to a cohomology theory. However, this cohomology theory is not very interesting: in degree zero, it simply gives the following variant of classical cohomology:
$$\prod_{n \geq 0} \HH^{n}( X; \pi_{n} Z).$$
In particular, complex $K$-theory cannot be obtained in this way. In other words, the
$Z = \Z \times BU$ for stable vector bundles cannot be equipped with the structure of a topological commutative ring. This reflects the fact that complex vector bundles on a space $X$ form a {\it category}, rather than just a set. The direct sum and tensor product operation on complex vector bundles satisfy the ring axioms, such as the distributive law
$$ \calE \otimes( \calF \oplus \calF' ) \simeq ( \calE \otimes \calF) \oplus ( \calE \otimes \calF'),$$
but only up to isomorphism. However, although $\Z \times BU$ has less structure than a commutative ring, it has more structure than simply a commutative ring object in the homotopy category $\calH$, because the isomorphism displayed above is actually {\em canonical} and satisfies certain coherence conditions (see \cite{categoricalring} for a discussion). 

To describe the kind of structure which exists on the topological space $\Z \times BU$, it is convenient to introduce the language of {\it commutative ring spectra}, or, as we will call them, {\it $E_{\infty}$-rings}. Roughly speaking, an $E_{\infty}$-ring can be thought of as a space $Z$ which is equipped with an addition and a multiplication for which the axioms for a commutative ring hold not only up to homotopy, but up to {\it coherent} homotopy. The $E_{\infty}$-rings play a role in stable homotopy theory analogous to the role played by commutative rings in ordinary algebra. As such, they are the fundamental building blocks of {\it derived algebraic geometry}.

Our goal in this paper is to introduce the theory of $E_{\infty}$-ring spectra from an $\infty$-categorical point of view. An ordinary commutative ring $R$ can be viewed as a {\it commutative algebra object} in the category of abelian groups, which we view as endowed with a symmetric monoidal structure given by tensor product of abelian groups. To obtain the theory of $E_{\infty}$-rings 
we will use the same definition, except that we will replace the ordinary category of abelian groups by the $\infty$-category of {\it spectra}. In order to carry out the details, we will need to develop an $\infty$-categorical version of the theory of symmetric monoidal categories. The construction of this theory will occupy our attention throughout most of this paper.

The notion of a symmetric monoidal $\infty$-category can be regarded as a special case of the more general notion of an {\it $\infty$-operad}. The theory of
$\infty$-operads bears the same relationship to the theory of $\infty$-categories that the classical
theory of colored operads (or multicategories) bears to ordinary category theory. Roughly speaking,
an $\infty$-operad $\calC^{\otimes}$ consists of an underlying $\infty$-category $\calC$
equipped with some additional data, given by ``multi-input'' morphism spaces
$\Mul_{\calC}( \{ X_i \}_{i \in I}, Y)$ for $Y \in \calC$, where $\{ X_i \}_{i \in I}$ is any finite collection of objects of $\calC$; we recover the usual morphism spaces in $\calC$ by taking
the set $I$ to contain a single element. We will give a more precise definition in
\S \ref{comm1}, and construct a number of examples.

\begin{remark}
An alternate approach to the theory of $\infty$-operads has been proposed by
Cisinski and Moerdijk, based on the formalism of {\em dendroidal sets} (see \cite{cismoer} for details). It seems very likely that their theory is equivalent to the one presented here. More precisely, there should be a Quillen equivalence between the category of dendroidal sets and the category $\PreOp$ of $\infty$-preoperads which we describe in \S \ref{comm1.8}. 

Nevertheless, the two perspectives differ somewhat at the level of technical detail. Our approach has the advantage of being phrased entirely in the language of simplicial sets, which allows us to
draw heavily upon the preexisting theory of $\infty$-categories (as described in \cite{topoi}) and to avoid direct contact with the combinatorics of trees (which play an essential role in the definition of a dendroidal set). However, the advantage comes at a cost: though our $\infty$-operads can be
described using the relatively pedestrian language of simplicial sets, the actual simplicial sets which
we need are somewhat unwieldy and tend to encode information in a inefficient way. Consequently, some of the proofs presented here are perhaps more difficult than they need to be: for example,
we suspect that Theorems \ref{oplk} and \ref{colmodd} admit much shorter proofs in the dendroidal setting.
\end{remark}

If $\calC^{\otimes}$ is a symmetric monoidal $\infty$-category (or, more generally, an
$\infty$-operad), then we can define an $\infty$-category
$\CAlg(\calC)$ of {\it commutative algebra objects of $\calC$}.
In \S \ref{comm4}, we will study these $\infty$-categories in some detail: in particular, we will
show that, under reasonable hypotheses, there is a good theory of {\em free algebras}
which can be used to construct colimits of diagrams in $\CAlg(\calC)$. If
$A \in \CAlg(\calC)$ is a commutative algebra object of $\calC$, we can also associate
to $A$ an $\infty$-category $\Mod_{A}(\calC)$ of {\it $A$-module objects of $\calC$}.
We will define and study these $\infty$-categories of modules in \S \ref{comm6}, and show
that (under some mild hypotheses) they inherit the structure of symmetric monoidal $\infty$-categories.

The symmetric monoidal $\infty$-category of greatest interest to us is the $\infty$-category
$\Spectra$ of {\it spectra}. In \S \ref{crung}, we will show that the smash product monoidal structure on
$\Spectra$ (constructed in \S \monoidref{hummingburg}) can be refined, in an essentially unique way, to a symmetric monoidal structure on $\Spectra$. 
We will define an {\it $\OpE{\infty}$-ring} to be a commutative algebra object in the $\infty$-category
of $\Spectra$. In \S \ref{comm8}, we will study the $\infty$-category $\CAlg(\Spectra)$ of $\OpE{\infty}$-rings and establish some of its basic properties. In particular, we will show that the $\infty$-category $\CAlg( \Spectra)$ contains, as a full subcategory, the (nerve of the) ordinary category of commutative rings. In other words, the theory of $\OpE{\infty}$-rings really can be viewed as a generalization of classical commutative algebra. We will also introduce the notion of a {\it Noetherian} $E_{\infty}$-ring, and prove a version of the Hilbert basis theorem.

\begin{remark}
There are many other definitions of commutative ring spectrum available in the literature. In \S \ref{comm9}, we will introduce the tools needed to compare our theory of $E_{\infty}$-rings with some of these other approaches. More precisely, we will show that if $\bfA$ is a simplicial model category
equipped with a symmetric monoidal structure satisfying suitable hypotheses, then the $\infty$-category $\CAlg( \Nerve(\bfA^{\degree} ) )$ is equivalent to the $\infty$-category underlying the category of {\em strictly} commutative algebras in $\bfA$ itself. Though the hypotheses needed for our argument are somewhat restrictive, they will apply in particular when $\bfA$ is the category of {\it symmetric spectra}; this will imply that our theory of $\OpE{\infty}$-rings is equivalent to the (more classical) theory of {\it commutative ring spectra}.
\end{remark}

\subsection*{Notation}

Throughout this paper, we will freely use the theory of $\infty$-categories developed in
\cite{topoi}. We will also use \cite{DAGStable} as a reference for the theory of stable $\infty$-categories, and \cite{monoidal} as reference for the theory of monoidal $\infty$-categories.
References to \cite{topoi} will be indicated by use of the letter T, references to \cite{DAGStable} will be indicated by use of the letter S, and references to \cite{monoidal} will be indicated by use of the letter M. For example, Theorem \toposref{mainchar} refers to Theorem \ref{HTT-mainchar} of \cite{topoi}.

If $p: X \rightarrow S$ is a map of simplicial sets and $s$ is a vertex of $S$, we will typically write
$X_{s}$ to denote the fiber $X \times_{S} \{s\}$. If $f: s \rightarrow s'$ is an edge of $S$
such that the induced map $X \times_{S} \Delta^1 \rightarrow \Delta^1$ is a coCartesian fibration,
then $p$ determines a functor from $X_{s}$ to $X_{s'}$ which we will denote by $f_{!}$ (this functor is well-defined up to canonical equivalence).

Given maps of simplicial sets $X \stackrel{f}{\rightarrow} S \leftarrow Y$, we let
$\Fun_{S}( X, Y)$ denote the simplicial set $Y^{X} \times_{ S^{X} } \{f \}$. 
 
 \begin{warning}
 Throughout this paper, we will encounter mild variations of the ideas presented in \cite{monoidal}. We will often recycle the notation of \cite{monoidal}, even in cases when this notation conflicts with its previous usage. For example, a symmetric monoidal structure on an $\infty$-category $\calC$ is encoded by a coCartesian fibration $\calC^{\otimes} \rightarrow \Nerve(\FinSeg)$ (see \S \ref{comm1}); in \cite{monoidal}, we used the same notation $\calC^{\otimes}$ to denote a simplicial set
equipped with a coCartesian fibration $\calC^{\otimes} \rightarrow \Nerve(\cDelta)^{op}$ which
encodes a (possibly noncommutative) monoidal structure on $\calC$.
\end{warning}
 
\begin{notation}\label{coll}
If $I$ is a finite set, we will often let $I_{\ast}$ denote the pointed set $I \coprod \{ \ast \}$ obtained
by adjoining a disjoint base point to $I$. For each element $i \in I$, we let $\Colp{i}: I_{\ast} 
\rightarrow \{ 1 \}_{\ast}$ be the map described by the formula
$$ \Colp{i}(j) = \begin{cases} 1 & \text{ if } i = j \\
\ast & \text{ if } i \neq j. \end{cases}$$
\end{notation}

\section{$\infty$-Operads}\label{comm1}

Our goal in this section is to develop an $\infty$-categorical version of the theory of
{\em colored operads}. We will begin in \S \ref{comm1.1} by reviewing the classical definition of a colored operad. The structure of an arbitrary colored operad $\calO$ is completely encoded
by a category which we will denote by $\calO^{\otimes}$, together with a forgetful functor
from $\calO^{\otimes}$ to Segal's category $\FinSeg$ of pointed finite sets. Motivated by this observation, we will define an {\it $\infty$-operad} to be an $\infty$-category
$\calO^{\otimes}$ equipped with a map $\calO^{\otimes} \rightarrow \Nerve(\FinSeg)$,
satisfying an appropriate set of axioms (Definition \ref{colopi}). 

Every symmetric monoidal category $\calC$ gives rise to a colored operad, from which we can recover the original symmetric monoidal category (up to symmetric monoidal equivalence). Consequently, we can view the theory of symmetric monoidal categories as a specialization of the theory of
colored operads. There is a corresponding specialization in the $\infty$-categorical case:
we will define a {\it symmetric monoidal $\infty$-category} to be an $\infty$-operad
$p: \calC^{\otimes} \rightarrow \Nerve(\FinSeg)$ for which the map $p$ is a coCartesian fibration.
More generally, we can associate to any $\infty$-operad $\calO^{\otimes}$ a theory
of {\it $\calO$-monoidal $\infty$-categories}, which are encoded by coCartesian fibrations
$\calC^{\otimes} \rightarrow \calO^{\otimes}$. The relevant definitions will be given in
\S \ref{comm1.2}. In the special case where $\calO^{\otimes}$
is the associative $\infty$-operad (Example \ref{xe2}), the theory of $\calO^{\otimes}$-monoidal $\infty$-categories is essentially equivalent to the theory of monoidal $\infty$-categories
introduced in \cite{monoidal} (Corollary \ref{stix}).

If $\calC^{\otimes}$ is a symmetric monoidal $\infty$-category, then there is an associated
theory of {\it commutative algebra objects} of $\calC^{\otimes}$: these can be defined as 
sections of the fibration $\calC^{\otimes} \rightarrow \Nerve(\FinSeg)$ which satisfy a mild additional condition. More generally, to any fibration of $\infty$-operads $\calC^{\otimes} \rightarrow \calO^{\otimes}$, we can associate an $\infty$-category $\Alg_{\calO}(\calC)$ of $\calO$-algebras in $\calC$.
We will give the basic definitions in \S \ref{algobj}, and undertake a much more comprehensive study in
\S \ref{comm4}.

The next few sections are devoted to giving examples of symmetric monoidal $\infty$-categories. In
\S \ref{comm2}, we will show that if $\calC$ is an $\infty$-category which admits finite products, then
the formation of products endows $\calC$ with a symmetric monoidal structure, which we will call
the {\it Cartesian symmetric monoidal structure} on $\calC$. There is a dual theory of 
{\em coCartesian symmetric monoidal structures} (on $\infty$-categories $\calC$ which admit finite coproducts) which we will describe in \S \ref{jilun}. In \S \ref{monenv}, we will show that
any $\infty$-operad $\calO^{\otimes}$ has a canonical enlargement to a symmetric monoidal
$\infty$-category $\Envv(\calO)^{\otimes}$, which we will call the {\it symmetric monoidal envelope}
of $\calO^{\otimes}$. Finally, in \S \ref{subcatmon}, we will study conditions under which
a subcategory of a symmetric monoidal $\infty$-category $\calC$ (such as a localization of $\calC$) inherits a symmetric monoidal structure.

In \S \ref{comm1.8}, we will show that the collection of $\infty$-operads can be identified
with the fibrant-cofibrant objects of a simplicial model category $\PreOp$, which we call
the {\it category of $\infty$-preoperads}. The category $\PreOp$ is in fact a monoidal model category:
the monoidal structure underlies the classical {\em tensor product} of operads. In \S \ref{comm1.9}, we will study a particular $\infty$-operad $\OpE{0}$ which is idempotent with respect to this tensor product:
tensor product with $\OpE{0}$ therefore determines a localization from the $\infty$-category
of $\infty$-operads to itself, whose essential image we will describe as the class of
{\it unital} $\infty$-operads. Finally, in \S \ref{comm1.10} we will introduce the notion of an
{\it $\infty$-operad family}: this is a mild generalization of the notion of an $\infty$-operad, which
will play an important role in \S \ref{comm6}.

\subsection{Basic Definitions}\label{comm1.1}

Our goal in this section is to introduce the definition of an {\it $\infty$-operad} (Definition \ref{colopi}), laying the foundations for the remainder of the paper. We begin by reviewing the classical theory of colored operads.

\begin{definition}\label{colop}
A {\it colored operad} $\calO$ consists of the following data:
\begin{itemize}
\item[$(1)$] A collection $\{ X, Y, Z, \ldots \}$ which we will refer to as the collection
of {\it objects} or {\it colors} of $\calO$. We will indicate that $X$ is an object of $\calO$ by
writing $X \in \calO$. 
\item[$(2)$] For every finite set of colors $\{ X_i \}_{i \in I}$ in $\calO$ and every
color $Y \in \calO$, a set $\Mul_{\calO}( \{ X_i \}_{i \in I}, Y)$, which we call the
set of {\it morphisms} from $\{ X_i \}_{i \in I}$ to $Y$.
\item[$(3)$] For every map of finite sets $I \rightarrow J$ having fibers
$\{ I_{j} \}_{j \in J}$, a finite collection of objects $\{ X_i \}_{i \in I}$, every finite
collection of objects $\{ Y_j \}_{j \in J}$, and every object $Z \in \calO$, a composition map
$$ \prod_{j \in J} \Mul_{\calO}( \{ X_i \}_{i \in I_j}, Y_j) \times
\Mul_{\calO}( \{Y_j \}_{j \in J}, Z) \rightarrow \Mul_{\calO}( \{ X_i \}_{i \in I}, Z).$$

\item[$(4)$] A collection of morphisms $\{ \id_{X} \in \Mul_{\calO}( \{X\}, X) \}_{X \in \calO}$ which are both left and right units for the composition on $\calO$ in the following sense: for every finite collection of colors $\{ X_i \}_{i \in I}$ and every color $Y \in \calO$, the compositions
\begin{eqnarray*} \Mul_{\calO}( \{ X_i \}_{i \in I}, Y) & \simeq & 
\Mul_{\calO}( \{ X_i \}_{i \in I}, Y) \times \{ \id_Y \} \\
&  \subseteq &
\Mul_{\calO}( \{ X_i \}_{i \in I}, Y) \times \Mul_{\calO}( \{Y \}, Y) \\
&  \rightarrow & \Mul_{\calO}( \{X_i \}_{i \in I}, Y) \end{eqnarray*}
\begin{eqnarray*} \Mul_{\calO}( \{X_i\}_{i \in I}, Y) & \simeq & 
(\prod_{i \in I} \{ \id_{X_{i}} \})  \times \Mul_{\calO}( \{ X_i \}_{i \in I}, Y) \\
& \subseteq & (\prod_{i \in I} \Mul_{\calO}( \{X_i\}, X_i) ) \times \Mul_{\calO}( \{X_i\}_{i \in I}, Y) \\
& \rightarrow & \Mul_{\calO}( \{X_i\}_{i \in I}, Y) \end{eqnarray*}
both coincide with the identity map from $\Mul_{\calO}( \{X_i \}_{i \in I}, Y)$ to itself.

\item[$(5)$] Composition is required to be associative in the following sense: for every
sequence of maps $I \rightarrow J \rightarrow K$ of finite sets together with collections of
colors $\{ W_i \}_{i \in I}$, $\{ X_{j} \}_{j \in J}$, $\{ Y_k \}_{k \in K}$, and every color $Z \in \calO$, the
diagram
$$ \xymatrix@!C=88pt{ & \prod_{j \in J} \Mul_{\calO}( \{ W_i \}_{i \in I_j}, X_j))
\times \prod_{k \in K} \Mul_{\calO}( \{ X_j \}_{j \in J_k}, Y_k)
\times \Mul_{\calO}( \{ Y_k \}_{k \in K}, Z) \ar[dl] \ar[dr] & \\
 \prod_{k \in K} \Mul_{\calO}( \{ W_i \}_{i \in I_k}, Y_k)
\times \Mul_{\calO}( \{ Y_k \}_{k \in K}, Z) \ar[dr] & & 
 \prod_{j \in J} \Mul_{\calO}( \{ W_i \}_{i \in I_j}, X_j)
\times \Mul_{\calO}( \{ X_j \}_{j \in J}, Z) \ar[dl] \\
& \Mul_{ \calO }( \{ W_{i} \}_{i \in I}, Z) & }$$
is commutative.
\end{itemize}
\end{definition}

\begin{remark}\label{inder}
Every colored operad $\calO$ has an underlying category (which we will also denote by $\calO$)
whose objects are the colors of $\calO$ and whose morphisms are defined by the formula
$\Hom_{\calO}(X, Y) = \Mul_{\calO}( \{ X\}, Y)$. Consequently, we can view a colored operad
as a category with some additional structure: namely, the collections of ``multilinear'' maps
$\Mul_{\calO}( \{ X_i \}_{i \in I}, Y)$. Some authors choose to emphasize this analogy by using the term {\it multicategory} for what we refer to here as a colored operad.
\end{remark}

\begin{variant}\label{col}
We can also define the notion of a {\it simplicial colored operad} by replacing the morphism
sets $\Mul_{\calO}( \{ X_i \}_{i \in I}, Y)$ by simplicial sets in Definition \ref{colop}.
\end{variant}

\begin{remark}
Remark \ref{inder} describes a forgetful functor from the category $\Cat$ of (small) categories
to the category $\MCat$ of (small) colored operads. This functor has a left adjoint: every category
$\calC$ can be regarded as a colored operad by defining
$$ \Mul_{\calC}( \{ X_i \}_{1 \leq i \leq n}, Y) = \begin{cases} \Hom(X_1, Y) & \text{if } $n=1$ \\
\emptyset & \text{otherwise.}  \end{cases}$$
This construction exhibits $\Cat$ as a full subcategory of $\MCat$: namely, the full subcategory spanned by those colored operads $\calO$ for which the sets $\Mul_{\calO}( \{ X_i \}_{i \in I}, Y)$ are empty
unless $I$ has exactly one element.
\end{remark}

\begin{example}\label{singu}
Let $\calC$ be a symmetric monoidal category: that is, a category equipped with a tensor product functor $\otimes: \calC \times \calC \rightarrow \calC$ which is commutative and associative up to coherent isomorphism (see \cite{maclane} for a careful definition). Then we can regard $\calC$ as a colored operad by setting
$$\Mul_{\calC}( \{ X_i \}_{i \in I}, Y) = \Hom_{\calC}( \otimes_{i \in I} X_i, Y).$$
We can recover the symmetric monoidal structure on $\calC$ (up to canonical isomorphism)
via Yoneda's lemma. For example, for every pair of objects $X,Y \in \calC$, the tensor product
$X \otimes Y$ is characterized up to canonical isomorphism by the fact that it corepresents the functor
$Z \mapsto \Mul_{\calC}( \{X, Y\}, Z)$. It follows from this analysis that we can regard the notion of a symmetric monoidal category as a special case of the notion of a colored operad.
\end{example} 

\begin{example}
An {\it operad} is a colored operad $\calO$ having only a single color ${\bf 1}$. For every nonnegative integer $n$ we let $\calO_{n} = \Mul( \{ {\bf 1} \}_{1 \leq i \leq n}, {\bf 1})$. We sometimes refer to
$\calO_n$ as the {\it set of $n$-ary operations of $\calO$}. The structure of $\calO$ as a colored
operad is determined by the sets $\{ \calO_n \}_{n \geq 0}$ together with the actions of the symmetric groups $\Sigma_n$ on $\calO_n$ and the ``substitution maps''
$$ \calO_{m} \times (\prod_{1 \leq i \leq m} \calO_{n_i}) \rightarrow \calO_{n_1 + \cdots + n_m}.$$
\end{example}

Definition \ref{colop} is phrased in a somewhat complicated way because the notion of 
morphism in a colored operad $\calO$ is assymetrical: the domain of a morphism consists of
a finite collection of object, while the codomain consists of only a single object. We can correct this assymetry by repackaging the definition of a colored operad in a different way. First, we recall the
definition of Segal's category $\FinSeg$ of pointed finite sets:

\begin{definition}
For each $n \geq 0$, we let $\nostar{n}$ denote the set
$\{ 1 , 2 , \ldots , n-1 , n\}$ and $\seg{n} = \seg{n}^{\circ}_{\ast} = \{ \ast, 1, \ldots, n \}$ the pointed
set obtained by adjoining a disjoint base point $\ast$ to $\nostar{n}$. We define
a category $\FinSeg$ as follows:
\begin{itemize}
\item[$(1)$] The objects of $\FinSeg$ are the sets $\seg{n}$, where $n \geq 0$.
\item[$(2)$] Given a pair of objects $\seg{m}, \seg{n} \in \FinSeg$, a morphism
from $\seg{m}$ to $\seg{n}$ in $\FinSeg$ is a map $\alpha: \seg{m} \rightarrow \seg{n}$ such that
$\alpha(\ast) = \ast$.
\end{itemize}
\end{definition}

\begin{remark}
The category $\FinSeg$ is equivalent to the category of {\em all} finite sets equipped with a distinguished point $\ast$. We will often invoke this equivalence implicitly, using the following device. Let $I$ be a finite linearly ordered set. Then there is a canonical bijection
$\alpha: I_{\ast} \simeq \seg{n}$, where $n$ is the cardinality of $I$; the bijection
$\alpha$ is determined uniquely by the requirement that the restriction of
$\alpha$ determines an order-preserving bijection of $I$ with $\nostar{n}$. We will generally
identify $I_{\ast}$ with the object $\seg{n} \in \FinSeg$ via this isomorphism.
\end{remark}

\begin{construction}\label{coco}
Let $\calO$ be a colored operad. We define a category $\calO^{\otimes}$ as follows:
\begin{itemize}
\item[$(1)$] The objects of $\calO^{\otimes}$ are finite sequences of colors $X_1, \ldots, X_n \in \calO$.
\item[$(2)$] Given two sequences of objects
$$ X_1, \ldots, X_m \in \calO \quad \quad Y_1, \ldots, Y_n \in \calO,$$
a morphism from $\{ X_i \}_{1 \leq i \leq m}$ to $\{ Y_j \}_{1 \leq j \leq n}$ is given by a map
$\alpha: \seg{m} \rightarrow \seg{n}$ in $\FinSeg$ together with a collection of morphisms
$$\{ \phi_{j} \in \Mul_{\calO}( \{ X_i \}_{i \in \alpha^{-1} \{j\} }, Y_j ) \}_{1 \leq j \leq n}$$ in $\calO$.
\item[$(3)$] Composition of morphisms in $\calO^{\otimes}$ is determined by the composition laws
on $\FinSeg$ and on the colored operad $\calO$.
\end{itemize}
\end{construction}

Let $\calO$ be a colored operad. By construction, the category $\calO^{\otimes}$ comes equipped with a forgetful functor $\pi: \calO^{\otimes} \rightarrow \FinSeg$. Using the functor $\pi$, we can reconstruct
the colored operad $\calO$ up to canonical equivalence. For example, the underlying category
of $\calO$ can be identified with the fiber $\calO^{\otimes}_{ \seg{1} } = \pi^{-1} \{ \seg{1} \}$. 
More generally, suppose that $n \geq 0$. For $1 \leq i \leq n$, let 
$\Colp{i}: \seg{n} \rightarrow \seg{1}$ be as in Notation \ref{coll}, so that
$\Colp{i}$ induces a functor $\Colll{i}{!}: \calO^{\otimes}_{\seg{n}} \rightarrow \calO^{\otimes}_{\seg{1}} \simeq \calO$ and these functors together determine an equivalence of categories
$\calO^{\otimes}_{\seg{n}} \simeq \calO^{n}$. Given a finite sequence of colors $X_1, \ldots, X_n \in \calO$, let $\vec{X}$ denote the corresponding object of $\calO^{\otimes}_{\seg{n}}$ (which is well-defined up to equivalence). For every color $Y \in \calO$, the set $\Mul_{\calO}( \{ X_i \}_{1 \leq i \leq n}, Y)$ can be identified with the set of morphisms $f: \vec{X} \rightarrow Y$ in $\calO^{\otimes}$ such that
$\pi(f): \seg{n} \rightarrow \seg{1}$ satisfies $\pi(f)^{-1} \{\ast \} = \{ \ast \}$. This shows that
$\pi: \calO^{\otimes} \rightarrow \FinSeg$ determines the morphism sets
$\Mul_{\calO}( \{ X_i \}_{1 \leq i \leq n}, Y)$ in the colored operad $\calO$; a more elaborate argument shows that the composition law for morphisms in the colored operad $\calO$ can be recovered from the composition law for morphisms in the category $\calO^{\otimes}$.

The above construction suggests that it is possible to give an alternate version of
Definition \ref{colop}: rather than thinking of a colored operad as a category-like structure $\calO$ equipped with an elaborate notion of morphism, we can think of a colored operad as an
ordinary category $\calO^{\otimes}$ equipped with a forgetful functor $\pi: \calO^{\otimes} \rightarrow \FinSeg$. Of course, we do not want to consider an {\em arbitrary} functor $\pi$: we only want to consider those functors which induce equivalences $\calO^{\otimes}_{\seg{n}} \simeq (\calO^{\otimes}_{\seg{1}})^{n}$, so that the category $\calO = \calO^{\otimes}_{\seg{1}}$ inherits the structure described in Definition \ref{colop}. This is one drawback of the second approach: it requires us to formulate a somewhat complicated-looking assumption on the functor $\pi$. The virtue of the second approach is that it can be phrased entirely in the language of category theory. This allows us to generalize the theory of colored operads to the $\infty$-categorical setting. First, we need to introduce a bit of terminology.

\begin{definition}\label{indef}
We will say that a morphism $f: \seg{m} \rightarrow \seg{n}$ in $\FinSeg$ is {\it inert} if, for each
element $i \in \nostar{n}$, the inverse image $f^{-1} \{i\}$ has exactly one element.
\end{definition}

\begin{remark}
Every inert morphism $f: \seg{m} \rightarrow \seg{n}$ in $\FinSeg$ induces an injective map of
sets $\alpha: \nostar{n} \rightarrow \nostar{m}$, characterized by the formula
$f^{-1} \{i \} = \{ \alpha(i) \}$.
\end{remark}

\begin{definition}\label{colopi}
An {\it $\infty$-operad} is a functor $\calO^{\otimes} \rightarrow \Nerve(\FinSeg)$ between
$\infty$-categories which satisfies the following conditions:
\begin{itemize}
\item[$(1)$] For every inert morphism $f: \seg{m} \rightarrow \seg{n}$ in $\Nerve(\FinSeg)$ and every
object $C \in \calO^{\otimes}_{\seg{m}}$, there exists a $p$-coCartesian morphism $\overline{f}: C \rightarrow C'$ in $\calO^{\otimes}$ lifting $f$. In particular, $f$ induces a functor
$f_{!}: \calO^{\otimes}_{\seg{m}} \rightarrow \calO^{\otimes}_{\seg{n}}$. 
\item[$(2)$] Let $C \in \calO^{\otimes}_{\seg{m}}$ and $C' \in \calO^{\otimes}_{\seg{n}}$ be objects,
let $f: \seg{m} \rightarrow \seg{n}$ be a morphism in $\FinSeg$, and let
$\bHom^{f}_{\calO^{\otimes}}(C,C')$ be the union of those connected components of
$\bHom_{\calO^{\otimes}}(C,C')$ which lie over $f \in \Hom_{\FinSeg}( \seg{m}, \seg{n} )$. 
Choose $p$-coCartesian morphisms $C' \rightarrow C'_{i}$ lying over the inert morphisms
$\Colp{i}: \seg{n} \rightarrow \seg{1}$ for $1 \leq i \leq n$.
Then the induced map 
$$ \bHom^{f}_{\calO^{\otimes}}(C,C') \rightarrow \prod_{1 \leq i \leq n} \bHom^{\Colp{i} \circ f}_{\calO^{\otimes}}(C, C'_i)$$ is a homotopy equivalence.
\item[$(3)$] For every finite collection of objects $C_1, \ldots, C_n \in \calO^{\otimes}_{\seg{1}}$, there
exists an object $C \in \calO^{\otimes}_{\seg{n}}$ and a collection of $p$-coCartesian morphisms
$C \rightarrow C_i$ covering $\Colp{i}: \seg{n} \rightarrow \seg{1}$.
\end{itemize}
\end{definition}

\begin{remark}
Definition \ref{colopi} is really an $\infty$-categorical generalization of the notion of a colored operad, rather than that of an operad. Our choice of terminology is motivated by a desire to avoid awkward language. To obtain an $\infty$-categorical analogue of the notion of an operad, we should consider
instead $\infty$-operads $\calO^{\otimes} \rightarrow \Nerve(\FinSeg)$ equipped with an essentially
surjective functor $\Delta^{0} \rightarrow \calO^{\otimes}_{\seg{1}}$. 
\end{remark}

\begin{remark}
Let $p: \calO^{\otimes} \rightarrow \Nerve(\FinSeg)$ be an $\infty$-operad. We will often abuse terminology by referring to $\calO^{\otimes}$ as an $\infty$-operad (in this case, it is implicitly assumed that we are supplied with a map $p$ satisfying the conditions listed in Definition \ref{colopi}).
We will usually denote the fiber $\calO^{\otimes}_{ \seg{1}} \simeq p^{-1} \{ \seg{1} \}$ by $\calO$. We will sometimes refer to $\calO$ as the {\it underlying $\infty$-category} of $\calO^{\otimes}$.
\end{remark}

\begin{remark}
Let $p: \calO^{\otimes} \rightarrow \Nerve(\FinSeg)$ be an $\infty$-operad. Then $p$ is a categorical fibration. To prove this, we first observe that $p$ is an inner fibration (this follows from Proposition \toposref{skawky} since $\calO^{\otimes}$ is an $\infty$-category and $\Nerve(\FinSeg)$ is the nerve of an ordinary category). In view of Corollary \toposref{gottaput}, it will suffice to show that if
$C \in \calO^{\otimes}_{\seg{m}}$ and $\alpha: \seg{m} \rightarrow \seg{n}$ is an isomorphism
in $\FinSeg$, then we can lift $\alpha$ to an equivalence $\overline{\alpha}: C \rightarrow C'$ in $\calO^{\otimes}$. Since $\alpha$ is inert, we can choose $\overline{\alpha}$ to be $p$-coCartesian; it then follows from Proposition \toposref{universalequiv} that $\overline{\alpha}$ is an equivalence.
\end{remark}

\begin{remark}\label{sil}
Let $p: \calO^{\otimes} \rightarrow \Nerve(\FinSeg)$ be a functor between $\infty$-categories which
satisfies conditions $(1)$ and $(2)$ of Definition \ref{colopi}. Then $(3)$ is equivalent to the following apparently stronger condition:
\begin{itemize}
\item[$(3')$] For each $n \geq 0$, the functors 
$\{ \Colll{i}{!} : \calO^{\otimes}_{ \seg{n}} \rightarrow \calO \}_{1 \leq i \leq n}$ induce an equivalence of $\infty$-categories $$\phi: \calO^{\otimes}_{\seg{n}} \rightarrow
\calO^{n}.$$
\end{itemize}
It follows easily from $(2)$ that $\phi$ is fully faithful, and condition $(3)$ guarantees that $\phi$ is an equivalence.
\end{remark}

\begin{remark}\label{jilkp}
Let $p: \calO^{\otimes} \rightarrow \Nerve(\FinSeg)$ be an $\infty$-operad and $\calO = \calO^{\otimes}_{\seg{1}}$ the underlying $\infty$-category. It follows from Remark \ref{sil} that we have
a canonical equivalence $\calO^{\otimes}_{\seg{n}} \simeq \calO^{n}$. Using this equivalence,
we can identify objects of $\calO^{\otimes}_{\seg{n}}$ with finite sequences
$(X_1, X_2, \ldots, X_n)$ of objects of $\calO$. We will sometimes denote the corresponding
object of $\calO^{\otimes}$ by $X_1 \oplus X_2 \oplus \cdots \oplus X_n$ (this object is well-defined up to equivalence). More generally, given objects $X \in \calO^{\otimes}_{\seg{m}}$ and
$Y \in \calO^{\otimes}_{\seg{n}}$ corresponding to sequences
$(X_1, \ldots, X_m)$ and $(Y_1, \ldots, Y_n)$, we let $X \oplus Y$ denote an object
of $\calO^{\otimes}_{\seg{m+n}}$ corresponding to the sequence
$(X_1, \ldots, X_m, Y_1, \ldots, Y_n)$. We will discuss the operation $\oplus$ more systematically in
\S \ref{monenv}.
\end{remark}

\begin{notation}\label{clink}
Let $\calO^{\otimes}$ be an $\infty$-operad, and suppose we are given a finite sequence
of objects $\{ X_i \}_{1 \leq i \leq n}$ of $\calO$ and another object $Y \in \calO$. We let
$\Mul_{\calO}( \{ X_i \}_{1 \leq i \leq n}, Y)$ denote the union of those components of
$\bHom_{\calO^{\otimes}}( X_1 \oplus \ldots \oplus X_n ,Y)$ which lie over the unique active morphism $\beta: \seg{n} \rightarrow \seg{1}$. We regard $\Mul_{\calO}( \{ X_i \}_{1 \leq i \leq n}, Y)$ as an object in
the homotopy category $\calH$ of spaces, which is well-defined up to canonical isomorphism.
\end{notation}

\begin{remark}
Let $\calO^{\otimes}$ be an $\infty$-operad. Then we should imagine that
$\calO^{\otimes}$ consists of an ordinary $\infty$-category $\calO$ which comes equipped
with a more elaborate notion of morphism supplied by the spaces
$\Mul_{\calO}( \bigdot, \bigdot)$ of Notation \ref{clink}. These morphism spaces are equipped with
a composition law satisfying the axiomatics of Definition \ref{colop} up to coherent homotopy.
A precise formulation of this idea is cumbersome, and Definition \ref{colopi} provides an efficient substitute.
\end{remark}

\begin{example}\label{xe1}
The identity map exhibits $\calO^{\otimes} = \Nerve(\FinSeg)$ as an $\infty$-operad, whose
underlying $\infty$-category $\calO$ is equivalent to $\Delta^0$. We will refer to this
$\infty$-operad as the {\it commutative $\infty$-operad}. When we wish to emphasize its
role as an $\infty$-operad, we will sometimes denote it by $\CommOp$.
\end{example}

\begin{example}\label{xe1.5}
Let $\InjSeg$ denote the subcategory of $\FinSeg$ spanned by all objects together
with those morphisms $f: \seg{m} \rightarrow \seg{n}$ such that $f^{-1} \{i\}$
has at most one element for $1 \leq i \leq n$. The nerve $\Nerve(\InjSeg)$
is an $\infty$-operad, which we will denote by $\OpE{0}$. 
\end{example}

\begin{example}\label{xe2}
We define a category $\CatAss$ as follows:
\begin{itemize}
\item The objects of $\CatAss$ are the objects of $\FinSeg$.
\item Given a pair of objects $\seg{m}, \seg{n} \in \FinSeg$, a morphism from
$\seg{m}$ to $\seg{n}$ in $\CatAss$ consists of a pair $(f, \{ \leq_i \}_{1 \leq i \leq n} )$, where
$f: \seg{m} \rightarrow \seg{n}$ is a morphism in $\FinSeg$ and $\leq_i$ is a linear ordering
on the inverse image $f^{-1} \{i\} \subseteq \nostar{m}$ for $1 \leq i \leq n$.
\item The composition of a pair of morphisms 
$$(f, \{ \leq_i \}_{1 \leq i \leq n} ): \seg{m} \rightarrow \seg{n}
\quad \quad (g, \{ \leq'_j \}_{1 \leq j \leq p} ): \seg{n} \rightarrow \seg{p}$$
is the pair $(g \circ f, \{ \leq''_{j} \}_{1 \leq j \leq p})$, where each
$\leq''_{j}$ is the lexicographical ordering characterized by the property that
for $a,b \in \nostar{m}$ such that $(g \circ f)(a) = (g \circ f)(b) = j$, we have
$a \leq''_{j} b$ if and only if $f(a) \leq'_{j} f(b)$ and $a \leq_{i} b$ if $f(a) = f(b) = i$.  
\end{itemize}
Let $\Ass$ denote the nerve of the category $\CatAss$. The forgetful functor
$\CatAss \rightarrow \FinSeg$ induces a map $\Ass \rightarrow \Nerve(\FinSeg)$ which
exhibits $\Ass$ as an $\infty$-operad; we will refer to $\Ass$ as the {\it associative $\infty$-operad}.
\end{example}

\begin{example}\label{trivop}
Let $\CatTriv$ be the subcategory of $\FinSeg$ spanned by the inert morphisms (together with
all objects of $\FinSeg$), and let $\Triv = \Nerve( \CatTriv)$. Then the inclusion
$\CatTriv \subseteq \FinSeg$ induces a functor $\Triv \rightarrow \Nerve(\FinSeg)$ which
exhibits $\Triv$ as an $\infty$-operad; we will refer to $\Triv$ as the {\it trivial $\infty$-operad}.
\end{example}

\begin{example}\label{xe3}
Each of the above examples is an instance of the following general observation. Let $\calO$ be a colored operad in the sense of Definition \ref{colop}, and let
$\calO^{\otimes}$ be the ordinary category given by Construction \ref{coco}. Then
the forgetful functor $\Nerve( \calO^{\otimes} ) \rightarrow \Nerve(\FinSeg)$ is an $\infty$-operad.
\end{example}

According to Definition \ref{indef}, a morphism $\gamma: \seg{m} \rightarrow \seg{n}$ in
$\FinSeg$ is {\it inert} if exhibits $\seg{n}$ as the quotient of $\seg{m}$ obtained by identifying
a subset of $\nostar{m}$ with the base point $\ast$. In this case, $\gamma^{-1}$ determines an injective map from $\nostar{n}$ to $\nostar{m}$. By design, for any $\infty$-operad
$\calO^{\otimes}$, the morphism $\gamma$ induces a functor $\gamma_{!}: \calO^{\otimes}_{\seg{m}} \rightarrow \calO^{\otimes}_{\seg{n}}$. This functor can be identified with the projection map 
$\calO^{m} \rightarrow \calO^{n}$ determined by $\gamma^{-1}$. Our choice of terminology is intended to emphasize the role of $\gamma_{!}$ as a forgetful functor, which is not really encoding the essential
structure of $\calO^{\otimes}$. We now introduce a class of morphisms which lies at the other extreme:

\begin{definition}
A morphism $f: \seg{m} \rightarrow \seg{n}$ in $\FinSeg$ is {\it active} if 
$f^{-1} \{\ast \} = \{ \ast \}$.
\end{definition}

\begin{remark}\label{singer}
Every morphism $f$ in $\FinSeg$ admits a factorization $f = f' \circ f''$, where $f''$ is inert and $f'$ is active; moreover, this factorization is unique up to (unique) isomorphism. In other words, the collections of inert and active morphisms determine a factorization system on $\FinSeg$.
\end{remark}

The classes of active and inert morphisms determine a factorization system on the category $\FinSeg$
which induces an analogous factorization system on any $\infty$-operad $\calO^{\otimes}$.

\begin{definition}
Let $p: \calO^{\otimes} \rightarrow \Nerve(\FinSeg)$ be an $\infty$-operad. We will say that a morphism
$f$ in $\calO^{\otimes}$ is {\it inert} if $p(f)$ is inert and $f$ is $p$-coCartesian. We will say that a morphism $f$ in $\calO^{\otimes}$ is {\it active} if $p(f)$ is active.
\end{definition}

\begin{proposition}\label{singwell}
Let $\calO^{\otimes}$ be an $\infty$-operad. Then the collections of active and inert morphisms
determine a factorization system on $\calO^{\otimes}$. 
\end{proposition}

Proposition \ref{singwell} is an immediate consequence of Remark \ref{singer} together with the following general result:

\begin{proposition}
Let $p: \calC \rightarrow \calD$ be an inner fibration of $\infty$-categories.
Suppose that $\calD$ admits a factorization system $(S_L, S_R)$ satisfying the following
condition:
\begin{itemize}
\item[$(\ast)$] For every object $C \in \calC$ and every morphism $\alpha: p(C) \rightarrow D$
in $\calD$ which belongs to $S_L$, there exists a $p$-coCartesian morphism
$\overline{\alpha}: C \rightarrow \overline{D}$ lifting $\alpha$.
\end{itemize}
Let $\overline{S}_{L}$ denote the collection of all $p$-coCartesian morphisms $\overline{\alpha}$ in
$\calC$ such that $p( \overline{\alpha} ) \in S_L$, and let $\overline{S}_{R} = p^{-1} S_{R}$.
Then $(\overline{S}_{L}, \overline{S}_{R} )$ is a factorization system on $\calC$.
\end{proposition}

\begin{proof}
We will prove that $( \overline{S}_{L}, \overline{S}_{R})$ satisfies conditions $(1)$ through
$(3)$ of Definition \toposref{spanhun}:
\begin{itemize}
\item[$(1)$] The collections $\overline{S}_{L}$ and $\overline{S}_{R}$ are stable under retracts. This follows from the stability of $S_{L}$ and $S_{R}$ under retracts, together with the observation that the collection of $p$-coCartesian morphisms is stable under retracts.

\item[$(2)$] Every morphism in $\overline{S}_L$ is left orthogonal to every morphism
in $\overline{S}_{R}$. To prove this, let $\overline{\alpha}: \overline{A} \rightarrow \overline{B}$ belong to $\overline{S}_{L}$ and $\overline{\beta}: \overline{X} \rightarrow \overline{Y}$ belong to $\overline{S}_R$. Let $\alpha: A \rightarrow B$ and $\beta: X \rightarrow Y$ denote the images
of $\overline{\alpha}$ and $\overline{\beta}$ under the functor $p$.
We wish to prove that the space $\bHom_{ \calC_{\overline{A}/ \, /\overline{Y}}}(\overline{B},\overline{X})$ is contractible. Using the fact that $\overline{\alpha}$ is $p$-coCartesian, we deduce that the map
$\bHom_{ \calC_{\overline{A}/ \, /\overline{Y}}}(\overline{B},\overline{X}) \rightarrow
\bHom_{ \calD_{A/ \, /Y} }( B, X)$ is a trivial Kan fibration. The desired result now follows from the fact that 
$\alpha \in S_L$ is left orthogonal to $\beta \in S_{R}$.

\item[$(3)$] Every morphism $\overline{\alpha}: \overline{X} \rightarrow \overline{Z}$ admits
a factorization $\overline{\alpha} = \overline{\beta} \circ \overline{\gamma}$, where
$\overline{\gamma} \in \overline{S}_L$ and $\overline{\beta} \in \overline{S}_R$.
To prove this, let $\alpha: X \rightarrow Z$ denote the image of $\alpha$ under $p$.
Using the fact that $(S_L, S_R)$ is a factorization system on $\calD$, we deduce that
$\alpha$ fits into a commutative diagram
$$ \xymatrix{ & Y \ar[dr]^{\beta} & \\
X \ar[ur]^{\gamma} \ar[rr]^{\alpha} & & Z }$$
where $\gamma \in S_L$ and $\beta \in S_{R}$. Using assumption $(\ast)$, we can lift
$\gamma$ to a $p$-coCartesian morphism $\overline{\gamma} \in \overline{S}_{R}$.
Using the fact that $\overline{\gamma}$ is $p$-coCartesian, we can lift the above diagram to a commutative triangle
$$ \xymatrix{ & \overline{Y} \ar[dr]^{\overline{\beta}} & \\
\overline{X} \ar[ur]^{\overline{\gamma}} \ar[rr]^{\overline{\alpha}} & & \overline{Z} }$$
in $\calC$ having the desired properties.

\end{itemize}
\end{proof}

We conclude this section with a simple observation that will prove useful throughout this paper:

\begin{remark}\label{talee}
Let $p: \calO^{\otimes} \rightarrow \Nerve(\FinSeg)$ be an $\infty$-operad, and suppose
we are given a collection of inert morphisms $\{ f_{i}: X \rightarrow X_i \}_{1 \leq i \leq m}$ in $\calO^{\otimes}$
covering maps $\seg{n} \rightarrow \seg{n_i}$ in $\FinSeg$ which induce a bijection
$\coprod_{1 \leq i \leq m} \nostar{n_i} \rightarrow \nostar{n}$. These morphisms determine a
$p$-limit diagram $q: (\nostar{m})^{\triangleleft} \rightarrow \calO^{\otimes}$. To prove this,
choose inert morphisms $g_{i,j}: X_{i} \rightarrow X_{i,j}$ covering the maps
$\Colp{j}: \seg{n_i} \rightarrow \seg{1}$ for $1 \leq j \leq n_i$. Using the
fact that $\calO^{\otimes}$ is an $\infty$-category, we obtain a diagram
$\overline{q}: ( \coprod_{1 \leq i \leq m} (\nostar{n_i})^{\triangleleft})^{\triangleleft} 
\rightarrow \calO^{\otimes}$. Since the inclusion $\nostar{m} \subseteq
\coprod_{1 \leq i \leq m} \nostar{n_i}^{\triangleright}$ is cofinal, it will suffice to show
that $\overline{q}$ is a $p$-limit diagram. Using the assumption that $\calO^{\otimes}$
is an $\infty$-operad, we deduce that $\overline{q} | ( \coprod_{1 \leq i \leq m} (\nostar{n_i})^{\triangleleft})$ is a $p$-right Kan extension of $\overline{q} | \nostar{n}$. According to Lemma \toposref{kan0}, it will suffice to show that $\overline{q} | \nostar{n}^{\triangleleft}$ is a $p$-limit diagram, which again follows from the assumption that $\calO^{\otimes}$ is an $\infty$-operad.
\end{remark}

\subsection{Fibrations of $\infty$-Operads}\label{comm1.2}

In order to make effective use of the theory of $\infty$-operads, we must understand them not only in isolation but also in relation to one another. To this end, we will introduce the notion of
{\it $\infty$-operad fibration} (Definition \ref{quas}). We begin with the more general notion of
an {\it $\infty$-operad map}, which we will study in more detail in \S \ref{algobj}.

\begin{definition}\label{kij}
Let $\calO^{\otimes}$ and ${\calO'}^{\otimes}$ be $\infty$-operads. An
{\it $\infty$-operad map} from $\calO^{\otimes}$ to ${\calO'}^{\otimes}$ is a map of simplicial sets
$f: \calO^{\otimes} \rightarrow {\calO'}^{\otimes}$ satisfying the following conditions:
\begin{itemize}
\item[$(1)$] The diagram
$$ \xymatrix{ \calO^{\otimes} \ar[rr]^{f} \ar[dr] & & {\calO'}^{\otimes} \ar[dl] \\
& \Nerve(\FinSeg) & }$$
commutes.
\item[$(2)$] The functor $f$ carries inert morphisms in $\calO^{\otimes}$ to inert morphisms
in ${\calO'}^{\otimes}$.
\end{itemize}
We let $\Fun^{\lax}( \calO^{\otimes}, {\calO'}^{\otimes})$ denote the full subcategory of
$\Fun_{ \Nerve(\FinSeg)}( \calO^{\otimes}, {\calO'}^{\otimes})$ spanned by the $\infty$-operad
maps.
\end{definition}

\begin{remark}\label{casper}
Let $\calO^{\otimes}$ and ${\calO'}^{\otimes}$ be $\infty$-operads, and let
$F: \calO^{\otimes} \rightarrow {\calO'}^{\otimes}$ be a functor satisfying condition
$(1)$ of Definition \ref{kij}. Then $F$ preserves all inert morphisms if and only if
it preserves inert morphisms in $\calO^{\otimes}$ lying over the maps
$\Colp{i}: \seg{n} \rightarrow \seg{1}$. For suppose that
$X \rightarrow Y$ is an inert morphism in $\calO^{\otimes}$ lying over
$\beta: \seg{m} \rightarrow \seg{n}$; we wish to prove that the induced map
$\beta_{!} F(X) \rightarrow F(Y)$ is an equivalence. Since ${\calO'}^{\otimes}$ is an $\infty$-operad, it suffices to show that the induced map $\Colll{i}{!} : \beta_{!} F(X) \rightarrow \Colll{i}{!}  F(Y)$ is an equivalence in $\calO'$ for $1 \leq i \leq n$. Our hypothesis allows us to identify this with the morphism $F( \Colll{i}{!}  \beta_{!} X) \rightarrow
F( \Colll{i}{!}  Y)$, which is an equivalence since $Y \simeq \beta_{!} X$.
\end{remark}

\begin{definition}\label{quas}
We will say that a map of $\infty$-operads $q: \calC^{\otimes} \rightarrow \calO^{\otimes}$ is
a {\it fibration of $\infty$-operads} if $q$ is a categorical fibration.
\end{definition}

\begin{remark}\label{cotallee}
Let $p: \calC^{\otimes} \rightarrow \calO^{\otimes}$ be a fibration of $\infty$-operads, and suppose we are given a collection of inert morphisms $\{ f_{i}: X \rightarrow X_i \}_{1 \leq i \leq m}$ in $\calC^{\otimes}$
covering maps $\seg{n} \rightarrow \seg{n_i}$ in $\FinSeg$ which induce a bijection
$\coprod_{1 \leq i \leq m} \nostar{n_i} \rightarrow \nostar{n}$. Then these morphisms determine a
$p$-limit diagram $q: \nostar{m}^{\triangleleft} \rightarrow \calC^{\otimes}$. This follows
from Remark \ref{talee} and Proposition \toposref{basrel}.
\end{remark}

The following result describes an important special class of $\infty$-operad fibrations:

\begin{proposition}\label{saber}
Let $\calO^{\otimes}$ be an $\infty$-operad, and let $p: \calC^{\otimes} \rightarrow \calO^{\otimes}$
be a coCartesian fibration. The following conditions are equivalent:
\begin{itemize}
\item[$(a)$] The composite map $q: \calC^{\otimes} \rightarrow \calO^{\otimes} \rightarrow
\Nerve(\FinSeg)$ exhibits $\calC^{\otimes}$ as an $\infty$-operad.
\item[$(b)$] For every object $T \simeq T_1 \oplus \cdots \oplus T_n \in \calO^{\otimes}_{\seg{n}}$, the
inert morphisms $T \rightarrow T_i$ induce an equivalence of
$\infty$-categories 
$\calC^{\otimes}_{T} \rightarrow \prod_{1 \leq i \leq n} \calC^{\otimes}_{T_i}$.
\end{itemize}
\end{proposition}

\begin{proof}
Suppose that $(a)$ is satisfied. We first claim that $p$ preserves inert morphisms.
To prove this, choose an inert morphism $f: C \rightarrow C'$ in $\calC^{\otimes}$.
Let $g: p(C) \rightarrow X$ be an inert morphism in $\calO^{\otimes}$ lifting
$q(f)$, and let $\overline{g}: C \rightarrow \overline{X}$ be a $p$-coCartesian lift
of $g$. It follows from Proposition \toposref{stuch} that $\overline{g}$ is
a $q$-coCartesian lift of $q(f)$ and therefore equivalent to $f$. We conclude that
$p(f)$ is equivalent to $p( \overline{g} ) = g$ and is therefore inert.

The above argument guarantees that for each $n \geq 0$, the maps
$\{ \Colp{i} \}_{1 \leq i \leq n}$ determine a homotopy commutative diagram
$$ \xymatrix{ \calC^{\otimes}_{\seg{n}} \ar[r] \ar[d] & \calO^{\otimes}_{\seg{n}} \ar[d] \\
\calC^{n} \ar[r] & \calO^{n} }$$
and the assumption that $\calC^{\otimes}$ and $\calO^{\otimes}$ are $\infty$-operads
guarantees that the vertical maps are categorical equivalences. Let $T$ be an object
of $\calO^{\otimes}_{\seg{n}}$. Passing to the homotopy fibers over the vertices $T$ and $( T_i )_{1 \leq i \leq n}$ (which are equivalent to the actual fibers by virtue of Corollary \toposref{basety}), we deduce that the canonical map $\calC^{\otimes}_{T} \rightarrow \prod_{1 \leq i \leq n} \calC^{\otimes}_{T_i}$ is an equivalence, which proves $(b)$.

Now suppose that $(b)$ is satisfied. We will prove that the functor
$q: \calC^{\otimes} \rightarrow \Nerve(\FinSeg)$ satisfies the conditions of
Definition \ref{colopi}. To prove $(1)$, consider an object $C \in \calC^{\otimes}$
and an inert morphism $\alpha: q(C) \rightarrow \seg{n}$ in $\FinSeg$. 
Since $\calO^{\otimes} \rightarrow \Nerve(\FinSeg)$ is an $\infty$-operad, there
exists an inert morphism $\widetilde{\alpha}: p(C) \rightarrow X$ in $\calO^{\otimes}$ lying
over $\alpha$. Since $p$ is a coCartesian fibration, we can lift $\widetilde{\alpha}$
to a $p$-coCartesian morphism $\overline{\alpha}: C \rightarrow \overline{X}$ in
$\calC^{\otimes}$. Proposition \toposref{stuch} implies that $\overline{\alpha}$ is
$p$-coCartesian, which proves $(1)$. 

Now let $C \in \calC^{\otimes}_{\seg{m}}$, $C' \in \calC^{\otimes}_{\seg{n}}$, and
$f: \seg{m} \rightarrow \seg{n}$ be as in condition $(2)$ of Definition \ref{colopi}.
Set $T= p(C)$, set $T' = p(C')$, choose inert morphisms $g_i: T' \rightarrow T'_i$ lying over
$\Colp{i}$ for $1 \leq i \leq n$, and choose $p$-coCartesian morphisms
$\overline{g}_i: C' \rightarrow C'_i$ lying over $g_i$ for $1 \leq i \leq n$.
Let $f_i = \Colp{i} \circ f$ for $1 \leq i \leq n$. We have a homotopy coherent diagram
$$ \xymatrix{ \bHom^{f}_{\calC^{\otimes}}( C, C') \ar[r] \ar[d] & \prod_{1 \leq i \leq n}
\bHom^{f_i}_{ \calC^{\otimes}}( C, C'_i) \ar[d] \\
\bHom^{f}_{\calO^{\otimes}}(T, T') \ar[r] & \prod_{1 \leq i \leq n} \bHom^{f_i}_{\calO^{\otimes}}(T, T'_i). }$$
Since $\calO^{\otimes}$ is an $\infty$-operad, the bottom horizontal map is a homotopy equivalence.
Consequently, to prove that the top map is a homotopy equivalence, it will suffice to show that it induces a homotopy equivalence after passing to the homotopy fiber over any 
$h: T \rightarrow T'$ lying over $f$. Let $D = h_{!} C$ and let $D_i = (g_i)_{!} D$ for
$1 \leq i \leq n$. Using the assumption that $p$ is a coCartesian fibration and
Proposition \toposref{compspaces}, we see that the map of homotopy fibers can be identified with
$$ \bHom_{ \calC^{\otimes}_{T'}}(D, C) \rightarrow \prod_{1 \leq i \leq n} \bHom_{\calC^{\otimes}_{T_i}}( D_i, C_i),$$
which is a homotopy equivalence by virtue of assumption $(b)$.

We now prove $(3)$. Fix a sequence of objects $\{ C_i \in \calC^{\otimes}_{\seg{1}} \}_{1 \leq i \leq n}$, and set $T_i = p(C_i)$. Since $\calO^{\otimes}$ is an $\infty$-operad, we can choose an object
$T \in \calO^{\otimes}$ equipped with inert morphisms $f_i: T \rightarrow T_i$ lifting
$\Colp{i}$ for $1 \leq i \leq n$. Invoking $(b)$, we conclude that there exists
an object $C$ together with $p$-coCartesian morphisms $\overline{f}_i: C \rightarrow C_i$
lifting $f_i$. It follows from Proposition \toposref{stuch} that each $\overline{f}_i$ is $q$-coCartesian, so that condition $(3)$ is satisfied and the proof is complete.
\end{proof}

\begin{definition}\label{omon}
Let $\calO^{\otimes}$ be an $\infty$-operad.
We will say that a map $p: \calC^{\otimes} \rightarrow \calO^{\otimes}$ is
a {\it coCartesian fibration of $\infty$-operads} if it satisfies the hypotheses
of Proposition \ref{saber}. In this case, we might also say that {\it $p$ exhibits
$\calC^{\otimes}$ as a $\calO$-monoidal $\infty$-category}.
\end{definition}

\begin{remark}
Let $\calO^{\otimes}$ be an $\infty$-operad and let $p: \calC^{\otimes} \rightarrow \calO^{\otimes}$
be a coCartesian fibration of $\infty$-operads. Then $p$ is a map of $\infty$-operads: this follows
from the first step of the proof of Proposition \ref{saber}. Combining this observation with
Proposition \toposref{funkyfibcatfib}, we conclude that $p$ is a fibration of $\infty$-operads in the sense of Definition \ref{quas}.
\end{remark}

\begin{remark}
In the situation of Definition \ref{omon}, we will generally denote the fiber product
$\calC^{\otimes} \times_{ \calO^{\otimes} } \calO$ by $\calC$, and abuse terminology by saying that
{\it $\calC$ is an $\calO$-monoidal $\infty$-category}. 
\end{remark}

\begin{remark}\label{swinn}
Let $\calO^{\otimes}$ be an $\infty$-operad and let $p: \calC^{\otimes} \rightarrow \calO^{\otimes}$ be a coCartesian fibration of $\infty$-operads. For every object $X \in \calO^{\otimes}$, we let $\calC^{\otimes}_{X}$ denote the inverse
image of $T$ under $p$. If $X \in \calO$, we will also denote this $\infty$-category by $\calC_{X}$.
Note that if $X \in \calO^{\otimes}_{\seg{n}}$ corresponds to a sequence of objects
$\{ X_i \}_{1 \leq i \leq n}$ in $\calO$, then we have a canonical equivalence
$\calC^{\otimes}_{X} \simeq \prod_{1 \leq i \leq n} \calC_{X_i}$.

Given a morphism $f \in \Mul_{\calO}( \{ X_i \}_{1 \leq i \leq n}, Y)$ in $\calO$, the 
coCartesian fibration $p$ determines a functor
$$ \prod_{1 \leq i \leq n} \calC_{X_i} \simeq \calC^{\otimes}_{X} \rightarrow
\calC_{Y}$$
which is well-defined up to equivalence; we will sometimes denote this functor by
$\otimes_{f}$.
\end{remark}

\begin{remark}
Let $\calO^{\otimes}$ be an $\infty$-operad, and let $\calC^{\otimes} \rightarrow \calO^{\otimes}$ be
a $\calO$-monoidal $\infty$-category. Then the underlying map $\calC \rightarrow \calO$ is
a coCartesian fibration of $\infty$-categories, which is classified by a functor
$\chi: \calO \rightarrow \Cat_{\infty}$. In other words, we can think of a $\calO$-monoidal
$\infty$-category as assigning to each color $X \in \calO$ an $\infty$-category
$\chi(X)$ of {\it $X$-colored objects}. We will later see that $\chi$ can be extended to a map
of $\infty$-operads, and that this map of $\infty$-operads determines $\calC^{\otimes}$ up to equivalence (Example \ref{interplace} and Proposition \ref{ungbatt}).
\end{remark}

\begin{definition}
A {\it symmetric monoidal $\infty$-category} is an
$\infty$-category $\calC^{\otimes}$ equipped with a coCartesian fibration
of $\infty$-operads $p: \calC^{\otimes} \rightarrow \Nerve(\FinSeg)$.
\end{definition}

\begin{remark}
In other words, a symmetric monoidal $\infty$-category is a coCartesian fibration
$p: \calC^{\otimes} \rightarrow \Nerve(\FinSeg)$ which induces equivalences
of $\infty$-categories $\calC^{\otimes}_{\seg{n}} \simeq \calC^{n}$ for each $n \geq 0$, where
$\calC$ denotes the $\infty$-category $\calC^{\otimes}_{\seg{1}}$. 
\end{remark}

\begin{remark}
Let $\calC^{\otimes} \rightarrow \Nerve(\FinSeg)$ be a symmetric monoidal $\infty$-category. We will refer to the fiber $\calC^{\otimes}_{\seg{1}}$ as the {\it underlying $\infty$-category} of $\calC^{\otimes}$; it will often be denoted by $\calC$. We will sometimes abuse terminology by referring to $\calC$ as a symmetric monoidal $\infty$-category, or say that $\calC^{\otimes}$ determines a symmetric monoidal structure on $\calC$. 

Using the constructions of Remark \ref{swinn}, we see that the active morphisms
$\alpha: \seg{0} \rightarrow \seg{1}$ and $\beta: \seg{2} \rightarrow \seg{1}$ determine
functors 
$$\Delta^0 \rightarrow \calC \quad \quad \calC \times \calC \rightarrow \calC,$$
which are well-defined up to a contractible space of choice. The first of these functors
determines an object of $\calC$ which we will denote by ${\bf 1}$ (or sometimes
${\bf 1}_{\calC}$ if we wish to emphasize the dependence on $\calC$) and refer to 
as the {\it unit object} of $\calC$. 

It is not difficult to verify that the unit object ${\bf 1} \in \calC$ and the tensor product $\otimes$ on
$\calC$ satisfy all of the usual axioms for a symmetric monoidal category up to homotopy.
In particular, these operations endow the homotopy category $\h{\calC}$ with a symmetric monoidal structure.
\end{remark}

\begin{example}\label{spudman}
Let $\calC$ be a symmetric monoidal category (see \cite{maclane}), so that we can
regard $\calC$ as a colored operad as in Example \ref{singu}. Then the $\infty$-operad
$\Nerve( \calC^{\otimes})$ of Example \ref{xe3} is a symmetric monoidal $\infty$-category.
(whose underlying $\infty$-category is the nerve $\Nerve(\calC)$).
We will discuss some generalizations of this example in \S \ref{comm9}.
\end{example}

We conclude this section with a useful criterion for detecting $\infty$-operad fibrations:

\begin{proposition}
Let $q: \calC^{\otimes} \rightarrow \calO^{\otimes}$ be a map of $\infty$-operads which is an
inner fibration. The following conditions are equivalent:
\begin{itemize}
\item[$(1)$] The map $q$ is a fibration of $\infty$-operads.
\item[$(2)$] For every object $C \in \calC^{\otimes}$ and every inert morphism
$f: q(C) \rightarrow X$ in $\calO^{\otimes}$, there exists an inert morphism
$\overline{f}: C \rightarrow \overline{X}$ in $\calC^{\otimes}$ such that $f = q( \overline{f} )$. 
\end{itemize}
Moreover, if these conditions are satisfied, then the inert morphisms of $\calC^{\otimes}$ are precisely
the $q$-coCartesian morphisms in $\calC^{\otimes}$ whose image in $\calO^{\otimes}$ is inert.
\end{proposition}

\begin{proof}
According to Corollary \toposref{gottaput}, condition $(1)$ is satisfied if and only if condition
$(2)$ is satisfied whenever $f$ is an equivalence; this proves that $(2) \Rightarrow (1)$.
For the reverse implication, suppose that $q$ is a categorical fibration and let
$f: q(C) \rightarrow X$ be as in $(2)$. Let $f_0: \seg{m} \rightarrow \seg{n}$ denote
the image of $f$ in $\Nerve(\FinSeg)$, and let $\overline{f}': C \rightarrow \overline{X}'$ be an
inert morphism in $\calC^{\otimes}$ lifting $f_0$. Since $f$ and $q( \overline{f}')$ are both
inert lifts of the morphism $f_0$, we can find a $2$-simplex
$$ \xymatrix{ & X' \ar[dr]^{g} & \\
q(C) \ar[ur]^{q( \overline{f}')} \ar[rr]^{f} & & X }$$
in $\calO^{\otimes}$ where $g$ is an equivalence. Since $q$ is a categorical fibration,
Corollary \toposref{gottaput} guarantees that we can choose an equivalence
$\overline{g}: \overline{X}' \rightarrow \overline{X}$ in $\calC^{\otimes}$ lifting $g$.
Using the fact that $q$ is an inner fibration, we can lift the above diagram
to a $2$-simplex
$$ \xymatrix{ & \overline{X}' \ar[dr]^{ \overline{g} } & \\
C \ar[ur]^{ \overline{f}' } \ar[rr]^{ \overline{f} } & & \overline{X} }$$
in $\calC^{\otimes}$. Since $\overline{f}$ is the composition of an inert morphism
with an equivalence, it is an inert morphism which lifts $f$, as desired. This completes the proof of the implication $(1) \Rightarrow (2)$; the final assertion follows immediately from Proposition \toposref{stuch} (and the assumption that $q$ preserves inert morphisms).
\end{proof}

\subsection{Algebra Objects}\label{algobj}

In this section, we will undertake a more systematic study of maps between
$\infty$-operads.

\begin{definition}\label{susqat}
Let $p: \calC^{\otimes} \rightarrow \calO^{\otimes}$ be a fibration of $\infty$-operads. A
{\it $\calO$-algebra object of $\calC$} is a map of $\infty$-operads 
$A: \calO^{\otimes} \rightarrow \calC^{\otimes}$ such that $p \circ A = \id_{ \calO^{\otimes}}$.
We let $\Alg_{\calO}(\calC)$ denote the full subcategory of 
$\Fun_{ \calO^{\otimes} }( \calO^{\otimes}, \calC^{\otimes} )$ spanned by the
$\calO$-algebra objects of $\calC$. In the special case where $\calO^{\otimes} = \Nerve(\FinSeg)$, we will denote
$\Alg_{\calO}(\calC)$ by $\CAlg(\calC)$, which we refer to as the {\it $\infty$-category of
commutative algebra objects of $\calC$}.

More generally, suppose that we are given a fibration of $\infty$-operads
$p: \calC^{\otimes} \rightarrow \calO^{\otimes}$ as above, and let 
$f: {\calO'}^{\otimes} \rightarrow \calO^{\otimes}$ be a map of $\infty$-operads.
A {\it $\calO'$-algebra object} of $\calC$ is an $\infty$-operad map
$A: {\calO'}^{\otimes} \rightarrow  \calC^{\otimes}$ such that
$p \circ A = f$. We let $\Alg_{\calO'}(\calC)$ denote the full subcategory of
$\Fun_{ \calO^{\otimes}}( {\calO'}^{\otimes} , \calC^{\otimes})$ spanned by 
the $\calO'$-algebra objects of $\calC$.
\end{definition}

\begin{remark}
The notation of Definition \ref{susqat} is somewhat abusive: note that
$\Alg_{\calO'}(\calC)$ depends strongly on the underlying $\infty$-operad $\calO^{\otimes}$, which
we generally assume to be clear from context.
\end{remark}

\begin{remark}
When $\calO^{\otimes} = \Nerve(\FinSeg)$, the $\infty$-category $\Alg_{\calO'}(\calC)$ coincides
with the $\infty$-category $\Fun^{\lax}( {\calO'}^{\otimes}, \calC^{\otimes})$ of Definition \ref{kij}. 
\end{remark}

\begin{example}
Let $\calC$ be a symmetric monoidal category, and let us regard $\Nerve(\calC)$ as a symmetric monoidal $\infty$-category via the construction described in Example \ref{spudman}. Then
$\CAlg( \Nerve(\calC) )$ can be identified with the nerve of the category of {\it commutative algebra objects} of $\calC$: that is, objects $A \in \calC$ equipped with a unit map
$1_{\calC} \rightarrow A$ and a multiplication $A \otimes A \rightarrow A$ which is commutative, unital, and associative.
\end{example}

\begin{remark}\label{funtime}
Let $\calO^{\otimes}$ be an $\infty$-operad, let $\calC^{\otimes} \rightarrow \calO^{\otimes}$ be a $\calO$-monoidal $\infty$-category, and let $K$ be an arbitrary simplicial set. Then
the induced map $$\Fun( K, \calC^{\otimes}) \times_{ \Fun( K, \calO^{\otimes})} \calO^{\otimes} \rightarrow \calO^{\otimes}$$
is another $\calO$-monoidal $\infty$-category $\calD^{\otimes}$. For every object
$X \in \calO$, we have a canonical isomorphism $\calD_{X} \simeq \Fun(K, \calC_{X})$, and the operations $\otimes_{f}$ of
Remark \ref{swinn} are computed pointwise. We have a canonical isomorphism
$$ \Alg_{\calO}( \calD) \simeq \Fun(K, \Alg_{\calO}( \calC) ).$$
\end{remark}

\begin{example}\label{algtriv}
Let $q: \calO^{\otimes} \rightarrow \Nerve(\FinSeg)$ be an $\infty$-operad such that the image
of $q$ is contained in $\Triv \subseteq \Nerve(\FinSeg)$, and let
$p: \calC^{\otimes} \rightarrow \calO^{\otimes}$ be a $\calO$-monoidal $\infty$-category.
We observe that a section $A: \calO^{\otimes} \rightarrow \calC^{\otimes}$ of $p$
is a $\calO$-algebra if and only if $A$ is a $p$-right Kan extension of its restriction
to $\calO \subseteq \calO^{\otimes}$. Using Proposition \toposref{lklk}, we deduce that the restriction map $\Alg_{\calO}(\calC) \rightarrow \bHom_{ \calO}( \calO, \calC)$ is a trivial Kan fibration.
In particular, taking $\calO^{\otimes}$ to be the trivial $\infty$-operad $\Triv$, we obtain
a trivial Kan fibration $\Alg_{\Triv}(\calC) \rightarrow \calC$.
\end{example}

\begin{remark}\label{coalgtriv}
Let $\calO^{\otimes}$ be an $\infty$-operad. Since $\Triv$ can be identified with
a simplicial subset of $\Nerve( \FinSeg)$, an $\infty$-operad map
$p: \calO^{\otimes} \rightarrow \Triv$ is unique if it exists; the existence of $p$ is equivalent to the requirement that every morphism in $\calO^{\otimes}$ is inert. If this condition is satisfied,
then $p$ is automatically a coCartesian fibration and exhibits $\calO^{\otimes}$ as
a $\Triv$-monoidal $\infty$-category.

Now let $\calO^{\otimes}$ be a general $\infty$-operad. The $\infty$-category
$\Alg_{\Triv}(\calO)$ can be identified with the collection of
$\Triv$-algebras in the $\Triv$-monoidal $\infty$-category
$\calO^{\otimes} \times_{ \Nerve(\FinSeg)} \Triv$. By virtue of Example
\ref{algtriv}, we deduce that evaluation at $\seg{1}$ induces a trivial Kan fibration
$\Alg_{\Triv}(\calO) \rightarrow \calO.$
\end{remark}

The following variant of Example \ref{algtriv} can be used to describe algebras over
the $\infty$-operad $\OpE{0}$ of Example \ref{xe1.5}:

\begin{proposition}\label{ezalg}
Let $p: \calO^{\otimes} \rightarrow \OpE{0}$ be a fibration of $\infty$-operads, and 
consider the map $\Delta^1 \rightarrow \OpE{0}$ corresponding to the
unique morphism $\alpha: \seg{0} \rightarrow \seg{1}$ of $\InjSeg \subseteq \FinSeg$.
Then the restriction functor
$$ \theta: \Alg_{ \OpE{0}}( \calO) \rightarrow \Fun_{ \OpE{0}}( \Delta^1, \calO^{\otimes})$$
is a trivial Kan fibration.
\end{proposition}

\begin{remark}
Suppose that $\calC^{\otimes}$ is a symmetric monoidal $\infty$-category with
unit object ${\bf 1}$. Proposition \ref{ezalg} implies that we can identify $\Alg_{\OpE{0}}( \calC)$
with the $\infty$-category $\calC_{ {\bf 1}/ }$ consisting of maps ${\bf 1} \rightarrow A$ in
$\calC$. In other words, an $\OpE{0}$-algebra object of $\calC$ is an object $A \in \calC$ equipped with a unit map ${\bf 1} \rightarrow A$, but no other additional structures.
\end{remark}

\begin{proof}
Since $p$ is a categorical fibration, the map $\theta$ is also a categorical fibration.
The morphism $\alpha$ determines a functor $s: [1] \rightarrow \InjSeg$.
Let $\calJ$ denote the categorical mapping cylinder of $s$. More precisely,
we define the category $\calJ$ as follows:
\begin{itemize}
\item An object of $\calJ$ is either an object $\seg{n} \in \InjSeg$ or
an object $i \in [1]$.
\item Morphisms in $\calJ$ are described by the formulas
$$ \Hom_{\calJ}( \seg{m}, \seg{n}) = \Hom_{\InjSeg}(\seg{m}, \seg{n}) \quad \quad
\Hom_{\calJ}( \seg{m}, i) = \Hom_{ \InjSeg}( \seg{m}, s(i))$$
$$ \Hom_{\calJ}( i,j) = \Hom_{ [1] }(i,j) \quad \quad \quad
\Hom_{\calJ}( i, \seg{m} ) = \emptyset.$$ 
\end{itemize}
Note that there is a canonical retraction $r$ of $\Nerve(\calJ)$ onto the simplicial
subset $\OpE{0} = \Nerve( \InjSeg)$. We let $\calC$ denote the full subcategory of
$\Fun_{ \OpE{0}}( \Nerve(\calJ), \calO^{\otimes})$ spanned by those functors
$F$ satisfying the following conditions:
\begin{itemize}
\item[$(a)$] For $i \in [1]$, the morphism $F( s(i) ) \rightarrow F(i)$ is an equivalence
in $\calO^{\otimes}$.
\item[$(b)$] The restriction $F | \OpE{0}$ belongs to $\Alg_{ \OpE{0}}( \calC)$. 
\end{itemize}
We observe that a functor $F$ satisfies condition $(\ast)$ if and only if
it is a left Kan extension of $F| \OpE{0}$. Using Proposition \toposref{lklk}, we deduce
that the restriction map $\phi: \calC \rightarrow \Alg_{ \OpE{0}}(\calC)$ is a trivial Kan fibration.
Let $\theta': \calC \rightarrow \Fun_{ \OpE{0}}( \Delta^1, \calO^{\otimes})$ be the map given
by restriction to $\Delta^1 \simeq \Nerve( [1] )$. The map $\theta$ can be written as a composition
$$ \Alg_{ \OpE{0}}(\calO) \rightarrow \calC \stackrel{\theta'}{\rightarrow}  \Fun_{ \OpE{0}}( \Delta^1, \calO^{\otimes}),$$
where the first map is the section of $\phi$ given by composition with $r$ (and therefore a categorical equivalence). By a two-out-of-three argument, we are reduced to proving that $\theta'$ is a categorical equivalence.

We will show that $\theta'$ is a trivial Kan fibration. By virtue of Proposition \toposref{lklk}, it will suffice to prove the following:
\begin{itemize}
\item[$(i)$] A functor $\Fun_{ \OpE{0}}( \Nerve(\calJ), \calO^{\otimes})$ belongs to
$\calC$ if and only if $F$ is a $p$-right Kan extension of $F | \Delta^1$.
\item[$(ii)$] Every functor $F_0 \in \Fun_{ \OpE{0} }( \Delta^1, \calO^{\otimes})$ admits
a $p$-right Kan extension $F \in \Fun_{ \OpE{0}}( \Nerve(\calJ), \calO^{\otimes})$.
\end{itemize}

To see this, fix an object $F \in \Fun_{ \OpE{0}}( \Nerve(\calJ), \calO^{\otimes})$
and consider an object $\seg{n} \in \InjSeg$. Let $\calI = \calJ_{ \seg{n} / } \times_{ \calJ} [1]$,
let $\calI_0$ denote the full subcategory of $\calI$ obtained by omitting the map
$\seg{n} \rightarrow 1$ corresponding to the nonsurjective map
$\seg{n} \rightarrow \seg{1}$ in $\InjSeg$, and let
$\calI_{1}$ denote the full subcategory of $\calI_0$ obtained by omitting the unique map
$\seg{n} \rightarrow 0$. We observe that the inclusion $\Nerve( \calI_0)^{op} \subseteq
\Nerve(\calI)^{op}$ is cofinal and that the restriction
$F | \Nerve( \calI_0)$ is a $p$-right Kan extension of $F| \Nerve(\calI_1)$. Using Lemma
\toposref{kan0}, we deduce that $F$ is a $p$-right Kan extension of $F| \Delta^1$
at $\seg{n}$ if and only if the induced map
$$ \Nerve( \calI_1)^{\triangleleft} \rightarrow \Nerve(\calJ) \stackrel{F}{\rightarrow} \calO^{\otimes}$$
is a $p$-limit diagram. Combining this with the observation that $\calI_1$ is a discrete
category (whose objects can be identified with the elements of $\nostar{n}$), we obtain
the following version of $(i)$:

\begin{itemize}
\item[$(i')$] A functor $F \in \Fun_{ \OpE{0}}( \Nerve(\calJ), \calO^{\otimes})$ is
a $p$-right Kan extension of $F | \Delta^1$ if and only if, for every nonnegative integer $n$,
the maps $F( \Colp{i}): F(\seg{n}) \rightarrow F(1)$ exhibit $F( \seg{n} )$ as a
$p$-product of the objects $\{ F(1) \}_{1 \leq i \leq n}$. 
\end{itemize}

Since $\calO^{\otimes}$ is an $\infty$-operad, condition $(i')$ is equivalent to the requirement
that each of the maps $F( \Colp{i}): F(\seg{n}) \rightarrow F(1)$ is inert.

We now prove $(i)$. Suppose first that $F$ is a $p$-right Kan extension
of $F | \Delta^1$. We will show that $F \in \calC$. We first show that
the map $F( s(i) ) \rightarrow F(i)$ is an equivalence in $\calO^{\otimes}$ for
$i \in [1]$. If $i = 0$ this is clear (the $\infty$-category $\calO^{\otimes}_{\seg{0}}$ is a contractible Kan complex, so every morphism in $\calO^{\otimes}_{\seg{0}}$ is an equivalence). If
$i = 1$, we apply condition $(i')$ in the case $n=1$. Using this observation,
we see that the condition of $(i')$ is equivalent to the requirement that the
maps $\Colp{i}: \seg{n} \rightarrow \seg{1}$ induce inert morphisms
$F( \seg{n} ) \rightarrow F(\seg{1})$ for $1 \leq i \leq n$, so that
$F| \OpE{0} \in \Alg_{ \OpE{0}}(\calO)$. This proves that $F \in \calC$ as desired.

Conversely, suppose that $F \in \calC$. If $1 \leq i \leq n$, then the map
$F( \Colp{i}): F( \seg{n} \rightarrow F(1)$ factors as a composition
$$ F( \seg{n} ) \rightarrow F( \seg{1} ) \rightarrow F(1).$$
The first map is inert by virtue of our assumption that
$F| \OpE{0}  \in \Alg_{ \OpE{0}}( \calO)$, and the second map is an
equivalence since $F$ satisfies condition $(a)$. It follows that
$F$ satisfies the hypothesis of $(i')$ and is therefore a $p$-right Kan extension
of $F| \Delta^1$. This completes the proof of $(i)$.

The proof of $(ii)$ is similar: using Lemma \toposref{kan0}, we can reduce to showing that
for each $n \geq 0$ the composite map
$$ \Nerve(\calI_1) \rightarrow \{1\} \subseteq \Delta^1 \stackrel{F_0}{\rightarrow} \calO^{\otimes}$$
can be extended to a $p$-limit diagram (lying over the canonical map
$\Nerve(\calI_1)^{\triangleleft} \rightarrow \Nerve(\calJ) \rightarrow \OpE{0}$). 
The existence of such a diagram follows immediately from Proposition \ref{saber}.
\end{proof}

Suppose that $p: \calC^{\otimes} \rightarrow \Nerve(\FinSeg)$ and $q: \calD^{\otimes} \rightarrow \Nerve(\FinSeg)$ are symmetric monoidal $\infty$-categories. An $\infty$-operad map
$F \in \Alg_{\calC}(\calD) = \Fun^{\lax}( \calC^{\otimes}, \calD^{\otimes})$ can be thought of as a functor $F: \calC \rightarrow \calD$ which is compatible with the symmetric monoidal structures in the sense that we are given maps
$$ F(C) \otimes F(C') \rightarrow F(C \otimes C')$$
$$ { \bf 1} \rightarrow F( {\bf 1} )$$
which are compatible with the commutativity and associativity properties of the tensor
products on $\calC$ and $\calD$. Our next definition singles out a special class of $\infty$-category maps for which the morphisms above are required to be equivalences.

\begin{definition}
Let $\calO^{\otimes}$ be an $\infty$-operad, and let $p: \calC^{\otimes} \rightarrow \calO^{\otimes}$ and $q: \calD^{\otimes} \rightarrow \calO^{\otimes}$ be coCartesian fibrations of $\infty$-operads. 
We will say that an $\infty$-operad map $f \in \Alg_{\calC}(\calD)$ is a
{\it $\calO$-monoidal functor} if it carries $p$-coCartesian morphisms to $q$-coCartesian morphisms.
We let $\Fun^{\otimes}_{\calO}( \calC, \calD)$ denote the full subcategories
$\Fun_{ \calO^{\otimes}}( \calC^{\otimes}, \calD^{\otimes})$ spanned by the
$\calO$-monoidal functors.

In the special case where $\calO^{\otimes} = \Nerve(\FinSeg)$, we will denote
$\Fun^{\otimes}_{\calO}( \calC, \calD)$ by $\Fun^{\otimes}(\calC, \calD)$; we will refer
to objects of $\Fun^{\otimes}(\calC, \calD)$ as {\it symmetric monoidal functors} from
$\calC^{\otimes}$ to $\calD^{\otimes}$.
\end{definition}

\begin{remark}\label{slame}
Let $\calO^{\otimes}$ be an $\infty$-operad. We define a simplicial category
$\Cat_{\infty}^{\cDelta, \lax, \calO}$ whose objects are fibrations of $\infty$-operads
$p: \calC^{\otimes} \rightarrow \calO^{\otimes}$, where
$\bHom_{ \Cat_{\infty}^{\cDelta, \lax, \calO}}( \calC^{\otimes}, \calD^{\otimes})$ is the largest
Kan complex contained in the $\infty$-category $\Alg_{\calC}(\calD)$. 
We let $\Cat_{\infty}^{\lax, \calO}$ denote the simplicial nerve $\sNerve( \Cat_{\infty}^{\cDelta, \lax, \calO})$; we will refer to $\Cat_{\infty}^{\lax, \calO}$ as the {\it $\infty$-category of $\calO$-operads}.
Let $\Cat_{\infty}^{\calO}$ denote the subcategory of $\Cat_{\infty}^{\lax, \calO}$ spanned by the
$\calO$-monoidal $\infty$-categories and $\calO$-monoidal functors between them.
\end{remark}

\begin{remark}\label{symmetricegg}
Let $F: \calC^{\otimes} \rightarrow \calD^{\otimes}$ be a $\calO$-monoidal functor between
$\calO$-monoidal $\infty$-categories, where $\calO^{\otimes}$ is an $\infty$-operad.
Using Corollary \toposref{usefir}, we deduce that the following conditions are equivalent:
\begin{itemize}
\item[$(1)$] The functor $F$ is an equivalence.
\item[$(2)$] The underlying map of $\infty$-categories $\calC \rightarrow \calD$ is an equivalence.
\item[$(3)$] For every color $X \in \calO$, the induced map of fibers
$\calC_{X} \rightarrow \calD_{X}$ is an equivalence.
\end{itemize}
The analogous statement for lax $\calO$-monoidal functors is false.
\end{remark}

Our goal for the remainder of this section is to describe the $\infty$-category of
$\Ass$-algebra objects.

\begin{construction}\label{urpas}
Let $\cDelta$ denote the category of {\it combinatorial simplices}: the objects of
$\cDelta$ are the linearly ordered sets $[n] = \{ 0, \ldots, n\}$, and the morphisms are (nonstrictly) order-preserving maps between linearly ordered sets. We define a functor $\phi:
\cDelta^{op} \rightarrow \CatAss$ as follows:
\begin{itemize}
\item[$(1)$] For each $n \geq 0$, we have $\phi( [n] ) = \seg{n}$.
\item[$(2)$] Given a morphism $\alpha: [n] \rightarrow [m]$ in $\cDelta$, the associated morphism
$\phi(\alpha): \seg{m} \rightarrow \seg{n}$ is given by the formula
$$\phi(\alpha)(i) = \begin{cases} j & \text{if } (\exists j) [\alpha(j-1) < i \leq \alpha(j)] \\
\ast & \text{otherwise.} \end{cases}$$
where we endow each $\phi(\alpha)^{-1} \{ j\}$ with the linear ordering induces by its inclusion into
$\nostar{n}$.
\end{itemize}
More informally, $\phi$ assigns to a nonempty linearly ordered set $[n]$ the set of
all ``gaps'' between adjacent elements of $[n]$.

Suppose that $p: \calC^{\otimes} \rightarrow \Ass$ is an $\Ass$-monoidal $\infty$-category.
Pulling back via the functor $\phi$, we obtain a coCartesian fibration
$q: \calC^{\otimes} \times_{ \Ass} \Nerve( \cDelta^{op}) \rightarrow \Nerve(\cDelta^{op})$.
Unwinding the definitions, we deduce that $q$ exhibits $\calC^{\otimes} \times_{\Ass} \Nerve(\cDelta)^{op}$ as a monoidal $\infty$-category, in the sense of Definition \monoidref{mainef}. 
We will refer to $\calC^{\otimes} \times_{\Ass} \Nerve(\cDelta)^{op}$ as the
{\it underlying monoidal $\infty$-category of $\calC^{\otimes}$}. 
\end{construction}

Algebras over the associative $\infty$-operad $\Ass$ can be described as follows:

\begin{proposition}\label{algass}
Let $q: \calC^{\otimes} \rightarrow \Ass$ be a fibration of $\infty$-operads, and let
$\calD^{\otimes} \rightarrow \Nerve( \cDelta^{op} )$ be the pullback of $q$
along the nerve of the functor $\phi: \cDelta^{op} \rightarrow \CatAss$ of Construction
\ref{urpas}. Then composition with $\phi$ induces an equivalence of $\infty$-categories
$\theta: \Alg_{\Ass}(\calC) \rightarrow \Alg(\calD)$, where $\Alg(\calD) \subseteq \Fun_{ \Nerve(\cDelta^{op}}( \Nerve(\cDelta)^{op}, \calD^{\otimes}) \simeq \Fun_{ \Ass}(\Nerve(\cDelta^{op}, \calC^{\otimes})$ is the full subcategory spanned by those sections which carry every convex morphism in $\cDelta^{op}$ to inert morphisms in $\calC^{\otimes}$.
\end{proposition}

\begin{remark}
In the situation of Proposition \ref{algass}, if $q$ is a coCartesian fibration of $\infty$-operads,
then $\calD^{\otimes}$ is a monoidal $\infty$-category and our definition of $\Alg(\calD)$ agrees with the one given in Definition \monoidref{suskin}. 
\end{remark}

\begin{proof}
The proof is similar to that of Proposition \monoidref{algcompare}.
We define a category $\calI$ as follows:
\begin{itemize}
\item[$(1)$] An object of $\calI$ is either an object of $\cDelta^{op}$ or an object of $\CatAss$.
\item[$(2)$] Morphisms in $\calI$ are give by the formulas
$$\Hom_{\calI}( [m], [n] ) = \Hom_{ \cDelta^{op} }( [m], [n]) \quad \Hom_{\calI}( \seg{m}, \seg{n} ) = \Hom_{\CatAss}( \seg{m}, \seg{n} )$$
$$ \bHom_{ \calI}( \seg{m}, [n] ) = \bHom_{ \LinSeg}( \seg{m}, \phi( [n] ) ) \quad
 \bHom_{ \calI}( [n], \seg{m} ) = \emptyset.$$ 
\end{itemize}
where $\phi: \cDelta^{op} \rightarrow \CatAss$ is the functor defined in Construction
\ref{urpas}. We observe that $\phi$ extends to a retraction $\overline{\phi}: \calI \rightarrow \CatAss$.
Let $\overline{\Alg}(\calC)$ denote the full subcategory of 
$\bHom_{ \Ass}( \Nerve(\calI), \calC^{\otimes})$ consisting of those functors $f: \Nerve(\calI) \rightarrow \calC^{\otimes}$ 
such that $q \circ f = \overline{\psi}$ and the following additional conditions are satisfied:
\begin{itemize}
\item[$(i)$] For each $n \geq 0$, $f$ carries the canonical map $\seg{n} \rightarrow [n]$ in
$\calI$ to an equivalence in $\calC^{\otimes}$.
\item[$(ii)$] The restriction $f | \Nerve(\cDelta)^{op}$ belongs to $\Alg(\calD)$.
\item[$(ii')$] The restriction $f | \Ass$ is an $\Ass$-algebra object of $\calC$.
\end{itemize}

If $(i)$ is satisfied, then $(ii)$ and $(ii')$ are equivalent to one another. Moreover, $(i)$ is equivalent to the assertion that $f$ is a $q$-left Kan extension of $f | \Ass$.
Since every functor $f_0: \Ass \rightarrow \calC$ admits a $q$-left Kan extension (given, for example, by $f_0 \circ \overline{\psi}$), Proposition \toposref{lklk} implies that the restriction map
$p: \overline{\Alg}(\calC) \rightarrow \Alg_{\Ass}(\calC)$ is a trivial Kan fibration. The map
$\theta$ is the composition of a section to $p$ (given by composition with $\overline{\psi}$) 
and the restriction map $p': \overline{\Alg}(\calC) \rightarrow \Alg(\calD)$. It will therefore suffice to show that $p'$ is a trivial fibration. In view of Proposition \toposref{lklk}, this will follow from the following pair of assertions:

\begin{itemize}
\item[$(a)$] Every $f_0 \in \Alg(\calD)$ admits a $q$-right Kan extension
$f \in \bHom_{ \Ass}( \Nerve(\calI), \calC^{\otimes})$.
\item[$(b)$] Given $f \in \bHom_{ \Ass}( \Nerve(\calI), \calC^{\otimes})$
such that $f_0 = f | \Nerve(\cDelta)^{op}$ belongs to $\Alg(\calD)$, $f$ is a $q$-right Kan extension of $f_0$ if and only if $f$ satisfies condition $(i)$ above.
\end{itemize}

To prove $(a)$, we fix an object $\seg{n} \in \CatAss$. Let $\calJ$ denote the category
$( \cDelta)^{op} \times_{ \CatAss} ( \CatAss)_{\seg{n}/ },$
and let $g$ denote the composition
$ \Nerve(\calJ) \rightarrow \Nerve(\cDelta)^{op} \rightarrow \calC^{\otimes}.$
According to Lemma \toposref{kan2}, it will suffice to show that $g$ admits a $q$-limit in $\calC^{\otimes}$ (for each $n \geq 0$). The objects of $\calJ$ can be identified with
morphisms $\alpha: \seg{n} \rightarrow \phi([m])$ in $\CatAss$. Let $\calJ_0 \subseteq \calJ$
denote the full subcategory spanned by those objects for which $\alpha$ is inert.
The inclusion $\calJ_0 \subseteq \calJ$ has a right adjoint, so that $\Nerve(\calJ_0)^{op} \rightarrow \Nerve(\calJ)^{op}$ is cofinal. Consequently, it will suffice to show that $g_0 = g | \Nerve(\calJ_0)$ admits a $q$-limit in $\calC$.

Let $\calJ_1$ denote the full subcategory of $\calJ_0$ spanned by the morphisms
$\Colp{j}: \seg{n} \rightarrow \phi( [1] )$. Using our assumption that
$f_0$ is an algebra object of $\calD$, we deduce that $g_0$ is a $q$-right Kan extension of
$g_1 = g_0 | \Nerve(\calJ_1)$. In view of Lemma \toposref{kan0}, it will suffice to show that
the map $g_1$ has a $q$-limit in $\calC$. But this is clear; our assumption that $f_0$ belongs to $\Alg(\calD)$ guarantees that $f_0$ exhibits $f_0([n])$ as a $q$-limit of $g_1$. This proves $(a)$. Moreover, the proof shows that $f$ is a $q$-right Kan extension of $f_0$ at $\seg{n}$ if and only if $f$ induces an equivalence $f( \seg{n} ) \rightarrow f( [n] )$; this immediately implies $(b)$ as well.
\end{proof}

\subsection{Cartesian Monoidal Structures}\label{comm2}

Our goal in this section is to prove that if $\calC$ is an $\infty$-category which admits finite products, then the formation of products endows $\calC$ with the structure of a symmetric monoidal $\infty$-category. The collection of symmetric monoidal structures which arise in this way is easy to characterize:

\begin{definition}
Let $\calC$ be an $\infty$-category. We will say that a monoidal structure on $\calC$ is {\it Cartesian} if the following conditions are satisfied:
\begin{itemize}
\item[$(1)$] The unit object $1_{\calC} \in \calC$ is final.
\item[$(2)$] For every pair of objects $C, D \in \calC$, the canonical maps
$$ C \simeq C \otimes 1_{\calC} \leftarrow C \otimes D \rightarrow 1_{\calC} \otimes D \simeq D$$
exhibit $C \otimes D$ as a product of $C$ and $D$ in the $\infty$-category $\calC$.
\end{itemize}
We will say that a symmetric monoidal structure on $\calC$ is {\it Cartesian} if the underlying monoidal structure (see Construction \ref{urpas}) on $\calC$ is Cartesian.
\end{definition}

Our first goal is to show that if $\calC$ admits finite products, then there is an essentially unique
Cartesian symmetric monoidal structure on $\calC$. To prove this, we need a criterion for detecting when a symmetric monoidal structure is Cartesian.

\begin{definition}\label{puggybear}
Let $p: \calC^{\otimes} \rightarrow \Nerve(\FinSeg)$ be an $\infty$-operad.
A {\it lax Cartesian structure} on $\calC^{\otimes}$ is a functor
$\pi: \calC^{\otimes} \rightarrow \calD$ satisfying the following condition:
\begin{itemize}
\item[$(\ast)$] Let $C$ be an object of $\calC^{\otimes}_{\seg{n}}$, and write
$C \simeq C_1 \oplus \cdots \oplus C_n$. Then the canonical maps
$\pi(C) \rightarrow \pi(C_i)$ exhibit $\pi(C)$ as a product 
$\prod_{1 \leq j \leq n} \pi(C_j)$ in the $\infty$-category $\calD$.
\end{itemize}

We will say that $\pi$ is a {\it weak Cartesian structure} if it is a lax Cartesian structure, 
$\calC^{\otimes}$ is a symmetric monoidal $\infty$-category, and the following additional condition is satisfied:
\begin{itemize}
\item[$(\ast')$] Let $f: C \rightarrow C'$ be a $p$-coCartesian morphism covering an active morphism $\alpha: \seg{n} \rightarrow \seg{1}$ in $\FinSeg$. Then $\pi(f)$ is an equivalence in
$\calD$.
\end{itemize}
We will say that a weak Cartesian structure $\pi$ is a {\it Cartesian structure} if $\pi$ induces an equivalence $\calC^{\otimes}_{\seg{1}} \rightarrow \calD$. 
\end{definition}

It follows immediately from the definition that if $\calC$ is a symmetric monoidal $\infty$-category and there exists a Cartesian structure $\calC^{\otimes} \rightarrow \calD$, then the symmetric monoidal structure on $\calC$ is Cartesian. Our first result is a converse: if $\calC$ is a Cartesian symmetric monoidal $\infty$-category, then there exists an essentially unique Cartesian structure on $\calC$.

\begin{proposition}\label{cardbade2}
Let $p: \calC^{\otimes} \rightarrow \Nerve(\FinSeg)$ be a Cartesian symmetric monoidal $\infty$-category and let $\calD$ be an $\infty$-category which admits finite products.
Let $\Fun^{\times}( \calC^{\otimes}, \calD)$ denote the full subcategory of $\Fun( \calC^{\otimes}, \calD)$ spanned by the weak Cartesian structures, and let $\Fun^{\times}(\calC, \calD)$ be the full subcategory of $\Fun(\calC, \calD)$ spanned by those functors which preserve finite products. The restriction map
$ \Fun^{\times}( \calC^{\otimes}, \calD) \rightarrow \Fun^{\times}(\calC, \calD)$ is an equivalence
of $\infty$-categories.
\end{proposition}

The proof will be given at the end of this section.

\begin{corollary}
Let $\calC^{\otimes}$ be a Cartesian symmetric monoidal $\infty$-category. Then
there exists an (essentially unique) Cartesian structure $\calC^{\otimes} \rightarrow \calC$
extending the identity map $\id_{\calC}$.
\end{corollary}

We now treat the problem of {\it existence} for Cartesian structures. We wish to show that, if
$\calC$ is an $\infty$-category which admits finite products, then there exists a Cartesian symmetric monoidal structure on $\calC$. We will prove this by giving an explicit construction of
a symmetric monoidal $\infty$-category $\calC^{\times}$ equipped with a Cartesian
structure $\calC^{\times} \rightarrow \calC$.

\begin{notation}\label{hugegrinn}
The category $\Fintimes$ is defined as follows:
\begin{itemize}
\item[$(1)$] An object of $\Fintimes$ consists of an ordered pair $(\seg{n}, S)$, where
$\seg{n}$ is an object of $\FinSeg$ and $S$ is a subset of $\nostar{n}$.

\item[$(2)$] A morphism from $(\seg{n}, S)$ to $(\seg{n'}, S')$ in $\Fintimes$
consists of a map $\alpha: \seg{n} \rightarrow \seg{n'}$ in $\FinSeg$ with the property that
$\alpha^{-1} S' \subseteq S$.
\end{itemize}
\end{notation}

We observe that the forgetful functor $\Fintimes \rightarrow \FinSeg$ is a Grothendieck fibration, so that the induced map of $\infty$-categories $\Nerve(\Fintimes) \rightarrow \Nerve(\FinSeg)$ is a Cartesian fibration (Remark \toposref{gcart}).

\begin{remark}\label{hungerfope}
The forgetful functor $\Fintimes \rightarrow \FinSeg$ has a canonical section $s$, given by
$s( \seg{n} )  = ( \seg{n}, \nostar{n})$.
\end{remark}

\begin{construction}\label{ilkur}
Let $\calC$ be an $\infty$-category. We define a simplicial set $\widetilde{\calC}^{\times}$ equipped with a map $\widetilde{\calC}^{\times} \rightarrow \Nerve(\FinSeg)$ by the following universal property: for every map of simplicial sets $K \rightarrow \Nerve(\FinSeg)$, we have a bijection
$$ \Hom_{ \Nerve(\cDelta)^{op} }(K, \widetilde{\calC}^{\times} )
\simeq \Hom_{ \sSet }( K \times_{ \Nerve(\FinSeg) } \Nerve(\Fintimes), \calC).$$

Fix $\seg{n} \in \FinSeg$. We observe that the fiber $\widetilde{\calC}^{\times}_{\seg{n}}$ can be identified with the $\infty$-category of functors $f: \Nerve(P)^{op} \rightarrow \calC$, where $P$ is the partially ordered set of subsets of $\nostar{n}$. We let $\calC^{\times}$ be the full simplicial subset
of $\widetilde{\calC}^{\times}$ spanned by those vertices which correspond to those functors $f$
with the property that for every $S \subseteq \nostar{n}$, the maps
$f(S) \rightarrow f( \{ j \})$ exhibit $f(S)$ as a product of the objects
$\{ f( \{ j \}  \}_{j \in S}$ in the $\infty$-category $\calC$.
\end{construction}

The fundamental properties of Construction \ref{ilkur} are summarized in the following
pair of results, whose proofs we defer until the end of this section:

\begin{proposition}\label{commonton2}
Let $\calC$ be an $\infty$-category.
\begin{itemize}
\item[$(1)$] The projection $p: \widetilde{\calC}^{\times} \rightarrow \Nerve(\FinSeg)$ is a coCartesian fibration.
\item[$(2)$] Let $\overline{\alpha}: f \rightarrow f'$ be a morphism of $\widetilde{\calC}^{\times}$ whose image in $\Nerve(\FinSeg)$ corresponds to a map $\alpha: \seg{n} \rightarrow \seg{n'}$. Then $\overline{\alpha}$ is $p$-coCartesian if and only if, for every $S \subseteq \nostar{n'}$, the induced map $f( \alpha^{-1} S ) \rightarrow f(S)$ is an equivalence in $\calC$.

\item[$(3)$] The projection $p$ restricts to a coCartesian fibration $\calC^{\times} \rightarrow \Nerve(\FinSeg)$ $($with the same class of coCartesian morphisms$)$.
\item[$(4)$] The projection $\calC^{\times} \rightarrow \Nerve(\FinSeg)$ is a symmetric monoidal $\infty$-category if and only if $\calC$ admits finite products.
\item[$(5)$] Suppose that $\calC$ admits finite products. 
Let $\pi: \calC^{\times} \rightarrow \calC$ be the map given by composition with
the section $s: \Nerve(\FinSeg) \rightarrow \Nerve( \Fintimes)$ defined in Remark \ref{hungerfope}. Then $\pi$ is a Cartesian structure on $\calC^{\times}$.
\end{itemize}
\end{proposition}
 
\begin{proposition}\label{sungtoo}
Let $\calO^{\otimes}$ be an $\infty$-operad, $\calD$ an $\infty$-category which admits finite products, and $\pi: \calD^{\times} \rightarrow \calD$ the Cartesian structure of Proposition \ref{commonton2}. Then
composition with $\pi$ induces a trivial Kan fibration
$$ \theta: \Fun^{\lax}( \calO^{\otimes}, \calD^{\times}) \rightarrow \Fun^{\lax}( \calO^{\otimes}, \calD)$$
where $\Fun^{\lax}( \calO^{\otimes}, \calD)$ denotes the full subcategory spanned by the lax
Cartesian structures. If $\calO^{\otimes}$ is a symmetric monoidal $\infty$-category, then
composition with $\pi$ induces a trivial Kan fibration
$$ \theta_0: \Fun^{\otimes}( \calO^{\otimes}, \calD^{\times}) \rightarrow \Fun^{\times}( \calO^{\otimes}, \calD)$$
where $\Fun^{\times}( \calO^{\otimes}, \calD)$ denotes the full subcategory of
$\Fun( \calO^{\otimes}, \calD)$ spanned by the weak Cartesian structures.
\end{proposition}

We are now in a position to establish the uniqueness of Cartesian symmetric monoidal structures. 
Let $\Cat_{\infty}^{\CommOp, \otimes}$ denote the $\infty$-category of symmetric monoidal $\infty$-categories (Remark \ref{slame}), $\Cat_{\infty}^{\CommOp, \times}$ the full subcategory spanned by the
Cartesian symmetric monoidal $\infty$-categories, and $\Cat_{\infty}^{\Cart}$ the subcategory
of $\Cat_{\infty}$ whose objects are $\infty$-categories which admit finite products, and whose morphisms are functors which preserve finite products. 

\begin{corollary}\label{unisim}
The forgetful functor
$\theta: \Cat_{\infty}^{\CommOp, \times} \rightarrow \Cat_{\infty}^{\Cart}$ is an equivalence of $\infty$-categories.
\end{corollary}

\begin{proof}
We first observe that if $\calC$ is an $\infty$-category which admits finite products, then
Proposition \ref{commonton2} implies that the symmetric monoidal $\infty$-category $\calC^{\times}$ is a preimage of $\calC$ under the forgetful functor $\theta$. It follows that $\theta$ is essentially surjective. Moreover, if $\calC^{\otimes}$ is {\em any} Cartesian symmetric monoidal structure on $\calC$, then Proposition \ref{cardbade2} guarantees the existence of a Cartesian structure
$\pi: \calC^{\otimes} \rightarrow \calC$. Applying Proposition \ref{sungtoo}, we can lift
$\pi$ to a symmetric monoidal functor $F: \calC^{\otimes} \rightarrow \calC^{\times}$. Remark \monoidref{eggsal} guarantees that $F$ is an equivalence.

We now show that $\theta$ is fully faithful. Let $\calC^{\otimes}$ and $\calD^{\otimes}$ be Cartesian symmetric monoidal structures on $\infty$-categories $\calC$ and $\calD$. We wish to show that the restriction map
$$ \bHom_{ \Cat_{\infty}^{\CommOp} }( \calC^{\otimes}, \calD^{\otimes} )
\rightarrow \bHom_{ \Cat_{\infty}^{\Cart} }( \calC, \calD )$$
is a homotopy equivalence. We will prove a slightly stronger assertion: namely, that the restriction map
$\psi: \Fun^{\otimes}( \calC^{\otimes}, \calD^{\otimes} ) \rightarrow \Fun^{\times}(\calC, \calD)$
is a categorical equivalence, where $\Fun^{\times}(\calC, \calD)$ denotes the full subcategory of
$\Fun(\calC, \calD)$ spanned by those functors which preserve finite products. In view of the above remarks, it suffices to prove this in the case where $\calD^{\otimes} = \calD^{\times}$. 
In this case, the map $\psi$ factors as a composition
$$ \Fun^{\otimes}(\calC^{\otimes}, \calD^{\otimes} )
\stackrel{\psi'}{\rightarrow} \Fun^{\times}( \calC^{\otimes}, \calD)
\stackrel{\psi''}{\rightarrow} \Fun^{\times}(\calC, \calD).$$
Proposition \ref{sungtoo} implies that $\psi'$ is a categorical equivalence, and 
Proposition \ref{cardbade2} implies that $\psi''$ is a categorical equivalence.
\end{proof}

We now study algebra objects in a Cartesian symmetric monoidal $\infty$-category.

\begin{definition}\label{gammaobj}
Let $\calC$ be an $\infty$-category and let $\calO^{\otimes}$ be an $\infty$-operad.
A {\it $\calO$-monoid} in $\calC$ is a functor
$M: \calO^{\otimes} \rightarrow \calC$ with the following property: for
every object $X \in \calO^{\otimes}_{\seg{n}}$ corresponding to a sequence of objects
$\{ X_i \in \calO \}_{1 \leq i \leq n}$, the canonical maps
$M(X) \rightarrow M(X_i)$ exhibit $M(X)$ as a product $\prod_{1 \leq i \leq n} M(X_i)$
in the $\infty$-category $\calC$. We let $\Mon_{\calO}( \calC)$ denote the full subcategory
of $\Fun( \calO^{\otimes}, \calC)$ spanned by the $\calO$-monoids in $\calC$.
\end{definition}

\begin{remark}
In the special case where $\calO^{\otimes}$ is the commutative $\infty$-operad, we will
refer to $\calO$-monoids in an $\infty$-category $\calC$ as {\it commutative monoid objects}
of $\calC$. These objects might also be referred to as {\it $\Gamma$-objects}
of $\calC$; in the special case where $\calC$ is the $\infty$-category of spaces, the theory of $\Gamma$ objects is essentially equivalent to Segal's theory of $\Gamma$-spaces.
\end{remark}


\begin{example}\label{interplace}
Let $\calO^{\otimes}$ be an $\infty$-operad. A functor
$M: \calO^{\otimes} \rightarrow \Cat_{\infty}$ is a $\calO$-monoid in $\Cat_{\infty}$ if and only if 
the coCartesian fibration $\calC^{\otimes} \rightarrow \calO^{\otimes}$ classified by $M$ is a $\calO$-monoidal $\infty$-category. Arguing as in Remark \monoidref{otherlander}, we deduce that there is a canonical equivalence of $\Cat_{\infty}^{\calO, \otimes}$ with the $\infty$-category of $\calO$-monoids
$\Mon_{\calO}( \Cat_{\infty})$.
\end{example}

\begin{proposition}\label{ungbatt}
Let $\calC^{\otimes}$ be a symmetric monoidal $\infty$-category,
$\pi: \calC^{\otimes} \rightarrow \calD$ a Cartesian structure, and
$\calO^{\otimes}$ an $\infty$-operad. Then composition with $\pi$ induces an equivalence of $\infty$-categories
$\Alg_{\calO}(\calC) \rightarrow \Mon_{\calO}( \calD)$.
\end{proposition}

\begin{proof}
As in the proof of Corollary \ref{unisim}, we may assume without loss of generality that
$\calC^{\otimes} = \calD^{\times}$. We now apply Proposition \ref{sungtoo} again to deduce that the map 
$$\Alg_{\calO}(\calC) = \Fun^{\lax}( \calO^{\otimes}, \calC^{\otimes}) \rightarrow
\Fun^{\lax}( \calO^{\otimes}, \calD) = \Mon_{\calO}(\calD)$$
is a trivial Kan fibration.
\end{proof}

\begin{corollary}\label{co1}
Let $\calO^{\otimes}$ be an $\infty$-operad, and regard $\Cat_{\infty}$ as endowed
with the Cartesian symmetric monoidal structure. Then we have a canonical
equivalence of $\infty$-categories
$$ \Cat_{\infty}^{\calO} \simeq \Alg_{ \calO}( \Cat_{\infty}).$$
\end{corollary}

\begin{proof}
Combine Proposition \ref{ungbatt} with Example \ref{interplace}.
\end{proof}

\begin{corollary}
Let $\Triv$ denote the trivial $\infty$-operad of Example \ref{trivop}.
Then passing to the fiber over the object $\seg{1} \in \Triv$ induces an
equivalence of $\infty$-categories $\Cat_{\infty}^{\Triv} \rightarrow \Cat_{\infty}$.
\end{corollary}

\begin{proof}
Combine Corollary \ref{co1} with Example \ref{algtriv}.
\end{proof}

\begin{corollary}\label{stix}
Let $\Ass$ denote the associative $\infty$-operad of Example \ref{xe2}. Then
Construction \ref{urpas} induces an equivalence of
$\infty$-categories
$$ \Cat_{\infty}^{\Ass} \rightarrow \Cat_{\infty}^{\Mon},$$
where $\Cat_{\infty}^{\Mon}$ is the $\infty$-category of monoidal $\infty$-categories $($Definition \monoidref{monfunc2}$)$. 
\end{corollary}

\begin{proof}
Combine Remark \monoidref{otherlander}, Corllary \ref{co1}, and 
Proposition \ref{algass}.
\end{proof}

The following criterion is useful for establishing that a symmetric monoidal structure on an $\infty$-category $\calC$ is Cartesian:

\begin{proposition}\label{fingertone}
Let $\calC$ be a symmetric monoidal $\infty$-category. The following conditions are equivalent:
\begin{itemize}
\item[$(1)$] The symmetric monoidal structure on $\calC$ is Cartesian.
\item[$(2)$] The induced symmetric monoidal structure on the homotopy category
$\h{\calC}$ is Cartesian.
\item[$(3)$] The unit object $1_{\calC}$ is final, and for each object $C \in \calC$ there
exists a {\it diagonal map} $\delta_{C}: C \rightarrow C \otimes C$ satisfying the following conditions:
\begin{itemize}
\item[$(i)$] Let $C$ be an object of $\calC$ and let $u: C \rightarrow 1_{\calC}$ be a map $(${}automatically unique up to homotopy{}$)$. Then the composition
$$ C \stackrel{\delta_{C}}{\rightarrow} C \otimes C \stackrel{\id \otimes u}{\longrightarrow} C \otimes 1_{\calC} \rightarrow C$$
is homotopic to the identity.
\item[$(ii)$] For every morphism $f: C \rightarrow D$ in $\calC$, the diagram
$$ \xymatrix{ C \ar[r]^{f} \ar[d]^{\delta_{C} } & D \ar[d]^{\delta_{D} } \\
C \otimes C \ar[r]^{f \otimes f} & D \otimes D }$$
commutes up to homotopy.
\item[$(iii)$] Let $C$ and $D$ be objects of $\calC$. Then the diagram
$$ \xymatrix{ &  C \otimes D \ar[dl]_{\delta_{C} \otimes \delta_{D} } \ar[dr]^{\delta_{C \otimes D}} \\
(C \otimes C) \otimes (D \otimes D) \ar[rr]^{\sim} & & (C \otimes D) \otimes (C \otimes D) }$$
commutes up to homotopy.
\end{itemize}
\end{itemize}
\end{proposition}

\begin{proof}
The implications $(1) \Rightarrow (2) \Rightarrow (3)$ are obvious. Let us suppose that $(3)$ is satisfied. We wish to show that, for every pair of objects $C, D \in \calC$, the maps
$$ C \simeq C \otimes 1_{\calC} \leftarrow C \otimes D \rightarrow 1_{\calC} \otimes D \simeq D$$
exhibit $C \otimes D$ as a product of $C$ and $D$ in $\calC$. In other words, we must show that
for every object $A \in \calC$, the induced map
$$ \phi: \bHom_{\calC}( A, C \otimes D) \rightarrow \bHom_{\calC}(A,C) \times \bHom_{\calC}(A,D)$$
is a homotopy equivalence. Let $\psi$ denote the composition
$$ \bHom_{\calC}(A,C) \otimes \bHom_{\calC}(A,D)
\stackrel{\otimes}{\rightarrow} \bHom_{\calC}(A \otimes A, C \otimes D)
\stackrel{\delta_{A}}{\rightarrow} \bHom_{\calC}(A, C \otimes D).$$
We claim that $\psi$ is a homotopy inverse to $\phi$. The existence of a homotopy
$\psi \circ \phi \simeq \id$ follows from $(i)$. We will show that $\psi \circ \phi$ is homotopic to the identity. In view of condition $(ii)$, $\psi \circ \phi$ is homotopic to the map defined by composition with
$$ C \otimes D \stackrel{ \delta_{C \otimes D} }{\rightarrow} (C \otimes D) \otimes (C \otimes D)
\rightarrow (C \otimes 1_{\calC} ) \otimes (1_{\calC} \otimes D)
\simeq C \otimes D.$$
Invoking $(iii)$, we deduce that this composition is homotopic with the map
$$ C \otimes D \stackrel{ \delta_{C} \otimes \delta_{D} }{\rightarrow}
(C \otimes C) \otimes (D \otimes D) \rightarrow ( C \otimes 1_{\calC} ) \otimes ( 1_{\calC} \otimes D)
\simeq C \otimes D,$$
which is equivalent to the identity in virtue of $(i)$.
\end{proof}

We conclude this section with the proofs of Propositions \ref{cardbade2}, \ref{commonton2},
and \ref{sungtoo}.

\begin{proof}[Proof of Proposition \ref{cardbade2}]
We define a subcategory $\calI \subseteq \FinSeg \times [1]$ as follows:

\begin{itemize}
\item[$(a)$] Every object of $\FinSeg \times [1]$ belongs to $\calI$.
\item[$(b)$] A morphism $( \seg{n}, i ) \rightarrow (\seg{n'}, i')$ in $\FinSeg \times [1]$ belongs
to $\calI$ if and only if either $i'=1$ or the induced map $\alpha: \seg{n} \rightarrow \seg{n'}$ satisfies
$\alpha^{-1} \{ \ast \} = \ast$.
\end{itemize}

Let $\calC'$ denote the fiber product $\calC^{\otimes} \times_{ \Nerve(\FinSeg) } \Nerve(\calI)$, which we regard as a subcategory of $\calC^{\otimes} \times \Delta^1$, and let
$p': \calC' \rightarrow \Nerve(\calI)$ denote the projection.
Let $\calC'_0$ and $\calC'_{1}$ denote the intersections of $\calC'$ with
$\calC^{\otimes} \times \{0\}$ and $\calC^{\otimes} \times \{1\}$, respectively.
We note that there is a canonical isomorphism $\calC'_{1} \simeq \calC^{\otimes}$.

Let $\calE$ denote the full subcategory of $\Fun( \calC', \calD)$ spanned by those functors $F$ which satisfy the following conditions:
\begin{itemize}
\item[$(i)$] For every object $C \in \calC^{\otimes}$, the induced map
$F(C,0) \rightarrow F(C,1)$ is an equivalence in $\calD$.
\item[$(ii)$] The restriction $F | \calC'_{1}$ is a weak Cartesian structure on $\calC^{\otimes}$.
\end{itemize}

It is clear that if $(i)$ and $(ii)$ are satisfied, then the restriction 
$F_0 = F | \calC'_0$ satisfies the following additional conditions:

\begin{itemize}
\item[$(iii)$] The restriction $F_0 | \calC^{\otimes}_{\seg{1}} \times \{0\}$ is a functor from $\calC$ to $\calD$ which preserves finite products.

\item[$(iv)$] For every $p'$-coCartesian morphism $\alpha$ in $\calC'_{0}$, the induced map
$F_0(\alpha)$ is an equivalence in $\calD$. 
\end{itemize}

Moreover, $(i)$ is equivalent to the assertion that $F$ is a right Kan extension of $F| \calC'_{1}$. Proposition \toposref{lklk} implies that the restriction map $r: \calE \rightarrow \Fun^{\times}( \calC^{\otimes}, \calD)$ induces a trivial Kan fibration onto its essential image. The map $r$ has a section $s$, given by composition with the projection map
$\calC' \rightarrow \calC^{\otimes}$. The restriction map $\Fun^{\times}( \calC^{\otimes}, \calD) \rightarrow \Fun^{\times}(\calC, \calD)$
factors as a composition
$$ \Fun^{\times}( \calC^{\otimes}, \calD) \stackrel{s}{\rightarrow} \calE \stackrel{e}{\rightarrow} \Fun^{\times}(\calC, \calD),$$
where $e$ is induced by composition with the inclusion
$\calC \subseteq \calC'_0 \subseteq \calC'$. Consequently, it will suffice to prove that $e$ is an equivalence of $\infty$-categories.

Let $\calE_0 \subseteq \Fun( \calC'_0, \calD)$ be the full subcategory spanned by those functors which satisfy conditions $(iii)$ and $(iv)$. The map $e$ factors as a composition 
$$ \calE \stackrel{e'}{\rightarrow} \calE_0 \stackrel{e''}{\rightarrow} \Fun^{\times}(\calC, \calD).$$
Consequently, it will suffice to show that $e'$ and $e''$ are trivial Kan fibrations.

Let $f: \calC'_0 \rightarrow \calD$ be an arbitrary functor, and let
$C \in \calC^{\otimes}_{\seg{n}} \subseteq \calC'_{0}$. There exists a unique map
$\alpha: ( \seg{n}, 0) \rightarrow ( \seg{1}, 0)$ in $\calI$; choose a $p'$-coCartesian morphism
$\overline{\alpha}: C \rightarrow C'$ lifting $\alpha$. We observe that $C'$ is
an initial object of $\calC \times ( \calC'_0 )_{/C'} \times_{ \calC'_0} \calC$. Consequently,
$f$ is a right Kan extension of $f | \calC$ at $C$ if and only if $f( \overline{\alpha})$ is an equivalence.
It follows that $f$ satisfies $(iv)$ if and only if $f$ is a right Kan extension of $f | \calC$. The same argument (and Lemma \toposref{kan0}) shows that every functor $f_0: \calC \rightarrow \calD$ admits a right Kan extension to $\calC'_0$. Applying Proposition \toposref{lklk}, we deduce that $e''$ is a trivial Kan fibration.

It remains to show that $e'$ is a trivial Kan fibration. In view of Proposition \toposref{lklk}, it will suffice to prove the following pair of assertions, for every functor $f \in \calE_0$:

\begin{itemize}
\item[$(1)$] There exist a functor $F: \calC' \rightarrow \calD$ which is a left Kan extension of
$f = F | \calC'_0$.
\item[$(2)$] An arbitrary functor $F: \calC' \rightarrow \calD$ which extends $f$ is a
left Kan extension of $f$ if and only if $F$ belongs to $\calE$.
\end{itemize}

For every finite linearly ordered set $J$, let $J^{+}$ denote the disjoint union $J \coprod \{ \infty \}$,
where $\infty$ is a new element larger than every element of $J$. Let $(C,1) \in \calC^{\otimes}_{J_{\ast}} \times \{1\} \subseteq \calC'$. Since there exists a final object $1_{\calC} \in \calC$, the $\infty$-category $\calC'_0 \times_{ \calC' } \calC'_{/C}$ also has a final object, given by the map $\alpha: (C',0) \rightarrow (C,1)$, where $C' \in \calC^{\otimes}_{J^{+}_{\ast}}$
corresponds, under the equivalence
$$ \calC^{\otimes}_{ J^{+}_{\ast} } \simeq \calC \times \calC^{\otimes}_{J_{\ast}},$$
to the pair $(1_{\calC}, C)$. We now apply Lemma \toposref{kan2} to deduce $(1)$, together with the following analogue of $(2)$:

\begin{itemize}
\item[$(2')$] An arbitrary functor $F: \calC' \rightarrow \calD$ which extends $f$ is a left
Kan extension of $f$ if and only if, for every morphism $\alpha: (C',0) \rightarrow (C,1)$ as above,
the induced map $F(C',0) \rightarrow F(C,1)$ is an equivalence in $\calD$.
\end{itemize}

To complete the proof, it will suffice to show that $F$ satisfies the conditions stated in $(2')$ if and only if $F \in \calE$. We first prove the ``if'' direction. Let $\alpha: (C', 0) \rightarrow (C,1)$ be as above; we wish to prove that $F(\alpha): F(C',0) \rightarrow F(C,1)$ is an equivalence in $\calD$. 
The map $\alpha$ factors as a composition
$$ (C',0 ) \stackrel{\alpha'}{\rightarrow} (C', 1) \stackrel{\alpha''}{\rightarrow} (C,1).$$
Condition $(i)$ guarantees that $F(\alpha')$ is an equivalence. Condition $(ii)$ guarantees that
$F( C', 1)$ is equivalent to a product $F(1_{\calC}, 1) \times F(C,1)$, and that
$F(\alpha'')$ can be identified with the projection onto the second factor. Moreover, since
$1_{\calC}$ is a final object of $\calC$, condition $(ii)$ also guarantees that
$F(1_{\calC}, 1)$ is a final object of $\calD$. It follows that $F(\alpha'')$ is an equivalence, so that
$F(\alpha)$ is an equivalence as desired.

Now let us suppose that $F$ satisfies the condition stated in $(2')$. We wish to prove that
$F \in \calE$. Here we must invoke our assumption that the monoidal structure on $\calC$ is Cartesian.
We begin by verifying condition $(i)$.
Let $C \in \calC^{\otimes}_{J_{\ast}}$ for some finite linearly ordered set $J$, and let $\alpha: ( C', 0) \rightarrow (C, 1)$ be defined as above. Let $\beta: ( J_{\ast}, 0 ) \rightarrow 
(J^{+}_{\ast}, 0)$ be the morphism in $\calI$ induced by the inclusion $J \subseteq J^{+}$. Choose
a $p'$-coCartesian morphism $\overline{\beta}: (C,0) \rightarrow (C'', 0)$ lifting $\beta$. Since the final object $1_{\calC} \in \calC$ is also the unit object of $\calC$, we can identify $C''$ with $C'$. The composition
$ (C, 0) \stackrel{\overline{\beta}}{\rightarrow} (C', 1) \stackrel{\alpha}{\rightarrow} (C,1)$
is homotopic to the canonical map $\gamma: (C,0) \rightarrow (C,1)$ appearing in the statement of $(i)$. Condition $(iv)$ guarantees that $F( \overline{\beta} )$ is an equivalence, and
$(2')$ guarantees that $F(\alpha)$ is an equivalence. Using the two-out-of-three property, we deduce that $F(\gamma)$ is an equivalence, so that $F$ satisfies $(i)$.

To prove that $F$ satisfies $(ii)$, we must verify two conditions:
\begin{itemize}
\item[$(ii_0)$] If $\overline{\beta}: (C, 1) \rightarrow (D,1)$ is a $p'$-coCartesian morphism
in $\calC'$, and the underlying morphism $\beta: \seg{m} \rightarrow \seg{n}$ satisfies
$\beta^{-1} \{\ast \} = \{ \ast\}$, then
$F(\overline{\beta})$ is an equivalence.
\item[$(ii_1)$] Let $C \in \calC^{\otimes}_{\seg{n}}$, and choose $p$-coCartesian morphisms
$\gamma_i: C \rightarrow C_{j}$ covering the maps $\Colp{j}: \seg{n} \rightarrow \seg{1}$.
Then the maps $\gamma_i$ exhibit $F(C,1)$ as a product
$\prod_{1 \leq j \leq n} F(C_j, 1)$ in the $\infty$-category $\calD$. 
\end{itemize}

Condition $(ii_0)$ follows immediately from $(i)$ and $(iv)$.
To prove $(ii_1)$, we consider the maps $\alpha: (C',0) \rightarrow (C,1)$ and
$\alpha_j: (C'_j, 0) \rightarrow (C_j,1)$ which appear in the statement of $(2')$. For each $1 \leq j \leq n$, we have a commutative diagram
$$ \xymatrix{ (C',0) \ar[r]^{\alpha} \ar[d]^{\gamma'_j} & (C,1) \ar[d]^{\gamma_j} \\
(C'_{j}, 0) \ar[r]^{\alpha_j} & (C_{j}, 1). }$$
Condition $(2')$ guarantees that the maps $F(\alpha)$ and 
$F(\alpha_i)$ are equivalences in $\calD$. Consequently, it will suffice to show that
the maps $f(\gamma'_i)$ exhibit $f(C',0)$ as a product 
$\prod_{j \in J} f( C'_{j}, 0)$ in $\calD$.
Let $f_0 = f | \calC$. Using condition $(iv)$, we obtain canonical equivalences
$$ f(C',0) \simeq f_0( 1_{\calC} \otimes \bigotimes_{j \in J} C_j) \quad \quad f(C'_j, 0) \simeq f_0 ( 1_{\calC} \otimes C_j )$$
Since condition $(iii)$ guarantees that $f_0$ preserves products, it will suffice to show that
the canonical map
$$ 1_{\calC} \otimes (\bigotimes_{1 \leq j \leq n} C_j)
\rightarrow \bigotimes_{1 \leq j \leq n} ( 1_{\calC} \otimes C_j )$$
is an equivalence in the $\infty$-category $\calC$. This follows easily from our assumption that
the symmetric monoidal structure on $\calC$ is Cartesian, using induction on $n$.
\end{proof}

\begin{proof}[Proof of Proposition \ref{commonton2}]
Assertions $(1)$ and $(2)$ follow immediately from Corollary \toposref{skinnysalad}, and $(3)$ follows from $(2)$ (since $\calC^{\times}$ is stable under the pushforward functors associated
to the coCartesian fibration $p$). 
We now prove $(4)$. If $\calC$ has no final object, then $\calC^{\times}_{\seg{0}}$ is empty; consequently, we may assume without loss of generality that $\calC$ has a final object.
Then $\calC^{\times}_{\seg{1}}$ is isomorphic to the $\infty$-category of diagrams
$X \rightarrow Y$ in $\calC$, where $Y$ is final. It follows that $\pi$ induces an equivalence $\calC^{\times}_{\seg{1}} \simeq \calC$. Consequently,
$\calC^{\times}$ is a symmetric monoidal $\infty$-category if and only if, for each $n \geq 0$, the
functors $\Colll{j}{!} $ determine an equivalence
$\phi: \calC^{\times}_{\seg{n}} \rightarrow \calC^{n}$.
Let $P$ denote the partially ordered set of subsets of $\nostar{n}$, and let
$P_0 \subseteq P$ be the partially ordered set consisting of subsets which consist of a single element.
Then $\calC^{\times}_{\seg{n}}$ can be identified with the
set of functors $f: \Nerve(P)^{op} \rightarrow \calC$ which are right Kan extensions of
$f| \Nerve(P_0)^{op}$, and $\phi$ can be identified with the restriction map determined by the inclusion $P_0 \subseteq P$. According to Proposition \toposref{lklk}, $\phi$ is fully faithful, and is essentially surjective if and only if every functor $f_0: \Nerve(P_0)^{op} \rightarrow \calC$ admits a right Kan extension to
$\Nerve(P)^{op}$. Unwinding the definitions, we see that this is equivalent to the assertion that every finite collection of objects of $\calC$ admits a product in $\calC$. This completes the proof of $(4)$. Assertion $(5)$ follows immediately from $(2)$ and the construction of $\calC^{\times}$.
\end{proof}

\begin{proof}[Proof of Proposition \ref{sungtoo}]
Unwinding the definitions, we can identify 
$\Fun^{\lax}( \calO^{\otimes}, \calD^{\times} )$ with the full subcategory of
$$\bHom( \calO^{\otimes} \times_{ \Nerve(\FinSeg) } \Nerve(\Fintimes), \calD)$$
spanned by those functors $F$ which satisfy the following conditions:
\begin{itemize}
\item[$(1)$] For every object $C \in \calO^{\otimes}_{\seg{n}}$ and every subset $S \subseteq \nostar{n}$, 
the functor $F$ induces an equivalence
$$F(C, S ) \rightarrow \prod_{j \in S} F(C, \{j\})$$
in the $\infty$-category $\calD$.
\item[$(2)$] For every inert morphism $C \rightarrow C'$ in $\calO^{\otimes}$ which
covers $\seg{n} \rightarrow \seg{n'}$ and every subset $S \subseteq \nostar{n'}$, the induced map
$F( C, \alpha^{-1} S ) \rightarrow F( C', S)$ is an equivalence in $\calD$.
\end{itemize} 
The functor $F' = \pi \circ F$ can be described by the formula $F'( C ) = F( C, \nostar{n} )$, for each
$C \in \calO^{\otimes}_{\seg{n}}$. In other words, $F'$ can be identified with the restriction of $F$ to the full subcategory $\calC \subseteq \calO^{\otimes} \times_{ \Nerve(\FinSeg) } \Nerve(\Fintimes)$ spanned by objects of the form $(C, \nostar{n})$. 

Let $X = (C, S)$ be an object of the fiber product $\calO^{\otimes} \times_{ \Nerve(\FinSeg) } \Nerve(\Fintimes)$. Here $C \in \calO^{\otimes}_{\seg{n}}$ and $S \subseteq \nostar{n}$. We claim that the $\infty$-category $\calC_{X/}$ has an initial object. More precisely, if we choose a $p$-coCartesian morphism $\overline{\alpha}: C \rightarrow C'$ covering the map $\alpha: \seg{n} \rightarrow S_{\ast}$ given by the formula
$$ \alpha( j )  =
\begin{cases} j & \text{if } j \in S \\
\ast & \text{otherwise,} \end{cases}$$
then the induced map $\widetilde{\alpha}: (C, S) \rightarrow ( C', S)$ is an initial object of $\calC_{X/}$. It follows that every functor $F': \calC \rightarrow \calD$
admits a right Kan extension to $\calO^{\otimes} \times_{ \Nerve(\FinSeg) } \Nerve(\Fintimes)$, and that an arbitrary functor
$F: \calO^{\otimes} \times_{ \Nerve(\FinSeg) } \Nerve(\Fintimes) \rightarrow
\calD$ is a right Kan extension of $F| \calC$ if and only if
$F( \widetilde{\alpha})$ is an equivalence, for every $\widetilde{\alpha}$ defined as above.

Let $\calE$ be the full subcategory of $\Fun( \calO^{\otimes} \times_{ \Nerve(\FinSeg) } \Nerve(\Fintimes), \calD)$ spanned by those functors $F$ which satisfy the following conditions:
\begin{itemize}
\item[$(1')$] The restriction $F' = F | \calC$ is a lax Cartesian structure on $\calO^{\otimes} \simeq \calC$. 
\item[$(2')$] The functor $F$ is a right Kan extension of $F'$. 
\end{itemize}
Using Proposition \toposref{lklk}, we conclude that the restriction map
$\calE \rightarrow \Fun^{\lax}( \calO^{\otimes}, \calD)$ is a trivial fibration of simplicial sets.
To prove that $\theta$ is a trivial Kan fibration, it will suffice to show that conditions $(1)$ and $(2)$ are equivalent to conditions $(1')$ and $(2')$. 

Suppose first that $(1')$ and $(2')$ are satisfied by a functor $F$. Condition then $(1)$ follows easily; we will prove $(2)$. Choose a map $C \rightarrow C'$ covering a inert morphism $\seg{n} \rightarrow \seg{n'}$ in $\FinSeg$, and let $S \subseteq \nostar{n'}$. Define another inert morphism
$\alpha: \seg{n'} \rightarrow S_{\ast}$ by the formula
$$ \alpha( j )  =
\begin{cases} j & \text{if } j \in S \\
\ast & \text{otherwise,} \end{cases}$$
and choose a $p$-coCartesian morphism $C' \rightarrow C''$ lifting $\alpha$. Condition
$(2')$ implies that the maps $F( C, \alpha^{-1} S) \rightarrow F( C'', S)$ and
$F( C', S) \rightarrow F(C'', S)$ are equivalences in $\calD$. Using the two-out-of-three property, we deduce that the map $F(C, \alpha^{-1} S) \rightarrow F(C',S)$ is likewise an equivalence in $\calD$. This proves $(2)$.

Now suppose that $(1)$ and $(2)$ are satisfied by $F$. The implication $(2) \Rightarrow (2')$ is obvious; it will therefore suffice to verify $(1')$. Let $C$ be an object of $\calO^{\otimes}_{\seg{n}}$, and choose $p$-coCartesian morphisms $g_j: C \rightarrow C_{j}$ covering the inert morphisms
$\Colp{j}: \seg{n} \rightarrow \seg{1}$ for $1 \leq j \leq n$.
We wish to show that the induced map
$ F( C, \nostar{n} ) \rightarrow \prod_{1 \leq j \leq n} F(C_j, \nostar{1} )$ is an equivalence in $\calD$, which follows immediately from $(1)$ and $(2')$. This completes the proof that $\theta$ is a trivial Kan fibration.

Now suppose that $\calO^{\otimes}$ is a symmetric monoidal $\infty$-category.
To prove that $\theta_0$ is a trivial Kan fibration, it will suffice to show that $\theta_0$ is a pullback of $\theta$. In other words, it will suffice to show that if 
$F: \calO^{\otimes} \times_{ \Nerve(\FinSeg) } \Nerve(\Fintimes) \rightarrow
\calD$ is a functor satisfying conditions $(1)$ and $(2)$, then $F| \calC$ is a weak Cartesian structure on $\calO^{\otimes} \simeq \calC$ if and only if $F$ determines a symmetric monoidal functor from
$\calO^{\otimes}$ into $\calD^{\times}$. Let $q: \calD^{\times} \rightarrow \Nerve(\FinSeg)$ denote the projection. Using the description of the class of $q$-coCartesian morphisms provided by Proposition \ref{commonton2}, we see that the latter condition is equivalent to

\begin{itemize}
\item[$(2_{+})$] For every $p$-coCartesian morphism $\overline{\alpha}: C \rightarrow C'$ in $\calO^{\otimes}$ covering a map $\alpha: \seg{n} \rightarrow \seg{n'}$ in $\FinSeg$, and every $S \subseteq \nostar{n'}$, 
the induced map
$F( C, \alpha^{-1}(S) ) \rightarrow F( C, S)$ is an equivalence in $\calD$.
\end{itemize}

Moreover, $F| \calC$ is a weak Cartesian structure if and only if $F$ satisfies the following:

\begin{itemize}
\item[$(3')$] For each $n \geq 0$ and every $p$-coCartesian morphism
$\overline{\beta}: C \rightarrow C'$ in $\calO^{\otimes}$ lifting the map
$\beta: \seg{n} \rightarrow \seg{1}$ such that $\beta^{-1} \{ \ast \} = \{ \ast \}$,
the induced map $F( C, \nostar{n}) \rightarrow F( C', \nostar{1})$ is an equivalence in $\calD$.
\end{itemize}

It is clear that $(2_{+})$ implies $(3')$. Conversely, suppose that $(3')$ is satisfied, and let
$\overline{\alpha}$ and $S \subseteq \nostar{n}$ be as in the statement of $(2_+)$. 
Choose $p$-coCartesian morphisms $\overline{\gamma}: C \rightarrow C_0$, $\overline{\gamma'}: C' \rightarrow C'_0$, $\overline{\beta}: C'_0 \rightarrow C''_0$ covering the maps
$\gamma: \seg{n} \rightarrow ( \alpha^{-1} S)_{\ast}$, $\gamma': \seg{n'} \rightarrow S_{\ast}$, 
$\beta: \seg{n'} \rightarrow \seg{1}$ described by the formulas
$$ \beta( j )  =
\begin{cases} 1 & \text{if } j \in S \\
\ast & \text{if } j = \ast\end{cases}$$
$$ \gamma( j )  =
\begin{cases} j & \text{if } j \in \alpha^{-1} S \\
\ast & \text{otherwise}\end{cases}$$
$$ \gamma'( j )  =
\begin{cases} 1 & \text{if } j \in S \\
\ast & \text{otherwise.}\end{cases}$$
We have a commutative diagram
$$ \xymatrix{ F( C, \alpha^{-1} S) \ar[r] \ar[d] & F(C', S) \ar[d] & \\
F( C_0, \alpha^{-1} S) \ar[r]^{g} & F( C'_0, S) \ar[r]^{g'} & F(C''_0, \{0\} ) }$$
Condition $(3')$ implies that $g'$ and $g' \circ g$ are equivalences, so that
$g$ is an equivalce by the two-out-of-three property. Condition $(2')$ implies that the vertical maps are equivalences, so that the upper horizontal map is also an equivalence, as desired.
\end{proof}

\subsection{CoCartesian Monoidal Structures}\label{jilun}

Let $\calO^{\otimes}$ be an $\infty$-operad. The theory of $\calO$-monoidal $\infty$-categories is not manifestly self-dual. However, using Example \ref{interplace} and the argument of Remark \monoidref{selfop2}, we see that a $\calO$-monoidal structure on an $\infty$-category $\calC$ determines a $\calO$-monoidal structure on $\calC^{op}$, which is well-defined up to equivalence.

In particular, if $\calC$ admits finite coproducts, then the Cartesian symmetric monoidal structure on $\calC^{op}$ determines a symmetric monoidal structure on $\calC$, which we will call the {\it coCartesian symmetric monoidal structure}. It is characterized up to equivalence by the following properties:
\begin{itemize}
\item[$(1)$] The unit object $1_{\calC} \in \calC$ is initial.
\item[$(2)$] For every pair of objects $C,D \in \calC$, the canonical maps
$$ C \simeq C \otimes 1_{\calC} \rightarrow C \otimes D
\leftarrow 1_{\calC} \otimes D \simeq D$$
exhibit $C \otimes D$ as a coproduct of $C$ and $D$ in the $\infty$-category $\calC$.
\end{itemize}

The coCartesian symmetric monoidal on an $\infty$-category $\calC$ can be generalized
to the case where $\calC$ does not admit finite coproducts: however, in this case, we
obtain an $\infty$-operad rather than a symmetric monoidal $\infty$-category.

\begin{construction}\label{jiu}
We define a category $\PFinSeg$ as follows:
\begin{itemize}
\item[$(1)$] The objects of $\PFinSeg$ are pairs $( \seg{n}, i)$ where
$i \in \nostar{n}$.
\item[$(2)$] A morphism in $\PFinSeg$ from $(\seg{m},i)$ to $(\seg{n}, j)$ is
a map of pointed sets $\alpha: \seg{m} \rightarrow \seg{n}$ such that 
$\alpha(i) = j$.
\end{itemize}
Let $\calC$ be any simplicial set. We define a new simplicial set
$\calC^{\amalg}$ equipped with a map $\calC^{\amalg} \rightarrow \Nerve(\FinSeg)$ so that the
following universal property is satisfied: for every map of simplicial sets $K \rightarrow \Nerve(\FinSeg)$, we have a canonical bijection
$$ \Hom_{\Nerve(\FinSeg)}(K, \calC^{\amalg}) \simeq \Hom_{\sSet}( K \times_{ \Nerve(\FinSeg)} \Nerve(\PFinSeg), \calC).$$
\end{construction}

\begin{remark}
If $\calC$ is a simplicial set, then each fiber $\calC^{\amalg}_{\seg{n}} = \calC^{\amalg}
\times_{ \Nerve(\FinSeg)} \{ \seg{n} \}$ can be identified with $\calC^{n}$; we will henceforth invoke
these identifications implicitly.
\end{remark}

\begin{proposition}
Let $\calC$ be an $\infty$-category. Then the map $p: \calC^{\amalg} \rightarrow \Nerve(\FinSeg)$
of Construction \ref{jiu} is an $\infty$-operad.
\end{proposition}

\begin{proof}
We first show that $p$ is an inner fibration of simplicial sets. Suppose we are given a lifting problem
$$ \xymatrix{ \Lambda^n_i \ar[r]^{f_0} \ar[d] & \calC^{\amalg} \ar[d] \\
\Delta^n \ar[r] \ar@{-->}[ur]^{f} & \Nerve(\FinSeg) }$$
where $0 < i < n$. The lower horizontal map determines a sequence of maps
$\seg{k_0} \rightarrow \ldots \rightarrow \seg{k_n}$ in $\FinSeg$. Unwinding the definitions,
we see that finding the desired extension $f$ of $f_0$ is equivalent to the problem of solving
a series of extension problems
$$ \xymatrix{ \Lambda^n_{i} \ar[r]^{f^{j}_0} \ar[d] & \calC \\
\Delta^n \ar@{-->}[ur]^{f^{j}} & }$$
indexed by those elements $j \in \nostar{k_0}$ whose image in $\seg{k_n}$ belongs to $\nostar{k_n}$. 
These extensions exist by virtue of the assumption that $\calC$ is an $\infty$-category.
If $i=0$ and the map $\seg{k_0} \rightarrow \seg{k_1}$ is inert, then the same argument applies:
we conclude that the desired extension of $f$ exists provided that $n \geq 2$ and
$f^{j}_0$ carries $\Delta^{ \{0,1\} }$ to an equivalence in $\calC$.

Unwinding the definitions, we see that an object of $\calC^{\amalg}$ consists of an object
$\seg{n} \in \FinSeg$ together with a sequence of objects $(C_1, \ldots, C_n)$ in $\calC$.
A morphism $f$ from $(C_1, \ldots, C_m)$ to $(C'_1, \ldots, C'_n)$ in $\calC^{\amalg}$ consists
of a map of pointed sets $\alpha: \seg{m} \rightarrow \seg{n}$ together with a sequence
of morphisms $\{ f_{i}: C_{i} \rightarrow C'_{\alpha(i)} \}_{ i \in \alpha^{-1} \nostar{n}}$. 
The above argument shows that $f$ is $p$-coCartesian if $\alpha$ is inert and each
of the maps $f_{i}$ is an equivalence in $\calC$. In particular (taking each $f_i$ to be the
identity map), we deduce that for every object $C \in \calC^{\amalg}_{\seg{m}}$ and every
inert morphism $\alpha: \seg{m} \rightarrow \seg{n}$ in $\FinSeg$, there exists a
$p$-coCartesian morphism $C \rightarrow C'$ in $\calC^{\amalg}$ lifting $\sigma$.

Let $C = (C_1, \ldots, C_n)$ be an object of $\calC^{\amalg}_{\seg{n}}$ and choose
$p$-coCartesian morphisms $C \rightarrow C'_i$ covering $\Colp{i}$ for $1 \leq i \leq n$, corresponding to equivalences $g_i: C_i \simeq C'_i$ in $\calC$.
These morphisms determine a diagram $\overline{q}: \nostar{n}^{\triangleleft} \rightarrow \calC^{\amalg}$; we must show that $\overline{q}$ is a $p$-limit diagram. To prove this, we must show that it is possible to solve lifting problems of the form 
$$ \xymatrix{ \bd \Delta^m \star \nostar{n} \ar[r]^{f_0} \ar[d] & \calC^{\amalg} \ar[d] \\
\Delta^m \star \nostar{n} \ar[r] \ar@{-->}[ur]^{f} & \Nerve(\FinSeg) }$$
provided that $f_0|(\{m\} \star \nostar{n})$ is given by $\overline{q}$. 
Unwinding the definitions, we see that this is equivalent to solving a collection of extension problems of the form
$$ \xymatrix{ \Lambda^{m+1}_{m+1} \ar[r]^{f'_0} \ar[d] & \calC \\
\Delta^{m+1}, \ar@{-->}[ur]^{f'} & }$$
where $f'_0$ carries the final edge of $\Delta^{m+1}$ to one of the morphisms
$g_i$. This is possible by virtue of our assumption that each $g_i$ is an equivalence.

To complete the proof that $\calC^{\amalg}$ is an $\infty$-operad, it suffices to show that
for each $n \geq 0$, the functors $\Colll{i}{!} : \calC^{\amalg}_{\seg{n}} \rightarrow
\calC^{\amalg}_{\seg{1}}$ induce an equivalence $\theta: \calC^{\amalg}_{\seg{n}} \rightarrow
\prod_{1 \leq i \leq n} \calC^{\amalg}_{\seg{1}}$. In fact, we have canonical isomorphisms
of simplicial sets $\calC^{\amalg}_{\seg{n}} \simeq \calC^{n}$ which allow us to identify
$\theta$ with $\id_{\calC^{n}}$. 
\end{proof}

\begin{remark}
Unwinding the definitions, we deduce that a map
$(C_1, \ldots, C_m) \rightarrow (C'_1, \ldots, C'_n)$ in $\calC^{\amalg}$ covering
a map $\alpha: \seg{m} \rightarrow \seg{n}$ is $p$-coCartesian if and only if for
each $j \in \nostar{n}$, the underlying maps $\{ f_i: C_{i} \rightarrow C'_j \}_{\alpha(i) = j }$
exhibit $C'_{j}$ as a coproduct of $\{ C_i \}_{\alpha(i) = j}$ in the $\infty$-category $\calC$.
It follows that $\calC^{\amalg}$ is a symmetric monoidal $\infty$-category if and only if
$\calC$ admits finite coproducts. If this condition is satisfied, then $\calC^{\amalg}$ determines a
coCartesian symmetric monoidal structure on $\calC$ and is therefore determined by $\calC$ up to essentially unique equivalence. We will see that the situation is similar even if $\calC$ does not admit finite coproducts.
\end{remark}

\begin{example}\label{stuly}
The projection map $\Nerve(\PFinSeg) \rightarrow \Nerve(\FinSeg)$ induces a
canonical map $\calC \times \Nerve(\FinSeg) \rightarrow \calC^{\amalg}$.
If $\calC = \Delta^0$, this map is an isomorphism (so that $\calC^{\amalg}$ is the
commutative $\infty$-operad $\Nerve(\FinSeg)$). For any $\infty$-operad $\calO^{\otimes}$,
we obtain a map $\calC \times \calO^{\otimes} \rightarrow \calC^{\amalg} \times_{ \Nerve(\FinSeg)} \calO^{\otimes}$ which determines a functor
$\calC \rightarrow \Alg_{\calO}(\calA)$, where $\calA^{\otimes}$ is the $\infty$-operad
$\calC^{\amalg} \times_{ \Nerve(\FinSeg)} \calO^{\otimes}$.
\end{example}

Our main goal in this section is to prove the following result, which characterizes
$\infty$-operads of the form $\calC^{\amalg}$ by a universal mapping property.

\begin{theorem}\label{kinema}
Let $\calC$ be an $\infty$-category and let $\calO^{\otimes}$ and $\calD^{\otimes}$ be $\infty$-operads. The construction of Example \ref{stuly} induces a trivial Kan fibration
$$ \theta: \Fun^{\lax}( \calC^{\amalg} \times_{ \Nerve(\FinSeg)} \calO^{\otimes}, \calD^{\otimes})
\rightarrow \Fun( \calC, \Alg_{\calO}(\calD)).$$
In particular (taking $\calO^{\otimes}$ to be the commutative operad), we have a trivial Kan fibration 
$\Fun^{\lax}( \calC^{\amalg}, \calD^{\otimes}) \rightarrow \Fun( \calC, \CAlg(\calD))$. 
\end{theorem}

\begin{corollary}\label{slamm}
Let $\calC$ and $\calD$ be $\infty$-categories. Then the restriction map
$\theta: \Fun^{\lax}( \calC^{\amalg}, \calD^{\amalg}) \rightarrow \Fun( \calC, \calD)$
is a trivial Kan fibration.
\end{corollary}

\begin{proof}
We can factor the map $\theta$ as a composition
$$ \Fun^{\lax}( \calC^{\amalg}, \calD^{\amalg}) \stackrel{\theta'}{\rightarrow}
\Fun(\calC, \CAlg( \calD)) \stackrel{\theta''}{\rightarrow} \Fun(\calC, \calD)$$
there $\theta'$ is a trivial Kan fibration by Theorem \ref{kinema} and
$\theta''$ is a trivial Kan fibration by Corollary \ref{kime2}.
\end{proof}

\begin{definition}
We will say that an $\infty$-operad $\calO^{\otimes}$ is {\it coCartesian} if it is
equivalent to $\calC^{\amalg}$, for some $\infty$-category $\calC$.
\end{definition}

\begin{corollary}
The forgetful functor $\phi: \Cat_{\infty}^{\lax, \CommOp} \rightarrow \Cat_{\infty}$ induces
an equivalence from the full subcategory of $\Cat_{\infty}^{\lax, \CommOp}$ spanned by the
coCartesian $\infty$-operads to $\Cat_{\infty}$.
\end{corollary}

\begin{proof}
Let $\phi_0$ denote the restriction of $\phi$ to the full subcategory spanned by the
coCartesian $\infty$-operads. Since every $\infty$-category $\calC$ can be identified with the underlying $\infty$-category of $\calC^{\amalg}$, we deduce that $\phi_0$ is essentially surjective.
Corollary \ref{slamm} guarantees that $\phi_0$ is fully faithful.
\end{proof}

We now turn to the proof of Theorem \ref{kinema}. We will need a few preliminaries.

\begin{lemma}\label{stealthwick}
Let $S$ be a finite set (regarded as a discrete simplicial set), let $v$ denote the cone
point of $S^{\triangleleft}$, and suppose
we are given coCartesian fibrations $p: X \rightarrow S^{\triangleleft}$
and $q: Y \rightarrow S^{\triangleleft}$ which induce categorical equivalences
$$ X_{v} \simeq \prod_{s \in S} X_{s} \quad \quad Y_{v} \simeq \prod_{s \in S} Y_{s}.$$
Let $\Fun'_{ S^{\triangleleft}}( X, Y)$ denote the full subcategory of
$\Fun_{ S^{\triangleleft}}( X,Y)$ spanned by those maps which carry
$p$-coCartesian morphisms to $q$-coCartesian morphisms. Then
the restriction functor
$$\Fun'_{ S^{\triangleleft}}( X, Y) \rightarrow \Fun_{S}( X \times_{S^{\triangleleft}} S,
Y \times_{ S^{\triangleleft}} S)$$
is a trivial Kan fibration.
\end{lemma}

\begin{proof}
In view of Proposition \toposref{lklk}, it will suffice to prove the following:
\begin{itemize}
\item[$(1)$] A functor $F \in \Fun_{ S^{\triangleleft}}( X,Y)$ belongs to 
$\Fun'_{ S^{\triangleleft}}( X, Y)$ if and only if $F$ is a $q$-right Kan extension of
$F | X$.
\item[$(2)$] Every map $F_0 \in \Fun_{S}( X \times_{S^{\triangleleft}} S,
Y \times_{ S^{\triangleleft}} S)$ can be extended to a map
$F \in \Fun_{ S^{\triangleleft}}( X,Y)$ satisfying the equivalent conditions of $(1)$.
\end{itemize}

To prove $(1)$, consider an arbitrary object $x$ of $X_{v}$, and choose
$p$-coCartesian morphisms $f_{s}: x \rightarrow x_{s}$ to objects $x_{s} \in X_{s}$ for
$s \in S$. We note that
the inclusion $\{ f_s \}_{s \in S} \hookrightarrow ((X \times_{S^{\triangleleft}} S)_{x/})^{op}$ is cofinal. It follows that a functor $F$ as in $(1)$
is a $q$-right Kan extension of $F_0$ at $x$ if and only if the maps
$F(f_s)$ exhibit $F(x)$ as a $q$-product of the objects $F_0(x_{s})$.
This is equivalent to the requirement that each $F(f_s)$ is $q$-coCartesian.
This proves the ``if'' direction of $(1)$; the converse follows from same argument
together with the observation that every $p$-coCartesian morphism $f: x \rightarrow x_{s}$ in
$X$ can be completed to a collection of $p$-coCartesian morphisms $\{ f_{s'}: x \rightarrow x_{s'} \}_{s' \in S}$. 

To prove $(2)$, it suffices (by Lemma \toposref{kan2}) to show that for every $x \in X_{v}$, 
the diagram
$(X \times_{ S^{\triangleleft} } S)_{x/} \rightarrow Y$
induced by $F_0$ can be extended to a $q$-limit diagram covering the projection map
$(X \times_{S^{\triangleleft}} S)_{x/}^{\triangleleft} \rightarrow S^{\triangleleft}$. Let
$\{ f_{s}: x \rightarrow x_{s} \}_{s \in S}$ be as above. We must show that there
exists an object $y \in Y_{v}$ equipped with morphisms
$y \rightarrow F_0( x_s)$ for $s \in S$, which exhibit $y$ as a $q$-product of
$\{ F_0(x_s) \}_{s \in S}$. It suffices to choose $y$ to be any preimage of
$\{ F_0(x_s) \}_{s \in S}$ under the equivalence $Y_{v} \simeq \prod_{s \in S} Y_{s}$.
\end{proof}

\begin{lemma}\label{sleuth}
Let $S$ be a finite set (regarded as a discrete simplicial set), let $n > 0$, and
suppose we are given inner fibrations $p: X \rightarrow \Delta^n \star S$ and
$q: Y \rightarrow \Delta^n \star S$. For every simplicial subset $K \subseteq \Delta^n$, let $X_{K}$ denote the fiber product $X \times_{ \Delta^n \star S} (K \star S)$, and define
$Y_{K}$ similarly. Assume that the maps $X_{ \{n\} } \rightarrow S^{\triangleleft}$ and
$Y_{ \{n\} } \rightarrow S^{\triangleleft}$ satisfy the hypotheses of Lemma \ref{stealthwick},
and for $\{n\} \subseteq K$ define
$\Fun'_{ K \star S}( X_{K}, Y_{K} )$ to be the fiber product $$\Fun_{ K \star S}( X_K, Y_K) \times_{
\Fun_{ \{n\} \star S}( X_{ \{n\} }, Y_{ \{n\} })} \Fun'_{ \{n\} \star S}( X_{ \{n\}}, Y_{ \{n\} }).$$
Then the map $$\theta: \Fun'_{ \Delta^n \star S}( X, Y) \rightarrow
\Fun'_{ \bd \Delta^n \star S}( X_{ \bd \Delta^n}, Y_{ \bd \Delta^n} )$$
is a trivial Kan fibration.
\end{lemma}

\begin{proof}
The proof proceeds by induction on $n$. We observe that $\theta$ is evidently a categorical fibration;
to prove that it is a trivial Kan fibration, it will suffice to show that $\theta$ is a categorical equivalence.
Let $\theta''$ denote the composition
$$ \Fun'_{ \Delta^n \star S}( X,Y) \stackrel{\theta}{\rightarrow}
\Fun'_{ \bd \Delta^n \star S}( X_{ \bd \Delta^n}, Y_{\bd \Delta^n} )
\stackrel{\theta'}{\rightarrow} \Fun_{\Delta^{n-1} \star S}( X_{ \Delta^{n-1}}, Y_{ \Delta^{n-1}}).$$
By a two-out-of-three argument, it will suffice to show that $\theta'$ and $\theta''$ are trivial
Kan fibrations. The map $\theta'$ is a pullback of the composition
$$ \Fun'_{K \star S}( X_{K}, Y_{K})
\stackrel{\phi}{\rightarrow} \Fun'_{ \{n\} \star S}( X_{ \{n\}}, Y_{ \{n\} }) \stackrel{\phi'}{\rightarrow}
\Fun_{S}( X_{\emptyset}, Y_{\emptyset}),$$
where $K = \Delta^{n-1} \coprod \{n\} \subseteq \Delta^n$. It follows from iterated application
of the inductive hypothesis that $\phi$ is a trivial Kan fibration, and it follows from Lemma
\ref{stealthwick} that $\phi'$ is a trivial Kan fibration. Consequently, to complete the proof, it will
suffice to show that $\theta''$ is a trivial Kan fibration.
In view of Proposition \toposref{lklk}, it will suffice to prove the following:
\begin{itemize}
\item[$(1)$] A map $F \in \Fun_{\Delta^n \star S}( X, Y)$
is a $q$-right Kan extension of $F_0 = F | X_{ \Delta^{n-1} }$ if and only if
it belongs to $\Fun'_{ \Delta^n \star S}(X, Y)$.
\item[$(2)$] Every map $F_0 \in \Fun_{ \Delta^{n-1} \star S}( X_{ \Delta^{n-1}}, Y_{ \Delta^{n-1} })$
admits an extension $F \in \Fun_{ \Delta^n \star S}(X,Y)$ satisfying the equivalent conditions
of $(1)$.
\end{itemize}
These assertions follow exactly as in the proof of Lemma \ref{stealthwick}.
\end{proof}

\begin{proof}[Proof of Theorem \ref{kinema}]
The map $\theta$ is evidently a categorical fibration; to prove that it is a trivial Kan fibration, it suffices to show that it is a categorical equivalence. We are therefore free to replace $\calC$ by any
$\infty$-category equivalent to $\calC$; using Proposition \toposref{minimod} we may reduce to the
case where the $\infty$-category $\calC$ is minimal.

We now proceed with the proof by defining a rather elaborate filtration of
$\calC^{\amalg}$. To describe this filtration, we need to introduce a bit of terminology.
Recall that a morphism $\overline{\alpha}$ from $(C_1, \ldots, C_m)$ to
$(C'_1, \ldots, C'_n)$ in $\calC^{\amalg}$ consists of a map of pointed sets
$\alpha: \seg{m} \rightarrow \seg{n}$ together with a morphism $f_i: C_i \rightarrow
C'_{\alpha(i)}$ for each $i \in \alpha^{-1} \nostar{n}$. We will say that $\overline{\alpha}$ is
{\it quasidegenerate} if each of the morphisms $f_i$ is a degenerate edge of $\calC$.

Let $\sigma$ be an $n$-simplex of $\calC^{\amalg}$ given by a sequence of morphisms
$$ \sigma(0) \stackrel{ \overline{\alpha}(1)}{\rightarrow} \sigma(1)
\rightarrow \cdots \stackrel{ \overline{\alpha}(n)}{\rightarrow} \sigma(n)$$
and let
$$ \seg{k_0} \stackrel{\alpha(1)}{\rightarrow} \seg{k_1} \rightarrow \cdots \stackrel{\alpha(n)}{\rightarrow} \seg{k_n}$$
be the underlying $n$-simplex of $\Nerve(\FinSeg)$. We will say that $\sigma$ is
{\it closed} if $k_n = 1$, and {\it open} otherwise. If $\sigma$ is closed, we define the
{\it tail length} of $\sigma$ to be the largest integer $m$ such that the maps 
$\alpha(k)$ are isomorphisms for $n-m < k \leq n$.
We will denote the tail length of $\sigma$ by $t(\sigma)$.
We define the {\it break point} of a closed simplex $\sigma$ to be smallest nonnegative integer $m$ such that the maps $\overline{\alpha}(k)$ are
active and quasidegenerate for $m < k \leq n- t(\sigma)$. We will denote the break point of $\sigma$ by
$b(\sigma)$. Let $S = \coprod_{0 \leq i \leq n} \nostar{k_i}$. We will say that
an element $j \in \nostar{k_i} \subseteq S$ is a {\it leaf} if $i=0$ or if $j$ does not lie in the image
of the map $\alpha(i)$, and we will say that $j$ is a {\it root} if $i=n$ or if $\alpha(i+1)(j) = \ast$.
We define the {\it complexity} $c(\sigma)$ of $\sigma$ to be $2l - r$, where $l$ is the number of
leaves of $\sigma$ and $r$ is the number of roots of $\sigma$. We will say that $\sigma$
is {\it flat} if it belongs to the image of the embedding $\Nerve(\FinSeg) \times \calC \rightarrow \calC^{\amalg}$. Since $\calC$ is minimal, Proposition \toposref{minstrict} implies that if
$\sigma$ is closed and $b(\sigma) = 0$, then $\sigma$ is flat.

We now partition the nondegenerate, nonflat simplices of $\calC^{\amalg}$ into six groups:

\begin{itemize}
\item[$(A)$] An $n$-dimensional nonflat nondegenerate simplex $\sigma$ of
$\calC^{\amalg}$ belongs to $A$ if $\sigma$ is closed and
the map $\alpha(b(\sigma))$ is not inert.

\item[$(A')$] An $n$-dimensional nonflat nondegenerate simplex $\sigma$ of
$\calC^{\amalg}$ belongs to $A'$ if $\sigma$ is closed, $b(\sigma) < n - t(\sigma)$,
and the map $\alpha(b(\sigma))$ is inert.

\item[$(B)$] An $n$-dimensional nonflat nondegenerate simplex $\sigma$ of
$\calC^{\amalg}$ belongs to $B$ if $\sigma$ is closed, $b(\sigma) = n - t(\sigma)$,
the map $\alpha( b(\sigma) )$ is inert, and $\overline{\alpha}( b(\sigma) )$ is 
not quasidegenerate.

\item[$(B')$] An $n$-dimensional nonflat nondegenerate simplex $\sigma$ of
$\calC^{\amalg}$ belongs to $B$ if $\sigma$ is closed, $b(\sigma) = n - t(\sigma) < n$,
the map $\alpha(b(\sigma))$ is inert, and $\overline{\alpha}( b(\sigma))$ is quasidegenerate.
\item[$(C)$] An $n$-dimensional nonflat nondegenerate simplex $\sigma$ of
$\calC^{\amalg}$ belongs to $C$ if it is open.

\item[$(C')$] An $n$-dimensional nonflat nondegenerate simplex $\sigma$ of
$\calC^{\amalg}$ belongs to $C'$ is it is closed, $b(\sigma) = n-t(\sigma) = n$, the
map $\alpha(b(\sigma))$ is inert, and $\overline{\alpha}(b(\sigma))$ is quasidegenerate.
\end{itemize}

If $\sigma$ belongs to $A'$, $B'$, or $C'$, then we define the {\it associate} $a(\sigma)$ of
$\sigma$ to be the face of $\sigma$ opposite the $b(\sigma)$th vertex. It follows from
Proposition \toposref{minstrict} that $a(\sigma)$ belongs to $A$ if $\sigma \in A'$,
$B$ if $\sigma \in B'$, and $C$ if $\sigma \in C'$. In this case, we will say that
$\sigma$ is an {\it associate} of $a(\sigma)$. We note that every simplex belonging to
$A$ or $B$ has a unique associate, while a simplex $\sigma$ of $C$ has precisely $k$ associates,
where $\seg{k}$ is the image of the final vertex of $\sigma$ in $\Nerve(\FinSeg)$.

For each $n \geq 0$, let $K(n) \subseteq \calC^{\amalg}$ be the simplicial subset
generated by those nondegenerate simplices which are either flat, have dimension $\leq n$, 
or have dimension $n+1$ and belong to either $A'$, $B'$, or $C'$. We observe that 
$K(0)$ is generated by $\calC \times \Nerve(\FinSeg)$ together with the collection of
$1$-simplices belonging to $C'$. Let
$\calE(0)$ denote the full subcategory of $\bHom_{\Nerve(\FinSeg)}( K(0) \times_{\Nerve(\FinSeg)}
\calO^{\otimes}, \calD^{\otimes} )$ spanned by those maps $F$ with the following properties:
\begin{itemize}
\item[$(i)$] The restriction of $F$ to $\calC \times \Nerve(\FinSeg)$ determines an object of
$\Fun( \calC, \Alg_{\calA}(\calD))$.
\item[$(ii)$] Let $f$ be an edge of $K(0) \times_{ \Nerve(\FinSeg)} \calO^{\otimes}$ whose
image in $\calO^{\otimes}$ is inert and whose image in $K(0)$ belongs to $C'$. Then
$F(f)$ is an inert morphism in $\calD^{\otimes}$.
\end{itemize}

To complete the proof, it will suffice to show that the restriction maps
$$  \Fun^{\lax}( \calC^{\amalg} \times_{ \Nerve(\FinSeg)} \calO^{\otimes}, \calD^{\otimes})
\stackrel{\theta'}{\rightarrow} \calE(0) \stackrel{\theta''}{\rightarrow} \Fun( \calC, \Alg_{\calO}(\calD))$$
are trivial Kan fibrations. For the map $\theta''$, this follows from repeated application
of Lemma \ref{stealthwick}. To prove that $\theta'$ is a trivial Kan fibration, we define 
$\calE(n)$ to be the full subcategory of 
$\bHom_{\Nerve(\FinSeg)}( K(n) \times_{\Nerve(\FinSeg)}
\calO^{\otimes}, \calD^{\otimes} )$ spanned by those functors $F$ whose restriction
to $K(0) \times_{ \Nerve(\FinSeg)} \calO^{\otimes}$ belongs to $\calE(0)$.
It now suffices to prove the following:

\begin{itemize}
\item[$(a)$] A functor $F \in \bHom_{ \Nerve(\FinSeg)}( \calC^{\amalg} \times_{ \Nerve(\FinSeg)}
\calO^{\otimes}, \calD^{\otimes})$ is a map of $\infty$-operads if and only if
$F$ satisfies conditions $(i)$ and $(ii)$. Consequently,
$\Fun^{\lax}( \calC^{\amalg} \times_{ \Nerve(\FinSeg)}
\calO^{\otimes}, \calD^{\otimes})$ can be identified with the inverse limit of the tower
$$ \cdots \rightarrow \calE(2) \rightarrow \calE(1) \rightarrow \calE(0).$$

\item[$(b)$] For $n > 0$, the restriction map $\calE(n) \rightarrow \calE(n-1)$ is a trivial Kan fibration.
\end{itemize}

We now prove $(a)$. The ``only if'' direction is obvious. For the converse, suppose that
$$F \in \bHom_{ \Nerve(\FinSeg)}( \calC^{\amalg} \times_{ \Nerve(\FinSeg)}
\calO^{\otimes}, \calD^{\otimes})$$ satisfies conditions $(i)$ and $(ii)$ above.
We wish to prove that $F$ preserves inert morphisms. Let
$f: X \rightarrow X'$ be an inert morphism in $\calC^{\amalg} \times_{ \Nerve(\FinSeg)} \calO^{\otimes}$
covering the map $\alpha: \seg{m} \rightarrow \seg{m'}$ in $\FinSeg$; we must show that
$F(f)$ is inert. If $m' = 1$, then we can identify the image of $X$ in $\calC^{\amalg}$ with a sequence of objects $(C_1, \ldots, C_m)$ in $\calC$, the image of $X'$ in $\calC^{\amalg}$ with an object $C' \in \calC$, and the image of $f$ in $\calC^{\amalg}$ with a morphism
$\overline{f}: C_i \rightarrow C'$ for some $i \in \nostar{m}$. Replacing $f$ by an equivalent morphism
if necessary, we can assume that $\overline{f}$ is degenerate. Then the image of
$f$ in $\calC^{\amalg}$ belongs to $K(0)$ and the desired result follows from 
either $(i)$ (if $C_1 = \cdots = C_m$) or $(ii)$ (otherwise).

In the general case, we can choose inert morphisms $f_i: X' \rightarrow X'_i$ in
$\calC^{\amalg} \times_{ \Nerve(\FinSeg)} \calO^{\otimes}$ covering
$\Colp{i}$ for $1 \leq i \leq m$. The argument above shows that
the morphisms $F( f_i \circ f) \simeq F(f_i) \circ F(f)$ and $F(f_i)$ are inert. 
Since $\calD^{\otimes}$ is an $\infty$-operad, it follows that $F(f)$ is inert, which completes the proof of $(a)$.

We now prove $(b)$. For each integer $c \geq 0$, let $K(n,c)$ denote the simplicial subset
$K(n)$ spanned by those simplices which either belong to $K(n-1)$ or have complexity $\leq c$.
Let $\calE(n,c)$ denote the full subcategory of
$\bHom_{\Nerve(\FinSeg)}( K(n,c) \times_{\Nerve(\FinSeg)}
\calO^{\otimes}, \calD^{\otimes} )$
spanned by those maps $F$ whose restriction to $K(0)$ satisfies conditions $(i)$ and $(ii)$.
We have a tower of simplicial sets
$$ \cdots \rightarrow \calE(n,2) \rightarrow \calE(n,1) \rightarrow \calE(n,0) \simeq \calE(n-1)$$
with whose inverse limit can be identified with $\calE(n)$. It will therefore suffice to show that for
each $c > 0$, the restriction map $\calE(n,c) \rightarrow \calE(n,c-1)$ is a trivial Kan fibration.

We now further refine our filtration as follows. Let $K(n,c)_{A}$ denote the 
simplicial subset of $K(n,c)$ spanned by $K(n,c-1)$ together with those simplices
of $K(n,c)$ which belong to $A$ or $A'$ and let $K(n,c)_{B}$ denote the simplicial subset
of $K(n,c)$ spanned by $K(n,c-1)$ together with those simplices which belong to
$A$, $A'$, $B$, or $B'$. Let
$\calE(n,c)_{A}$ denote the full subcategory of
$\bHom_{\Nerve(\FinSeg)}( K(n,c)_{A} \times_{\Nerve(\FinSeg)}
\calO^{\otimes}, \calD^{\otimes} )$
spanned by those functors which satisfy $(i)$ and $(ii)$, and define
$\calE(n,c)_{B}$ similarly. To complete the proof, it will suffice to prove the following:

\begin{itemize}
\item[$(A)$] The restriction map $\calE(n,c)_{A} \rightarrow \calE(n,c-1)$ is a trivial Kan fibration.
To prove this, it suffices to show that the inclusion
$K(n,c-1) \times_{ \Nerve(\FinSeg)} \calO^{\otimes} \rightarrow K(n,c)_{A} \times_{ \Nerve(\FinSeg)} \calO^{\otimes}$ is a categorical equivalence.
Let $A_{n,c}$ denote the collection of all $n$-simplices belonging to $A$ having
complexity $c$. Choose a well-ordering of $A_{n,c}$ with the following properties:
\begin{itemize}
\item If $\sigma, \sigma' \in A_{n,c}$ and $t(\sigma) < t(\sigma')$, then $\sigma < \sigma'$.
\item If $\sigma, \sigma' \in A_{n,c}$, $t(\sigma) = t(\sigma')$, and $b(\sigma) < b(\sigma')$, then
$\sigma < \sigma'$.
\end{itemize}
For each $\sigma \in A_{n,c}$, let $K(n,c)_{\leq \sigma}$ denote the simplicial subset
of $K(n,c)$ generated by $K(n,c-1)$, all simplices $\tau \leq \sigma$ in $A_{n,c}$, and
all of the simplices in $A'$ which are associated to simplices of the form $\tau \leq \sigma$. Define $K(n,c)_{< \sigma}$ similarly. Using transfinite induction on $A_{n,c}$, we are reduced to proving
that for each $\sigma \in A_{n,c}$, the inclusion
$$ i: K(n,c)_{< \sigma} \times_{ \Nerve(\FinSeg)} \calO^{\otimes} \rightarrow K(n,c)_{\leq \sigma} \times_{ \Nerve(\FinSeg)} \calO^{\otimes}$$
is a categorical equivalence. Let $\sigma': \Delta^{n+1} \rightarrow \calC^{\amalg}$ be the
unique $(n+1)$-simplex of $A'$ associated to $\sigma$. We observe that $\sigma'$ determines
a pushout diagram
$$ \xymatrix{ \Lambda^{n+1}_{b(\sigma')} \ar[r] \ar[d] & K(n,c)_{< \sigma} \ar[d] \\
\Delta^{n+1} \ar[r] & K(n,c)_{\leq \sigma}. }$$
Consequently, the map $i$ is a pushout of an inclusion
$$ i': \Lambda^{n+1}_{b(\sigma')} \times_{ \Nerve(\FinSeg)} \calO^{\otimes}
\rightarrow \Delta^{n+1}_{b(\sigma)} \times_{ \Nerve(\FinSeg)} \calO^{\otimes}.$$
Since the Joyal model structure is left proper, it suffices to show that $i'$ is a
categorical equivalence, which follows from Lemma \ref{skope}.

\item[$(B)$] The map $\calE(n,c)_{B} \rightarrow \calE(n,c)_{A}$ is a trivial Kan fibration.
To prove this, it suffices to show that the inclusion 
$K(n,c)_{A} \times_{ \Nerve(\FinSeg)} \calO^{\otimes}
\subseteq K(n,c)_{B} \times_{ \Nerve(\FinSeg)} \calO^{\otimes}$
is a categorical equivalence of simplicial sets. Let $B_{n,c}$ denote the collection
of all $n$-simplices belonging to $B$ having complexity $c$. Choose
a well-ordering of $B_{n,c}$ such that the function $\sigma \mapsto t(\sigma)$ is nonstrictly
decreasing. For each $\sigma \in B_{n,c}$, we let
$K(n,c)_{\leq \sigma}$ be the simplicial subset of $K(n,c)$ generated by
$K(n,c)_{A}$, those simplices $\tau$ of $B_{n,c}$ such that $\tau \leq \sigma$, and
those simplices of $B'$ which are associated to $\tau \leq \sigma \in B_{n,c}$.
Let $K(n,c)_{< \sigma}$ be defined similarly. Using a induction on $B_{n,c}$, we
can reduce to the problem of showing that each of the inclusions
$$ K(n,c)_{< \sigma} \times_{ \Nerve(\FinSeg)} \calO^{\otimes} \rightarrow
K(n,c)_{\leq \sigma} \times_{ \Nerve(\FinSeg)} \calO^{\otimes}$$
is a categorical equivalence. 
Let $\sigma': \Delta^{n+1} \rightarrow \calC^{\amalg}$ be the
unique $(n+1)$-simplex of $B'$ associated to $\sigma$. We observe that $\sigma'$ determines
a pushout diagram
$$ \xymatrix{ \Lambda^{n+1}_{b(\sigma')} \ar[r] \ar[d] & K(n,c)_{< \sigma} \ar[d] \\
\Delta^{n+1} \ar[r] & K(n,c)_{\leq \sigma}. }$$
Consequently, the map $i$ is a pushout of an inclusion
$$ i': \Lambda^{n+1}_{b(\sigma')} \times_{ \Nerve(\FinSeg)} \calO^{\otimes}
\rightarrow \Delta^{n+1}_{b(\sigma)} \times_{ \Nerve(\FinSeg)} \calO^{\otimes}.$$
Since the Joyal model structure is left proper, it suffices to show that $i'$ is a
categorical equivalence, which follows from Lemma \ref{skope}.

\item[$(C)$] The map $\calE(n,c) \rightarrow \calE(n,c)_{B}$ is a trivial
Kan fibration. To prove this, let $C_{n,c}$ denote the subset of
$C$ consisting of $n$-dimensional simplices of complexity $c$, and choose
a well-ordering of $C_{n,c}$. For each $\sigma \in C_{n,c}$, let
$K(n,c)_{\leq \sigma}$ denote the simplicial subset of $K(n,c)$ generated
by $K(n,c)_{B}$, those simplices $\tau \in C_{n,c}$ such that $\tau \leq \sigma$,
and those simplices of $C'$ which are associated to $\tau \in C_{n,c}$ with
$\tau \leq \sigma$. Let $\calE(n,c)_{\leq \sigma}$ be the
full subcategory of $\bHom_{\Nerve(\FinSeg)}( K(n,c)_{\leq \sigma} \times_{\Nerve(\FinSeg)}
\calO^{\otimes}, \calD^{\otimes} )$ spanned by those maps $F$ satisfying
$(i)$ and $(ii)$. We define $K(n,c)_{< \sigma}$ and
$\calE(n,c)_{< \sigma}$ similarly. Using transfinite induction on $C_{n,c}$, we
are reduced to the problem of showing that for each $\sigma \in C_{n,c}$, the map
$\psi: \calE(n,c)_{\leq \sigma} \rightarrow \calE(n,c)_{< \sigma}$ is a trivial Kan fibration.

Let $\seg{k}$ denote the image of the final vertex of $\sigma$ in $\FinSeg$.
For $1 \leq i \leq k$, let $\sigma_i \in C'$ denote the unique $(n+1)$-simplex
associated to $\sigma$ such that $\sigma_{i}$ carries $\Delta^{ \{n, n+1\} }$ to
the morphism $\Colp{i}$ in $\FinSeg$. The simplices
$\{ \sigma_{i} \}_{1 \leq i \leq k}$ determine a map of simplicial sets
$\Delta^n \star \nostar{k} \rightarrow K(n,c)_{\leq \sigma}$. We have a pushout diagram of simplicial sets
$$ \xymatrix{ (\bd \Delta^n) \star \nostar{k} \ar[r] \ar[d] & K(n,c)_{< \sigma} \ar[d] \\
\Delta^n \star \nostar{k} \ar[r] & K(n,c)_{\leq \sigma}. }$$
The map $\psi$ fits into a pullback diagram
$$ \xymatrix{ \calE(n,c)_{\leq \sigma} \ar[r] \ar[d]^{\psi} &
\bHom'_{ \Delta^n \star \nostar{k}}( (\Delta^n \star \nostar{k}) \times_{ \Nerve(\FinSeg)} \calO^{\otimes},
\calD^{\otimes} \times_{ \Nerve(\FinSeg)} ( \Delta^n \star \nostar{k})) \ar[d]^{\psi'} \\
 \calE(n,c)_{< \sigma} \ar[r] & \bHom'_{\Nerve(\FinSeg)}( ( \bd \Delta^n \star \nostar{k})
\times_{ \bd \Delta^n \star \nostar{k}} \calO^{\otimes}, \calD^{\otimes} 
\times_{ \Nerve(\FinSeg)} ( \bd \Delta^n \star \nostar{k}))}$$
where $\psi'$ is the trivial fibration of Lemma \ref{sleuth}. It follows that $\psi$
is a trivial fibration, as desired.
\end{itemize}
\end{proof}

\subsection{Monoidal Envelopes}\label{monenv}

By definition, every symmetric monoidal $\infty$-category can be regarded
as an $\infty$-operad, and every symmetric monoidal functor as a map of $\infty$-operads.
This observation provides a forgetful functor from the $\infty$-category
$\Cat_{\infty}^{\CommOp, \otimes}$ of symmetric monoidal $\infty$-categories
to the $\infty$-category $\Cat_{\infty}^{\CommOp, \lax}$ of $\infty$-operads.
Our goal in this section is to construct a left adjoint to this forgetful functor. More
generally, we will give a construction which converts an arbitrary
fibration of $\infty$-operads $\calC^{\otimes} \rightarrow \calO^{\otimes}$ into
a coCartesian fibration of $\infty$-operads $\Envv_{\calO}(\calC)^{\otimes} \rightarrow \calO^{\otimes}$.

\begin{definition}\label{kapit}
Let $\calO^{\otimes}$ be an $\infty$-operad. We let $\Act(\calO^{\otimes})$ denote the full subcategory of
$\Fun( \Delta^1, \calO^{\otimes})$ spanned by the active morphisms.
Suppose that $p: \calC^{\otimes} \rightarrow \calO^{\otimes}$ is a fibration of $\infty$-operads. We let
$\Envv_{\calO}(\calC)^{\otimes}$ denote the fiber product
$$ \calC^{\otimes} \times_{ \Fun( \{0\}, \calO^{\otimes}) } \Act( \calO^{\otimes}).$$
We will refer to $\Envv_{\calO}(\calC)^{\otimes}$ as the {\it $\calO$-monoidal envelope} of $\calC^{\otimes}$. In the special case where $\calO^{\otimes}$ is the commutative $\infty$-operad, we will
denote $\Envv_{\calO}(\calC)^{\otimes}$ simply by $\Envv(\calC)^{\otimes}$.
\end{definition}

\begin{remark}
Let $\calC^{\otimes}$ be an $\infty$-operad. 
Evaluation at $\{1\} \subseteq \Delta^1$ induces a map
$\Envv(\calC)^{\otimes} \rightarrow \Nerve(\FinSeg)$. We let $\Envv(\calC)$ denote the fiber
$\Envv(\calC)^{\otimes}_{\seg{1}}$. Unwinding the definitions, we deduce that
$\Envv(\calC)$ can be identified with the subcategory $\calC^{\otimes}_{\acti} \subseteq \calC^{\otimes}$ spanned by all objects and active morphisms between them.
\end{remark}

We will defer the proof of the following basic result until the end of this section:

\begin{proposition}\label{luckma}
Let $p: \calC^{\otimes} \rightarrow \calO^{\otimes}$ be a fibration of $\infty$-operads. The functor
$q: \Envv_{\calO}( \calC)^{\otimes} \rightarrow \calO^{\otimes}$ given by evaluation at $\{1\} \in \Delta^1$
exhibits $\Envv_{\calO}(\calC)^{\otimes}$ as a $\calO$-monoidal $\infty$-category.
\end{proposition}

For any $\infty$-operad $\calO^{\otimes}$, the diagonal embedding
$\calO^{\otimes} \rightarrow \Fun( \Delta^1, \calO^{\otimes})$ factors through $\Act(\calO^{\otimes})$.
Pullback along this embedding induces an inclusion
$\calC^{\otimes} \subseteq \Envv_{\calO}(\calC)^{\otimes}$ for any fibration of $\infty$-operads
$\calC^{\otimes} \rightarrow \calO^{\otimes}$. It follows from Proposition \ref{luckma} and Lemma \ref{placeman} (below) that this inclusion is a map of $\infty$-operads.
The terminology ``$\calO$-monoidal envelope'' is justified by the following result, whose proof we will also defer:

\begin{proposition}\label{tull}
Let $p: \calC^{\otimes} \rightarrow \calO^{\otimes}$ be a fibration of $\infty$-operads, and let
$q: \calD^{\otimes} \rightarrow \calO^{\otimes}$ be a $\calO$-monoidal $\infty$-category.
The inclusion $i: \calC^{\otimes} \subseteq \Envv_{\calO}( \calC)^{\otimes}$ induces 
an equivalence of $\infty$-categories
$$ \Fun^{\otimes}_{\calO}( \Envv_{\calO}(\calC)^{\otimes}, \calD^{\otimes} ) \rightarrow \Fun^{\lax}_{\calO}( \calC^{\otimes}, \calD^{\otimes}).$$
\end{proposition}

Taking $\calO^{\otimes}$ to be the commutative $\infty$-operad, we deduce the following:

\begin{corollary}\label{letm}
Let $\calC^{\otimes}$ be an $\infty$-operad. Then there is a canonical symmetric monoidal structure on the $\infty$-category $\calC^{\otimes}_{\acti}$ of active morphisms.
\end{corollary}

\begin{remark}\label{opl}
If $\calC^{\otimes}$ is an $\infty$-operad, we let
$\oplus: \calC^{\otimes}_{\acti} \times \calC^{\otimes}_{\acti} \rightarrow \calC^{\otimes}_{\acti}$
denote the functor induced by the symmetric monoidal structure described in Corollary
\ref{letm}. This operation can be described informally as follows: if
$X \in \calC^{\otimes}_{\seg{m}}$ classifies a sequence of objects $\{ X_i \in \calC \}_{1 \leq i \leq m}$
and $Y \in \calC^{\otimes}_{\seg{n}}$ classifies a sequence $\{ Y_j \in \calC \}_{1 \leq j \leq n}$, then
$X \oplus Y \in \calC^{\otimes}_{ \seg{m+n} }$ corresponds to the sequence of objects
$\{ X_i \in \calC \}_{1 \leq i \leq m} \cup \{ Y_j \in \calC \}_{1 \leq j \leq n}$ obtained by concatenation.
\end{remark}

\begin{remark}
The symmetric monoidal structure on $\calC^{\otimes}_{\acti}$ described in
Corollary \ref{letm} can actually be extended to a symmetric monoidal structure on
$\calC^{\otimes}$ itself, but we will not need this.
\end{remark}

\begin{remark}\label{qact}
Let $\calO^{\otimes}$ be an $\infty$-operad, and let $q: \calC^{\otimes} \rightarrow \calO^{\otimes}$ be a
$\calO$-monoidal $\infty$-category. Then the collection of active $q$-coCartesian morphisms
in $\calC^{\otimes}$ is stable under the operation $\oplus: \calC^{\otimes}_{\acti} \times \calC^{\otimes}_{\acti} \rightarrow \calC^{\otimes}_{\acti}$. To see this, let
$\alpha: C \rightarrow C'$ and $\beta: D \rightarrow D'$ be active $q$-coCartesian morphisms in
$\calC^{\otimes}$. Let $\gamma: C \oplus D \rightarrow E$ be a $q$-coCartesian morphism
lifting $q( \alpha \oplus \beta)$. We have a commutative diagram in $\calO^{\otimes}$
$$ \xymatrix{ p(C) \ar[d] & p(C \oplus D) \ar[r] \ar[l] \ar[d] & p(D) \ar[d] \\
p(C') & p(C' \oplus D') \ar[l] \ar[r] & p(D') }$$ 
We can lift this to a diagram of $q$-coCartesian morphisms
$$ \xymatrix{ C \ar[d] & C \oplus D \ar[r] \ar[l] \ar[d] & D \ar[d] \\
C' & E \ar[l] \ar[r] & D'. }$$
Let $\delta: E \rightarrow C' \oplus D'$ be the canonical map in
$\calC^{\otimes}_{p(C' \oplus D')}$; the above diagram shows that
the image of $\delta$ is an equivalence in both $\calC^{\otimes}_{p(C')}$ and $\calC^{\otimes}_{p(D')}$.
Since $\calC^{\otimes}_{p(C' \oplus D')} \simeq \calC^{\otimes}_{p(C')} \times \calC^{\otimes}_{p(D')}$,
it follows that $\delta$ is an equivalence, so that $\alpha \oplus \beta$ is $q$-coCartesian as desired.
\end{remark}

\begin{remark}\label{lazier}
Since every diagonal embedding $\calO^{\otimes} \rightarrow \Act( \calO^{\otimes})$
is fully faithful, we conclude that for every fibration of $\infty$-operads
$\calC^{\otimes} \rightarrow \calO^{\otimes}$ the inclusion
$\calC^{\otimes} \hookrightarrow \Envv_{\calO}(\calC)^{\otimes}$ is fully faithful.
In particular, we deduce that for every $\infty$-operad $\calC^{\otimes}$ there
exists a fully faithful $\infty$-operad map $\calC^{\otimes} \rightarrow {\calD}^{\otimes}$,
where ${\calD}^{\otimes}$ is a symmetric monoidal $\infty$-category.
\end{remark}

We now turn to the proofs of Propositions \ref{luckma} and \ref{tull}. We will need several preliminary results.

\begin{lemma}\label{caraway}
Let $p: \calC \rightarrow \calD$ be a coCartesian fibration of $\infty$-categories.
Let $\calC'$ be a full subcategory of $\calC$ satisfying the following conditions:
\begin{itemize}
\item[$(i)$] For each $D \in \calD$, the inclusion $\calC'_{D} \subseteq \calC_{D}$ admits a
left adjoint $L_{D}$.
\item[$(ii)$] For every morphism $f: D \rightarrow D'$ in $\calD$, the associated functor
$f_{!}: \calC_{D} \rightarrow \calC_{D'}$ carries $L_{D}$-equivalences to $L_{D'}$-equivalences.
\end{itemize}
Then:
\begin{itemize}
\item[$(1)$] The restriction $p' = p | \calC'$ is a coCartesian fibration.
\item[$(2)$] Let $f: C \rightarrow C'$ be a morphism in $\calC'$
lying over $g: D \rightarrow D'$ in $\calD$, and let $g_{!}: \calC_{D} \rightarrow \calC_{D'}$
be the functor induced by the coCartesian fibration $p$. Then $f$ is $p'$-coCartesian if and only if
the induced map $\alpha: g_{!} C \rightarrow C'$ is an $L_{D'}$-equivalence.
\end{itemize}
\end{lemma}

\begin{proof}
We first prove the ``if'' direction of $(2)$. According to Proposition \toposref{charCart}, it will
suffice to show that for every object $C'' \in \calC'$ lying over $D'' \in \calD$, the outer square
in the homotopy coherent diagram
$$ \xymatrix{ \bHom_{\calC}( C', C'') \ar[r]^{\theta} \ar[d] & \bHom_{\calC}( g_{!} C, C'') \ar[r] \ar[d] & \bHom_{\calC}(C, C'') \ar[d] \\
\bHom_{\calD}( D', D'') \ar[r] & \bHom_{\calD}( D', D'') \ar[r] & \bHom_{\calD}( D, D'') }$$
is a homotopy pullback square. Since the right square is a homotopy pullback (Proposition \toposref{charCart}), it will suffice to show that $\theta$ induces a homotopy equivalence after passing the homotopy fiber over any map $h: D' \rightarrow D''$. Using Proposition
\toposref{compspaces}, we see that this is equivalent to the assertion that the map
$h_{!}(\alpha)$ is an $L_{D''}$-equivalence. This follows from $(ii)$, since
$\alpha$ is an $L_{D'}$-equivalence.

To prove $(1)$, it will suffice to show that for every $C \in \calC'$ and every $g: p(C) \rightarrow D'$
in $\calD$, there exists a morphism $f: C \rightarrow C'$ lying over $g$ satisfying the criterion
of $(2)$. We can construct $f$ as a composition
$ C \stackrel{f'}{\rightarrow} C'' \stackrel{f''}{\rightarrow} C'$, where $f'$ is a 
$p$-coCartesian lift of $g$ in $\calC$, and $f'': C'' \rightarrow C'$ exhibits
$C'$ as an $\calC'_{D'}$-localization of $C''$.

We conclude by proving the ``only if'' direction of $(2)$. Let $f: C \rightarrow C'$ be a 
$p'$-coCartesian morphism in $\calC'$ lying over $g: D \rightarrow D'$. Choose a
factorization of $f$ as a composition
$$ C \stackrel{f'}{\rightarrow} C'' \stackrel{f''}{\rightarrow} C'$$
where $f'$ is $p$-coCartesian and $f''$ is a morphism in $\calC_{D'}$.
The map $f''$ admits a factorization as a composition
$$ C'' \stackrel{h}{\rightarrow} C''' \stackrel{h'}{\rightarrow} C'$$
where $h$ exhibits $C'''$ as a $\calC'_{D'}$-localization of $C''$. 
The first part of the proof shows that the composition $h \circ f'$ is a $p'$-coCartesian lift
of $g$. Since $f$ is also a $p'$-coCartesian lift of $g$, we deduce that $h'$ is an equivalence.
It follows that $f'' = h' \circ h$ exhibits $C'$ as a $\calC'_{D'}$-localization of $C'' \simeq g_{!} C$, so that $f$ satisfies the criterion of $(2)$.
\end{proof}

\begin{remark}\label{poke}
Let $\calC' \subseteq \calC$ and $p: \calC \rightarrow \calD$ be as in Lemma \ref{caraway}.
Hypotheses $(i)$ and $(ii)$ are equivalent to the following:
\begin{itemize}
\item[$(i')$] The inclusion $\calC' \subseteq \calC$ admits a left adjoint $L$.
\item[$(ii')$] The functor $p$ carries each $L$-equivalence in $\calC$ to an equivalence in $\calD$.
\end{itemize}
Suppose first that $(i)$ and $(ii)$ are satisfied. To prove $(i')$ and $(ii')$, it will suffice
(by virtue of Proposition \toposref{testreflect}) to show that for each object
$C \in \calC_{D}$, if $f: C \rightarrow C'$ exhibits $C'$ as a $\calC'_{D}$-localization of
$C$, then $f$ exhibits $C'$ as a $\calC'$-localization of $C$. In other words, we must show that
for each $C'' \in \calC'$ lying over $D'' \in \calD$, composition with
$f$ induces a homotopy equivalence
$\bHom_{ \calC}( C', C'') \rightarrow \bHom_{\calC}(C, C'')$. Using Proposition \toposref{compspaces}, we can reduce to showing that for every morphism $g: D \rightarrow D''$ in $\calD$, the induced map
$\bHom_{ \calC_{D''} }( g_{!} C', C'') \rightarrow \bHom_{ \calC_{D''}}( g_{!} C, C'')$ is a homotopy equivalence. For this, it suffices to show that $g_{!}(f)$ is an $L_{D''}$-equivalence, which follows
immediately from $(ii)$.

Conversely, suppose that $(i')$ and $(ii')$ are satisfied. We first prove $(i)$. Fix
$D \in \calD$. To prove that the inclusion $\calC'_{D} \subseteq \calC_{D}$ admits a left
adjoint, it will suffice to show that for each object $C \in \calC_{D}$ there exists
a $\calC'_{D}$-localization of $C$ (Proposition \toposref{testreflect}). Fix a map $f: C \rightarrow C'$ in $\calC$ which exhibits $C'$ as a $\calC'$-localization of $C$. Assumption $(ii')$ guarantees that
$p(f)$ is an equivalence in $\calD$. Replacing $f$ by an equivalent morphism if necessary, we may suppose that $p(f) = \id_{D}$ so that $f$ is a morphism in $\calC_{D}$. We claim that
$f$ exhibits $C'$ as a $\calC'_{D}$-localization of $C$. To prove this, it suffices to show that
for each $C'' \in \calC'_{D}$, composition with $f$ induces a homotopy equivalence
$\bHom_{ \calC_{D}}( C', C'') \rightarrow \bHom_{ \calC_{D}}( C, C'')$. Using 
Proposition \toposref{compspaces}, we can reduce to showing that $f$ induces a homotopy
equivalence $\bHom_{\calC}(C', C'') \rightarrow \bHom_{\calC}(C,C'')$, which follows from the assumption that $f$ exhibits $C'$ as a $\calC'$-localization of $C$. 

We now prove $(ii)$. Let $f: C \rightarrow C'$ be an $L_{D}$-equivalence in $\calC_{D}$, and let
$g: D \rightarrow D''$ be a morphism in $\calD$. We wish to show that $g_{!}(f)$ is an
$L_{D''}$-equivalence in $\calC_{D''}$. In other words, we wish to show that for each
object $C'' \in \calC'_{D''}$, the map
$\bHom_{\calC_{D''}}( g_{!} C', C'') \rightarrow \bHom_{ \calC_{D''}}( g_{!} C, C'')$ is a homotopy
equivalence. This follows from Proposition \toposref{compspaces} and the fact that
$\bHom_{\calC}(C', C'') \rightarrow \bHom_{\calC}(C,C'')$ is a homotopy equivalence.
\end{remark}

\begin{lemma}\label{kape}
Let $p: \calC \rightarrow \calD$ be an inner fibration of $\infty$-categories.
Let $\calD' \subseteq \calD$ be a full subcategory and set $\calC' = \calD' \times_{\calD} \calC$.
Assume that:
\begin{itemize}
\item[$(i)$] The inclusion $\calD' \subseteq \calD$ admits a left adjoint.
\item[$(ii)$] Let $C \in \calC$ be an object and let $g: p(C) \rightarrow D$ be a morphism
which exhibits $D$ as a $\calD'$-localization of $p(C)$. Then $g$ can be lifted to
a $p$-coCartesian morphism $C \rightarrow C'$.
\end{itemize}
Then:
\begin{itemize}
\item[$(1)$] A morphism
$f: C \rightarrow C'$ exhibits $C'$ as a $\calC'$-localization of $C$ if and only if
$p(f)$ exhibits $p(C')$ as a $\calD'$-localization of $p(C)$ and $f$ is $p$-coCartesian.
\item[$(2)$] The inclusion $\calC' \subseteq \calC$ admits a left adjoint. 
\end{itemize}
\end{lemma}

\begin{proof}
We first prove the ``if'' direction of $(1)$. Fix an object $C'' \in \calC'$; we wish to prove that
$f$ induces a homotopy equivalence $\bHom_{\calC}( C', C'') \rightarrow \bHom_{\calC}( C, C'')$.
Using Proposition \toposref{charCart}, we deduce that the homotopy coherent diagram
$$ \xymatrix{ \bHom_{\calC}( C', C'') \ar[r] \ar[d] & \bHom_{\calC}( C, C'') \ar[d] \\
\bHom_{\calD}( p(C'), p(C'') ) \ar[r] & \bHom_{\calD}( p(C), p(C'') ) }$$
is a homotopy pullback square. It therefore suffices to show that the bottom horizontal map is a homotopy equivalence, which follows from the assumptions that $p(C'') \in \calD'$ and
$p(f)$ exhibits $p(C')$ as a $\calD'$-localization of $p(C)$.

Assertion $(2)$ now follows from $(ii)$ together with Proposition \toposref{testreflect}.
To complete the proof, we verify the ``only if'' direction of $(1)$. Let
$f: C \rightarrow C'$ be a map which exhibits $C'$ as a $\calC'$-localization of
$C$, and let $g: D \rightarrow D'$ be the image of $f$ in $\calD$. Then $f$ factors as a composition
$$ C \stackrel{f'}{\rightarrow} g_{!} C \stackrel{f''}{\rightarrow} C';$$
we wish to prove that $f''$ is an equivalence. This follows from the first part of the proof, which shows
that $f'$ exhibits $g_{!} C$ as a $\calC'$-localization of $C$.
\end{proof}

\begin{lemma}\label{swade}
Let $p: \calC^{\otimes} \rightarrow \calO^{\otimes}$ be a fibration of $\infty$-operads, and 
let $\calD = \calC^{\otimes} \times_{ \Fun( \{0\}, \calO^{\otimes} ) } \Fun( \Delta^1, \calO^{\otimes})$.
Then the inclusion $\Envv_{\calO}(\calC)^{\otimes} \subseteq \calD$ admits a left adjoint.
Moreover, a morphism $\alpha: D \rightarrow D'$ in $\calD$ exhibits $D'$ as an
$\Envv_{\calO}(\calC)^{\otimes}$-localization of $D$ if and only if $D'$ is active, the image of
$\alpha$ in $\calC^{\otimes}$ is inert, and the image of $\alpha$ in $\calO^{\otimes}$ is an equivalence.
\end{lemma}

\begin{proof}
According to Proposition \ref{singwell}, the active and inert morphisms determine
a factorization system on $\calO^{\otimes}$. It follows from Lemma \toposref{prefukt} that
the inclusion $\Act(\calO^{\otimes}) \subseteq \Fun( \Delta^1, \calO^{\otimes})$ admits a left adjoint, and that a morphism $\alpha: g \rightarrow g'$ in $\Fun( \Delta^1, \calO^{\otimes})$ corresponding to a commutative diagram
$$ \xymatrix{ X \ar[r]^{f} \ar[d]^{g} & X' \ar[d]^{g'} \\
Y \ar[r]^{f'} & Y' }$$
in $\calO^{\otimes}$ exhibits $g'$ as an $\Act(\calO^{\otimes})$-localization of
$g$ if and only if $g'$ is active, $f$ is inert, and $f'$ is an equivalence. The desired result
now follow from Lemma \ref{kape}.
\end{proof}

\begin{lemma}\label{iba}
Let $p: \calC^{\otimes} \rightarrow \calO^{\otimes}$ be a fibration of $\infty$-operads, and 
let $\calD = \calC^{\otimes} \times_{ \Fun( \{0\}, \calO^{\otimes} ) } \Fun( \Delta^1, \calO^{\otimes})$.
Then:
\begin{itemize}
\item[$(1)$] Evaluation at $\{1\}$ induces a coCartesian fibration
$q': \calD \rightarrow \calO^{\otimes}$.
\item[$(2)$] A morphism in $\calD$ is $q'$-coCartesian if and only if its
image in $\calC^{\otimes}$ is an equivalence.
\item[$(3)$] The map $q'$ restricts to a coCartesian fibration
$q: \Envv_{\calO}(\calC)^{\otimes} \rightarrow \calO^{\otimes}$.
\item[$(4)$] A morphism $f$ in $\Envv_{\calO}(\calC)^{\otimes}$
is $q$-coCartesian if and only if its image in $\calC^{\otimes}$ is inert.
\end{itemize}
\end{lemma}

\begin{proof}
Assertions $(1)$ and $(2)$ follow from Corollary \toposref{tweezegork}. 
Assertions $(3)$ and $(4)$ follow by combining Lemma \ref{swade}, Remark \ref{poke}, and Lemma
\ref{caraway}.
\end{proof}

\begin{lemma}\label{placeman}
Let $\calC$ denote the full subcategory of $\Fun( \Delta^1, \Nerve(\FinSeg))$ spanned by the
active morphisms, and let $p: \calC \rightarrow \Nerve(\FinSeg)$ be given by evaluation on the vertex 
$1$. Let $X \in \calC$ be an object with $p(X) = \seg{n}$, and choose $p$-coCartesian morphisms
$f_i: X \rightarrow X_i$ covering the maps $\Colp{i}: \seg{n} \rightarrow \seg{1}$
for $1 \leq i \leq n$.
These morphisms determine a $p$-limit diagram $\nostar{n}^{\triangleleft} \rightarrow \calC$.
\end{lemma}

\begin{proof}
Let $X$ be given by an active morphism $\beta: \seg{m} \rightarrow \seg{n}$.
Each of the maps $f_i$ can be identified with a commutative diagram
$$ \xymatrix{ \seg{m} \ar[r]^{\gamma_i} \ar[d]^{\beta} & \seg{m_i} \ar[d]^{\beta_i} \\
\seg{n} \ar[r]^{ \Colp{i} } & \seg{1} }$$
where $\beta_i$ is active and $\gamma_i$ is inert. Unwinding the definitions, we must show the following:
\begin{itemize}
\item[$(\ast)$] Given an active morphism $\beta': \seg{m'} \rightarrow \seg{n'}$ in 
$\FinSeg$, a map $\delta: \seg{n'} \rightarrow \seg{n}$, and a collection of commutative
diagrams
$$ \xymatrix{ \seg{m'} \ar[r]^{\epsilon_i} \ar[d] & \seg{m_i} \ar[d]^{\beta_i} \\
\seg{n'} \ar[r]^{ \Colp{i} \circ \delta} & \seg{1}, }$$
there is a unique morphism $\epsilon: \seg{m'} \rightarrow \seg{m}$ such that
$\epsilon_i = \gamma_i \circ \epsilon$.
\end{itemize}
For each $j \in \seg{m'}$, let $j' = (\delta \circ \beta')(j) \in \seg{n}$. Then
$\epsilon$ is given by the formula
$$ \epsilon(j) = \begin{cases} \ast & \text{ if } j' = \ast \\
\gamma_{j'}^{-1}( \epsilon_{i}(j) ) & \text{ if } j' \neq \ast. \end{cases}$$
\end{proof}

\begin{proof}[Proof of Proposition \ref{luckma}]
Lemma \ref{iba} implies that $q$ is a coCartesian fibration. It will therefore suffice to show that
$\Envv_{\calO}(\calC)^{\otimes}$ is an $\infty$-operad. Let $r$ denote the composition
$$\Envv_{\calO}(\calC)^{\otimes} \stackrel{q}{\rightarrow} \calO^{\otimes} \rightarrow \Nerve(\FinSeg).$$
Let $X \in \Envv_{\calO}(\calC)^{\otimes}_{\seg{n}}$, and choose $r$-coCartesian morphisms
$X \rightarrow X_i$ covering $\Colp{i}: \seg{n} \rightarrow \seg{1}$ for
$1 \leq i \leq n$. We wish to prove that these morphisms determine an $r$-limit diagram $\alpha: \nostar{n}^{\triangleleft} \rightarrow \Envv_{\calO}(\calC)^{\otimes}$. 

Let $\calD = \calC^{\otimes} \times_{ \Fun( \{0\}, \calO^{\otimes} ) } \Fun( \Delta^1, \calO^{\otimes})$.
Then $r$ extends to a map $r': \calD \rightarrow \Nerve(\FinSeg)$. To show that $\alpha$ is an
$r$-limit diagram, it will suffice to show that $\alpha$ is an $r'$-limit diagram. Write
$r'$ as a composition
$$ \calD \stackrel{r'_0}{\rightarrow} \Fun( \Delta^1, \calO^{\otimes}) \stackrel{r'_1}{\rightarrow}
\Fun(\Delta^1, \Nerve(\FinSeg)) \stackrel{r'_2}{\rightarrow} \Nerve(\FinSeg).$$
In view of Proposition \toposref{basrel}, it will suffice to show that
$\alpha$ is an $r'_0$-limit diagram, that $r'_0 \circ \alpha$ is an $r'_1$-limit diagram, and that
$r'_1 \circ r'_0 \circ \alpha$ is an $r'_2$-limit diagram. The second of these assertions follows
from Remark \ref{talee} and Lemma \monoidref{surtybove}, and the last from Lemma \ref{placeman}. 
To prove that $\alpha$ is an $r'_0$-limit diagram, we consider the pullback diagram
$$ \xymatrix{ \calD \ar[r] \ar[d] & \Fun( \Delta^1, \calO^{\otimes}) \ar[d] \\
\calC^{\otimes} \ar[r] & \calO^{\otimes}. }$$
Using Proposition \toposref{basrel}, we are reduced to the problem of showing that
the induced map $\nostar{n}^{\triangleleft} \rightarrow \calC^{\otimes}$ is a $p$-limit diagram; this follows from Remark \ref{cotallee}.

To complete the proof, it will suffice to show for any finite collection of objects
$X_{i} \in \Envv_{\calO}(\calC)^{\otimes}_{\seg{1}}$ (parametrized by $1 \leq i \leq n$), there
exists an object $X \in \Envv_{\calO}(\calC)^{\otimes}_{\seg{n}}$ and a collection
of $r$-coCartesian morphisms $X \rightarrow X_i$ covering the maps
$\Colp{i}: \seg{n} \rightarrow \seg{1}$. Each
$X_i$ can be identified with an object $C_i \in \calC^{\otimes}$ and an active morphism
$\beta_{i}: p(C_i) \rightarrow Y_i$ in $\calO^{\otimes}$, where $Y_i \in \calO$. 
The objects $Y_i$ determine a diagram $g: \nostar{n} \rightarrow \calO^{\otimes}$.
Using the assumption that $\calO^{\otimes}$ is an $\infty$-operad, we deduce
the existence of an object $Y \in \calO^{\otimes}_{\seg{n}}$ and a collection of inert morphisms
$Y \rightarrow Y_i$ covering the maps $\Colp{i}: \seg{n} \rightarrow \seg{1}$.
We regard these morphisms as providing an object $\overline{Y} \in \calO^{\otimes}_{/g}$ lifting $Y$.
Let $C_i \in \calC^{\otimes}$ lie over $\seg{m_i}$ in $\FinSeg$. Since
$\calC^{\otimes}$ is an $\infty$-operad, there exists an object
$C \in \calC^{\otimes}_{\seg{m}}$ and a collection of inert morphisms
$C \rightarrow C_i$, where $m = m_1 + \cdots + m_n$. Composing these maps with
the $\beta_{i}$, we can lift $p(C)$ to an object $Z \in \calO^{\otimes}_{/g}$.
To construct the object $X$ and the maps $X \rightarrow X_i$, it suffices to select
a morphism $Z \rightarrow \overline{Y}$ in $\calO^{\otimes}_{/g}$. The existence of such a morphism follows from the observation that $\overline{Y}$ is a final object of $\calO^{\otimes}_{/g}$.
\end{proof}

\begin{proof}[Proof of Proposition \ref{tull}]
Let $\calE^{\otimes}$ denote the essential image of $i$. We can identify
$\calE^{\otimes}$ with the full subcategory of $\Envv_{\calO}(\calC)^{\otimes}$ spanned by
those objects $(X, \alpha: p(X) \rightarrow Y)$ for which $X \in \calC^{\otimes}$ and $\alpha$ is an equivalence in $\calO^{\otimes}$. We observe that $i$ induces an equivalence of $\infty$-operads $\calC^{\otimes} \rightarrow \calE^{\otimes}$. It will therefore suffice to prove that the restriction functor
$$ \Fun^{\otimes}_{\calO}( \Envv_{\calO}(\calC)^{\otimes}, \calD^{\otimes} ) \rightarrow \Fun^{\lax}_{\calO}( \calE^{\otimes}, \calD^{\otimes})$$
is an equivalence of $\infty$-categories. In view of Proposition \toposref{lklk}, it will suffice to show the following:
\begin{itemize}
\item[$(a)$] Every $\infty$-operad map $\theta_0: \calE^{\otimes} \rightarrow \calD^{\otimes}$
admits a $q$-left Kan extension $\theta: \Envv_{\calO}(\calC)^{\otimes} \rightarrow \calD^{\otimes}$.
\item[$(b)$] An arbitrary map $\theta: \Envv_{\calO}(\calC)^{\otimes} \rightarrow \calD^{\otimes}$
in $(\sSet)_{/ \calO^{\otimes}}$ is a $\calO$-monoidal functor if and only if
it is a $q$-left Kan extension of $\theta_0 = \theta | \calE^{\otimes}$ and 
$\theta_0$ is an $\infty$-operad map.
\end{itemize}
To prove $(a)$, we will use criterion of Lemma \toposref{kan2}: it suffices to show that for
every object $\overline{X} = (X, \alpha: p(X) \rightarrow Y)$ in
$\Envv_{\calO}(\calC)^{\otimes}$, the induced diagram
$$\calE^{\otimes} \times_{ \Envv_{\calO}(\calC)^{\otimes}} \Envv_{\calO}(\calC)^{\otimes}_{/ \overline{X}} \rightarrow 
\calD^{\otimes}$$
admits a $q$-colimit covering the natural map
$$(\calE^{\otimes} \times_{ \Envv_{\calO}(\calC)^{\otimes}} \End_{\calO}(\calC)^{\otimes}_{/ \overline{X}} )^{\triangleright}
\rightarrow \calD^{\otimes}.$$
To see this, we observe that the $\infty$-category
$\calE^{\otimes} \times_{ \Envv_{\calO}(\calC)^{\otimes}} \Envv_{\calO}(\calC)^{\otimes}_{/ \overline{X}}$ has a final object, given by the pair $(X, \id_{p(X)})$. It therefore suffices to show that
there exists a $q$-coCartesian morphism $\theta_0(X, \id_{p(X)}) \rightarrow
C$ lifting $\alpha: p(X) \rightarrow Y$, which follows from the assumption that
$q$ is a coCartesian fibration. This completes the proof of $(a)$ and yields the following version of $(b)$:
\begin{itemize}
\item[$(b')$] Let $\theta: \Envv_{\calO}(\calC)^{\otimes} \rightarrow \calD^{\otimes}$ be a morphism
in $(\sSet)_{/ \calO^{\otimes}}$  such that
the restriction $\theta_0 = \theta| \calE^{\otimes}$ is an $\infty$-operad map.
Then $\theta$ is a $q$-left Kan extension of $\theta_0$ if and only if, for every object
$(X, \alpha: p(X) \rightarrow Y) \in \Envv(\calO)^{\otimes}$, the canonical map
$\theta(X, \id_{p(X)}) \rightarrow \theta(X,\alpha)$ is $q$-coCartesian.
\end{itemize}

We now prove $(b)$. Let $\theta: \Envv_{\calO}(\calC)^{\otimes} \rightarrow \calD^{\otimes}$
be such that the restriction $\theta_0 = \theta | \calE^{\otimes}$ is an
$\infty$-operad map. In view of $(b')$, it will suffice to show that
$\theta$ is a $\calO$-monoidal functor if and only if
$\theta( X, \id_{p(X)} \rightarrow \theta(X, \alpha)$ is $q$-coCartesian, for each
$(X, \alpha) \in \Envv_{\calO}(\calC)^{\otimes}$. The ``only if'' direction is clear, since Lemma \ref{placeman} implies that the morphism $(X, \id_{p(X)} \rightarrow (X, \alpha)$ in $\Envv_{\calO}(\calC)^{\otimes}$ is $p'$-coCartesian, where $p': \Envv_{\calO}(\calC)^{\otimes} \rightarrow \calO^{\otimes}$ denotes the projection. 
For the converse, suppose that we are given a $p'$-coCartesian morphism
$f: (X, \alpha) \rightarrow (Y, \alpha')$ in $\Envv_{\calO}(\calC)^{\otimes}$. Let
$\beta: p(X) \rightarrow p(Y)$ be the induced map in $\calO^{\otimes}$, and choose 
a factorization $\beta \simeq \beta'' \circ \beta'$ where $\beta'$ is inert and
$\beta''$ is active. Choose a $p$-coCartesian morphism $\overline{\beta}': X \rightarrow X''$
lifting $\beta'$. We then have a commutative diagram
$$ \xymatrix{ (X, \id_{p(X)}) \ar[r]^{f'} \ar[d]^{g} & (X'', \id_{p(X'')}) \ar[d]^{g'} \\
(X, \alpha) \ar[r]^{f} & (X', \alpha'). }$$
The description of $p'$-coCartesian morphisms supplied by Lemma
\ref{placeman} shows that the map $X'' \rightarrow X'$ is an equivalence in $\calO^{\otimes}$.
If $\theta$ satisfies the hypotheses of $(b')$, then $\theta(g)$ and $\theta(g')$ are
$q$-coCartesian. The assumption that $\theta_0$ is an $\infty$-operad map
guarantees that $\theta(f')$ is $q$-coCartesian. It follows from Proposition
\toposref{protohermes} that $\theta(f)$ is $q$-coCartesian. By allowing $f$ to
range over all morphisms in $\Envv_{\calO}(\calC)^{\otimes}$ we deduce that  
$\theta$ is a $\calO$-monoidal functor, as desired.
\end{proof}

\subsection{Subcategories of $\calO$-Monoidal $\infty$-Categories}\label{subcatmon}

Let $\calC^{\otimes}$ be an $\infty$-operad, and let $\calD \subseteq \calC$ be a full subcategory of
the underlying $\infty$-category of $\calC^{\otimes}$ which is stable under equivalence.
In this case, we let $\calD^{\otimes}$ denote the full subcategory of $\calC^{\otimes}$ spanned by those
objects $D \in \calC^{\otimes}$ having the form $D_1 \oplus \cdots \oplus D_{n}$, where
each object $D_i$ belongs to $\calD$. It follows immediately from the definitions that
$\calD^{\otimes}$ is again an $\infty$-operad, and that the inclusion $\calD^{\otimes} \subseteq \calC^{\otimes}$ is a map of $\infty$-operads. In this section, we will consider some refinements
of the preceding statement: namely, we will suppose that there exists a coCartesian fibration
of $\infty$-operads $p: \calC^{\otimes} \rightarrow \calO^{\otimes}$, and obtain criteria
which guarantee that $p | \calD^{\otimes}$ is again a coCartesian fibration of $\infty$-operads.
The most obvious case to consider is that in which $\calD$ is stable under the relevant tensor product operations on $\calC$.


\begin{proposition}\label{yak2}
Let $p: \calC^{\otimes} \rightarrow \calO^{\otimes}$ be a coCartesian fibration of $\infty$-operads, let $\calD \subseteq \calC$ be a full subcategory which is stable under equivalence, and let $\calD^{\otimes} \subseteq \calC^{\otimes}$ be defined as above. Suppose that, for every operation $f \in \Mul_{\calO}( \{ X_i \}, Y)$, the functor
$\otimes_{f}: \prod_{1 \leq i \leq n} \calC_{X_i} \rightarrow \calC_Y$ defined in Remark \ref{swinn}
carries $\prod_{1 \leq i \leq n} \calD_{X_i}$ into $\calD_Y$. Then:
\begin{itemize}

\item[$(1)$] The restricted map $\calD^{\otimes} \rightarrow \calO^{\otimes}$ is
a coCartesian fibration of $\infty$-operads.

\item[$(2)$] The inclusion $\calD^{\otimes} \subseteq \calC^{\otimes}$ is a $\calO$-monoidal functor.

\item[$(3)$] Suppose that, for every object $X \in \calO$, the inclusion
$\calD_{X} \subseteq \calC_{X}$ admits a right adjoint $L_X$ $($so that $\calD_X$ is a {\it colocalization} of $\calC_X${}$)$. Then there exists a commutative diagram
$$ \xymatrix{ \calC^{\otimes} \ar[dr]^{p} \ar[rr]^{L^{\otimes}} & & \calD^{\otimes} \ar[dl] \\
& \calO^{\otimes} & }$$
and a natural transformation of functors $\alpha: L^{\otimes} \rightarrow \id_{\calC^{\otimes}}$
which exhibits $L^{\otimes}$ as a colocalization functor $($see Proposition \toposref{recloc}$)$
and such that, for every object $C \in \calC^{\otimes}$, the image $p( \alpha(C))$ is a degenerate edge
of $\calO^{\otimes}$.
\item[$(4)$] Under the hypothesis of $(3)$, the functor $L^{\otimes}$ is a map of
$\infty$-operads.
\end{itemize}
\end{proposition}

\begin{remark}\label{yak3}
Assume that $\calO^{\otimes}$ is the commutative $\infty$-operad, and let
$\calC^{\otimes} \rightarrow \calO^{\otimes}$ be a symmetric monoidal $\infty$-category. 
A full subcategory $\calD \subseteq \calC$ (stable under equivalence) satisfies the hypotheses of Proposition \ref{yak2} if and only if $\calD$ contains the unit object of $\calC$ and is closed under the
tensor product functor $\otimes: \calC \times \calC \rightarrow \calC$.
\end{remark}

\begin{proof}
Assertions $(1)$ and $(2)$ follow immediately from the definitions. Now suppose that the hypotheses of $(3)$ are satisfied. Let us say that a morphism $f: D \rightarrow C$ in $\calC^{\otimes}$ is
{\it colocalizing} if $D \in \calD^{\otimes}$, and the projection map
$$ \psi: \calD^{\otimes} \times_{ \calC^{\otimes} } \calC^{\otimes}_{/f} 
\rightarrow \calD^{\otimes} \times_{ \calC^{\otimes} } \calC^{\otimes}_{/C}$$
is a trivial Kan fibration. Since $\psi$ is automatically a right fibration (Proposition \toposref{sharpen}), the map $f$ is colocalizing if and only if the fibers of $\psi$ are contractible (Lemma \toposref{toothie}). 
For this, it suffices to show that for each $D' \in \calD^{\otimes}$, the induced map of fibers
$$\psi_{D'}: \{ D' \} \times_{ \calC^{\otimes} } \calC^{\otimes}_{/f}
\rightarrow \{ D' \} \times_{ \calC^{\otimes} }  \calC^{\otimes}_{/C}$$
has contractible fibers. Since $\psi_{D'}$ is a right fibration between Kan complexes, it is a Kan fibration (Lemma \toposref{toothie2}). Consequently, $\psi_{D'}$ has contractible fibers if and only if
it is a homotopy equivalence. In other words, $f$ is colocalizing if and only if composition
with $f$ induces a homotopy equivalence
$$ \bHom_{\calC^{\otimes}}( D', D) \rightarrow \bHom_{\calC^{\otimes}}( D', C)$$
for all $D' \in \calD^{\otimes}$.

Suppose now that $f: D \rightarrow C$ is a morphism belonging to a particular fiber
$\calC^{\otimes}_{X}$ of $p$. Then, for every object $D' \in \calC^{\otimes}_{Y}$, we have canonical maps
$$\bHom_{\calC^{\otimes}}( D', D) \rightarrow \bHom_{\calO^{\otimes}}(Y, X) \leftarrow
\bHom_{\calC^{\otimes}}(D', C)$$
whose homotopy fibers over a morphisms $g: Y \rightarrow X$ in $\calO^{\otimes}$ can be identified
with $\bHom_{\calC^{\otimes}_{X}}( g_{!} D', D)$ and $\bHom_{\calC^{\otimes}_{X}}( g_{!} D', C)$
(Proposition \toposref{compspaces}). Our hypothesis that $\calD$ is stable under the tensor operations on $\calC$ implies that $g_{!} D' \in \calD^{\otimes}_{X}$. It follows that $f$ is colocalizing in $\calC^{\otimes}$ if and only if it is colocalizing in $\calC^{\otimes}_{X}$.

To prove $(3)$, we need to construct a commutative diagram
$$ \xymatrix{ \calC^{\otimes} \times \Delta^1 \ar[rr]^{\alpha} \ar[dr] & & \calC^{\otimes} \ar[dl]^{p} \\
& \calO^{\otimes} & }$$
with $\alpha | \calC^{\otimes} \times \{1\}$ equal to the identity, and having 
the property that, for every object $C \in \calC^{\otimes}$, the restriction
$\alpha| \{ C\} \times \Delta^1$ is colocalizing. For this, we construct $\alpha$ one simplex at a time. To define $\alpha$ on a vertex $C \in \calC^{\otimes}_{X}$ where
$X \in \calO^{\otimes}_{\seg{n}}$ corresponds to a sequence of objects
$\{ X_i \}_{1 \leq i \leq n}$ in $\calO$, we use the equivalence
$\calC^{\otimes}_{X} \simeq \prod_{1 \leq i \leq X} \calC_{X_i}$ and the assumption that 
each $\calD_{X_i} \subseteq \calC_{X_i}$
is a colocalization. For simplices of higher dimension, we need to solve extensions problems of the type indicated below:
$$ \xymatrix{ (\bd \Delta^n \times \Delta^1) \coprod_{ \bd \Delta^n \times \{1\} } ( \Delta^n \times \{1\}) \ar[r]^-{\alpha_0} \ar@{^{(}->}[d] & \calC^{\otimes} \ar[d]^{p} \\
\Delta^n \times \Delta^1 \ar[r] \ar@{-->}[ur] & \calO^{\otimes}, }$$
where $n > 0$ and $\alpha_0$ carries each of the edges $\{ i \} \times \Delta^1$ to a colocalizing morphism in $\calC^{\otimes}$. We now observe that $\Delta^n \times \Delta^1$ admits a filtration
$$ X_0 \subseteq X_1 \subseteq \ldots \subseteq X_n \subseteq X_{n+1} = \Delta^n \times \Delta^1,$$
where $X_0 = (\bd \Delta^n \times \Delta^1) \coprod_{ \bd \Delta^n \times \{1\} } ( \Delta^n \times \{1\})$ and there exist pushout diagrams
$$ \xymatrix{ \Lambda^{n+1}_{i+1} \ar@{^{(}->}[r] \ar[d] & \Delta^{n+1} \ar[d] \\
X_{i} \ar[r] & X_{i+1}. }$$
We now argue, by induction on $i$, that the map $\alpha_0$ admits an extension to $X_i$ (compatible with the projection $p$). For $i \leq n$, this follows from the fact that $p$ is an inner fibration. For $i = n+1$, it follows from the definition of a colocalizing morphism. This completes the proof of $(3)$. Assertion $(4)$ follows from the construction.
\end{proof}


\begin{remark}\label{yukan}
Let $\calD^{\otimes} \subseteq \calC^{\otimes}$ be as in the statement of Proposition \ref{yak2}, and suppose that each inclusion $\calD_{X} \subseteq \calC_{X}$ admits a right adjoint $L_X$.
The inclusion $i: \calD^{\otimes} \subseteq \calC^{\otimes}$ is a $\calO$-monoidal functor, and therefore induces a fully faithful embedding $\Alg_{\calO}(\calD) \subseteq \Alg_{\calO}( \calC)$. Let $L^{\otimes}$ be the functor constructed in Proposition \ref{yak2}. Since $L^{\otimes}$ is a map of
$\infty$-operads, it also induces a
functor $f: \Alg_{\calO}(\calC) \rightarrow \Alg_{\calO}(\calD)$. It is not difficult to see that $f$ is
right adjoint to the inclusion $\Alg_{\calO}(\calD) \subseteq \Alg_{\calO}(\calC)$. Moreover, if $\theta: \Alg_{\calO}(\calC) \rightarrow \calC_{X}$ denotes the evaluation functor associated to an object
$X \in \calO$, then for each $A \in \Alg_{\calO}(\calC)$ the map $\theta(f(A)) \rightarrow \theta(A)$
induced by the colocalization map $f(A) \rightarrow A$ determines an equivalence
$\theta(f(A)) \simeq L_{X} \theta(A)$.
\end{remark}

Proposition \ref{yak2} has an obvious converse: if $i: \calD^{\otimes} \rightarrow \calC^{\otimes}$
is a fully faithful $\calO$-monoidal functor between $\calO$-monoidal $\infty$-categories, then
then the essential image of $\calD$ in $\calC$ is stable under the tensor operations
of Remark \ref{swinn}. However, it is possible for this converse to fail if we only assume that
$i$ belongs to $\Alg_{\calD}(\calC)$. We now discuss a general class of examples where $\calD$ is {\em not} stable under tensor products, yet $\calD$ nonetheless inherits a $\calO$-monoidal structure from that of $\calC$.

\begin{definition}\label{compatsymon}
Let $\calC^{\otimes} \rightarrow \calO^{\otimes}$
be a coCartesian fibration of $\infty$-operads and suppose we are given a family of
localization functors $L_{X}: \calC_{X} \rightarrow \calC_{X}$ for $X \in \calO$. 
We will say that the family $\{ L_{X} \}_{X \in \calO}$ is {\it compatible with the
$\calO$-monoidal structure on $\calC$} if the following condition
is satisfied:
\begin{itemize}
\item[$(\ast)$] Let $f \in \Mul_{\calO}( \{ X_i \}_{1 \leq i \leq n}, Y)$ be an operation in
$\calO$, and suppose we are given morphisms $g_i$ in $\calC_{X_i}$ for $1 \leq i \leq n$.
If each $g_i$ is a $L_{X_i}$-equivalence, then the morphism 
$\otimes_{f} \{ g_i \}_{1 \leq i \leq n}$ in $\calC_Y$ is an $L_Y$-equivalence.
\end{itemize}
\end{definition}

\begin{example}
In the situation where $\calO^{\otimes}$ is the commutative $\infty$-operad, the condition of
Definition \ref{compatsymon} can be simplified: a localization functor
$L: \calC \rightarrow \calC$ is compatible with a symmetric monoidal structure on
$\calC$ if and only if, for every $L$-equivalence $X \rightarrow Y$ in $\calC$ and
every object $Z \in \calC$, the induced map $X \otimes Z \rightarrow Y \otimes Z$ is again an
$L$-equivalence. This is equivalent to the requirement that $L$ be compatible with
the induced monoidal structure on $\calC$ (Definition \monoidref{compatmon}).
\end{example}


\begin{proposition}\label{localjerk2}
Let $p: \calC^{\otimes} \rightarrow \calO^{\otimes}$ be a
be a coCartesian fibration of $\infty$-operads and suppose we are given a family of localization functors
$\{ L_X: \calC_X \rightarrow \calC_{X} \}_{X \in \calO}$ which are compatible with
the $\calO$-monoidal structure on $\calC$. Let $\calD$ denote the collection
of all objects of $\calC$ which lie in the image of some $L_X$, and let $\calD^{\otimes}$ be
defined as above. Then:

\begin{itemize}
\item[$(1)$] There exists a commutative diagram
$$ \xymatrix{ \calC^{\otimes} \ar[rr]^{L^{\otimes}} \ar[dr]^{p} & & \calD^{\otimes} \ar[dl] \\
& \calO^{\otimes} & }$$
and a natural transformation $\alpha: \id_{\calC^{\otimes}} \rightarrow L^{\otimes}$ which
exhibits $L^{\otimes}$ as a left adjoint to the inclusion $\calD^{\otimes} \subseteq \calC^{\otimes}$ 
and such that $p( \alpha)$ is the identity natural transformation from $p$ to itself.

\item[$(2)$] The restriction $p| \calD^{\otimes}: \calD^{\otimes} \rightarrow \Nerve(\FinSeg)$ is a coCartesian fibration of $\infty$-operads.

\item[$(3)$] The inclusion functor $\calD^{\otimes} \subseteq \calC^{\otimes}$ is a map of $\infty$-operads and $L^{\otimes}: \calC^{\otimes} \rightarrow \calD^{\otimes}$ is a $\calO$-monoidal functor.
\end{itemize}
\end{proposition}

\begin{proof}
We first prove $(1)$. Let us say that a map $f: C \rightarrow D$ in $\calC^{\otimes}$ is
{\it localizing} if it induces a trivial Kan fibration
$$ \calD^{\otimes} \times_{ \calC^{\otimes} } \calC^{\otimes}_{f/} 
\rightarrow \calD^{\otimes} \times_{ \calC^{\otimes} } \calC^{\otimes}_{D/}.$$
Arguing as in the proof of Proposition \ref{yak2}, we see that $f$ is localizing if and only if,
for every object $D' \in \calD^{\otimes}$, composition with $f$ induces a homotopy equivalence
$\bHom_{ \calC^{\otimes}} (D, D') \rightarrow \bHom_{ \calC^{\otimes}}(C,D').$
Suppose that $f$ lies in a fiber $\calC^{\otimes}_{X}$ of $p$, and let
$D' \in \calC^{\otimes}_{Y}$ of $p$. We have canonical maps
$$ \bHom_{\calC^{\otimes}}(D, D') \rightarrow \bHom_{\calO^{\otimes}}( X, Y) \leftarrow \bHom_{\calC^{\otimes}}(C, D')$$
and Proposition \toposref{compspaces} allows us to identify the homotopy fibers of these maps over a morphism $g: X \rightarrow Y$ in $\calO^{\otimes}$ with
$\bHom_{\calC^{\otimes}_{Y}}( g_{!} D, D')$ and $\bHom_{\calC^{\otimes}_{Y}}( g_{!} C, D')$.
It follows that $f$ is localizing in $\calC^{\otimes}$ if and only if 
$g_{!}(f)$ is localizing in $\calC^{\otimes}_{Y}$ for every map $g: X \rightarrow Y$ in
$\calO^{\otimes}$. Let $X \in \calO^{\otimes}_{\seg{n}}$ and $Y \in \calO^{\otimes}_{\seg{m}}$ correspond to sequences of objects $\{ X_i \}_{1 \leq i \leq n}$ and $\{ Y_j \}_{1 \leq j \leq m}$ in
$\calO$. Using the equivalence $\calC^{\otimes}_{Y} \simeq \prod_{1 \leq j \leq m} \calC_{Y_i}$ we see that it suffices to check this in the case $m=1$. Invoking the assumption that the localization functors
Invoking the assumption that the localization functors $\{ L_Z \}_{Z \in \calO}$ are compatible
with the $\calO$-monoidal structure on $\calC^{\otimes}$, we see that $f$ is localizing if
and only if it induces a $L_{X_i}$-localizing morphism in $\calC_{X_i}$ for $1 \leq i \leq n$. 

We now argue as in the proof of Proposition \ref{yak2}. To prove $(1)$, we need to construct
a commutative diagram
$$ \xymatrix{ \calC^{\otimes} \times \Delta^1 \ar[rr]^{\alpha} \ar[dr] & & \calC^{\otimes} \ar[dl]^{p} \\
& \Nerve(\FinSeg) & }$$
with $\alpha | \calC^{\otimes} \times \{0\}$ equal to the identity, and having 
the property that, for every object $C \in \calC^{\otimes}$, the restriction
$\alpha| \{ C\} \times \Delta^1$ is localizing. For this, we construct $\alpha$ one simplex at a time. To define $\alpha$ on a vertex $C \in \calC^{\otimes}_{X}$ for $X \in \calO^{\otimes}_{\seg{n}}$, we choose a morphism $f: C \rightarrow C'$ corresponding to the localization morphisms 
$C_i \rightarrow L_{X_i} C_i$ under the equivalence $\calC^{\otimes}_{X} \simeq \prod_{1 \leq i \leq n} \calC_{X_i}$ and apply the argument above.
For simplices of higher dimension, we need to solve extensions problems of the type indicated below:
$$ \xymatrix{ (\bd \Delta^n \times \Delta^1) \coprod_{ \bd \Delta^n \times \{0\} } ( \Delta^n \times \{0\}) \ar[r]^-{\alpha_0} \ar@{^{(}->}[d] & \calC^{\otimes} \ar[d]^{p} \\
\Delta^n \times \Delta^1 \ar[r] \ar@{-->}[ur] & \calO^{\otimes}, }$$
where $n > 0$ and $\alpha_0$ carries each of the edges $\{ i \} \times \Delta^1$ to a localizing morphism in $\calC^{\otimes}$. We now observe that $\Delta^n \times \Delta^1$ admits a filtration
$$ X_0 \subseteq X_1 \subseteq \ldots \subseteq X_n \subseteq X_{n+1} = \Delta^n \times \Delta^1,$$
where $X_0 = (\bd \Delta^n \times \Delta^1) \coprod_{ \bd \Delta^n \times \{0\} } ( \Delta^n \times \{0\})$ and there exist pushout diagrams
$$ \xymatrix{ \Lambda^{n+1}_{n-i} \ar@{^{(}->}[r] \ar[d] & \Delta^{n+1} \ar[d] \\
X_{i} \ar[r] & X_{i+1}. }$$
We now argue, by induction on $i$, that the map $\alpha_0$ admits an extension to $X_i$ (compatible with the projection $p$). For $i \leq n$, this follows from the fact that $p$ is an inner fibration. For $i = n+1$, it follows from the definition of a localizing morphism. This completes the proof of $(1)$.

Lemma \monoidref{surine} implies that $p' = p | \calD^{\otimes}$ is a coCartesian fibration. It follows immediately from the definition that for every object
$X \in \calO^{\otimes}_{\seg{m}}$ corresponding to $\{ X_i \in \calO \}_{1 \leq i \leq m}$, the equivalence
$\calC^{\otimes}_{X} \simeq \prod_{1 \leq i \leq m} \calC_{X_i}$ restricts to an equivalence
$\calD^{\otimes}_{X} \simeq \prod_{1 \leq i \leq m} \calD_{X_i}$. This proves that $p': \calD^{\otimes} \rightarrow \calO^{\otimes}$ is a coCartesian fibration of $\infty$-operads and that the inclusion
$\calD^{\otimes} \subseteq \calC^{\otimes}$ is a map of $\infty$-operads. Lemma \monoidref{surine} implies that $L^{\otimes}$ carries $p$-coCartesian edges to $p'$-coCartesian edges and is therefore a $\calO$-monoidal functor.
\end{proof}

\begin{example}\label{swoon}
Let $\calC$ be a symmetric monoidal $\infty$-category which is stable. We will say that
a t-structure on $\calC$ is {\it compatible} with the symmetric monoidal structure on $\calC$ if it is compatible with the underlying monoidal structure on $\calC$, in the sense of Definition \monoidref{tuppa}. In other words, a t-structure on $\calC$ is compatible with the symmetric monoidal structure if
\begin{itemize}
\item[$(1)$] For every object $C \in \calC$, the functor $C \otimes \bigdot$ is exact.
\item[$(2)$] The full subcategory $\calC_{\geq 0}$ contains the unit object $1_{\calC}$ and is closed under tensor products.
\end{itemize}
Proposition \ref{yak2} implies that $\calC_{\geq 0}$ inherits a the structure of a symmetric monoidal $\infty$-category. Applying Proposition \monoidref{jumperr}, we deduce that the truncation functors
$\tau_{ \leq n}: \calC_{\geq 0} \rightarrow ( \calC_{\geq 0} )_{\leq n}$
are compatible with the symmetric monoidal structure on $\calC_{\geq 0}$. Consequently, 
Proposition \ref{localjerk2} implies that $( \calC_{\geq 0} )_{\leq n}$ inherits the structure of 
a symmetric monoidal $\infty$-category. Taking $n=0$, we deduce that the heart of $\calC$ inherits the structure of a symmetric monoidal category (in the sense of classical category theory).
\end{example}

\subsection{$\infty$-Preoperads}\label{comm1.8}

In \S \ref{algobj}, we saw that the collection of all
$\infty$-operads could be organized into an $\infty$-category
$\Cat_{\infty}^{\lax}$. In this section, we will show that $\Cat_{\infty}^{\lax}$ can
be identified with the underlying $\infty$-category of a simplicial model category
$\PreOp$.

\begin{definition}
An {\it $\infty$-preoperad} is a marked simplicial set $(X,M)$ equipped with
a map of simplicial sets $f: X \rightarrow \Nerve(\FinSeg)$ with the following property:
for each edge $e$ of $X$ which belongs to $M$, the image $f(e)$ is an inert morphism
in $\FinSeg$. A morphism from an $\infty$-preoperad $(X,M)$ to an
$\infty$-preoperad $(Y,N)$ is a map of marked simplicial sets $(X, M) \rightarrow (Y,N)$ such that the diagram
$$ \xymatrix{ X \ar[rr] \ar[dr] & & Y \ar[dl] \\
& \Nerve(\FinSeg). & }$$
The collection of $\infty$-preoperads (and $\infty$-preoperad morphisms) forms a category,
which we will denote by $\PreOp$.
\end{definition}

\begin{remark}\label{rej}
The category $\PreOp$ of $\infty$-preoperads is naturally tensored over simplicial sets:
if $\overline{X}$ is an $\infty$-preoperad and $K$ is a simplicial set, we let
$\overline{X} \otimes K = \overline{X} \times K^{\sharp}$. This construction endows
$\PreOp$ with the structure of a simplicial category.
\end{remark}

\begin{notation}
Let $\calO^{\otimes}$ be an $\infty$-operad. We let $\calO^{\otimes, \natural}$ denote the
$\infty$-preoperad $( \calO^{\otimes}, M)$, where $M$ is the collection of all inert morphisms
in $\calO^{\otimes}$.
\end{notation}

\begin{proposition}\label{cannwell}
There exists a left proper combinatorial model structure on $\PreOp$ which may be characterized as follows:
\begin{itemize}
\item[$(C)$] A morphism $f: \overline{X} \rightarrow \overline{Y}$ in $\PreOp$ is a cofibration if and only if it induces a monomorphism between the underlying simplicial sets of $\overline{X}$ and $\overline{Y}$.

\item[$(W)$] A morphism $f: \overline{X} \rightarrow \overline{Y}$ in $\PreOp$ is a weak equivalence if and only if, for every $\infty$-operad $\calO^{\otimes}$, the induced map
$$ \bHom_{ \PreOp}( \overline{Y}, \calO^{\otimes, \natural}) \rightarrow \bHom_{ \PreOp}( \overline{X}, \calO^{\otimes, \natural} )$$
is a homotopy equivalence of simplicial sets.
\end{itemize}

This model structure is compatible with the simplicial structure of Remark \ref{rej}. Moreover,
an object of $\PreOp$ is fibrant if and only if it is of the form $\calO^{\otimes, \natural}$, for some
$\infty$-operad $\calO^{\otimes}$.
\end{proposition}

\begin{proof}
Let $\CatP = (M, T, \{ p_{\alpha}: \Lambda^2_0 \rightarrow \Nerve(\FinSeg) \}_{\alpha \in A} )$
be the categorical pattern on $\Nerve(\FinSeg)$ where $M$ consists of all
inert morphisms, $T$ the collection of all $2$-simplices, and $A$ ranges over all diagrams
$$ \seg{p} \leftarrow \seg{n} \rightarrow \seg{q}$$
where the maps are inert and induce a bijection $\nostar{n} \simeq \nostar{p} \coprod \nostar{q}$.
We now apply Theorem \bicatref{theo}.
\end{proof}

\begin{remark}
We will refer to the model structure of Proposition \ref{cannwell} as
the {\it $\infty$-operadic model structure} on $\PreOp$.
\end{remark}

The following result is a special case of Proposition \bicatref{fibraman}:

\begin{proposition}\label{sayid}
Let $\calO^{\otimes}$ be an $\infty$-operad. A map
$\overline{X} \rightarrow \calO^{\otimes, \natural}$ in $\PreOp$ is a fibration (with respect to the model structure of Proposition \ref{cannwell}) if and only if
$\overline{X} \simeq \calC^{\otimes, \natural}$, where the underlying map
$\calC^{\otimes} \rightarrow \calO^{\otimes}$ is a fibration of $\infty$-operads (in the sense of
Definition \ref{quas}).
\end{proposition}

\begin{example}\label{kapf}
The inclusion $\{ \seg{1} \}^{\flat} \subseteq \Triv^{\natural}$ is a weak
equivalence of $\infty$-preoperads. This is an immediate
consequence of Example \ref{algtriv}. In other words, we can view
$\Triv^{\natural}$ as a fibrant replacement for the object $\{ \seg{1} \}^{\flat} \in \PreOp$.
\end{example}

\begin{example}\label{capf}
Using Proposition \ref{ezalg}, we deduce that the morphism $\seg{0} \rightarrow \seg{1}$
in $\FinSeg$ induces a weak equivalence of $\infty$-preoperads
$( \Delta^1)^{\flat} \subseteq \OpE{0}^{\natural}$, so that
$\OpE{0}^{\natural}$ can be viewed as a fibrant replacement for
the $\infty$-preoperad $(\Delta^1)^{\flat}$.
\end{example}

\begin{example}
Let $\overline{ \Nerve(\cDelta) } = ( \Nerve(\cDelta), M)$, where $M$ is the collection
of edges of $\Nerve(\cDelta)$ corresponding to convex morphisms in $\cDelta$. Proposition
\ref{algass} implies that the functor $\phi$ of Construction \ref{urpas} induces a
weak equivalence of $\infty$-preoperads $\overline{ \Nerve(\cDelta)}^{op} \rightarrow
\Ass^{\natural}$. 
\end{example}

Our next goal is to show that $\PreOp$ has the structure of a {\em monoidal} model category, underlying the operation of {\em tensor product} of $\infty$-operads.

\begin{notation}\label{crumb}
We define a functor $\wedge: \FinSeg \times \FinSeg \rightarrow \FinSeg$ as follows:
\begin{itemize}
\item[$(i)$] On objects, $\wedge$ is given by the formula
$\seg{m} \wedge \seg{n} = \seg{mn}$.
\item[$(ii)$] If $f: \seg{m} \rightarrow \seg{m'}$ and $g: \seg{n} \rightarrow \seg{n'}$ are morphisms in $\FinSeg$, then $f \wedge g$ is given by the formula
$$ (f \wedge g)( an+b - n) = \begin{cases}
\ast & \text{if } f(a) = \ast \text{ or } g(b) = \ast \\
f(a) n' + g(b) - n' & \text{otherwise.} \end{cases}$$
In other words, $\wedge$ is given by the formula
$\seg{m} \wedge \seg{n} = ( \nostar{m} \times \nostar{n})_{\ast}$, where we identify
$\nostar{m} \times \nostar{n}$ with $\nostar{mn}$ via the lexicographical ordering.
\end{itemize}
\end{notation}

\begin{definition}
Let $\overline{X} = (X,M)$ and $\overline{Y} = (Y, N)$ be $\infty$-preoperads. We let
$\overline{X} \odot \overline{Y}$ denote the marked simplicial set
$\overline{X} \times \overline{Y}$, regarded as an $\infty$-preoperad by the composite map
$$ X \times Y \rightarrow \Nerve(\FinSeg) \times \Nerve(\FinSeg) 
\stackrel{\wedge}{\rightarrow} \Nerve(\FinSeg).$$
The construction $(\overline{X}, \overline{Y}) \mapsto \overline{X} \odot \overline{Y}$ endows
the category $\PreOp$ of $\infty$-preoperads with a monoidal structure.
\end{definition}

\begin{warning}
The monoidal structure on the category $\PreOp$ of $\infty$-preoperads can be regarded
as obtained from the monoidal structure on the category $\FinSeg$ via the functor
$\wedge$. In order to guarantee that the operation $\odot$ on $\PreOp$ is associative up to isomorphism (rather than merely associative up to weak equivalence), it is necessary to arrange that the operation $\wedge$ on $\FinSeg$ is {\em strictly} associative (the construction of
Notation \ref{crumb} has this property). We note that $\FinSeg$ is actually a symmetric monoidal category with respect to $\wedge$. However, this does not lead to a symmetric monoidal 
structure on $\PreOp$: if $\overline{X}$ and $\overline{Y}$ are $\infty$-preoperads, then
$\overline{X} \odot \overline{Y}$ and $\overline{Y} \odot \overline{X}$ are generally not isomorphic
in $\PreOp$ (though there is a canonical isomorphism between them in the homotopy category
$\h{ \PreOp}$).
\end{warning}

\begin{proposition}\label{candlestick}
The functor $\odot$ endows $\PreOp$ with the structure of a monoidal model category.
\end{proposition}

\begin{proof}
Since every object of $\PreOp$ is cofibrant, it will suffice to show that the functor
$$\odot: \PreOp \times \PreOp \rightarrow \PreOp$$ is a left Quillen bifunctor.
Let $\CatP$ be as in the proof of Proposition \ref{cannwell}. Using Remark \bicatref{slaver},
we deduce that the Cartesian product functor
$\PreOp \times \PreOp \rightarrow \mset{ \CatP \times \CatP}$ is a left Quillen bifunctor.
The desired result now follows by applying Proposition \bicatref{supper} to the product
functor $\Nerve(\FinSeg) \times \Nerve(\FinSeg) \rightarrow \Nerve(\FinSeg)$.
\end{proof}

\begin{remark}\label{stoss}
Since every object of $\PreOp$ is cofibrant, we observe that the functor
$\odot$ preserves weak equivalences separately in each variable.
\end{remark}

\begin{notation}\label{scambag}
Let $\calO^{\otimes}$ and ${\calO'}^{\otimes}$ be maps of $\infty$-operads, let
$\overline{X}$ be an $\infty$-preoperad, and suppose we are given a map
$f: \overline{X} \odot \calO^{\otimes, \natural} \rightarrow {\calO'}^{\otimes, \natural}$.
For every fibration of $\infty$-operads $\calC^{\otimes} \rightarrow \calO^{\otimes}$, we let
$\overline{\Alg}_{\calO}( \calC)^{\otimes}$ denote an $\infty$-preoperad with a map
$\overline{\Alg}_{\calO}( \calC)^{\otimes} \rightarrow \overline{X}$ with the following universal property: for every map of $\infty$-preoperads $\overline{Y} \rightarrow \overline{X}$, we have a canonical bijection
$$ \Hom_{( \PreOp)_{/ \overline{X}}}( \overline{Y}, \overline{\Alg}_{\calO}(\calC)^{\otimes})
\simeq \Hom_{(\PreOp)_{/ {\calO'}^{\otimes} }}( \overline{Y} \odot {\calO}^{\otimes, \natural},
\calC^{\otimes, \natural}).$$

We will denote the underlying simplicial set of $\overline{\Alg}_{\calO}(\calC)^{\otimes}$ by
$\Alg_{\calO}(\calC)^{\otimes}$.
\end{notation}

\begin{remark}
It follows immediately from Proposition \ref{candlestick} (and Proposition \ref{sayid}) that the map
$\overline{\Alg}_{\calO}( \calC)^{\otimes} \rightarrow \overline{X}$ is a fibration
in $\PreOp$. In particular, if $\overline{X} = \calD^{\otimes, \natural}$ for some $\infty$-operad $\calD^{\otimes}$, then the map $\Alg_{\calO}(\calC)^{\otimes} \rightarrow \calD^{\otimes}$ is a fibration of $\infty$-operads.
\end{remark}

\begin{example}\label{clise}
Suppose that $\overline{X} = ( \Delta^0)^{\sharp}$, regarded as an $\infty$-preoperad by the map
$\Delta^0 \simeq \{ \seg{1} \} \subseteq \Nerve(\FinSeg)$. Then
$\overline{X} \odot \calO^{\otimes, \natural} \simeq \calO^{\otimes, \natural}$, so that we can identify
a map $f$ as in Notation \ref{scambag} with a map of $\infty$-operads
$\calO^{\otimes} \rightarrow {\calO'}^{\otimes}$. In this case, we can identify
$\Alg_{\calO}(\calC)^{\otimes}$ with the $\infty$-category $\Alg_{\calO}(\calC)$ of Definition
\ref{susqat}. Moreover, an edge of $\overline{\Alg}_{\calO}(\calC)^{\otimes}$ is marked if and only if it is an equivalence in $\Alg_{\calO}(\calC)$.

More generally, if $\overline{X} \in \PreOp$ is an arbitrary object containing a vertex $x$
lying over $\seg{1} \in \FinSeg$, then a map $f: \overline{X} \odot \calO^{\otimes, \natural} \rightarrow {\calO'}^{\otimes, \natural}$ restricts on $\{x\}^{\sharp} \odot \calO^{\otimes, \natural}$ to a map
of $\infty$-operads $\calO^{\otimes} \rightarrow {\calO'}^{\otimes}$, and the fiber
$\Alg_{\calO}(\calC)^{\otimes}_{x}$ is canonically isomorphic to $\Alg_{\calO}(\calC)$.
\end{example}

\begin{definition}\label{rocks}
Let $\calO^{\otimes}$, ${\calO'}^{\otimes}$, and ${\calO''}^{\otimes}$ be $\infty$-operads.
We will say that a map of simplicial sets
$$ \calO^{\otimes} \times {\calO'}^{\otimes} \rightarrow {\calO''}^{\otimes}$$
is a {\it bifunctor of $\infty$-operads} if it induces a map of $\infty$-preoperads
$\calO^{\otimes, \natural} \odot {\calO'}^{\otimes, \natural} \rightarrow
{\calO''}^{\otimes, \natural}$.
\end{definition}

More concretely, a map $ \calO^{\otimes} \times {\calO'}^{\otimes} \rightarrow {\calO''}^{\otimes}$
as in Definition \ref{rocks} is a bifunctor of $\infty$-operads if the following conditions are satisfied:

\begin{itemize}
\item[$(1)$] The diagram
$$ \xymatrix{ \calO^{\otimes} \times {\calO'}^{\otimes} \ar[r] \ar[d] & {\calO''}^{\otimes} \ar[d] \\
\Nerve(\FinSeg) \times \Nerve(\FinSeg) \ar[r]^-{\wedge} & \Nerve(\FinSeg) }$$
commutes.
\item[$(2)$] If $f$ is a morphism in $\calO^{\otimes} \times {\calO'}^{\otimes}$ which
projects to an inert morphism in each factor, then the image of
$f$ in ${\calO''}^{\otimes}$ is inert.
\end{itemize}

\begin{proposition}\label{slape}
Suppose given a bifunctor of
$\infty$-operads $\calO^{\otimes} \times {\calO'}^{\otimes} \rightarrow {\calO''}^{\otimes}$. Let $p: \calC^{\otimes} \rightarrow {\calO''}^{\otimes}$ be a
$\calO''$-monoidal $\infty$-category. Then:
\begin{itemize}
\item[$(1)$] Define a simplicial set $\calD^{\otimes}$ equipped with a map
$q: \calD^{\otimes} \rightarrow \calO^{\otimes}$ so that the following condition is satisfied:
for every map of simplicial sets $K \rightarrow \calO^{\otimes}$, we have a canonical
bijection
$$ \Hom_{ \calO^{\otimes}}(K, \calD^{\otimes}) \simeq
\Hom_{ {\calO''}^{\otimes}}( K \times {\calO'}^{\otimes}, \calC^{\otimes}).$$
Then $q$ is a coCartesian fibration. In particular, $\calD^{\otimes}$ is an $\infty$-category.

\item[$(2)$] A morphism $f: D \rightarrow D'$ in $\calD^{\otimes}$ is $q$-coCartesian if and only if,
for each object $X \in {\calO'}^{\otimes}$, the induced map $D(X) \rightarrow D'(X)$ in
$\calC^{\otimes}$ is $p$-coCartesian.

\item[$(3)$] The full subcategory $\Alg_{\calO'}(\calC)^{\otimes} \subseteq \calD^{\otimes}$
has the following property: if $f: D \rightarrow D'$ is a coCartesian morphism
in $\calD^{\otimes}$ such that $D \in \Alg_{\calO'}(\calC)^{\otimes}$, then
$D' \in \Alg_{\calO'}(\calC)^{\otimes}$.

\item[$(4)$] The restriction $q_0$ of $q$ to $\Alg_{\calO'}(\calC)^{\otimes}$ is a coCartesian fibration, which therefore exhibits $\Alg_{\calO'}(\calC)$ as a $\calO$-monoidal $\infty$-category.

\item[$(5)$] A morphism $f: D \rightarrow D'$ in $\Alg_{\calO'}(\calC)^{\otimes}$
is $q_0$-coCartesian if and only if, for every $X \in \calO'$, the induced map
$D(X) \rightarrow D'(X)$ in $\calC^{\otimes}$ is $q$-coCartesian.
\end{itemize}
\end{proposition}

\begin{proof}
Assertions $(1)$ and $(2)$ follow immediately from Proposition \toposref{doog}, and
$(4)$ is an obvious consequence of $(1)$, $(2)$, and $(3)$. To prove $(3)$, suppose
we are given a $q$-coCartesian morphism $f: D \rightarrow D'$ in $\calD^{\otimes}$
where $D \in \Alg_{\calO'}(\calC)^{\otimes}$. We can view $f$ as a natural transformation of functors $D,D': {\calO'}^{\otimes} \rightarrow \calC^{\otimes}$.
Let $\alpha: X \rightarrow Y$ be an inert morphism of ${\calO'}^{\otimes}$, so we have a commutative diagram
$$ \xymatrix{ D(X) \ar[r] \ar[d] & D(Y) \ar[d] \\
D'(X) \ar[r] & D'(Y). }$$
We wish to prove that the bottom horizontal morphism is $p$-coCartesian. 
Invoking Proposition \toposref{protohermes}, we may reduce to showing that the upper
horizontal morphism and both of the vertical morphisms are $p$-coCartesian.
The first follows from the assumption that $D \in  \Alg_{\calO'}(\calC)^{\otimes}$, and the
second from the characterization of $q$-coCartesian morphisms given by $(2)$.

We now prove $(5)$. The ``only if'' assertion follows immediately from $(2)$.
To prove the converse, suppose that $f: D \rightarrow D'$ is a morphism in
$\Alg_{\calO'}(\calC)^{\otimes}$ such that $D(X) \rightarrow D'(X)$ is
$q$-coCartesian for every $X \in \calO'$. Choose a factorization of $D$ as a composition
$$ D \stackrel{f'}{\rightarrow} D'' \stackrel{f''}{\rightarrow} D'$$
where $f'$ is $q_0$-coCartesian and $q_0(f'')$ is degenerate. To prove that
$f$ is $q_0$ coCartesian, it will suffice to show that $f''$ is an equivalence.
Using Proposition \toposref{protohermes}, we deduce that $f''$ induces an
equivalence $D''(X) \rightarrow D'(X)$ for each $X \in \calO'$.
Fix any object $Y \in {\calO'}^{\otimes}$, and choose inert morphisms
$Y \rightarrow Y_i$ in ${\calO'}^{\otimes}$ covering the maps $\Colp{i}: \seg{n} \rightarrow \seg{1}$ in $\Nerve(\FinSeg)$. Since $D''$ and $D'$ preserve inert morphisms, we deduce that
the vertical maps in the diagram
$$ \xymatrix{ D''(Y) \ar[r]^{g} \ar[d] & D'(Y) \ar[d] \\
D''(Y_i) \ar[r]^{g_i} & D'(Y_i) }$$
are inert. Since $\calC^{\otimes}$ is an $\infty$-operad, we conclude that
$g$ is an equivalence if and only if each $g_i$ is an equivalence. The desired
result now follows from the observation that each $Y_i \in \calO'$.
\end{proof}

\begin{example}\label{slabber}
Let $\calO^{\otimes}$ and $\calC^{\otimes}$ be any $\infty$-operads.
We have a canonical $\infty$-operad bifunctor
$$ \calO^{\otimes} \times \CommOp
\rightarrow \CommOp \times \CommOp \rightarrow \CommOp,$$
where the second map is induced by the functor $\wedge: \FinSeg \times \FinSeg \rightarrow \FinSeg$.
Applying Example \ref{clise}, we obtain an $\infty$-operad $\Alg_{\calO}(\calC)^{\otimes}$ whose underlying $\infty$-category agrees with the $\infty$-category $\Alg_{\calO}(\calC)$ of
$\infty$-operad maps from $\calO^{\otimes}$ to $\calC^{\otimes}$. If
$\calC^{\otimes}$ is a symmetric monoidal $\infty$-category, then Proposition \ref{slape}
ensures that $\Alg_{\calO}(\calC)^{\otimes}$ is again symmetric monoidal, and that for
each $X \in \calO$ the evaluation functor
the forgetful functor $\Alg_{\calO}(\calC)^{\otimes} \rightarrow \calC^{\otimes}$ is symmetric monoidal.

In particular, if $\calC^{\otimes}$ is an $\infty$-operad, then we can define a new
$\infty$-operad $\CAlg(\calC)^{\otimes} = \Alg_{\CommOp}(\calC)^{\otimes}$. If
$\calC^{\otimes}$ is a symmetric monoidal $\infty$-category, then $\CAlg(\calC)^{\otimes}$ has the same property.
\end{example}

\begin{example}
Let $\calO^{\otimes}$ be any $\infty$-operad. Combining Example \ref{kapf} and
Remark \ref{stoss}, we deduce that the inclusion $\{ \seg{1} \} \subseteq \Triv$
induces equivalences of preoperads
$$ \calO^{\otimes, \natural} \simeq \calO^{\otimes, \natural} \odot \{ \seg{1} \}^{\flat}
\rightarrow \calO^{\otimes, \natural} \odot \Triv^{\natural}.$$
$$ \calO^{\otimes, \natural} \simeq \{ \seg{1} \}^{\flat} \odot \calO^{\otimes, \natural}
\rightarrow \Triv^{\natural} \odot  \calO^{\otimes, \natural}.$$
In other words, we can regard the trivial $\infty$-operad $\Triv$ as the unit object
of $\Cat_{\infty}^{\lax}$ with respect to the monoidal structure determined by
the operation $\odot$ of Proposition \ref{candlestick}.
\end{example}

\begin{example}\label{kilm}
The functor $\wedge: \Nerve(\FinSeg) \times \Nerve(\FinSeg) \rightarrow \Nerve(\FinSeg)$
is a bifunctor of $\infty$-operads. Moreover, $\wedge$ induces a weak equivalence
of $\infty$-preoperads $\Nerve(\FinSeg)^{\natural} \odot \Nerve(\FinSeg)^{\natural} \rightarrow
\Nerve(\FinSeg)^{\natural}$. To prove this, it suffices to show that for any
$\infty$-operad $\calC$, the functor $\wedge$ induces an equivalence
of $\infty$-categories $\CAlg(\calC) \rightarrow \CAlg( \CAlg(\calC) )$.
This follows from Corollary \ref{kime2}, since the $\infty$-operad $\CAlg(\calC)^{\otimes}$ is coCartesian
(Proposition \ref{socarte}).
\end{example}

Let $V: \Cat_{\infty}^{\lax} \rightarrow \Cat_{\infty}^{\lax}$ be induced by the left Quillen
functor $\bigdot \odot \Nerve(\FinSeg)^{\natural}$ from $\PreOp$ to itself. It follows from Example
\ref{kilm} and Proposition \toposref{recloc} that $V$ is a localization functor on $\Cat_{\infty}^{\lax}$.
The following result characterizes the essential image of $V$:

\begin{proposition}
Let $\calC^{\otimes}$ be an $\infty$-operad. The following conditions are equivalent:
\begin{itemize}
\item[$(1)$] The $\infty$-operad $\calC^{\otimes}$ lies in the essential image of the localization
functor $V$ defined above: in other words, there exists another $\infty$-operad
$\calD^{\otimes}$ and a bifunctor of $\infty$-operads $\theta: \Nerve(\FinSeg) \times
\calD^{\otimes} \rightarrow \calC^{\otimes}$ which induces a weak equivalence of $\infty$-preoperads.

\item[$(2)$] The $\infty$-operad $\calC^{\otimes}$ is coCartesian.
\end{itemize}
\end{proposition}

\begin{proof}
Suppose first that $(1)$ is satisfied. We will prove that $\calC^{\otimes}$ is equivalent
to $\calD^{\coprod}$. Let $\theta_0: \Nerve(\FinSeg) \times \calD \rightarrow \calC^{\otimes}$
be the restriction of $\theta$. In view of Theorem \ref{kinema}, it will suffice to show that
for every $\infty$-operad $\calE^{\otimes}$, composition with $\theta_0$ induces an equivalence
of $\infty$-categories $\Alg_{\calC}(\calE) \rightarrow \Fun( \calD, \CAlg(\calE) )$.
This map factors as a composition
$$ \Alg_{\calC}(\calE) \stackrel{\phi}{\rightarrow} \Alg_{\calD}( \CAlg(\calE) ) \stackrel{\phi'}{\rightarrow} \Fun( \calD, \CAlg(\calE)).$$  
Our assumption that $\theta$ induces a weak equivalence of $\infty$-preoperads guarantees
that $\phi$ is an equivalence of $\infty$-categories. To prove that $\phi'$ is a weak equivalence,
it suffices to show that $\CAlg(\calE)^{\otimes}$ is a coCartesian $\infty$-operad
(Proposition \ref{kime1}), which follows from Proposition \ref{socarte}.

Conversely, suppose that $\calC^{\otimes}$ is coCartesian. We will prove that the canonical
map $\calC^{\otimes} \times \{ \seg{1} \} \rightarrow \calC^{\otimes} \times \Nerve(\FinSeg)$
induces a weak equivalence of $\infty$-preoperads $\calC^{\otimes, \natural}
\simeq \calC^{\otimes, \natural} \odot \Nerve(\FinSeg)^{\natural}$. Unwinding the definitions,
it suffices to show that for every $\infty$-operad $\calD^{\otimes}$, the canonical map
$$ \Alg_{\calC}( \CAlg(\calD)) \rightarrow \Alg_{\calC}( \calD )$$
is an equivalence of $\infty$-categories. Invoking Theorem \ref{kinema}, we can reduce
to proving that the forgetful functor
$\psi: \Fun( \calC, \CAlg(\CAlg(\calD)) \rightarrow \Fun( \calC, \CAlg(\calD))$ is
an equivalence of $\infty$-categories. This is clear, since $\psi$ is a left
inverse to the functor induced by the equivalence $\CAlg(\calD) \rightarrow \CAlg(\CAlg(\calD))$ of
Example \ref{kilm}.
\end{proof}

\subsection{Unital $\infty$-Operads}\label{comm1.9}

Our goal in this section is to introduce a special class of $\infty$-operads, which we will
call {\it unital} $\infty$-operads. We will also establish a basic result about
$\calO$-algebra objects in the case where $\calO$ is unital (Proposition \ref{kime1}), which was already used in \S \ref{jilun}. 

We begin with an elementary observation about the
$\infty$-operad $\OpE{0}$. Note that the functor $\wedge: \FinSeg \times \FinSeg \rightarrow \FinSeg$
of Notation \ref{crumb} restricts to a functor $\InjSeg \times \InjSeg \rightarrow \InjSeg$. Passing to nerves, we obtain a bifunctor of $\infty$-operads $\OpE{0} \times \OpE{0} \rightarrow \OpE{0}$.

\begin{proposition}\label{jaw}
The above construction determines an equivalence of $\infty$-preoperads
$$ \OpE{0}^{\natural} \odot \OpE{0}^{\natural} \rightarrow \OpE{0}^{\natural}.$$
\end{proposition}

\begin{proof}
Consider the map $f: \Delta^1 \rightarrow \OpE{0}$ determined by the morphism
$\seg{0} \rightarrow \seg{1}$ in $\InjSeg$. Example \ref{capf} asserts that
$f$ induces a weak equivalence of $\infty$-preoperads
$(\Delta^1)^{\flat} \rightarrow \OpE{0}^{\natural}$. It follows from Remark \ref{stoss}
that the induced map
$$ (\Delta^1)^{\flat} \odot ( \Delta^1)^{\flat} \rightarrow \OpE{0}^{\natural}
\odot \OpE{0}^{\natural}$$
is a weak equivalence of $\infty$-preoperads. By a two-out-of-three argument, we are reduced
to showing that the induced map $( \Delta^1)^{\flat} \odot (\Delta^1)^{\flat} \rightarrow \OpE{0}^{\natural}$ is a weak equivalence of $\infty$-preoperads. Example \ref{capf} asserts that the composition
of this map with the diagonal $\delta: (\Delta^1)^{\flat} \rightarrow ( \Delta^1)^{\flat} \odot (\Delta^1)^{\flat}$ is a weak equivalence of $\infty$-preoperads. Using a two-out-of-three argument again,
we are reduced to proving that $\delta$ is a weak equivalence. 

Unwinding the definitions, it suffices to show the following: for every $\infty$-operad
$p: \calO^{\otimes} \rightarrow \Nerve(\FinSeg)$, composition with $\delta$ induces a trivial Kan fibration
$$\Fun_{ \Nerve(\FinSeg)}( \Delta^1 \times \Delta^1, \calO^{\otimes}) 
\rightarrow \Fun_{ \Nerve(\FinSeg)}( \Delta^1, \calO^{\otimes}).$$
This follows from Proposition \toposref{lklk}, since every functor $F \in \Fun_{ \Nerve(\FinSeg)}( \Delta^1 \times \Delta^1, \calO^{\otimes})$ is a $p$-left Kan extension of $F \circ \delta$ (because
every morphism in $\calO^{\otimes}_{\seg{0}}$ is an equivalence).
\end{proof}

\begin{corollary}\label{stass}
Let $i: \{ \seg{0} \}^{\flat} \rightarrow \OpE{0}^{\natural}$ denote the inclusion. Then
composition with $i$ induces weak equivalences of $\infty$-preoperads
$$ \OpE{0}^{\natural} \simeq \OpE{0}^{\natural} \odot \{ \seg{0} \}^{\flat} \rightarrow
\OpE{0}^{\natural} \odot \OpE{0}^{\natural}$$
$$ \OpE{0}^{\natural} \simeq \{ \seg{0} \}^{\flat} \odot \OpE{0}^{\natural} \rightarrow
\OpE{0}^{\natural} \odot \OpE{0}^{\natural}.$$
\end{corollary}

\begin{corollary}\label{kuj}
Let $\Cat_{\infty}^{\lax}$ denote the $\infty$-category of (small) $\infty$-operads,
which we identify with the underlying $\infty$-category $\Nerve( \PreOp^{\degree})$ of
the monoidal simplicial model category $\PreOp$. Let $U: \Cat_{\infty}^{\lax} \rightarrow
\Cat_{\infty}^{\lax}$ be induced by the left Quillen functor
$\overline{X} \mapsto \overline{X} \odot \OpE{0}^{\natural}$. Then $U$
is a localization functor from $\Cat_{\infty}^{\lax}$ to itself.
\end{corollary}

\begin{proof}
Combine Corollary \ref{stass} with Proposition \toposref{recloc}.
\end{proof}

Our next goal is to describe the essential image of the localization functor
$U: \Cat_{\infty}^{\lax} \rightarrow \Cat_{\infty}^{\lax}$ of Corollary \ref{kuj}.

\begin{definition}
We will say that an $\infty$-operad $\calO^{\otimes}$ is {\it unital} if, for every object
$X \in \calO$, the space $\Mul_{\calO}( \emptyset, X)$ is contractible.
\end{definition}

\begin{example}
The commutative and associative $\infty$-operads are unital. The $\infty$-operad
$\OpE{0}$ is unital. The trivial $\infty$-operad is not unital.
\end{example}

\begin{proposition}\label{skymall}
Let $\calO^{\otimes}$ be an $\infty$-operad. The following conditions are equivalent:
\begin{itemize}
\item[$(1)$] The $\infty$-category $\calO^{\otimes}$ is pointed (that is, there exists an
object of $\calO^{\otimes}$ which is both initial and final).
\item[$(2)$] The $\infty$-operad $\calO^{\otimes}$ is unital.
\item[$(3)$] The $\infty$-operad $\calO^{\otimes}$ lies in the essential image of the localization
functor $U:  \Cat_{\infty}^{\lax} \rightarrow \Cat_{\infty}^{\lax}$ of Corollary \ref{kuj}.
\end{itemize}
\end{proposition}

The proof depends on the following preliminary results:

\begin{lemma}\label{sad}
Let $p: \calO^{\otimes} \rightarrow \Nerve(\FinSeg)$ be an $\infty$-operad.
Then an object $X$ of the $\infty$-category $\calO^{\otimes}$ is final if and only if
$p(X) = \seg{0}$. Moreover, there exists an object of $\calO^{\otimes}$ satisfying this condition.
\end{lemma}

\begin{proof}
Since $\calO^{\otimes}$ is an $\infty$-operad, we have an equivalence
$\calO^{\otimes}_{\seg{0}} \simeq \calO^{0} \simeq \Delta^0$; this proves the
existence of an object $X \in \calO^{\otimes}$ such that $p(X) = \seg{0}$.
Since $\calO^{\otimes}$ is an $\infty$-operad, the object $X \in \calO^{\otimes}$ is
$p$-final. Since $p(X) = \seg{0}$ is a final object of $\Nerve(\FinSeg)$, it follows that
$X$ is a final object of $\calO^{\otimes}$. 

To prove the converse, suppose that $X'$ is any final object of $\calO^{\otimes}$.
Then $X' \simeq X$ so that $p(X') \simeq \seg{0}$. It follows that $p(X') = \seg{0}$, as desired.
\end{proof}

\begin{proof}[Proof of Proposition \ref{skymall}]
We first show that $(1) \Leftrightarrow (2)$. According to Lemma \ref{sad}, the $\infty$-category $\calO^{\otimes}$ admits a final object
$Y$. Assertion $(1)$ is equivalent to the requirement that $Y$ is also initial: that is, that
the space $\bHom_{ \calO^{\otimes}}( Y, X)$ is contractible for every $X \in \calO^{\otimes}$.
Let $\seg{n}$ denote the image of $X$ in $\Nerve(\FinSeg)$, and choose inert morphisms
$X \rightarrow X_i$ covering $\Colp{i}: \seg{n} \rightarrow \seg{1}$ for $1 \leq i \leq n$. 
Since $X$ is a $p$-limit of the diagram $\{ X_i \}_{1 \leq i \leq n} \rightarrow \calO^{\otimes}$, we conclude that $\bHom_{\calO^{\otimes}}( Y, X) \simeq \prod_{1 \leq i \leq n} \bHom_{\calO^{\otimes}}( Y, X_i )$. Assertion $(1)$ is therefore equivalent to the requirement that
$\bHom_{ \calO^{\otimes}}( Y, X)$ is contractible for $X \in \calO$, which is a rewording of condition $(2)$.

We now show that $(1) \Rightarrow (3)$. Let $\alpha: \calO^{\otimes} \rightarrow U \calO^{\otimes}$
be a morphism in $\Cat_{\infty}^{\lax}$ which exhibits $U \calO^{\otimes}$ as a $U$-localization
of $\calO^{\otimes}$. We will prove that there exists a morphism $\beta: U \calO^{\otimes} \rightarrow \calO^{\otimes}$ such that $\beta \circ \alpha$ is equivalent to $\id_{ \calO^{\otimes} }$. We claim
that $\beta$ is a homotopy inverse to $\alpha$: to prove this, it suffices to show that
$\alpha \circ \beta$ is homotopic to the identity $\id_{ U \calO^{\otimes} }$. Since
$U \calO^{\otimes}$ is $U$-local and $\alpha$ is a $U$-equivalence, it suffices to show that
$\alpha \circ \beta \circ \alpha$ is homotopic to $\alpha$, which is clear.

To construct the map $\beta$, we observe that $U \calO^{\otimes}$ can be identified
with a fibrant replacement for the object $\calO^{\otimes, \natural} \odot (\Delta^1)^{\flat} \in \PreOp$
(here we regard $( \Delta^1)^{\flat}$ as an $\infty$-preoperad as in the proof of Proposition \ref{jaw}). 
It will therefore suffice to construct a map $h: \calO^{\otimes} \times \Delta^1 \rightarrow \calO^{\otimes}$
such that $h| (\calO^{\otimes} \times \{1\}) = \id_{\calO^{\otimes}}$ and $h | (\calO^{\otimes} \times \{0\} )$
factors through $\calO^{\otimes}_{\seg{0}}$. The existence of $h$ follows immediately from assumption $(1)$.

To show that $(3) \Rightarrow (1)$, we reverse the above reasoning: if $\calO^{\otimes}$ is $U$-local, then there exists a morphism $\beta: U \calO^{\otimes} \rightarrow \calO^{\otimes}$ which is
right inverse to $\alpha$, which is equivalent to the existence of a map $h: \calO^{\otimes} \times \Delta^1 \rightarrow \calO^{\otimes}$ as above. We may assume without loss of generality
that $h| ( \calO^{\otimes} \times \{0\})$ is the constant map taking some value $X \in \calO^{\otimes}_{\seg{0}}$. Then $h$ can be regarded as a section of the left fibration
$( \calO^{\otimes} )^{X/} \rightarrow \calO^{\otimes}$. This proves that $X$ is an initial
object of $\calO^{\otimes}$. Since $X$ is also a final object of $\calO^{\otimes}$, we deduce that
$\calO^{\otimes}$ is pointed as desired.
\end{proof}

\begin{proposition}\label{kime1}
Let $\calO^{\otimes}$ be a unital $\infty$-operad. Let $\calC$ be an $\infty$-category, and
regard $\calC$ as the $\infty$-category underlying the $\infty$-operad $\calC^{\amalg}$.
Then the restriction functor $\Alg_{\calO}(\calC) \rightarrow \Fun( \calO, \calC)$
is an equivalence of $\infty$-categories.
\end{proposition}

\begin{corollary}\label{kime2}
Let $\calC$ be an $\infty$-category, and regard
$\calC$ as the underlying $\infty$-category of the $\infty$-operad $\calC^{\amalg}$. Then the forgetful functor $\CAlg( \calC) \rightarrow \calC$ is a trivial Kan fibration.
\end{corollary}

\begin{proof}[Proof of Proposition \ref{kime1}]
Let $\calD$ denote the fiber product $$\calO^{\otimes} \times_{ \Nerve(\FinSeg)} \Nerve(\PFinSeg).$$
By definition, a map $\calO^{\otimes} \rightarrow \calC^{\amalg}$ in
$(\sSet)_{/ \Nerve(\FinSeg)}$ can be identified with a functor
$A: \calD \rightarrow \calC$.
Such a functor determines a map of $\infty$-operads if and only if the following
condition is satisfied:
\begin{itemize}
\item[$(\ast)$] Let $\alpha$ be a morphism in $\calD$ whose image in $\Nerve(\FinSeg)$ is inert. 
Then $A(\alpha)$ is an equivalence in $\calC$.
\end{itemize}
We can identify $\Alg_{\calO}(\calC)$ with the full subcategory of $\Fun( \calD, \calC)$
spanned by those functors which satisfy $(\ast)$.

We observe that the inverse image in $\calD$ of $\seg{1} \in \Nerve(\FinSeg)$ is canonically isomorphic to $\calO$. Via this isomorphism, we will regard
$\calO$ as a full subcategory of $\calD$.
In view of Proposition \toposref{lklk}, it will suffice to prove the following:

\begin{itemize}
\item[$(a)$] A functor $A: \calD \rightarrow \calC$
is a left Kan extension of $A | \calO$ if and only if it satisfies condition $(\ast)$.
\item[$(b)$] Every functor $A_0: \calO \rightarrow \calC$ admits an extension
$A:  \calD \rightarrow \calC$ satisfying the equivalent conditions of $(a)$.
\end{itemize}

We can identify objects of $\calD$
with pairs $(X,i)$, where $X \in \calO^{\otimes}_{\seg{n}}$ and $1 \leq i \leq n$.
For any such pair, we can choose an inert morphism $X \rightarrow X_i$ lying over
$\Colp{i}: \seg{n} \rightarrow \seg{1}$. We then have a morphism
$f: (X, i) \rightarrow (X_i, 1)$. Using the assumption that $\calO^{\otimes}$ is unital, we deduce that the map
$$\bHom_{  \calD}( Y, (X,i) )
\rightarrow \bHom_{\calD}( Y, (X_i,1))
\simeq \bHom_{\calO}( Y, X_i)$$ is a homotopy equivalence for
each $Y \in \calO \subseteq \calD$.
In particular, we conclude that $f$ admits a right homotopy inverse
$g: (X_i,1) \rightarrow (X,i)$. It follows that composition with $g$ induces a homotopy
equivalence
$$\bHom_{ \calD}( Y, (X_i,1))
\rightarrow \bHom_{  \calD }( Y, (X,i) )
$$
for each $Y \in \calO$. This implies that the inclusion
$\calO \subseteq \calD$ admits a right adjoint $G$, given by
$(X,i) \mapsto (X_i, 1)$. This immediately implies $(b)$ (we can take $A = A_0 \circ G$) together with the following version of $(a)$: 
\begin{itemize}
\item[$(a')$] A functor $A: \calD \rightarrow \calC$
is a left Kan extension of $A | \calO$ if and only if, for every
object $(X,i) \in \calD$, the map $A(g)$ is an equivalence in $\calC$, where
$g: (X_i, 1) \rightarrow (X,i)$ is defined as above.
\end{itemize}
Since $g$ is a right homotopy inverse to the inert morphism $f: (X,i) \rightarrow (X_i, 1)$, 
assertion $(a')$ can be reformulated as follows: a functor $A: \calD \rightarrow \calC$ is
a left Kan extension of $A | \calO$ if and only if the following condition is satisfied:
\begin{itemize}
\item[$(\ast')$] Let $(X,i) \in \calD$ be an object, and let $f: (X,i) \rightarrow (X_i,1)$ be
defined as above. Then $A(f)$ is an equivalence in $\calC$.
\end{itemize}

To complete the proof, it will suffice to show that conditions $(\ast)$ and
$(\ast')$ are equivalent. The implication $(\ast) \Rightarrow (\ast')$ is obvious.
For the converse, suppose that $h: (Y,j) \rightarrow (X,i)$ is an arbitrary morphism in
$\calD$ whose image in $\calO^{\otimes}$ is inert. We then have a commutative diagram
$$ \xymatrix{ (Y,j) \ar[dr]^{f'} \ar[rr]^{f} & & (X,i) \ar[dl]^{f''} \\
& (X_i, 1) & }$$
Condition $(\ast')$ guarantees that $A(f')$ and $A(f'')$ are equivalences, so that
$A(f)$ is an equivalence by the two-out-of-three property.
\end{proof}

\subsection{$\calO$-Operad Families}\label{comm1.10}

Let $\calO^{\otimes}$ be an $\infty$-operad. Then, for each $n \geq 0$, we have a canonical
equivalence of $\infty$-categories $\calO^{\otimes}_{\seg{n}} \simeq \calO^{n}$. 
Our goal in this section is to introduce a mild generalization of the notion of an
$\infty$-operad, where we replace the absolute $n$th power $\calO^{n}$ with
the $n$th fiber power over some base simplicial set $S \simeq \calO^{\otimes}_{\seg{0}}$. 
The discussion in this section is of a somewhat technical nature and we recommend
that it be skipped at a first reading (these ideas will be needed in \S \ref{comm6}). 

It will be convenient for us to employ the formalism of {\it categorical patterns}
developed in \S \bicatref{bisec3.2}; in what follows, we will assume that the reader is familiar with
the terminology described there.

\begin{definition}\label{suskup}
Let $\calO^{\otimes}$ be an $\infty$-operad. We define a categorical pattern
(see Definition \bicatref{catpat}) $\CatP_{\calO}^{\famm} = (M, T, \{ p_{\alpha}: \Delta^1 \times \Delta^1 \rightarrow \calO^{\otimes}\}_{\alpha \in A} )$ on $\calO^{\otimes}$ as follows:
\begin{itemize}
\item[$(1)$] The set $M$ consists of all inert morphisms in $\calO^{\otimes}$.
\item[$(2)$] The set $T$ consists of all $2$-simplices in $\calO^{\otimes}$.
\item[$(3)$] The set $A$ parametrizes all diagrams $\Delta^1 \times \Delta^1 \rightarrow \calO^{\otimes}$
of inert morphisms having the property that the induced diagram of inert morphisms in
$\FinSeg$
$$ \xymatrix{ \seg{m} \ar[r] \ar[d] & \seg{n} \ar[d] \\
\seg{m'} \ar[r] & \seg{n'} }$$
induces a bijection of finite sets $\nostar{m'} \coprod_{ \nostar{n'}} \nostar{n} \rightarrow \nostar{m}$.
\end{itemize}

We will say that a map of simplicial sets $\calC^{\otimes} \rightarrow \calO^{\otimes}$
exhibits $\calC^{\otimes}$ as a {\it $\calO$-operad family} it it is $\CatP_{\calO}^{\famm}$-fibered,
in the sese of Remark \bicatref{caroline}. In the special case $\calO^{\otimes} = \Nerve(\FinSeg)$, we will simply say that $\calC^{\otimes}$ is an {\it $\infty$-operad family}.
\end{definition}

\begin{definition}
Let $\calO^{\otimes}$ be an $\infty$-operad, and let $p: \calC^{\otimes} \rightarrow \calO^{\otimes}$
be an $\infty$-operad family. We will say that a morphism $f$ in $\calC^{\otimes}$ is
{\it inert} if $p(f)$ is an inert morphism of $\calO^{\otimes}$ and $f$ is $p$-coCartesian.
\end{definition}

The following basic result guarantees that $\infty$-operad families are not far from
being $\infty$-operads themselves. 

\begin{proposition}\label{prejule}
Let $p: \calC^{\otimes} \rightarrow \Nerve(\FinSeg)$ be a map of simplicial sets.
The following conditions are equivalent:
\begin{itemize}
\item[$(1)$] The map $p$ exhibits $\calC^{\otimes}$ as an $\infty$-operad family, and
the fiber $\calC^{\otimes}_{ \seg{0}}$ is a contractible Kan complex.
\item[$(2)$] The map $p$ exhibits $\calC^{\otimes}$ as an $\infty$-operad.
\end{itemize}
\end{proposition}

We will give the proof of Proposition \ref{prejule} at the end of this section.

\begin{remark}\label{casta}
Let $\calO^{\otimes}$ be an $\infty$-operad, and let $\calD$ be any $\infty$-category.
Then the projection map $\calD \times \calO^{\otimes} \rightarrow \calO^{\otimes}$
exhibits $\calD \times \calO^{\otimes}$ as a $\calO$-operad family.
Moreover, a morphism in $\calD \times \calO^{\otimes}$ is inert if and only if
its image in $\calO^{\otimes}$ is inert and its image in $\calD$ is an equivalence. To prove this,
it suffices to show that the product functor $\overline{X} \mapsto \overline{X} \times ( \calO^{\otimes}, M)$
is a right Quillen functor from the category $\mSet$ of marked simplicial sets to
$\mset{ \CatP^{\famm}_{\calO} }$; here $M$ denotes the collection of all inert morphisms in $\calO^{\otimes}$.
Equivalently, it suffices to show that the forgetful functor
$\mset{ \CatP^{\famm}_{\calO}} \rightarrow \mSet$ is a left Quillen functor. This follows from
Example \bicatref{megaparnum}.
\end{remark}

\begin{proposition}\label{supercae}
Let $\calO^{\otimes}$ be an $\infty$-operad, let
$p: \calC^{\otimes} \rightarrow \calO^{\otimes}$ be an $\infty$-operad family.
Suppose we are given a pullback diagram of $\infty$-categories
$$ \xymatrix{ {\calC'}^{\otimes} \ar[r] \ar[d] & \calC^{\otimes} \ar[d]^{q} \\
\calD' \ar[r] & \calD }$$
where $q \times p: \calC^{\otimes} \rightarrow \calD \times \calO^{\otimes}$ is a categorical fibration and $q$ carries inert morphisms in
$\calC^{\otimes}$ to equivalences in $\calD$.
Then:
\begin{itemize}
\item[$(1)$] The $\infty$-category ${\calC'}^{\otimes}$ is an $\calO$-operad family, and a morphism in
${\calC'}^{\otimes}$ is inert if and only if its image in $\calC^{\otimes}$ is inert and its image
in $\calD'$ is an equivalence.
\item[$(2)$] Suppose that $\calD'$ is a contractible Kan complex and that $q$ induces an equivalence
of $\infty$-categories $\calC^{\otimes}_{\seg{0}} \rightarrow \calD$. Then
the functor ${ \calC'}^{\otimes} \rightarrow \calO^{\otimes}$ is a fibration of $\infty$-operads.
\end{itemize}
\end{proposition}

\begin{proof}
Let $M_{\calC}$ denote the collection of inert morphisms in $\calC^{\otimes}$, let
$M_{\calO}$ denote the collection of inert morphisms in $\calO^{\otimes}$.
For each $\infty$-category $\calE$, let $\calE^{\natural}$ denote the marked simplicial
set $(\calE, M)$, where $M$ denotes the collection of all equivalences in $\calE$.
Assertion $(1)$ is equivalent to the requirement that the fiber product
$$ ( {\calD'}^{\natural} \times (\calO^{\otimes}, M_{\calO}))
\times_{ \calD^{\natural} \times ( \calO^{\otimes}, M_{\calO})}
(\calC^{\otimes}, M_{\calC})$$
is a fibrant object of $\mset{ \CatP^{\famm}_{\calO}}$. Since
$(\calD')^{\natural} \times (\calO^{\otimes}, M_{\calO})$ is fibrant
(Remark \ref{casta}), it will suffice to show that $q \times p$ induces a fibration
$$ ( \calC^{\otimes}, M_{\calC}) \rightarrow \calD^{\natural} \times (\calO^{\otimes}, M_{\calO}).$$
Since the domain and codomain of this morphism are both fibrant (using Remark \ref{casta} again),
this is equivalent to our assumption that $q \times p$ is a categorical fibration
(Proposition \bicatref{fibraman}). This proves $(1)$. 

Using Proposition \bicatref{fibraman} again, we deduce that the underlying map of simplicial sets
${\calC'}^{\otimes} \rightarrow \calD' \times \calO^{\otimes}$ is a categorical fibration, so the map
${\calC'}^{\otimes} \rightarrow \calO^{\otimes}$ is a categorical fibration. To prove
$(2)$, it will suffice to show that ${\calC'}^{\otimes}$ is an $\infty$-operad. We have a pullback diagram
$$ \xymatrix{ {\calC'}^{\otimes}_{\seg{0}} \ar[r] \ar[d] & \calC^{\otimes}_{\seg{0}} \ar[d]^{q_0} \\
\calD' \ar[r] & \calD. }$$
Since $p \times q$ is a categorical fibration, the map $q_0$ is a categorical fibration.
If $q_0$ is also an equivalence of $\infty$-categories, then it is a trivial Kan fibration.
If $\calD'$ is a contractible Kan complex, we conclude that ${\calC'}^{\otimes}_{\seg{0}}$ is a contactible Kan complex, so that ${\calC'}^{\otimes}$ is an $\infty$-operad by virtue of Proposition \ref{prejule}.
\end{proof}

\begin{definition}\label{hunker}
Let $\calO^{\otimes}$ be an $\infty$-operad, and let $\calC^{\otimes} \rightarrow \calO^{\otimes}$
be an $\calO$-operad family. We let $\Alg_{\calO}(\calC)$ denote the full subcategory of
$\Fun_{ \calO^{\otimes}}( \calO^{\otimes}, \calC^{\otimes})$ spanned by those functors which preserve inert morphisms. 
\end{definition}

\begin{remark}
In the special case where $\calC^{\otimes} \rightarrow \calO^{\otimes}$ is a
fibration of $\infty$-operads, Definition \ref{hunker} agrees with Definition \ref{susqat}.
\end{remark}

\begin{example}\label{insert}
Let $\calO^{\otimes}$ be an $\infty$-operad and let $\calD$ be an $\infty$-category,
and regard $\calC^{\otimes} = \calD \times \calO^{\otimes}$ as a $\calO$-operad family (Remark \ref{casta}). Unwinding the definitions, we see that $\Alg_{\calO}(\calC)$ can be identified
with the full subcategory of $\Fun( \calO^{\otimes}, \calD)$ spanned by those functors
$F$ which carry inert morphisms in $\calO^{\otimes}$ to equivalences in $\calD$. 
This condition is equivalent to the requirement that $F$ be a right Kan extension
of its restriction to the contractible Kan complex $\calO^{\otimes}_{\seg{0}}$. Using
Proposition \toposref{lklk}, we deduce that evaluation at any object
$X \in \calO^{\otimes}_{\seg{0}}$ induces a trivial Kan fibration
$\Alg_{\calO}( \calC) \rightarrow \calD$. This Kan fibration has a section $s: \calD \rightarrow \Alg_{\calO}(\calC)$, given by the diagonal embedding $\calD \rightarrow \Fun( \calO^{\otimes}, \calD)$.
It follows that $s$ is an equivalence of $\infty$-categories.
\end{example}

\begin{proof}[Proof of Proposition \ref{prejule}]
We first prove that $(1) \Rightarrow (2)$. Fix an object
$X \in \calC^{\otimes}_{\seg{n}}$, and choose $p$-coCartesian morphisms
$f_i: X \rightarrow X_i$ covering the maps $\Colp{i}: \seg{n} \rightarrow \seg{1}$
in $\Nerve(\FinSeg)$. These maps determine a diagram 
$q: \nostar{n}^{\triangleleft} \rightarrow \calC^{\otimes}$; we wish to prove that
$q$ is a $p$-limit diagram. The proof proceeds by induction on $n$. 

If $n = 0$,
then we must show that every object $X \in \calC^{\otimes}_{\seg{0}}$ is $p$-final.
In other words, we must show that for every object $Y \in \calC^{\otimes}_{\seg{m}}$, 
the homotopy fiber of the map
$$ \bHom_{ \calC^{\otimes}}( Y, X) \rightarrow \bHom_{ \Nerve(\FinSeg)}( \seg{m}, \seg{0})$$
is contractible. Fix a map $\alpha: \seg{m} \rightarrow \seg{0}$ (in fact, there is a unique such morphism); since $\alpha$ is inert, we can choose a $p$-coCartesian morphism $Y \rightarrow Y'$ lifting
$\alpha$. In view of Proposition \toposref{compspaces}, we can identify the relevant homotopy fiber with
$\bHom_{ \calC^{\otimes}_{\seg{0}}}( Y', X)$. This space is contractible by virtue of our assumption that $\calC^{\otimes}_{\seg{0}}$ is a contractible Kan complex.

If $n = 1$, there is nothing to prove. Assume that $n > 1$. Let $\beta: \seg{n} \rightarrow \seg{n-1}$
be defined by the formula
$$ \beta(i) = \begin{cases} i & \text{if } 1 \leq i \leq n-1 \\
\ast & \text{otherwise,} \end{cases}$$
and choose a $p$-coCartesian morphism $g: X \rightarrow X'$ lying over $\beta$.
Using the assumption that $g$ is $p$-coCartesian, we obtain factorizations of
$f_i$ as a composition
$$ X \stackrel{g}{\rightarrow} X' \stackrel{f'_{i}}{\rightarrow} X_i $$ 
for $1 \leq i \leq n$. These factorizations determine a diagram
$$q': ( \nostar{n-1}^{\triangleleft} \coprod \{n\})^{\triangleleft} \rightarrow \calC^{\otimes}$$
extending $q$. Fix an object $X_0 \in \calC^{\otimes}_{\seg{0}}$. We have seen
that $X_0$ is a $p$-final object of $\calC^{\otimes}$. Since $\seg{0}$ is also a final
object of $\Nerve(\FinSeg)$, we deduce that $X_0$ is a final object of $\calC^{\otimes}$.
We may therefore extend $q'$ to a diagram
$$ q'': \{x\} \star (\nostar{n-1}^{\triangleleft} \coprod \{n\} ) \star \{x_0 \} \rightarrow \calC^{\otimes}$$
carrying $x$ to $X$ and $x_0$ to $X_0$.

In view of Lemma \toposref{kan0}, to prove that $q$ is a $p$-limit diagram it will suffice to show the following:
\begin{itemize}
\item[$(a)$] The restriction $q''| ( \nostar{ n-1 }^{\triangleleft} \coprod \{n\}) \star \{x_0 \}$ is
a $p$-right Kan extension of $q'' | \nostar{n}$.
\item[$(b)$] The diagram $q''$ is a $p$-limit.
\end{itemize}

Using Proposition \toposref{acekan}, we can break the proof of $(a)$ into two parts:
\begin{itemize}
\item[$(a')$] The restriction $q''| ( \nostar{n-1}^{\triangleleft} \coprod \{n\}) \star \{x_0\}$ is a $p$-right Kan extension of $q''| ( \nostar{n-1}^{\triangleleft} \coprod \{n\})$.
\item[$(a'')$] The restriction $q''| ( \nostar{n-1}^{\triangleleft} \coprod \{n\})$ is a $p$-right Kan
extension of $q'' | \nostar{n}$.
\end{itemize}
Assertion $(a')$ follows from the observation that $X_0$ is a $p$-final object of
$\calC^{\otimes}$, and $(a'')$ follows from the inductive hypothesis.

To prove $(b)$, we observe that the inclusion $( \emptyset^{\triangleleft} \coprod \{n\}) \star \{x_0\} \subseteq ( \nostar{n-1}^{\triangleleft} \coprod \{n\}) \star \{x_0\}$ is cofinal (for example, using Theorem \toposref{hollowtt}). Consequently, it suffices to show that the restriction of
$q''$ to $\{ x \} \star ( \emptyset^{\triangleleft} \coprod \{n\} ) \star \{x_0 \}$ is a $p$-limit diagram,
which follows from assumption $(2)$ (since this restriction of $q''$ is a $p$-coCartesian lift of the inert diagram
$$ \xymatrix{ \seg{n} \ar[r]^{\Colp{i}} \ar[d]^{\beta} & \seg{1} \ar[d] \\
\seg{n-1} \ar[r] & \seg{0}. }$$

To complete the proof that $(1) \Rightarrow (2)$, we must show that for
each $n \geq 0$, the maps $\Colll{i}{!} : \calC^{\otimes}_{ \seg{n}} \rightarrow
\calC^{\otimes}_{ \seg{1}}$ induce an equivalence of $\infty$-categories
$\theta_n: \calC^{\otimes}_{ \seg{n} } \rightarrow ( \calC^{\otimes}_{\seg{1}})^{n}$. The proof
again proceeds by induction on $n$. When $n=0$, this follows immediately
from the contractibility of $\calC^{\otimes}_{\seg{0}}$, and when $n=1$ there is nothing to prove. Assume therefore that $n \geq 2$ and observe that $\theta$ is equivalent to the composition
$$ \calC^{\otimes}_{\seg{n}} \stackrel{ \beta_{!} \times \alpha_{!} }{\rightarrow} \calC^{\otimes}_{ \seg{n-1}} \times \calC^{\otimes}_{\seg{1}} \stackrel{ \theta_{n-1} \times \id}{\rightarrow}
(\calC^{\otimes}_{\seg{1}})^{n},$$
where $\beta: \seg{n} \rightarrow \seg{n-1}$ is defined as above and
$\alpha = \Colp{n}$. By virtue of the inductive hypothesis, it suffices to show that
the map $\beta_{!} \times \alpha_{!}$ is an equivalence of $\infty$-categories.
We have a homotopy coherent diagram of $\infty$-categories
$$ \xymatrix{ \calC^{\otimes}_{\seg{n}} \ar[r]^{\beta_{!} } \ar[d]^{\alpha_{!}} & \calC^{\otimes}_{ \seg{n-1}} \ar[d] \\
\calC^{\otimes}_{ \seg{1} } \ar[r] & \calC^{\otimes}_{ \seg{0} }.}$$
Because $\calC^{\otimes}$ is an $\infty$-operad family, this square is a homotopy pullback.
Since $\calC^{\otimes}_{\seg{0}}$ is a contractible Kan complex, we conclude that
$\beta_{!} \times \alpha_{!}$ is a categorical equivalence as desired. This completes the proof
that $(1) \Rightarrow (2)$.

We now prove that $(2) \Rightarrow (1)$. Assume that $\calC^{\otimes}$ is an
$\infty$-operad. Then we have a categorical equivalence $\calC^{\otimes}_{\seg{0}} \rightarrow
(\calC^{\otimes}_{\seg{1}})^{0}$, so that $\calC^{\otimes}_{\seg{0}}$ is a contractible Kan complex.
It will therefore suffice to show that $\calC^{\otimes}$ is an $\infty$-operad family.
Fix a diagram $\sigma: \Delta^1 \times \Delta^1 \rightarrow \Nerve(\FinSeg)$
of inert morphisms
$$ \xymatrix{ \seg{n} \ar[r]^{\alpha} \ar[d]^{\beta} \ar[dr]^{\gamma} & \seg{m} \ar[d] \\
\seg{m'} \ar[r] & \seg{k} }$$
which induces a bijection $\nostar{m} \coprod_{ \nostar{k} } \nostar{m'}$.
We wish to prove the following:
\begin{itemize}
\item[$(i)$] Every map $\overline{\sigma}: \Delta^1 \times \Delta^1 \rightarrow \calC^{\otimes}$
lifting $\sigma$ which carries every morphism in $\Delta^1 \times \Delta^1$ to an inert
morphism in $\calC^{\otimes}$ is a $p$-limit diagram.

\item[$(ii)$] Let $\sigma_0$ denote the restriction of $\sigma$ to the full subcategory $K$ of
$\Delta^1 \times \Delta^1$ obtained by omitting the initial object. If
$\overline{\sigma}_0: K \rightarrow \calC^{\otimes}$ is a map lifting
$\sigma_0$ which carries every edge of $K$ to an inert morphism in $\calC^{\otimes}$, then
$\overline{\sigma}_0$ can be extended to a map $\overline{\sigma}: \Delta^1 \times \Delta^1 \rightarrow \calC^{\otimes}$ satisfying the hypothesis of $(i)$.
\end{itemize}

To prove these claims, consider the $\infty$-category
$\widetilde{A} = ( \Delta^1 \times \Delta^1) \star \nostar{n}$, and let
$A$ denote the subcategory obtained by removing those morphisms
of the form $(1,1) \rightarrow i$ where $i \in \gamma^{-1} \{ \ast \}$,
$(0,1) \rightarrow i$ where $i \in \beta^{-1} \{ \ast \}$, and $(1,0) \rightarrow i$ where
$i \in \alpha^{-1} \{ \ast \}$. We observe that $\sigma$ can be extended uniquely
to a diagram $\tau: A \rightarrow \Nerve(\FinSeg)$ such that
$\tau(i) = \seg{1}$ for $i \in \nostar{n}$, and $\tau$ carries
the morphism $(0,0) \rightarrow i$ to the map $\Colp{i}: \seg{n} \rightarrow \seg{1}$.
The assumption that $\nostar{n} \simeq \nostar{m} \coprod_{ \nostar{k} } \nostar{m'}$. 
guarantees that for each $i \in \nostar{n}$, the $\infty$-category
$(\Delta^1 \times \Delta^1) \times_{A} A_{/i}$ contains a final object
corresponding to a morphism $(j,j') \rightarrow i$ in $A$, where
$(j,j') \neq (0,0)$. Note that the image of this morphism in $\Nerve(\FinSeg)$ is inert.

Let $\overline{\sigma}: \Delta^1 \times \Delta^1$ be as $(i)$.
Using Lemma \toposref{kan2}, we can choose a $p$-left Kan extension
$\overline{\tau}: A \rightarrow \calC^{\otimes}$ of $\overline{\sigma}$ such that
$p \circ \overline{\tau} = \tau$. Let $A^{0}$ denote the full subcategory of
$A$ obtained by removing the object $(0,0)$. We observe that
the inclusion $K^{op} \subseteq (A^{0})^{op}$
is cofinal (Theorem \toposref{hollowtt}). Consequently, to prove that
$\overline{\sigma}$ is a $q$-limit diagram, it suffices to show that
$\overline{\tau}$ is a $q$-limit diagram. Since $\calC^{\otimes}$ is an
$\infty$-operad, the restriction of $\overline{\tau}$ to
$\{ (0,0) \} \star \nostar{n}$ is a $p$-limit diagram. To complete the proof,
it will suffice (by virtue of Lemma \toposref{kan0}) to show that
$\overline{\tau} | A^0$ is a $p$-right Kan extension of $\overline{\tau} | \nostar{n}$.
This again follows immediately from our assumption that $\calC^{\otimes}$ is an $\infty$-operad.

We now prove $(ii)$. Let $\overline{\sigma}_0: K \rightarrow \calC^{\otimes}$ be
as in $(ii)$. Using Lemma \toposref{kan2}, we can choose
a $p$-left Kan extension $\overline{\tau}_0: A^0 \rightarrow \calC^{\otimes}$ of
$\overline{\sigma}_0$ covering the map $\tau_0 = \tau| A^0$.
Using the assumption that $\calC^{\otimes}$ is an $\infty$-operad, we deduce
that $\overline{\tau}_0$ is a $p$-right Kan extension of
$\overline{\tau}_0 | \nostar{n}$, and that $\overline{\tau}_0 | \nostar{n}$
can be extended to a $p$-limit diagram $\widetilde{\tau}_0: \{ (0,0) \} \star \nostar{n} \rightarrow \calC^{\otimes}$ lifting $\tau | ( \{(0,0) \} \star \nostar{n} )$; moreover, any such diagram
carries each edge of $\{ (0,0) \} \star \nostar{n}$ to an inert morphism in $\calC^{\otimes}$.
Invoking Lemma \toposref{kan0}, we can amalgamate $\widetilde{\tau}_0$ and
$\overline{\tau}_0$ to obtain a diagram $\overline{\tau}: A \rightarrow \calC^{\otimes}$ covering
$\tau$. We claim that $\overline{\sigma} = \overline{\tau} | \Delta^1 \times \Delta^1$
is the desired extension of $\overline{\sigma}_0$. To prove this, it suffices to show that
$\overline{\tau}$ carries each morphism of $\Delta^1 \times \Delta^1$ to an inert
morphism of $\calC^{\otimes}$. Since the composition of inert morphisms
in $\calC^{\otimes}$ is inert, it will suffice to show that the maps
$$\overline{\tau}(0,1) \stackrel{\overline{\beta}}\leftarrow \overline{\tau}(0,0) \stackrel{\overline{\alpha}}{\rightarrow} \overline{\tau}(1,0)$$
are inert, where $\overline{\alpha}$ and $\overline{\beta}$ are the morphisms lying over
$\alpha$ and $\beta$ determined by $\tau$. We will prove that $\overline{\alpha}$ is inert;
the case of $\overline{\beta}$ follows by the same argument. We can factor
$\overline{\alpha}$ as a composition
$$ \overline{\tau}(0,0) \stackrel{ \overline{\alpha}'}{\rightarrow} \alpha_{!} \overline{\tau}(0,0)
\stackrel{ \overline{\alpha}''}{\rightarrow} \overline{\tau}(1,0).$$
We wish to prove that $\overline{\alpha}''$ is an equivalence in the $\infty$-category
$\calC^{\otimes}_{\seg{m}}$. Since $\calC^{\otimes}$ is an $\infty$-operad, it will suffice to show that for $1 \leq j \leq m$, the functor $\Colll{j}{!} : \calC^{\otimes}_{\seg{m}} \rightarrow \calC$
carries $\overline{\alpha}''$ to an equivalence in $\calC$. Unwinding the definitions, this is
equivalent to the requirement that the map $\overline{\tau}(0,0) \rightarrow \overline{\tau}(i)$
is inert, where $i = \alpha^{-1}(j) \in \nostar{n}$, which follows immediately from our construction.
\end{proof}

\section{Algebras}\label{comm4}

Let $\calC^{\otimes} \rightarrow \calO^{\otimes}$ be a fibration of $\infty$-operads.
In \S \ref{algobj}, we introduced the $\infty$-category $\Alg_{\calO}(\calC)$ of
$\calO$-algebra objects of $\calC$. Our goal in this section is to study these $\infty$-categories in more detail. In particular, we will study conditions which guarantee the existence (and allow for the computation of) limits and colimits in $\Alg_{\calO}(\calC)$.
We begin in \S \ref{limco1} with the study of limits in $\Alg_{\calO}(\calC)$. This is fairly straightforward: the basic result is that limits in $\Alg_{\calO}(\calC)$ can be computed in the underlying $\infty$-category $\calC$ (Proposition \ref{limycomps2}). 

The study of colimits is much more involved. First of all, we do not expect colimits
in $\Alg_{\calO}(\calC)$ to be computed in the underlying $\infty$-category $\calC$ in general.
This is sometimes true for colimits of a special type (namely, colimits of diagrams indexed by
{\em sifted} simplicial sets). However, even in this case the existence of colimits in $\Alg_{\calO}(\calC)$ requires strong assumptions than the existence of the relevant colimits in $\calC$: we must also assume that these colimits are {\em compatible} with the $\infty$-operad structure on $\calC$. In \S \ref{comm35} we will axiomatize the relevant compatibility by introducing the theory of {\em operadic} colimit diagrams.
As with the usual theory of colimits, there is a relative version of the notion of an operadic colimit diagram, which we will refer to as an {\it operadic left Kan extension}. We will study operadic
left Kan extensions in \S \ref{iljest}, and prove a basic existence and uniqueness result
(Theorem \ref{oplk}). The proof is somewhat technical, and involves a mild generalization of the theory of minimal $\infty$-categories which we review in \S \ref{strictun}.

Suppose we are given a map of $\infty$-operads ${\calO'}^{\otimes} \rightarrow \calO^{\otimes}$.
Composition with this map induces a forgetful functor $\theta: \Alg_{\calO}(\calC) \rightarrow \Alg_{\calO'}(\calC)$. In \S \ref{comm33}, we will apply the theory of operadic left Kan extensions to construct
a left adjoint to the forgetful functor $\theta$ (assuming that $\calC$ satisfies some reasonable hypotheses). In the special case where ${\calO'}^{\otimes}$ is the trivial $\infty$-operad, we can think of this left adjoint as a ``free algebra'' functor $F$. Though it is difficult to construct colimits of algebras in
general, it is often much easier to construct colimits of free algebras, since the functor
$F$ preserves colimits (because it is a left adjoint). In \S \ref{limco2} we will exploit this observation to construct general colimits in $\Alg_{\calO}(\calC)$: the basic idea is to resolve arbitrary algebras with free algebras. Although this strategy leads to a fairly general existence result (Corollary \ref{smurf}), 
it is somewhat unsatisfying because it does not yield a formula which {\em describes} the colimit
of a diagram $K \rightarrow \Alg_{\calO}(\calC)$ except in the special case where $K$ Is sifted.
However, we can sometimes do better for specific choices of the underlying $\infty$-operad
$\calO^{\otimes}$. For example, if $\calO^{\otimes}$ is the commutative $\infty$-operad
$\CommOp = \Nerve(\FinSeg)$, then it is easy to construct finite coproducts in the
$\infty$-category $\Alg_{\calO}(\calC) = \CAlg(\calC)$: these are simply given by tensor products
in the underlying $\infty$-category $\calC$ (Proposition \ref{cocarten}). For coproducts of
{\em empty} collections, we can be more general: there is an easy construction for initial objects
of $\Alg_{\calO}(\calC)$ whenever the $\infty$-operad $\calO^{\otimes}$ is unital. We will
describe this construction in \S \ref{unitr}.

\subsection{Limits of Algebras}\label{limco1}

In this section, we will study the formation of limits in $\infty$-categories of algebra
objects. Our main result is the following:

\begin{proposition}\label{limycomps2}
Suppose we are given maps of $\infty$-operads
$\calC^{\otimes} \stackrel{p}{\rightarrow} \calD^{\otimes} \rightarrow \calO^{\otimes}$.
and a commutative diagram
$$ \xymatrix{ K \ar[r]^{f} \ar[d] & \Alg_{\calO}(\calC) \ar[d]^{q} \\
K^{\triangleleft} \ar[r]^{g} & \Alg_{\calO}(\calD). }$$
Assume that for every object $X \in \calO$, the induced diagram
$$ \xymatrix{ K \ar[r]^{f_X} \ar[d] & \calC_{X} \ar[d] \\
K^{\triangleleft} \ar[r] \ar@{-->}[ur]^{ \overline{f}_X} & \calD_{X}  }$$
admits an extension as indicated, where $\overline{f}_X$ is a 
$p$-limit diagram. Then:

\begin{itemize}
\item[$(1)$] There exists an extension $\overline{f}: K^{\triangleleft} \rightarrow \Alg_{\calO}(\calC)$
of $f$ which is compatible with $g$, such that $\overline{f}$ is a $q$-limit diagram.

\item[$(2)$] Let $\overline{f}: K^{\triangleleft} \rightarrow \Alg_{\calO}(\calC)$
be an arbitrary extension of $f$ which is compatible with $g$. Then $\overline{f}$ is a $q$-limit
diagram if and only if for every object $X \in \calO$, the induced map $\overline{f}_{X}: K^{\triangleleft} \rightarrow \calC_{X}$ is a $p$-limit diagram.
\end{itemize}
\end{proposition}

\begin{warning}\label{keps}
Let $\overline{f}: K^{\triangleleft} \rightarrow \Alg_{\calO}(\calC)$ be as in part $(2)$ of
Proposition \ref{limycomps2}, and let $X \in \calO$. Although the map
$\overline{f}_{X}$ takes values in $\calC \subseteq \calC^{\otimes}$, the condition that
$\overline{f}_{X}$ be a $p$-limit diagram is generally {\em stronger} than the condition that
$\overline{f}_{X}$ be a $p_0$-limit diagram, where $p_0: \calC \rightarrow \calD$ denotes
the restriction of $p$.
\end{warning}

In spite of Warning \ref{keps}, the criterion of Proposition \ref{limycomps2} can be simplified
if we are willing to restrict our attention to the setting of $\calO$-monoidal $\infty$-categories.

\begin{corollary}
Let $\calO^{\otimes}$ be an $\infty$-operad and let
$p: \calC^{\otimes} \rightarrow \calD^{\otimes}$ be a $\calO$-monoidal functor
between $\calO$-monoidal $\infty$-categories. 
Suppose given a commutative diagram
$$ \xymatrix{ K \ar[r]^{f} \ar[d] & \Alg_{\calO}(\calC) \ar[d]^{q} \\
K^{\triangleleft} \ar[r]^{g} & \Alg_{\calO}(\calD) }$$
such that, for every object $X \in \calO$, the induced diagram
$$ \xymatrix{ K \ar[r]^{f_X} \ar[d] & \calC_{X} \ar[d]^{p_X} \\
K^{\triangleleft} \ar[r] \ar@{-->}[ur]^{ \overline{f}_X} & \calD_{X}  }$$
admits an extension as indicated, where $\overline{f}_X$ is a 
$p_X$-limit diagram. Then:
\begin{itemize}
\item[$(1)$] There exists an extension $\overline{f}: K^{\triangleright} \rightarrow \Alg_{\calO}(\calC)$
of $f$ which is compatible with $g$, such that $\overline{f}$ is a $q$-limit diagram.

\item[$(2)$] Let $\overline{f}: K^{\triangleleft} \rightarrow \Alg_{\calO}(\calC)$
be an arbitrary extension of $f$ which is compatible with $g$. Then $\overline{f}$ is a $q$-limit
diagram if and only if for every object $X \in \calO$, the induced map $K^{\triangleleft} \rightarrow \calC$ is a $p_X$-limit diagram.
\end{itemize}
\end{corollary}

\begin{proof}
Combine Proposition \ref{limycomps2} with Corollary \toposref{pannaheave}.
\end{proof}

Passing to the case $\calD^{\otimes} = \calO^{\otimes}$, we obtain the following result:

\begin{corollary}
Let $p: \calC^{\otimes} \rightarrow \calO^{\otimes}$ be a coCartesian fibration of 
$\infty$-operads, and let $q: K \rightarrow \Alg_{\calO}(\calC)$ be a diagram.
Suppose that, for every object $X \in \calO$, the induced diagram $q_{X}: K \rightarrow \calC_{X}$
admits a limit. Then:
\begin{itemize}
\item[$(1)$] The diagram $q: K \rightarrow \Alg_{\calO}(\calC)$ admits a limit.
\item[$(2)$] An extension $\overline{q}: K^{\triangleleft} \rightarrow \Alg_{\calO}(\calC)$ of
$q$ is a limit diagram if and only if the induced map $\overline{q}_{X}: K^{\triangleleft} \rightarrow \calC_{X}$ is a limit diagram for each $X \in \calO$.
\end{itemize}
\end{corollary}

\begin{corollary}\label{slimycamp2}
Let $p: \calC^{\otimes} \rightarrow \calO^{\otimes}$ be a coCartesian fibration of $\infty$-operads and let $K$ be a simplicial set. Suppose that, for every object $X \in \calO^{\otimes}$, the fiber
$\calC_{X}$ admits $K$-indexed limits. Then:
\begin{itemize}
\item[$(1)$] The $\infty$-category $\Alg_{\calO}(\calC)$ admits $K$-indexed limits.
\item[$(2)$] An arbitrary diagram $K^{\triangleleft} \rightarrow \Alg_{\calO}(\calC)$ is a limit
diagram if and only the composite diagram
$$ K^{\triangleleft} \rightarrow \Alg_{\calO}(\calC) \rightarrow \calC_{X}$$
is a limit, for each $X \in \calO$.
\end{itemize}
In particular, for each $X \in \calO$ the evaluation functor $\Alg_{\calO}(\calC) \rightarrow \calC_{X}$ preserves $K$-indexed limits.
\end{corollary}

\begin{corollary}\label{jumunj22}
Let $p \calC^{\otimes} \rightarrow \calO^{\otimes}$ be a coCartesian fibration of
$\infty$-operads. Then a morphism $f$ in $\Alg_{\calO}(\calC)$ is an equivalence
if and only if its image in $\calC_{X}$ is an equivalence, for each $X \in \calO$.
\end{corollary}

\begin{proof}
Apply Corollary \ref{slimycamp2} in the case $K = \Delta^0$.
\end{proof}

We now turn to the proof of Proposition \ref{limycomps2}. We will need the the following simple observation:

\begin{lemma}\label{surinore2}
Let $p: \calC^{\otimes} \rightarrow \calO^{\otimes}$ be a fibration of $\infty$-operads, and let
$\gamma: A \rightarrow A'$ be a morphism in $\Alg_{\calO}(\calC)$. The following conditions are equivalent:
\begin{itemize}
\item[$(1)$] The morphism $\gamma$ is an equivalence in $\Alg_{\calO}(\calC)$.
\item[$(2)$] For every object $X \in \calO$, the morphism $\gamma(X)$ is an equivalence in $\calC$.
\end{itemize}
\end{lemma}

\begin{proof}
The implication $(1) \Rightarrow (2)$ is obvious. Conversely, suppose that $(2)$ is satisfied. Let
$X \in \calO^{\otimes}_{\seg{n}}$; we wish to prove that $\gamma(X)$ is an equivalence in
$\calC^{\otimes}_{\seg{n}}$. Since $\calC^{\otimes}$ is a symmetric monoidal $\infty$-category, it suffices to show that for every $1 \leq j \leq n$, the image of $\gamma(X)$ under the functor
$\Colll{j}{!} : \calC^{\otimes}_{\seg{n}} \rightarrow \calC$ is an equivalence.
Since $A$ and $A'$ are lax monoidal functors, this morphism can be identified with $\gamma(X_j)$ where $X_j$ is the image of $X$ under the corresponding functor $\calO^{\otimes}_{\seg{n}} \rightarrow \calO$. The desired result now follows immediately from $(2)$.
\end{proof}

\begin{proof}[Proof of Proposition \ref{limycomps2}]
We first establish the following:
\begin{itemize}
\item[$(\ast)$] Let $h: K^{\triangleleft} \rightarrow \calC^{\otimes}_{\seg{n}}$ be a diagram, and
for $1 \leq i \leq n$ let $h_i$ denote the image of $h$ under the functor
$\Colll{i}{!} : \calC^{\otimes}_{\seg{n}} \rightarrow \calC$. If each
$h_i$ is a $p$-limit diagram, then $h$ is a $p$-limit diagram.
\end{itemize}
To prove $(\ast)$, we observe that there are natural transformations
$h \rightarrow h_i$ which together determine a map $H: K^{\triangleleft} \times \nostar{n}^{\triangleleft} \rightarrow \calC^{\otimes}$. Since the restriction of $H$ to $K^{\triangleleft} \times \{i\}$ is a $p$-limit
diagram for $1 \leq i \leq n$ and the restriction of $H$ to
$\{v\} \times \nostar{n}^{\triangleleft}$ is a $p$-limit diagram for all vertices
$v$ in $K^{\triangleleft}$ (Remark \ref{cotallee}), we deduce from Lemma \toposref{limitscommute} that
$h$ is a $p$-limit diagram as desired.

We now prove the ``if'' direction of $(2)$. Suppose that $\overline{f}: K^{\triangleleft} \rightarrow
\Alg_{\calO}(\calC)$ induces a $p$-limit diagram $\overline{f}_{X}: K^{\triangleleft} \rightarrow
\calC^{\otimes}$ for each $X \in \calO$. We claim that the same assertion holds for each
$X \in \calO^{\otimes}$. To prove this, let $\seg{n}$ denote the image of $X$ in
$\Nerve(\FinSeg)$, and choose inert morphisms $X \rightarrow X_i$ in $\calO^{\otimes}$
lifting $\Colp{i}: \seg{n} \rightarrow \seg{1}$ for $1 \leq i \leq n$.
The image of $\overline{f}_{X}$ under $\Colll{i}{!} $ can be identified
with $\overline{f}_{X_i}$, and is therefore a $p$-limit diagram; the desired result now follows from $(\ast)$. Applying Lemma \monoidref{surtybove}, we deduce that the diagram $\overline{f}$ is a $q'$-limit diagram, where $q'$ denotes the projection
$\bHom_{ \calO^{\otimes}}( \calO^{\otimes}, \calC^{\otimes}) \rightarrow \bHom_{ \calO^{\otimes}}( \calO^{\otimes}, \calD^{\otimes})$. Passing to the full subcategories
$$\Alg_{\calO}(\calC) \subseteq \bHom_{ \calO^{\otimes}}( \calO^{\otimes}, \calC^{\otimes}) \quad \quad 
\Alg_{\calO}(\calD) \subseteq \bHom_{ \calO^{\otimes}}( \calO^{\otimes}, \calD^{\otimes})$$
we deduce that $\overline{f}$ is also a $q$-limit diagram, as desired.

We now prove $(1)$. Choose a categorical equivalence $i: K \rightarrow K'$, where
$K'$ is an $\infty$-category and $i$ is a monomorphism. Since $\Alg_{\calO}(\calD)$ is
an $\infty$-category and $q$ is a categorical fibration, we can assume that $f$ and $g$
factor through compatible maps $f': K' \rightarrow \Alg_{\calO}(\calC)$ and
${K'}^{\triangleleft} \rightarrow \Alg_{\calO}(\calD)$. Using Proposition \toposref{princex}, we deduce that
the extension $\overline{f}$ exists if and only if there is an analogous extension of $f'$. We are therefore free to replace $K$ by $K'$ and reduce to the case where $K$ is an $\infty$-category.
The map $f$ classifies a functor $F: K \times \calO^{\otimes} \rightarrow \calC^{\otimes}$
and the map $g$ classifies a functor $G: K^{\triangleleft} \times \calO^{\otimes} \rightarrow \calC^{\otimes}$. We first claim:
\begin{itemize}
\item[$(\ast')$] There exists an extension $\overline{F}: K^{\triangleleft} \times \calO^{\otimes} \rightarrow \calC^{\otimes}$ of $F$ lying over $G$ such that $\overline{F}$ is a $p$-right Kan extension of $F$.
\end{itemize}
Let $v$ denote the cone point of $K^{\triangleleft}$.
According to Lemma \toposref{kan2}, it suffices to show that for each object $X \in \calO^{\otimes}$, the induced diagram
$(K \times \calO^{\otimes})_{(v, X)/} \simeq K \times \calO^{\otimes}_{X/}
\rightarrow \calC^{\otimes}$ can be extended to a $p$-limit diagram compatible with $G$.
Since the inclusion $K \times \{ \id_{X} \} \hookrightarrow K \times \calO^{\otimes}_{X/}$ is
the opposite of a cofinal morphism, it suffices to show that we can complete the diagram
$$ \xymatrix{ K \ar[r]^{f_X} \ar[d] & \calC^{\otimes} \ar[d] \\
K^{\triangleleft} \ar[r]^{g_{X}} \ar@{-->}[ur]^{ \overline{f}_X} & \calD^{\otimes} }$$
so that $\overline{f}_{X}$ is a $p$-limit diagram. Let $\seg{n}$ denote the image
of $X$ in $\Nerve(\FinSeg)$. Since $\calC^{\otimes}$ and $\calD^{\otimes}$ are $\infty$-operads, 
we have a homotopy commutative diagram
$$ \xymatrix{ \calC^{\otimes}_{\seg{n}} \ar[r] \ar[d] & \calC^{n} \ar[d] \\
\calD^{\otimes}_{\seg{n}} \ar[r] & \calD^{n} }$$
where the horizontal morphisms are categorical equivalences. Using Proposition \toposref{princex}
and our assumption on $f$, we deduce that $f_{X}$ admits an extension $\overline{f}_{X}$ (compatible with $g_X$) whose
images under the functors $\Colll{i}{!} $ are $p$-limit diagrams in $\calC$.
It follows from $(\ast)$ that $\overline{f}_{X}$ is a $p$-limit diagram, which completes the verification
of $(\ast')$. Moreover, the proof shows that an extension $\overline{F}$ of $F$ (compatible with $G$)
is a $p$-right Kan extension if and only if, for each $X \in \calO^{\otimes}_{\seg{n}}$,
the functor $\Colll{i}{!} : \calC^{\otimes}_{\seg{n}} \rightarrow \calC$ carries
$\overline{F} | (K^{\triangleleft} \times \{X\})$ to a $p$-limit diagram in $\calC^{\otimes}$.

Let $s: \calO^{\otimes} \rightarrow \calC^{\otimes}$ be the functor obtained by restricting
$\overline{F}$ to the cone point of $K^{\triangleleft}$. Then $s$ preserves inert morphisms
lying over the maps $\Colp{i}: \seg{n} \rightarrow \seg{1}$ and is therefore a map of
$\infty$-operads (Remark \ref{casper}). It follows that $\overline{F}$ determines an extension
$\overline{f}: K^{\triangleleft} \rightarrow \Alg_{\calO}(\calC)$ of $f$ lifting $g$.
The first part of the proof shows that $\overline{f}$ is a $q$-limit diagram, which completes the proof of $(1)$.

We conclude by proving the ``only if'' direction of $(2)$. Let $\overline{f}$ be a $q$-left Kan extension of $f$ lying over $g$. The proof of $(1)$ shows that
we can choose another $q$-left Kan extension $\overline{f}'$ of $f$ lying over $g$ such that
$\overline{f}'_{X}: K^{\triangleleft} \rightarrow \calC^{\otimes}$ is a $p$-limit diagram for
each $X \in \calO$. It follows that $\overline{f}$ and $\overline{f}'$ are equivalent, so that
each $\overline{f}_{X}$ is also a $p$-limit diagram.
\end{proof}

\subsection{Operadic Colimit Diagrams}\label{comm35}

Let $\calC^{\otimes}$ be a symmetric monoidal $\infty$-category. 
If the underlying $\infty$-category $\calC$ of $\calC^{\otimes}$ admits limits, then
the $\infty$-category of commutative algebra objects $\CAlg(\calC)$ admits limits:
moreover, these limits can be computed at the level of the underlying objects of
$\calC$ (Corollary \ref{slimycamp2}). The analogous statement for colimits is false:
colimits of commutative algebra objects generally cannot be computed
at the level of the underlying $\infty$-category $\calC$ (though this is often
true for sifted colimits: see Proposition \ref{fillfemme}). Moreover, to guarantee
the existence of colimits in $\CAlg(\calC)$ it is not enough to assume that $\calC$
admits all (small) colimits: we must also assume that the formation of colimits
in $\calC$ is compatible with the symmetric monoidal structure on $\calC$.

\begin{itemize}
\item[$(\ast)$] For every object $C \in \calC$, the functor $\bigdot \mapsto C \otimes \bigdot$
preserves colimits.
\end{itemize}

Our goal in this section is to study a generalization of condition $(\ast)$, which is adapted
to fibrations of $\infty$-operads $q: \calC^{\otimes} \rightarrow \calO^{\otimes}$ which
are not necessarily coCartesian fibrations, where $\calO^{\otimes}$ is not necessarily
the commutative $\infty$-operad, and for which the fibers $\calC_{X}$ are not necessarily
assumed to admit colimits in general. To this end, we will introduce the rather technical
notion of an {\it operadic $q$-colimit diagram} (Definition \ref{swm}). Though this notion is somewhat complicated, it often reduces (under the assumption of a condition like $(\ast)$) to the usual theory
of colimits: see Example \ref{coltens}.

\begin{notation}
Let $\calO^{\otimes}$ be an $\infty$-operad. We let $\calO^{\otimes}_{\acti}$ denote the subcategory
of $\calO^{\otimes}$ spanned by the active morphisms.
\end{notation}

\begin{definition}\label{swm}
Let $\calO^{\otimes}$ be an $\infty$-operad and let $p: K \rightarrow \calO^{\otimes}$ be a diagram.
We let $\calO_{p/}$ denote the $\infty$-category $\calO \times_{ \calO^{\otimes}} \calO^{\otimes}_{p/}$.
If $p$ factors through $\calO^{\otimes}_{\acti}$, we let
$\calO^{\acti}_{p/} \subseteq \calO_{p/}$ denote the $\infty$-category
$\calO \times_{ \calO^{\otimes}_{\acti} } (\calO^{\otimes}_{\acti})_{p/}$.

Let $q: \calC^{\otimes} \rightarrow \calO^{\otimes}$ be a fibration of $\infty$-operads,
let $\overline{p}: K^{\triangleright} \rightarrow \calC^{\otimes}_{\acti}$ be a diagram, and
let $p = \overline{p} | K$. We will say that $\overline{p}$ is a {\it weak operadic $q$-colimit diagram} if the evident map
$$ \psi: \calC^{\acti}_{ \overline{p}/ } \rightarrow
\calC^{\acti}_{ p/} \times_{ \calO^{\acti}_{qp/} } \calO^{\acti}_{q \overline{p}/}$$
is a trivial Kan fibration.

We say that an active diagram $p: K \rightarrow \calC^{\otimes}$ is an
{\it operadic $q$-colimit diagram} if the composite functor
$$ K \stackrel{p}{\rightarrow} \calC^{\otimes}_{\acti} \stackrel{\bigdot \oplus C}{\longrightarrow}
\calC^{\otimes}_{\acti}$$
is a weak operadic $q$-colimit diagram, for every object $C \in \calC^{\otimes}$
(here the functor $\oplus$ is defined as in Remark \ref{opl}).
\end{definition}

\begin{remark}\label{woperad}
Let $q: \calC^{\otimes} \rightarrow \calO^{\otimes}$ be as in Definition \ref{swm}, let
$\overline{p}: K^{\triangleright} \rightarrow \calC^{\otimes}_{\acti}$ be a weak
operadic $q$-colimit diagram, and let $L \rightarrow K$ be a cofinal map of simplicial sets.
Then the induced map $L^{\triangleright} \rightarrow \calC^{\otimes}_{\acti}$ is a
weak operadic $q$-colimit diagram.
\end{remark}

\begin{remark}\label{coheath}
The map $\psi$ appearing in Definition \ref{swm} is always a left fibration
(Proposition \toposref{sharpen}); consequently, it is a trivial Kan fibration if and only if
it is a categorical equivalence (Corollary \toposref{heath}).
\end{remark}

\begin{example}\label{jim}
Let $q: \calC^{\otimes} \rightarrow \calO^{\otimes}$ be as in Definition \ref{swm}, let
$K = \Delta^0$, and suppose that $\overline{p}: K^{\triangleright} \simeq \Delta^1 \rightarrow \calC^{\otimes}_{\acti}$ corresponds to an equivalence in $\calC^{\otimes}$. Then
$\overline{p}$ is an operadic $q$-colimit diagram.
\end{example}

The remainder of this section is devoted to a series of somewhat technical results which are useful for working with operadic colimit diagrams.

\begin{proposition}\label{quadly}
Let $q: \calC^{\otimes} \rightarrow \calO^{\otimes}$ be a fibration of $\infty$-operads
let $\overline{p}: K^{\triangleright} \rightarrow \calC^{\otimes}_{\acti}$ be a diagram.
The following conditions are equivalent:
\begin{itemize}
\item[$(1)$] The map $\overline{p}$ is a weak operadic $q$-colimit diagram.
\item[$(2)$] For every $n > 0$ and every diagram
$$ \xymatrix{ K \star \bd \Delta^n \ar[r]^{f_0} \ar[d] & \calC^{\otimes}_{\acti} \ar[d] \\
K \star \Delta^n \ar@{-->}[ur]^{f} \ar[r]^{\overline{f}} & \calO^{\otimes}_{\acti} }$$
such that the restriction of $f_0$ to $K \star \{0\}$ coincides with $\overline{p}$
and $f_0(n) \in \calC$, there exists a dotted arrow $f$ as indicated above.
\end{itemize}
\end{proposition}

\begin{proof}
The implication $(2) \Rightarrow (1)$ follows immediately from the definition.
For the converse, suppose that $(1)$ is satisfied. Let
$r: \Delta^n \times \Delta^1 \rightarrow \Delta^n$ be the map which is
given on vertices by the formula
$$ r(i,j) = \begin{cases} i & \text{if } j = 0 \\
n & \text{if } j = 1. \end{cases}$$
Set $$X = \bd \Delta^n \coprod_{ \Lambda^{n}_{n} \times \{0\} }
(\Lambda^n_{n} \times \Delta^1) \coprod_{ \Lambda^n_{n} \times \{1\} } \Delta^n,$$
and regard $X$ as a simplicial subset of $\Delta^n \times \Delta^1$; we observe that
$r$ carries $X$ into $\bd \Delta^n$. Composing $f_0$ and $\overline{f}$ with
$r$, we obtain a diagram
$$ \xymatrix{ X \ar[r]^{g_0} \ar[d] & (\calC^{\otimes}_{\acti})_{f_0|K /} \ar[d]^{q'} \\
\Delta^n \times \Delta^1) \ar[r] \ar@{-->}[ur]^{g} & (\calO^{\otimes}_{\acti})_{\overline{f}|K/}. }$$
To complete the proof, it will suffice to show that we can supply the dotted arrow in this
diagram (we can then obtain $f$ by restricting $g$ to $\Delta^n \times \{0\}$).  
To prove this, we define a sequence of simplicial subsets
$$ X = X(0) \subset X(1) \subset X(2) \subset \cdots \subset X(2n+1) = \Delta^n \times \Delta^1$$
and extend $g_0$ to a sequence of maps $g_i: X(i) \rightarrow (\calC^{\otimes}_{\acti})_{f_0|K/}$
compatible with the projection $q'$. The analysis proceeds as follows:

\begin{itemize}
\item[$(i)$] For $1 \leq i \leq n$, we let $X(i)$ be the simplicial subset
of $\Delta^n \times \Delta^1$ generated by $X(i-1)$ and the
$n$-simplex $\sigma: \Delta^n \rightarrow \Delta^n \times \Delta^1$
given on vertices by the formula
$$ \sigma(j) = \begin{cases} (j, 0) & \text{if } j \leq n-i \\
(j-1, 1) & \text{if } j > n-i. \end{cases}$$
We have a pushout diagram of simplicial sets
$$ \xymatrix{ \Lambda^{n}_{n-i} \ar[r] \ar[d] & X(i-1) \ar[d] \\
\Delta^n \ar[r] & X(i). }$$
If $1 \leq i < n$, then the extension $g_i$ of $g_{i-1}$ exists because $q'$ is an inner fibration. If $i = n$, then the desired extension exists by virtue of assumption $(1)$.

\item[$(ii)$] If $i = n+ i'$ for $0 < i' \leq n+1$, we let $X(i)$ be the simplicial subset
of $\Delta^n \times \Delta^1$ generated by $X(i-1)$ and the $(n+1)$-simplex
$\sigma: \Delta^{n+1} \rightarrow \Delta^{n} \times \Delta^1$ given on vertices by the formula
$$ \sigma(j) = \begin{cases} (j, 0) & \text{if } j < i' \\
(j-1, 1) & \text{if } j \geq i'. \end{cases}$$
We have a pushout diagram of simplicial sets
$$ \xymatrix{ \Lambda^{n+1}_{i'} \ar[r] \ar[d] & X(i-1) \ar[d] \\
\Delta^{n+1} \ar[r] & X(i). }$$
If $1 \leq i' \leq n$, then the extension $g_i$ of $g_{i-1}$ exists because $q'$ is an inner fibration. 
If $i = n+1$, then the desired extension exists because $g_{i-1}$ carries
$\Delta^{ \{n, n+1\}}$ to an equivalence in $(\calC^{\otimes}_{\acti})_{f_0|K/}$.
\end{itemize}
\end{proof}

\begin{proposition}\label{swingup}
Let $q: \calC^{\otimes} \rightarrow \calO^{\otimes}$ be a fibration of $\infty$-operads, and suppose we are given a finite collection 
$\{ \overline{p}_{i}: K^{\triangleright}_{i} \rightarrow \calC^{\otimes}_{\acti} \}_{i \in I}$ of
operadic $q$-colimit diagrams. Let $K = \prod_{i \in I} K_i$, and let
$\overline{p}$ denote the composition
$$ K^{\triangleright} \rightarrow
\prod_{i \in I} K_{i}^{\triangleright} \rightarrow \prod_{i \in I} \calC^{\otimes}_{\acti}
\stackrel{ \oplus_{I} }{\rightarrow} \calC^{\otimes}_{\acti}.$$
Then $\overline{p}$ is an operadic $q$-colimit diagram.
\end{proposition}

The proof of Proposition \ref{swingup} depends on the following lemma:

\begin{lemma}\label{fatwise}
Let $q: \calC^{\times} \rightarrow \calO^{\otimes}$ be a fibration of $\infty$-operads, let $p': K_0 \times K_1^{\triangleright} \rightarrow \calC^{\otimes}_{\acti}$, and let
$p = p' | K_0 \times K_1$. Suppose that, for each vertex $v$ of $K_0$, the induced map
$p'_{v}: K_1^{\triangleright} \rightarrow \calC^{\otimes}_{\acti}$ is a weak operadic
$q$-colimit diagram. Then the map
$$ \theta: \calC^{\acti}_{p'/} \rightarrow \calC^{\acti}_{p/} \times_{ \calO^{\acti}_{qp/}}
\calO^{\acti}_{qp'/}$$
is a trivial Kan fibration.
\end{lemma}

\begin{proof}
For each simplicial subset $A \subseteq K_0$, let $L_{A}$ denote the pushout
$(K_0 \times K_1) \coprod_{ (A \times K_1) } (A \times K_1^{\triangleright} )$, and let
$p_{A} = p' | L_A$. For $A' \subseteq A$, let $\theta^{A}_{A'}$ denote the map
$$ \calC^{\acti}_{p_A/} \rightarrow \calC^{\acti}_{p_{A'}/} \times_{ \calO^{\acti}_{qp_{A'}/}}
\calO^{\acti}_{qp_{A}/}.$$
We will prove that each of the maps $\theta^{A}_{A'}$ is a trivial Kan fibration.
Taking $A = K_0$ and $A' = \emptyset$, this will imply the desired result.

Working simplex-by-simplex on $A$, we may assume without loss of generality that
$A$ is obtained from $A'$ by adjoining a single nondegenerate $n$-simplex whose boundary already belongs to $A'$. Replacing $K_0$ by $A$, we may assume that $K_0 = A = \Delta^n$ and
$A' = \bd \Delta^n$. Working by induction on $n$, we may assume that the map
$\theta^{A'}_{\emptyset}$ is a trivial Kan fibration. Since the map
$\theta^{A}_{A'}$ is a left fibration (Proposition \toposref{sharpen}), it is
a trivial Kan fibration if and only if it is a categorical fibration (Corollary \toposref{heath}).
Since $\theta^{A}_{\emptyset} = \theta^{A'}_{\emptyset} \circ \theta^{A}_{A'}$, we can use
the two-out-of-three property to reduce to the problem of showing that
$\theta^{A}_{\emptyset}$ is a categorical equivalence. Since $A = \Delta^n$, the inclusion
$\{ n\} \subseteq A$ is cofinal so that $\theta^{A}_{ \{n\} }$ is a trivial Kan fibration.
It therefore suffices to show that $\theta^{ \{n\} }_{\emptyset}$ is a trivial Kan fibration, which follows
from our assumption that each $p'_{v}$ is a weak operadic $q$-colimit diagram.
\end{proof}

\begin{proof}[Proof of Proposition \ref{swingup}]
If $I$ is empty, the desired result follows from Example \ref{jim}. If $I$ is a singleton, the result is obvious. To handle the general case, we can use induction on the cardinality of $I$ to reduce to the case
where $I$ consists of two elements. Let us therefore suppose that we are given
operadic $q$-colimit diagrams $\overline{p}_0: K_0^{\triangleright} \rightarrow \calC^{\otimes}$
and $\overline{p}_1: K_1^{\triangleright} \rightarrow \calC^{\otimes}$; we wish to prove that the induced map $\overline{p}: (K_0 \times K_1)^{\triangleright} \rightarrow \calC^{\otimes}$ is an operadic
$q$-colimit diagram. Let $X \in \calC^{\otimes}$; we must show that the composition
$$ (K_0 \times K_1)^{\triangleright} \rightarrow \calC^{\otimes} \stackrel{ \oplus X}{\rightarrow} \calC^{\otimes}$$
is a weak operadic $q$-colimit diagram. Replacing $\overline{p}_0$ by the composition
$$ K_0^{\triangleright} \rightarrow \calC^{\otimes} \stackrel{ \oplus X}{\rightarrow} \calC^{\otimes},$$
we can reduce to the case where $X$ is trivial; it will therefore suffice to prove that
$\overline{p}$ is a weak operadic $q$-colimit diagram.

Let $L$ denote the simplicial set $K_0 \times K_1^{\triangleright}$. We have a commutative diagram of simplicial sets
$$ \xymatrix{ K \ar[r] \ar[d] & L \ar[d] & K_0 \times \{v\} \ar[l] \ar[d] \\
K^{\triangleright} \ar[r] & L^{\triangleright} & K^{\triangleright}_0 \times \{v\} \ar[l] }$$
where $v$ denotes the cone point of $K_1^{\triangleright}$. 
Let $\overline{p}'$ denote the composition
$$ L^{\triangleright} \rightarrow K_0^{\triangleright} \times K_1^{\triangleright} 
\rightarrow \calC^{\otimes}_{\acti} \times \calC^{\otimes}_{\acti} \stackrel{\oplus}{\rightarrow} \calC^{\otimes}_{\acti}.$$
Since $\overline{p}_0$ is an operadic $q$-colimit diagram, we deduce that
$\overline{p}' | (K_0^{\triangleright} \times \{v\} )$ is a weak operadic $q$-colimit diagram.
Since the inclusion $K_0 \times \{v\} \subseteq L$ is cofinal, we deduce that
$\overline{p}'$ is a weak operadic $q$-colimit diagram (Remark \ref{woperad}).
The inclusion $K^{\triangleright} \subseteq L^{\triangleright}$ is cofinal
(since the cone point is cofinal in both). Consequently, to prove that $\overline{p}$ is
a weak operadic $q$-colimit diagram, it will suffice to show that the map
$$ \calC^{\acti}_{\overline{p}'/} \rightarrow \calC^{\acti}_{p/} \times_{ \calO^{\acti}_{qp/}}
\calO^{\acti}_{ q\overline{p}'/},$$
where $p = \overline{p} | K$. Since $\overline{p}'$ is a weak operadic $q$-colimit diagram,
it will suffice to show that the map
$$ \calC^{\acti}_{p'/} \times_{ \calO^{\acti}_{qp'/} }
\calO^{\acti}_{ q\overline{p}'/} \rightarrow \calC^{\acti}_{p/}
\times_{ \calO^{\acti}_{qp/} } \calO^{\acti}_{ q\overline{p}'/ }$$
is a trivial Kan fibration, where $p' = \overline{p}' | L$. This map is a pullback of
$$ \calC^{\acti}_{p'/} \rightarrow \calC^{\acti}_{p/} \times_{ \calO^{\acti}_{qp/} }
\calO^{\acti}_{qp'/},$$
which is a trivial Kan fibration by Lemma \ref{fatwise} (since $\overline{p}_1$ 
is an operadic $q$-colimit diagram).
\end{proof}

\begin{proposition}\label{copan}
Let $q: \calC^{\otimes} \rightarrow \calO^{\otimes}$ be a fibration of $\infty$-operads. Suppose we are given a finite collection of operadic $q$-colimit diagrams
$\{ \overline{p}_{i}: K_{i}^{\triangleright} \rightarrow \calC^{\otimes}_{\acti} \}_{i \in I}$, where each $\overline{p}_i$ carries
the cone point of $K_i^{\triangleright}$ into $\calC \subseteq \calC^{\otimes}$. 
Let $K = \prod_{i \in I} K_i$, and let $\overline{p}: K^{\triangleright} \rightarrow \calC^{\otimes}$ be defined
as in Proposition \ref{swingup}. If each of the simplicial sets $K_i$ is weakly contractible, then
$\overline{p}$ is a $q$-colimit diagram.
\end{proposition}

We will need a few preliminaries:

\begin{lemma}\label{sidewell}
Let $q: \calC^{\otimes} \rightarrow \calO^{\otimes}$ be a fibration of $\infty$-operads, let $X \in \calC^{\otimes}$ be an object lying over $\seg{n} \in \FinSeg$, and choose inert morphisms $f_i: X \rightarrow X_i$ in $\calC^{\otimes}$ covering $\Colp{i}: \seg{n} \rightarrow \seg{1}$ for $1 \leq i \leq n$. Then the maps $\{ f_i: X \rightarrow X_i \}_{1 \leq i \leq n}$
determine a $q$-limit diagram $\chi: \nostar{n}^{\triangleleft} \rightarrow \calC^{\otimes}$.
\end{lemma}

\begin{proof}
Let $p: \calO^{\otimes} \rightarrow \Nerve(\FinSeg)$ be the map which exhibits
$\calO^{\otimes}$ as an $\infty$-operad. Since each $q(f_i)$ is an inert morphism in
$\calO^{\otimes}$, we can invoke the definition of an $\infty$-operad to deduce that
$q \circ \chi$ is a $p$-limit diagram in $\calO^{\otimes}$. Similarly, we conclude that
$\chi$ is a $p \circ q$-limit diagram in $\calC^{\otimes}$. The desired result now follows from 
Proposition \toposref{basrel}.
\end{proof}

\begin{lemma}\label{pin}
Let $q: \calC^{\otimes} \rightarrow \calO^{\otimes}$ be a fibration of $\infty$-operads, let $\overline{p}: K^{\triangleright} \rightarrow \calC^{\otimes}$ be a diagram, and let
$p = \overline{p} | K$. Then $\overline{p}$ is a $q$-colimit diagram if and only if
the induced map
$$\phi_0: \calC_{ \overline{p}/ } \rightarrow \calC_{p/} \times_{ \calO_{qp/}} \calO_{q \overline{p}/}$$
is a trivial Kan fibration.
\end{lemma}

\begin{proof}
The ``only if'' direction is obvious. To prove the converse, suppose that $\phi_0$ is a trivial Kan fibration, and consider the map
$$\phi: \calC^{\otimes}_{ \overline{p}/ } \rightarrow \calC^{\otimes}_{p/} \times_{ \calO^{\otimes}_{qp/}} \calO^{\otimes}_{q \overline{p}/}.$$
We wish to prove that $\phi$ is a trivial Kan fibration. Since $\phi$ is a left fibration
(Proposition \toposref{sharpen}), it suffices to show that the fibers of $\phi$ are contractible
(Lemma \toposref{toothie}). Let $\overline{X}$ be an object in the codomain of $\phi$ having
image $X \in \calC^{\otimes}$ and image $\seg{n} \in \FinSeg$. For $1 \leq i \leq n$, choose
an inert morphism $X \rightarrow X_{i}$ in $\calC^{\otimes}$ lying over
$\Colp{i}: \seg{n} \rightarrow \seg{1}$. Since the projection
$$\calC^{\otimes}_{p/} \times_{ \calO^{\otimes}_{qp/}} \calO^{\otimes}_{q \overline{p}/} \rightarrow \calC$$ (this follows from repeated application of Proposition \toposref{sharpen}), we can lift each
of these inert morphisms in an essentially unique way to a map $\overline{X} \rightarrow \overline{X}_i$.
Using Lemma \ref{sidewell}, we deduce that the fiber
$\phi^{-1}(X)$ is homotopy equivalent to the product $\prod_{1 \leq i \leq n} \phi^{-1}( X_i )$, and
therefore contractible (since $\phi^{-1}(X_i) = \phi_0^{-1}(X_i)$ is a fiber of the trivial Kan fibration
$\phi_0$ for $1 \leq i \leq n$), as desired.
\end{proof}

\begin{proof}[Proof of Proposition \ref{copan}]
In view of Lemma \ref{pin}, it will suffice to show that the map
$$\phi_0: \calC_{ \overline{p}/ } \rightarrow \calC_{p/} \times_{ \calO_{qp/}} \calO_{q \overline{p}/}$$
is a trivial Kan fibration. Since $\phi_0$ is a left fibration (Proposition \toposref{sharpen}), it will suffice
to show that $\phi_0^{-1} \{ \overline{X} \}$ is contractible for each vertex $\overline{X} \in \calC_{p/} \times_{ \calO_{qp/}} \calO_{q \overline{p}/}$ (Lemma \toposref{toothie}). Such a vertex determines
a map $\alpha: I_{\ast} \rightarrow \seg{1}$; let $I_0 = \alpha^{-1} \{1\} \subseteq I$. Let
$K' = \prod_{i \in I_0} K_i$, let $\overline{p}'$ denote the composite map
$$ {K'}^{\triangleright} \rightarrow
\prod_{i \in I_0} K_i^{\triangleright} \rightarrow \prod_{i \in I_0} \calC^{\otimes}
\stackrel{ \oplus}{\rightarrow} \calC^{\otimes},$$
and let $\overline{p}''$ denote the composition
$K^{\triangleright} \rightarrow {K'}^{\triangleright} \stackrel{ \overline{p}'}{\rightarrow} \calC^{\otimes}$.

There is a canonical map $\overline{P}: \Delta^1 \times K^{\triangleright} \rightarrow \calC^{\otimes}$
which is an inert natural transformation from $\overline{p}$ to $\overline{p}''$. Let
$p'' = \overline{p}'' | K$ and $P = \overline{P} | \Delta^1 \times K$, so we have a commutative diagram
$$ \xymatrix{ \calC_{\overline{p}/ } \ar[r]^{\phi_0} &  \calC_{p/} \times_{ \calO_{qp/}} \calO_{q \overline{p}/} \\
\calC_{ \overline{P}/} \ar[r] \ar[u] \ar[d] & \calC_{P/} \times_{ \calO_{qP/}} \calO_{q \overline{P}/} \ar[u] \ar[d] \\
\calC_{\overline{p}''/ } \ar[r] &  \calC_{p''/} \times_{ \calO_{qp''/}} \calO_{q \overline{p}''/} }$$
which determines a homotopy equivalence $\phi_0^{-1} \overline{X}$ with a fiber of the map
$\psi: \calC_{ \overline{p}''/}^{\acti} \rightarrow  \calC^{\acti}_{p''/} \times_{ \calO^{\acti}_{qp''/}} \calO^{\acti}_{q \overline{p}''/}.$
It will therefore suffice to show that the fibers of $\psi$ are contractible. Since
$\psi$ is a left fibration (Proposition \toposref{sharpen}), this is equivalent to the assertion
that $\psi$ is a categorical equivalence (Corollary \toposref{heath}).

Because the simplicial sets $\{ K_i \}_{i \in I - I_0}$ are weakly contractible, the projection map
$K \rightarrow K'$ is cofinal. Consequently, we have a commutative diagram
$$ \xymatrix{ \calC^{\acti}_{\overline{p}'/ } \ar[r]^-{\psi'} \ar[d]  &  \calC^{\acti}_{p'/} \times_{ \calO^{\acti}_{qp'/}} \calO^{\acti}_{q \overline{p}'/} \ar[d] \\ 
\calC^{\acti}_{\overline{p}''/ } \ar[r]^-{\psi} &  \calC^{\acti}_{p''/} \times_{ \calO^{\acti}_{qp''/}} \calO^{\acti}_{q \overline{p}''/} }$$
where the vertical maps are categorical equivalences. Using a two-out-of-three argument, we
are reduced to proving that $\psi'$ is a categorical equivalence. We conclude by observing that
Proposition \ref{swingup} guarantees that $\overline{p}'$ is an operadic $q$-colimit diagram,
so that $\psi'$ is a trivial Kan fibration.
\end{proof}

\begin{proposition}\label{justica}
Let $q: \calC^{\otimes} \rightarrow \calO^{\otimes}$ be a fibration of $\infty$-operads and let $K = \Delta^0$. The following conditions are equivalent:
\begin{itemize}
\item[$(1)$] The map $q$ exhibits $\calC^{\otimes}$ as a $\calO$-monoidal $\infty$-category
(in other words, the map $q$ is a coCartesian fibration).
\item[$(2)$] For every map $p: K \rightarrow \calC^{\otimes}_{\acti}$ and every
extension $\overline{p}_0: K^{\triangleright} \rightarrow \calO^{\otimes}_{\acti}$ of
$q \circ p$ carrying the cone point of $K^{\triangleright}$ into $\calO$, there exists
an operadic $q$-colimit diagram $\overline{p}: K^{\triangleright}
\rightarrow \calC^{\otimes}_{\acti}$ which extends $p$ and lifts $\overline{p}_0$.
\end{itemize}
\end{proposition}

\begin{lemma}\label{starr}
Let $\calO^{\otimes}$ be an $\infty$-operad and let $q: \calC^{\otimes} \rightarrow \calO^{\otimes}$ be a $\calO$-monoidal $\infty$-category. Let $\overline{p}: \Delta^1 \rightarrow \calC^{\otimes}_{\acti}$ classify an active $q$-coCartesian morphism $X \rightarrow Y$ in $\calC^{\otimes}$. Then $\overline{p}$ is an operadic $q$-colimit diagram.
\end{lemma}

\begin{proof}
Let $Z \in \calC^{\otimes}$. Replacing $p$ by its composition with the functor
$\calC^{\otimes} \stackrel{ \oplus Z}{\rightarrow} \calC^{\otimes}$ (and using Remark \ref{qact}), we can reduce to showing that $\overline{p}$ is a weak operadic $q$-colimit diagram, which is clear.
\end{proof}

\begin{proof}[Proof of Proposition \ref{justica}]
To prove that $(1)$ implies $(2)$, we observe that we can take $\overline{p}$ to be
a $q$-coCartesian lift of $\overline{p}_0$ (by virtue of Lemma \ref{starr}). Conversely,
suppose that $(2)$ is satisfied. Choose an object $\overline{X} \in \calC^{\otimes}$, 
let $X = q( \overline{X})$, and suppose we are given a morphism
$\alpha: X \rightarrow Z$ in $\calO^{\otimes}$; we wish to prove that we can lift
$\alpha$ to a $q$-coCartesian morphism $\overline{X} \rightarrow \overline{Z}$ in
$\calC^{\otimes}$. Using Proposition \ref{singwell}, we can factor $\alpha$ as a composition
$$ X \stackrel{ \alpha'}{\rightarrow} Y \stackrel{\alpha''}{\rightarrow} Z$$
where $\alpha'$ is inert and $\alpha''$ is active. Let $\alpha'_0$ denote the image of
$\alpha'$ in $\FinSeg$, and choose an inert morphism $\overline{\alpha}': \overline{X} \rightarrow
\overline{Y}$ in $\calC^{\otimes}$ lifting $\alpha'_0$. Then $q( \overline{\alpha}')$ is an inert
lift of $\alpha'_0$, so we may assume without loss of generality that $\overline{\alpha}'$ lifts
$\alpha'$. Proposition \toposref{stuch} guarantees that $\overline{\alpha}'$ is $q$-coCartesian.
Since the collection of $q$-coCartesian morphisms in $\calC^{\otimes}$ is stable under composition
(Proposition \toposref{protohermes}), we may replace $\overline{X}$ by $\overline{Y}$ and
thereby reduce to the case where the map $\alpha$ is active.

Let $\seg{n}$ denote the image of $Z$ in $\FinSeg$. Since $\calC^{\otimes}$ and
$\calC^{\otimes}$ are $\infty$-operads, we can identify
$\overline{X}$ with a concatenation $\oplus_{1 \leq i \leq n} \overline{X}_i$ and
$\alpha$ with $\oplus_{1 \leq i \leq n} \alpha_i$, where
each $\alpha_i: X_i = q( \overline{X}_i) \rightarrow Z_i$ is an active morphism in
$\calO^{\otimes}$ where $Z_i$ lies over $\seg{1} \in \FinSeg$. 
Assumption $(2)$ guarantees that each $\alpha_i$ can be lifted to
an operadic $q$-colimit diagram $\overline{\alpha}_{i}: \overline{X}_i \rightarrow \overline{Z}_i$
in $\calC^{\otimes}$. It follows from Proposition \ref{copan} that the concatenation
$\overline{\alpha} = \oplus_{1 \leq i \leq n} \overline{\alpha}_i$ is $q$-coCartesian.
The map $q( \overline{\alpha})$ is equivalent to $\oplus_{1 \leq i \leq n} q( \overline{\alpha}_i) \simeq \alpha$. Since $q$ is a categorical fibration, we can replace $\overline{\alpha}$ by an
equivalent morphism if necessary to guarantee that $q( \overline{\alpha}) = \alpha$, which completes the proof of $(1)$.
\end{proof}

Our next two results, which are counterparts of Propositions \toposref{chocolatelast} and \toposref{relcolfibtest}, are useful for detecting the existence of operadic
colimit diagrams:

\begin{proposition}\label{chocolateoperad}
Let $q: \calC^{\otimes} \rightarrow \calO^{\otimes}$ be a fibration of $\infty$-operads, let $K$ be a simplicial set, and let 
$\overline{h}: \Delta^1 \times K^{\triangleright} \rightarrow \calC^{\otimes}_{\acti}$ be a natural transformation from
$\overline{h}_0 = \overline{h} | \{0\} \times K^{\triangleright}$ to $\overline{h}_1 = \overline{h} | \{1\} \times K^{\triangleright}$.
Suppose that
\begin{itemize}
\item[$(a)$] For every vertex $x$ of $K^{\triangleright}$, the restriction
$\overline{h} | \Delta^1 \times \{x\}$ is a $q$-coCartesian edge of $\calC^{\otimes}$.
\item[$(b)$] The composition
$$ \Delta^1 \times \{v\} \subseteq \Delta^1 \times K^{\triangleright}
\stackrel{\overline{h}}{\rightarrow} \calC^{\otimes} \stackrel{q}{\rightarrow} \calO^{\otimes}$$
is an equivalence in $\calO^{\otimes}$, where $v$ denotes the cone point of
$K^{\triangleright}$.
\end{itemize}
Then:
\begin{itemize}
\item[$(1)$] 
The map $\overline{h}_0$ is a weak operadic $q$-colimit diagram if and only if $\overline{h}_1$ is a weak operadic $q$-colimit diagram.
\item[$(2)$] Assume that $q$ is a coCartesian fibration. Then $\overline{h}_0$ is an operadic
$q$-colimit diagram if and only if $\overline{h}_1$ is an operadic $q$-colimit diagram.
\end{itemize}
\end{proposition}

\begin{proof}
Assertion $(2)$ follows from $(1)$ and Remark \ref{qact} (after composing with the functor
$\calC^{\otimes} \stackrel{ \oplus X}{\rightarrow} \calC^{\otimes}$ determined by an arbitrary object
$X \in \calC^{\otimes}$). It will therefore suffice to prove $(1)$. 
Let $h = \overline{h} | \Delta^1 \times K$, $h_0 = h | \{0\} \times K$, and $h_1 = h | \{1\} \times K$.
We have a commutative diagram
$$ \xymatrix{ \calC^{\acti}_{\overline{h}_0/} \ar[d] & \calC^{\acti}_{\overline{h}/} \ar[l]_{\phi} \ar[r] \ar[d] & \calC^{\acti}_{\overline{h}_1/} \ar[d] \\
\calC^{\acti}_{h_0/} \times_{\calO^{\acti}_{q  h_0/} } \calO^{\acti}_{q  \overline{h}_0/} &
\calC^{\acti}_{h/} \times_{\calO^{\acti}_{q  h/}} \calO^{\acti}_{q  \overline{h}/} \ar[r] \ar[l]_{\psi} &
\calC^{\acti}_{h_1/} \times_{\calO^{\acti}_{q  h_1/}} \calO^{\acti}_{q  \overline{h}_1/}. }$$
According to Remark \ref{coheath}, it suffices to show that the left vertical map is a
categorical equivalence if and only if the right vertical map is a categorical equivalence.
Because the inclusions $\{1\} \times K \subseteq \Delta^1 \times K$ and $\{1\} \times K^{\triangleright} \subseteq \Delta^1 \times K^{\triangleright}$ are right anodyne, the horizontal maps on the right are trivial fibrations. Using a diagram chase, we are reduced to proving that the maps $\phi$
and $\psi$ are categorical equivalences.

Let $f: x \rightarrow y$ denote the morphism $\calC^{\otimes}$ obtained by restricting $\overline{h}$ to the cone point of $K^{\triangleright}$. The map $\phi$ fits into a commutative diagram
$$ \xymatrix{ \calC^{\acti}_{\overline{h}/} \ar[r]^{\phi} \ar[d] & \calC^{\acti}_{h_0/} \ar[d] \\
\calC^{\acti}_{f/} \ar[r] & \calC^{\acti}_{x/}. }$$
Since the inclusion of the cone point into $K^{\triangleright}$ is right anodyne, the vertical arrows are trivial fibrations. Moreover, hypotheses $(1)$ and $(2)$ guarantee that $f$ is an equivalence in $\calC^{\otimes}$, so that the map $\calC^{\acti}_{f/} \rightarrow \calC^{\acti}_{x/}$ is a trivial fibration. This proves that $\phi$ is a categorical equivalence.

The map $\psi$ admits a factorization
$$ \calC^{\acti}_{h/} \times_{\calO^{\acti}_{q h/}} \calO^{\acti}_{q  \overline{h}/}
\stackrel{\psi'}{\rightarrow} 
\calC^{\acti}_{h_0/} \times_{ \calO^{\acti}_{q  h_0/}} \calO^{\acti}_{q  \overline{h}/}
\stackrel{\psi''}{\rightarrow}
\calO^{\acti}_{h_0} \times_{ \calO^{\acti}_{q  h_0/}} \calO^{\acti}_{ q  \overline{h}_0/}.$$
To complete the proof, it will suffice to show that $\psi'$ and $\psi''$ are trivial fibrations of simplicial sets. We first observe that $\psi'$ is a pullback of the map
$$\calC^{\acti}_{h/} \rightarrow \calC^{\acti}_{h_0/} \times_{\calO^{\acti}_{q  h_0/} } \calO^{\acti}_{q  h/},$$
which is a trivial fibration (by Proposition \toposref{eggwhite}). The map $\psi''$
is a pullback of the left fibration $\psi''_0: \calO^{\acti}_{q  \overline{h}/} \rightarrow \calO^{\acti}_{q  \overline{h}_0/}$. It therefore suffices to show that $\psi''_0$ is a categorical equivalence.
To prove this, we consider the diagram
$$ \xymatrix{ \calO^{\acti}_{q  \overline{h}/} \ar[r]^{\psi''_0} \ar[d] & \calO^{\acti}_{q  \overline{h}_0/} \ar[d] \\
\calO^{\acti}_{q(f)/} \ar[r]^{\psi''_1} & \calO^{\acti}_{q(x)/}. }$$
As above, we observe that the vertical arrows are trivial fibrations and that $\psi''_1$ is a trivial fibration (because the morphism $q(f)$ is an equivalence in $\calO^{\otimes}$). It follows that $\psi''_0$ is a categorical equivalence, as desired.
\end{proof}

In the case where $q: \calC^{\otimes} \rightarrow \calO^{\otimes}$ is a categorical fibration,
Proposition \ref{chocolateoperad} allows us to reduce the problem of testing whether
a diagram $\overline{p}: K^{\triangleright} \rightarrow \calC^{\otimes}_{\acti}$ is an
operadic $q$-colimit to the special case where $\overline{p}$ factors through
$\calC^{\otimes}_{X}$, for some $X \in \calO^{\otimes}$. In this case, we can apply
the following criterion:

\begin{proposition}\label{optest}
Let $\calO^{\otimes}$ be an $\infty$-operad, let $q: \calC^{\otimes} \rightarrow \calO^{\otimes}$
be a $\calO$-monoidal $\infty$-category, and let $\overline{p}: K^{\triangleright}
\rightarrow \calC^{\otimes}_{X}$ be a diagram for some $X \in \calO^{\otimes}$. The
following conditions are equivalent:
\begin{itemize}
\item[$(1)$] The map $\overline{p}$ is an operadic $q$-colimit diagram.
\item[$(2)$] For every object $\overline{Y} \in \calC^{\otimes}$ with image
$Y \in \calO^{\otimes}$, every object $Z \in \calO$, and every active morphism
$m: X \oplus Y \rightarrow Z$ in $\calO$, the induced functor
$$ \calC^{\otimes}_{X} \stackrel{ \oplus \overline{Y} }{\rightarrow}
\calC^{\otimes}_{X \oplus Y} \stackrel{m_{!}}{\rightarrow} \calC_{Z}$$
carries $\overline{p}$ to a colimit diagram in $\calC_{Z}$.
\end{itemize}
\end{proposition}

\begin{example}\label{coltens}
Let $q: \calC^{\otimes} \rightarrow \Nerve(\FinSeg)$ exhibit $\calC$ as a symmetric monoidal
$\infty$-category. Then a diagram $\overline{p}: K^{\triangleright} \rightarrow \calC$
is an operadic $q$-colimit diagram if and only if, for every object $C \in \calC$, the induced diagram
$$ K^{\triangleright} \stackrel{ \overline{p}}{\rightarrow} \calC \stackrel{ \otimes C}{\rightarrow} \calC$$
is a colimit diagram in $\calC$. More informally: an operadic colimit diagram in $\calC$ is a colimit diagram which remains a colimit diagram after tensoring with any object of $\calC$.
\end{example}

\begin{proof}[Proof of Proposition \ref{optest}]
Replacing $\overline{p}$ by its image under the functor $\calC^{\otimes}
\stackrel{ \oplus \overline{Y} }{\rightarrow} \calC^{\otimes}$, we are reduced to proving
the equivalence of the following pair of assertions:
\begin{itemize}
\item[$(1')$] The map $\overline{p}$ is a weak operadic $q$-colimit diagram.
\item[$(2')$] For every object $Z \in \calO$ and every active morphism
$m: X \rightarrow Z$ in $\calO$, the induced functor
$$ \calC^{\otimes}_{X} \stackrel{m_!}{\rightarrow} \calC_{Z}$$
carries $\overline{p}$ to a colimit diagram in $\calC_{Z}$.
\end{itemize}

Assertion $(1')$ is equivalent to the statement that the map
$$ \theta: \calC^{\acti}_{\overline{p}/} \rightarrow \calC^{\acti}_{p/} \times_{ \calO^{\acti}_{qp/} } \calO^{\acti}_{q \overline{p}/}$$
is a trivial fibration of simplicial sets. Since $\theta$ is a left fibration (Proposition \toposref{sharpen}), it will suffice to show that the fibers of $\theta$ are contractible. Consider an arbitrary vertex of 
$\calO^{\acti}_{q \overline{p}/}$ corresponding to a diagram $t: K \star \Delta^1 \rightarrow \calO^{\otimes}_{\acti}$. Since $K \star \Delta^1$
is categorically equivalent to $( K \star \{0\} ) \coprod_{ \{0\} } \Delta^1$ and
$t | K \star \{0\}$ is constant, we may assume without loss of generality that
$t$ factors as a composition
$$ K \star \Delta^1 \rightarrow \Delta^1 \stackrel{m}{\rightarrow} \calO^{\otimes}_{\acti}.$$
Here $m: X \rightarrow Z$ can be identified with an active morphism in $\calO^{\otimes}$. 
It will therefore suffice to show that the following assertions are equivalent, where
$m: X \rightarrow Z$ is {\em fixed}:

\begin{itemize}
\item[$(1'')$] The map $\calC^{\acti}_{ \overline{p}/} \times_{ \calO^{\acti}_{q \overline{p}/} } \{ t \}
\rightarrow \calC^{\acti}_{p/} \times_{ \calO^{\acti}_{qp/} } \{ t_0 \}$
where $t \in \calO^{\otimes}_{q \overline{p}/}$ is determined by $m$ as above, and $t_0$ is
the image of $t$ in $\calO^{\acti}_{qp/}$. 

\item[$(2'')$] The map $m_{!}: \calC^{\otimes}_{X} \rightarrow \calC_{Z}$ carries
$\overline{p}$ to a colimit diagram in $\calC_{Z}$.
\end{itemize}

This equivalence results from the observation that the fibers of the maps
$$ \calC^{\acti}_{ \overline{p}/ } \rightarrow \calO^{\acti}_{ q \overline{p}/}
\quad \quad \calC^{\acti}_{p/} \rightarrow \calO^{\acti}_{q p/}$$
can be identified with the fibers of the maps
$$ \calC^{\acti}_{ m_{!} \overline{p}/} \rightarrow \calO^{\acti}_{q m_{!} \overline{p}/}
\quad \quad \calC^{\acti}_{m_{!} p/} \rightarrow \calO^{\acti}_{q m_{!} p/}.$$
\end{proof}

We now establish a general criterion for the existence of operadic colimit diagrams.

\begin{definition}\label{compatK}
Let $\calO^{\otimes}$ be an $\infty$-operad, let $q: \calC^{\otimes} \rightarrow \calO^{\otimes}$ be
a coCartesian fibration of $\infty$-operads and let $K$ be a simplicial set. We will say that
the {\it $\calO$-monoidal structure on $\calC$ is compatible with $K$-indexed colimits}
if the following conditions are satisfied:
\begin{itemize}
\item[$(1)$] For every object $X \in \calO$, the $\infty$-category $\calC_{X}$ admits $K$-indexed colimits.
\item[$(2)$] For every operation $f \in \Mul_{\calO}( \{ X_i \}_{1 \leq i \leq n}, Y)$, the functor
$\otimes_{t}: \prod_{1 \leq i \leq n} \calC_{X_i} \rightarrow \calC_{Y}$ of Remark \ref{swinn}
preserves $K$-indexed colimits separately in each variable.
\end{itemize}
\end{definition}

\begin{proposition}\label{slavewell}
Let $\calO^{\otimes}$ be an $\infty$-operad, let $K$ be a simplicial set, and let
$q: \calC^{\otimes} \rightarrow \calO^{\otimes}$ be a $\calO$-monoidal structure on
$\calC$ which is compatible with $K$-indexed colimits. Let $p: K \rightarrow \calC^{\otimes}_{\acti}$ be a diagram, and let $\overline{p}_0: K^{\triangleright} \rightarrow \calO^{\otimes}_{\acti}$ be an
extension of $qp$ which carries the cone point of $K^{\triangleright}$ to an object
$X \in \calO$. Then there exists an operadic $q$-colimit diagram
$\overline{p}: K^{\triangleright} \rightarrow \calC^{\otimes}$ which extends $p$ and lifts
$\overline{p}_0$.
\end{proposition}

\begin{proof}
Let $p_0 = qp$. The map $\overline{p}_0$ determines a natural transformation
$\alpha: p_0 \rightarrow \underline{X}$ of diagrams $K \rightarrow \calO^{\otimes}$, where
$\underline{X}$ denotes the constant diagram taking the value $X$. Choose a 
$q$-coCartesian natural transformation $\overline{\alpha}: p \rightarrow p'$ lifting
$\alpha$. Since $\calC_{X}$ admits $K$-indexed colimits, we can extend
$\overline{p}'$ to a colimit diagram $\overline{p}': K^{\triangleright} \rightarrow \calC_{X}$.
The compatibility of the $\calO$-monoidal structure on $\calC$ with $K$-indexed colimits
and Proposition \ref{optest} imply that $\overline{p}'$ is an operadic $q$-colimit diagram.
Let $C \in \calC_{X}$ be the image under $\overline{p}'$ of the cone point of
$K^{\triangleright}$, so we can regard $\overline{p}'$ as a diagram
$K \rightarrow \calC^{\otimes}_{/C}$ which lifts $p'$. Using the assumption that
$\calC^{\otimes}_{/C} \rightarrow 
\calC^{\otimes} \times_{ \calO^{\otimes} } \calO^{\otimes}_{/X}$ is a right fibration, we can choose
a transformation $\overline{p} \rightarrow \overline{p}'$ lifting $\overline{\alpha}: p \rightarrow p'$.
It follows from Proposition \ref{chocolateoperad} that $\overline{p}$ is an extension of $p$
with the desired properties.
\end{proof}

\begin{corollary}\label{cyster}
Let $q: \calC^{\otimes} \rightarrow \calO^{\otimes}$ be a fibration of $\infty$-operads and
let $\kappa$ be a regular cardinal. The following conditions are equivalent:
\begin{itemize}
\item[$(1)$] The map $q$ is a coCartesian fibration of $\infty$-operads, and the
$\calO$-monoidal structure on $\calC$ is compatible with $\kappa$-small colimits.
\item[$(2)$] For every $\kappa$-small simplicial set $K$ and every diagram
$$ \xymatrix{ K \ar[r]^{p} \ar[d] & \calC^{\otimes}_{\acti} \ar[d] \\
K^{\triangleright} \ar@{-->}[ur]^{ \overline{p} } \ar[r]^{\overline{p}_0} & \calO^{\otimes} }$$
such that $\overline{p}_0$ carries the cone point of $K^{\triangleright}$ into $\calO$, there
exists an extension $\overline{p}$ of $p$ as indicated in the diagram, which is an
operadic $q$-colimit diagram.
\end{itemize}
\end{corollary}

\begin{proof}
The implication $(1) \Rightarrow (2)$ is immediate from Proposition \ref{slavewell}.
Assume that $(2)$ is satisfied. Taking $K = \Delta^0$ and applying Proposition
\ref{justica}, we deduce that $q$ is a coCartesian fibration of $\infty$-operads.
Applying $(2)$ in the case where the map $\overline{p}_0$ takes some constant value
$X \in \calO$, we deduce that every $\kappa$-small diagram $K \rightarrow \calC^{\otimes}_{X}$
can be extended to an operadic $q$-colimit diagram. This operadic $q$-colimit diagram is in particular
a colimit diagram, so that $\calC^{\otimes}_{X}$ admits $\kappa$-small colimits. Moreover, the uniqueness properties of colimits show that every $\kappa$-small colimit diagram
$\overline{p}: K^{\triangleright} \rightarrow \calC^{\otimes}_{X}$ is a $q$-operadic colimit diagram. 
Combining this with the criterion of Proposition \ref{optest}, we conclude that
the $\calO$-monoidal structure on $\calC$ is compatible with $\kappa$-small colimits.
\end{proof}

\subsection{Unit Objects and Trivial Algebras}\label{unitr}

Let $\calO^{\otimes}$ be a unital $\infty$-operad, and let $p: \calC^{\otimes} \rightarrow \calO^{\otimes}$ be a coCartesian fibration of $\infty$-operads. For each object $X \in \calO$, there is an essentially
unique morphism $0 \rightarrow X$ (where $0$ denotes the zero object of $\calO^{\otimes}$), which
determines a functor $\calC^{\otimes}_0 \rightarrow \calC_{X}$. Since $\calC^{\otimes}_0$
is a contractible Kan complex, we can identify this functor with an object of ${\bf 1}_{X} \in \calC_{X}$.
We will refer to the object ${\bf 1}_{X}$ as the {\it unit object} of $\calC$ (of color $X$).
Our goal in this section is to study the basic features of these unit objects. We begin by formulating a slightly more general definition, which makes sense even if we do not assume that
$p$ is a coCartesian fibration.

\begin{definition}\label{cuspe}
Let $p: \calC^{\otimes} \rightarrow \calO^{\otimes}$ be a fibration of $\infty$-operads,
let $X \in \calO^{\otimes}$, and let $f: C_0 \rightarrow {\bf 1}_{X}$ be a morphism in 
$\calC^{\otimes}$, where ${\bf 1}_{X} \in \calC_{X}$. 
We will say that {\it $f$ exhibits ${\bf 1}_{X}$ as an $X$-unit object} if the following conditions are satisfied:
\begin{itemize}
\item[$(1)$] The object $C_0$ belongs to $\calC^{\otimes}_{\seg{0}}$. 
\item[$(2)$] The morphism $f$ determines an operadic $p$-colimit diagram
$(\Delta^0)^{\triangleright} \rightarrow \calC^{\otimes}$.
\end{itemize}
More generally, we will say that an arbitrary morphism $C_0 \rightarrow C$
in $\calC^{\otimes}$ {\it exhibits $C$ as a unit object} if, for every inert morphism
$C \rightarrow C'$ with $C' \in \calC$, the composite map
$C_0 \rightarrow C'$ exhibits $C'$ as a $p(C')$-unit object.

Suppose that $\calO^{\otimes}$ is a unital $\infty$-operad. We will say that a fibration
$\calC^{\otimes} \rightarrow \calO^{\otimes}$ {\it has unit objects} if,
for every object $X \in \calO^{\otimes}$, there
exists a morphism $f: C_0 \rightarrow {\bf 1}_{X}$ in $\calC^{\otimes}$ which exhibits
${\bf 1}_{X}$ as an $X$-unit object.
\end{definition}

\begin{remark}
The terminology of Definition \ref{cuspe} is slightly abusive: the condition that a morphism
$f: C_0 \rightarrow C$ exhibits $C$ as a unit object of $\calC^{\otimes}$ depends not only
on the $\infty$-operad $\calC^{\otimes}$, but also on the $\infty$-operad fibration
$\calC^{\otimes} \rightarrow \calO^{\otimes}$.
\end{remark}

\begin{example}\label{swagg}
Corollary \ref{swink} implies that if $p: \calC^{\otimes} \rightarrow \calO^{\otimes}$ is a fibration of $\infty$-operads and $\calO^{\otimes}$ is coherent, then $\Mod^{\calO}(\calC)^{\otimes} \rightarrow \calO^{\otimes}$ has units.
\end{example}

\begin{remark}\label{coope}
Let $p: \calC^{\otimes} \rightarrow \calO^{\otimes}$ be a fibration of $\infty$-operads.
If $X \in \calO$ and $f: C_0 \rightarrow C$ is a morphism in $\calC^{\otimes}$ which exhibits 
$C$ as a unit object, then $f$ is $p$-coCartesian: this follows immediately from
Proposition \ref{copan}.
\end{remark}

\begin{remark}\label{toopse}
Let $p: \calC^{\otimes} \rightarrow \calO^{\otimes}$ be a fibration of $\infty$-operads, and assume that $\calO^{\otimes}$ is unital. Let $\calX$ denote the full subcategory of $\Fun( \Delta^1, \calC^{\otimes})$ spanned by those morphisms $f: C_0 \rightarrow C$ which exhibit $C$ as a unit object.
Then the composite map
$$ \phi: \calX \subseteq \Fun( \Delta^1, \calC^{\otimes}) \rightarrow
\Fun( \{1\}, \calC^{\otimes}) \stackrel{p}{\rightarrow} \calO^{\otimes}$$
induces a trivial Kan fibration onto a full subcategory of $\calO^{\otimes}$. To see
this, we let $\calX'$ denote the full subcategory of $\Fun( \Delta^1, \calC^{\otimes})$ spanned by those morphisms $f: C_0 \rightarrow C$ which are $p$-coCartesian and satisfy $C_0 \in \calC^{\otimes}_{\seg{0}}$. Proposition \toposref{lklk} implies that the map
$$ \calX' \rightarrow \calC^{\otimes}_{\seg{0}} \times_{ \Fun( \{0\}, \calO^{\otimes}) }
\Fun( \Delta^1, \calO^{\otimes})$$ 
is a trivial Kan fibration onto a full subcategory of the fiber product
$\calX'' = \calC^{\otimes}_{\seg{0}} \times_{ \Fun( \{0\}, \calO^{\otimes}) }
\Fun( \Delta^1, \calO^{\otimes})$.
Since $\calC^{\otimes}_{\seg{0}} \rightarrow \calO^{\otimes}_{\seg{0}}$ is a
categorical fibration between contractible Kan complexes, it is a trivial Kan fibration; it follows
that the induced map
$$ \calX'' \rightarrow \calO^{\otimes}_{ \seg{0}} \times_{ \Fun( \{0\}, \calO^{\otimes} )} \Fun( \Delta^1, \calO^{\otimes})$$
is also a trivial Kan fibration. The target of this map can be identified with the $\infty$-category
$\calO^{\otimes}_{\ast}$ of pointed objects of $\calO^{\otimes}$, and the forgetful functor
$\calO^{\otimes}_{\ast} \rightarrow \calO^{\otimes}$ is a trivial Kan fibration because
$\calO^{\otimes}$ is unital. The desired result now follows by observing that $\phi$ is given by the composition
$$ \calX \subseteq \calX' \rightarrow \calX'' \rightarrow \calO^{\otimes}_{\ast} \rightarrow \calO^{\otimes}.$$

We can summarize our discussion as follows: if $f: C_0 \rightarrow C$ is a morphism in
$\calC^{\otimes}$ which exhibits $C$ as a unit object, then $f$ is determined (up to canonical equivalence) by the object $p(C) \in \calO^{\otimes}$. We observe that
$\calC^{\otimes} \rightarrow \calO^{\otimes}$ has unit objects if and only if the functor $\phi$ above is essentially surjective. In this case
the inclusion $\calX \rightarrow \calX'$ must also be essentially surjective: in other words,
a morphism $f: C_0 \rightarrow C$ exhibits $C$ as a unit object if and only if $f$ is $p$-coCartesian
and $C_0 \in \calC^{\otimes}_{\seg{0}}$. 
\end{remark}

\begin{remark}\label{carver}.
It is easy to see that the essential image of the functor $\phi$ of 
Remark \ref{toopse} is stable under the concatenation functor
$\oplus$ of Remark \ref{opl}. Consequently, to show that a fibration of $\infty$-operads $\calC^{\otimes} \rightarrow
\calO^{\otimes}$ has units objects, it suffices to show that for every object $X \in \calO$ 
there exists a map $f: C_0 \rightarrow {\bf 1}_{X}$ which exhibits ${\bf 1}_{X}$ as an
$X$-unit object of $\calC$.
\end{remark}

\begin{example}
Let $p: \calC^{\otimes} \rightarrow \calO^{\otimes}$ be a coCartesian fibration of
$\infty$-operads, where $\calO^{\otimes}$ is unital. Then $p$ has unit objects. This
follows immediately from Proposition \ref{slavewell} and Remark \ref{carver}.
\end{example}

\begin{definition}
Let $p: \calC^{\otimes} \rightarrow \calO^{\otimes}$ be a fibration of $\infty$-operads, and assume that $\calO^{\otimes}$ is unital. We will say that an algebra object $A \in \Alg_{\calO}(\calC)$ is
{\it trivial} if, for every object $X \in \calO^{\otimes}$, the induced map $A(0) \rightarrow A(X)$ exhibits
$A(X)$ as a unit object; here $0$ denotes the zero object of $\calO^{\otimes}$.
\end{definition}

When trivial algebra objects exist, they are precisely the initial objects of $\Alg_{\calO}(\calC)$:

\begin{proposition}\label{gargle1}
Let $p: \calC^{\otimes} \rightarrow \calO^{\otimes}$ be a fibration of $\infty$-operads, where
$\calO^{\otimes}$ is unital. Assume that $p$ has unit objects. The following conditions on a $\calO$-algebra object $A \in \Alg_{\calO}(\calC)$ are equivalent:
\begin{itemize}
\item[$(1)$] The object $A$ is an initial object of $\Alg_{\calO}(\calC)$.
\item[$(2)$] The functor $A$ is a $p$-left Kan extension of $A| \calO^{\otimes}_{\seg{0}}$.
\item[$(3)$] The algebra object $A$ is trivial.
\end{itemize}
\end{proposition}

\begin{corollary}\label{gargle2}
Let $\calC^{\otimes}$ be a symmetric monoidal $\infty$-category. Then
the $\infty$-category $\CAlg(\calC)$ has an initial object. Moreover, a commutative
algebra object $A$ of $\calC$ is an initial object of $\CAlg(\calC)$ if and only if the
unit map ${\bf 1} \rightarrow A$ is an equivalence in $\calC$.
\end{corollary}

The proof of Proposition \ref{gargle1} depends on the following:

\begin{lemma}\label{kanus}
Let $p: \calC^{\otimes} \rightarrow \calO^{\otimes}$ be a fibration of $\infty$-operads, where
$\calO^{\otimes}$ is unital. The following conditions are equivalent:
\begin{itemize}
\item[$(1)$] The fibration $p$ has unit objects.
\item[$(2)$] There exists a trivial $\calO$-algebra object $A \in \Alg_{\calO}(\calC)$.
\end{itemize}
\end{lemma}

\begin{proof}
The implication $(2) \Rightarrow (1)$ follows immediately from Remark \ref{carver}.
Conversely, suppose that $(1)$ is satisfied. Choose an arbitrary section
$s$ of the trivial Kan fibration $\calC^{\otimes}_{\seg{0}} \rightarrow \calO^{\otimes}_{\seg{0}}$.
Since $p$ has units, Lemma \toposref{kan2} (and Remark \ref{coope}) imply that
$s$ admits a $p$-left Kan extension $A: \calO^{\otimes} \rightarrow \calC^{\otimes}$.
For every morphism $f: X \rightarrow Y$ in $\calO^{\otimes}$, we have a commutative diagram
$$ \xymatrix{ & X \ar[dr]^{f} & \\
0 \ar[ur]^{g} \ar[rr]^{h} & & Y}$$
in $\calO^{\otimes}$. Since $A(g)$ and $A(h)$ are $p$-coCartesian morphisms of
$\calC^{\otimes}$, Proposition \toposref{protohermes} implies that $f$ is $p$-coCartesian.
It follows in particular that $A$ preserves inert morphisms, so that $A \in \Alg_{\calO}(\calC)$. 
Using Remark \ref{toopse}, we conclude that $A$ is trivial.
\end{proof}

\begin{proof}[Proof of Proposition \ref{gargle1}]
The implication $(2) \Rightarrow (3)$ follows from the proof of Lemma \ref{kanus}, and
$(3) \Rightarrow (2)$ follows immediately from Remark \ref{coope}. The implication
$(2) \Rightarrow (1)$ follows from Proposition \toposref{leftkanadj} (since
the $\infty$-category of sections of the trivial Kan fibration $\calC^{\otimes}_{\seg{0}} \rightarrow
\calO^{\otimes}_{\seg{0}}$ is a contractible Kan complex). The reverse implication
$(1) \Rightarrow (2)$ follows from the uniqueness of initial objects up to equivalence, since
Lemma \ref{kanus} guarantees the existence of an object $A \in \Alg_{\calO}(\calC)$ satisfying
$(3)$ (and therefore $(2)$ and $(1)$ as well).
\end{proof}

For the remainder of this section, we will assume familiarity with the theory of
module objects developed in \S \ref{comm6}. Our goal is to prove the following result, which
asserts that trivial algebras have equally trivial module categories:

\begin{proposition}\label{stupe}
Let $p: \calC^{\otimes} \rightarrow \calO^{\otimes}$ be a fibration of $\infty$-operads, and assume that
$\calO^{\otimes}$ is coherent. Let $A$ be a trivial $\calO$-algebra object of $\calC$.
Then the forgetful functor $\theta: \Mod^{\calO}_{A}(\calC)^{\otimes} \rightarrow \calC^{\otimes}$ is an equivalence of $\infty$-categories.
\end{proposition}

The proof will require a few preliminary results.

\begin{lemma}\label{loopse}
Let $p: X \rightarrow S$ be an inner fibration of simplicial sets.
Let $X^{0}$ be a full simplicial subset of $X$, and assume that the restriction
map $p^0 = p | X^0$ is a coCartesian fibration. Let $q: Y \rightarrow Z$ be a categorical fibration of simplicial sets. Define a simplicial sets $A$ and $B$ equipped with a maps $A,B \rightarrow S$ so that the following universal property is satisfied: for every map of simplicial sets $K \rightarrow S$, we have
bijections 
$$\Hom_{ (\sSet)_{/S}}( K, A) \simeq \Fun(X \times_{S} K, Y)$$
$$\Hom_{ ( \sSet)_{/S}}( K, B) \simeq  \Fun(X^{0} \times_{S} K, Y)
\times_{  \Fun(X^{0} \times_{S} K, Z) }  \Fun(X \times_{S} K, Z).$$
Let $\phi: A \rightarrow B$ be the restriction map.
Let $A'$ denote the full simplicial subset of $A$ spanned by those vertices corresponding
to maps $f: X_{s} \rightarrow Y$ such that $f$ is a $q$-left Kan extension of $f| X^{0}_{s}$, and
let $B'$ denote the full simplicial subset of $B$ spanned by those vertices of the form
$\phi(f)$ where $f \in A'$. Then $\phi$ induces a trivial Kan fibration $\phi': A' \rightarrow B'$.
\end{lemma}

\begin{proof}
For every map of simplicial sets $T \rightarrow S$, let $F(T) = \bHom_{ S}( T, A')$
and let $G(T) = \bHom_{S}( T, B')$. If $T_0 \subseteq T$ is a simplicial subset, we
have a restriction map $\theta_{T_0,T}: F(T) \rightarrow G(T) \times_{ G(T_0)} F(T_0)$. To prove
that $\phi'$ is a trivial Kan fibration, it will suffice to show that
$\theta_{T_0,T}$ is surjective on vertices whenever $T = \Delta^n$ and $T_0 = \bd \Delta^n$.We will complete the proof by showing that $\theta_{T_0, T}$ is a trivial Kan fibration whenever $T$ has only finitely many simplices.

The proof proceeds by induction on the dimension of $T$ (if $T$ is empty, the result is trivial). Assume first that $T = \Delta^n$. If $T_0 = T$ there is nothing to prove. Otherwise, we may assume that
$T_0$ has dimension smaller than $n$. Using
the fact that $q$ is a categorical fibration, we deduce that $\theta_{T_0, T}$ is a categorical
fibration. It therefore suffices to show that $\theta_{T_0,T}$ is a categorical equivalence.
We have a commutative diagram
$$ \xymatrix{ F(T) \ar[dr]^{\psi} \ar[rr] & & G(T) \times_{ G(T_0)} F(T_0) \ar[dl]^{\psi'} \\
& G(T). & }$$ 
The inductive hypothesis (applied to the inclusion $\emptyset \subseteq T_0$) guarantees
that $\psi'$ is a trivial Kan fibration. It will therefore suffice to show that $\psi$ is a trivial
Kan fibration. In view of Proposition \toposref{lklk}, it will suffice to prove the following assertions:
\begin{itemize}
\item[$(a)$] Let $F: \Delta^n \times_{S} X \rightarrow Y$ be a functor. Then
$F$ is a $q$-left Kan extension of $F|(\Delta^n \times_{S} X^0)$ if and only if
$F|(\{i\} \times_{S} X)$ is a $q$-left Kan extension of $F| ( \{i\} \times_{S} X^{0} )$ for
$0 \leq i \leq n$. 

\item[$(b)$] Suppose given a commutative diagram 
$$ \xymatrix{ \Delta^n \times_{S} X^0 \ar[r]^{f} \ar[d] & Y \ar[d]^{q} \\
\Delta^n \times_{S} X \ar[r]^{g} \ar@{-->}[ur]^{F} & Z }$$
such that for $0 \leq i \leq n$, there exists a $q$-left Kan extension
$F_i: \{i\} \times_{S} X \rightarrow Y$ of $f | \{i\} \times_{S} X^0$
which is compatible with $g$. Then there exists a dotted arrow $F$
as indicated satisfying the condition described in $(a)$. 
\end{itemize}

Assertion $(a)$ follows from the observation that for $x \in \{i\} \times_{S} X$, the
assumption that $p^{0}$ is a coCartesian fibration guarantees that
$X^{0} \times_{X} ( \{i\} \times_{S} X)_{/x}$ is cofinal in
$X^{0} \times_{X} ( \Delta^n \times_{S} X)_{/x}$. Assertion $(b)$ follows from
the same observation together with Lemma \toposref{kan2}.

We now complete the proof by considering the case where $T$ is not a simplex.
We use induction on the number $k$ of simplices of $T'$ which do not belong to $T$.
If $k = 0$, then $T' = T$ and there is nothing to prove. If $k = 1$, then there is a pushout diagram
$$ \xymatrix{ \bd \Delta^n \ar[r] \ar[d] & T_0 \ar[d] \\
\Delta^n \ar[r] & T. }$$
It follows that $\theta_{T_0, T}$ is a pullback of the map $\theta_{ \bd \Delta^n, \Delta^n}$, and we
are reduced to the case where $T$ is a simplex. If $k > 1$, then we have nontrivial inclusions
$T_0 \subset T_1 \subset T$. Using the inductive hypothesis, we conclude that
$\theta_{T_1, T}$ and $\theta_{T_0, T_1}$ are trivial Kan fibrations. The desired result follows from the observation that $\theta_{T_0,T}$ can be obtained by composing $\theta_{T_1,T}$ with a pullback
of the morphism $\theta_{T_0, T_1}$.
\end{proof}

\begin{lemma}\label{kolmat}
Let $p: \calC^{\otimes} \rightarrow \calO^{\otimes}$ be a fibration of $\infty$-operads.
Assume that $\calO^{\otimes}$ is unital and that $p$ has unit objects. Let
$C \in \calC^{\otimes}$ and let $\alpha: p(C) \rightarrow Y$ be a semi-inert morphism
in $\calO^{\otimes}$. Then $\alpha$ can be lifted to a $p$-coCartesian morphism
$\overline{\alpha}: C \rightarrow \overline{Y}$ in $\calC^{\otimes}$.
\end{lemma}

\begin{proof}
The map $\alpha$ can be factored as the composition of an inert morphism and an active morphism.
We may therefore reduce to the case where $\alpha$ is either active or inert. If $\alpha$ is inert, we can choose $\overline{\alpha}$ to be an inert morphism lifting $\alpha$. Assume therefore that $\alpha$ is active. Write $C \simeq C_1 \oplus \ldots \oplus C_m$ (using the notation of Remark \ref{opl}),
and write $Y = p(C_1) \oplus \ldots \oplus C_m \oplus Y_1 \oplus \ldots \oplus Y_n$. 
Since $\calO^{\otimes}$ is unital, we may assume that $\alpha$ has the form
$\id_{p(C_1)} \oplus \ldots \oplus \id_{ p(C_m)} \oplus \alpha_1 \oplus \ldots \oplus \alpha_n$,
where each $\alpha_i: 0 \rightarrow Y_i$ is a morphism with $0 \in \calO^{\otimes}_{\seg{0}}$.
Since $p$ has units, we can lift each $\alpha_{i}$ to a morphism
$\overline{\alpha}_{i}: \overline{0}_i \rightarrow \overline{Y}_i$ which exhibits
$\overline{Y}_i$ as a $p$-unit. Let $\overline{\alpha} = \id_{C_1} \oplus \ldots
\oplus \id_{C_m} \oplus \overline{\alpha}_1 \oplus \ldots \oplus \overline{\alpha}_{n}$.
It follows from Proposition \ref{copan} that $\overline{\alpha}$ is $p$-coCartesian.
\end{proof}

\begin{proof}[Proof of Proposition \ref{stupe}]
We may assume that $p$ has unit objects (otherwise the assertion is vacuous). Let
$\phi: \calO^{\otimes} \times \Alg_{\calO}(\calC) \rightarrow \PrAlg_{\calO}(\calC)$
be the equivalence of Remark \ref{user2}, and let $\calX \subseteq \PrAlg_{\calO}(\calC)$
denote the essential image of the full subcategory spanned by those pairs $(X, A)$ where
$A$ is trivial. Let $\calX'$ denote the fiber product $\calX \times_{ \PrAlg_{\calO}(\calC)} \MMod^{\calO}(\calC)^{\otimes}$. Since Proposition \ref{gargle1} implies that trivial $\calO$-algebras form a contractible Kan complex, the inclusion $\Mod^{\calO}_{A}(\calC)^{\otimes} \subseteq \calX'$ is a categorical equivalence. It will therefore suffice to show that composition with the diagonal
map $\delta: \calO^{\otimes} \rightarrow \calK_{\calO}$ induces a categorical equivalence
$\calX' \rightarrow \calC^{\otimes}$.

Let $\calK^{1}_{\calO}$ denote the essential image of $\delta$, and define a simplicial
set $\calY$ equipped with a map $\calY \rightarrow \calO^{\otimes}$ so that the following
universal property is satisfied: for every map of simplicial sets $K \rightarrow \calO^{\otimes}$, we have a canonical bijection
$$ \Hom_{ (\sSet)_{/ \calO^{\otimes}}}( K, \calY)
\simeq \Hom_{ (\sSet)_{/ \calO^{\otimes}}}( K \times_{ \Fun( \{0\}, \calO^{\otimes})} \calK^{1}_{\calO}, \calC^{\otimes}).$$
Since $\delta$ is fully faithful, it induces a categorical equivalence $\calO^{\otimes} \rightarrow \calK^{1}_{\calO}$. It follows that the canonical map $\calY \rightarrow \calC^{\otimes}$ is a categorical equivalence.

We have a commutative diagram
$$ \xymatrix{ \calX' \ar[rr]^{\theta'} \ar[dr] & & \calY \ar[dl] \\
& \calC^{\otimes}. & }$$
Consequently, it will suffice to show that $\theta'$ is a categorical equivalence.
We will prove that $\theta'$ is a trivial Kan fibration. 

Define a simplicial set $\calD$ equipped with a map $\calD \rightarrow \calO^{\otimes}$ so that the following universal property is satisfied: for every map of simplicial sets $K \rightarrow \calO^{\otimes}$, we have a canonical bijection
$$ \Hom_{ (\sSet)_{/ \calO^{\otimes}}}( K, \calD)
\simeq \Hom_{ (\sSet)_{/ \calO^{\otimes}}}( K \times_{ \Fun( \{0\}, \calO^{\otimes})} \calK_{\calO}, \calC^{\otimes}).$$
For each $X \in \calO^{\otimes}$, let $A_{X}$ denote the full subcategory of
$(\calO^{\otimes})^{X/}$ spanned by the semi-inert morphisms $X \rightarrow Y$ in
$\calO^{\otimes}$, and let $A^{1}_{X}$ denote the full subcategory of
$(\calO^{\otimes})^{X/}$ spanned by the equivalences $X \rightarrow Y$ in
$\calO^{\otimes}$. An object of $\calD$ can be identified with a pair $(X, F)$, where $X \in \calO^{\otimes}$ and
$F: A_{X} \rightarrow \calC^{\otimes}$ is a functor. We will prove the following:
\begin{itemize}
\item[$(a)$] The full subcategory $\calX' \subseteq \calD$ is spanned by those pairs
$(X,F)$ where $F: A_{X} \rightarrow \calC^{\otimes}$ is a $q$-left Kan extension of
$F|A_{X}^{1}$.
\item[$(b)$] For every $X \in \calO^{\otimes}$ and every functor $f \in \Fun_{ \calO^{\otimes}}( A_{X}^{1}, \calC^{\otimes})$, there
exists a $q$-left Kan extension $F \in \Fun_{ \calO^{\otimes}}( A_{X}, \calC^{\otimes})$ of $f$.
\end{itemize}
Assuming that $(a)$ and $(b)$ are satisfied, the fact that the restriction functor
$\calX' \rightarrow \calY$ is a trivial Kan fibration will follow immediately from Lemma \ref{loopse}.

Note that for $X \in \calO^{\otimes}$, we can identify $A^{1}_{X}$ with the full subcategory
of $A_{X}$ spanned by the initial objects. Consequently, a functor $f \in \Fun_{ \calO^{\otimes}}( A_{X}^{1}, \calC^{\otimes})$ as in $(b)$ is determined up to equivalence by $f( \id_{X}) \in 
\calC^{\otimes}_{X}$. Using Lemma \toposref{kan2}, we deduce that $f$ admits a
$q$-left Kan extension $F \in \Fun_{ \calO^{\otimes}}( A_{X}, \calC^{\otimes})$ if and only if
every semi-inert morphism $X \rightarrow Y$ in $\calO^{\otimes}$ can be lifted to
a $q$-coCartesian morphism $f(\id_{X}) \rightarrow \overline{Y}$ in $\calC^{\otimes}$.
Assertion $(b)$ now follows from Lemma \ref{kolmat}.

We now prove $(a)$. Suppose first that $F$ is a $q$-left Kan extension of
$F| A_{X}^{1}$. The proof of $(b)$ shows that $F(u)$ is $q$-coCartesian for
every morphism $u: Y \rightarrow Z$ in $A_{X}$ such that $Y$ is an initial object
of $A_{X}$. Since every morphism $u$ in $A_{X}$ fits into a commutative diagram
$$ \xymatrix{ & Y \ar[dr]^{u} & \\
\id_{X} \ar[ur] \ar[rr] & & Z, }$$
Proposition \toposref{protohermes} guarantees that $F$ carries {\em every} morphism
in $A_{X}$ to a $q$-coCartesian morphism in $\calC^{\otimes}$. In particular,
$F$ carries inert morphisms in $A_{X}$ to inert morphisms in $\calC^{\otimes}$, and
therefore belongs to $\MMod^{\calO}(\calC)^{\otimes} \times_{ \calO^{\otimes} } \{ X \}$.
Let $A_{X}^{0}$ denote the full subcategory of $A_{X}$ spanned by the null morphisms
$X \rightarrow Y$ in $\calO^{\otimes}$, and let $s: \calO^{\otimes} \rightarrow A_{X}^{0}$ denote a section to the trivial Kan fibration $A_{X}^{0} \rightarrow \calO^{\otimes}$. 
To prove that $F \in \calX'$, it suffices to show that the composition
$$\calO^{\otimes} \stackrel{s}{\rightarrow} A_{X}^{0} \subseteq A_{X} \stackrel{F}{\rightarrow} \calC^{\otimes}$$
is a trivial $\calO$-algebra. Since this composition carries every morphism in $\calO^{\otimes}$ to
a $q$-coCartesian morphism in $\calC^{\otimes}$, it is an $\calO$-algebra: the triviality now follows from Remark \ref{toopse}.

To complete the proof of $(a)$, let us suppose that $F \in \calX'$; we wish to show that
$F$ is a $q$-left Kan extension of $f = F | \calA^{1}_{X}$. Using $(b)$, we deduce that
$f$ admits a $q$-left Kan extension $F' \in \Fun_{ \calO^{\otimes}}( A_{X}, \calC^{\otimes})$.
Let $\alpha: F' \rightarrow F$ be a natural transformation which is the identity on $f$;
we wish to prove that $\alpha$ is an equivalence. Fix an object $\overline{Y} \in A_{X}$,
corresponding to a semi-inert morphism $X \rightarrow Y$ in $\calO^{\otimes}$.
Let $\seg{n} \in \Nerve(\FinSeg)$ denote the image of $Y$, and choose inert morphisms
$Y \rightarrow Y_i$ lifting the maps $\Colp{i}: \seg{n} \rightarrow \seg{1}$.
Let $\overline{Y}_i$ denote the composition of $\overline{Y}$ with the map
$Y \rightarrow Y_i$ for $1 \leq i \leq n$.
Since $F$ and $F'$ both preserve inert morphisms and $\calC^{\otimes}$ is an $\infty$-operad,
it suffices to prove that $\alpha_{\overline{Y}_i}: F'( \overline{Y}_i) \rightarrow F( \overline{Y}_i)$ is an equivalence for $1 \leq i \leq n$. We may therefore replace $Y$ by $Y_i$ and reduce to the case where $Y \in \calO$. In this case, the semi-inert morphism $\overline{Y}$ is either null or inert.

If the map $\overline{Y}: X \rightarrow Y$ is null, then $\overline{Y} \in A_{X}^{0}$.
Since $F \circ s$ and $F' \circ s$ both determine trivial $\calO$-algebra objects, the induced natural transformation $F' \circ s \rightarrow F \circ s$ is an equivalence (Proposition \ref{gargle1}).
It follows that the natural transformation $F' | A^{0}_{X} \rightarrow F| A^{0}_{X}$ is an equivalence,
so that $F'( \overline{Y}) \simeq F( \overline{Y})$.

If the map $\overline{Y}: X \rightarrow Y$ is inert, then we have an inert morphism
$u: \id_{X} \rightarrow \overline{Y}$ in $A_{X}$. Since $F$ and $F'$ both preserve inert morphisms, it suffices to show that the map $F'( \id_{X}) \rightarrow F(\id_{X})$ is an equivalence. This is clear,
since $\id_{X} \in A^{1}_{X}$.
\end{proof}

\subsection{Strictly Unital $\infty$-Categories}\label{strictun}

Let $\calC$ be an $\infty$-category. In \S \toposref{minin}, we saw that $\calC$ admits
a {\it minimal model}: that is, there exists an $\infty$-category $\calC_0 \subseteq \calC$
which is categorically equivalent to $\calC$, which is in some sense as small as possible
(and is determined up to noncanonical isomorphism). The theory of minimal models is a useful technical device: it allows us to to reduce problems about arbitrary $\infty$-categories to problems about
minimal $\infty$-categories. Our goal in this section is to expand the usefulness of this device: we will introduce a slightly larger class of $\infty$-categories, which we call {\it strictly unital $\infty$-categories}.
These $\infty$-categories possess some of the crucial features of minimal $\infty$-categories
(Proposition \ref{appstrict}), but also enjoy closure properties which are not shared by minimal
$\infty$-categories (Remark \ref{tope2}).

\begin{definition}
Let $\calC$ be an $\infty$-category. We will say that $\calC$ is {\it strictly unital} if the following
condition is satisfied:
\begin{itemize}
\item[$(\ast)$] Let $f: \Delta^n \rightarrow \calC$ be a morphism with the property that
$f | \Delta^{ \{i, i+1 \} }$ is a degenerate edge of $\calC$, for $0 \leq i < n$. Then $f$ factors as a composition 
$$ \Delta^n \stackrel{f'}{\rightarrow} \Delta^{n-1} \stackrel{f''}{\rightarrow} \calC,$$
where $f'$ is given on vertices by the formula
$$ f'( j) = \begin{cases} j & \text{if } j \leq i \\
j-1 & \text{if } j > i. \end{cases}$$
\end{itemize}
More generally, suppose that $p: X \rightarrow S$ is an inner fibration of simplicial sets.
We will say that $p$ is {\it strictly unital} if, for every $k$-simplex $\Delta^{k} \rightarrow S$, the
fiber product $\Delta^{k} \times_{S} X$ is a strictly unital $\infty$-category.
\end{definition}

\begin{remark}\label{conse}
An inner fibration of simplicial sets $p: X \rightarrow S$ is strictly unital if and only if it satisfies
the following condition:
\begin{itemize}
\item[$(\ast)$] Let $0 \leq i < n$, and let $f': \Delta^n \rightarrow \Delta^{n-1}$ be the surjective
map which collapses the edge $\Delta^{ \{i, i+1\} }$ to the vertex $\{i\} \subseteq \Delta^{n-1}$.
Suppose that $f: \Delta^{n} \rightarrow X$ is a map whose restriction to
$\Delta^{ \{i, i+1\} }$ is a degenerate edge of $X$. Then $f$ factors through $f'$
if and only if $p \circ f$ factors through $f'$. 
\end{itemize}
\end{remark}

\begin{remark}\label{tope1}
An $\infty$-category $\calC$ is strictly unital if and only if the inner fibration
$\calC \rightarrow \Delta^0$ is strictly unital.
\end{remark}

\begin{remark}\label{tope2}
Every isomorphism of simplicial sets is a strictly unital inner fibration. Moreover, the collection 
of strictly unital inner fibrations is closed under composition (this follows immediately from
Remark \ref{conse}).
\end{remark}

\begin{example}\label{cose}
Let $\calC$ be a minimal $\infty$-category. Then $\calC$ is strictly unital. This is just a
translation of Proposition \toposref{minstrict}.
\end{example}

\begin{example}\label{telos}
Every minimal inner fibration $p: X \rightarrow S$ is strictly unital. This follows
from Example \ref{cose} and Remark \toposref{tablor}. 
\end{example}

\begin{proposition}\label{appstrict}
Let $\calC$ be a strictly unital $\infty$-category equipped with a
map $\calC \rightarrow \Delta^1$, having fibers $\calC_0$ and $\calC_1$.
For every simplicial subset $K \subseteq \calC_1$, let $\calC[K]$ denote the simplicial subset of $\calC$ spanned those simplices whose intersection with $\calC_1$ belongs to $K$.
Suppose that $K' \subseteq \calC_1$ is obtained from
$K$ by adjoining a single nondegenerate $m$-simplex $\tau: \Delta^m \rightarrow \calC$.
Then the diagram
$$ \xymatrix{ (\calC_0 \times_{\calC} \calC_{/ \tau}) \star \bd \Delta^m \ar[r] \ar[d] & \calC[K] \ar[d] \\
(\calC_{0} \times_{\calC} \calC_{/ \tau}) \star \Delta^m \ar[r] & \calC[K'] }$$
is a pushout square of simplicial sets.
\end{proposition}

\begin{proof}
It is easy to see that we have an injection of simplicial sets
$$ j: ((\calC_{0} \times_{\calC} \calC_{/ \tau}) \star \Delta^m)
\coprod_{  (\calC_0 \times_{\calC} \calC_{/ \tau}) \star \bd \Delta^m } \calC[K] \rightarrow \calC[K'].$$
To show that this map is surjective, consider an arbitrary nondegenerate simplex $\sigma: \Delta^n \rightarrow \calC[K']$ of $\calC[K']$. If $\sigma$ does not belong to $\calC[K]$, then
there exists $0 \leq i \leq n$ such that $\sigma | \Delta^{ \{0, \ldots, i \} }$ factors through
$\calC_0$ and $\sigma| \Delta^{ \{ i+1, \ldots, n \} }$ is a surjection onto the nondegenerate
simplex $\tau$ of $K'$. Since $\sigma$ is nondegenerate and $\calC$ is strictly unital,
no restriction $\sigma | \Delta^{ \{j, j+1\} }$ is degenerate; it follows that
$\sigma$ induces an isomorphism $\Delta^{ \{ i+1, \ldots, n \} } \simeq \Delta^m$, so that
$\sigma$ lies in the image of the map
$(\calC_0 \times_{ \calC} \calC_{/ \tau} ) \star \Delta^m \rightarrow \calC[K']$.
\end{proof}

\subsection{Operadic Left Kan Extensions}\label{iljest}

The theory of colimits can be regarded as a special case of the theory of left Kan extensions:
a functor between $\infty$-categories $\overline{p}: K^{\triangleright} \rightarrow \calC$ is a colimit diagram if and only $\overline{p}$ is a left Kan extension of $p = \overline{p} | K$. The theory of left Kan extensions has a counterpart in the setting of $\infty$-operads, which we will call
{\it operadic left Kan extensions}. Our goal in this section is to introduce this
counterpart (Definition \ref{campe}) and to prove a basic existence result (Theorem \ref{oplk}).
The results of this section will play a crucial role in our discussion of free algebras
in \S \ref{comm33}.

We begin by introducing a bit of terminology.

\begin{definition}\label{sike}
Let $\calC$ be an $\infty$-category. A {\it $\calC$-family of $\infty$-operads} is
a categorical fibration $p: \calO^{\otimes} \rightarrow \calC \times \Nerve(\FinSeg)$ with the following properties:
\begin{itemize}
\item[$(1)$] Let $C \in \calC$ be an object, let
$X \in \calO^{\otimes}_{C}$ have image $\seg{m} \in \FinSeg$, and let
$\alpha: \seg{m} \rightarrow \seg{n}$ be an inert morphism. Then
there exists a $p$-coCartesian morphism $\overline{\alpha}: X \rightarrow Y$
in $\calO^{\otimes}_{C}$.
\end{itemize}
We will say that a morphism $\overline{\alpha}$ of $\calO^{\otimes}$ is {\it inert} if
$\overline{\alpha}$ is $p$-coCartesian, the image of $\overline{\alpha}$ in $\Nerve(\FinSeg)$ is inert, and the image of $\overline{\alpha}$ in $\calC$ is an equivalence.
\begin{itemize}
\item[$(2)$] Let $X \in \calO^{\otimes}$ have images $C \in \calC$ and
$\seg{n} \in \Nerve(\FinSeg)$. For $1 \leq i \leq n$, let 
$f_i: X \rightarrow X_i$ be an inert morphism in $\calO^{\otimes}_{C}$ which
covers $\Colp{i}: \seg{n} \rightarrow \seg{1}$. Then the collection of morphisms
$\{ f_i \}_{1 \leq i \leq n}$ determines a $p$-limit diagram $\nostar{n}^{\triangleleft} \rightarrow \calO^{\otimes}$.
\item[$(3)$] For each object $C \in \calC$, the induced map $\calO^{\otimes}_{C} \rightarrow \Nerve(\FinSeg)$ is an $\infty$-operad.
\end{itemize}

Let ${\calO'}^{\otimes} \rightarrow \Nerve(\FinSeg)$ be an $\infty$-operad and let
$\calO^{\otimes} \rightarrow \calC \times \Nerve(\FinSeg)$ be a $\calC$-family of $\infty$-operads. An {\it $\infty$-operad family map} from $\calO^{\otimes}$ to ${\calO'}^{\otimes}$
is a map $F: \calO^{\otimes} \rightarrow {\calO'}^{\otimes}$ with the following properties:
\begin{itemize}
\item[$(i)$] The diagram
$$ \xymatrix{ \calO^{\otimes} \ar[rr] \ar[dr] & & {\calO'}^{\otimes} \ar[dl] \\
& \Nerve(\FinSeg) & }$$
is commutative.
\item[$(ii)$] The functor $F$ carries inert morphisms in $\calO^{\otimes}$ to inert morphisms in
${\calO'}^{\otimes}$.
\end{itemize}
We let $\Alg_{\calO}(\calO')$ denote the full subcategory
of $\bHom_{ \Nerve(\FinSeg)}( \calO^{\otimes}, {\calO'}^{\otimes})$ spanned by the
$\infty$-operad family maps.
\end{definition}

\begin{remark}
Let $\calO^{\otimes} \rightarrow \calC \times \Nerve(\FinSeg)$ be a $\calC$-family of
$\infty$-operads. Then the underlying map $\calO^{\otimes} \rightarrow \Nerve(\FinSeg)$
exhibits $\calO^{\otimes}$ as an $\infty$-operad family, in the sense of Definition
\ref{suskup}. Conversely, suppose that $p: \calO^{\otimes} \rightarrow \Nerve(\FinSeg)$ is
an $\infty$-operad family, and set $\calC = \calO^{\otimes}_{\seg{0}}$. 
Choose a $p$-coCartesian natural transformation $\id_{ \calO^{\otimes}} \rightarrow q$,
where $q: \calO^{\otimes} \rightarrow \calO^{\otimes}_{\seg{0}} = \calC$, and choose
a factorization of $q \times p$ as a composition
$$ \calO^{\otimes} \stackrel{\theta}{\rightarrow} {\calO'}^{\otimes} \stackrel{\theta'}{\rightarrow} \calC \times \Nerve(\FinSeg)$$
where $\theta$ is a categorical equivalence and $\theta'$ is a categorical fibration. Then
$\theta'$ exhibits ${\calO'}^{\otimes}$ as a $\calC$-family of $\infty$-operads.
In other words, we can regard Definition \ref{sike} as an slight variation on
Definition \ref{suskup}: an $\infty$-operad family is essentially the same thing
as a $\calC$-family of $\infty$-operads, where the underlying $\infty$-category
$\calC$ has not been specified.
\end{remark}

\begin{notation}
Let $\calO^{\otimes} \rightarrow \calC \times \Nerve(\FinSeg)$ be a $\calC$-family of $\infty$-operads.
We let $\calO^{\otimes}_{\acti}$ denote the subcategory of $\calO^{\otimes}$ spanned by those
morphisms which induce active morphisms in $\Nerve(\FinSeg)$.
\end{notation}

\begin{definition}
A {\it correspondence of $\infty$-operads} is a $\Delta^1$-family of $\infty$-operads
$p: \calM^{\otimes} \rightarrow \Nerve(\FinSeg) \times \Delta^1$. In this case, we will say that
$\calM^{\otimes}$ is a {\it correspondence} from the $\infty$-operad
$\calA^{\otimes} = \calM^{\otimes} \times_{ \Delta^1} \{0\}$ to the
$\infty$-operad $\calB^{\otimes} = \calM^{\otimes} \times_{ \Delta^1} \{1\}$. 
\end{definition}

\begin{definition}\label{campe}
Let $\calM^{\otimes} \rightarrow \Nerve(\FinSeg) \times \Delta^1$ be a correspondence
from an $\infty$-operad $\calA^{\otimes}$ to another $\infty$-operad $\calB^{\otimes}$, let
$q: \calC^{\otimes} \rightarrow \calO^{\otimes}$ be a fibration of $\infty$-operads, and let
$\overline{F}: \calM^{\otimes} \rightarrow \calC^{\otimes}$ be an $\infty$-operad family map.
We will say that $\overline{F}$ is an {\it operadic $q$-left Kan extension}
of $F = \overline{F} | \calA^{\otimes}$ if the following condition is satisfied, for every object 
$B \in \calB^{\otimes}$: 

\begin{itemize}
\item[$(\ast)$] Let $K = (\calM^{\otimes}_{\acti})_{/B} \times_{ \calM^{\otimes}} \calA^{\otimes}$.
Then the composite map
$$ K^{\triangleright} \rightarrow (\calM^{\otimes})^{\triangleright}_{/B} \rightarrow
\calM^{\otimes} \stackrel{ \overline{F}}{\rightarrow} \calC^{\otimes}$$
is an operadic $q$-colimit diagram.
\end{itemize}
\end{definition}

The remainder of this section will be devoted to the proof of the following result:

\begin{theorem}\label{oplk}
Let $n \geq 1$, let
$p: \calM^{\otimes} \rightarrow \Nerve(\FinSeg) \times \Delta^n$ be
a $\Delta^n$-family of $\infty$-operads, and let $q: \calC^{\otimes} \rightarrow \calO^{\otimes}$ be a fibration of $\infty$-operads. Suppose we are given a commutative diagram
of $\infty$-operad family maps
$$ \xymatrix{ \calM^{\otimes} \times_{\Delta^n} \Lambda^{n}_0 \ar[d] \ar[r]^-{f_0} & \calC^{\otimes} \ar[d]^{q} \\
\calM^{\otimes} \ar@{-->}[ur]^{f} \ar[r]^{g} & \calO^{\otimes}. }$$
Then
\begin{itemize}
\item[$(A)$] Suppose that $n=1$. The following conditions are equivalent:
\begin{itemize}
\item[$(a)$] There exists a dotted arrow $f$ as indicated in the diagram
which is an operadic $q$-left Kan extension of $f_0$ (in particular, such that $f$ is
a map of $\infty$-operad families).

\item[$(b)$] For every object $B \in \calM^{\otimes} \times_{ \Delta^n } \{1\}$,
the diagram
$$ (\calM^{\otimes}_{\acti})_{/B} \times_{ \Delta^n } \{0\}
\rightarrow \calM^{\otimes} \times_{ \Delta^n} \{0\} \stackrel{f_0}{\rightarrow} \calC^{\otimes}$$
can be extended to an operadic $q$-colimit diagram lifting the map
$$ (( \calM^{\otimes}_{\acti})_{/B} \times_{ \Delta^n } \{0\})^{\triangleright}
\rightarrow \calM^{\otimes} \stackrel{g}{\rightarrow} \calO^{\otimes}.$$
\end{itemize}
\item[$(B)$] Suppose that $n > 1$, and that the restriction of $f_0$ to
$\calM^{\otimes} \times_{ \Delta^n} \Delta^{ \{0,1\} }$ is an operadic $q$-left Kan extension
of $f_0 | ( \calM^{\otimes} \times_{ \Delta^n} \{0\})$. Then there exists a dotted arrow $f$ as indicated in the diagram (automatically a map of $\infty$-operad families).
\end{itemize}
\end{theorem}


The proof of Theorem \ref{oplk} will requires some preliminaries.

\begin{lemma}\label{hookman}
Suppose we are given an inner fibration of simplicial sets $p: \calC \rightarrow \Lambda^3_2$ satisfying the following condition: for every object $X \in \calC_1$, there exists a $p$-coCartesian morphism
$f: X \rightarrow Y$ where $Y \in \calC_2$. Then there exists a pullback diagram
$$ \xymatrix{ \calC \ar[r] \ar[d]^{p} & \calC' \ar[d]^{p'} \\
\Lambda^3_2 \ar[r] & \Delta^3 }$$
where $\calC'$ is an $\infty$-category and the map $\calC \hookrightarrow \calC'$ is a categorical equivalence.
\end{lemma}

\begin{proof}
We will construct a sequence of categorical equivalences
$$ \calC = \calC(0) \subseteq \calC(1) \subseteq \cdots$$
in $(\sSet)_{/ \Delta^3}$ with the following properties: 
\begin{itemize}
\item[$(i)$] If $0 < i < n$ and
$\Lambda^{n}_{i} \rightarrow \calC(m)$ is a map, then the induced map
$\Lambda^{n}_{i} \rightarrow \calC(m+1)$ can be extended to an $n$-simplex of
$\calC(m+1)$.
\item[$(ii)$] For each $m \geq 0$, the inclusion $\calC \subseteq \calC(m) \times_{ \Delta^3} \Lambda^3_2$ is an isomorphism of simplicial sets.
\end{itemize}
Our construction proceeds by induction on $m$. Assume that $\calC(m)$ has been constructed,
and let $A$ denote the collection of all maps $\alpha: \Lambda^n_i \rightarrow \calC(m)$ for
$0 < i < n$. We will prove that, for each $\alpha \in A$, there exists a categorical equivalence
$\calC(m) \subseteq \calC(m,\alpha)$ in $(\sSet)_{/ \Delta^3}$ with the following properties:
\begin{itemize}
\item[$(i')$] The composite map $\Lambda^n_i \stackrel{\alpha}{\rightarrow} \calC(m)
\subseteq \calC(m,\alpha)$ can be extended to an $n$-simplex of $\calC(m,\alpha)$.
\item[$(ii')$] The map $\calC \rightarrow \calC(m,\alpha) \times_{ \Delta^3} \Lambda^3_2$
is an isomorphism of simplicial sets. 
\end{itemize}
Assuming that we can satisfy these conditions, we can complete the proof by defining
$\calC(m+1)$ to be the amalgamation $\coprod_{\alpha \in A} \calC(m,\alpha)$ in the category
$( \sSet)_{\calC(m)/ \, / \Delta^3}$.

Suppose now that $\alpha: \Lambda^n_i \rightarrow \calC(m)$ is given. If the composite map
$\alpha_0: \Lambda^n_i \stackrel{\alpha}{\rightarrow} \calC(m) \rightarrow \Delta^3$ factors through
$\Lambda^3_2$, then we can use the assumption that $p$ is an inner fibration to
extend $\alpha$ to an $n$-simplex of $\calC$. In this case, we can satisfy requirements
$(i')$ and $(ii')$ by setting $\calC(m,\alpha) = \calC(m)$. We may therefore assume without loss
of generality that the image of $\alpha_0$ contains $\Delta^{ \{0,1,3\} } \subseteq \Delta^3$.
In particular, we see that $\alpha_0(0)$ and $\alpha_0(n) = 3$.

We observe that $\alpha_0$ extends uniquely to a map $\overline{\alpha}_0: \Delta^n \rightarrow \Delta^3$. The pushout $\calC(m) \coprod_{ \Lambda^{n}_{i} } \Delta^n$ evidently satisfies
condition $(i')$. It satisfies condition $(ii')$ unless $\overline{\alpha}_0$ carries
$\Delta^{ \{0, \ldots, i-1, i+1, \ldots, n \}}$ into $\Lambda^3_2$. If condition $(ii')$ is satisfied,
we can set $\calC(m,\alpha) = \calC(m) \coprod_{ \Lambda^n_i} \Delta^n$; otherwise,
we have $\overline{\alpha}_0^{-1} \{1\} = \{ i \}$. 

Let $A = \overline{\alpha}_{0}^{-1} \{0\} \subseteq \Delta^n$, 
$B = \overline{\alpha}_0^{-1} \{2\} \subseteq \Delta^n$, and $C = \overline{\alpha}_0^{-1} \{3\}
\subseteq \Delta^n$, so we have a canonical isomorphism
$\Delta^n \simeq A \star \{x\} \star B \star C$ where $x$ corresponds to the $i$th vertex of $\Delta^n$.
Let $\beta$ denote the restriction of $\alpha$ to $A \star \{x\} \star B$ and $\gamma$ the restriction of $\alpha$ to $\{x \} \star B \star C$. 
Let $X \in \calC_1$ denote the image of
$x$ in $\calC_1$, and choose a $p$-coCartesian morphism $f: X \rightarrow Y$
in $\calC$ where $Y \in \calC_2$. Since $f$ is $p$-coCartesian, we can choose a map
$\overline{\gamma}: \{x\} \star \{y\} \star B \star C \rightarrow \calC$ which is compatible
with $f$ and $\gamma$. Let $\overline{\gamma}_0 = \overline{\gamma} | \{x\} \star \{y\} \star B$, and
choose a map $\overline{\beta}: A \star \{x\} \star \{y\} \star B \rightarrow \calC$ compatible with
$\overline{\gamma}_0$ and $\beta$. The restrictions of $\overline{\beta}$ and
$\overline{\gamma}$ determine a map
$\delta: (A \star \{y\} \star B) \coprod_{ \{y\} \star B} (\{y\} \star B \star C) \rightarrow \calC$. Using the fact that $p$ is an inner fibration, we can extend $\delta$ to a map
$\overline{\delta}: A \star \{y\} \star B \star C \rightarrow \calC$. 

We define simplicial subsets
$$ K_0 \subseteq K_1 \subseteq K_2 \subseteq A \star \{x\} \star \{y\} \star B \star C$$
as follows:
\begin{itemize}
\item[$(a)$] The simplicial subset $K_0 \subseteq A \star \{x\} \star \{y\} \star B \star C$ is generated by $A \star \{x \} \star \{y\} \star B$ and
$\{x\} \star \{y\} \star B \star C$. 
\item[$(b)$] The simplicial subset $K_1 \subseteq A \star \{x\} \star \{y\} \star B \star C$ is generated by
$K_0$ together with $\Lambda^n_i \subseteq \Delta^n \simeq A \star \{x\} \star B \star C$.
Since the inclusions $\Lambda^n_i \subseteq \Delta^n$ and
$(A \star \{x\} \star B) \coprod_{ \{x\} \star B} (\{x\} \star B \star C) \subseteq
A \star \{x\} \star B \star C$ are inner anodyne, the inclusion
$$i: (A \star \{x\} \star B) \coprod_{ \{x\} \star B} (\{x\} \star B \star C) \subseteq \Lambda^n_i$$ is a categorical equivalence.
The inclusion $K_0 \subseteq K_1$ is a pushout of $i$, and therefore also a categorical equivalence.

\item[$(c)$] The simplicial subset $K_2 \subseteq A \star \{x\} \star \{y\} \star B \star C$ is generated
by $K_1$ together with $A \star \{y\} \star B \star C$. The inclusion $K_1 \subseteq K_2$ is a pushout of the
inclusion $(A \star \{y\} \star B) \coprod_{ \{y\} \star B} (\{y\} \star B \star C) \subseteq 
A \star \{y\} \star B \star C$, and therefore a categorical equivalence.
\end{itemize}
The inclusion $K_0 \subseteq  A \star \{x\} \star \{y\} \star B \star C$ is evidently a categorical equivalence. It follows from a two-out-of-three argument that the inclusion
$K_2 \subseteq A \star \{x\} \star \{y\} \star B \star C$ is also a categorical equivalence.
The maps $\overline{\delta}$, $\overline{\beta}$, $\overline{\beta}'$, and $\alpha$
determine a map $K_2 \rightarrow \calC(m)$. We now define $\calC(m,\alpha)$ to be the pushout
$\calC(m) \coprod_{K_2} ( A \star \{x\} \star \{y\} \star B \star C)$. It is not difficult to verify that  
$\calC(m,\alpha)$ has the desired properties.
\end{proof}

\begin{lemma}\label{skope}
Let $p: \calC \rightarrow \Delta^n$ be a map of $\infty$-categories, let
$0 < i < n$, and assume that for every object $X \in \calC_{i-1}$ there
exists a $p$-coCartesian morphism $f: X \rightarrow Y$, where $Y \in \calC_{i}$. 
Then the inclusion $p^{-1} \Lambda^{n}_{i} \subseteq \calC$ is a categorical equivalence.
\end{lemma}

\begin{proof}
We first treat the special case $i=1$. The proof proceeds by induction on $n$.
Let $S$ be the collection of all nondegenerate
simplices in $\Delta^n$ which contain the vertices $0$, $1$, and at least one other vertex.
For each $\sigma \in S$, let $\sigma'$ be the simplex obtained from $\sigma$ by removing
the vertex $1$. Choose an ordering $S = \{ \sigma_1, \ldots, \sigma_m \}$ of $S$ where
the dimensions of the simplices $\sigma_j$ are nonstrictly decreasing as a function of $i$
(so that $\sigma_1 = \Delta^n$). For $0 \leq j \leq m$, let 
$K_j$ denote the simplicial subset of $\Delta^n$ obtained by removing
the simplices $\sigma_{k}$ and $\sigma'_{k}$ for $k \leq j$. If we let
$n_{j}$ denote the dimension of $\sigma_{j}$, then we have a pushout diagram
$$\xymatrix{ \Lambda^{n_j}_{1} \ar[r] \ar[d] & K_{j} \ar[d] \\
\Delta^{n_j} \ar[r] & K_{j-1} }$$
Applying the inductive hypothesis (and the left properness of the Joyal model structure), we
deduce that the inclusion $K_{j} \times_{\Delta^n} \calC \rightarrow K_{j-1} \times_{\Delta^n} \calC$
is a categorical equivalence for $1 < j \leq m$. Combining these facts, we deduce that the map
$K_{m} \times_{ \Delta^n} \calC \subseteq K_{1} \times_{\Delta^n} \calC \simeq p^{-1} \Lambda^{n}_1$
is a categorical equivalence. By a two-out-of-three argument, we are reduced to proving that
the inclusion $K_m \times_{\Delta^n} \calC \rightarrow \calC$ is a categorical equivalence.
Let $q: \Delta^{n} \rightarrow \Delta^2$ be the map given on vertices by the formula
$$ q(k) = \begin{cases} 0 & \text{if } k=0 \\
1 & \text{if } k = 1  \\
2 & \text{otherwise} \end{cases}$$
and observe that $K_m = q^{-1} \Lambda^2_1$. We may therefore replace 
$p$ by $q \circ p$ and thereby reduce to the case $n=2$.

Choose a map $h: \Delta^1 \times \calC_0 \rightarrow \calC$ which is a natural
transformation from the identity map $\id_{\calC_0}$ to a functor $F: \calC_0 \rightarrow \calC_1$,
such that $h$ carries $\Delta^1 \times \{X\}$ to a $p$-coCartesian morphism in $\calC$ for
each $X \in \calC_0$. Let $\calD = p^{-1} \Delta{ \{1,2\} }$.
The natural transformation $h$ induces maps
$$ (\Delta^1 \times \calC_0) \coprod_{ \{1\} \times \calC_0} \calC_1 \rightarrow p^{-1} \Delta^{ \{0,1\} }$$
$$ ( \Delta^1 \times \calC_0) \coprod_{ \{1\} \times \calC_0} \calD \rightarrow \calC$$
and it follows from Proposition \toposref{simplexplay} that these maps are categorical equivalences.
Consider the diagram
$$ \xymatrix{ & (\Delta^1 \times \calC_0) \coprod_{ \{1\} \times \calC_0} \calD \ar[dl] \ar[dr] & \\
p^{-1} \Lambda^2_1 \ar[rr] & & \calC. }$$
It follows from the above arguments (and the left properness of the Joyal model structure) that
the diagonal maps are categorical equivalences, so that the horizontal map is a categorical equivalence by the two-out-of-three property. This completes the proof in the case $i=1$.

We now treat the case $i > 1$. The proof again proceeds by induction on $n$. Let $q: \Delta^n \rightarrow \Delta^3$ be the map defined by the formula
$$ q(k) = \begin{cases} 0 & \text{if } k < i-1 \\
1 & \text{if } k = i-1  \\
2 & \text{if } k = i \\
3 & \text{otherwise.}\end{cases}$$
Let $S$ denote the collection of all nondegenerate simplices $\sigma$ of $\Delta^n$ such that
the restriction $q| \sigma$ is surjective. For each $\sigma \in S$, let $\sigma'$ denote the simplex
obtained from $\sigma$ by deleting the vertex $i$. Choosing an ordering
$S = \{ \sigma_1, \ldots, \sigma_m \}$ of $S$ where the dimension of the simplex
$\sigma_{j}$ is a nondecreasing function of $j$ (so that $\sigma_1 = \Delta^n$), and for
$0 \leq j \leq m$ let $K_{j}$ be the simplicial subset of $\Delta^n$ obtained by deleting 
$\sigma_{k}$ and $\sigma'_{k}$ for $k \leq j$. If we let $n_{j}$ denote the dimension of $\sigma_{j}$, then
we have pushout diagrams
$$\xymatrix{ \Lambda^{n_j}_{p} \ar[r] \ar[d] & K_{j} \ar[d] \\
\Delta^{n_j} \ar[r] & K_{j-1} }$$
where $1 < p < n_j$. Applying the inductive hypothesis and the left properness of the Joyal model structure, we deduce that $K_{j} \times_{ \Delta^n} \calC \rightarrow K_{j-1} \times_{ \Delta^n} \calC$
is a categorical equivalence for $1 < j \leq m$. It follows that the map
$$K_{m} \times_{ \Delta^n} \calC \rightarrow K_{1} \times_{ \Delta^n} \calC = p^{-1} \Lambda^n_i$$
is a categorical equivalence. To complete the proof, it will suffice to show that the inclusion
$K_{m} \times_{ \Delta^n} \calC \rightarrow \calC$ is a categorical equivalence. We observe
that $K_{m} = q^{-1} \Lambda^{3}_2$. We may therefore replace $p$ by $q \circ p$ and reduce to the
case where $n=3$ and $i=2$.

Applying Lemma \ref{hookman}, we can factor the map $p^{-1} \Lambda^3_2 \rightarrow \Delta^3$
as a composition
$$ p^{-1} \Lambda^3_{2} \stackrel{i}{\rightarrow} \calC' \stackrel{p'}{\rightarrow} \Delta^3,$$
where $\calC'$ is an $\infty$-category and $i$ is a categorical equivalence which induces
an isomorphism $p^{-1} \Lambda^3_{2} \simeq {p'}^{-1} \Lambda^3_{2}$. In particular, 
$i$ is a trivial cofibration with respect to the Joyal model structure, so there exists a solution to
the following lifting property:
$$ \xymatrix{ p^{-1} \Lambda^3_{2} \ar[d]^{i} \ar[r]^{j} & \calC \ar[d] \\
\calC' \ar[r] \ar@{-->}[ur]^{g} & \Delta^3. }$$
Since the map $g$ induces an isomorphism $\Lambda^3_{2} \times_{ \Delta^3} \calC'
\rightarrow \Lambda^{3}_{2} \times_{ \Delta^3} \calC$, it is a categorical equivalence
(it is bijective on vertices and induces isomorphisms $\Hom^{\rght}_{\calC'}(x,y)
\rightarrow \Hom^{\rght}_{\calC}(g(x), g(y) )$ for every pair of vertices $x,y \in \calC'$,
since $\Lambda^3_{2}$ contains every edge of $\Delta^3$). It follows that
$j = g \circ i$ is a categorical equivalence as well, which completes the proof.
\end{proof}

\begin{proof}[Proof of Theorem \ref{oplk}]
The implication $(a) \Rightarrow (b)$ in the case $(A)$ is obvious. Let us therefore assume
either that $(b)$ is satisfied (if $n=1$) or that $f_0 | (\calM^{\otimes} \times_{ \Delta^n} \Delta^{ \{0,1\} }$ is an operadic $q$-left Kan extension of $f_0 | (\calM^{\otimes} \times_{ \Delta^n} \{0\})$ (if $n > 1$). To complete the proof, we need to construct a functor $f: \calM^{\otimes} \rightarrow \calC^{\otimes}$ satisfying some natural conditions. The construction is somewhat elaborate, and will require several steps.


Using Proposition \toposref{minimod}, we can choose a simplicial subset 
$\calN^{\otimes} \subseteq \calM^{\otimes}$ such that $\calN^{\otimes}$ is categorically
equivalent to $\calM^{\otimes}$, the induced map $\calN^{\otimes} \rightarrow \Nerve(\FinSeg) \times \Delta^n$ is a minimal inner fibration, and there is a retraction of $\calM^{\otimes}$ onto $\calN^{\otimes}$ which is compatible with the projection to $\Nerve(\FinSeg) \times \Delta^n$. 
Consider the associated lifting problem 
$$ \xymatrix{ {\calN}^{\otimes} \times_{ \Delta^n} \Lambda^n_0 \ar[r] \ar[d] & \calC^{\otimes} \ar[d] \\
{\calN}^{\otimes} \ar[r] \ar@{-->}[ur]^{f'} & \calO^{\otimes}. }$$
Suppose this lifting problem admits a solution (where $f'$ is an operadic left Kan extension
of $f' | \calN^{\otimes} \times_{ \Delta^n } \{0\}$ if $n=1$). Applying Proposition \toposref{princex}, we deduce that the original lifting problem admits a solution $f$ such that $f| \calN^{\otimes}$ is equivalent to $f'$. If $n=1$, it will then follow that $f$ is also an operadic $q$-left Kan extension of
$f_0$. Consequently, we are free to replace $\calM^{\otimes}$ by $\calN^{\otimes}$ and thereby reduce to the case where $\calM^{\otimes} \rightarrow \Nerve(\FinSeg) \times \Delta^n$ is a minimal inner fibration. Since $\Nerve(\FinSeg) \times \Delta^n$ is a minimal $\infty$-category
(this follows from the observation that every isomorphism class in the category $\FinSeg \times [n]$
has a unique representative), we deduce that $\calM^{\otimes}$ is a strictly unital $\infty$-category
(see \S \ref{strictun}).

Our next step is to decompose the simplices of $\Nerve(\FinSeg) \times \Delta^{ \{1, \ldots, n\} }$
into groups. We will say that a morphism $\alpha$ in $\Nerve(\FinSeg) \times \Delta^{ \{1, \ldots, n\} }$
is {\it active} if its image in $\Nerve(\FinSeg)$ is active, {\it inert} if its image in
$\Nerve(\FinSeg)$ is inert and its image in $\Delta^{ \{1, \ldots, n\} }$ is degenerate, and
{\it neutral} if it is neither active nor inert. 
Let $\sigma$ be an $m$-simplex of $\Nerve(\FinSeg) \times \Delta^{ \{1, \ldots, n\} }$, corresponding to a chain of morphisms $$ (\seg{k_0}, e_0) \stackrel{\alpha(1)}{\rightarrow} (\seg{k_1}, e_1) \stackrel{\alpha(2)}{\rightarrow}
\cdots \stackrel{\alpha(m)}{\rightarrow} (\seg{k_m}, e_m).$$
Let $J_{\sigma}$ denote the collection of integers $j \in \{1, \ldots, m\}$ for which
the map $\alpha(j)$ is not an isomorphism. We will denote the cardinality of
$J_{\sigma}$ by $l(\sigma)$ and refer to it as the {\it length} of $\sigma$. 
For $1 \leq d \leq l(\sigma)$, we let $j^{\sigma}_{d}$ denote the $d$th element of
$J_{\sigma}$ and set $\alpha^{\sigma}_{d} = \alpha( j^{\sigma}_{d})$.
We will say that $\sigma$ is {\it closed} if $k_{m} = 1$; otherwise we will say that
$\sigma$ is {\it open}. We will say that $\sigma$ is {\it new} if the projection map
$\sigma \rightarrow \Delta^{ \{1, \ldots, n\} }$ is surjective.

We now partition the collection of new simplices $\sigma$ of $\Nerve(\FinSeg) \times
\Delta^{ \{1, \ldots, n\} }$ into eleven groups:

\begin{itemize}
\item[$(G'_{(1)})$] A new simplex $\sigma$ of $\Nerve(\FinSeg) \times \Delta^{ \{1, \ldots, n\} }$ belongs to $G'_{(1)}$ if
it is a closed and the maps $\alpha^{\sigma}_{i}$ are active for $1 \leq i \leq l(\sigma)$.

\item[$(G_{(2)})$] A new simplex $\sigma$ of $\Nerve(\FinSeg) \times \Delta^{ \{1, \ldots, n\} }$ belongs
to $G_{(2)}$ if $\sigma$ is closed and there exists
$1 \leq k < l(\sigma)$ such that $\alpha^{\sigma}_{k}$ is inert, while
$\alpha^{\sigma}_{j}$ is active for $k < j \leq l(\sigma)$.

\item[$(G'_{(2)})$]  A new simplex $\sigma$ of $\Nerve(\FinSeg) \times \Delta^{ \{1, \ldots, n\} }$ belongs
to $G'_{(2)}$ if $\sigma$ is closed and there exists $1 \leq k \leq l(\sigma)$ such that
$\alpha^{\sigma}_{k}$ is neutral while the maps $\alpha^{\sigma}_{j}$ are active for $k < j \leq l(\sigma)$.

\item[$(G_{(3)})$] A new simplex $\sigma$ of $\Nerve(\FinSeg) \times \Delta^{ \{1, \ldots, n\} }$ belongs
to $G_{(3)}$ if $\sigma$ is closed and there exists $1 \leq k < l(\sigma)-1$ such that
$\alpha^{\sigma}_{k}$ is inert, the maps $\alpha^{\sigma}_{j}$ are active for
$k < j < l(\sigma)$, and $\alpha^{\sigma}_{l(\sigma)}$ is inert.

\item[$(G'_{(3)})$] A new simplex $\sigma$ of $\Nerve(\FinSeg) \times \Delta^{ \{1, \ldots, n\} }$ belongs
to $G'_{(3)}$ if $\sigma$ is closed and there exists $1 \leq k < l(\sigma)$ such that the map $\alpha^{\sigma}_{k}$ is neutral, the maps $\alpha^{\sigma}_{j}$ are active for $k < j < l(\sigma)$, and
$\alpha^{\sigma}_{l(\sigma)}$ is inert. 

\item[$(G_{(4)})$] A new simplex $\sigma$ of $\Nerve(\FinSeg) \times \Delta^{ \{1, \ldots, n\} }$ belongs to $G_{(4)}$ if
it is a closed, the maps $\alpha^{\sigma}_{i}$ are active for $1 \leq i < l(\sigma)$, and
the map $\alpha^{\sigma}_{l(\sigma)}$ is inert.

\item[$(G'_{(4)})$] A new simplex $\sigma$ of $\Nerve(\FinSeg) \times \Delta^{ \{1, \ldots, n\} }$ belongs to $G'_{(4)}$ if
it is an open and the maps $\alpha^{\sigma}_{i}$ are active for $1 \leq i \leq l(\sigma)$.

\item[$(G_{(5)})$] A new simplex $\sigma$ of $\Nerve(\FinSeg) \times \Delta^{ \{1, \ldots, n\} }$ belongs
to $G_{(5)}$ if $\sigma$ is an open and there exists $1 \leq k < l(\sigma)$ such that
$\alpha^{\sigma}_{k}$ is inert and $\alpha^{\sigma}_{j}$ is active for
$k < j \leq l(\sigma)$. 

\item[$(G'_{(5)})$] A new simplex $\sigma$ of $\Nerve(\FinSeg) \times \Delta^{ \{1, \ldots, n\} }$ belongs
to $G'_{(5)}$ if $\sigma$ is an open and there exists $1 \leq k \leq l(\sigma)$ such that
$\alpha^{\sigma}_{k}$
is neutral and the maps $\alpha^{\sigma}_{j}$ are active for $k < j \leq l(\sigma)$.

\item[$(G_{(6)})$] A new simplex $\sigma$ of $\Nerve(\FinSeg) \times \Delta^{ \{1, \ldots, n\} }$ belongs to $G_{(6)}$ if
it is closed, has length $\geq 2$, and the maps $\alpha^{\sigma}_{l(\sigma)-1}$ and
$\alpha^{\sigma}_{l(\sigma)}$ are both inert.

\item[$(G'_{(6)})$] A new simplex $\sigma$ of $\Nerve(\FinSeg) \times \Delta^{ \{1, \ldots, n\} }$ belongs to $G'_{(6)}$ if
it is open, has length at least $1$, and the map $\alpha^{\sigma}_{l(\sigma)}$ is inert.
\end{itemize}

For each integer $m \geq 0$, let $F(m)$ denote simplicial subset of
$\Nerve(\FinSeg) \times \Delta^{ \{1, \ldots, n\} }$ spanned by the simplices which
either belong to $\Nerve(\FinSeg) \times \bd \Delta^{ \{1, \ldots, n\} }$, have
length less than $m$, or have length $m$ and belong to one of the groups $G_{(i)}$
for $2 \leq i \leq 6$. Let $\overline{F}(m)$ be the
simplicial subset of $\Nerve(\FinSeg) \times \Delta^n$ spanned by those simplices
whose intersection with $\Nerve(\FinSeg) \times \Delta^{ \{1, \ldots, n\} }$ belongs to $F(m)$, and let
$\calM^{\otimes}_{(m)}$ denote the inverse image of $\overline{F}(m)$ in
$\calM^{\otimes}$. We observe that $\calM^{\otimes}_{(0)} = \calM^{\otimes}
\times_{\Delta^n} \Lambda^n_0$. To complete the proof, it will suffice to show that
$f_0: \calM^{\otimes}_{(0)} \rightarrow \calC^{\otimes}$ can be extended to a sequence of maps
$f_m: \calM^{\otimes}_{(m)} \rightarrow \calC^{\otimes}$ which fit into commutative diagrams
$$ \xymatrix{ \calM^{\otimes}_{(m-1)} \ar[d] \ar[r]^{f_{m-1}} & \calC^{\otimes} \ar[d]^{q} \\
\calM^{\otimes}_{(m)} \ar@{-->}[ur]^{f_m} \ar[r] & \calO^{\otimes} }$$
and such that $f_1$ has the following special properties:
\begin{itemize}
\item[$(i)$] For each object $B \in \calM \times_{ \Delta^n} \{1\}$, the map
$$ ((\calM^{\otimes}_{\acti})_{/B} \times_{ \Delta^n} \{0\})^{\triangleright}
\rightarrow \calM^{\otimes}_{(1)} \stackrel{f_1}{\rightarrow} \calC^{\otimes}$$
is an operadic $q$-colimit diagram.
\item[$(ii)$] For every inert morphism $e: M' \rightarrow M$ in $\calM^{\otimes} \times_{ \Delta^n} \{1\}$ such that $M \in \calM = \calM^{\otimes}_{\seg{1}}$, the functor
$f_1$ carries $e$ to an inert morphism in $\calC^{\otimes}$.
\end{itemize}

Fix $m > 0$, and suppose that $f_{m-1}$ has already been constructed.
We define a filtration
$$ F(m-1) = K(0) \subseteq K(1) \subseteq K(2) \subseteq K(3) \subseteq K(4) \subseteq K(5) \subseteq K(6) = F(m)$$
as follows:
\begin{itemize}
\item We let $K(1)$ denote the simplicial subset of $\Nerve(\FinSeg) \times \Delta^{ \{1, \ldots, n\} }$ spanned by those simplices which either belong to $K(0)$ or have length $(m-1)$ and belong
to $G'_{(1)}$. 
\item For $2 \leq i \leq 6$, we let $K(i)$ be the simplicial subset of $\Nerve(\FinSeg) \times \Delta^{ \{1, \ldots, n\} }$ spanned by those simplices which either belong to $K(i-1)$, have length $m$ and belong to $G_{(i)}$, or have length $m-1$ and belong to $G'_{(i)}$. 
\end{itemize}
For $0 \leq i \leq 6$, we let $\overline{K}(i)$ denote the simplicial subset of
$\Nerve(\FinSeg) \times \Delta^n$ spanned by those simplices whose intersection
with $\Nerve(\FinSeg) \times \Delta^{ \{1, \ldots, n\} }$ belongs to $K(i)$, and $\calM^{\otimes}_{i}$ the
inert image of $\overline{K}(i)$ in $\calM^{\otimes}$. We will define maps
$f^{i}: \calM^{\otimes}_{i} \rightarrow \calC^{\otimes}$ with $f^{0} = f_{m-1}$.
The construction now proceeds in six steps:

\begin{itemize}
\item[$(1)$] Assume that $f^{0} = f_{m-1}$ has been constructed; we wish to define
$f^{1}$. Let $\{ \sigma_{a} \}_{a \in A}$ be the collection of all nondegenerate
simplices of $\calM^{\otimes} \times_{ \Delta^n} \Delta^{ \{1, \ldots, n\} }$ whose image in $\Nerve(\FinSeg) \times \Delta^{ \{1, \ldots, n\} }$ has length $(m-1)$ and belongs to the class $G'_{(1)}$.
Choose a well-ordering of the set $A$ such that the dimensions of the simplices
$\sigma_{a}$ form a (nonstrictly) increasing function of $a$. For each $a \in A$, let $\calM^{\otimes}_{\leq a}$ denote the
simplicial subset of $\calM^{\otimes}$ spanned by those simplices which
either belong to $\calM^{\otimes}_{0}$ or whose intersection with
$\calM^{\otimes} \times_{ \Delta^n} \Delta^{ \{1, \ldots, n\} }$ is contained in $\sigma_{a'}$ for some $a' \leq a$, and define $\calM^{\otimes}_{< a}$ similarly. We construct a compatible family of maps
$f^{\leq a}: \calM^{\otimes}_{\leq a} \rightarrow \calC^{\otimes}$ extending $f^{0}$, using transfinite
induction on $a$. Assume that $f^{\leq a'}$ has been constructed for $a' < a$; these
maps can be amalgamated to obtain a map $f^{<a}: \calM^{\otimes}_{< a} \rightarrow \calC^{\otimes}$. 
Let $X = \{0\} \times_{ \Delta^n } \calM^{\otimes}_{/ \sigma_{a}}$. Since
$\calM^{\otimes}$ is strictly unital, Proposition \ref{appstrict} guarantees that we have
a pushout diagram of simplicial sets
$$ \xymatrix{ X \star \bd \sigma_a \ar[r] \ar[d] & \calM^{\otimes}_{< a} \ar[d] \\
X \star \sigma_{a} \ar[r] & \calM^{\otimes}_{\leq a}. }$$
It will therefore suffice to extend the composition
$$ g_0: X \star \bd \sigma_{a} \rightarrow \calM^{\otimes}_{< a} 
\stackrel{f^{<a}}{\rightarrow} \calC^{\otimes}$$
to a map $g: X \star \star \sigma_{a} \rightarrow \calC^{\otimes}$ which is compatible
with the projection to $\calO^{\otimes}$.

We first treat the special case where the simplex $\sigma_{a}$ is zero-dimensional
(in which case we must have $m=1$). We can identify $\sigma_{a}$ with
an object $B \in \calM^{\otimes}$.
Let $X_0 = (\calM^{\otimes}_{\acti})_{/B} \times_{ \Delta^n} \{0\} \subseteq X$.
Using assumption $(b)$, we can choose a map $g_{1}: X_0^{\triangleright} \rightarrow \calC^{\otimes}$
compatible with the projection to $\calO^{\otimes}$ which is an operadic $q$-colimit diagram.
We note that $X_0$ is a localization of
$X$ so that the inclusion $X_0 \subseteq X$ is cofinal; it follows that
$g_{1}$ can be extended to a map $X^{\triangleright} \rightarrow \calC^{\otimes}$
with the desired properties. Note that our particular construction of $g$ guarantees
that the map $f_{1}$ will satisfy condition $(i)$ above for $M \in \calM \times_{ K \star \Delta^n} \{0\}$.

Now suppose that $\sigma_{a}$ is a simplex of positive dimension. We again let
$X_0$ denote the simplicial subset of $X$ spanned by those vertices of
$X$ which correspond to diagrams $\sigma_{a}^{\triangleleft} \rightarrow
\calC^{\otimes}$ which project to a sequence of active morphisms in
$\Nerve(\FinSeg)$. The inclusion $X_0 \subseteq X$ admits a left adjoint and is therefore
cofinal; it follows that the induced map $\calC^{\otimes}_{(f_{0} | X)/} \rightarrow
\calC^{\otimes}_{(f_0|X_0)/} \times_{ \calO^{\otimes}_{ (qf_0|X_0)/}}
\calO^{\otimes}_{(q f_0|X)/}$ is a trivial Kan fibration. It therefore suffices to show that
the restriction $g'_0 = g_0 | (X_0 \star \bd \sigma_{a})$ can be extended to a map
$g': X_0 \star \sigma_{a} \rightarrow \calC^{\otimes}$ compatible with the projection
to $q$. In view of Proposition \ref{quadly}, it will suffice to show that the restriction
$g'_0 | (X_0 \star \{B\})$ is a $q$-operadic colimit diagram.
Since the inclusion $\{B\} \subseteq \sigma_{a}$ is left anodyne, the projection
map $X_0 \rightarrow (\calA^{\otimes}_{\acti})_{/B}$ is a trivial Kan fibration.
It will therefore suffice to show that $f_{m-1}$ induces an operadic
$q$-colimit diagram $\delta: (\calA^{\otimes}_{\acti})^{\triangleright}_{/B} \rightarrow
\calC^{\otimes}$. 

Let $\seg{p}$ be the image of $B$ in $\Nerve(\FinSeg)$. For $1 \leq i \leq p$, choose
an inert morphism $B \rightarrow B_{i}$ in $\calB^{\otimes}$ covering the map
$\Colp{i}: \seg{p} \rightarrow \seg{1}$, and let
$\delta_{i}: (\calA^{\otimes}_{\acti})^{\triangleright}_{/B_i} \rightarrow \calC^{\otimes}$
be the map induced by $f_{n-1}$. The above construction guarantees that
each $\delta_{i}$ is an operadic $q$-colimit diagram. Using the fact that
$\calM^{\otimes}$ is a $\Delta^n$-family of $\infty$-operads, we deduce that
the maps $B \rightarrow B_i$ induce an equivalence
$$ ( \calA^{\otimes}_{\acti})_{/B} \simeq \prod_{1 \leq i \leq p} ( \calA^{\otimes}_{\acti})_{/B_i}.$$
Using $(ii)$, we see that under this equivalence we can identify $\delta$ with the diagram
$$ ( \prod_{1 \leq i \leq p} ( \calA^{\otimes}_{\acti})_{/B_i}^{\triangleright}
\rightarrow \prod_{1 \leq i \leq p} ( \calA^{\otimes}_{\acti})_{/B_i}^{\triangleright}
\stackrel{ \oplus \delta_i }{\rightarrow} \calC^{\otimes}.$$ 
The desired result now follows from Proposition \ref{swingup}.

\item[$(2)$] We now assume that $f^{1}$ has been constructed. 
Since $q$ is a categorical fibration, to produce the desired extension
$f^{2}$ of $f^{1}$ it is sufficient to show that the inclusion
$\calM^{\otimes}_{1} \subseteq \calM^{\otimes}_{2}$ is a categorical equivalence.
For each simplex $\sigma$ of $\calM^{\otimes} \times_{ \Delta^n} \Delta^{ \{1, \ldots, n\} }$
which belongs to $G_{(2)}$, let $k(\sigma) < l(\sigma)$ be the integer such that
$\alpha^{\sigma}_{k(\sigma)}$ is inert while
$\alpha^{\sigma}_{j}$ is active for $k(\sigma) < j \leq l(\sigma)$. We will say that
$\sigma$ is {\it good} if $\alpha^{\sigma}_{k(\sigma)}$ induces a map
$\seg{p} \rightarrow \seg{p'}$ in $\FinSeg$ whose restriction to $(\alpha^{\sigma}_{k})^{-1} \nostar{p'}$
is order preserving. Let $\{ \sigma_{a} \}_{a \in A}$ be the collection of all nondegenerate simplices
of $\Nerve(\FinSeg) \times \Delta^n$ such that the intersection
$\sigma'_{a} = \sigma_{a} \cap ( \Nerve(\FinSeg) \times \Delta^{ \{1, \ldots, n \} })$ is a
(nonempty) good simplex belonging to
$G_{(2)}$. For each $a \in A$, let $k_{a} = k(\sigma'_{a})$ 
and $j_{a} = j^{\sigma'_{a}}_{k}$. Choose a well-ordering of $A$ with the following
properties:
\begin{itemize}
\item The map $a \mapsto k_{a}$ is a (nonstrictly) increasing function of $a \in A$.
\item For each integer $k$, the dimension of the simplex $\sigma_{a}$ is a nonstrictly
increasing function of $a \in A_{k} = \{ a \in A: k_a = k \}$. 
\item Fix integers $k,d \geq 0$, and let $A_{k,d}$
be the collection of elements $a \in A_{k}$ such that $\sigma_{a}$ has dimension $d$.
The map $a \mapsto j_{a}$ is a nonstrictly increasing function on
$A_{k,d}$.
\end{itemize}
For each $a \in A$, let $\overline{K}_{\leq a}$ denote the simplicial subset of
$\Nerve(\FinSeg) \times \Delta^n$ generated by $\overline{K}(1)$ and
the simplices $\{ \sigma_{a'} \}_{a' \leq a}$, and define $\overline{K}_{< a}$ similarly.
The inclusion $\calM^{\otimes}_{1} \subseteq \calM^{\otimes}_{2}$ can be
obtained as a transfinite composition of inclusions
$$ i_{a}: \calM^{\otimes} \times_{ \Nerve(\FinSeg) \times \Delta^n}
\overline{K}_{< a} \rightarrow
\calM^{\otimes} \times_{ \Nerve(\FinSeg) \times \Delta^n} \overline{K}_{\leq a}.$$
Each $i_{a}$ is a pushout of an inclusion
$$ i'_{a}: \calM^{\otimes} \times_{ \Nerve(\FinSeg) \times \Delta^n} \sigma^{0}_{a}
\calM^{\otimes} \times_{ \Nerve(\FinSeg) \times \Delta^n} \sigma_{a},$$
where $\sigma^{0}_{a} \subseteq \sigma_{a}$ denotes the inner horn
obtained by removing the interior of $\sigma_{a}$ together with the
face opposite the $j_{a}$th vertex of $\sigma'_{a}$. It follows from
Lemma \ref{skope} that $i'_{a}$ is a categorical equivalence. Since
the Joyal model structure is left proper, we deduce that each $i_{a}$ is an
equivalence, so that the inclusion $\calM^{\otimes}_{1} \subseteq \calM^{\otimes}_{2}$
is a categorical equivalence as desired.

\item[$(3)$] To find the desired extension $f^{3}$ of $f^{2}$, it suffices to show that
the inclusion $\calM^{\otimes}_{2} \subseteq \calM^{\otimes}_{3}$ is a categorical equivalence. This follows from the argument given in step $(2)$.

\item[$(4)$] Let $\{ \sigma_{a} \}_{a \in A}$ denote the collection of all nondegenerate simplices
$\sigma$ of $\calM^{\otimes}$ with the property that $\sigma \times_{ \Delta^n} \Delta^{ \{1, \ldots, n\}}$
projects to a simplex of $\Nerve(\FinSeg) \times \Delta^{ \{1, \ldots, n \} }$ which belongs to $G'_{4}$. Choose a well-ordering
of $A$ having the property that the dimensions of the simplices $\sigma_{a}$ form a (nonstrictly)
increasing function of $a$. Let $\calD^{\otimes} = \calM^{\otimes} \times_{ \Delta^n} \{n\}$. For 
each $a \in A$ let $D_{a} \in \calD^{\otimes}$ denote the final vertex of $\sigma_{a}$
and let $X_{a}$ denote the full subcategory of $\calD \times_{ \calM^{\otimes} } \calM^{\otimes}_{\sigma_a/}$ spanned by those objects for which the underlying map
$D_{a} \rightarrow D$ is inert. We have a canonical map
$t_{a}: \sigma_{a} \star X_{a} \rightarrow \calM^{\otimes}$. For each
$a \in A$, let $\calM^{\otimes}_{\leq a} \subseteq \calM^{\otimes}$ denote the simplicial subset generated by $\calM^{\otimes}_{3}$ together with the image of $t_{a'}$ for all $a' \leq a$, and define
$\calM^{\otimes}_{< a}$ similarly. Using our assumption that $\calM^{\otimes}$ is strictly unital
(and a variation on Proposition \ref{appstrict}), we deduce that $\calM^{\otimes}_{4} = \bigcup_{a \in A} \calM^{\otimes}_{n,\leq a}$ and that for each $a \in A$ we have a pushout diagram of simplicial sets
$$ \xymatrix{ \bd \sigma_{a} \star X_{a} \ar[r] \ar[d] & \calM^{\otimes}_{< a} \ar[d] \\
\sigma_{a} \star X_{a} \ar[r] & \calM^{\otimes}_{\leq a}. }$$
To construct $f^{4}$, we are reduced to the problem of solving a sequence of extension
problems of the form
$$ \xymatrix{ \bd \sigma_{a} \star X_{a} \ar[r]^{g_0} \ar[d] & \calC^{\otimes} \ar[d]^{q} \\
\sigma_{a} \star X_{a} \ar[r] \ar@{-->}[ur]^{g} & \calO^{\otimes}. }$$
Let $\seg{p}$ denote the image of $D_{a}$ in $\Nerve(\FinSeg)$; note that we have a canonical decomposition $$X_{a} \simeq \coprod_{1 \leq i \leq p} X_{a,i}.$$
Each of the $\infty$-categories $X_{a,i}$ has an initial object $\overline{B}_i$, given by any
map $\sigma_{a} \star \{ D_i \} \rightarrow \calC^{\otimes}$ which induces an inert
morphism $D_{a} \rightarrow D_{i}$ covering $\Colp{i}: \seg{p} \rightarrow \seg{1}$.
Let $h: X_{a} \rightarrow \calC^{\otimes}$ be the map induced by $g_0$, and let
$h'$ be the restriction of $h$ to the discrete simplicial set $X'_{a} = \{ \overline{B}_i \}_{1 \leq i \leq p}$.
Since inclusion $X'_{a} \rightarrow X_a$ is left anodyne, we have
a trivial Kan fibration $\calC^{\otimes}_{/h} \rightarrow \calC^{\otimes}_{/h'}
\times_{ \calO^{\otimes}_{/qh'}} \calO^{\otimes}_{/qh}$. Unwinding the definitions, we are reduced to the lifting probelm depicted in the diagram
$$ \xymatrix{ \bd \sigma_{a} \star X'_{a} \ar[r]^{g'_0} \ar[d] & \calC^{\otimes} \ar[d]^{q} \\
\sigma_{a} \star X'_{a} \ar[r]^{\overline{g}'} \ar@{-->}[ur]^{g'} & \calO^{\otimes}. }$$
If the dimension of $\sigma_{a}$ is positive, then it suffices to show that
$g'_0$ carries $\{D_a \} \star X'_{a}$ to a $q$-limit diagram in $\calC^{\otimes}$.
Let $q'$ denote the canonical map $\calO^{\otimes} \rightarrow \Nerve(\FinSeg)$.
In view of Proposition \toposref{basrel}, it suffices to show that
$g'_0$ carries $\{D_a \} \star X'_{a}$ to a $(q' \circ q)$-limit diagram and that
$q \circ g'_0$ carries $\{D_a \} \star X'_{a}$ to a $q'$-limit diagram. The first of these
assertions follows from $(ii)$ and from the fact that $\calC^{\otimes}$ is an $\infty$-operad,
and the second follows by the same argument (since $q$ is a map of $\infty$-operads).

It remains to treat the case where $\sigma_{a}$ is zero dimensional (in which case
we must have $m=1$). Since $\calC^{\otimes}$ is an $\infty$-operad, we can solve the lifting problem depicted in the diagram
$$ \xymatrix{ \bd \sigma_{a} \star X'_{a} \ar[r]^{g'_0} \ar[d] & \calC^{\otimes} \ar[d]^{q' \circ q} \\
\sigma_{a} \star X'_{a} \ar[r] \ar@{-->}[ur]^{g''} & \Nerve(\FinSeg) }$$
in such a way that $g''$ carries edges of $\sigma_{a} \star X'_{a}$ to inert morphisms
in $\calC^{\otimes}$. Since $q$ is an $\infty$-operad map, it follows that
$q \circ g''$ has the same property. Since $\calO^{\otimes}$ is an $\infty$-operad, we conclude that
$q \circ g''$ and $\overline{g}'$ are both $q'$-limit diagrams in $\calO^{\otimes}$
extending $q \circ g'_0$, and therefore equivalent to one another via
an equivalence which is fixed on $X'_{a}$ and compatible with the projection to
$\Nerve(\FinSeg)$. Since $q$ is a categorical fibration, we can lift this equivalence to an
equivalence $g'' \simeq g'$, where $g': \sigma_{a} \star X'_{a} \rightarrow \calC^{\otimes}$
is the desired extension of $g'_0$. We note that this construction ensures that condition
$(ii)$ is satisfied.

\item[$(5)$] To find the desired extension $f^5$ of $f^4$, it suffices to show that
the inclusion $\calM^{\otimes}_{4} \subseteq \calM^{\otimes}_{5}$ is a categorical equivalence. This again follows from the argument given in step $(2)$.

\item[$(6)$] The verification that $f^{5}$ can be extended to a map
$f^{6}: \calM^{\otimes}_{6} \rightarrow \calC^{\otimes}$ with the desired properties
proceeds as in step $(4)$, but is slightly easier (because $G'_{(6)}$ contains no simplices
of length zero).
\end{itemize}
\end{proof}

\subsection{Free Algebras}\label{comm33}

Let $q: \calC^{\otimes} \rightarrow \calO^{\otimes}$ be a fibration of $\infty$-operads, and suppose
we are given $\infty$-operad maps
$$ \calA^{\otimes} \stackrel{i}{\rightarrow} \calB^{\otimes} \stackrel{j}{\rightarrow} \calO^{\otimes}.$$
Composition with $i$ induces a forgetful functor $\theta: \Alg_{\calB}(\calC) \rightarrow
\Alg_{\calA}(\calC)$. Our goal in this section is to show that, under suitable conditions, the
functor $\theta$ admits a left adjoint. We can think of this left adjoint as carrying an object
$F \in \Alg_{\calA}(\calC)$ to the free $\calB$-algebra object of $\calC$ generated by $F$.
Our first step is to make this idea more precise:

\begin{definition}\label{dfree}
Let $q: \calC^{\otimes} \rightarrow \calO^{\otimes}$ and $\calA^{\otimes} \stackrel{i}{\rightarrow} \calB^{\otimes} \stackrel{j}{\rightarrow} \calO^{\otimes}$ be as above, let
$F \in \Alg_{\calA}(\calC)$ and $F \in \Alg_{\calB}(\calC)$, and let
$f: F \rightarrow \theta(F')$ be a morphism of $\calA$-algebra objects in $\calC$.

For every object $B \in \calB$, we let  $(\calA^{\otimes}_{\acti})_{/B}$ denote the fiber product
$\calA^{\otimes} \times_{ \calB^{\otimes} } (\calB^{\otimes}_{\acti})_{/B}$.
The maps $F$ and $F'$ induce maps
$\alpha, \alpha': (\calA^{\otimes}_{\acti})_{/B} \rightarrow \calC^{\otimes}_{\acti}$, $f$
determines a natural transformation $g: \alpha \rightarrow \alpha'$. We note that
$\alpha'$ lifts to a map $\overline{\alpha}': (\calA^{\otimes}_{\acti})_{/B} \rightarrow (\calC^{\otimes}_{\acti})_{/F'(B)}$;
since the projection 
$(\calC^{\otimes}_{\acti})_{/F'(B)} \rightarrow \calC^{\otimes}_{\acti} \times_{ \calO^{\otimes}_{\acti}} (\calO^{\otimes}_{\acti})_{/qF'(B)}$ is a right fibration, we can lift $g$ (in an essentially unique fashion) to a natural transformation 
$\overline{g}: \overline{\alpha} \rightarrow \overline{\alpha}'$ which is compatible with the projection
to $\calO^{\otimes}$.

We will say that $f$ {\it exhibits $F'$ as a $q$-free $\calB$-algebra generated by $F$} if
the following condition is satisfied, for every object $B \in \calB$:
\begin{itemize}
\item[$(\ast)$] The map $\overline{\alpha}$ above determines an operadic $q$-colimit diagram
$(\calA^{\otimes}_{\acti})_{/B}^{\triangleright} \rightarrow \calC^{\otimes}$.
\end{itemize}
\end{definition}

The following result guarantees that free algebras are unique (up to equivalence) whenever they exist:

\begin{proposition}\label{wax2}
Let $q: \calC^{\otimes} \rightarrow \calO^{\otimes}$ be a fibration of $\infty$-operads, and suppose we are given maps of $\infty$-operads $\calA^{\otimes} \rightarrow \calB^{\otimes} \rightarrow \calO^{\otimes}$. Let $\theta: \Alg_{\calB}(\calC) \rightarrow \Alg_{\calA}(\calC)$ denote the forgetful functor, let $F \in \Alg_{\calA}(\calC)$, let $F' \in \Alg_{\calB}(\calC)$, and let
$f: F \rightarrow \theta(F')$ be a map which exhibits $F'$ as a $q$-free $\calB$-algebra generated by $F$. For every $F'' \in \Alg_{\calB}(\calC)$, composition with $f$ induces a homotopy equivalence
$$ \gamma: \bHom_{ \Alg_{\calB}(\calC)}( F', F'') \rightarrow \bHom_{ \Alg_{\calA}(\calC)}( F, \theta(F'') ).$$
\end{proposition}

The proof of Proposition \ref{wax2} will require some preliminaries.

\begin{lemma}\label{siwd}
Let $\calM^{\otimes} \rightarrow \Nerve(\FinSeg) \times \Delta^1$ be a correspondence
between $\infty$-operads $\calA^{\otimes}$ and $\calB^{\otimes}$, and let
$q: \calC^{\otimes} \rightarrow \calO^{\otimes}$ be a fibration of $\infty$-operads.
Suppose that $n > 0$ and we are given maps of $\infty$-operad families
$$ \xymatrix{ (\calA^{\otimes} \times \Delta^n) \coprod_{ \calA^{\otimes} \times \bd \Delta^n}
(\calM^{\otimes} \times \bd \Delta^n) \ar[r]^-{f_0} \ar[d] & \calC^{\otimes} \ar[d] \\
\calM^{\otimes} \times \Delta^n \ar[r] \ar@{-->}[ur]^{f} & \calO^{\otimes} }$$
where the restriction of $f_0$ to $\calM^{\otimes} \times \{0\}$ is an operadic
$q$-left Kan extension. Then there exists an $\infty$-operad family map $f$ as indicated in the diagram. 
\end{lemma}

\begin{proof}
Let $K(0) = (\Delta^n \times \{0\}) \coprod_{ \bd \Delta^n \times \{0\} } ( \bd \Delta^n \times \Delta^1)$,
which we identify with a simplicial subset of $\Delta^{n} \times \Delta^1$. 
We define a sequence of simplicial subsets
$$ K(0) \subseteq K(1) \subseteq \ldots \subseteq K(n+1) = \Delta^n \times \Delta^1$$
so that each $K(i+1)$ is obtained from $K(i)$ by adjoining the image of the simplex
$\sigma_{i}: \Delta^{n+1} \rightarrow \Delta^n \times \Delta^1$ which is given on vertices by the formula
$$ \sigma_{i}(j) = \begin{cases} (j,0) & \text{if } j \leq n-i \\
(j-1,1) & \text{otherwise.} \end{cases}$$
We construct a compatible family of maps $f_{i}: K(i) \times_{ \Delta^1} \calM^{\otimes} \rightarrow \calC^{\otimes}$ extending $f_0$ using induction on $i$. Assuming $f_{i}$ has been constructed,
to build $f_{i+1}$ it suffices to solve the lifting problem presented in the following diagram:
$$ \xymatrix{ \Lambda^{n+1}_{n-i} \times_{ \Delta^1} \calM^{\otimes} \ar[r]^{j} \ar[d] & \calC^{\otimes} \ar[d] \\
\Delta^{n+1} \times_{ \Delta^1} \calM^{\otimes} \ar[r] \ar@{-->}[ur] & \calO^{\otimes}. }$$
If $i < n$, then Lemma \ref{skope} guarantees that $j$ is a categorical equivalence, so the dotted
arrow exists by virtue of the fact that $q$ is a categorical fibration. If $i = n$, then the
lifting problem admits a solution by Theorem \ref{oplk}.
\end{proof}

\begin{remark}\label{sweed}
Let $q: \calC^{\otimes} \rightarrow \calO^{\otimes}$, $\calA^{\otimes} \stackrel{i}{\rightarrow}
\calB^{\otimes} \stackrel{j}{\rightarrow} \calO^{\otimes}$, $\theta: \Alg_{\calB}(\calC) \rightarrow
\Alg_{\calA}(\calC)$, and $f: F \rightarrow \theta(F')$ be as in Definition \ref{dfree}.
The maps $f$, $F$, and $F'$ determine a map
$$ h: ( \calA^{\otimes} \times \Delta^1) \coprod_{ \calA^{\otimes} \times \{1\} }
\calB^{\otimes} \rightarrow \calC^{\otimes} \times \Delta^1.$$
Choose a factorization of $h$ as a composition
$$  ( \calA^{\otimes} \times \Delta^1) \coprod_{ \calA^{\otimes} \times \{1\} }
\calB^{\otimes} \stackrel{h'}{\rightarrow} \calM^{\otimes} \stackrel{h''}{\rightarrow} \calC^{\otimes} \times \Delta^1$$
where $h'$ is a categorical equivalence and $\calM^{\otimes}$ is an $\infty$-category;
we note that the composite map $\calM^{\otimes} \rightarrow \Nerve(\FinSeg) \times \Delta^1$
is exhibits $\calM^{\otimes}$ as a correspondence of $\infty$-operads. Unwinding
the definitions, we see that $f$ exhibits $F'$ as a $q$-free $\calB$-algebra generated by $F$ if and only if the map $h''$ is an operadic $q$-left Kan extension.
\end{remark}

\begin{notation}
Let $q: \calC^{\otimes} \rightarrow \calO^{\otimes}$ be a fibration of $\infty$-operads,
let $\calM^{\otimes} \rightarrow \Nerve(\FinSeg) \times \Delta^1$ be a correspondence between
$\infty$-operads, and suppose we are given a map of $\infty$-operad families
$\calM^{\otimes} \rightarrow \calO^{\otimes}$. We let $\Alg_{\calM}( \calC)$ denote
the full subcategory of $\bHom_{ \calO^{\otimes}}( \calM^{\otimes}, \calC^{\otimes})$ spanned
by those functors which preserve inert morphisms.
\end{notation}

\begin{proof}[Proof of Proposition \ref{wax2}]
The triple $(F, F', f)$ determines a map
$h: ( \calA^{\otimes} \times \Delta^1) \coprod_{ \calA^{\otimes} \times \{1\} }
\calB^{\otimes} \rightarrow \calC^{\otimes} \times \Delta^1$. Choose
a factorization of $h$ as a composition
$$  ( \calA^{\otimes} \times \Delta^1) \coprod_{ \calA^{\otimes} \times \{1\} }
\calB^{\otimes} \stackrel{h'}{\rightarrow} \calM^{\otimes} \stackrel{h''}{\rightarrow} \calC^{\otimes} \times \Delta^1$$
where $h'$ is a categorical equivalence and $\calM^{\otimes}$ is a correspondence of
$\infty$-operads from $\calA^{\otimes}$ to $\calB^{\otimes}$. Let $\Alg'_{\calM}(\calC)$ denote the full subcategory of
$\Alg_{\calM}(\calC)$ spanned by the functors $\calM^{\otimes} \rightarrow \calC^{\otimes}$ which
are operadic $q$-left Kan extensions, and let $\Alg''_{\calM}(\calC)$ denote the full subcategory of
$\Alg_{\calM}(\calC)$ spanned by those functors which are $q$-right Kan extensions of their
restrictions to $\calB^{\otimes}$. Proposition \toposref{lklk} guarantees that
the restriction map $\Alg''_{\calM}(\calC) \rightarrow \Alg_{\calB}(\calC)$ is a trivial Kan fibration;
let $s: \Alg_{\calB}(\calC) \rightarrow \Alg''_{\calM}(\calC)$ be a section. We observe that
the functor $\theta$ is equivalent to the composition of $s$ with the restriction map
$\Alg_{\calM}(\calC) \rightarrow \Alg_{\calA}(\calC)$. Let $\overline{F} \in \Alg_{\calM}(\calC)$ be
the object determined by $h''$ so that we have a homotopy commutative diagram
$$ \xymatrix{ & \Hom^{\rght}_{ \Alg_{\calM}(\calC)}( \overline{F}, s(F'') ) \ar[dr]^{\beta} \ar[dl]^{\beta'} & \\
\Hom^{\rght}_{ \Alg_{\calB}(\calC) }( F', F'') \ar[rr] & & \Hom^{\rght}_{ \Alg_{\calA}(\calC)}( F, s(F'') | \calA^{\otimes}); }$$
where $\gamma$ can be identified with the bottom horizontal map. By the two-out-of-three property,
it will suffice to show that $\beta$ and $\beta'$ are trivial Kan fibrations. For the map
$\beta$, this follows from the observation that $s(F'')$ is a $q$-right Kan extension of $F''$.
For the map $\beta'$, we apply Lemma \ref{siwd} (together with the observation that
$\overline{F}$ is an operadic $q$-left Kan extension, by virtue of Remark \ref{sweed}).
\end{proof}

We record a few other easy consequences of Lemma \ref{siwd}:

\begin{corollary}\label{kidj}
Let $q: \calC^{\otimes} \rightarrow \calO^{\otimes}$ be a fibration of $\infty$-operads,
let $\calM^{\otimes} \rightarrow \Nerve(\FinSeg) \times \Delta^1$ be a correspondence between
$\infty$-operads $\calA^{\otimes}$ and $\calB^{\otimes}$, and suppose we are given an $\infty$-operad family map
$\calM^{\otimes} \rightarrow \calO^{\otimes}$. Let $\Alg'_{\calM}(\calC)$ denote the full
subcategory of $\Alg_{\calM}(\calC)$ spanned by those maps $f: \calM^{\otimes} \rightarrow \calC^{\otimes}$ which are operadic $q$-left Kan extensions of $f | \calA^{\otimes}$
and let $\Alg'_{\calA}(\calC)$
spanned by the functors $f_0: \calA^{\otimes} \rightarrow \calC^{\otimes}$ with the 
property that for each $B \in \calB$, the map $(\calA^{\otimes}_{\acti})_{/B} \rightarrow \calC^{\otimes}_{\acti}$ can be extended to an operadic $q$-colimit diagram lifting
the composite map
$$ (\calA^{\otimes}_{\acti})_{/B}^{\triangleright} \rightarrow \calM^{\otimes}
\rightarrow \calO^{\otimes}.$$
Then the restriction map $\theta: \Alg'_{\calM}(\calC) \rightarrow \Alg'_{\calA}(\calC)$ is a
trivial Kan fibration.
\end{corollary}

\begin{proof}
We must show that $\theta$ has the right lifting property with respect to every inclusion
$\bd \Delta^n \subseteq \Delta^n$. If $n =0$, this follows from Theorem \ref{oplk}; if
$n > 0$, it follows from Lemma \ref{siwd}.
\end{proof}

\begin{corollary}
Using the notation of Corollary \ref{kidj}, suppose that $\Alg'_{\calA}(\calC) = \Alg_{\calA}(\calC)$.
Then the restriction functor $\theta: \Alg_{\calM}(\calC) \rightarrow \Alg_{\calA}(\calC)$ admits a fully faithful left adjoint, given by any section $s$ of $\theta': \Alg'_{\calM}(\calC) \rightarrow \Alg_{\calA}(\calC)$. In particular, $\Alg'_{\calM}(\calC)$ is a colocalization of $\Alg_{\calM}(\calC)$.
\end{corollary}

\begin{proof}
Corollary \ref{kidj} implies that $\theta'$ admits a section $s$. We claim that the identity transformation
$\id_{ \Alg_{\calA}(\calC)} \rightarrow \theta \circ s$ exhibits $s$ as a left adjoint to $\theta$.
To prove this, it suffices to show that for every $A \in \Alg_{\calA}(\calC)$ and every
$B \in \Alg_{\calM}(\calC)$, the restriction map
$$ \Hom^{\rght}_{\Alg_{\calM}(\calC)}( s(A), B) \rightarrow \Hom^{\rght}_{ \Alg_{\calA}(\calC)}(A, \theta(B))$$
is a homotopy equivalence. Lemma \ref{siwd} guarantees that this map is a trivial Kan fibration. 
\end{proof}

We now consider the question of existence for free algebra objects.

\begin{proposition}\label{wax1}
Let $q: \calC^{\otimes} \rightarrow \calO^{\otimes}$ be a fibration of $\infty$-operads, and suppose we are given maps of $\infty$-operads $\calA^{\otimes} \rightarrow \calB^{\otimes} \rightarrow \calO^{\otimes}$. Let $\theta: \Alg_{\calB}(\calC) \rightarrow \Alg_{\calA}(\calC)$ denote the forgetful functor and let $F \in \Alg_{\calA}(\calC)$. The following conditions are equivalent:
\begin{itemize}
\item[$(1)$] There exists $F' \in \Alg_{\calB}(\calC)$ and a map $F \rightarrow \theta(F')$ which
exhibits $F'$ as a $q$-free $\calB$-algebra generated by $F$.

\item[$(2)$] For every object $B \in \calB$, the induced map
$$(\calA^{\otimes}_{\acti})_{/B} \rightarrow \calA^{\otimes}_{\acti} \stackrel{F}{\rightarrow} \calC^{\otimes}_{\acti}$$
can be extended to an operadic $q$-colimit diagram lying over the composition
$$(\calA^{\otimes}_{\acti})_{/B}^{\triangleright} \rightarrow \calB^{\otimes}_{\acti}
\rightarrow \calO^{\otimes}_{\acti}.$$
\end{itemize}
\end{proposition}

\begin{proof}
The implication $(1) \Rightarrow (2)$ is obvious. For the converse, we use the small object argument to choose a factorization of the map $h: ( \calA^{\otimes} \times \Delta^1) \coprod_{ \calA^{\otimes} \times \{1\} }
\calB^{\otimes} \rightarrow \calO^{\otimes} \times \Delta^1$ as a composition
$$ ( \calA^{\otimes} \times \Delta^1) \coprod_{ \calA^{\otimes} \times \{1\} }
\calB^{\otimes} \stackrel{h'}{\rightarrow} \calM^{\otimes} \stackrel{h''}{\rightarrow} \calO^{\otimes} \times \Delta^1$$
where $h'$ is inner anodyne and $\calM^{\otimes}$ is a correspondence of $\infty$-operads from
$\calA^{\otimes}$ to $\calB^{\otimes}$. Using assumption $(2)$ and Theorem \ref{oplk}, we can
solve the lifting problem depicted in the diagram
$$ \xymatrix{ \calA^{\otimes} \ar[d] \ar[r] & \calC^{\otimes} \ar[d] \\
\calM^{\otimes} \ar[r] \ar@{-->}[ur]^{G} & \calO^{\otimes} }.$$
in such a way that $G$ is an operadic $q$-left Kan extension. Composing
$G$ with $h'$, we obtain an object $F' \in \Alg_{\calB}(\calC)$ and a natural transformation
$f: F \rightarrow \theta(F')$. It follows from Remark \ref{sweed} that $f$ exhibits
$F'$ as a $q$-free $\calB$-algebra generated by $F$.
\end{proof}

\begin{corollary}\label{slavetime}
Let $q: \calC^{\otimes} \rightarrow \calO^{\otimes}$ be a fibration of $\infty$-operads, and suppose
we are given maps of $\infty$-operads $\calA^{\otimes} \rightarrow \calB^{\otimes} \rightarrow \calO^{\otimes}$. Assume that the following condition is satisfied:
\begin{itemize}
\item[$(\ast)$] For every object $B \in \calB$ and every $F \in \Alg_{\calA}(\calC)$, the diagram
$$ (\calA^{\otimes}_{\acti})_{/B} \rightarrow \calA^{\otimes}_{\acti} \stackrel{F}{\rightarrow} \calC^{\otimes}_{\acti}$$
can be extended to an operadic $q$-colimit diagram lifting the composition
$$ (\calA^{\otimes}_{\acti})^{\triangleright}_{/B} \rightarrow (\calB^{\otimes}_{\acti})_{/B}^{\triangleright} \rightarrow \calB^{\otimes}_{\acti} \rightarrow \calO^{\otimes}_{\acti}.$$
\end{itemize}
Then the forgetful functor $\theta: \Alg_{\calB}(\calC) \rightarrow \Alg_{\calA}(\calC)$ admits
a left adjoint, which carries each $F \in \Alg_{\calA}(\calC)$ to a $q$-free $\calB$-algebra generated by $F$.
\end{corollary}

\begin{proof}
Combine Propositions \ref{wax1}, \ref{wax2}, and \toposref{sumpytump}.
\end{proof}

\begin{corollary}\label{spaltwell}
Let $\kappa$ be an uncountable regular cardinal, let $\calO^{\otimes}$ be an $\infty$-operad,
and let $\calC^{\otimes} \rightarrow \calO^{\otimes}$ be a $\calO$-monoidal $\infty$-category which is
compatible with $K$-indexed colimits for every $\kappa$-small simplicial set $K$ (see Definition \ref{compatK}). Let ${\calO'}^{\otimes} \rightarrow \calO^{\otimes}$ be a map of $\infty$-operads, and suppose that both $\calO^{\otimes}$ and ${\calO'}^{\otimes}$ are $\kappa$-small. Then
the forgetful functor $\Alg_{\calO}(\calC) \rightarrow \Alg_{\calO'}(\calC)$ admits a left adjoint,
which carries each $\calO'$-algebra $F$ in $\calC$ to a $q$-free $\calO$-algebra generated by $F$.
\end{corollary}

\begin{proof}
Combine Proposition \ref{slavewell} with Corollary \ref{slavetime}.
\end{proof}

\begin{example}\label{lad}
Let $\calO^{\otimes}$ be a $\kappa$-small $\infty$-operad, let $p: \calC^{\otimes} \rightarrow \calO^{\otimes}$ be a coCartesian fibration of $\infty$-operads which is compatible with
$\kappa$-indexed colimits. Applying Corollary \ref{spaltwell} in the
case ${\calO'}^{\otimes} = \calO^{\otimes} \times_{ \Nerve(\FinSeg)} \Triv$ and applying Example
\ref{algtriv}, we deduce that the forgetful functor
$\Alg_{\calO}(\calC) \rightarrow \bHom_{ \calO}( \calO, \calC)$ admits a left adjoint.
\end{example}

To apply Corollary \ref{spaltwell} in practice, it it convenient to have a more explicit description of
free algebras. To obtain such a description, we will specialize to the case where the $\infty$-operad
${\calO'}^{\otimes}$ is trivial.

\begin{construction}\label{aball}
Let $\calO^{\otimes}$ be an $\infty$-operad, let $X,Y \in \calO$ be objects, and let
$f: \Triv \rightarrow \calO^{\otimes}$ be a map of $\infty$-operads such that
$f( \seg{1} ) = X$ (such a map exists an is unique up to equivalence, by Remark \ref{coalgtriv}).

Let $p: \calC^{\otimes} \rightarrow \calO^{\otimes}$ be a $\calO$-monoidal $\infty$-category
and let $C \in \calC_{X}$ be an object. Using Example \ref{algtriv}, we conclude that
$C$ determines an essentially unique $\Triv$-algebra $\overline{C} \in \Alg_{\Triv}( \calC)$
such that $\overline{C}( \seg{1}) = C$.

For each $n \geq 0$, let $\calP(n)$ denote the full subcategory of
$\Triv \times_{ \calO^{\otimes} } \calO^{\otimes}_{/Y}$ spanned by the active morphisms $f(\seg{n}) \rightarrow Y$; we observe that $\calP(n)$ is an $\infty$-groupoid and that the fibers
of the canonical map $q: \calP(n) \rightarrow \Nerve( \Sigma_n)$ can be identified with
the space of $n$-ary operations $\Mul_{\calO}( \{ X \}_{1 \leq i \leq n}, Y)$. 
By construction, we have a canonical map $h: \calP(n) \times \Delta^1 \rightarrow \calO^{\otimes}$
which we regard as a natural transformation from $h_0 = f \circ q$ to
the constant map $h_1: \calP(n) \rightarrow \{ Y \}$. Since $p$ is a coCartesian
fibration, we can choose a $p$-coCartesian natural transformation
$\overline{h}: \overline{C} \circ q \rightarrow \overline{h}_1$, for some map
$\overline{h}_1: \calP(n) \rightarrow \calC_{Y}$. We let
$\Sym^{n}_{\calO,Y}(C)$ denote a colimit of the diagram
$\overline{h}_1$, if such a homotopy colimit exists. We observe that
$\Sym^{n}_{\calO,Y}(C) \in \calC_{Y}$ is well-defined up to equivalence (in fact, up to a contractible space of choices).
\end{construction}

\begin{notation}
In the special case $X = Y$, we will simply write $\Sym^{n}_{\calO}(C)$ for
$\Sym^{n}_{\calO,Y}(C)$. When $\calO^{\otimes}$ is the commutative $\infty$-operad $\Nerve(\FinSeg)$
(so that $X = Y = \seg{1}$) we will denote $\Sym^{n}_{\calO}(C)$ by $\Sym^{n}(C)$.
\end{notation}

\begin{remark}\label{spann}
In the situation of Construction \ref{aball}, suppose that $A \in \Alg_{\calO}(\calC)$
is a $\calO$-algebra and that we are equipped with a map
$\overline{C} \rightarrow A \circ f$ of $\Triv$-algebras (in view of Example
\ref{algtriv}, this is equivalent to giving a map $C \rightarrow A(X)$ in $\calC_{X}$). 
This map induces a natural transformation from $\overline{h}_1$ to the constant map
$\calP(n) \rightarrow \{ A(Y) \}$, which we can identify with a map
$\Sym^{n}_{\calO,Y}(C) \rightarrow A(Y)$ provided that the left side is defined.
\end{remark}

\begin{definition}
Let $q: \calC^{\otimes} \rightarrow \calO^{\otimes}$ be a fibration of $\infty$-operads,
let $X \in \calO$, let $A \in \Alg_{\calO}(\calC)$, and suppose we are given a morphism
$f: C \rightarrow A(X)$ in $\calC_{X}$. Using Remark \ref{coalgtriv} and Example \ref{algtriv}, we can extend (in an essentially unique way) $X$ to an $\infty$-operad map $\Triv \rightarrow \calO$, $C$ to
an object $\overline{C} \in \Alg_{\Triv}(\calC)$, and $f$ to a map of $\Triv$-algebras
$\overline{f}: \overline{C} \rightarrow A| \Triv$. We will say that {\it $f$ exhibits
$A$ as a free $\calO$-algebra generated by $C$} if $\overline{f}$ exhibits
$A$ as a $q$-free $\calO$-algebra generated by $\overline{C}$.
\end{definition}

\begin{remark}\label{tappus}
Let $\kappa$ be an uncountable regular cardinal, $\calO^{\otimes}$ a $\kappa$-small
$\infty$-operad, and $q: \calC^{\otimes} \rightarrow \calO^{\otimes}$ a $\calO$-monoidal
$\infty$-category which is compatible with $\kappa$-small colimits. Let $X \in \calO$ and
let $C \in \calC_{X}$. From the above discussion we deduce the following:
\begin{itemize}
\item[$(i)$] There exists an algebra $A \in \Alg_{\calO}(\calC)$ and a map
$C \rightarrow A(X)$ which exhibits $A$ as a free $\calO$-algebra generated by
$X$.
\item[$(ii)$] An arbitrary map $f: C \rightarrow A(X)$ as in $(i)$ exhibits
$A$ as a free $\calO$-algebra generated by $X$ if and only if, for
every object $Y \in \calO$, the maps
$\Sym^{n}_{\calO,Y}(C) \rightarrow A(Y)$ of Remark \ref{spann} exhibit
$A(Y)$ as a coproduct $\coprod_{n \geq 0} \Sym^{n}_{\calO,Y}(C)$.
\end{itemize}
Assertion $(ii)$ follows by combining Corollary \ref{spaltwell} with Propositions
\ref{chocolateoperad} and \ref{optest}.
\end{remark}

\begin{example}
Let $\calC^{\otimes}$ be a symmetric monoidal $\infty$-category. Assume that
the underlying $\infty$-category $\calC$ admits countable colimits, and that for
each $X \in \calC$ the functor $Y \mapsto X \otimes Y$ preserves countable colimits.
Then the forgetful functor $\CAlg(\calC) \rightarrow \calC$ admits a left adjoint, which
is given informally by the formula
$$ C \mapsto \coprod_{n \geq 0} \Sym^{n}(C).$$
\end{example}

\subsection{Colimits of Algebras}\label{limco2}

Let $p: \calC^{\otimes} \rightarrow \calO^{\otimes}$ be a fibration of $\infty$-operads.
Our goal in this section is to construct colimits in the $\infty$-category
$\Alg_{\calO}(\calC)$, given suitable hypotheses on $p$. Our strategy is as follows:
we begin by constructing {\em sifted} colimits in $\Alg_{\calO}(\calC)$ (Proposition \ref{fillfemme}), which are given by forming colimits of the underlying objects in $\calC$. We will then combine
this construction with existence results for free algebras (see \S \ref{comm33}) to deduce
the existence of general colimits (Corollary \ref{smurf}). 

Our main result can be stated as follows:

\begin{proposition}\label{fillfemme}
Let $K$ be a sifted simplicial set and let $p: \calC^{\otimes} \rightarrow \calO^{\otimes}$ be a coCartesian fibration of $\infty$-categories which is compatible with $K$-indexed colimits.
Then:
\begin{itemize}
\item[$(1)$] The $\infty$-category $\bHom_{ \calO^{\otimes}}( \calO^{\otimes}, \calC^{\otimes})$ of
sections of $p$ admits $K$-indexed colimits.

\item[$(2)$] A map $\overline{f}: K^{\triangleright} \rightarrow \bHom_{ \calO^{\otimes}}( \calO^{\otimes}, \calC^{\otimes})$ is a colimit diagram if and only if, for each $X \in \calO^{\otimes}$, the induced diagram
$\overline{f}_{X}: K^{\triangleright} \rightarrow \calC^{\otimes}_{X}$ is a colimit diagram.

\item[$(3)$] The full subcategories 
$$\Fun^{\otimes}_{\calO}(\calO, \calC) \subseteq \Alg_{\calO}(\calC)
\subseteq \bHom_{ \calO^{\otimes}}( \calO^{\otimes}, \calC^{\otimes})$$
are stable under $K$-indexed colimits.

\item[$(4)$] A map $\overline{f}: K^{\triangleright} \rightarrow \Alg_{\calO}(\calC)$ is a colimit diagram if and only if, for each $X \in \calO$, the induced diagram
$\overline{f}_{X}: K^{\triangleright} \rightarrow \calC_{X}$ is a colimit diagram.
\end{itemize}
\end{proposition}

We will give the proof of Proposition \ref{fillfemme} at the end of this section.

\begin{corollary}\label{filtfemme}
Let $\calC$ be a symmetric monoidal $\infty$-category, and let $K$ be a sifted simplicial set such that the symmetric monoidal structure on $\calC$ is compatible with $K$-indexed colimits.
Then:
\begin{itemize}
\item[$(1)$] The $\infty$-category $\CAlg(\calC)$ admits $K$-indexed colimits.

\item[$(2)$] The forgetful functor $\theta: \CAlg(\calC) \rightarrow \calC$ detects $K$-indexed colimits. More precisely, a map $\overline{q}: K^{\triangleright} \rightarrow \CAlg(\calC)$ is a colimit diagram if and only if $\theta \circ \overline{q}$ is a colimit diagram.
\end{itemize}
\end{corollary}

\begin{proof}
Apply Proposition \ref{fillfemme} in the case where $\calO^{\otimes}$ is the commutative $\infty$-operad.
\end{proof}

We now treat the case of more general colimits.

\begin{corollary}\label{smurf}
Let $\kappa$ be an uncountable regular cardinal and let
$\calO^{\otimes}$ be an essentially $\kappa$-small $\infty$-operad. 
Let $p: \calC^{\otimes} \rightarrow \calO^{\otimes}$ be a coCartesian fibration of $\infty$-operads which is compatible with $\kappa$-small colimits. Then $\Alg_{\calO}(\calC)$ admits $\kappa$-small colimits.
\end{corollary}

\begin{proof}
Without loss of generality we may assume that $\calO^{\otimes}$ is $\kappa$-small.
In view of Corollary \toposref{uterrr} and Lemma \stableref{simpenough}, arbitrary $\kappa$-small colimits in $\Alg_{\calO}(\calC)$ can be built from $\kappa$-small filtered colimits, geometric realizations of simplicial objects, and finite coproducts. Since $\Alg_{\calO}(\calC)$ admits $\kappa$-small sifted colimits (Corollary \ref{filtfemme}), it will suffice to prove that $\Alg_{\calO}(\calC)$ admits finite coproducts.

According to Example \ref{lad}, the forgetful functor $G: \Alg_{\calO}(\calC) \rightarrow
\bHom_{ \calO}( \calO, \calC)$ admits a left adjoint which we will denote by $F$.
Let us say that an object of $\Alg_{\calO}(\calC)$ is {\it free} if it belongs to the essential image of
$F$. Since $\bHom_{\calO}(\calO, \calC)$ admits $\kappa$-small colimits (Lemma \monoidref{hunkerin}) and $F$ preserves $\kappa$-small colimits (Proposition \toposref{adjointcol}), a finite collection of objects $\{ A^i \}$ of $\Alg_{\calO}(\calC)$ admits a coproduct whenever each $A_i$ is free.

We now claim the following:
\begin{itemize}
\item[$(\ast)$] For every object $A \in \Alg_{\calO}(\calC)$, there exists a simplicial object
$A_{\bigdot}$ of $\Alg_{\calO}(\calC)$ such that $A$ is equivalent to the geometric realization
$| A_{\bigdot} |$ and each $A_n$ is free.
\end{itemize}

This follows from Proposition \monoidref{littlebeck} since the forgetful functor $G$ is
conservative (Corollary \ref{jumunj22}) and commutes with geometric realizations
(Corollary \ref{filtfemme}). 

Now let $\{ A^i \}$ be an arbitrary finite collection of objects of $\Alg_{\calO}(\calC)$;
we wish to prove that $A^{i}$ admits a coproduct in $\Alg_{\calO}(\calC)$.
According to $(\ast)$, each $A^{i}$ can be obtained as the geometric realization of a
simplicial object $A^{i}_{\bigdot}$ of $\Alg_{\calO}(\calC)$, where
each $A^{i}_{n}$ is free. It follows that the diagrams $A^{i}_{\bigdot}$ admit
a coproduct $A_{\bigdot}: \Nerve(\cDelta^{op}) \rightarrow \Alg_{\calO}(\calC)$.
Using Lemma \toposref{limitscommute}, we conclude that a geometric realization of
$A_{\bigdot}$ (which exists by virtue of Corollary \ref{filtfemme}) is a coproduct
for the collection $A^{i}$.
\end{proof}

We next give a criterion for establishing that an $\infty$-category of algebras
$\Alg_{\calO}(\calC)$ is presentable.

\begin{lemma}\label{sillyfunn}
Let $\calO^{\otimes}$ be a small $\infty$-operad, and let
$p: \calC^{\otimes} \rightarrow \calO^{\otimes}$ be a coCartesian fibration of $\infty$-operads.
The following conditions are equivalent:
\begin{itemize}
\item[$(1)$] For each $X \in \calO$, the fiber $\calC_{X}$ is accessible, and $p$
is compatible with $\kappa$-filtered colimits for $\kappa$ sufficiently large.

\item[$(2)$] Each fiber of $p$ is accessible, and for every morphism $f: X \rightarrow Y$
in $\calO^{\otimes}$, the associated functor $f_{!}: \calC^{\otimes}_{X} \rightarrow \calC^{\otimes}_{Y}$ is accessible.
\end{itemize}
\end{lemma}

\begin{proof}
The implication $(2) \Rightarrow (1)$ is obvious. Conversely, suppose that $(1)$ is satisfied.
For each $X \in \calO^{\otimes}_{\seg{m}}$, choose inert morphisms
$X \rightarrow X_i$ lying over $\Colp{i}: \seg{m} \rightarrow \seg{1}$ for
$1 \leq i \leq m$. Proposition \ref{saber} implies that
$\calC^{\otimes}_{X}$ is equivalent to the product $\prod_{1 \leq i \leq m} \calC_{X_i}$, which is accessible by virtue of $(1)$ and Lemma \toposref{complexhorse}. 
Now suppose that $f: X \rightarrow Y$ is a morphism in $\calO^{\otimes}$; we wish to prove
that the functor $f_{!}: \calC^{\otimes}_{X} \rightarrow \calC^{\otimes}_{Y}$ is accessible.
This follows from $(1)$ and Lemma \ref{hungerine}.
\end{proof}

\begin{corollary}\label{algprice1}
Let $p: \calC^{\otimes} \rightarrow \calO^{\otimes}$ be a coCartesian fibration of $\infty$-operads,
where $\calO^{\otimes}$ is small. Assume that for each $X \in \calO$, the fiber $\calC_{X}$ is accessible.

\begin{itemize}
\item[$(1)$] If the $\calO$-monoidal structure on $\calC$ is compatible with $\kappa$-filtered colimits for $\kappa$ sufficiently large,
$\Alg_{\calO}(\calC)$ is an accessible $\infty$-category.
\item[$(2)$] Suppose that each fiber
$\calC_{X}$ is presentable and that the $\calO$-monoidal structure on $\calC$ is compatible with small colimits. Then $\Alg_{\calO}(\calC)$ is a presentable $\infty$-category.
\end{itemize}
\end{corollary}

In particular, if $\calC$ is a presentable $\infty$-category equipped with a symmetric monoidal structure, and the tensor product $\otimes: \calC \times \calC \rightarrow \calC$ preserves colimits separately in each variable, then $\CAlg(\calC)$ is again presentable.

\begin{proof}
Assertion $(1)$ follows from Lemma \ref{sillyfunn} and Proposition \toposref{prestorkus}. To prove $(2)$, we combine $(1)$ with Corollary \ref{smurf}.
\end{proof}

The proof of the existence for colimits in $\Alg_{\calO}(\calC)$ given in Corollary \ref{smurf} is rather indirect. In general, this seems to be necessary: there is no simple formula which describes
the coproduct of a pair of associative algebras, for example. However, for commutative algebras we can be more explicit: coproducts can be computed by forming the tensor product of the underlying algebras. This assertion can be formulated more precisely as follows:

\begin{proposition}\label{cocarten}
Let $\calC$ be a symmetric monoidal $\infty$-category. Then the symmetric monoidal
structure on $\CAlg(\calC)$ provided by Example \ref{slabber} is coCartesian.
\end{proposition}

\begin{proof}
We will show that the symmetric monoidal structure on $\CAlg(\calC)$ satisfies (the dual of) criterion
$(3)$ of Proposition \ref{fingertone}. We first show that the unit object of $\CAlg(\calC)$ is initial. It follows from Example \ref{slabber} that the forgetful functor $\theta: \CAlg(\calC) \rightarrow \calC$ can be promoted to a symmetric monoidal functor. In particular, $\theta$ carries the unit object $A \in \CAlg(\calC)$ to the unit object of $\calC$. It follows from Corollary \ref{gargle2} that $A$ is initial in $\CAlg(\calC)$.

It remains to show that we can produce a collection of codiagonal
map $\{ \delta_{A}: A \otimes A \rightarrow A\}_{A \in \CAlg(\calC)}$, satisfying the axioms $(i)$, $(ii)$, and $(iii)$ of 
Proposition \ref{fingertone}. For this, we observe that composition with the multiplication functor
functor $\wedge$ of Notation \ref{crumb} induces a functor
$l: \CAlg(\calC) \rightarrow \CAlg( \CAlg(\calC) )$. In particular, for every object $A \in \CAlg(\calC)$, the multiplication on $l(A)$ induces a map $\delta_{A}: A \otimes A \rightarrow A$ in
the $\infty$-category $\CAlg(\calC)$. It is readily verified that these maps possess the desired properties.
\end{proof}

\begin{corollary}\label{exi}
Let $\calC$ be a symmetric monoidal $\infty$-category. Then the $\infty$-category
$\CAlg(\calC)$ admits finite coproducts.
\end{corollary}

\begin{remark}\label{slaphap}
Let $\calC$ be an $\infty$-category which admits finite coproducts, and let
$\calD^{\otimes}$ be a symmetric monoidal $\infty$-category. According to Theorem
\ref{kinema}, an $\infty$-operad map $F: \calC^{\amalg} \rightarrow \calD^{\otimes}$
is classified up to equivalence by the induced functor $f: \calC \rightarrow \CAlg(\calD)$.
We note that $F$ is a symmetric monoidal functor if and only if for every finite
collection of objects $C_{i} \in \calC$, the induced map
$\theta: \otimes f( C_i ) \rightarrow f( \coprod C_i )$ is an equivalence in $\CAlg(\calD)$. 
According to Proposition \ref{cocarten}, we can identify the domain of $\theta$ with the
coproduct of the commutative algebra obejcts $f(C_i)$. Consequently, we conclude that
$F$ is symmetric monoidal if and only if $f$ preserves finite coproducts.
\end{remark}

The following generalization of Proposition \ref{cocarten} is valid for an arbitrary $\infty$-operad $\calC^{\otimes}$.

\begin{proposition}\label{socarte}
Let $\calC^{\otimes}$ be an $\infty$-operad. Then the $\infty$-operad
$\CAlg(\calC)^{\otimes}$ of Example \ref{slabber} is coCartesian.
\end{proposition}

\begin{proof}
Consider the functor $\CAlg(\calC) \rightarrow \CAlg( \CAlg(\calC) )$ appearing
in the proof of Proposition \ref{cocarten}. According to Theorem \ref{kinema}, this map
determines an $\infty$-operad map $\theta: \CAlg(\calC)^{\amalg} \rightarrow \CAlg(\calC)^{\otimes}$ which is well-defined up to homotopy. We wish to prove that $\theta$ is an equivalence.
We first observe that $\theta$ induces an equivalence between the underlying $\infty$-categories, and is therefore essentially surjective. To prove that $\theta$ is fully faithful, choose a fully faithful
map of $\infty$-operads $\calC^{\otimes} \rightarrow \calD^{\otimes}$ where $\calD^{\otimes}$ is a symmetric monoidal $\infty$-category (Remark \ref{lazier}). We have a homotopy commutative diagram
$$ \xymatrix{ \CAlg(\calC)^{\amalg} \ar[r]^{\theta} \ar[d] & \CAlg(\calC)^{\otimes} \ar[d] \\
\CAlg(\calD)^{\amalg} \ar[r]^{\theta'} & \CAlg(\calD)^{\otimes} }$$
where the vertical maps are fully faithful. Consequently, it will suffice to show that
$\theta'$ is an equivalence of symmetric monoidal $\infty$-categories.
In view of Remark \ref{symmetricegg}, it will suffice to show that
$\theta'$ is a symmetric monoidal functor. This is equivalent to the assertion that the composite functor
$\CAlg(\calD)^{\amalg} \rightarrow \CAlg(\calD)^{\otimes} \rightarrow \calD^{\otimes}$ is symmetric monoidal (Proposition \ref{slape}), which follows from Remark \ref{slaphap} because the underlying functor $\CAlg(\calD) \rightarrow \CAlg(\calD)$ is equivalent to the identity (and therefore preserves finite coproducts).
\end{proof}

We now turn to the proof of Proposition \ref{fillfemme}.
We will need the following generalization of Proposition \toposref{urbil}:

\begin{lemma}\label{posturbil}
Let $K$ be a sifted simplicial set, let $\{ \calC_i \}_{1 \leq i \leq n}$ be a finite collection
of $\infty$-categories which admit $K$-indexed colimits, let $\calD$ be another $\infty$-category which admits $K$-indexed colimits, and let $F: \prod_{1 \leq i \leq n} \calC_i \rightarrow \calD$ be a functor. Suppose that $F$ preserves $K$-indexed colimits
separately in each variable: that is, if we are given $1 \leq i \leq n$ and objects
$\{ C_j \in \calC_{j} \}_{j \neq i}$, then the restriction of $F$ to
$\{ C_1 \} \times \cdots \times \{ C_{i-1} \} \times \calC_{i} \times \{ C_{i+1} \}
\times \cdots \times \{ C_n \}$ preserves $K$-indexed colimits.
Then $F$ preserves $K$-indexed colimits.
\end{lemma}

\begin{proof}
Choose $S \subseteq \{1, \ldots, n \}$, and consider the following assertion:
\begin{itemize}
\item[$(\ast_{S})$] Fix objects $\{ C_i \in \calC_i \}_{i \notin S}$. Then the restriction of
$F$ to
$\prod_{i \in S} \calC_i \times \prod_{i \notin S} \{ C_i \}$ preserves $K$-indexed colimits.
\end{itemize}
We will prove $(\ast_{S})$ using induction on the cardinality of $S$.
Taking $S = \{ 1, \ldots, n \}$, we can deduce that $F$ preserves $K$-indexed colimits and complete the proof.

If $S$ is empty, then $(\ast_S)$ follows from Corollary \toposref{silt}, since
$K$ is weakly contractible (Proposition \toposref{siftcont}). If $S$ contains a single element, then
$(\ast_{S})$ follows from the assumption that $F$ preserves $K$-indexed colimits separately in each variable. We may therefore assume that $S$ has at least two elements, so we can write
$S = S' \cup S''$ where $S'$ and $S''$ are disjoint nonempty subsets of $S$. We now observe
that $(\ast_{S})$ follows from $(\ast_{S'})$, $(\ast_{S''})$, and Proposition \toposref{urbil}.
\end{proof}

\begin{lemma}\label{hungerine}
Let $K$ be a sifted simplicial set, let $\calO^{\otimes}$ be an $\infty$-operad, and let
$p: \calC^{\otimes} \rightarrow \calO^{\otimes}$ be a coCartesian fibration of $\infty$-operads which is compatible with $K$-indexed colimits (Definition \ref{compatK}). For each morphism
$f: X \rightarrow Y$ in $\calO^{\otimes}$, the associated functor $f_{!}: \calC^{\otimes}_{X}
\rightarrow \calC^{\otimes}_{Y}$ preserves $K$-indexed colimits.
\end{lemma}

\begin{proof}
Factor $f$ as a composition $X \stackrel{f'}{\rightarrow} Z \stackrel{f''}{\rightarrow} Y$, where
$f'$ is inert and $f''$ is active. Using Proposition \ref{saber}, we deduce that
$\calC^{\otimes}_{X}$ is equivalent to a product $\calC^{\otimes}_{Z} \times \calC^{\otimes}_{Z'}$ and
that $f'_{!}$ can be identified with the projection onto the first factor. Since $\calC^{\otimes}_{Z}$ and
$\calC^{\otimes}_{Z'}$ both admit $K$-indexed colimits, this projection preserves $K$-indexed colimits. We may therefore replace $X$ by $Z$ and thereby reduce to the case where $f$ is active.

Let $\seg{n}$ denote the image of $Y$ in $\Nerve(\FinSeg)$, and choose inert morphisms
$g_i: Y \rightarrow Y_i$ lying over $\Colp{i}: \seg{n} \rightarrow \seg{1}$ for
$1 \leq i \leq n$. Applying Proposition \ref{saber} again, we deduce that the functors
$(g_{i})_{!}$ exhibit $\calC^{\otimes}_{Y}$ as equivalent to the product
$\prod_{1 \leq i \leq n} \calC^{\otimes}_{Y_i}$. It will therefore suffice to show that each of the
functors $(g_i \circ f)_{!}$ preserve $K$-indexed colimits. Replacing $Y$ by $Y_i$, we can
reduce to the case $n=1$.

Let $\seg{m}$ denote the image of $X$ in $\Nerve(\FinSeg)$ and choose inert morphisms
$h_{j}: X \rightarrow X_j$ lying over $\Colp{j}: \seg{m} \rightarrow \seg{1}$ for
$1 \leq j \leq m$. Invoking Proposition \ref{saber} again, we obtain an equivalence
$\calC^{\otimes}_{X} \simeq \prod_{1 \leq j \leq m} \calC^{\otimes}_{X_i}$. It will therefore
suffice to show that the composite map
$$ \phi: \prod_{1 \leq j \leq m} \calC^{\otimes}_{X_i} \simeq \calC^{\otimes}_{X} \stackrel{f_{!}}{\rightarrow} \calC^{\otimes}_{Y}$$
preserves $K$-indexed colimits. Since $p$ is compatible with $K$-indexed colimits, 
we conclude that $\phi$ preserves $K$-indexed colmits separately in each variable.
Since $K$ is sifted, the desired result now follows from Lemma \ref{posturbil}.
\end{proof}

\begin{proof}[Proof of Proposition \ref{fillfemme}]
Assertions $(1)$, $(2)$, and $(3)$ follow immediately from Lemmas \ref{hungerine} and
\monoidref{hunkerin}. The ``only if'' direction of $(4)$ follows immediately from $(2)$.
To prove the converse, suppose that $\overline{f}: K^{\triangleleft} \rightarrow
\Alg_{\calO}(\calC)$ has the property that $\overline{f}_{X}$ is a colimit diagram in
$\calC_{X}$ for each $X \in \calO$. We wish to prove that the analogous assertion holds
for any $X \in \calO^{\otimes}_{\seg{n}}$. Choose inert morphisms
$g(i): X \rightarrow X_i$ lying over $\Colp{i}: \seg{n} \rightarrow \seg{1}$ for
$1 \leq i \leq n$. According to Proposition \ref{saber}, the functors
$g(i)_{!}$ induce an equivalence $\calC^{\otimes}_{X} \simeq \prod_{1 \leq i \leq n} \calC_{X_i}$. 
It will therefore suffice to prove that each composition $g(i)_{!} \circ \overline{f}_{X}$
is a colimit diagram $K^{\triangleright} \rightarrow \calC_{X_i}$. This follows
from the observation that $g(i)_{!} \circ \overline{f}_{X} \simeq \overline{f}_{X_i}$.
\end{proof}

\section{Modules}\label{comm6}

Let $\calC$ be a symmetric monoidal category, and suppose that
$A$ is an algebra object of $\calC$: that is, an object of $\calC$ equipped with an
associative (and unital) multiplication $m: A \otimes A \rightarrow A$. It then makes sense to consider left
$A$-modules in $\calC$: that is, objects $M \in \calC$ equipped with a map $a: A \otimes M \rightarrow
M$ such that the following diagrams commute
$$ \xymatrix{ A \otimes A \otimes M \ar[r]^-{ m \otimes \id_{M} } \ar[d]^{ \id_{A} \otimes a} & A \otimes M \ar[d]^{a} & {\bf 1} \otimes M \ar[rr] \ar[dr] & & A \otimes M \ar[dl]^{a} \\
A \otimes M \ar[r]^{a} & M & & M. & & }$$
There is also a corresponding theory of right $A$-modules in $\calC$. If the multiplication on
$A$ is commutative, then we can identify left modules with right modules and speak simply of
{\it $A$-modules in $\calC$}. Moreover, a new phenomenon occurs: in many cases, the category
$\Mod_{A}(\calC)$ of $A$-modules in $\calC$ inherits the structure of a symmetric monoidal category
with respect to the relative tensor product over $A$: namely, we can define
$M \otimes_{A} N$ to be the coequalizer of the diagram
$$\xymatrix{ M \otimes A \otimes N \ar@<.4ex>[r] \ar@<-.4ex>[r] & M \otimes N}$$
determined by the actions of $A$ on $M$ and $N$ (this construction is not always sensible: we must assume that the relevant coequalizers exist, and that they behave well with respect to the symmetric monoidal structure on $\calC$). 

If we do not assume that $A$ is commutative, then we can often still make sense of the relative tensor product $M \otimes_{A} N$ provided that $M$ is a right $A$-module and $N$ is a left $A$-module.
However, this tensor product is merely an object of $\calC$: it does not inherit any action of the algebra $A$. We can remedy this situation by assuming that $M$ and $N$ are {\em bimodules} for the
algebra $A$: in this case, we can use the right module structure on $M$ and the left module structure on $N$ to define the tensor product $M \otimes_{A} N$, and this tensor product inherits the structure of a bimodule using the left module structure on $M$ and the right module structure on $N$. In good cases, this definition endows the category of $A$-bimodules in $\calC$ with the structure of a monoidal
category (which is not symmetric in general). 

We can summarize the above discussion as follows: if $A$ is an object of $\calC$ equipped with some algebraic structure (in this case, the structure of either a commutative or an associative algebra), then
we can construct a new category of $A$-modules (in the associative case, we consider
bimodules rather than left or right modules). This category is equipped with the same sort of algebraic structure as the original algebra $A$ (when $A$ is a commutative algebra, the category of $A$-modules inherits a tensor product which is commutative and associative up to isomorphism; in the
case where $A$ is associative, the tensor product of bimodules $\otimes_{A}$ is merely associative up to isomorphism).

Our goal in this section is to obtain an analogous picture in the case where $\calC$ is a symmetric monoidal $\infty$-category, and we consider algebras of a more general nature: namely,
$\calO$-algebra objects of $\calC$ for an arbitrary $\infty$-operad $\calO^{\otimes}$.
We have seen that the collection of such algebras can itself be organized into an
$\infty$-category $\Alg_{\calO}(\calC)$. If $\calO$ is unital, then we can associate
to every object $A \in \Alg_{\calO}(\calC)$ a new $\infty$-category
$\Mod^{\calO}_{A}( \calC)$, whose objects we will refer to as {\it $A$-module objects of
$\calC$}. The exact sense in which an object $M \in \Mod^{\calO}_{A}(\calC)$ depends strongly
on the $\infty$-operad $\calO^{\otimes}$: if $\calO^{\otimes}$ is the commutative $\infty$-operad,
we recover the usual theory of (left) $A$-modules (see \S \ref{modcom}), while
if $\calO^{\otimes}$ is the associative $\infty$-operad we will obtain instead a theory of
bimodules (see \S \ref{bimid}).

There is very little we can say about the $\infty$-categories $\Mod^{\calO}_{A}(\calC)$ of
module objects in the case where $\calO^{\otimes}$ is a general unital $\infty$-operad.
We have therefore adopted to work instead in the more restrictive setting of
{\em coherent} $\infty$-operads. We will introduce the relevant definitions in \S \ref{moduldef}, and
show that they give rise to fibration of $\infty$-operads $\Mod^{\calO}_{A}(\calC)^{\otimes} \rightarrow \calO^{\otimes}$. In \S \ref{almod}, we will study this $\infty$-operad fibration and obtain
a precise description of the $\calO$-algebra objects of $\Mod^{\calO}_{A}(\calC)$ (these
can be identified with $\calO$-algebra objects $B$ of $\calC$ equipped with a map $A \rightarrow
B$: see Corollary \ref{skoke2}). 

Our other major goal in this section is to establish criteria which guarantee the existence of
limits and colimits in $\infty$-categories of the form $\Mod_{A}^{\calO}(\calC)$. The case
of limits is fairly straightforward: in \S \ref{com34}, we will show that under some mild hypotheses, these limits exist and can be computed in the underlying $\infty$-category $\calC$ (Corollary \ref{splasher}).
There is an analogous statement for colimits in \S \ref{com35}, which we will prove in \S \ref{com35}. However, both the statement and the proof are considerably more involved: we must assume not only that the relevant colimits exist in the underlying $\infty$-category $\calC$, but that they are
{\em operadic} colimits in the sense of \S \ref{comm35}.

\subsection{Coherent $\infty$-Operads and Modules}\label{moduldef}

Let $\calC^{\otimes} \rightarrow \calO^{\otimes}$ be a fibration of $\infty$-operads
and let $A \in \Alg_{\calO}(\calC)$. Our goal in this section is to define the simplicial set
$\Mod^{\calO}_{A}(\calC)$ of $A$-module objects of $\calC$. We will then introduce
the notion of a {\em coherent $\infty$-operad} (Definition \ref{koopa}) and show that
if $\calO^{\otimes}$ is coherent, then the construction of $\Mod^{\calO}_{A}(\calC)$ is well-behaved
(Theorem \ref{qwise}). We begin by establishing some terminology.

\begin{definition}
We will say that a morphism $\alpha: \seg{m} \rightarrow \seg{n}$ in $\Nerve(\FinSeg)$ is
{\it semi-inert} if $\alpha^{-1} \{i \}$ has at most one element, for each $i \in \nostar{n}$.
We will say that $\alpha$ is {\it null} if $\alpha$ carries $\seg{m}$ to the base point of
$\seg{n}$.

Let $p: \calO^{\otimes} \rightarrow \Nerve(\FinSeg)$ be an $\infty$-operad, and let $f: X \rightarrow Y$ be a morphism in $\calO^{\otimes}$. We will say that $f$ is {\it semi-inert} if the following conditions
are satisfied:
\begin{itemize}
\item[$(1)$] The image $p(f)$ is a semi-inert morphism in $\Nerve(\FinSeg)$.
\item[$(2)$] For every inert morphism $g: Y \rightarrow Z$ in $\calO^{\otimes}$, if
$p(g \circ f)$ is an inert morphism in $\Nerve(\FinSeg)$, then $g \circ f$ is an inert morphism
in $\calO^{\otimes}$.
\end{itemize}
We will say that $f$ is {\it null} if its image in $\Nerve(\FinSeg)$ is null.
\end{definition}

\begin{remark}\label{ijl}
Let $\calO^{\otimes}$ be an $\infty$-operad. We can think of $\calO^{\otimes}$
as an $\infty$-category whose objects are finite sequences of objects $(X_1, \ldots, X_m)$
of $\calO$. Let $f: (X_1, \ldots, X_m) \rightarrow (Y_1, \ldots, Y_n)$ be a morphism in
$\calO^{\otimes}$, corresponding to a map $\alpha: \seg{m} \rightarrow \seg{n}$ and a collection of maps
$$ \phi_j \in \Mult_{\calO}( \{ X_i \}_{ \alpha(i) = j}, Y_j ).$$
Then $f$ is semi-inert if and only for every $j \in \nostar{n}$, exactly one of the following conditions holds:
\begin{itemize}
\item[$(a)$] The set $\alpha^{-1} \{j\}$ is empty.
\item[$(b)$] The set $\alpha^{-1} \{j\}$ contains exactly one element $i$, and the map
$\phi_i$ is an equivalence between $X_i$ and $Y_j$ in the $\infty$-category $\calO$.
\end{itemize}
The morphism $f$ is inert if $(b)$ holds for each $j \in \nostar{n}$, and null if
$(a)$ holds for each $j \in \nostar{n}$. It follows that every null morphism in $\calO^{\otimes}$ is semi-inert.

Note that $\calO^{\otimes}$ is a unital $\infty$-operad, then the maps $\phi_{j}$ are unique up to
homotopy (in fact, up to a contractible space of choices) when $\alpha^{-1} \{j\}$ is empty.
\end{remark}

\begin{remark}\label{sinkus}
Let $\calO^{\otimes}$ be an $\infty$-operad, and suppose we are given a commutative diagram
$$ \xymatrix{ & Y \ar[dr]^{g} & \\
X \ar[ur]^{f} \ar[rr]^{h} & & Z }$$
in $\calO^{\otimes}$, where $f$ is inert. Then $g$ is semi-inert if and only if
$h$ is semi-inert; this follows immediately from Remark \ref{ijl}.
\end{remark}

\begin{notation}\label{salvi}
Let $\calO^{\otimes}$ be a unital $\infty$-operad. We let $\calK_{\calO}$ denote
the full subcategory of $\Fun( \Delta^1, \calO^{\otimes})$ spanned by the semi-inert morphisms.
For $i \in \{0,1\}$ we let $e_i: \calK_{\calO} \rightarrow \calO^{\otimes}$ be
the map given by evaluation on $i$. We will say that a morphism in $\calK_{\calO}$ is
{\it inert} if its images under $e_0$ and $e_1$ are inert morphisms in $\calO^{\otimes}$.
\end{notation}

\begin{construction}
Let $\calO^{\otimes}$ be an $\infty$-operad and let
$\calC^{\otimes} \rightarrow \calO^{\otimes}_1$ be a $\calO$-operad family.
We define a simplicial set $\tildeMod^{\calO}(\calC)^{\otimes}$ equipped
with a map $\tildeMod^{\calO}(\calC)^{\otimes} \rightarrow \calO^{\otimes}$ so that the following
universal property is satisfied: for every map of simplicial sets $X \rightarrow \calO^{\otimes}$,
there is a canonical bijection
$$ \Fun_{\calO^{\otimes}}( X, \tildeMod^{\calO}(\calC)^{\otimes})
\simeq \Fun_{ \Fun( \{1\}, \calO^{\otimes})}( X \times_{\Fun( \{0\}, \calO^{\otimes})}
\calK_{\calO}, \calC^{\otimes}).$$
We let $\MMod^{\calO}(\calC)^{\otimes}$ denote the full simplicial subset of
$\tildeMod^{\calO}(\calC)^{\otimes}$ spanned by those vertices $\overline{X}$ lying
over an object $X \in \calO^{\otimes}$ which classify functors
$$ \calK_{\calO} \times_{ \calO^{\otimes} } \{X\} \rightarrow \calC^{\otimes}$$
which carry inert morphisms to inert morphisms.
\end{construction}

\begin{notation}\label{splay}
Let $\calO^{\otimes}$ be an $\infty$-operad. We let $\calK^{0}_{\calO}$ denote the full subcategory
of $\calK_{\calO}$ spanned by the null morphisms $f: X \rightarrow Y$ in $\calO^{\otimes}$.
For $i = 0,1$, we let $e_i^{0}$ denote the restriction to $\calK^{0}_{\calO}$ of the evaluation map
$e_{i}: \calK_{\calO} \rightarrow \calO^{\otimes}$.
\end{notation}

\begin{notation}
Let $\calO^{\otimes}$ be an $\infty$-operad and let $\calC^{\otimes} \rightarrow \calO^{\otimes}$ be a $\calO$-operad family. We define a simplicial set $\tildePAlg_{\calO}(\calC)$ equipped
with a map $\tildePAlg_{\calO}(\calC) \rightarrow \calO^{\otimes}_{0}$ so that the following universal property is satisfied: for every simplicial set $X$ with a map $X \rightarrow \calO^{\otimes}_0$, we have a canonical bijection
$$ \Fun_{ \calO^{\otimes}}( X, \tildePAlg_{\calO}(\calC) ) \simeq
\bHom_{\Fun( \{1\}, \calO^{\otimes})}( X \times_{ \Fun( \{0\}, \calO^{\otimes})} \calK^{0}_{\calO}, \calC^{\otimes}).$$

We observe that for every object $Y \in \calO^{\otimes}$, a vertex of the fiber
$\tildePAlg_{\calO}(\calC) \times_{ \calO^{\otimes} } \{Y\}$ can be identified with
a functor $F: \{Y\} \times_{ \Fun( \{0\}, \calO^{\otimes})} \calK^{0}_{\calO} \rightarrow \calC^{\otimes}$.
We let $\PrAlg_{\calO}(\calC)$ denote the full simplicial subset of
$\tildePAlg_{\calO}(\calC)$ spanned by those vertices for which the
functor $F$ preserves inert morphisms.
\end{notation}

\begin{remark}\label{user}
Let $\calO^{\otimes}$ be a unital $\infty$-operad and $\calC^{\otimes} \rightarrow \calO^{\otimes}$
a $\calO$-operad family. The map $\theta: \calK_{\calO}^{0} \rightarrow \calO^{\otimes} \times \calO^{\otimes}$ of Lemma \ref{preswope} induces a map
$\phi: \calO^{\otimes} \times \Fun_{ \calO^{\otimes}}( \calO^{\otimes}, \calC^{\otimes})
\rightarrow \tildePAlg_{\calO}(\calC)$. 
Unwinding the definitions, we see that $\phi^{-1} \PrAlg_{\calO}(\calC)$ can be identified with
the product $\calO^{\otimes} \times \Alg_{\calO}(\calC)$. 
\end{remark}

\begin{definition}\label{forse}
Let $\calO^{\otimes}$ be a $\infty$-operad and $\calC^{\otimes} \rightarrow \calO^{\otimes}$ a $\calO^{\otimes}$-operad family. We let $\Mod^{\calO}(\calC)^{\otimes}$ denote the fiber product
$$\MMod^{\calO}(\calC)^{\otimes} \times_{ \PrAlg_{\calO}(\calC)} (\calO^{\otimes} \times \Alg_{\calO}(\calC)).$$
For every algebra object $A \in \Alg_{\calO}(\calC)$, we let $\Mod^{\calO}_{A}(\calC)^{\otimes}$
denote the fiber product
$$ \MMod^{\calO}(\calC)^{\otimes} \times_{ \PrAlg_{\calO}(\calC)} \{A\}
\simeq \Mod^{\calO}(\calC)^{\otimes} \times_{ \Alg_{\calO}(\calC)} \{A\}.$$
\end{definition}

For a general fibration of $\infty$-operads $\calC^{\otimes} \rightarrow \calO^{\otimes}$, the
simplicial set $\Mod^{\calO}(\calC)^{\otimes}$ defined above need not be an $\infty$-category.
To guarantee that our constructions are well-behaved, we need to make an assumption on $\calO^{\otimes}$.

\begin{definition}\label{koopa}
We will say that an $\infty$-operad $\calO^{\otimes}$ is
{\it coherent} if the following conditions are satisfied:
\begin{itemize}
\item[$(1)$] The $\infty$-operad $\calO^{\otimes}$ is unital.
\item[$(2)$] The evaluation map $e_0: \calK_{\calO} \rightarrow \calO^{\otimes}$
is a flat categorical fibration.
\end{itemize}
\end{definition}

\begin{example}
Let $\calO^{\otimes}$ be the $\infty$-operad $\OpE{0}$ of Example \ref{xe1.5}. 
Then every morphism in $\calO^{\otimes}$ is semi-inert, so that
$\calK_{\calO} = \Fun( \Delta^1, \calO^{\otimes})$. The evaluation functor
$e_0: \Fun( \Delta^1, \calO^{\otimes}) \rightarrow \calO^{\otimes}$ is
a Cartesian fibration and therefore a flat categorical fibration (Example \bicatref{cartflat}), so
that $\OpE{0}$ is a coherent $\infty$-operad.
\end{example}

\begin{example}
We will show later that the commutative and associative $\infty$-operads are coherent
(Propositions \ref{asscoh} and \ref{commcoh}). 
\end{example}

We can now state the main result of this section as follows:

\begin{theorem}\label{qwise}
Let $p: \calC^{\otimes} \rightarrow \calO^{\otimes}$ be a fibration of $\infty$-operads,
where $\calO^{\otimes}$ is coherent, and let $A \in \Alg_{\calO}(\calC)$. Then the induced map
$\Mod^{\calO}_{A}(\calC)^{\otimes} \rightarrow \calO^{\otimes}$ is a fibration of $\infty$-operads.
\end{theorem}

The main ingredient in the proof of Theorem \ref{qwise} is the following
result, which we will prove at the end of this section.

\begin{proposition}\label{cammite}
Let $\calO^{\otimes}$ be a coherent $\infty$-operad, and
$\calC^{\otimes} \rightarrow \calO^{\otimes}$ be a $\calO$-operad family. Then
the induced map $q: \MMod^{\calO}(\calC)^{\otimes} \rightarrow \calO^{\otimes}$
is again a $\calO$-operad family. Moreover, a morphism
$f$ in $\MMod^{\calO}(\calC)^{\otimes}$ is inert if and only if it satisfies the following
conditions:
\begin{itemize}
\item[$(1)$] The morphism $f_0 = q(f): X \rightarrow Y$ is inert in $\calO^{\otimes}$.
\item[$(2)$] Let $F: \calK_{\calO} \times_{ \calO^{\otimes}_{0}} \Delta^1 \rightarrow
\calC^{\otimes}$ be the functor classified by $f$. For every
$e_0$-Cartesian morphism $\widetilde{f}$ of $\calK_{\calO}$ lifting
$f_0$ (in other words, for every lift $\widetilde{f}$ of $f_0$ such that
$e_1( \widetilde{f})$ is an equivalence: see part $(4)$ of the proof of
Proposition \ref{sapless}), the morphism $F( \widetilde{f} )$ of
$\calC^{\otimes}$ is inert.
\end{itemize}
\end{proposition}

We will obtain Theorem \ref{qwise} by combining Proposition \ref{cammite} with
a few elementary observations.

\begin{lemma}\label{kaj}
Let $p: \calO^{\otimes}$ be a unital $\infty$-operad. Let $X$ and $Y$ be objects of
$\calO^{\otimes}$ lying over objects $\seg{m}, \seg{n} \in \FinSeg$, and let $\bHom_{\calO^{\otimes}}^{0}(X,Y)$ denote the summand of $\bHom_{\calO^{\otimes}}(X,Y)$ corresponding to those maps which cover the null morphism $\seg{m} \rightarrow \{ \ast \} \rightarrow \seg{n}$ in $\FinSeg$.
Then $\bHom_{ \calO^{\otimes} }^{0}(X, Y)$ is contractible.
\end{lemma}

\begin{proof}
Choose an inert morphism $f: X \rightarrow Z$ covering the map $\seg{m} \rightarrow \{ \ast \}$ in
$\FinSeg$. We observe that composition with $f$ induces a homotopy equivalence
$$ \bHom_{\calO^{\otimes}}(Z, Y) \rightarrow \bHom_{\calO^{\otimes}}^{0}( X, Y).$$
Since $\calO^{\otimes}$ is unital, the space $\bHom_{\calO^{\otimes}}(Z, Y)$ is contractible.
\end{proof}

\begin{lemma}\label{jinlo}
Let $\calO^{\otimes}$ be a unital $\infty$-operad and let
$\calC^{\otimes} \rightarrow \calO^{\otimes}$ be a $\calO$-operad family. Then the map
$\tildePAlg_{\calO}(\calC) \rightarrow \calO^{\otimes}_{0}$ is a categorical fibration.
In particular, $\tildePAlg_{\calO}(\calC)$ is an $\infty$-category.
\end{lemma}

\begin{proof}
The simplicial set $\tildePAlg_{\calO}(\calC)$ can be described as
$(e_0^{0})_{\ast}( \calK_{\calO}^{0} \times_{ \Fun( \{1\}, \calO^{\otimes})} \calC^{\otimes})$.
The desired result now follows from Proposition \bicatref{saeve} and Corollary \ref{swope}.
\end{proof}

\begin{remark}
In the situation of Lemma \ref{jinlo}, the simplicial set $\PrAlg_{\calO}(\calC)$ is a full
simplicial subset of an $\infty$-category, and therefore also an $\infty$-category. Since
$\PrAlg_{\calO}(\calC)$ is evidently stable under equivalence in
$\tildePAlg_{\calO}(\calC)$, we conclude that the map
$\PrAlg_{\calO}(\calC) \rightarrow \calO^{\otimes}$ is also a categorical fibration.
\end{remark}

\begin{lemma}\label{preswope}
Let $\calO^{\otimes}$ be a unital $\infty$-operad. Then the maps $e_0^{0}$ and
$e_{1}^{0}$ determine a trivial Kan fibration
$\theta: \calK_{\calO}^{0} \rightarrow \calO^{\otimes} \times \calO^{\otimes}$.
\end{lemma}

\begin{proof}
Since $\theta$ is evidently a categorical fibration, it will suffice to show that $\theta$ is a categorical equivalence. 
Corollary \toposref{tweezegork} implies that evaluation at $\{0\}$ induces a 
Cartesian fibration $p: \Fun( \Delta^1, \calO^{\otimes}) \rightarrow \calO^{\otimes}$.
If $f$ is a null morphism in $\calO^{\otimes}$, then so is $f \circ g$ for every morphism
$g$ in $\calO^{\otimes}$; it follows that if $\alpha: f \rightarrow f'$ is a $p$-Cartesian
morphism in $\Fun( \Delta^1, \calO^{\otimes})$ and $f' \in \calK^{0}_{\calO}$, then
$f \in \calK^{0}_{\calO}$. Consequently, $p$ restricts to a Cartesian fibration
$e_0^{0}: \calK^{0}_{\calO} \rightarrow \calO^{\otimes}$. Moreover, a morphism
$\alpha$ in $\calK^{0}_{\calO}$ is $e_0^{0}$-Cartesian if and only if $e_1^{0}(\alpha)$ is an equivalence. Consequently, $\theta$ fits into a commutative diagram
$$ \xymatrix{ \calK^{0}_{\calO} \ar[rr]^{\theta} \ar[dr]^{ e_0^{0}} & & \calO^{\otimes} \times \calO^{\otimes} \ar[dl]^{\pi} \\
& \calO^{\otimes} & }$$
and carries $e_0^{0}$-Cartesian morphisms to $\pi$-Cartesian morphisms. According to
Corollary \toposref{usefir}, it will suffice to show that for every object $X \in \calO^{\otimes}$, the
map $\theta$ induces a categorical equivalence $\theta_{X}: ( \calO^{\otimes})^{X/}_{0} \rightarrow \calO^{\otimes}$, where $(\calO^{\otimes})^{X/}_{0}$ is the full subcategory of
$(\calO^{\otimes})^{X/}$ spanned by the null morphisms $X \rightarrow Y$.

The map $\theta_{X}$ is obtained by restricting the left fibration $q: (\calO^{\otimes})^{X/} \rightarrow \calO^{\otimes}$. We observe that if $f$ is a null morphism in $\calO^{\otimes}$, then
so is every composition of the form $g \circ f$. It follows that if $\alpha: f \rightarrow f'$ is
a morphism in $(\calO^{\otimes})^{X/}$ such that $f \in (\calO^{\otimes})^{X/}_0$, then
$f' \in (\calO^{\otimes})^{X/}_0$, so the map $q$ restricts to a left fibration
$\theta_{X}: (\calO^{\otimes})^{X/}_{0} \rightarrow \calO^{\otimes}$. Since the fibers
of $\theta_{X}$ are contractible (Lemma \ref{kaj}), Lemma \toposref{toothie} implies that $\theta_{X}$ is a trivial Kan fibration.
\end{proof}

\begin{corollary}\label{swope}
Let $\calO^{\otimes}$ be a unital $\infty$-operad.
The projection map $e^{0}_{0}: \calK^{0}_{\calO} \rightarrow \calO^{\otimes}$ is a flat categorical fibration.
\end{corollary}

\begin{remark}\label{user2}
Let $\calO^{\otimes}$ be a unital $\infty$-operad.
Since the map $e_{1}^{0}: \calK_{\calO}^{0} \rightarrow \calO^{\otimes}$
(Lemma \ref{preswope}), we conclude that 
$\phi: \calO^{\otimes} \times \Fun_{ \calO^{\otimes}}( \calO^{\otimes}, \calC^{\otimes})
\rightarrow \tildePAlg_{\calO}(\calC)$ and induces a categorical equivalence
$\calO^{\otimes} \times \Alg_{\calO}(\calC) \rightarrow \PrAlg_{\calO}(\calC).$
\end{remark}

\begin{lemma}\label{sqeak}
Let $\calO^{\otimes}$ be a coherent $\infty$-operad, and let $\calC^{\otimes} \rightarrow \calO^{\otimes}$ be a $\calO$-operad family. Then the inclusion $\calK_{\calO}^{0} \rightarrow \calK_{\calO}$
induces a categorical fibration
$\MMod^{\calO}( \calC)^{\otimes} \rightarrow \PrAlg_{\calO}(\calC)$.
\end{lemma}

\begin{proof}
Since $\MMod^{\calO}(\calC)^{\otimes}$ and $\PrAlg_{\calO}(\calC)$ are full subcategories
stable under equivalence in $\widetilde{ \MMod}^{\calO}(\calC)^{\otimes}$ and
$\tildePAlg_{\calO}(\calC)$, it suffices to show that the map
$\widetilde{ \MMod}^{\calO}(\calC)^{\otimes} \rightarrow \tildePAlg_{\calO}(\calC)$
has the right lifting property with respect to every trivial cofibration $A \rightarrow B$
with respect to the Joyal model structure on $(\sSet)_{/ \calO^{\otimes}}$. Unwinding
the definitions, we are required to provide solutions to lifting problems of the form
$$ \xymatrix{ (A \times_{ \Fun( \{0\}, \calO^{\otimes})} \calK_{\calO})
\coprod_{ A \times_{ \Fun( \{0\}, \calO^{\otimes})} \calK_{\calO}^{0} } ( B \times_{ \Fun( \{0\}, \calO^{\otimes}) }\calK^{0}_{\calO}) \ar[d]^{i} \ar[r] & \calC^{\otimes} \ar[d]^{q} \\
B \times_{ \Fun( \{0\}, \calO^{\otimes})} \calK_{\calO} \ar[r] & \Fun( \{1\}, \calO^{\otimes}). }$$
Since $q$ is a categorical fibration, it suffices to prove that the monomorphism $i$ is a categorical equivalence of simplicial sets. In other words, we need to show that the diagram
$$ \xymatrix{ A \times_{ \Fun( \{0\}, \calO^{\otimes})} \calK_{\calO}^{0} \ar[r] \ar[d] & A \times_{ \Fun( \{0\}, \calO^{\otimes})} \calK_{\calO} \ar[d] \\
B \times_{ \Fun( \{0\}, \calO^{\otimes})} \calK_{\calO}^{0} \ar[r] & B \times_{ \Fun( \{0\}, \calO^{\otimes})} \calK_{\calO} }$$
is a homotopy pushout square (with respect to the Joyal model structure). To prove this, it suffices
to show that the vertical maps are categorical equivalences. Since $A \rightarrow B$ is a
categorical equivalence, this follows from Corollary \bicatref{silman}, since
the restriction maps $e_0^{0}: \calK_{\calO}^{0} \rightarrow \calO^{\otimes}$ and
$e_0: \calK_{\calO} \rightarrow \calO^{\otimes}_0$ are flat categorical fibrations
(Corollary \ref{swope}).
\end{proof}

\begin{remark}\label{saef}
Let $\calC^{\otimes} \rightarrow \calO^{\otimes}$ be as in Lemma \ref{sqeak}. Since
every morphism $f: X \rightarrow Y$ in $\calO^{\otimes}$ with $X \in \calO^{\otimes}_{\seg{0}}$ is automatically null, the categorical fibration
$\MMod^{\calO}( \calC)^{\otimes} \rightarrow \PrAlg_{\calO}(\calC)$ induces an isomorphism
$\MMod^{\calO}(\calC)^{\otimes}_{\seg{0}} \rightarrow \PrAlg_{\calO}(\calC)_{\seg{0}}$.
\end{remark}

\begin{remark}\label{capp}
In the situation of Definition \ref{forse}, consider the pullback diagram
$$ \xymatrix{ \Mod^{\calO}(\calC)^{\otimes} \ar[r]^{j} \ar[d]^{\theta'} & \MMod^{\calO}(\calC)^{\otimes} \ar[d]^{\theta} \\
\calO^{\otimes} \times \Alg_{\calO}(\calC) \ar[r]^{j'} & \PrAlg_{\calO}(\calC). }$$
Since $\theta$ is a categorical fibration and $\calO^{\otimes} \times \Alg_{\calO}(\calC)$ is an $\infty$-category, this diagram is a homotopy pullback square (with respect to the Joyal model structure).
Remark \ref{user2} implies that $j'$ is a categorical equivalence, so $j$ is also a categorical
equivalence. Using Proposition \ref{cammite}, we deduce that
the map $\Mod^{\calO}(\calC)^{\otimes} \rightarrow \calO^{\otimes}$ exhibits
$\Mod^{\calO}(\calC)^{\otimes}$ as a $\calO$-operad family. Remark 
\ref{saef} implies that $\theta'$ induces an isomorphism
$\Mod^{\calO}(\calC)^{\otimes}_{\seg{0}} \rightarrow \Alg_{\calO}(\calC)$.
\end{remark}

\begin{proof}[Proof of Theorem \ref{qwise}]
Combine Proposition \ref{supercae} with Remark \ref{saef}.
\end{proof}

We now turn to the proof of Proposition \ref{cammite}.

\begin{lemma}\label{susker}
Let $\calO^{\otimes}$ be an $\infty$-operad, let $\overline{X}: X \rightarrow X'$
be an object of $\Fun( \Delta^1, \calO^{\otimes})$, and let $e_0: \Fun( \Delta^1, \calO^{\otimes}) \rightarrow \calO^{\otimes}$ be given by evaluation at $0$.
\begin{itemize}
\item[$(1)$] For every inert morphism $f_0: X \rightarrow Y$ in $\calO^{\otimes}$, there exists an $e_0$-coCartesian morphism $f: \overline{X} \rightarrow \overline{Y}$ in $\Fun(\Delta^1, \calO^{\otimes})$ lifting $f_0$. 

\item[$(2)$] An arbitrary morphism $f: \overline{X} \rightarrow \overline{Y}$ in $\Fun(\Delta^1,\calO^{\otimes})$ lifting
$f_0$ is $e_0$-coCartesian if and only if the following conditions are satisfied:
\begin{itemize}
\item[$(i)$] The image of $f$ in $\Fun( \{1\}, \calO^{\otimes})$ is inert.
\item[$(ii)$] Let $\alpha: \seg{m} \rightarrow \seg{m'}$ be the morphism in
$\FinSeg$ determined by $\overline{X}$, and let $\beta: \seg{n} \rightarrow \seg{n'}$ be the morphism in $\FinSeg$ determined by $\overline{Y}$. Then $f$ induces a diagram
$$ \xymatrix{ \seg{m} \ar[r]^{\delta} \ar[d]^{\alpha} & \seg{n} \ar[d]^{\beta} \\
\seg{m'} \ar[r]^{\gamma} & \seg{n'} }$$
with the property that $\gamma^{-1} \{ \ast \} = \alpha( \delta^{-1} \{\ast \} )$. 
\end{itemize}

\item[$(3)$] If $\overline{X} \in \calK_{\calO}$ and $f$ is $e_0$-coCarteian, then $\overline{Y} \in \calK_{\calO}$ (the map $f$ is then automatically coCartesian with respect to the restriction $e_0 | \calK_{\calO}$).
\end{itemize}
\end{lemma}

\begin{proof}
We note that the ``only if'' direction of $(2)$ follows from the ``if'' direction together with $(1)$,
since a coCartesian lift $f$ of $f_0$ is determined uniquely up to homotopy.
We first treat the case where $\calO^{\otimes} = \Nerve(\FinSeg)$. 
Then we can identify $\overline{X}$ with $\alpha$ with $\beta$, and
$f_0$ with an inert morphism $\delta: \seg{m} \rightarrow \seg{n}$.
Choose a map $\gamma: \seg{m'} \rightarrow \seg{n'}$ which identifies
$\seg{n'}$ with the finite pointed set obtained from $\seg{m'}$ by collapsing
$\alpha( \delta^{-1} \{ \ast \})$ to a point, so that we have a commutative diagram
$$ \xymatrix{ \seg{m} \ar[r]^{\delta} \ar[d]^{\alpha} & \seg{n} \ar[d]^{\beta} \\
\seg{m'} \ar[r]^{\gamma} & \seg{n'} }$$
as in the statement of $(ii)$. By construction, the map $\gamma$ is inert.
The above diagram
is evidently a pushout square in $\FinSeg$ and so therefore corresponds to a morphism
morphism $f$ in $\Fun( \Delta^1, \Nerve(\FinSeg))$ which is coCartesian with respect to the projection $\Fun( \Delta^1, \Nerve(\FinSeg))$, and therefore with respect to the projection
$e_0: \calK_{ \Nerve(\FinSeg)} \rightarrow \Nerve(\FinSeg)$. This completes the proof
of $(1)$ and the ``if'' direction of $(2)$; assertion $(3)$ follows from the observation that
$\beta$ is semi-inert whenever $\alpha$ is semi-inert.

We now prove the ``if'' direction of $(2)$ in the general case. Since $f_0$ is inert, 
we conclude (Proposition \toposref{stuch}) that $f$ is $e_0$-coCartesian if and only if it is coCartesian with respect to the composition
$$ \Fun(\Delta^1, \calO^{\otimes}) \stackrel{e_0}{\rightarrow} \calO^{\otimes} \rightarrow \Nerve(\FinSeg).$$
This map admits another factorization
$$ \Fun( \Delta^1, \calO^{\otimes}) \stackrel{p}{\rightarrow} \Fun( \Delta^1, \Nerve(\FinSeg)) \stackrel{q}{\rightarrow} \Nerve(\FinSeg).$$
The first part of the proof shows that $p(f)$ is $q$-coCartesian. Since $f_0$ is inert,
condition $(i)$ and Lemma \monoidref{surtybove} guarantee that $f$ is $p$-coCartesian.
Applying Proposition \toposref{stuch}, we deduce that $f$ is $(q \circ p)$-coCartesian as desired.

To prove $(1)$, we first choose a diagram
$$ \xymatrix{ \seg{m} \ar[r]^{\delta} \ar[d]^{\alpha} & \seg{n} \ar[d]^{\beta} \\
\seg{m'} \ar[r]^{\gamma} & \seg{n'} }$$
satisfying the hypotheses of $(2)$. We can identify this diagram with a morphism
$\overline{f}$ in $\Fun( \Delta^1, \Nerve(\FinSeg))$. We will show that it is possible to
choose a map $f: \overline{X} \rightarrow \overline{Y}$ satisfying $(i)$ and $(ii)$, which lifts
both $f_0$ and $\overline{f}$. Using the assumption that $\calO^{\otimes}$ is
an $\infty$-operad, we can choose an inert morphism $f_1: X' \rightarrow Y'$
lifting $\gamma$. Using the fact that the projection $\calO^{\otimes} \rightarrow \Nerve(\FinSeg)$ is an inner fibration, we obtain a commutative diagram
$$ \xymatrix{ X \ar[d]^{ \overline{X} } \ar[r]^{f_0} \ar[dr] & Y \ar@{-->}[d] \\
X' \ar[r]^{f_1} & Y'} $$
in $\calO^{\otimes}$. To complete the proof, it suffices to show that we can complete the diagram
by filling an appropriate horn $\Lambda^2_0 \subseteq \Delta^2$ to supply the dotted arrow.
The existence of this arrow follows from the fact that $f_0$ is inert, and therefore coCartesian with respect to the projection $\calO^{\otimes} \rightarrow \Nerve(\FinSeg)$.

We now prove $(3)$. Let $f: \overline{X} \rightarrow \overline{Y}$ be a morphism corresponding
to a diagram as above, where $\overline{X}$ is semi-inert; we wish to show that $\overline{Y}$ is semi-inert. The first part of the proof shows that the image of $\overline{Y}$ in $\Nerve(\FinSeg)$ is semi-inert. It will therefore suffice to show that for any inert morphism $g: Y' \rightarrow Z$, if the image
of $g \circ \overline{Y}$ in $\Nerve(\FinSeg)$ is inert, then $g \circ \overline{Y}$ is inert.
We note that $g \circ \overline{Y} \circ f_0 \simeq g \circ f_1 \circ \overline{X}$ has inert image
in $\Nerve(\FinSeg)$; since $\overline{X}$ is semi-inert we deduce that
$g \circ \overline{Y} \circ f_0$ is inert. Since $f_0$ is inert, we deduce from Proposition
\toposref{protohermes} that $g \circ \overline{Y}$ is inert as desired.
\end{proof}

\begin{proposition}\label{sapless}
Let $\calO^{\otimes}$ be a coherent $\infty$-operad, let $M_0$ denote the collection of all inert morphisms in $\calO^{\otimes}$, and let $M$ denote the collection of all inert morphisms in $\calK_{\calO}$. 
Then the construction
$$ \overline{X} \mapsto \overline{X} \times_{ (\Fun( \{0\}, \calO^{\otimes}), M_0)} (\calK_{\calO},M)$$
determines a left Quillen functor from $\mset{\CatP^{\famm}_{\Fun( \{0\}, \calO^{\otimes})}}$
to $\mset{\CatP^{\famm}_{\Fun( \{1\}, \calO^{\otimes})}}$. 
\end{proposition}

\begin{proof}
It will suffice to show that the map $\calK_{\calO} \rightarrow \calO^{\otimes} \times \calO^{\otimes}$ satisfies the hypotheses of Theorem \bicatref{makera}:
\begin{itemize}
\item[$(1)$] The evaluation map $e_0: \calK_{\calO} \rightarrow \calO^{\otimes}$ is
a flat categorical fibration. This follows from our assumption that $\calO^{\otimes}$ is coherent.

\item[$(2)$] The collections of inert morphisms in $\calO^{\otimes}$ and $\calK_{\calO}$ contain
all equivalences and are closed under composition. This assertion is clear from the definitions.

\item[$(3)$] For every $2$-simplex $\sigma$ of $\calK_{\calO}$, if $e_0( \sigma)$ is thin,
then $e_1(\sigma)$ thin. Moreover, a $2$-simplex $\Delta^2 \rightarrow \calK_{\calO}$ is thin
if its restriction to $\Delta^{ \{0,1\} }$ is thin. These assertions are clear, since the categorical pattern $\CatP^{\famm}_{\calO^{\otimes}}$ designates every $2$-simplex as thin.

\item[$(4)$] For every inert edge $\Delta^1 \rightarrow \Fun( \{0\}, \calO^{\otimes})$, the induced map
$\calK_{\calO} \times_{ \Fun( \{0\}, \calO^{\otimes}) } \Delta^1 \rightarrow \Delta^1$ is a Cartesian
fibration. To prove this, it suffices to show that if $\overline{X}$ is an object of
$\calK_{\calO}$ corresponding to a semi-inert morphism $X \rightarrow X'$ in
$\calO^{\otimes}$, and $f: Y \rightarrow X$ is an inert morphism in $\calO^{\otimes}$, then
we can lift $f$ to an $e_0$-Cartesian morphism $\overline{f}: \overline{Y} \rightarrow \overline{X}$.
According to Corollary \toposref{tweezegork}, we can lift $f$ to a
morphism $\overline{f}: \overline{Y} \rightarrow \overline{X}$ in $\Fun( \Delta^1, \calO^{\otimes})$
which is $e'_0$-coCartesian, where $e'_0: \Fun( \Delta^1, \calO^{\otimes}) \rightarrow \calO^{\otimes}_0$ is given by evaluation at $0$. To complete the proof, it suffices
to show that $\overline{Y}$ belongs to $\calK_{\calO}$. We can identify $\overline{f}$ with a commutative diagram
$$ \xymatrix{ Y \ar[r]^{f} \ar[d]^{ \overline{Y} } & X \ar[d]^{\overline{X} } \\
Y' \ar[r]^{f'} & X' }$$
in $\calO^{\otimes}$; we wish to prove that the morphism $\overline{Y}$ is semi-inert.
Since $\overline{f}$ is $e'_0$-Cartesian, Corollary \toposref{tweezegork} implies that
$f'$ is an equivalence. We can therefore identify $\overline{Y}$ with
$f' \circ \overline{Y} \simeq \overline{X} \circ f$, which is inert by virtue of
Remark \ref{sinkus}. This completes the verification of condition $(4)$ and establishes
the following criterion: if $\overline{f}$ is a morphism in $\calK_{\calO}$ such that
$f = e_0( \overline{f})$ is inert, the following conditions are equivalent:
\begin{itemize}
\item[$(i)$] The map $\overline{f}$ is locally $e_0$-Cartesian.
\item[$(ii)$] The map $\overline{f}$ is $e_0$-Cartesian.
\item[$(iii)$] The morphism $e_1( \overline{f} )$ is an equivalence in $\calO^{\otimes}_{1}$.
\end{itemize}

\item[$(5)$] Let $p: \Delta^1 \times \Delta^1 \rightarrow\Fun( \{0\}, \calO^{\otimes})$ be one of the diagrams
specified in the definition of the categorical pattern $\CatP^{\famm}_{\Fun( \{0\}, \calO^{\otimes})}$. Then
the restriction of $e_0$ determines a coCartesian fibration $(\Delta^1 \times \Delta^1) \times_{ \Fun( \{0\}, \calO^{\otimes})} \calK_{\calO} \rightarrow \Delta^1 \times \Delta^1$. This follows immediately from Lemma \ref{susker}.

\item[$(6)$] Let $p: \Delta^1 \times \Delta^1 \rightarrow \Fun( \{0\}, \calO^{\otimes})$ be one of the diagrams
specified in the definition of the categorical pattern $\CatP^{\famm}_{\Fun( \{0\}, \calO^{\otimes})}$ and
let $s$ be a coCartesian section of the projection $(\Delta^1 \times \Delta^1) \times_{ \Fun( \{0\}, \calO^{\otimes})} \calK_{\calO} \rightarrow \Delta^1 \times \Delta^1$. Then the composite map
$$q: \Delta^1 \times \Delta^1 \stackrel{s}{\rightarrow} (\Delta^1 \times \Delta^1) \times_{ \Fun( \{0\}, \calO^{\otimes})} \calK_{\calO} \rightarrow \calK_{\calO} \stackrel{e_1}{\rightarrow} \Fun( \{1\},\calO^{\otimes})$$
is one of the diagrams specified in the definition of the categorical patterm $\CatP^{\famm}_{\Fun( \{1\}, \calO^{\otimes})}$. It follows from Lemma \ref{susker} that
$q$ carries each morphism in $\Delta^1 \times \Delta^1$ to an inert morphism in $\Fun( \{1\}, \calO^{\otimes})$.
Unwinding the definitions (and using the criterion provided by Lemma \ref{susker}), we are reduced to verifying the following simple combinatorial fact: given a semi-inert morphism $\seg{m} \rightarrow \seg{k}$ in
$\FinSeg$ and a commutative diagram of inert morphisms
$$ \xymatrix{ \seg{m} \ar[r] \ar[d] & \seg{n} \ar[d] \\
\seg{m'} \ar[r] & \seg{n'} }$$
which induces a bijection $\nostar{n} \coprod_{ \nostar{n'} } \nostar{m'} \rightarrow \nostar{m}$, the induced diagram
$$ \xymatrix{ \seg{k} \ar[r] \ar[d] & \seg{n} \coprod_{ \seg{m} } \seg{k} \ar[d] \\
\seg{m'} \coprod_{ \seg{m} } \seg{k} \ar[r] & \seg{n'} \coprod_{ \seg{m} } \seg{k} }$$ 
has the same property. 

\item[$(7)$] Suppose we are given a commutative diagram
$$\xymatrix{ & Y \ar[dr]^{g} & \\
X \ar[ur]^{f} \ar[rr]^{h} & & Z }$$
in $\calK_{\calO}$, where $g$ is locally $e_0$-Cartesian, $e_0(g)$ is inert, and 
$e_0(f)$ is an equivalence. We must show that $f$ is inert if and only if $h$ is inert.
Consider the underlying diagram in $\calO^{\otimes}$
$$ \xymatrix{ X_0 \ar[r]^{f_0} \ar[d] & Y_0 \ar[r]^{g_0} \ar[d] & Z_0 \ar[d] \\
X_1 \ar[r]^{f_1} & Y_1 \ar[r]^{g_1} & Z_1. }$$ 
Since $f_0$ is an equivalence and $g_0$ is inert, $h_0 = g_0 \circ f_0$ is inert.
It will therefore suffice to prove that $f_1$ is inert if and only if $h_1 = g_1 \circ f_1$ is inert.
For this, it sufficess to show that $g_1$ is an equivalence in $\calO^{\otimes}$.
This follows from the proof of $(4)$, since $g$ is assumed to be locally
$e_0$-Cartesian.

\item[$(8)$] Suppose we are given a commutative diagram
$$ \xymatrix{ & Y \ar[dr]^{g} & \\
X \ar[ur]^{f} \ar[rr]^{h} & & Z }$$
in $\calK_{\calO}$ where $f$ is $e_0$-coCartesian, $e_0(f)$ is inert, and $e_0(g)$ is an equivalence.
We must show that $g$ is inert if and only if $h$ is inert. Consider the underlying
diagram in $\calO^{\otimes}$
$$ \xymatrix{ X_0 \ar[r]^{f_0} \ar[d] & Y_0 \ar[r]^{g_0} \ar[d] & Z_0 \ar[d] \\
X_1 \ar[r]^{f_1} & Y_1 \ar[r]^{g_1} & Z_1. }$$ 
Since $f_0$ and $f_1$ are inert (Lemma \ref{susker}), Propositions \toposref{swimmm} and \ref{singwell} guarantee that $g_0$ is inert if and only if $h_0$ is inert, and that
$g_1$ is inert if and only if $h_1$ is inert. Combining these facts, we conclude that
$g$ is inert if and only if $h$ is inert.
\end{itemize}
\end{proof}

\begin{proof}[Proof of Proposition \ref{cammite}]
Let $M$ denote the collection of inert morphisms in $\calK_{\calO}$, and $M_0$ the collection of inert morphisms in $\calO^{\otimes}$. It follows
from Proposition \ref{sapless} that the construction
$$\overline{X} \mapsto \overline{X} \times_{ (\calO^{\otimes},M_0)} (\calK_{\calO},M)$$
is a left Quillen functor from $\CatP_{\calO}$ to itself. Let $G$ denote the right
adjoint of this functor, so that $G$ preserves fibrant objects. It follows that
if $\calC^{\otimes}$ is a $\calO$-operad family and $M'$ denotes the collection of
inert morphisms in $\calC^{\otimes}$, then $G( \calC^{\otimes}, M')$ is a fibrant
object of $\CatP_{\calO}$. Unwinding the definitions, we see that
$G( \calC^{\otimes}, M')$ is given by the pair
$(\MMod^{\calO}(\calC)^{\otimes}, M'')$, where $M''$ is the collection of all morphisms
$f$ in $\MMod^{\calO}(\calC)^{\otimes}$ satisfying condition $(1)$ together with the following analogue of $(2)$:
\begin{itemize}
\item[$(2')$] Let $F: \calK_{\calO} \times_{ \Fun( \{0\}, \calO^{\otimes})} \Delta^1 \rightarrow
\calC^{\otimes}$ be the functor classified by $f$. For every
inert morphism $\widetilde{f}$ of $\calK_{\calO}$ lifting
the morphism $e_0(f): X \rightarrow Y$ in $\calO^{\otimes}$, the morphism $F( \widetilde{f} )$ of
$\calC^{\otimes}$ is inert.
\end{itemize}
Since every $e_0$-Cartesian morphism $\widetilde{f}$ of $\calK_{\calO}$ 
such that $e_0( \widetilde{f}) = f$ is inert, it is clear that $(2') \Rightarrow (2)$.
To prove the converse, consider an arbitrary inert morphism
$\widetilde{f}$ lifting $f$, and factor $\widetilde{f}$ as a composition
$\widetilde{f} \simeq \widetilde{f}' \circ \widetilde{f}''$, where $\widetilde{f}'$
is $e_0$-Cartesian (so that $F( \widetilde{f}')$ is inert in $\calC^{\otimes}$ by virtue of $(2)$)
and $\widetilde{f}''$ is a morphism in the fiber $\calK_{\calO} \times_{ \Fun( \{0\}, \calO^{\otimes})} \{X\}$.
Part $(7)$ of the proof of Proposition \ref{sapless} shows that $\widetilde{f}''$ is inert, so that
$F( \widetilde{f}'')$ is inert in $\calC^{\otimes}$ (since the domain of the morphism $f$ belongs
to $\MMod^{\calO}(\calC)^{\otimes}$). Since the collection of inert morphisms in $\calC^{\otimes}$ is stable under composition, we conclude that $F( \widetilde{f} )$ is inert, as desired. 
\end{proof}

\subsection{Algebra Objects of $\infty$-Categories of Modules}\label{almod}

Let $\calC$ be a symmetric monoidal category, let $A$ be a commutative algebra object of
$\calC$, and let $\calD = \Mod_{A}(\calC)$ be the category of $A$-modules in $\calC$.
Under some mild hypotheses, the category $\calD$ inherits the structure of a symmetric monoidal category. Moreover, one can show the following:

\begin{itemize}
\item[$(1)$] The forgetful functor $\theta: \calD \rightarrow \calC$ induces an equivalence of categories from the category of commutative algebra objects $\CAlg(\calD)$ to the category 
$\CAlg(\calC)_{A/}$ of commutative algebra objects $A' \in \CAlg(\calD)$ equipped with a map
$A \rightarrow A'$.
\item[$(2)$] Given a commutative algebra object $B \in \CAlg(\calD)$, the category of
$B$-modules in $\calD$ is equivalent to the category of $\theta(B)$-modules in $\calC$.
\end{itemize}

Our goal in this section is to obtain $\infty$-categorical analogues of the above statements
for algebras over an arbitrary coherent $\infty$-operad (Corollaries \ref{skoke} and \ref{skoke2}).
Before we can state our results, we need to introduce a bit of terminology.

\begin{notation}\label{quay}
Let $p: \calO^{\otimes} \rightarrow \Nerve(\FinSeg)$ be an $\infty$-operad, and let
$\calK_{\calO} \subseteq \Fun( \Delta^1, \calO^{\otimes})$ be defined as in Notation \ref{salvi}.
We let $\calO^{\otimes}_{\ast}$ denote the $\infty$-category of pointed objects of
$\calO^{\otimes}$: that is, the full subcategory of $\Fun( \Delta^1, \calO^{\otimes})$ spanned
by those morphisms $X \rightarrow Y$ such that $X$ is a final object of
$\calO^{\otimes}$ (which is equivalent to the requirement that $p(X) = \seg{0}$).
If $\calO^{\otimes}$ is unital, then a diagram $\Delta^1 \rightarrow \calO^{\otimes}$ belongs
to $\calO^{\otimes}_{\ast} \subseteq \Fun( \Delta^1, \calO^{\otimes})$ is and only if it is
a left Kan extension of its restriction to $\{1\}$. In this case, Proposition \toposref{lklk} implies that
evaluation at $\{1\}$ induces a trivial Kan fibration $e: \calO^{\otimes}_{\ast} \rightarrow \calO^{\otimes}$.
We let $s: \calO^{\otimes} \rightarrow \calO^{\otimes}_{\ast}$ denote a section of $e$, and
regard $s$ as a functor from $\calO^{\otimes}$ to $\Fun(\Delta^1, \calO^{\otimes}$. We observe that there is a canonical natural transformation $s \rightarrow \delta$, where
$\delta$ is the diagonal embedding $\calO^{\otimes} \rightarrow \Fun( \Delta^1, \calO^{\otimes})$. 
We regard this natural transformation as defining a map
$$ \gamma_{ \calO^{\otimes}}: \calO^{\otimes} \times \Delta^1 \rightarrow \calK_{\calO}.$$
\end{notation}

\begin{remark}
More informally, the map $\gamma_{\calO^{\otimes}}$ can be described as follows.
If $X \in \calO^{\otimes}$, then 
$$\gamma_{\calO^{\otimes}}(X,i) = \begin{cases} \id_{X} \in \calK_{\calO} & \text{if } i = 1 \\
(f: 0 \rightarrow X) \in \calK_{\calO} & \text{if } i = 0. \end{cases}$$
Here $0$ denotes a zero object of $\calO^{\otimes}$.
\end{remark}

Let $\calO^{\otimes}$ be a coherent $\infty$-operad and let
$\calC^{\otimes} \rightarrow \calO^{\otimes}$ be a $\calO$-operad family.
Unwinding the definitions, we see that giving an $\calO$-algebra in
$\MMod^{\calO}(\calC)^{\otimes}$ is equivalent to giving a commutative diagram of simplicial sets
$$ \xymatrix{ \calK_{\calO} \ar[dr]^{e_1} \ar[rr]^{f} & & \calC^{\otimes} \ar[dl] \\
& \calO^{\otimes}. & }$$ 
such that $f$ preserves inert morphisms. Composing with the map $\gamma_{\calO}$ of
Notation \ref{quay}, we obtain a map
$$ \Alg_{\calO}( \MMod^{\calO}(\calC)) \rightarrow \Fun_{ \calO^{\otimes}}( \calO^{\otimes} \times \Delta^1, \calC^{\otimes}).$$

Our main results can be stated as follows:

\begin{proposition}\label{urs}
Let $\calO^{\otimes}$ be a coherent $\infty$-operad, and let
$q: \calC^{\otimes} \rightarrow \calO^{\otimes}$ be a fibration of $\infty$-operads.
Then the construction above determines a categorical equivalence
$$ \Alg_{\calO}( \MMod^{\calO}(\calC)) \rightarrow \Fun( \Delta^1, \Alg_{\calO}(\calC)).$$
\end{proposition}

\begin{proposition}\label{stabs}
Let $\calO^{\otimes}$ be a coherent $\infty$-operad, let
$p: \calC^{\otimes} \rightarrow \calO^{\otimes}$ be a fibration of $\infty$-operads,
and let $A \in \Alg_{\calO}(\calC)$. Then the forgetful functor
$\Mod^{\calO}_{A}(\calC)^{\otimes} \rightarrow \calC^{\otimes}$ induces a homotopy pullback diagram
of $\infty$-categories
$$ \xymatrix{ \Mod^{\calO}( \Mod^{\calO}_{A}(\calC))^{\otimes} \ar[r] \ar[d] & \Mod^{\calO}(\calC)^{\otimes} \ar[d] \\
\Alg_{\calO}( \Mod^{\calO}_{A}(\calC) ) \times \calO^{\otimes} \ar[r] & \Alg_{\calO}(\calC) \times \calO^{\otimes}. }$$
\end{proposition}

We defer the proofs of Proposition \ref{urs} and \ref{stabs} until the end of this section.

\begin{corollary}\label{usr}
Let $\calO^{\otimes}$ be a coherent $\infty$-operad and let
$\calC^{\otimes} \rightarrow \calO^{\otimes}$ be a fibration of $\infty$-operads.
The composition
$$ \theta: \Alg_{\calO}( \Mod^{\calO}(\calC)) \rightarrow
\Alg_{\calO}(\MMod^{\calO}(\calC)) \rightarrow \Fun( \Delta^1, \Alg_{\calO}(\calC) )$$
is an equivalence of $\infty$-categories.
\end{corollary}

\begin{proof}
Combine Proposition \ref{urs} with Remark \ref{capp}.
\end{proof}

\begin{remark}\label{infl}
Let $\theta: \Alg_{\calO}( \Mod^{\calO}(\calC) ) \rightarrow \Fun( \Delta^1, \Alg_{\calO}(\calC))$ be
the categorical equivalence of Corollary \ref{usr}. Composing with the map
$\Fun( \Delta^1, \Alg_{\calO}(\calC)) \rightarrow \Alg_{\calO}(\calC)$ given by evaluation
at $\{0\}$, we obtain a map $\theta_0: \Alg_{\calO}( \Mod^{\calO}(\calC)) \rightarrow \Alg_{\calO}(\calC)$.
Unwinding the definitions, we see that $\theta_0$ factors as a composition
$$ \Alg_{\calO}( \Mod^{\calO}(\calC)) \stackrel{\theta'_0}{\rightarrow}
\Alg_{\calO}( \Alg_{\calO}(\calC) \times \calO) \stackrel{\theta''_0}{\rightarrow} \Alg_{\calO}(\calC),$$
where $\theta''_0$ is induced by composition with the diagonal embedding
$\calO^{\otimes} \rightarrow \calO^{\otimes} \times \calO^{\otimes}$. In particular, if
$A$ is a $\calO$-algebra object of $\calC$, then the
the restriction of $\theta_0$ to $\Alg_{\calO}( \Mod^{\calO}_{A}(\calC) )$
is a constant map taking the value $A \in \Alg_{\calO}(\calC)$.
\end{remark}

\begin{corollary}\label{skoke}
Let $\calO^{\otimes}$ be a coherent $\infty$-operad, let $\calC^{\otimes} \rightarrow \calO^{\otimes}$ be a fibration of $\infty$-operads, and let $A \in \Alg_{\calO}(\calC)$ be a $\calO$-algebra object
of $\calC$. Then the categorical equivalence $\theta$ of Corollary \ref{usr} restricts to a categorical equivalence
$$ \theta_{A}: \Alg_{\calO}( \Mod^{\calO}_{A}(\calC) ) \rightarrow \Alg_{\calO}(\calC)^{A/}.$$
\end{corollary}

\begin{proof}
Remark \ref{infl} guarantees that the restriction of $\theta$ carries 
$\Alg_{\calO}( \Mod^{\calO}_{A}(\calC))$ into $$\Alg_{\calO}(\calC)^{A/} \subseteq
\Fun( \Delta^1, \Alg_{\calO}(\calC)^{A/}).$$ Consider the diagram
$$ \xymatrix{ \Alg_{\calO}( \Mod^{\calO}_{A}(\calC)) \ar[r] \ar[d] & \Alg_{\calO}( \Mod^{\calO}(\calC)) \ar[r]^{\theta} \ar[d] & \Fun( \Delta^1, \Alg_{\calO}( \calC) ) \ar[d] \\
\Delta^0 \ar[r] & \Alg_{\calO}( \Alg_{\calO}(\calC) \times \calO) \ar[r]^{\theta'} & \Fun( \{0\}, \Alg_{\calO}(\calC) ). }$$
The left square is a homotopy pullback, since it is a pullback square between fibrant objects in which the vertical maps are categorical fibrations. The right square is a homotopy pullback since both of the horizontal arrows are categorical equivalences ($\theta$ is a categorical equivalence by virtue of Corollary \ref{usr}, and $\theta'$ is a categorical equivalence since it is left inverse to the
categorical equivalence described in Example \ref{insert}). It follows that the outer square is
also a homotopy pullback, which is equivalent to the assertion that $\theta_{A}$ is a categorical equivalence (Proposition \toposref{basechangefunky}). 
\end{proof}

\begin{corollary}\label{costabs}
Let $\calO^{\otimes}$ be a coherent $\infty$-operad, $p: \calC^{\otimes} \rightarrow \calO^{\otimes}$ a fibration of $\infty$-operads, and $A \in \Alg_{\calO}(\calC)$ an algebra object. Then
there is a canonical equivalence of $\infty$-categories
$$ \Mod^{\calO}( \Mod^{\calO}_{A}( \calC))^{\otimes}
\rightarrow \Mod^{\calO}(\calC)^{\otimes} \times_{ \Alg_{\calO}(\calC)} \Alg_{\calO}(\calC)^{A/}.$$
\end{corollary}

\begin{proof}
Combine Proposition \ref{stabs} with Corollary \ref{skoke}.
\end{proof}

\begin{corollary}\label{skoke2}
Let $\calO^{\otimes}$ be a coherent $\infty$-operad, $p: \calC^{\otimes} \rightarrow \calO^{\otimes}$ a fibration of $\infty$-operads. Let $A \in \Alg_{\calO}(\calC)$, let
$\overline{B} \in \Alg_{\calO}( \Mod^{\calO}_{A}(\calC))$, and let $B \in \Alg_{\calO}(\calC)$
be the algebra object determined by $\overline{B}$. Then there is a canonical
equivalence of $\infty$-operads
$$ \Mod^{\calO}_{\overline{B}}( \Mod^{\calO}_{A}(\calC))^{\otimes} \rightarrow \Mod^{\calO}_{B}(\calC)^{\otimes}.$$
\end{corollary}

We now turn to the proof of Propositions \ref{urs} and \ref{stabs}. We will need a few preliminary results.

\begin{lemma}\label{ilker}
Suppose we are given a commutative diagram of $\infty$-categories
$$ \xymatrix{ \calC \ar[rr]^{F} \ar[dr]^{p} & & \calD \ar[dl]^{q} \\
& \calE & }$$
where $p$ and $q$ are Cartesian fibrations and the map $F$ carries
$p$-Cartesian morphisms to $q$-Cartesian morphisms. Let
$D \in \calD$ be an object, let $E = q(D)$, and let $\calC_{D/} = \calC \times_{ \calD} \calD_{D/}$. Then:
\begin{itemize}
\item[$(1)$] The induced map $p': \calC_{D/} \rightarrow \calE_{E/}$
is a Cartesian fibration.
\item[$(2)$] A morphism $f$ in $\calC_{D/}$ is $p'$-Cartesian
if and only if its image in $\calC$ is $p$-Cartesian.
\end{itemize}
\end{lemma}

\begin{proof}
Let us say that a morphism in $\calC_{D/}$ is {\it special} if
its image in $\calC$ is $p$-Cartesian. We first prove the ``if'' direction of $(2)$ by showing
that every special morphism of $\calC_{D/}$ is $p$-Cartesian; we will simultaneously show that $p'$ is an inner fibration. For this, we must show that every lifting problem of the form
$$ \xymatrix{ \Lambda^n_i \ar[d] \ar[r]^{g_0} & \calC_{D/} \\
\Delta^n \ar[r] \ar@{-->}[ur]^{g} & \calE_{E/ } }$$
admits a solution, provided that $n \geq 2$ and either $0 < i < n$ or $i = n$ and
$g_0$ carries $\Delta^{ \{n-1, n\} }$ to a special morphism $e$ in $\calC_{D/}$.
To prove this, we first use the fact that $p$ is an inner fibration (together with the
observation that the image $e_0$ of $e$ in $\calC$ is $p$-Cartesian when $i=n$) to solve the associated lifting problem
$$ \xymatrix{ \Lambda^n_i \ar[d] \ar[r]^{g'_0} & \calC \ar[d] \\
\Delta^n \ar[r] \ar@{-->}[ur]^{g'} & \calE. }$$
To extend this to solution of our original lifting problem, we are required to solve
another lifting problem of the form
$$ \xymatrix{ \Lambda^{n+1}_{i+1} \ar[r]^{g''_0} \ar[d] & \calD \ar[d]^{q} \\
\Delta^{n+1} \ar@{-->}[ur]^{g''} \ar[r] & \calE. }$$
If $i < n$, the desired solution exists by virtue of our assumption that $q$ is
an inner fibration. If $i = n$, then it suffices to observe that
$g''_0( \Delta^{ \{n, n+1\} }) = F(e_0)$ is a $q$-Cartesian morphism in $\calD$.

To prove $(1)$, it will suffice to show that for every object $\overline{C} \in \calC_{D/}$ and
every morphism $\overline{f}_0: \overline{E}' \rightarrow p'( \overline{C})$ in $\calE_{E/}$, there
exists a special morphism $\overline{f}$ in $\calC_{D/}$ with $p'(\overline{f}) = \overline{f}_0$. We can identify
$\overline{C}$ with an object $C \in \calC$ together with a morphism $\alpha: D \rightarrow F(C)$
in $\calD$, and we can identify $\overline{f}_0$ with a $2$-simplex
$$ \xymatrix{ & E' \ar[dr]^{f_0} & \\
E \ar[ur] \ar[rr]^{q(\alpha)} & & p(C) }$$
in $\calE$. Since $p$ is a Cartesian fibration, we can choose a $p$-Cartesian morphism
$f: C' \rightarrow C$ with $p(f) = f_0$. In order to lift $f$ to a special morphism
$\overline{f}: \overline{C}' \rightarrow \overline{C}$, it suffices to complete the diagram
$$ \xymatrix{ & F(C') \ar[dr]^{ F(f_0) } \ar[dr] & \\
D \ar@{-->}[ur] \ar[rr]^{\alpha} & & F(C) }$$
to a $2$-simplex of $\calD$. This is possibly by virtue of our assumption that $F(f_0)$
is $q$-Cartesian.

To complete the proof of $(2)$, it will suffice to show that every $p'$-Cartesian morphism
$\overline{f}': \overline{C}'' \rightarrow \overline{C}$ of $\calC_{D/}$ is special. The proof
of $(1)$ shows that there exists a special morphism $\overline{f}: \overline{C}' \rightarrow \overline{C}$
with $p'( \overline{f}) = p'( \overline{f}')$. Since $\overline{f}$ is also $p'$-Cartesian, it is equivalent to
$\overline{f}'$, so that $\overline{f}'$ is also special.
\end{proof}

\begin{lemma}\label{cambus}
Suppose given a commutative diagram of $\infty$-categories
$$ \xymatrix{ \calC \ar[rr]^{F} \ar[dr]^{p} & & \calD \ar[dl]^{q} \\
& \calE & }$$
where $p$ and $q$ are Cartesian fibrations, and the map $F$ carries
$p$-Cartesian morphisms to $q$-Cartesian morphisms. Suppose
furthermore that for every object $E \in \calE$, the induced functor
$F_{E}: \calC_{E} \rightarrow \calD_{E}$ is cofinal. Then
$F$ is cofinal.
\end{lemma}

\begin{proof}
In view of Theorem \toposref{hollowtt}, it suffices to show that for each object
$D \in \calD$, the simplicial set $\calC_{D/} = \calC \times_{ \calD} \calD_{D/}$ is weakly contractible.
Let $E$ denote the image of $D$ in $\calE$; we observe that $\calC_{D/}$ comes equipped with a map $p': \calC_{D/} \rightarrow \calE_{E/}$. Moreover, the fiber of $p'$ over the initial object
$\id_{E} \in \calE_{E/}$ can be identified with $\calC_{E} \times_{ \calD_{E} } (\calD_{E})_{D/}$, which
is weakly contractible by virtue of our assumption that $F_{E}$ is cofinal (Theorem \toposref{hollowtt}).
To prove that $\calC_{D/}$ is contractible, it will suffice to show that the inclusion
$i: {p'}^{-1} \{ \id_{E} \} \hookrightarrow \calC_{D/}$ is a weak homotopy equivalence.
We will prove something slightly stronger: the inclusion $i^{op}$ is cofinal.
Since the inclusion $\{ \id_{E} \}^{op} \rightarrow \calE_{E/}^{op}$ is evidently cofinal, it will
suffice to show that $p'$ is a Cartesian fibration (Proposition \toposref{strokhop}). This follows
from Lemma \ref{ilker}.
\end{proof}

\begin{lemma}\label{qlum}
Let $q: \calC^{\otimes} \rightarrow \calO^{\otimes}$ be a fibration of $\infty$-operads,
and let $\alpha: \seg{n} \rightarrow \seg{m}$ be an inert morphism in $\FinSeg$.
Let $\sigma$ denote the diagram
$$ \xymatrix{ \seg{n} \ar[r]^{\alpha} \ar[d]^{\id} & \seg{m} \ar[d]^{\id} \\
\seg{n} \ar[r]^{\alpha} & \seg{m} }$$
in $\Nerve(\FinSeg)$, let $K \simeq \Lambda^2_2$ be the full subcategory of
$\Delta^1 \times \Delta^1$ obtained by omitting the initial vertex, and let
$\sigma_0 = \sigma|K$. Suppose that $\overline{ \sigma_0}: K \rightarrow  \calC^{\otimes}$
is a diagram lifting $\sigma_0$, corresponding to a commutative diagram
$$ \xymatrix{ & X \ar[d] \\
Y' \ar[r]^{\overline{\alpha}'} & X' }$$
where $\overline{\alpha}'$ is inert.
Then:
\begin{itemize}
\item[$(1)$] Let $\overline{\sigma}: \Delta^1 \times \Delta^1 \rightarrow \calC^{\otimes}$
be an extension of $\overline{\sigma}_0$ lifting $\sigma$, corresponding to a commutative
diagram
$$ \xymatrix{ Y \ar[r]^{\overline{\alpha}} \ar[d]^{\overline{\beta}} & X \ar[d] \\
Y' \ar[r]^{\overline{\alpha}'} & X' }$$
in $\calC^{\otimes}$. Then $\sigma$ is a $q$-limit diagram if and only if
it satisfies the following conditions:
\begin{itemize}
\item[$(i)$] The map $\overline{\alpha}$ is inert. 
\item[$(ii)$] Let $\gamma: \seg{n} \rightarrow \seg{k}$ be an inert morphism in
$\Nerve(\FinSeg)$ such that $\alpha^{-1} \nostar{m} \subseteq \gamma^{-1} \{ \ast \}$, and
let $\overline{\gamma}: Y' \rightarrow Z$ be an inert morphism in $\calC^{\otimes}$ lifting
$\gamma$. Then $\overline{\gamma} \circ \overline{\beta}: Y \rightarrow Z$ is an inert morphism in $\calC^{\otimes}$.
\end{itemize}

\item[$(2)$] There exists an extension $\overline{\sigma}$ of $\overline{\sigma}_0$
lying over $\sigma$ which satisfies conditions $(i)$ and $(ii)$ of $(1)$.
\end{itemize}
\end{lemma}

\begin{proof}
This is a special case of Lemma \ref{helpless}.
\end{proof}

\begin{proof}[Proof of Proposition \ref{urs}]
We define a simplicial set $\calM$ equipped with a map
$p: \calM \rightarrow \Delta^1$ so that the following universal property is satisfied:
for every map of simplicial sets $K \rightarrow \Delta^1$, the set
$\Hom_{ (\sSet)_{/\Delta^1}}( K, \calM)$ can be identified with the collection of all
commutative diagrams
$$ \xymatrix{ K \times_{ \Delta^1} \{1\} \rightarrow \ar[r] \ar[d] & \calO^{\otimes} \times \Delta^1 \ar[d]^{\gamma_{\calO}} \\
K \ar[r] & \calK_{\calO}. }$$
The map $p$ is a Cartesian fibration, associated to the functor
$\gamma_{\calO}$ from $\calM_{1} \simeq \calO^{\otimes} \times \Delta^1$ to
$\calM_0 \simeq \calK_{\calO}$. We observe that $\calM$ is equipped with a functor
$\calM \rightarrow \calO^{\otimes}$, whose restriction to $\calM_{1} \simeq \calO^{\otimes} \times \Delta^1$ is given by projection onto the first factor and whose restriction
to $\calM_0 \simeq \calK_{\calO}$ is given by evaluation at $\{1\}$. 

We let $\calX$ denote the full subcategory of $\Fun_{ \calO^{\otimes}}( \calM, \calC^{\otimes})$ spanned
by those functors $F$ satisfying the following pair of conditions:
\begin{itemize}
\item[$(i)$] The functor $F$ is a $q$-left Kan extension of $F | \calM_0$.
\item[$(ii)$] The restriction $F | \calM_0 \in \Fun_{ \calO^{\otimes}}( \calK_{\calO}, \calC^{\otimes})$
belongs to $\Alg_{ \calO}( \MMod^{\calO}(\calC))$. 
\end{itemize}
Since
$p$ is a Cartesian fibration, condition $(i)$ can be reformulated as follows:
\begin{itemize}
\item[$(i')$] For every $p$-Cartesian morphism $f$ in $\calM$, the image
$F(f)$ is a $q$-coCartesian morphism in $\calC^{\otimes}$. Since the image of
$f$ in $\calO^{\otimes}$ is an equivalence, this is equivalent to the requirement
that $F(f)$ is an equivalence in $\calC^{\otimes}$.
\end{itemize}
Using Proposition \toposref{lklk}, we deduce that the restriction map
$\calX \rightarrow \Alg_{\calO}( \MMod^{\calO}(\calC))$ is a trivial Kan fibration.
This restriction map has a section $s$, given by composition with the natural retraction
$r: \calM \rightarrow \calM_0$. It follows that $s$ is a categorical equivalence, and that
every object $F \in \calX$ is equivalent $(F | \calM_0) \circ r$. We deduce that
restriction to $\calM_1 \subseteq \calM$ induces a functor
$\theta': \calX \rightarrow \Fun( \Delta^1, \Alg_{\calO}(\calC) )$. 
We have a commutative diagram $$ \xymatrix{ & \calX \ar[dr]^{\theta'} & \\
\Alg_{ \calO}( \MMod^{\calO}(\calC)) \ar[ur]^{s} \ar[rr]^{\theta} & & \Fun( \Delta^1, \Alg_{\calO}(\calC)). }$$
To complete the proof, it will suffice to show that $\theta'$ is a categorical equivalence. We will
show that $\theta'$ is a trivial Kan fibration. In view of Proposition \toposref{lklk}, it will suffice to prove the following:

\begin{itemize}
\item[$(a)$] A functor $F \in \Fun_{ \calO^{\otimes}}( \calM, \calC^{\otimes})$ belongs to
$\calX$ if and only if $F_1 = F| \calM_1 \in \Fun( \Delta^1, \Alg_{\calO}(\calC))$ and
$F$ is a $q$-right Kan extension of $F_1$.

\item[$(b)$] Every object $F_1 \in \Fun_{\calO^{\otimes}}(\calM_1, \calC^{\otimes})$
belonging to $\Fun( \Delta^1, \Alg_{\calO}(\calC))$ admits a $q$-right Kan extension
of $F \in \Fun_{ \calO^{\otimes}}( \calM, \calC^{\otimes})$.
\end{itemize}

To prove these claims, we will need a criterion for detecting whether a functor
$F \in \Fun_{ \calO^{\otimes}}( \calM, \calC^{\otimes})$ is a $q$-right Kan extension
of $F_1 = F| \calM_1 \in \Fun( \Delta^1, \Alg_{\calO}(\calC))$ at an object
$\overline{X} \in \calM_0$. Let $\overline{X}$ correspond to a semi-inert
morphism $\overline{\alpha}: X' \rightarrow X$ in $\calO^{\otimes}$, covering a morphism
$\alpha: \seg{m} \rightarrow \seg{n}$ in $\Nerve(\FinSeg)$. Let
$\calD$ denote the $\infty$-category $( \calO^{\otimes} \times \Delta^1) \times_{ \calM} \calM_{\overline{X}/}$, so that $\calD$ is equipped with a projection $\calD \rightarrow \Delta^1$; we let $\calD_0$ and $\calD_1$ denote the fibers of this map. Form a pushout diagram
$$ \xymatrix{ \seg{m} \ar[r]^{\alpha} \ar[d] & \seg{n} \ar[d]^{\beta} \\
\seg{0} \ar[r] & \seg{k} }$$
in $\FinSeg$, and choose an inert morphism $\overline{\beta}: X' \rightarrow Y$ in $\calO^{\otimes}$
lying over $\beta$. Let $\overline{Y}$ denote the image of $Y$ under our chosen section
$s: \calO^{\otimes} \rightarrow \calO^{\otimes}_{\ast}$, so we can identify $\overline{Y}$ with
a morphism $0 \rightarrow Y$ in $\calO^{\otimes}$, where $0$ is a zero object of $\calO^{\otimes}$.
We can therefore lift $\overline{\beta}$ to a morphism $\widetilde{\beta}:\overline{X} \rightarrow \overline{Y}$ in $\calK_{\calO}$, corresponding to a commutative diagram
$$ \xymatrix{ X' \ar[r]^{\overline{\alpha}} \ar[d] & X \ar[d]^{\overline{\beta}} \\
0 \ar[r] & Y. }$$
The pair $(Y, \widetilde{\beta})$ can be identified with an object of $\calN_{0}$, which we
will denote by $\widetilde{Y}$. We claim that $\widetilde{Y}$ is an initial object of
$\calN_0$. Unwinding the definitions, this is equivalent to the following assertion:
for every object $A \in \calO^{\otimes}$, composition with $\widetilde{\beta}$ induces a homotopy
equivalence $$ \phi: \bHom_{ \calO^{\otimes}}(Y, A) \rightarrow \bHom_{ \calK_{\calO}}(\overline{X}, s(A)).$$
To prove this, we observe that $\phi$ factors as a composition
$$ \bHom_{ \calO^{\otimes}}( Y,A) \stackrel{\phi'}{\rightarrow}
\bHom_{\calK_{\calO}}(\overline{Y}, s(A)) \stackrel{\phi''}{\rightarrow}
\bHom_{\calK_{\calO}}( \overline{X}, s(A) ).$$
The map $\phi'$ is a homotopy equivalence because $s$ is a categorical equivalence,
and the map $\phi''$ is a homotopy equivalence because $\overline{\beta}$ is
coCartesian with respect to the projection $\calK_{\calO} \rightarrow \calK_{\Nerve(\FinSeg)}$.

We have an evident natural transformation $\widetilde{\gamma}: \overline{X} \rightarrow \id_{X}$
in $\calK_{\calO}$. The pair $(X, \widetilde{\gamma})$ determines an object
$\widetilde{Z} \in \calN_1$. We claim that $\widetilde{Z}$ is an initial object of $\calN_1$.
Unwinding the definitions, we see that this is equivalent to the assertion that for every
object $A \in \calO^{\otimes}$, composition with $\widetilde{\gamma}$ induces a homotopy equivalence
$$ \psi: \bHom_{ \calO^{\otimes}}( X, A) \rightarrow \bHom_{ \calK_{\calO}}( \overline{X}, \delta(A) ),$$
where $\delta: \calO^{\otimes} \rightarrow \calK_{\calO}$ is the diagonal embedding.
To prove this, we factor $\psi$ as a composition
$$ \bHom_{ \calO^{\otimes}}( X,A) \stackrel{\psi'}{\rightarrow} \bHom_{ \calK_{\calO}}( \delta(X), \delta(A))
\stackrel{\psi''}{\rightarrow} \bHom_{ \calK_{\calO}}( \overline{X}, \delta(A)).$$
The map $\psi'$ is a homotopy equivalence since $\delta$ is fully faithful, and the map
$\psi''$ is a homotopy equivalence by virtue of Corollary \toposref{funcsys} (applied to the
trivial factorization system on $\calO^{\otimes}$).

Since $\calN \rightarrow \calO^{\otimes} \times \Delta^1$ is a left fibration,
we can lift the map $(Y,0) \rightarrow (Y,1)$ to a map $e: \widetilde{Y} \rightarrow \widetilde{Y}'$
in $\calN$. Since $\widetilde{Z}$ is an initial object of $\calN_{1}$, we can choose a map
$e': \widetilde{Z} \rightarrow \widetilde{Y}$ in $\calN_1$. Let $C \simeq \Lambda^2_{2}$
denote the full subcategory of $\Delta^1 \times \Delta^1$ obtained by omitting the final vertex,
so that $e$ and $e'$ together determine a map of simplicial sets $C \rightarrow \calN$. Applying
the dual of Lemma \ref{cambus} to the diagram
$$ \xymatrix{ C \ar[rr] \ar[dr] & & \calN \ar[dl] \\
& \Delta^1, & }$$
we deduce that $C^{op} \rightarrow \calN^{op}$ is cofinal. We therefore arrive at the following:

\begin{itemize}
\item[$(\ast)$] A functor $F \in \Fun_{ \calO^{\otimes}}( \calM, \calC^{\otimes})$ is a $q$-right Kan extension of $F_1 = F| \calM_1 \in \Fun( \Delta^1, \Alg_{\calO}(\calC))$ at an object
$\overline{X} \in \calK_{\calO}$ if and only if the induced diagram
$$ \xymatrix{ F( \overline{X}) \ar[r] \ar[d] & F_1(X,1) \ar[d] \\
F_1(Y,0) \ar[r] & F_1(Y,1) }$$
is a $q$-limit diagram.
\end{itemize}

Moreover, Lemma \toposref{kan2} yields the following:
\begin{itemize}
\item[$(\ast')$] A functor $F_1 \in \Fun( \Delta^1, \Alg_{\calO}(\calC))$ admits
a $q$-right Kan extension $F \in \Fun_{ \calO^{\otimes}}( \calM, \calC^{\otimes})$ if and only if,
for every object $\overline{X} \in \calK_{\calO}$, the diagram
$$ \xymatrix{ & F_1(X,1) \ar[d] \\
F_1(Y,0) \ar[r] & F_1(Y,1) }$$
can be extended to a $q$-limit diagram lying over the diagram
$$ \xymatrix{ X \ar[r]^{\id} \ar[d] & X \ar[d] \\
Y \ar[r]^{\id} & Y. }$$
\end{itemize}

Assertion $(b)$ follows immediately from $(\ast')$ together with Lemma \ref{qlum}.
Combining assertion $(\ast)$ with Lemma \ref{qlum}, we deduce that a functor
$F \in \Fun_{ \calO^{\otimes}}( \calM, \calC^{\otimes})$ is a 
$q$-right Kan extension of $F_1 =  F| \calM_1 \in \Fun( \Delta^1, \Alg_{\calO}(\calC))$
if and only if the following conditions are satisfied:

\begin{itemize}
\item[$(i')$] The restriction $F_1$ belongs to $\Fun( \Delta^1, \Alg_{\calO}(\calC))$.
That is, $F$ carries every inert morphism in $\calO^{\otimes} \times \{j\} \subseteq \calM_1$ to an inert
morphism in $\calC^{\otimes}$, for $j \in \{0,1\}$.

\item[$(ii')$] Let $\overline{X}$ be as above.
Then the induced morphism $F( \overline{X} ) \rightarrow F( Y, 0)$ is inert.

\item[$(iii')$] Let $\overline{X}$ be as above, and suppose that we are given
an inert morphism $\alpha': \seg{n} \rightarrow \seg{l}$ such that
the composite map $\alpha' \circ \alpha: \seg{m} \rightarrow \seg{n'}$ is surjective
together with an inert morphism $\overline{\alpha}': X \rightarrow X''$ lifting $\alpha'$.
Then the composite map $F( \overline{X} ) \rightarrow F(X,1) \rightarrow F(X'',1)$ is an inert morphism in $\calC^{\otimes}$.
\end{itemize}
To complete the proof, it will suffice to show that a functor
$F \in \Fun_{\calO^{\otimes}}( \calM, \calC^{\otimes})$ satisfies conditions $(i)$ and
$(ii)$ if and only if it satisfies conditions $(i')$, $(ii')$, and $(iii')$. 

Suppose first that $F$ satisfies $(i)$ and $(ii)$. We have already seen that
$F$ must also satisfy $(i')$. To prove $(ii')$, we observe that the map $F( \overline{X}) \rightarrow F(Y,0)$ factors as a composition
$$ F( \overline{X}) \rightarrow F( \overline{Y} ) \rightarrow F(Y,0).$$
The first map is inert because $F$ satisfies $(ii)$ and $\widetilde{\beta}: \overline{X} \rightarrow \overline{Y}$ is an inert morphism in $\calK_{\calO}$. The second morphism is inert by virtue of assumption $(i)$. Now suppose that we are given an inert morphism $X \rightarrow X''$
as in $(iii')$. We have a commutative diagram
$$ \xymatrix{ X \ar[r] \ar[d] & X' \ar[d] \\
X'' \ar[r]^{\id} & X'' }$$
corresponding to an inert morphism $\overline{X} \rightarrow \delta( X'')$ in $\calK_{\calO}$, and
the map $F( \overline{X}) \rightarrow F(X'',1)$ factors as a composition
$$ F( \overline{X} ) \rightarrow F( \delta(X'') ) \rightarrow F( X'', 1).$$
The first of these maps is inert by virtue of assumption $(ii)$, and the second by virtue of assumption
$(i)$.

Now suppose that $F$ satisfies $(i')$, $(ii')$, and $(iii')$; we wish to show that $F$ satisfies
$(i)$ and $(ii)$. To prove $(i)$, we must show that for every object $X \in \calO^{\otimes}$, the
morphisms $F( s(X) ) \rightarrow F(X,0)$ and $F( \delta(X) ) \rightarrow F(X,1)$ are
inert in $\calC$. The first of these assertions is a special case of $(ii')$, and the second
is a special case of $(iii')$. To prove $(ii)$, consider an arbitrary inert morphism
$\overline{\beta}: \overline{X} \rightarrow \overline{Y}$ in $\calO^{\otimes}$, corresponding
to a commutative diagram $\sigma$:
$$ \xymatrix{ X' \ar[d]^{\beta'} \ar[r] & X \ar[d]^{\beta} \\
Y' \ar[r] & Y }$$
in the $\infty$-category $\calO^{\otimes}$. We wish to show that $F( \overline{\beta})$ is an inert morphism
in $\calC^{\otimes}$. Let $\beta_0: \seg{n} \rightarrow \seg{k}$ be the image of
$\beta$ in $\Nerve(\FinSeg)$, and let $(\beta_0)_{!}: \calC^{\otimes}_{\seg{n}} \rightarrow \calC^{\otimes}_{\seg{k}}$ denote the induced functor. Then $F( \overline{\beta})$ factors as a composition $F( \overline{X}) \stackrel{\epsilon}{\rightarrow} (\beta_0)_{!} F(\overline{X}) \stackrel{\epsilon'}{\rightarrow} F(\overline{Y})$, where $\epsilon$ is inert; we wish to prove that $\epsilon'$ is an equivalence in $\calC^{\otimes}_{\seg{k}}$. Since $\calC^{\otimes}$ is an $\infty$-operad, it will suffice to show that $\Colll{j}{!} \epsilon'$ is an equivalence in
$\calC^{\otimes}_{\seg{1}}$ for $1 \leq j \leq k$. Let
$$ \xymatrix{ \seg{n} \ar[r] \ar[d] & \seg{m} \ar[d]^{\beta_0} \\
\seg{k'} \ar[r] & \seg{k} \\ }$$
denote the image of $\sigma$ in $\Nerve(\FinSeg)$. This diagram admits a unique
extension
$$ \xymatrix{ \seg{n} \ar[r] \ar[d] & \seg{m} \ar[d]^{\beta_0} \\
\seg{k'} \ar[d] \ar[r]^{\chi} & \seg{k} \ar[d]^{\Colp{j}} \\
\seg{t} \ar[r]^{\chi'} & \seg{1}} $$
where the vertical morphisms are inert, the integer $t$ is equal to
$1$ and $\chi'$ is an isomorphism if $j$ lies in the image of
$\seg{n} \rightarrow \seg{k}$, and $t=0$ otherwise.
We can lift this diagram to a commutative triangle
$$ \xymatrix{ \overline{X} \ar[rr]^{\overline{\beta}} \ar[dr]^{\overline{\beta}'} & & \overline{Y} \ar[dl]^{\overline{\beta}''} \\
& \overline{Z} & }$$
of inert morphisms in $\calK_{\calO}$. If $F( \overline{\beta}'')$ is inert, then
we can identify $F( \overline{\beta}')$ with the composition
$$ F( \overline{X} ) \rightarrow (\Colp{j} \circ \beta_0)_{!}
F( \overline{X}) \stackrel{ \Colll{j}{!}  \epsilon'}{\rightarrow} F( \overline{Z} ),$$
so that $\Colll{j}{!}  \epsilon'$ is an equivalence in $\calC^{\otimes}_{\seg{1}}$ if
and only if $F( \overline{\beta}')$ is inert. We are therefore reduced to proving
that $F( \overline{\beta}')$ and $F( \overline{\beta}'')$ are inert. Replacing
$\overline{\beta}$ by $\overline{\beta}'$ or $\overline{\beta}''$, we may reduce to the where
either $\chi$ is an isomorphism or $k' = 0$.

If $k' = 0$, then we can identify $\overline{Y}$ with $s(Y)$ and condition $(i)$ guarantees that $F( \overline{Y} ) \rightarrow F(Y,0)$ is an equivalence. It therefore suffices to show that the composite map $F( \overline{X} ) \rightarrow F(Y,0)$ is inert. Since the collection of inert morphisms
in $\calC^{\otimes}$ is stable under composition, this follows from $(i')$ and $(ii')$.
If $\chi$ is an isomorphism, then we can identify $\overline{Y}$ with $s(Y)$, and condition
$(i)$ guarantees that $F( \overline{Y}) \rightarrow F(Y,1)$ is an equivalence. 
It therefore suffices to show that the composite map $F( \overline{X}) \rightarrow F(Y,1)$ is inert. This
follows from $(i')$ and $(iii')$, again because the collection of inert morphisms in $\calC^{\otimes}$ is stable under composition. This completes the verification of condition $(ii)$ and the proof of Proposition \ref{urs}
\end{proof}

\begin{proof}[Proof of Proposition \ref{stabs}]
Let $\calK \subseteq \Fun( \Lambda^2_1, \calO^{\otimes})$ be the full subcategory spanned by
those diagrams $X \stackrel{\alpha}{\rightarrow} Y \stackrel{\beta}{\rightarrow} Z$ in $\calO^{\otimes}$ where $\alpha$ and $\beta$ are semi-inert, and let $e_i: \calK \rightarrow \calO^{\otimes}$ be the map given by evaluation at the vertex $\{i\}$ for $0 \leq i \leq 2$. We will say that a morphism in
$\calK$ is {\it inert} if its image under each $e_i$ is an inert morphism in $\calO^{\otimes}$.
If $S$ is a full subcategory of $\calK$, we $\overline{\calX}(S)$ denote the simplicial set 
$(e_0 | S)_{\ast} (e_{2} | S)^{\ast} \calC^{\otimes}$: that is, $\overline{\calX}(S)$ is a simplicial set
equipped with a map $\overline{\calX}(S) \rightarrow \calO^{\otimes}$ characterized by the following universal property: for any map of simplicial sets $K \rightarrow \calO^{\otimes} \simeq \Fun( \{0\}, \calO^{\otimes})$, we have a canonical bijection
$$ \Hom_{ \Fun( \{0\}, \calO^{\otimes})}( K, \overline{\calX}(S)) = \Hom_{ \Fun( \{2\}, \calO^{\otimes})}( K \times_{ \Fun( \{0\}, \calO^{\otimes})} S, \calC^{\otimes}).$$
We let $\calX(S)$ denote the full simplicial subset of $\overline{\calX}(S)$ spanned by those 
vertices which classify functors carrying inert morphisms in $S$ to inert morphisms in $\calC^{\otimes}$.

Let $\calK_1$ denote the full subcategory of $\calK$ spanned by those diagrams
$X \stackrel{\alpha}{\rightarrow} Y \stackrel{\beta}{\rightarrow} Z$ where $\beta$ is an equivalence,
and let $\calK_{01}$ denote the full subcategory of $\calK_1$ spanned by those diagrams where
$\alpha$ is null. We have a canonical embedding $j: \calK_{\calO} \rightarrow \calK_{1}$
which carries $\alpha: X \rightarrow Y$ to the diagram $X \stackrel{\alpha}{\rightarrow} Y \stackrel{\id_{Y}}{\rightarrow} Y$. Note that this embedding restricts to an embedding $j_0: \calK_{\calO}^{0} \hookrightarrow \calK_{01}$. Composition with these embeddings gives rise to a commutative diagram
$$ \xymatrix{ \Mod^{\calO}(\calC)^{\otimes} \ar[d] \ar[r] & \MMod^{\calO}(\calC)^{\otimes} \ar[r] \ar[d] & \calX(\calK_1) \ar[d] \\
\Alg_{\calO}(\calC) \times \calO^{\otimes} \ar[r] & \PrAlg_{\calO}(\calC) \ar[r] & \calX( \calK_{01}). }$$
The left horizontal maps are categorical equivalences (Remark \ref{user2}).
The right horizontal maps are categorical equivalences because $j$ and $j_0$ admit simplicial
homotopy inverses (given by restriction along the inclusion $\Delta^{ \{0,1\}} \subseteq \Lambda^2_1$). 
Consequently, we are reduced to proving that the diagram
$$ \xymatrix{ \Mod^{\calO}( \Mod^{\calO}_{A}(\calC))^{\otimes} \ar[r] \ar[d] & \calX(\calK_1) 
\ar[d] \\
\Alg_{\calO}( \Mod^{\calO}_{A}(\calC) ) \times \calO^{\otimes} \ar[r] & \calX(\calK_{01}) }$$
is a homotopy pullback square.

Let $\calK_{0}$ denote the full subcategory of $\calK$ spanned by those diagrams 
$X \stackrel{\alpha}{\rightarrow} Y \stackrel{\beta}{\rightarrow} Z$ for which $\alpha$ is null, and consider the diagram
$$ \xymatrix{ \Mod^{\calO}( \Mod^{\calO}_{A}(\calC))^{\otimes} \ar[r] \ar[d] & \calX(\calK) \ar[d] \ar[r] &  \calX(\calK_1) 
\ar[d] \\
\Alg_{\calO}( \Mod^{\calO}_{A}(\calC) ) \times \calO^{\otimes} \ar[r] & \calX( \calK_0) \ar[r] & \calX(\calK_{01}). }$$
To complete the proof, it will suffice to show that both of the squares appearing in this diagram
are homotopy pullback squares.

We first treat the square on the left. Consider the diagram
$$ \xymatrix{ \Mod^{\calO}( \Mod^{\calO}_{A}(\calC)) \ar[r] \ar[d] & \MMod^{\calO}( \Mod^{\calO}_{A}(\calC)) \ar[r] \ar[d] & \calX(\calK) \ar[d] \\
\Alg_{\calO}( \Mod^{\calO}_{A}(\calC)) \times \calO^{\otimes} \ar[r] & \PrAlg_{\calO}( \Mod^{\calO}_{A}(\calC)) \ar[r] & 
\calX(\calK_0). }$$
Since the left horizontal maps are categorical equivalences (Remark \ref{user2}), it suffices to show that
the right square is homotopy Cartesian. Let $\calK_{2}$ be the full subcategory of $\calK$ spanned by those diagrams $X \stackrel{\alpha}{\rightarrow} Y \stackrel{\beta}{\rightarrow} Z$ for which $\beta$ is null, and let $\calK_{02}$ denote the full subcategory of $\calK$ spanned by those diagrams
where $\alpha$ and $\beta$ are both null. The algebra $A \in \Alg_{\calO}(\calC)$ determines
a vertex $v$ of $\calX( \calK_{2})$ (and therefore, by restriction, a vertex $v'$ of $\calX( \calK_{02})$).
We have a commutative diagram
$$ \xymatrix{ \MMod^{\calO}( \Mod^{\calO}_{A}(\calC)) \ar[r] \ar[d] & \calX(\calK) \ar[r]^{\theta} \ar[d] & \calX( \calK_2) \ar[d]^{\phi} \\
\PrAlg_{\calO}(\Mod^{\calO}_{A}(\calC) \ar[r] & \calX( \calK_0) \ar[r]^{\theta'} & \calX( \calK_{02}) }$$
where the horizontal maps are fiber sequences (where the fibers are taken over the vertices
$v$ and $v'$, respectively). To show that the left square is a homotopy pullback, it suffices to prove the following:
\begin{itemize}
\item[$(i)$] The maps $\theta$ and $\theta'$ are categorical fibrations of $\infty$-categories.
\item[$(ii)$] The map $\phi$ is a categorical equivalence. 
\end{itemize}

To prove $(i)$, we first show that the simplicial sets $\calX( \calK)$, $\calX( \calK_0)$, 
$\calX( \calK_2)$, and $\calX( \calK_{02})$ are $\infty$-categories. In view of
Proposition \bicatref{saeve}, it will suffice to show that the maps
$$ e_0: \calK \rightarrow \calO^{\otimes} \quad \quad e_0^{0}: \calK_0 \rightarrow \calO^{\otimes} $$
$$ e_0^{2}: \calK_2 \rightarrow \calO^{\otimes} \quad \quad e_0^{02}: \calK_{02} \rightarrow \calO^{\otimes}$$
are flat categorical fibrations. The map $e_0$ can be written as a
$$ \calK \stackrel{e'_0}{\rightarrow} \calK_{\calO} \stackrel{e''_0}{\rightarrow} \calO^{\otimes}$$
where $e'_0$ is given by restriction along the inclusion $\Delta^{ \{0,1\} } \subseteq \Lambda^2_1$ and
$e''_0$ is given by evaluation at $\{0\}$. The map $e''_0$ is a flat categorical fibration by virtue of
our assumption that $\calO^{\otimes}$ is coherent. The map $e'_0$ is a pullback of $e''_0$, and therefore also flat. Applying Corollary \bicatref{stiffum}, we deduce that $e_0$ is flat. The proofs
in the other cases three are similar: the only additional ingredient that is required is the observation
that evaluation at $0$ induces a flat categorical fibration $\calK^{0}_{\calO} \rightarrow \calO^{\otimes}$, which follows from Lemma \ref{preswope}.

To complete the proof of $(i)$, we will show that $\theta$ and $\theta'$ are categorical fibrations.
We will give the proof for the map $\theta$; the case of $\theta'$ is handled similarly. 
We wish to show that $\theta$
has the right lifting property with respect to every trivial cofibration $A \rightarrow B$ in $(\sSet)_{/ \calO^{\otimes}}$. Unwinding
the definitions, we are required to provide solutions to lifting problems of the form
$$ \xymatrix{ (A \times_{ \calO^{\otimes}} \calK)
\coprod_{ A \times_{ \calO^{\otimes}} \calK_2 } ( B \times_{ \calO^{\otimes} }\calK_{2}) \ar[d]^{i} \ar[r] & \calC^{\otimes} \ar[d]^{p} \\
B \times_{ \calO^{\otimes}} \calK \ar[r] & \calO^{\otimes}. }$$
Since $p$ is a categorical fibration, it suffices to prove that the monomorphism $i$ is a categorical equivalence of simplicial sets. In other words, we need to show that the diagram
$$ \xymatrix{ A \times_{ \calO^{\otimes}} \calK_2 \ar[r] \ar[d] & A \times_{ \calO^{\otimes}} \calK \ar[d] \\
B \times_{ \calO^{\otimes}} \calK_{2} \ar[r] & B \times_{ \calO^{\otimes}} \calK }$$
is a homotopy pushout square (with respect to the Joyal model structure). To prove this, it suffices
to show that the vertical maps are categorical equivalences. Since $A \rightarrow B$ is a
categorical equivalence, this follows from Corollary \bicatref{silman} (since the maps
$e_0$ and $e_0^{2}$ are flat categorical fibrations).

We now prove $(ii)$. Let $\calK_{3}$ denote the full subcategory of $\calK$ spanned by those diagrams
$X \stackrel{ \alpha}{\rightarrow} Y \stackrel{\alpha}{\rightarrow} Z$ in $\calO^{\otimes}$ where
$Y \in \calO^{\otimes}_{\seg{0}}$. We have a commutative diagram
$$ \xymatrix{ \calX( \calK_0) \ar[dr] \ar[rr]^{\theta} & & \calX( \calK_{02} ) \ar[dl] \\
& \calX( \calK_3). & }$$
Consequently, to show that $\theta$ is a categorical equivalence, it suffices to show that the diagonal
maps in this diagram are categorical equivalences. We will show that $\calX( \calK_0 ) \rightarrow \calX( \calK_3)$ is a categorical equivalence; the proof for $\calX( \calK_{02}) \rightarrow \calX( \calK_3)$ is similar. Let $X \stackrel{ \alpha}{ \rightarrow} Y \stackrel{ \beta}{\rightarrow } Z$
be an object $K \in \calK_0$, and choose a morphism $\gamma: Y \rightarrow Y_0$ where
$Y_0 \in \calO^{\otimes}$. Since $\beta$ is null, it factors through $\gamma$, and we obtain a commutative diagram
$$ \xymatrix{ X \ar[r]^{\alpha} \ar[d]^{\id} & Y \ar[d]^{ \gamma} \ar[r]^{\beta} & Z \ar[d]^{\id}  \\
X \ar[r]^{ \gamma \circ \alpha} & Y_0 \ar[r] & Z. }$$
We can interpret this diagram as a morphism $\overline{\gamma}: K \rightarrow K_0$ in
$\calK_{0}$. It is not difficult to see that this $\overline{\gamma}$ exhibits $K_0$ as a
$\calK_{3}$-localization of $K$. Consequently, the construction $K \mapsto K_0$ can
be made into a functor $L: \calK_0 \rightarrow \calK_{3}$, equipped with a natural transformation
$t: \id \rightarrow L$. Without loss of generality, we can assume that $L$ and $t$ commute
with the evaluation maps $e_0$ and $e_2$. Composition with $L$ determines a map
$\calX( \calK_{3}) \rightarrow \calX( \calK_0)$, and the transformation $t$ exhibits
this map as a homotopy inverse to the restriction map $\calX( \calK_0) \rightarrow \calX( \calK_{3})$. This completes the proof of $(ii)$.

It remains to show that the diagram
$$ \xymatrix{ \calX(\calK) \ar[d] \ar[r] & \calX(\calK_1) 
\ar[d]^{\psi} \\
\calX( \calK_0) \ar[r] & \calX(\calK_{01}) }$$
is a homotopy pullback square. We first claim that $\psi$ is a categorical equivalence
of $\infty$-categories. The proof is similar to the proof of $(i)$: the only nontrivial point is
to verify that the restriction maps
$$ e_0^{1}: \calK_{1} \rightarrow \calO^{\otimes} \quad \quad e_{0}^{01}: \calK_{01} \rightarrow \calO^{\otimes}$$
are flat categorical fibrations. We will give the proof for $e_{0}^{1}$; the proof for $e_{0}^{01}$ is similar.
We can write $e_0^{1}$ as a composition
$\calK_{1} \rightarrow \calK_{\calO} \rightarrow \calO^{\otimes}$, where the second map
is a flat categorical fibration by virtue of our assumption that $\calO^{\otimes}$ is coherent.
The first map is a pullback of the restriction map $\Fun^{0}( \Delta^1, \calO^{\otimes})
\Fun( \{0\}, \calO^{\otimes})$, where $\Fun^{0}( \Delta^1, \calO^{\otimes})$ denotes the full subcategory of
$\Fun( \Delta^1, \calO^{\otimes})$ spanned by the equivalences in $\calO^{\otimes}$, and therefore a trivial Kan fibration (and, in particular, a flat categorical fibration). Applying Corollary \bicatref{stiffum}, we conclude that $e_0^{1}$ is a flat categorical fibration, as desired.

Since $\psi$ is a categorical fibration of $\infty$-categories and $\calX( \calK_0)$ is an $\infty$-category, we have a homotopy pullback diagram
$$ \xymatrix{ \calX( \calK_0) \times_{ \calX( \calK_{01}) } \calX( \calK_1) \ar[r] \ar[d] & \calX( \calK_1) \ar[d] \\
\calX( \calK_0) \ar[r] & \calX( \calK_{01}). }$$
To complete the proof, it will suffice to show that the restriction map
$\tau: \calX( \calK) \rightarrow \calX( \calK_0) \times_{ \calX( \calK_{01}) } \calX( \calK_1)$
is a categorical equivalence. We will show that $\tau$ is a trivial Kan fibration. Note that
the evaluation map $e_0: \calK \rightarrow \calO^{\otimes}$ is a Cartesian fibration; moreover,
if $K \rightarrow K'$ is an $e_0$-Cartesian moprhism in $\calK$ and $K' \in \calK_0 \coprod_{ \calK_{01} } \calK_{1}$, then $K \in \calK_0 \coprod_{ \calK_{01} } \calK_{1}$. It follows that
$e_0$ restricts to a Cartesian fibration $\calK_0 \coprod_{ \calK_{01} } \calK_{1} \rightarrow \calO^{\otimes}$. In view of Lemma \ref{loopse}, the map $\tau$ will be a trivial Kan fibration provided that the following pair of assertions holds:

\begin{itemize}
\item[$(a)$] Let $F \in \overline{\calX}( \calK)$ be an object lying over $X \in \calO^{\otimes}$ which we will identify with a functor $\{ X \} \times_{ \calO^{\otimes} } \calK \rightarrow \calC^{\otimes}$.
Let $F_0 = F | ( \{X \} \times_{ \calO^{\otimes}} ( \calK_0 \coprod_{ \calK_{01}} \calK_{1})$, and assume
that $F_0 \in \calX( \calK_0) \times_{ \calX( \calK_{01}) } \calX( \calK_1)$. Then $F \in \calX(\calK)$ if and only if $F$ is a $p$-right Kan extension of $F_0$.
\item[$(b)$] Let $F_0 \in  \calX( \calK_0) \times_{ \calX( \calK_{01}) } \calX( \calK_1)$. Then there
exists an extension $F$ of $F_0$ which satisfies the equivalent conditions of $(a)$.
\end{itemize}

To prove these assertions, let us consider an object $X \in \calO^{\otimes}$ and an object
$F_0 \in \calX( \calK_0) \times_{ \calX( \calK_{01}) } \calX( \calK_1)$ lying over $X$. 
Let $\calD =\{ X\} \times_{ \calO^{\otimes}} \calK$ denote the $\infty$-category of diagrams
$$ X \stackrel{\alpha}{\rightarrow} Y \stackrel{\beta}{\rightarrow} Z $$ in
$\calO^{\otimes}$, and define full subcategories $\calD_0$, $\calD_1$, and $\calD_{01}$ similarly.
Let $K$ be an object of $\calD$, corresponding to a decomposition
$X \simeq X_0 \oplus X_1 \oplus X_2$ and a diagram
$$ X_0 \oplus X_1 \oplus X_2 \rightarrow X_0 \oplus X_1 \oplus Y_0 \oplus Y_1
\rightarrow X_0 \oplus Y_0 \oplus Z$$
of semi-inert morphisms in $\calO^{\otimes}$. We have a commutative diagram
$$ \xymatrix{ K \ar[r] \ar[d] & K_1 \ar[d] \\
K_0 \ar[r] & K_{01} }$$
in $\calD$, where $K_0 \in \calD_0$ represents the diagram
$X_0 \oplus X_1 \oplus X_2 \rightarrow Y_0 \oplus Y_1 \rightarrow Y_0 \oplus Z$,
$K_1 \in \calD_1$ represents the diagram $X_0 \oplus X_1 \oplus X_2 \rightarrow X_0 \oplus Y_0 \oplus Z \simeq X_0 \oplus Y_0 \oplus Z$, and $K_{01} \in \calD_{01}$ represents the diagram
$X_0 \oplus X_1 \oplus X_2 \rightarrow Y_0 \oplus Z \rightarrow Y_0 \oplus Z$.
This diagram exhibits $K_0$, $K_1$, and $K_{01}$ as initial objects of
$(\calD_{0})_{K/}$, $(\calD_1)_{K/}$, and $(\calD_{01})_{K/}$, respectively.
Applying Theorem \toposref{hollowtt}, we conclude that the induced map
$\Lambda^2_{2} \rightarrow ( \calD_0 \coprod_{ \calD_{01} } \calD_1)_{K/}$ is
the opposite of a cofinal map, so that an extension $F$ of $F_0$ is a $p$-right Kan
extension of $F_0$ at $K$ if and only if the diagram
$$ \xymatrix{ F(K) \ar[r] \ar[d] & F_0(K_1) \ar[d] \\
F_0( K_0) \ar[r] & F_0( K_{01} ) }$$
is a $p$-limit diagram in $\calC^{\otimes}$. Choose equivalences
$F_0( K_0) = y_0 \oplus z$, $F_0(K_1) = x'_0 \oplus y'_0 \oplus z'$, and
$F_0(K_{01}) = y''_0 \oplus z''$. If we let $K' \in \calD$ denote the object corresponding
to the diagram $X_0 \oplus X_1 \oplus X_2 \rightarrow Y_0 \simeq Y_0$ and apply
our assumption that $F_0 | \calD_{0} \in \calX( \calK_0)$
to the commutative diagram $$ \xymatrix{ K_0 \ar[rr] \ar[dr] & & K_{01} \ar[dl] \\
& K', & }$$
then we deduce that $F_0( K_0) \rightarrow F_0( K_{01})$ induces an equivalence
$y_0 \simeq y''_0$. Similarly, the assumption that $F_0 | \calD_1 \in \calX( \calK_1)$ guarantees that $F_0( K_1) \rightarrow F_0( K_{01})$ is inert, so that $y'_0 \simeq y''_0$ and $z' \simeq z''$.
It follows from Lemma \ref{helpless} that the diagram $F_0(K_1) \rightarrow F_0( K_{01}) \leftarrow F_0(K_0)$ admits a $p$-limit (covering the evident map $\Delta^1 \times \Delta^1 \rightarrow \calO^{\otimes}$), so that $F_0 | (\calD_0 \coprod_{ \calD_{01} } \calD_1)_{K/}$ also admits
a $p$-limit (covering the map $( \calD_0 \coprod_{ \calD_{01} } \calD_1)^{\triangleleft}_{K/} \rightarrow \calO^{\otimes})$); assertion $(b)$ now follows from Lemma \toposref{kan2}. Moreover,
the criterion of Lemma \ref{helpless} gives the following version of $(a)$:

\begin{itemize}
\item[$(a')$] An extension $F \in \overline{\calX}(\calK)$ of $F_0$ is a $p$-right Kan extension of $F_0$ at $K$ if and only if, for every object $K \in \calD$ as above, the maps
$F(K) \rightarrow F_0(K_0)$ and $F(K) \rightarrow F_0(K_1)$ induce an equivalence
$F(K) \simeq x'_0 \oplus y_0 \oplus z$ in $\calC^{\otimes}$.
\end{itemize} 

To complete the proof, it will suffice to show that the criterion of $(a')$ holds if and only if
$F \in \calX(\calK)$. We first prove the ``if'' direction. Fix $K \in \calD$, so that we have
an equivalence $F(K) \simeq \overline{x}_0 \oplus \overline{y}_0 \oplus \overline{z}$.
Since $K \rightarrow K_0$ is inert, the assumption that $F \in \calX(\calK)$ implies that
$F(K) \rightarrow F(K_0)$ is inert, so that $\overline{y}_0 \simeq y_0$ and
$\overline{z} \simeq z$. Let $K'' \in \calD$ be the diagram
$X_0 \oplus X_1 \oplus X_2 \rightarrow Y_0 \simeq Y_0$, so that we have a
commutative diagram
$$ \xymatrix{ K \ar[rr] \ar[dr] & & K_1 \ar[dl] \\
& K'' & }$$
in which the diagonal maps are inert. It follows that the morphisms
$F(K) \rightarrow F(K'') \leftarrow F(K_1)$ are inert, so that the map
$F(K) \rightarrow F(K_1)$ induces an equivalence $\overline{y}_0 \simeq y'_0$.

We now prove the ``only if'' direction. Assume that $F \in \overline{\calX}(\calK)$ is an
extension of $F_0$ which satisfies the criterion given in $(a')$; we will show that
$F$ carries inert morphisms in $\calD \subseteq \calK$ to inert morphisms in $\calC^{\otimes}$.
Let $K \rightarrow L$ be an inert morphism in $\calD$, where $K$ is as above and
$L$ corresponds to a diagram $X \rightarrow Y' \rightarrow Z'$; we wish to prove that the induced
map $F(K) \rightarrow F(L)$ is inert. Let $\seg{n}$ denote the image of $Z'$ in $\Nerve(\FinSeg)$,
and choose inert morphisms $Z' \rightarrow Z'_i$ lying over $\Colp{i}: \seg{n} \rightarrow \seg{1}$ for $1 \leq i \leq n$. Each of the induced maps $Y' \rightarrow Z_i$ factors as a composition
$Y' \rightarrow Y'_{i} \rightarrow Z_i$, where the first map is inert and the second is active.
Let $L_i \in \calD$ denote the diagram $X \rightarrow Y'_{i} \rightarrow Z'_{i}$. To show
that $F(K) \rightarrow F(L)$ is inert, it will suffice to show that the maps
$F(K) \rightarrow F(L_i) \leftarrow F(L)$ are inert for $1 \leq i \leq n$. Replacing
$L$ by $L_i$ (and possibly replacing $K$ by $L$), we may reduce to the case where
$n = 1$ and the map $Y' \rightarrow Z'$ is active. There are two cases to consider:

\begin{itemize}
\item The map $Y' \rightarrow Z'$ is an equivalence. In this case, the map $K \rightarrow L$
factors as a composition $K \rightarrow K_1 \rightarrow L$. Since $F_0 | \calD_1 \in
\calX(\calK_1)$, the map $F(K_1) \rightarrow F(L)$ is inert. Consequently, the assertion
that $F(K) \rightarrow F(L)$ is inert follows our assumption that $F(K) \rightarrow F(K_1)$
induces an equivalence $\overline{y}_0 \rightarrow y'_0$.

\item The map $Y' \rightarrow Z'$ is null (so that $Y' \in \calO^{\otimes}_{\seg{0}}$). In this case,
the map $K \rightarrow L$ factors as a composition $K \rightarrow K_0 \rightarrow L$.
Assumption $(a')$ guarantees that $F(K) \rightarrow F(K_0)$ is inert, and the assumption
that $F_0 | \calK_0 \in \calX(\calK_0)$ guarantees that $F(K_0) \rightarrow F(L)$ is inert.
\end{itemize}
\end{proof}

\subsection{Limits of Modules}\label{com34}

Let $\calC$ be a symmetric monoidal category and let $A$ be a commutative algebra object
of $\calC$. Suppose we are given a diagram $\{ M_{\alpha} \}$ in the category of $A$-modules,
and let $M = \varprojlim M_{\alpha}$ be the limit of this diagram in the category $\calC$.
The collection of maps
$$ A \otimes M \rightarrow A \otimes M_{\alpha} \rightarrow M_{\alpha}$$
determines a map $A \otimes M \rightarrow M$, which endows $M$ with the structure of an
$A$-module. Moreover, we can regard $M$ also as a limit of the diagram $\{ M_{\alpha} \}$ in the category of $A$-modules.

Our goal in this section is to prove an analogous result in the $\infty$-categorical setting, for algebras over an arbitrary coherent $\infty$-operad. We can state our main result as follows:

\begin{theorem}\label{maincand}
Let $q: \calC^{\otimes} \rightarrow \calO^{\otimes}$ be a fibration of $\infty$-operads, where
$\calO^{\otimes}$ is coherent. Suppose we are given a commutative diagram
$$ \xymatrix{ K \ar[r]^{p} \ar[d] & \Mod^{\calO}(\calC)^{\otimes} \ar[d]^{\psi} \\
K^{\triangleleft} \ar[r]^-{ \overline{p}_0} \ar@{-->}[ur]^{\overline{p}} & \Alg_{\calO}(\calC) \times \calO^{\otimes} }$$
such that the underlying map $K^{\triangleleft} \rightarrow \calO^{\otimes}$
takes some constant value $X \in \calO$, and the lifting problem 
$$ \xymatrix{ K \ar[d] \ar[r]^{p'} & \calC^{\otimes} \ar[d]^{q} \\
K^{\triangleleft} \ar[r] \ar@{-->}[ur]^{ \overline{p}'} & \calO^{\otimes} }$$
admits a solution, where $\overline{p}$ is a $q$-limit diagram. Then: 
\begin{itemize}
\item[$(1)$] There exists a map $\overline{p}$ making the original diagram commute, such that
$\delta \circ \overline{p}$ is a $q$-limit diagram in $\calC^{\otimes}$ (here
$\delta: \Mod^{\calO}(\calC)^{\otimes} \rightarrow \calC^{\otimes}$ denotes the map
given by composition with the diagonal embedding $\calO^{\otimes} \rightarrow
\calK_{\calO} \subseteq \Fun( \Delta^1, \calO^{\otimes})$.

\item[$(2)$] Let $\overline{p}$ be an arbitrary map making the above diagram commute.
Then $\overline{p}$ is a $\psi$-limit diagram if and only if $\delta \circ \overline{p}$ is a
$q$-limit diagram.
\end{itemize}
\end{theorem}

\begin{corollary}\label{sus}
Let $q: \calC^{\otimes} \rightarrow \calO^{\otimes}$ be a fibration of $\infty$-operads, where
$\calO^{\otimes}$ is coherent, and let
$A \in \Alg_{\calO}(\calC)$ and $X \in \calO$. Suppose we are given a diagram
$p: K \rightarrow \Mod^{\calO}_{A}(\calC)^{\otimes}_{X}$ such that the induced map
$p': K \rightarrow \calC^{\otimes}_{X}$ can be extended to a $q$-limit diagram
$\overline{p}': K^{\triangleleft} \rightarrow \calC^{\otimes}_{X}$. Then:
\begin{itemize}
\item[$(1)$] There exists an extension $\overline{p}: K^{\triangleleft} \rightarrow \Mod^{\calO}_{A}(\calC)^{\otimes}_{X}$ of $p$ such that the induced map
$K^{\triangleleft} \rightarrow \calC^{\otimes}$ is a $q$-limit diagram.

\item[$(2)$] Let $\overline{p}: K^{\triangleleft} \rightarrow \Mod^{\calO}_{A}(\calC)^{\otimes}_{X}$
be an arbitrary extension of $p$. Then $\overline{p}$ is a limit diagram if and only if
it induces a $q$-limit diagram $K^{\triangleleft} \rightarrow \calC^{\otimes}$.
\end{itemize}
\end{corollary}

\begin{corollary}\label{sabel}
Let $q: \calC^{\otimes} \rightarrow \calO^{\otimes}$ be a fibration of $\infty$-operads, where
$\calO^{\otimes}$ is coherent. Then for each $X \in \calO$, the functor
$\phi: \Mod^{\calO}(\calC)^{\otimes}_{X} \rightarrow \Alg^{\calO}(\calC)$ is a Cartesian fibration.
Moreover, a morphism $f$ in $\Mod^{\calO}(\calC)_{X}$ is $\phi$-Cartesian if and only if
its image in $\calC^{\otimes}_{X}$ is an equivalence.
\end{corollary}

\begin{proof}
Apply Corollary \ref{sus} in the case $K = \Delta^0$.
\end{proof}

\begin{corollary}\label{postsabel}
Let $q: \calC^{\otimes} \rightarrow \calO^{\otimes}$ be a fibration of $\infty$-operads, where
$\calO^{\otimes}$ is coherent. Then:
\begin{itemize}
\item[$(1)$] The functor
$\phi: \Mod^{\calO}(\calC)^{\otimes} \rightarrow \Alg^{\calO}(\calC)$ is a Cartesian fibration.
\item[$(2)$] A morphism $f \in \Mod^{\calO}(\calC)^{\otimes}$ is $\phi$-Cartesian if and only if
its image in $\calC^{\otimes}$ is an equivalence.
\end{itemize}
\end{corollary}

\begin{proof}
Let $M \in \Mod^{\calO}_{A}(\calC)^{\otimes}$ and let $\seg{n}$ denote its image in
$\Nerve(\FinSeg)$. Suppose we are given a morphism $f_0: A' \rightarrow A$ in
$\Alg^{\calO}(\calC)$; we will to construct a $\phi$-Cartesian morphism
$f: M' \rightarrow M$ lifting $f_0$. 

Choose inert morphisms $g_i: M \rightarrow M_i$ in
$\Mod^{\calO}(\calC)^{\otimes}$ lying over $\Colp{i}: \seg{n} \rightarrow \seg{1}$ for
$1 \leq i \leq n$. These maps determine a diagram $F: \nostar{n}^{\triangleleft} \rightarrow
\Mod^{\calO}_{A}( \calC)^{\otimes}$. Let $X_i$ denote the image of $M_i$ in $\calO$, and let
$\phi_i: \Mod^{\calO}(\calC)_{X_i} \rightarrow \Alg^{\calO}(\calC)$ be the restriction of $\phi$.
Using Corollary \ref{sabel}, we can choose $\phi_i$-Cartesian morphisms
$f_i: M'_i \rightarrow M_i$ in $\Mod^{\calO}(\calC)_{X_i}$ lying over $f_0$, whose
images in $\calC_{X_i}$ are equivalences. Since $q: \Mod^{\calO}(\calC)^{\otimes} \rightarrow
\Alg^{\calO}(\calC) \times \calO^{\otimes}$ is a $\Alg^{\calO}(\calC)$-family of $\infty$-operads, we can choose a $q$-limit diagram $F': \nostar{n}^{\triangleleft} \rightarrow \Mod^{\calO}_{A'}(\calC)^{\otimes}$
with $F'(i) = M'_i$ for $1 \leq i \leq n$, where $F'$ carries the cone point of $\nostar{n}^{\triangleleft}$ to
$M' \in \Mod^{\calO}_{A'}(\calC)^{\otimes}_{X}$. Using the fact that $F$ is a $q$-limit diagram,
we get a natural transformation of functors $F' \rightarrow F$, which we may view as
a diagram $H: \nostar{n}^{\triangleleft} \times \Delta^1 \rightarrow \Mod^{\calO}(\calC)^{\otimes}$. 

Let $v$ denote the cone point of $\nostar{n}^{\triangleleft}$, and let $f = H | \{v\} \times \Delta^1$. Since each composition
$$ \{i\} \times \Delta^1 \stackrel{H}{\rightarrow} \Mod^{\calO}(\calC)^{\otimes} \rightarrow
\calC^{\otimes}$$
is an equivalence for $1 \leq i \leq n$, the assumption that $\calC^{\otimes}$ is an $\infty$-operad guarantees also that the image of $f$ in $\calC^{\otimes}$ is an equivalence. We will prove
that $f$ is a $\phi$-coCartesian lift of $f_0$. In fact, we will prove the slightly stronger assertion that $f$ is
$q$-Cartesian. Since the inclusion $\{v\} \subseteq \nostar{n}^{\triangleleft}$ is the opposite of
a cofinal map, it will suffice to show that $H | ( \nostar{n}^{\triangleleft} \times \{1\})^{\triangleleft}$
is a $q$-limit diagram. Since $H | (\nostar{n}^{\triangleleft} \times \{1\})$ is a $q$-right Kan extension
of $H | ( \nostar{n} \times \{1\} )$, it will suffice to show that the restriction 
$H | ( \nostar{n} \times \{1\} )^{\triangleleft}$ is a $q$-limit diagram (Lemma \toposref{kan0}). 
Note that $H | ( \nostar{n} \times \Delta^1)$ is a $q$-right Kan extension
of $H | ( \nostar{n} \times \{1\})$ (this follows from the construction, since the maps
$f_i$ are $\phi_i$-Cartesian and therefore also $q$-Cartesian, by virtue of Theorem \ref{maincand}).
Using Lemma \toposref{kan0} again, we are reduced to showing that $H| ( \nostar{n} \times \Delta^1)^{\triangleleft}$ is a $q$-limit diagram. Since the inclusion $\nostar{n} \times \{0\} \subseteq \nostar{n} \times \Delta^1$ is the opposite of a cofinal map, it suffices to show that
$F' = H | ( \nostar{n}^{\triangleleft} \times \{0\} )$ is a $q$-limit diagram, which follows from our assumption.

The above argument shows that for every $M \in \Mod^{\calO}_{A}(\calC)^{\otimes}$ 
and every morphism $f_0: A' \rightarrow A$ in $\Alg^{\calO}(\calC)$, there exists
a $\phi$-Cartesian morphism $f: M' \rightarrow M$ lifting $f_0$ whose image
in $\calC^{\otimes}$ is an equivalence. This immediately implies $(1)$, and
the ``only if'' direction of $(2)$ follows from the uniqueness properties of
Cartesian morphisms. To prove the ``if'' direction of $(2)$, suppose that
$g: M'' \rightarrow M$ is a lift of $f_0$ whose image in $\calC^{\otimes}$ is an equivalence, and let
$f: M' \rightarrow M$ be as above. Since $f$ is $\phi$-Cartesian, we have a commutative diagram
$$ \xymatrix{ & M' \ar[dr]^{f} & \\
M'' \ar[ur]^{h} \ar[rr]^{g} & & M; }$$
to prove that $g$ is $\phi$-Cartesian it will suffice to show that $h$ is an equivalence.
Since $\Mod^{\calO}_{A'}(\calC)^{\otimes}$ is an $\infty$-operad, it suffices to show that
each of the maps $h_i = \Colll{i}{!} (h)$ is an equivalence in
$\Mod^{\calO}_{A'}(\calC)$, for $1 \leq i \leq n$. This follows from Corollary \ref{sabel}, since
each $h_i$ maps to an equivalence in $\calC$.
\end{proof}

\begin{corollary}
Let $\calO^{\otimes}$ be a coherent $\infty$-operad, and let
$q: \calC^{\otimes} \rightarrow \calO^{\otimes}$ be a $\calO$-monoidal $\infty$-category.
Let $X \in \calO$, and suppose we are given a commutative diagram
$$ \xymatrix{ K \ar[r]^-{p} \ar[d] & \Mod^{\calO}(\calC)^{\otimes}_{X} \ar[d]^{\psi_X} \\
K^{\triangleleft} \ar[r]^-{ \overline{p}_0} \ar@{-->}[ur]^{\overline{p}} & \Alg_{\calO}(\calC) }$$
such that the induced diagram $K \rightarrow \calC^{\otimes}_{X}$ admits a limit.
Then there extension $\overline{p}$ of $p$ (as indicated in the diagram) which is
a $\psi_{X}$-limit diagram. Moreover, an arbitrary extension $\overline{p}$ of $p$
(as in the diagram) is a $\psi_{X}$-limit if and only if it induces a limit diagram
$K^{\triangleleft} \rightarrow \calC^{\otimes}_{X}$. 
\end{corollary}

\begin{proof}
Combine Corollary \ref{sus} with Corollary \toposref{pannaheave}.
\end{proof}

\begin{corollary}\label{splasher}
Let $\calO^{\otimes}$ be a coherent $\infty$-operad, let
$q: \calC^{\otimes} \rightarrow \calO^{\otimes}$ be a $\calO$-monoidal $\infty$-category, and let
$X \in \calO$. Assume that the $\infty$-category $\calC^{\otimes}_{X}$ admits $K$-indexed limits, for some simplicial set $K$.
Then:
\begin{itemize}
\item[$(1)$] For every algebra object $A \in \Alg_{\calO}(\calC)$, the $\infty$-category
$\Mod^{\calO}_{A}(\calC)^{\otimes}_{X}$ admits $K$-indexed limits.
\item[$(2)$] A functor $\overline{p}: K^{\triangleleft} \rightarrow \Mod^{\calO}_{A}(\calC)^{\otimes}_{X}$ is a limit diagram if and only if it induces a limit diagram $K^{\triangleleft} \rightarrow \calC^{\otimes}_{X}$.
\end{itemize}
\end{corollary}

We now turn to the proof of Theorem \ref{maincand}. First, choose an inner anodyne map
$K \rightarrow K'$, where $K'$ is an $\infty$-category. Since $\Alg_{\calO}(\calC) \times \calO^{\otimes}$ is an $\infty$-category and $\psi$ is a categorical fibration, we can extend our commutative diagram as indicated:
$$ \xymatrix{ K \ar[r] \ar[d] & K' \ar[r] \ar[d] & \Mod^{\calO}(\calC)^{\otimes} \ar[d]^{\psi} \\
K^{\triangleleft} \ar[r]  & {K'}^{\triangleleft} \ar[r] & \Alg_{\calO}(\calC) \times \calO^{\otimes}. }$$
Using Proposition \toposref{princex}, we see that it suffices to prove Theorem \ref{maincand} after replacing $K$ by $K'$. We may therefore assume that $K$ is an $\infty$-category. In this case,
the desired result is a consequence of the following:

\begin{proposition}\label{canda}
Let $q: \calC^{\otimes} \rightarrow \calO^{\otimes}$ be a fibration of $\infty$-operads, where
$\calO^{\otimes}$ is coherent. Let $K$ be an $\infty$-category.
Suppose we are given a commutative diagram
$$ \xymatrix{ K \ar[r]^-{p} \ar[d] & \MMod^{\calO}(\calC)^{\otimes} \ar[d]^{\psi} \\
K^{\triangleleft} \ar[r]^-{ \overline{p}_0} \ar@{-->}[ur]^{\overline{p}} & \PrAlg_{\calO}(\calC), }$$
where the induced diagram $K^{\triangleleft} \rightarrow \calO^{\otimes}$ is the constant map taking some value $C \in \calO$. Assume that:
\begin{itemize}
\item[$(\ast)$] 
The induced lifting problem
$$ \xymatrix{ K \ar[d] \ar[r]^{p'} & \calC^{\otimes} \ar[d]^{q} \\
K^{\triangleleft} \ar[r] \ar@{-->}[ur]^{ \overline{p}'} & \calO^{\otimes} }$$
admits a solution, where $\overline{p}$ is a $q$-limit diagram.
\end{itemize}
 Then: 
\begin{itemize}
\item[$(1)$] There exists a map $\overline{p}$ making the original diagram commute, such that
$\delta \circ \overline{p}$ is a $q$-limit diagram in $\calC^{\otimes}$ (here
$\delta: \MMod^{\calO}(\calC)^{\otimes} \rightarrow \calC^{\otimes}$ denotes the map
given by composition with the diagonal embedding $\calO^{\otimes} \rightarrow
\calK_{\calO} \subseteq \Fun( \Delta^1, \calO^{\otimes})$.

\item[$(2)$] Let $\overline{p}$ be an arbitrary map making the above diagram commute.
Then $\overline{p}$ is a $\psi$-limit diagram if and only if $\delta \circ \overline{p}$ is a
$q$-limit diagram.
\end{itemize}
\end{proposition}

The proof of Proposition \ref{canda} will require some preliminaries. We first
need the following somewhat more elaborate version of Proposition \bicatref{skitter}:

\begin{proposition}\label{skitter2}
Suppose we are given a diagram of $\infty$-categories $X \stackrel{\phi}{\rightarrow} Y \stackrel{\pi}{\rightarrow} Z$ where $\pi$ is a flat categorical fibration and $\phi$ is a categorical fibration.
Let $Y' \subseteq Y$ be a full subcategory, let $X' = X \times_{Y} Y'$, let
$\pi' = \pi| Y'$, and let $\psi: \pi_{\ast} X \rightarrow \pi'_{\ast} X'$ be the canonical map.
(See Notation \bicatref{cabbel}.) Let $K$ be an $\infty$-category and
$\overline{p}_0: K^{\triangleleft} \rightarrow \pi'_{\ast} X'$ a diagram. 
Assume that the following conditions are satisfied:
\begin{itemize}
\item[$(i)$] The full subcategory $Y' \times_{Z} K^{\triangleleft} \subseteq Y \times_{Z} K^{\triangleleft}$ is a cosieve on $Y$.
\item[$(ii)$] For every object $y \in Y'$ and every morphism
$f: z \rightarrow \pi(y)$ in $Z$, there exists a $\pi$-Cartesian morphism
$\overline{f}: \overline{z} \rightarrow y$ in $Y'$ such that $\pi( \overline{f} ) = f$.
\item[$(iii)$] Let $\pi''$ denote the projection map $K^{\triangleleft} \times_{Z} Y
\rightarrow K^{\triangleleft}$. Then $\pi''$ is a coCartesian fibration.
\item[$(iv)$] Let $v$ denote the cone point of $K^{\triangleleft}$, let
$\calC = {\pi''}^{-1} \{v\}$, and let $\calC' = \calC \times_{Y} Y'$. Then
$\calC'$ is a localization of $\calC$.
\end{itemize}

Condition $(iii)$ implies that there is a map $\delta': 
K^{\triangleleft} \times \calC \rightarrow K^{\triangleleft} \times_{Z} Y$
which is the identity on $\{v\} \times \calC$ and carries 
carries $e \times \{C\}$ to a $\pi''$-coCartesian edge of
$K^{\triangleleft} \times_{Z} Y$, for each edge $e$ of $K^{\triangleleft}$ and each object $C$ of $\calC$.
Condition $(iv)$ implies that there is a map $\delta'': \calC \times \Delta^1 \rightarrow \calC$
such that $\delta'' | \calC \times \{0\} = \id_{\calC}$ and $\delta''| \{C\} \times \Delta^1$
exhibits $\delta''(C,1)$ as a $\calC'$-localization of $C$, for each $C \in \calC$.
Let $\delta$ denote the composition
$$ K^{\triangleleft} \times \calC \times \Delta^1 \stackrel{\delta''}{\rightarrow}
K^{\triangleleft} \times \calC \stackrel{\delta'}{\rightarrow} K^{\triangleleft}
\times_{Z} Y.$$
Then:
\begin{itemize}
\item[$(1)$] Let $\overline{p}: K^{\triangleleft} \rightarrow \pi_{\ast} X$ be a map lifting
$\overline{p}_0$, corresponding to a functor $\overline{F}: K^{\triangleleft} \times_{Z} Y \rightarrow X$.
Suppose that for each $C \in \calC$, the induced map
$$ K^{\triangleleft} \times \{C\} \times \Delta^1 \hookrightarrow
K^{\triangleleft} \times \calC \stackrel{\delta}{\rightarrow} K^{\triangleleft} \times_{Z} Y
\stackrel{\overline{F}}{\rightarrow} X$$
is a $\phi$-limit diagram. Then $\overline{p}$ is a $\psi$-limit diagram.

\item[$(2)$] Suppose that $p: K^{\triangleleft} \rightarrow \pi_{\ast} X$ is a map lifting
$p_0 = \overline{p}_0 | K$, corresponding to a functor
$F: (K^{\triangleleft} \times_{Z} Y') \coprod_{ K \times_{Z} Y'}
(K \times_{Z} Y) \rightarrow X$. Assume furthermore that for each
$C \in \calC$, the induced map
\begin{eqnarray*} (K^{\triangleleft} \times \{C\} \times \{1\})
\coprod_{ K \times \{C\} \times \{1\}} (K \times \{C\} \times \Delta^1) & \rightarrow &
(K^{\triangleleft} \times \calC \times \{1\})
\coprod_{ K \times \calC \times \{1\}} (K \times \calC \times \Delta^1) \\
& \stackrel{\delta}{\rightarrow} & (K^{\triangleleft} \times_{Z} Y') \coprod_{
K \times_{Z} Y'} (K \times_{Z} Y ) \\
& \stackrel{F}{\rightarrow} & X \end{eqnarray*}
can be extended to a $\psi$-limit diagram lifting the map
$$K^{\triangleleft} \times \{C\} \times \Delta^1 \hookrightarrow K^{\triangleleft} \times \calC \times \Delta^1
\stackrel{\delta}{\rightarrow} K^{\triangleleft} \times_{Z} Y \stackrel{\pi''}{\rightarrow} Y.$$ 
Then there exists an extension $\overline{p}: K^{\triangleleft} \rightarrow \pi_{\ast} X$
of $p$ lifting $\overline{p}_0$ which satisfies condition $(1)$.
\end{itemize}
\end{proposition}

\begin{proof}
Let $W = K^{\triangleleft} \times_{Y} Z$ and let
$W_0$ denote the coproduct
$( K^{\triangleleft} \times_{Z} Y') \coprod_{ K \times_{Z} Y'}
(K \times_{Z} Y)$; condition $(i)$ allows us to identify
$W_0$ with a full subcategory of $W$. Let $\overline{p}: K^{\triangleleft} \rightarrow \pi_{\ast} X$ satisfy the condition described
in $(1)$, corresponding to a functor $\overline{F}: W \rightarrow X$. 
In view of assumptions
$(i)$, $(ii)$, and Proposition \bicatref{critter}, it will suffice to show that
$\overline{F}$ is a $\phi$-right Kan extension of 
$F = \overline{F} | W_0$. 
Pick an object $C \in \calC$; we wish to show that $\overline{F}$ is a $\phi$-right Kan extension of
$F$ at $C$. In other words, we wish to show that the map
$$ (W_0 \times_{W} W_{C/})^{\triangleleft} \rightarrow W
\stackrel{\overline{F}}{\rightarrow} X$$
is a $\phi$-limit diagram. Restricting $\delta$, we obtain a map
$K^{\triangleleft} \times \{C\} \times \Delta^1 \rightarrow W$, which we can identify
with a map 
$$ s: (K^{\triangleleft} \times \{C\} \times \{1\})
\coprod_{ K \times \{C\} \times \{1\} } (K \times \{C\} \times \Delta^1) \rightarrow 
W' \times_{W} W_{C/}.$$
Since $\overline{p}$ satisfies $(1)$, it will suffice to show that
$s^{op}$ is cofinal. We have a commutative diagram
$$ \xymatrix{ (K^{\triangleleft} \times \{C\} \times \{1\})
\coprod_{ K \times \{C\} \times \{1\} } (K \times \{C\} \times \Delta^1) \ar[dr]^-{\theta} \ar[rr]^-{s} & & W' \times_{W} W_{C/} \ar[dl]^{\theta'} \\
& K^{\triangleleft} & }$$
The map $\theta$ is evidently a coCartesian fibration, and $\theta'$ is a coCartesian fibration
by virtue of assumptions $(i)$ and $(iii)$. Moreover, the map $s$ carries $\theta$-coCartesian
edges to $\theta'$-coCartesian edges. Invoking Lemma \stableref{camber}, we are reduced to showing
that for each vertex $k$ of $K^{\triangleleft}$, the map of fibers $s_{k}^{op}$ is cofinal.
If $k = v$ is the cone point of $K^{\triangleleft}$, then we are required to show that
$s$ carries $\{v\} \times \{C \} \times \{1\}$ to an initial object of $\calC'_{C/}$: this follows
from the definition of $\delta'$. If $k \neq v$, then we are required to show that
$s$ carries $K^{\triangleleft} \times \{C\} \times \{0\}$ to an initial object of
$W_{C/} \times_{K^{\triangleleft}} \{k\}$, which follows from our assumption that
$\delta$ carries $\{v\}^{\triangleleft} \times \{C\} \times \{0\}$ to a $\pi''$-coCartesian edge of $W$.
This completes the proof of $(1)$.

We now prove $(2)$. The diagram $p$ gives rise to a map $F: W_0 \rightarrow X$ fitting into a commutative diagram
$$ \xymatrix{ W_0 \ar[r]^{F} \ar[d] & X \ar[d]^{\phi} \\
W \ar[r] \ar@{-->}[ur]^{\overline{F} } & Y. }$$
The above argument shows that a dotted arrow $\overline{F}$ as indicated will correspond to a map $\overline{p}: K^{\triangleleft} \rightarrow \pi_{\ast} X$ satisfying $(1)$ if and only if
$\overline{F}$ is a $\phi$-right Kan extension of $F$. In view of Lemma \toposref{kan2}, the
existence of such an extension is equivalent to the requirement that for each $C \in \calC$, 
the diagram 
$$ W_0 \times_{W} W_{C/} \rightarrow W_0
\stackrel{F}{\rightarrow} X$$
can be extended to a $\phi$-limit diagram lifting the map
$$ (W_0 \times_{W} W_{C/})^{\triangleleft} \rightarrow W
\rightarrow Y.$$
This follows from the hypothesis of part $(2)$ together with the cofinality of the map
$s^{op}$ considered in the proof of $(1)$.
\end{proof}

\begin{definition}\label{capper}
Let $\seg{n}$ be an object of $\FinSeg$. A {\it splitting} of $\seg{n}$ is a pair of inert morphisms
$\alpha: \seg{n} \rightarrow \seg{n_0}$, $\beta: \seg{n} \rightarrow \seg{n_1}$ with the property
that the map $( \alpha^{-1} \coprod \beta^{-1}): \nostar{n_0} \coprod \nostar{n_1} \rightarrow \nostar{n}$ is a bijection.

More generally, let $K$ be a simplicial set. A {\it splitting} of a diagram
$p: K \rightarrow \Nerve(\FinSeg)$ is a pair of natural transformations
$\alpha: p \rightarrow p_0$, $\beta: p \rightarrow p_1$ with the following property:
for every vertex $k$ of $K$, the morphisms $\alpha_{k}: p(k) \rightarrow p_0(k)$
and $\beta_{k}: p(k) \rightarrow p_1(k)$ determine a splitting of $p(k)$.

We will say that a natural transformation $\alpha: p \rightarrow p_0$ of diagrams
$p,p_0: K \rightarrow \Nerve(\FinSeg)$ {\it splits} if there exists another natural
transformation $\beta: p \rightarrow p_1$ which gives a splitting of $p$.
\end{definition}

\begin{remark}
Let $\alpha: p \rightarrow p_0$ be a natural transformation of diagrams $p,p_0: K \rightarrow \Nerve(\FinSeg)$. If $\alpha$ splits, then the natural transformation $\beta: p \rightarrow p_1$ which provides the splitting of $p$ is well-defined up to (unique) equivalence. Moreover,
a bit of elementary combinatorics shows that $\alpha$ splits if and only if it satisfies
the following conditions:
\begin{itemize}
\item[$(1)$] The natural transformation $\alpha$ is {\it inert}: that is, for each vertex
$k \in K$, the map $\alpha_{k}: p(k) \rightarrow p_0(k)$ is an inert morphism in $\FinSeg$.
\item[$(2)$] For every edge $e: x \rightarrow x'$ in $K$, consider the diagram
$$ \xymatrix{ \seg{n} \ar[d]^{p(e)} \ar[r]^{\alpha_x} & \seg{n_0} \ar[d]^{p_0(e)} \ar[d] \\
\seg{m} \ar[r]^{\alpha_{x'}} & \seg{m_0} }$$
in $\FinSeg$ obtained by applying $\alpha$ to $e$. Then $p(e)$ carries $(\alpha^{-1}_{x} \nostar{n_0})_{\ast} \subseteq \seg{n}$
into $( \alpha_{x'}^{-1} \nostar{m_0})_{\ast} \subseteq \seg{m}$.
\end{itemize}
\end{remark}

\begin{definition}
Let $q: \calO^{\otimes} \rightarrow \Nerve(\FinSeg)$ be an $\infty$-operad. We will say that a natural transformation $\alpha: p \rightarrow p_0$ of diagrams $p,p_0: K \rightarrow \calO^{\otimes}$
is {\it inert} if the induced map $\alpha_{k}: p(k) \rightarrow p_0(k)$ is an inert morphism
in $\calO^{\otimes}$ for every vertex $k \in K$. 

A {\it splitting} of $p: K \rightarrow \calO^{\otimes}$ is a pair of inert natural transformations
$\alpha: p \rightarrow p_0$, $\beta: p \rightarrow p_1$ such that the induced transformations
$q \circ p_0 \leftarrow q \circ p \rightarrow q \circ p_1$ determine a splitting
of $q \circ p: K \rightarrow \Nerve(\FinSeg)$, in the sense of Definition \ref{capper}.

We will say that an inert natural transformation $\alpha: p \rightarrow p_0$ is {\it split}
if there exists another inert natural transformation $\beta: p \rightarrow p_1$ such
that $\alpha$ and $\beta$ are a splitting of $p$. In this case, we will say that
$\beta$ is a {\it complement} to $\alpha$.
\end{definition}

\begin{lemma}
Let $q: \calO^{\otimes} \rightarrow \Nerve(\FinSeg)$ and let $\alpha: p \rightarrow p_0$ be an
inert natural transformation of diagrams $p,p_0: K \rightarrow \calO^{\otimes}$. The following conditions are equivalent:
\begin{itemize}
\item[$(1)$] The natural transformation $\alpha$ is split: that is, there exists a complement
$\beta: p \rightarrow p_1$ to $\alpha$.
\item[$(2)$] The natural transformation $\alpha$ induces a split natural transformation 
$\overline{\alpha}: q \circ p \rightarrow q \circ p_0$.
\end{itemize}
Moreover, if these conditions are satisfied, then $\beta$ is determined uniquely up to equivalence.
\end{lemma}

\begin{proof}
The implication $(1) \Rightarrow (2)$ is clear: if $\beta: p \rightarrow p_1$ is a complement to
$\alpha$, then the induced transformations $q \circ p_0 \leftarrow q \circ p \rightarrow q \circ p_1$
form a splitting of $q \circ p: K \rightarrow \Nerve(\FinSeg)$. Conversely, suppose that
$q \circ p$ is split, and choose a complement $\overline{\beta}: q \circ p \rightarrow \overline{p}_1$
to $\overline{\alpha}$. Then $\overline{\beta}$ is inert, so we can choose a $q$-coCartesian lift
$\beta: p \rightarrow p_1$ of $\overline{\beta}$ which is a complement to $\alpha$. The uniqueness
of $\beta$ follows from the observation that $\overline{\beta}$ and its $q$-coCartesian lift are both well-defined up to equivalence.
\end{proof}

\begin{lemma}\label{squaat}
Let $q: \calC^{\otimes} \rightarrow \calO^{\otimes}$ be a fibration of $\infty$-operads.
Let $\alpha: X \rightarrow X_0$ and $\beta: X \rightarrow X_1$ be morphisms
in $\calO^{\otimes}$ which determine a splitting of $X$, and suppose that
$\overline{X}_0$ and $\overline{X}_1$ are objects of $\calC^{\otimes}$ lying
over $X_0$ and $X_1$, respectively. Then:
\begin{itemize}
\item[$(1)$] Let $\overline{\alpha}: \overline{X} \rightarrow \overline{X}_0$
and $\overline{\beta}: \overline{X} \rightarrow \overline{X}_1$ be morphisms in
$\calC^{\otimes}$ lying over $\alpha$ and $\beta$. Then $\overline{\alpha}$
and $\overline{\beta}$ determine a splitting of $\overline{X}$ if and only if they
exhibit $\overline{X}$ as a $q$-product of $\overline{X}_0$ and $\overline{X}_1$.

\item[$(2)$] There exist morphisms $\overline{\alpha}: \overline{X} \rightarrow \overline{X}_0$
and $\overline{\beta}: \overline{X} \rightarrow \overline{X}_1$ satisfying the equivalent conditions
of $(1)$.
\end{itemize}
\end{lemma}

\begin{proof}
We will prove $(2)$ and the ``if'' direction of $(1)$; the ``only if'' direction follows from $(2)$ together with the uniqueness properties of $q$-limit diagrams. We begin with $(1)$. Choose a diagram
$\sigma: \Delta^1 \times \Delta1 \rightarrow \calC^{\otimes}$
$$ \xymatrix{ \overline{X} \ar[r]^{\overline{\alpha}} \ar[d]^{\overline{\beta}} & \overline{X}_0 \ar[d] \\
\overline{X}_1 \ar[r] & 0, }$$
where $0$ is a final object of $\calC^{\otimes}$ (in other words, $0$ lies in
$\calC^{\otimes}_{\seg{0}}$). 
Let $K \simeq \Lambda^2_0$ denote the full subcategory of
$\Delta^1 \times \Delta^1$ obtained by removing the final object. 
Since $0$ is final in $\calC^{\otimes}$ and
$q(0)$ is final in $\calO^{\otimes}$, we deduce that $0$ is a a $q$-final object of
$\calC^{\otimes}$ (Proposition \toposref{basrel}), so that $\sigma$ is a $q$-right Kan extension of $\sigma|K$. Using Propositions \ref{prejule} and \bicatref{fibraman}, we deduce that
$\calC^{\otimes}$ is a $\calO$-operad family, so the diagram $\sigma$ is a $q$-limit.
Applying Lemma \toposref{kan0}, we deduce that $\sigma|K$ is a $q$-limit, so that
$\sigma$ exhibits $\overline{X}$ as a $q$-product of $\overline{X}_0$ and $\overline{X}_1$.

We now prove $(2)$.
Let $0$ be an object of $\calC^{\otimes}_{\seg{0}}$. Since $0$ is a final object of
$\calC^{\otimes}$, we can find morphisms (automatically inert) $\gamma: \overline{X}_0 \rightarrow
0$ and $\delta: \overline{X}_1 \rightarrow 0$ in $\calC^{\otimes}$. Since $q(0)$ is a
final object of $\calO^{\otimes}$, we can find a commutative square
$$ \xymatrix{ X \ar[r]^{\alpha} \ar[d]^{\beta} & X_0 \ar[d]^{q(\gamma)} \\
X_1 \ar[r]^{ q(\delta)} & q(0) }$$
in $\calO^{\otimes}$. Since $\calC^{\otimes}$ is a $\calO$-operad family
we can lift this to a $q$-limit diagram $\sigma$
$$ \xymatrix{ \overline{X} \ar[r]^{\overline{\alpha}} \ar[d]^{\overline{\beta}} & \overline{X}_0 \ar[d]^{\gamma} \\
\overline{X}_1 \ar[r]^{\delta} & 0 }$$
in $\calC^{\otimes}$, where $\overline{\alpha}$ and $\overline{\beta}$ are inert and therefore determine a splitting of $\overline{X}$.
\end{proof}

\begin{corollary}\label{spade}
Let $q: \calC^{\otimes} \rightarrow \calO^{\otimes}$ be a fibration of $\infty$-operads, let
$\calX$ be the full subcategory of $\Fun( \Lambda^2_0, \calC^{\otimes})$ be the full subcategory spanned by those diagrams $ X_0 \leftarrow X \rightarrow X_1$
which determine a splitting of $X$, and let $\calY \subseteq \Fun( \Lambda^2_0, \calO^{\otimes})$ be defined similarly. Then the canonical map
$$ \calX \rightarrow \calY \times_{ \Fun( \{1\}, \calO^{\otimes}) \times \Fun( \{2\}, \calO^{\otimes})
} ( \Fun( \{1\}, \calC^{\otimes}) \times \Fun( \{2\}, \calC^{\otimes}) )$$
is a trivial Kan fibration.
\end{corollary}

\begin{proof}
Combine Lemma \ref{squaat} with Proposition \toposref{lklk}.
\end{proof}

\begin{lemma}\label{helpless}
Let $q: \calC^{\otimes} \rightarrow \calO^{\otimes}$ be a fibration of $\infty$-operads.
Suppose we are given a split natural transformation
$\alpha: p \rightarrow p_0$ of diagrams $p, p_0: K^{\triangleleft} \rightarrow \calO^{\otimes}$.
Let $\overline{p}_0: K^{\triangleleft} \rightarrow \calC^{\otimes}$ be a diagram lifting
$p_0$, let $\overline{p}': K \rightarrow \calC^{\otimes}$ be a diagram lifting
$p' = p|K$, and let $\overline{\alpha}': \overline{p}' \rightarrow \overline{p}_0|K$ be a natural
transformation lifting $\alpha' = \alpha| ( \Delta^1 \times K)$. Suppose that the following
condition is satisfied:
\begin{itemize}
\item[$(\ast)$] Let $\beta: p \rightarrow p_1$ be a complement to $\alpha$, let 
$\beta' = \beta| ( \Delta^1 \times K)$, and let $\overline{\beta}': \overline{p}' \rightarrow \overline{p}'_1$ be
a $q$-coCartesian natural transformation lifting $\beta'$ (so that $\overline{\beta}'$ is a complement
to $\overline{\alpha}'$). Then $\overline{p}'_1$ can be extended to a $q$-limit diagram
$\overline{p}_1: K^{\triangleleft} \rightarrow \calC^{\otimes}$ such that $q \circ \overline{p}_1 = p_1$.
\end{itemize}
Then: 
\begin{itemize}
\item[$(1)$] Let $\overline{\alpha}: \overline{p} \rightarrow \overline{p}_0$ be a 
natural transformation of diagrams $\overline{p}, \overline{p}_0: K^{\triangleleft} \rightarrow \calC^{\otimes}$ which extends $\overline{\alpha}'$ and lies over $\alpha$. The following conditions are equivalent:
\begin{itemize}
\item[$(i)$] The map of simplicial sets
$\overline{\alpha}: \Delta^1 \times K^{\triangleleft} \rightarrow \calC^{\otimes}$
is a $q$-limit diagram.
\item[$(ii)$] The natural transformation $\overline{\alpha}$ is inert (and therefore split), and if
$\overline{\beta}: \overline{p} \rightarrow \overline{p}_1$ is a complement to 
$\overline{\alpha}$, then $\overline{p}_1$ is a $q$-limit diagram.
\end{itemize}
\item[$(2)$] There exists a natural transformation $\overline{\alpha}: \overline{p} \rightarrow \overline{p}_0$ satisfying the equivalent conditions of $(1)$.
\end{itemize}
\end{lemma}

\begin{proof}
We first prove the implication $(ii) \Rightarrow (i)$ of assertion $(1)$. Choose
a complement $\overline{\beta}: \overline{p} \rightarrow \overline{p}_{1}$ to $\overline{\alpha}$,
so that $\overline{\alpha}$ and $\overline{\beta}$ together determine a map
$\overline{F}: \Lambda^2_0 \times K^{\triangleleft} \rightarrow \calC^{\otimes}$
with $\overline{F} | \Delta^{ \{0,1\} } \times K^{\triangleleft} = \overline{\alpha}$ and
$\overline{F} | \Delta^{ \{0,2\} } \times K^{\triangleleft} = \overline{\beta}$.
Using the small object argument, we can choose an inner anodyne map
$K \rightarrow K'$ which is bijective on vertices, where $K'$ is an $\infty$-category.
Since $\calC^{\otimes}$ is an $\infty$-category, the map $\overline{F}$ factors as a composition
$K^{\triangleleft} \times \Lambda^2_0 \rightarrow {K'}^{\triangleleft} \times \Lambda^2_0 \rightarrow \calC^{\otimes}$. We may therefore replace $K$ by $K'$ and thereby reduce to the case where
$K$ is an $\infty$-category.

The inclusion $i: K \times \{0\} \subseteq K \times \Delta^{ \{0,2\} }$ is left anodyne, so that
$i^{op}$ is cofinal. It will therefore suffice to show that the restriction $\overline{F}_0$ of
$\overline{F}$ to 
$(K^{\triangleleft} \times \Delta^{ \{0,1\} }) \coprod_{ K \times \{0\} } (K \times \Delta^{ \{0,2\} })$
is a $q$-limit diagram. Since $\overline{p}_1$ is a $q$-limit diagram, $\overline{F}$ is a $q$-right Kan extension of $\overline{F}_0$; according to Lemma \toposref{kan0} it will suffice to prove that $\overline{F}$ is a $q$-limit diagram. 

Let $v$ denote the cone point of $K^{\triangleleft}$. 
Let $\calD$ be the full subcategory of $K^{\triangleleft} \times \Lambda^2_0$ spanned
by $K^{\triangleleft} \times \{1\}$, $K^{\triangleleft} \times \{2\}$, and $(v,0)$. Using Lemma \ref{squaat}, we deduce that $\overline{F}$ is a $q$-right Kan extension of $F| \calD$. Using Lemma \toposref{kan0} again, we are reduced to proving that $\overline{F}| \calD$ is a $q$-limit diagram. Since the inclusion
$$\{ (v,1) \} \coprod \{ (v,2) \} \subseteq (K^{\triangleleft} \times \{1\})
\coprod (K^{\triangleleft} \times \{2\})$$
is the opposite of a cofinal map, it suffices to show that
$F | \{v\} \times \Lambda^2_0$ is a $q$-limit diagram, which follows from Lemma \ref{squaat}. 
This completes the verification of condition $(i)$.

We now prove $(2)$. Choose a complement $\beta: p \rightarrow p_1$ to the split natural
transformation $\alpha$, let $\beta' = \beta | \Delta^1 \times K$, and choose an $q$-coCartesian
natural transformation $\overline{\beta}': \overline{p}' \rightarrow \overline{p}'_1$ lifting
$\beta'$. Invoking assumption $(\ast)$, we can extend $\overline{p}'_{1}$ to a 
$q$-limit diagram $\overline{p}_1: K^{\triangleleft} \rightarrow \calC^{\otimes}$ such that
$q \circ \overline{p}_1 = p_1$. The maps $\overline{p}_0$, $\overline{p}_1$, 
$\overline{\alpha}'$ and $\overline{\beta}'$ can be amalgamated to give a map
$$F: (\Lambda^2_0 \times K) \coprod_{ ( \{1, 2\}) \times K } ( \{1,2\} \times K^{\triangleleft})
\rightarrow \calC^{\otimes}.$$
Using Corollary \ref{spade}, we can extend $F$ to a map $\overline{F}: \Lambda^2_0 \times K^{\triangleleft} \rightarrow \calC^{\otimes}$ corresponding to a pair of morphisms
$\overline{\alpha}: \overline{p} \rightarrow \overline{p}_0$ and
$\overline{\beta}: \overline{p} \rightarrow \overline{p}_1$ having the desired properties.

The implication $(i) \Rightarrow (ii)$ of $(1)$ now follows from $(2)$, together with the uniqueness properties of $q$-limit diagrams.
\end{proof}

\begin{proof}[Proof of Proposition \ref{canda}]
We first treat the case where $K$ is an $\infty$-category.
Let $Y = \calK_{\calO}$ and $Y' = \calK^{0}_{\calO} \subseteq Y$. Let
$\pi: Y \rightarrow \calO^{\otimes}$ be the map given by evaluation at $\{0\}$, and let
$\pi' = \pi|Y'$. Our assumption that $\calO^{\otimes}$ is coherent guarantees that $\pi$ is a flat categorical fibration. Let $X = Y \times_{ \Fun( \{1\}, \calO^{\otimes})} \calC^{\otimes}$, and let
$X' = X \times_{Y} Y'$. The map $\psi: \MMod^{\calO}(\calC)^{\otimes} \rightarrow \PrAlg_{\calO}(\calC)$ can be identified with a restriction of the map $\pi_{\ast} X \rightarrow \pi'_{\ast} X'$.
We are given a diagram
$$ \xymatrix{ K \ar[r]^{p} \ar[d] & \pi_{\ast} X \ar[d] \\
K^{\triangleleft} \ar[r]^{ \overline{p}_0} \ar@{-->}[ur]^{\overline{p}} & \pi'_{\ast} X'. }$$
We claim that this situation satisfies the hypotheses of Proposition \ref{skitter2}:

\begin{itemize}
\item[$(i)$] The full subcategory $Y' \times_{ \calO^{\otimes}} K^{\triangleleft}$ is a cosieve on $Y
\times_{ \calO^{\otimes} } K^{\triangleleft}$. Since the map $K^{\triangleleft} \rightarrow
\calO^{\otimes}$ is constant taking some value $C \in \calO$, it will suffice to show that
the $Y' \times_{ \calO^{\otimes} } \{C\}$ is a cosieve on $Y \times_{ \calO^{\otimes}} \{C\}$. 
Unwinding the definitions, this amounts to the following assertion: given a commutative diagram
$$ \xymatrix{ C \ar[r] \ar[d]^{\id} & D \ar[d] \\
C \ar[r] & D' }$$
in $\calO^{\otimes}$, if the upper horizontal map is null then the lower horizontal map is null.
This is clear, since the collection of null morphisms in $\calO^{\otimes}$ is closed under composition with other morphisms. 

\item[$(ii)$] For every object $y \in Y'$ and every morphism
$f: z \rightarrow \pi(y)$ in $\calO$, there exists a $\pi$-Cartesian morphism
$\overline{f}: \overline{z} \rightarrow y$ in $Y'$ such that $\pi( \overline{f} ) = f$.
We can identify $y$ with a semi-inert morphism $y_0 \rightarrow y_1$ in $\calO^{\otimes}$, and $f$ with a morphism $z \rightarrow y_0$ in $\calO^{\otimes}$.
Using Corollary \toposref{tweezegork}, we see that the morphism $\overline{f}$ can be taken to correspond to the commutative diagram
$$ \xymatrix{ z \ar[r] \ar[d] & y_1 \ar[d]^{\id} \\
y_0 \ar[r] & y_1 }$$
in $\calO^{\otimes}$: our assumption that $y_0 \rightarrow y_1$ is null guarantees
that the composite map $z \rightarrow y$ is also null.

\item[$(iii)$] Let $\pi''$ denote the projection map
$K^{\triangleleft} \times_{\calO^{\otimes} } Y \rightarrow K^{\triangleleft}$. Then
$\pi''$ is a coCartesian fibration. This is clear, since $\pi''$ is a pullback of the
coCartesian fibration $( \calK_{\calO} \times_{ \calO^{\otimes} } \{C\}) \rightarrow \Delta^0$.

\item[$(iv)$] Let $v$ denote the cone point of $K^{\triangleleft}$ and let
$\calD = {\pi''}^{-1} \{v\}$. Then $\calD' = \calD \times_{Y} Y'$ is a localization of
$\calD$. We can identify an object of $\calD$ with a semi-inert morphism
$f: C \rightarrow C'$ in $\calO^{\otimes}$. We wish to prove that for any such object
$f$, there exists a morphism $f \rightarrow g$ in $\calD$ which exhibits
$g$ as a $\calD'$-localization of $f$. Let $f_0: \seg{1} \rightarrow \seg{k}$ denote
the underlying morphism in $\FinSeg$. If $f_0$ is null, then $f \in \calD'$ and there
is nothing to prove. Otherwise, $f_0(1) = i$ for some $1 \leq i \leq k$. Choose
an inert map $h_0: \seg{k} \rightarrow \seg{k-1}$ such that $h_0(i) = \ast$, and choose
an inert morphism $h: C' \rightarrow D$ in $\calO^{\otimes}$ lifting $h_0$.
We then have a commutative diagram
$$ \xymatrix{ C \ar[r]^{f} \ar[d]^{\id} & C' \ar[d]^{h} \\
C \ar[r]^{g} & D } $$
in $\calO^{\otimes}$, corresponding to a map $\alpha: f \rightarrow g$ in $\calD$; by construction,
$g$ is null so that $g \in \calD'$. We claim that $\alpha$ exhibits $g$ as a $\calD'$-localization
of $\calC$. To prove this, choose any object $g': C \rightarrow D'$ in $\calD'$; we wish to show that
composition with $\alpha$ induces a homotopy equivalence 
$$\bHom_{ \calD}( g, g') \simeq \bHom_{ (\calO^{\otimes})^{C/}}( g, g')
\rightarrow \bHom_{ (\calO^{\otimes})^{C/}}( f, g') \simeq \bHom_{ \calD}(f,g).$$
Since the projection map $$(\calO^{\otimes})^{C/} \rightarrow \calO^{\otimes}$$ is a left fibration, it will suffice to show that the map $\bHom_{ \calO^{\otimes}}( D, D') \rightarrow
\bHom^{0}_{\calO^{\otimes}}( C', D')$, where the superscript indicates that we consider
only morphisms $C' \rightarrow D'$ such that the underlying map $\seg{k} \rightarrow \seg{k'}$
carries $i$ to the base point $\ast \in \seg{k'}$. Since $h$ is inert, this follows from the observation that
composition with $h_0$ induces an injection
$\Hom_{ \FinSeg}( \seg{k-1}, \seg{k'}) \rightarrow \Hom_{ \FinSeg}( \seg{k}, \seg{k'})$ whose
image consists of those maps which carry $i$ to the base point.
\end{itemize}

Fix an object of $\calD$ corresponding to a semi-inert morphism $f:C \rightarrow C'$ in $\calO^{\otimes}$, and let $\alpha: f \rightarrow g$ be a map in $\calD$ which exhibits $g$ as a $\calD'$-localization of
$f$ (as in the proof of $(iv)$). Using the maps $p$ and $\overline{p}_0$, we get commutative
diagram
$$ \xymatrix{ (K \times \Delta^1) \coprod_{ K \times \{1\} } (K^{\triangleleft} \times \{1\} ) \ar[r] \ar[d] & \calC^{\otimes} \ar[d]^{q} \\
K^{\triangleleft} \times \Delta^1 \ar[r] \ar@{-->}[ur]^{\theta} & \calO^{\otimes}. }$$
To apply Proposition \ref{skitter2}, we must know that every such diagram admits
an extension as indicated, where $\theta$ is a $q$-limit. This follows from
Lemma \ref{helpless} and assumption $(\ast)$. Moreover, we obtain the following
criterion for testing whether $\theta$ is a $q$-limit diagram:
\begin{itemize}
\item[$(\ast')$] Let $\theta: K^{\triangleleft} \times \Delta^1 \rightarrow \calC^{\otimes}$ be
as above, and view $\theta$ as a natural transformation $d \rightarrow d_0$ of diagrams
$d, d_0: K^{\triangleleft} \rightarrow \calC^{\otimes}$. Then $\theta$ is a $q$-limit diagram
if and only if it is an inert (and therefore split) natural transformation, and admits a complement
$d \rightarrow d_1$ where $d_1: K^{\triangleleft} \rightarrow \calC^{\otimes}$ is a $q$-limit diagram.
\end{itemize} 

Applying Proposition \ref{skitter2}, we obtain the following:

\begin{itemize}
\item[$(a)$] There exists a solution to the lifting problem
$$ \xymatrix{ K \ar[r]^-{p} \ar[d] & \tildeMod^{\calO}(\calC)^{\otimes} \ar[d]^{\overline{\psi}} \\
K^{\triangleleft} \ar[r]^-{ \overline{p}_0} \ar@{-->}[ur]^{\overline{p}} & \tildePAlg_{\calO}(\calC), }$$ where $\overline{p}$ is an $\overline{\psi}$-limit diagram.

\item[$(b)$] An arbitrary extension $\overline{p}$ as above is an $\overline{\psi}$-limit diagram
if and only if the following condition is satisfied:
\begin{itemize}
\item[$(\ast'')$] For every object $f: C \rightarrow C'$ in $\calD$ and every
morphism $\alpha: f \rightarrow g$ in $\calD$ which exhibits $g$ as a $\calD'$-localization
of $f$, if $\theta: K^{\triangleleft} \times \Delta^1 \rightarrow \calC^{\otimes}$ is defined
as above, then $\theta$ is a split natural transformation of diagrams
$d, d_0: K^{\triangleleft} \rightarrow \calC^{\otimes}$ and admits a complement
$d \rightarrow d_1$ where $d_1: K^{\triangleleft} \rightarrow \calC^{\otimes}$ is a $q$-limit diagram.
\end{itemize}
\end{itemize}

To complete the proof, we must show that condition $(\ast'')$ is equivalent to the following
pair of assertions:
\begin{itemize}
\item[$(I)$] The map $\overline{p}$ carries $K^{\triangleleft}$ into the full subcategory
$\MMod^{\calO}(\calC)^{\otimes} \subseteq \tildeMod^{\calO}(\calC)^{\otimes}$.
(Since we know already that $p$ has this property, it suffices to check that
$\overline{p}$ carries the cone point $v$ of $K^{\triangleleft}$ into 
$\tildeMod^{\calO}(\calC)^{\otimes}$). 
\item[$(II)$] The composite map
$$ \overline{p}: K^{\triangleleft} \rightarrow  \tildeMod^{\calO}(\calC)^{\otimes} 
\rightarrow \calC^{\otimes}$$
is a $q$-limit diagram. Here the second map is induced by composition with the diagonal embedding $\calO^{\otimes} \hookrightarrow \calK_{\calO}$.
\end{itemize}

Assume first that condition $(\ast'')$ is satisfied by $\overline{p}: K^{\triangleleft}
\rightarrow \tildeMod^{\calO}(\calC)^{\otimes}$; we will prove that
$\overline{p}$ also satisfies $(I)$ and $(II)$. We can identify $\overline{p}$ with a map
$\overline{P}: K^{\triangleleft} \times \calD \rightarrow \calC^{\otimes}$. 
We first prove $(II)$. Fix an object $f \in \calD$ corresponding to an equivalence $C \rightarrow C'$ in
$\calO^{\otimes}$; we will show that $\overline{P}_{f} = \overline{P}| K^{\triangleleft} \times \{f\}$ is a $q$-limit diagram. Choose a morphism $f \rightarrow g$ which exhibits
$g$ as a $\calD'$-localization of $f$, and let $\theta: \overline{P}_{f} \rightarrow \overline{P}_{g}$
be the induced natural transformation as in $(\ast'')$. Let $\theta': \overline{P}_{f} \rightarrow d_1$
be a complement to $\theta$. Since $C \simeq C'$, 
$\overline{P}_{g}$ takes values in $\calC^{\otimes}_{\seg{0}}$, so $\theta'$ is an equivalence of diagrams. Condition $(\ast'')$ implies that $d_1$ is a $q$-limit diagram, so that
$\overline{P}_{f}$ is a $q$-limit diagram.

To prove $(I)$, we must show that
for every morphism $\alpha: f \rightarrow f'$ in $\calD$ whose
image in $\calK_{\calO}$ is inert, the induced map $\overline{P}(v,f) \rightarrow
\overline{P}(v,f')$ is inert in $\calC^{\otimes}$. There are several cases to consider:
\begin{itemize}
\item[$(I1)$] The map $f$ belongs to $\calD'$. Then
$f' \in \calD'$ and the desired result follows from our assumption that $\overline{p}_0$ factors through
$\PrAlg_{\calO}(\calC)$. 
\item[$(I2)$] The map $f$ does not belong to $\calD'$, but $f'$ does. Then $\alpha$
factors as a composition
$$ f \stackrel{\alpha'}{\rightarrow} g \stackrel{\alpha''}{\rightarrow} f',$$
where $\alpha'$ exhibits $g$ as a $\calD'$-localization of $f$. Since the composition of inert morphisms in $\calC^{\otimes}$ is inert and $g \in \calD'$, we can apply $(I1)$ to reduce to the case where
$\alpha = \alpha'$. In this case, the desired result follows immediately from $(\ast'')$.
\item[$(I3)$] The map $f'$ is an equivalence in $\calO^{\otimes}$. Let $\overline{\alpha}:
\overline{P}_{f} \rightarrow \overline{P}_{f'}$ be the natural transformation induced by $\alpha$; it will suffice
to show that this natural transformation in inert. Let
$\beta: f \rightarrow g$ be a map in $\calD$ which exhibits $g$ as a
$\calD'$-localization of $f$. Then $\beta$ induces a natural transformation
$\theta: \overline{P}_{f} \rightarrow \overline{P}_{g}$. Using $(\ast'')$, we can choose a complement
$\theta': \overline{P}_{f} \rightarrow d_1$ to $\theta$. Since $\theta'$ is a $q$-coCartesian transformation of diagrams, we obtain a factorization of $\overline{\alpha}$ as a composition
$$ \overline{P}_{f} \stackrel{\theta'}{\rightarrow} d_1 \stackrel{\gamma}{\rightarrow} \overline{P}_{f'}.$$
We wish to prove that $\gamma$ is an equivalence. Since $d_1$ is a $q$-limit diagram
(by virtue of $(\ast'')$) and $\overline{P}_{f'}$ is a $q$-limit diagram (by virtue of $(I)$),
it will suffice to show that $\gamma$ induces an equivalence $d_1 | K \rightarrow \overline{P}_{f'}|K$.
This follows from the fact that $p$ factors through $\MMod^{\calO}(\calC)^{\otimes}$.
\item[$(I4)$] The map $f'$ does not belong to $\calD'$. Let us identify
$f'$ with a semi-inert morphism $C \rightarrow C'$ in $\calO^{\otimes}$, lying over
an injective map $j: \seg{1} \hookrightarrow \seg{k}$ in $\FinSeg$. Choose a splitting
$C'_0 \leftarrow C' \rightarrow C'_1$ of $C'$ corresponding to the decomposition
$\nostar{k} \simeq \nostar{1} \coprod \nostar{k-1}$ induced by $j$. This splitting can be
lifted to a pair of morphisms $f' \rightarrow f'_0$ and $f' \rightarrow f'_1$ in
$\calD$. Using $(I2)$ and $(I3)$, we deduce that the maps $\overline{P}(v,f') \rightarrow
\overline{P}(v,f'_0)$ and $\overline{P}(v,f') \rightarrow \overline{P}(v,f'_{1})$ are inert.
Since $\calC^{\otimes}$ is an $\infty$-operad, to prove that the map
$\overline{P}(v,f) \rightarrow \overline{P}(v,f')$ is inert, it will suffices to show that the composite
maps $\overline{P}(v,f) \rightarrow \overline{P}(v,f'_0)$ and $\overline{P}(v,f) \rightarrow
\overline{P}(v, f'_{1})$ are inert. In other words, we may replace $f'$ by $f'_0$ or
$f'_1$ and thereby reduce to the cases $(I2)$ and $(I3)$.
\end{itemize}

Now suppose that conditions $(I)$ and $(II)$ are satisfied; we will prove
$(\ast'')$. Fix an object $f$ in $\calD$, let $\alpha: f \rightarrow g$ be a
map which exhibits $g$ as a $\calD'$-localization of $f$, let 
$\theta: \overline{P}_{f} \rightarrow \overline{P}_{g}$ be the induced natural transformation.
Our construction of $\alpha$ together with assumption $(I)$ guarantees that
$\theta$ is split; let $\theta': \overline{P}_{f} \rightarrow d_1$ be a complement to
$\theta$. We wish to prove that $d_1$ is a $q$-limit diagram. If
$f \in \calD'$, then $d_1$ takes values in $\calC^{\otimes}_{ \seg{0} }$
and the result is obvious. We may therefore assume that
$f: C \rightarrow C'$ induces an injective map $\seg{1} \rightarrow \seg{k}$ in
$\FinSeg$; choose a splitting $C'_0 \leftarrow C' \rightarrow C'_1$ corresponding
the decomposition $\nostar{k} \simeq \nostar{1} \coprod \nostar{k-1}$. 
This splitting lifts to a pair of maps $\beta_0: f \rightarrow f_0$,
$\beta_1: f \rightarrow f_1$ in $\calD$, and we can identify $\beta_1$ with
$\alpha: f \rightarrow g$. Using assumption $(I)$, we see that
$\beta_0$ induces a transformation $\overline{P}_{f} \rightarrow \overline{P}_{f_0}$
which is a complement to $\theta$. We are therefore reduced to showing that
$\overline{P}_{f_0}$ is a $q$-limit diagram. This follows from
$(II)$, since $f_0: C \rightarrow C'_0$ is an equivalence in $\calC$ and
therefore equivalent (in $\calD$) to the identity map $\id_{C}$.
\end{proof}

\subsection{Colimits of Modules}\label{com35}

Let $\calC$ be a symmetric monoidal category, let $A$ be a commutative algebra object
of $\calC$, let $\{ M_{\alpha} \}$ be a diagram in the category of $A$-modules, and let
$M = \varinjlim M_{\alpha}$ be a colimit of this diagram in the underling category $\calC$.
For each index $\alpha$, we have a canonical map
$$ A \otimes M_{\alpha} \rightarrow M_{\alpha} \rightarrow M.$$
These maps together determine a morphism $\varinjlim (A \otimes M_{\alpha}) \rightarrow M$.
If tensor product with $A$ preserves colimits, then we can identify the domain of this map
with $A \otimes M$, and the object $M \in \calC$ inherits the structure of an $A$-module
(which is then a colimit for the diagram $\{ M_{\alpha} \}$ in the category of $A$-modules).

Our goal in this section is to obtain an $\infty$-categorical generalization of the above discussion.
We first formalize the idea that ``tensor products commute with colimits''.

\begin{definition}
Let $\calO^{\otimes}$ be an $\infty$-operad. We will say that a fibration of $\infty$-operads
$q: \calC^{\otimes} \rightarrow \calO^{\otimes}$ is a {\it presentable $\calO$-monoidal $\infty$-category} if the following conditions are satisfied:
\begin{itemize}
\item[$(1)$] The functor $q$ is a coCartesian fibration of $\infty$-operads.
\item[$(2)$] The $\calO$-monoidal structure on $\calC$ is compatible with small colimits.
\item[$(3)$] For each $X \in \calO$, the fiber $\calC^{\otimes}_{X}$ is a presentable $\infty$-category.
\end{itemize}
\end{definition}

\begin{theorem}\label{turkwell}
Let $\calO^{\otimes}$ be a small coherent $\infty$-operad, and let $q: \calC^{\otimes} \rightarrow \calO^{\otimes}$ be a presentable $\calO$-monoidal $\infty$-category. Let
$A \in \Alg_{\calO}(\calC)$ be a $\calO$-algebra object of $\calC$. Then the induced map
$\psi: \Mod^{\calO}_{A}(\calC)^{\otimes} \rightarrow \calO^{\otimes}$ exhibits
$\Mod^{\calO}_{A}(\calC)^{\otimes}$ as a presentable $\calO$-monoidal $\infty$-category.
\end{theorem}

We will deduce Theorem \ref{turkwell} from the a more general result,
which can be used to construct colimits in $\infty$-categories of module objects in a wide variety of situations. The statement is somewhat complicated, since the idea that ``tensor product with
$A$ preserves colimits'' needs to be formulated using the theory of operadic colimit diagrams
described in \S \ref{comm33}.

\begin{theorem}\label{colmodd}
Let $q: \calC^{\otimes} \rightarrow \calO^{\otimes}$ be a fibration of $\infty$-operads, where $\calO^{\otimes}$ is coherent. Let $K$ be an $\infty$-category and let $A \in \Alg_{\calO}(\calC)$ be a 
$\calO$-algebra object of $\calC$. 
Suppose we are given a commutative diagram
$$ \xymatrix{ K \ar[r]^-{p} \ar[d] & \Mod^{\calO}_{A}(\calC)^{\otimes} \ar[d]^{\psi} \\
K^{\triangleright} \ar[r] \ar@{-->}[ur]^{\overline{p}} & \calO^{\otimes}. }$$
Let $\overline{\calD} = K^{\triangleright} \times_{ \calO^{\otimes}} \calK_{\calO}$
and let $\calD = K \times_{ \calO^{\otimes}} \calK_{\calO} \subseteq \overline{\calD}$, so that
$p$ classifies a diagram $F: \calD \rightarrow \calC^{\otimes}$. Assume the following:
\begin{itemize}
\item[$(i)$] The induced map $K^{\triangleright} \rightarrow \calO^{\otimes}$ factors
through $\calO^{\otimes}_{\acti}$, and carries the cone point of $K^{\triangleright}$ to an
object $X \in \calO$.
\item[$(ii)$] Let $D = (v, \id_{X}) \in \overline{\calD}$. Let $\calD_{/D}^{\acti}$ denote the full subcategory of
$\calD \times_{ \overline{\calD} } \overline{\calD}_{/D}$ spanned by those morphisms
$D' \rightarrow D$ in $\overline{\calD}$ which induce diagrams
$$ \xymatrix{ X' \ar[r] \ar[d] & Y' \ar[d]^{f} \\
X \ar[r]^{\id} & X }$$
in $\calO^{\otimes}$, where $f$ is active. Then the diagram
$$ \calD_{/D}^{\acti} \rightarrow \calD \stackrel{F}{\rightarrow} \calC^{\otimes}$$
can be extended to a $q$-operadic colimit diagram
$(\calD_{/D}^{\acti})^{\triangleright} \rightarrow \calC^{\otimes}$ lying over the composite map
$$ (\calD_{/D}^{\acti})^{\triangleright} \rightarrow \overline{\calD}^{\triangleright}_{/D}
\rightarrow \overline{\calD} \rightarrow K^{\triangleright} \rightarrow \calO^{\otimes}.$$
\end{itemize}

Then:
\begin{itemize}
\item[$(1)$] Let $\overline{p}$ be an extension of $p$ as indicated in the above diagram,
corresponding to a map $\overline{F}: \overline{\calD} \rightarrow \calC^{\otimes}$.
Then $\overline{p}$ is an operadic $\psi$-colimit diagram if and only if the following
condition is satisfied:
\begin{itemize}
\item[$(\ast)$] For every object $(v, \id_{X})$ as in $(ii)$, the map
$$ ( \calD_{/D}^{\acti} )^{\triangleright} \rightarrow 
\overline{\calD}_{/D}^{\triangleright} \rightarrow \overline{\calD} \stackrel{ \overline{F}}{\rightarrow}
\calC^{\otimes}$$
is an operadic $q$-colimit diagram.
\end{itemize}
\item[$(2)$] There exists an extension $\overline{p}$ of $p$ satisfying condition $(\ast)$.
\end{itemize}
\end{theorem}

The proof of Theorem \ref{colmodd} is rather technical, and will be given at the end of this section.

\begin{corollary}\label{swink}
Let $p: \calC^{\otimes} \rightarrow \calO^{\otimes}$ be a fibration of
$\infty$-operads, where $\calO^{\otimes}$ is coherent. Let
$A \in \Alg_{\calO}(\calC)$. Let $f: M_0 \rightarrow M$
be a morphism in $\Mod^{\calO}_{A}( \calC)^{\otimes}$ be a morphism
where $M_0 \in \Mod^{\calO}_{A}(\calC)^{\otimes}_{\seg{0}}$ and
$M \in \Mod^{\calO}_{A}(\calC)$. The following conditions are equivalent:
\begin{itemize}
\item[$(1)$] The morphism $f$ is classified by an operadic
$q$-colimit diagram $$\Delta^1 \rightarrow \Mod^{\calO}_{A}(\calC)^{\otimes},$$ where
$q: \Mod^{\calO}_{A}(\calC)^{\otimes} \rightarrow \calO^{\otimes}$ denotes the projection.
\item[$(2)$] Let $F: \calK_{\calO} \times_{ \calO^{\otimes}} \Delta^1 \rightarrow \calC^{\otimes}$
be the map corresponding to $f$. Then $F$ induces an equivalence
$F( q(f) ) \rightarrow F( \id_{q(M)} )$. 
\end{itemize}
Moreover, for every $X \in \calO^{\otimes}$, there exists a morphism
$f: M_0 \rightarrow M$ satisfying the above conditions, with $q(M) = X$.
\end{corollary}

\begin{proof}
Apply Theorem \ref{colmodd} together with the observation that
the $\infty$-category $\calD_{/D}^{\acti}$ has a final object.
\end{proof}

\begin{corollary}\label{sabus}
Let $\kappa$ be an uncountable regular cardinal. Let $\calO^{\otimes}$ be a $\kappa$-small coherent $\infty$-operad, and let $q: \calC^{\otimes} \rightarrow \calO^{\otimes}$ be a $\calO$-monoidal
$\infty$-category which is compatible with $\kappa$-small colimits.
Let $A \in \Alg_{\calO}(\calC)$ be a $\calO$-algebra object of $\calC$. Then:
\begin{itemize}
\item[$(1)$] The map $\psi: \Mod^{\calO}_{A}(\calC)^{\otimes} \rightarrow \calO^{\otimes}$ is a $\calO$-monoidal $\infty$-category, which is compatible with $\kappa$-small colimits. 

\item[$(2)$] For each object $X \in \calO$, consider the induced functor $\phi: \Mod^{\calO}_{A}(\calC)^{\otimes}_{X} \rightarrow \calC^{\otimes}_{X}$. Let $K$ be a $\kappa$-small simplicial set and
let $\overline{p}: K^{\triangleright} \rightarrow \Mod^{\calO}_{A}(\calC)^{\otimes}_{X}$ be a map. Then $\overline{p}$ is a colimit diagram if and only if $\phi \circ \overline{p}$ is a colimit diagram.
\end{itemize}
\end{corollary}

\begin{proof}
Assertion $(1)$ follows immediately from Theorem \ref{colmodd} and Corollary \ref{cyster}.
We will prove $(2)$. Without loss of generality, we may assume that $K$ is an $\infty$-category.
Let $\calD = \calK_{\calO} \times_{ \calO^{\otimes} } \{X\}$ denote the full subcategory of $( \calO^{\otimes})^{X/}$ spanned by the semi-inert morphisms
$X \rightarrow Y$ in $\calO^{\otimes}$, so that we can identify $\overline{p}$ with a functor
$F: \calD \times K^{\triangleright} \rightarrow \calC^{\otimes}$.
Let $p = \overline{p} | K$. It follows from $(1)$ (and Corollary \ref{cyster}) that
$p$ can be extended to an operadic $\psi$-colimit diagram in $\Mod^{\calO}_{A}(\calC)^{\otimes}_{X}$,
and any such diagram is automatically a colimit diagram. From the uniqueness properties
of colimit diagrams, we deduce that $\overline{p}$ is a colimit diagram if and only if it is an
operadic $\psi$-colimit diagram. In view of Theorem \ref{colmodd}, this is true if and only
if $F$ satisfies the following condition:
\begin{itemize}
\item[$(\ast)$] Let $D  = \id_{X} \in \calD$. Then the diagram
$$\calD_{/D} \times K^{\triangleright} \rightarrow \calD \times K^{\triangleright}
\stackrel{F}{\rightarrow} \calC^{\otimes}$$
is an operadic $q$-colimit diagram.
\end{itemize}
Since the inclusion $\{ \id_{D} \} \hookrightarrow \calD_{/D}$ is cofinal, condition
$(\ast)$ is equivalent to the requirement that $\phi \circ \overline{p} = F| \{D\} \times K^{\triangleright}$
is an operadic $q$-colimit diagram. Since $\phi \circ p$ can be extended to an
operadic $q$-colimit diagram in $\calC^{\otimes}_{X}$ (Corollary \ref{cyster}) and any such
diagram is automatically a colimit diagram in $\calC^{\otimes}_{X}$, the uniqueness properties
of colimit diagrams show that $(\ast)$ is equivalent to the requirement that
$\phi \circ \overline{p}$ is a colimit diagram in $\calC^{\otimes}_{X}$.
\end{proof}

\begin{example}
Let $q: \calC^{\otimes} \rightarrow \calO^{\otimes}$ be a fibration of $\infty$-operads, where
$\calO^{\otimes}$ is coherent. Let $A \in \Alg_{\calO}(\calC)$ be a $\calO$-algebra object of
$\calC$. Then the map
$$ \calK_{\calO} \stackrel{ e_1}{\rightarrow} \calO^{\otimes} \stackrel{A}{\rightarrow} \calC^{\otimes}$$
determines an object in $\Alg_{\calO} ( \Mod_{A}^{\calO}( \calC) )$, which we will denote by
$\overline{A}$ (it is a preimage of the identity map $\id_{A}$ under the equivalence
$\Alg_{\calO} ( \Mod_{A}^{\calO}( \calC) ) \simeq \Alg_{\calO}(\calC)^{A/}$ of Corollary \ref{skoke}). 
We can informally summarize the situation by saying that any algebra object
$A \in \Alg_{\calO}(\calC)$ can be viewed as a module over itself. 

Let $0 \in \calO^{\otimes}_{\seg{0}}$ be a zero object of $\calO^{\otimes}$ and let
$X \in \calO$ be any object. Then $\overline{A}(0)$ is a zero object of
$\Mod_{A}^{\calO}(\calC)^{\otimes}$, and $\overline{A}(X)$ is an object of
$\Mod_{A}^{\calO}(\calC)^{\otimes}_{X}$. Any choice of map $0 \rightarrow X$ in
$\calO^{\otimes}$ induces a map $\eta_X: \overline{A}(0) \rightarrow \overline{A}(X)$, which
is given by an edge $\overline{p}: \Delta^1 \rightarrow \Mod_{A}^{\calO}(\calC)^{\otimes}$. 
We claim that $\overline{p}$ is an operadic $\psi$-colimit diagram, where
$\psi: \Mod_{A}^{\calO}(\calC)^{\otimes} \rightarrow \calO^{\otimes}$ denotes the projection.
In view of Theorem \ref{colmodd}, it suffices to prove that $\overline{p}$ induces
an operadic $q$-colimit diagram
$$\theta: ((\calO^{\otimes})^{0/} \times_{\calD} \calD_{/D})^{\triangleright} \rightarrow \calC^{\otimes},$$
where $\calD = \calK_{\calO} \times_{ \calO^{\otimes}} \Delta^1$ and $D$ is the object of
$\calD$ determined by the pair $(\id_{X}, 1)$. We observe that
the fiber product $(\calO^{\otimes})^{0/} \times_{\calD} \calD_{/D}$ contains a final object
$C$, corresponding to the diagram
$$ \xymatrix{ 0 \ar[r] \ar[d] & X \ar[d]^{\id} \\
X \ar[r]^{\id} & X }$$
in $\calO^{\otimes}$. It therefore suffices to show that the restriction
$\theta_0 = \theta| \{ C \}^{\triangleright} \rightarrow \calC^{\otimes}$ is an operadic
$q$-colimit diagram (Remark \ref{woperad}). This is clear, since $\theta_0$ corresponds
to the identity morphism $\id: A(X) \rightarrow A(X)$ in $\calC^{\otimes}$.

We can summarize the situation informally as follows: for every $X \in \calO$, the map
$\eta_{X}$ exhibits $\overline{A}(X) \in \Mod_{A}^{\calO}(\calC)^{\otimes}_{X}$ as
a ``unit object'' with respect to the $\calO$-operad structure on
$\Mod_{A}^{\calO}(\calC)$.
\end{example}

We can deduce Theorem \ref{turkwell} from Theorem \ref{colmodd}:

\begin{proof}[Proof of Theorem \ref{turkwell}]
In view of Corollary \ref{sabus}, it will suffice to show that for each $X \in \calO$, the
fiber $\Mod^{\calO}_{A}(\calC)^{\otimes}_{X}$ is an accessible $\infty$-category.
Let $\calD$ denote the full subcategory of $(\calO^{\otimes})^{X/}$ spanned by the
semi-inert morphisms $f: X \rightarrow Y$, let $\calD_0 \subseteq \calD$ be the full subcategory spanned by those objects for which $f$ is an equivalence, and let $A': \calD_0 \rightarrow \calC^{\otimes}$. We will say that a morphism in $\calD$ is inert if its image in $\calO^{\otimes}$ is inert. We observe that $\Mod^{\calO}_{A}(\calC)^{\otimes}_{X}$ can be identified with a fiber of the
restriction functor
$$ \phi: \Fun^{0}_{ \calO^{\otimes}}( \calD, \calC^{\otimes}) \rightarrow
\Fun^{0}_{\calO^{\otimes}}( \calD_0, \calC^{\otimes} ),$$
where the superscript $0$ indicates that we consider only those functors which carry
inert morphisms in $\calD$ (or $\calD_0$) to inert morphisms in $\calC^{\otimes}$.
It follows from Corollary \toposref{storkus1} that the domain and codomain of $\phi$ are
accessible $\infty$-categories and that $\phi$ is an accessible functor.
Invoking Proposition \toposref{horse2}, we deduce that $\Mod^{\calO}_{A}(\calC)^{\otimes}_{X}$ is accessible as desired.
\end{proof}

We now turn to the proof of Theorem \ref{colmodd}. We will treat assertions
$(1)$ and $(2)$ separately. In both cases, our basic strategy is similar to that of
Theorem \ref{oplk}.

\begin{proof}[Proof of Part $(2)$ Theorem \ref{colmodd}]
Using Propositions \toposref{minimod} and \toposref{princex}, we can assume without loss of generality
that the $\infty$-category $K$ is minimal.
Let $\overline{\calD}^{0}$ denote the inverse image of $\calK_{\calO}^{0}$ in
$\overline{\calD}$, and let $\calD^{0} = \overline{\calD}^{0} \cap \calD$. 
Let $\calD'$ denote the full subcategory of $\overline{\calD}$ spanned by
$\calD$ together with $\overline{\calD}^{0}$. Note that there is a unique map
$z: \calD' \rightarrow \Delta^2$ such that $z^{-1} \Delta^{ \{0,1\} } \simeq \calD$
and $z^{-1} \Delta^{ \{1,2\} } = \overline{\calD}^{0}$. The map $z$ is a coCartesian fibration, and therefore flat. The algebra $A$ and the map $F$ determine a map $\Lambda^2_{1} \times_{ \Delta^2} \calD' \rightarrow \calC^{\otimes}$. Using the fact that $q$ is a categorical fibration and
that the inclusion $\Lambda^2_1 \times_{ \Delta^2} \calD' \subseteq \calD'$ is a categorical
equivalence (Proposition \bicatref{seda}), we can find a map
$F_0 \in \Fun_{\calO^{\otimes}}(\calD', \calC^{\otimes})$ compatible with $F$ and
$A$. To complete the proof, we wish to prove that $F_0$ can be extended to a map
$\overline{F} \in \Fun_{ \calO^{\otimes}}( \overline{\calD}, \calC^{\otimes})$ satisfying
$(\ast)$ together with the following condition (which guarantees that
$\overline{F}$ encodes a diagram $\overline{p}: K^{\triangleright} \rightarrow \Mod^{\calO}_{A}(\calC)^{\otimes}$):
\begin{itemize}
\item[$(\star)$] Let $\alpha: D \rightarrow D'$ be a morphism in $\overline{\calD}$ lying over the cone
point of $K^{\triangleleft}$ whose image in $\calK_{\calO}$ is inert. Then 
$\overline{F}'(\alpha)$ is an inert morphism of $\calC^{\otimes}$.
\end{itemize}
Note that, because $\calC^{\otimes}$ is an $\infty$-operad, it suffices to verify condition
$(\star)$ when the object $D' \in \overline{\calD}$ lies over $\seg{1} \in \Nerve(\FinSeg)$.

Using Proposition \toposref{minimod}, we can choose a
simplicial subset $\overline{\calE} \subseteq \overline{\calD}$ such that the projection
$\overline{\calE} \rightarrow K^{\triangleright}$ is a minimal inner fibration, where
$\overline{\calE}$ is a deformation retract
of $\overline{\calD}$ (fixed over $K^{\triangleright}$); in particular the inclusion
$\overline{\calE} \subseteq \overline{\calD}$ is a categorical equivalence.
It follows from Example \ref{telos} that $\overline{\calE} \rightarrow K^{\triangleright}$ is a strictly unital inner fibration and that $K^{\triangleright}$ is a strictly unital $\infty$-category, so that
$\overline{\calE}$ is strictly unital (Remarks \ref{tope1} and \ref{tope2}).
Let $\calE' = \overline{\calE} \cap \calD'$ and let $F'_0 = F | \calE'$. We will prove that there exists an extension $\overline{F}'
\in \Fun_{ \calO^{\otimes}}( \overline{\calE}, \calC^{\otimes})$ of $F'_0$ satisfying the following
analogues of conditions $(\ast)$ and $(\star)$: 
\begin{itemize}
\item[$(\ast')$] For every object $D \in \overline{\calE}$ whose image in $\overline{\calD}$
is equivalent to a pair $(v, \id_{X})$, the map 
$$ ( \calE_{/D}^{\acti} )^{\triangleright} \rightarrow 
\overline{\calE}_{/D}^{\triangleright} \rightarrow \overline{\calE} \stackrel{ \overline{F}'}{\rightarrow}
\calC^{\otimes}$$
is an operadic $q$-colimit diagram. Here $\calE_{/D}^{\acti}$ denotes the fiber product
$\calE \times_{ \overline{\calD} } \overline{\calD}^{\acti}_{/D}$. 
\item[$(\star')$] Let $\alpha: D \rightarrow D'$ be a morphism in $\overline{\calE}$ lying over the cone
point of $K^{\triangleleft}$ whose image in $\calK_{\calO}$ is inert, such that the image of $D'$
in $\Nerve(\FinSeg)$ is $\seg{1}$. Then $\overline{F}'(\alpha)$ is an inert morphism of $\calC^{\otimes}$.
\end{itemize}
It will then follow from Proposition \toposref{princex} that $F$ admits an extension
$\overline{F} \in \Fun_{ \calO^{\otimes}}( \overline{\calD}, \calC^{\otimes})$ satisfying $(\ast)$ and
$(\star')$, and the proof will be complete.

Let $E$ denote the full subcategory of $\overline{\calE}$ spanned by those objects which lie
over the cone point of $K^{\triangleright}$, and let $E^0 = E \times_{ \calK_{\calO} } \calK_{\calO}^{0}$.
Let $\calJ$ denote the category $(\FinSeg)_{ \seg{1}/ }$ of pointed objects of $\FinSeg$.
There is an evident forgetful functor $E \rightarrow \Nerve(\calJ)$, given by the map
$$ E \subseteq \{X\} \times_{ \calO^{\otimes} } \calK_{\calO} 
\subseteq ( \calO^{\otimes})^{X/} \rightarrow \Nerve( \FinSeg)^{ \seg{1}/} \simeq \Nerve(\calJ).$$
We will say that a morphism $\alpha$ in $\calJ$ is {\it active} or {\it inert} if its
image in $\FinSeg$ is active or inert, respectively; otherwise, we will say that $\alpha$ is {\it neutral}. Let
$\sigma$ be an $m$-simplex of $\Nerve(\calJ)$, corresponding to a chain of morphisms
$$ \seg{1} \stackrel{\alpha(0)}{\rightarrow} (\seg{k_0}) \stackrel{\alpha(1)}{\rightarrow} (\seg{k_1}) \stackrel{\alpha(2)}{\rightarrow}
\cdots \stackrel{\alpha(m)}{\rightarrow} (\seg{k_m}).$$
in the category $\FinSeg$. We will say that $\sigma$ is {\it new} if it is nondegenerate and the map
$\alpha(0)$ is not null. We let $J_{\sigma}$ denote the collection of integers $j \in \{1, \ldots, m\}$ for which the map $\alpha(j)$ is not an isomorphism. We will denote the cardinality of
$J_{\sigma}$ by $l(\sigma)$ and refer to it as the {\it length} of $\sigma$ (note that this length
is generally smaller than $m$).
For $1 \leq d \leq l(\sigma)$, we let $j^{\sigma}_{d}$ denote the $d$th element of
$J_{\sigma}$ and set $\alpha^{\sigma}_{d} = \alpha( j^{\sigma}_{d})$.
We will say that $\sigma$ is {\it closed} if $k_{m} = 1$; otherwise we will say that
$\sigma$ is {\it open}.  We now partition the collection of new simplices $\sigma$ of $E$ into eleven groups, as in the proof of Theorem \ref{oplk}:

\begin{itemize}
\item[$(G'_{(1)})$] A new simplex $\sigma$ of $\Nerve(\calJ)$ belongs to $G'_{(1)}$ if
it is a closed and the maps $\alpha^{\sigma}_{i}$ are active for $1 \leq i \leq l(\sigma)$.

\item[$(G_{(2)})$] A new simplex $\sigma$ of $\Nerve(\calJ)$ belongs
to $G_{(2)}$ if $\sigma$ is closed and there exists
$1 \leq k < l(\sigma)$ such that $\alpha^{\sigma}_{k}$ is inert, while
$\alpha^{\sigma}_{j}$ is active for $k < j \leq l(\sigma)$.

\item[$(G'_{(2)})$]  A new simplex $\sigma$ of $\Nerve(\calJ)$ belongs
to $G'_{(2)}$ if $\sigma$ is closed and there exists $1 \leq k \leq l(\sigma)$ such that
$\alpha^{\sigma}_{k}$ is neutral while the maps $\alpha^{\sigma}_{j}$ are active for $k < j \leq l(\sigma)$.

\item[$(G_{(3)})$] A new simplex $\sigma$ of $\Nerve(\calJ)$ belongs
to $G_{(3)}$ if $\sigma$ is closed and there exists $1 \leq k < l(\sigma)-1$ such that
$\alpha^{\sigma}_{k}$ is inert, the maps $\alpha^{\sigma}_{j}$ are active for
$k < j < l(\sigma)$, and $\alpha^{\sigma}_{l(\sigma)}$ is inert.

\item[$(G'_{(3)})$] A new simplex $\sigma$ of $\Nerve(\calJ)$ belongs
to $G'_{(3)}$ if $\sigma$ is closed and there exists $1 \leq k < l(\sigma)$ such that the map $\alpha^{\sigma}_{k}$ is neutral, the maps $\alpha^{\sigma}_{j}$ are active for $k < j < l(\sigma)$, and
$\alpha^{\sigma}_{l(\sigma)}$ is inert. 

\item[$(G_{(4)})$] A new simplex $\sigma$ of $\Nerve(\calJ)$ belongs to $G_{(4)}$ if
it is a closed, the maps $\alpha^{\sigma}_{i}$ are active for $1 \leq i < l(\sigma)$, and
the map $\alpha^{\sigma}_{l(\sigma)}$ is inert.

\item[$(G'_{(4)})$] A new simplex $\sigma$ of $\Nerve(\calJ)$ belongs to $G'_{(4)}$ if
it is an open and the maps $\alpha^{\sigma}_{i}$ are active for $1 \leq i \leq l(\sigma)$.

\item[$(G_{(5)})$] A new simplex $\sigma$ of $\Nerve(\calJ)$ belongs
to $G_{(5)}$ if $\sigma$ is an open and there exists $1 \leq k < l(\sigma)$ such that
$\alpha^{\sigma}_{k}$ is inert and $\alpha^{\sigma}_{j}$ is active for
$k < j \leq l(\sigma)$. 

\item[$(G'_{(5)})$] A new simplex $\sigma$ of $\Nerve(\calJ)$ belongs
to $G'_{(5)}$ if $\sigma$ is an open and there exists $1 \leq k \leq l(\sigma)$ such that
$\alpha^{\sigma}_{k}$
is neutral and the maps $\alpha^{\sigma}_{j}$ are active for $k < j \leq l(\sigma)$.

\item[$(G_{(6)})$] A new simplex $\sigma$ of $\Nerve(\calJ)$ belongs to $G_{(6)}$ if
it is closed, has length $\geq 2$, and the maps $\alpha^{\sigma}_{l(\sigma)-1}$ and
$\alpha^{\sigma}_{l(\sigma)}$ are both inert.

\item[$(G'_{(6)})$] A new simplex $\sigma$ of $\Nerve(\calJ)$ belongs to $G'_{(6)}$ if
it is open, has length at least $1$, and the map $\alpha^{\sigma}_{l(\sigma)}$ is inert.
\end{itemize}

For each integer $m \geq 0$, we let $\Nerve(\calJ)_{m}$ denote the simplicial
subset spanned by those simplices which are either not new, have length
$\leq m$, or have length $m$ and belong to one of the groups 
$G_{(i)}$ for $2 \leq i \leq 6$. Let $E(m)$ denote the inverse image
$E \times_{ \Nerve(\calJ) } \Nerve(\calJ)_{m}$ and let 
$\overline{\calE}(m)$ denote the simplicial subset
of $\overline{\calE}$ spanned by those simplices whose intersection with
$E$ belongs to $E(m)$. Then $\overline{\calE}(0) = \calE'$ is the domain of the map
$F'_0$. We will complete the proof by extending $F'_0$ to a compatible sequence of maps
$F'_{m} \in \Fun_{ \calO^{\otimes}}( \overline{\calE}(m), \calC^{\otimes})$, where
$F'_{1}$ satisfies conditions $(\ast')$ and $(\star')$.

Let us now fix $m > 0$ and assume that $F'_{m-1}$ has already been constructed.
We define a filtration
$$ \Nerve(\calJ)_{m-1} = K(0) \subseteq K(1) \subseteq K(2) \subseteq K(3) \subseteq K(4) \subseteq K(5) \subseteq K(6) = \Nerve(\calJ)_{m}$$
as follows:
\begin{itemize}
\item We let $K(1)$ denote the simplicial subset of $\Nerve(\calJ)$ spanned by those simplices which either belong to $K(0)$ or have length $(m-1)$ and belong
to $G'_{(1)}$. 
\item For $2 \leq i \leq 6$, we let $K(i)$ be the simplicial subset of $\Nerve(\calJ)$ spanned by those simplices which either belong to $K(i-1)$, have length $m$ and belong to $G_{(i)}$, or have length $m-1$ and belong to $G'_{(i)}$. 
\end{itemize}
For $0 \leq i \leq 6$, we let $\overline{K}(i)$ denote the simplicial subset of
$\overline{\calE}$ spanned by those simplices whose intersection
with $E$ belongs to the inverse image of $K(i)$. We will define maps
$f^{i}: \overline{K}(i) \rightarrow \calC^{\otimes}$ with $f^{0} = F'_{m-1}$.
The construction now proceeds in six steps:

\begin{itemize} 
\item[$(1)$] Assume that $f^{0} = F'_{m-1}$ has been constructed; we wish to define
$f^{1}$. Let $\{ \sigma_{a} \}_{a \in A}$ be the collection of all 
simplices of $E$ whose image in $\Nerve(\calJ)$ have length
$(m-1)$ and belong to $G'_{(1)}$.
Choose a well-ordering of the set $A$ such that the dimensions of the simplices
$\sigma_{a}$ form a (nonstrictly) increasing function of $a$. For each $a \in A$, let $\overline{\calE}_{\leq a}$ denote the
simplicial subset of $\overline{\calE}$ spanned by those simplices which
either belong to $\overline{K}(0)$ or whose intersection with
$E$ is contained in $\sigma_{a'}$ for some $a' \leq a$, and define $\overline{\calE}_{< a}$ similarly. We construct a compatible family of maps
$f^{\leq a} \in \Fun_{\calO^{\otimes}}(\overline{\calE}_{\leq a}, \calC^{\otimes})$ extending $f^{0}$, using transfinite induction on $a$. Assume that $f^{\leq a'}$ has been constructed for $a' < a$; these
maps can be amalgamated to obtain a map $f^{<a} \in \Fun_{\calO^{\otimes}}( \overline{\calE}_{< a},\calC^{\otimes})$. Let $Z = \calE \times_{ \overline{\calE} } \overline{\calE}_{/ \sigma_{a} }$. 
Since $\overline{\calE}$ is strictly unital, Proposition \ref{appstrict} guarantees that we have
a pushout diagram of simplicial sets
$$ \xymatrix{ Z \star \bd \sigma_a \ar[r] \ar[d] & \overline{\calE}_{< a} \ar[d] \\
Z \star \sigma_{a} \ar[r] & \overline{\calE}_{\leq a}. }$$
It will therefore suffice to extend the composition
$$ g_0: Z \star \bd \sigma_{a} \rightarrow \overline{\calE}_{< a} 
\stackrel{f^{<a}}{\rightarrow} \calC^{\otimes}$$
to a map $g \in \Fun_{\calO^{\otimes}}(Z \star \star \sigma_{a}, \calC^{\otimes})$.

We first treat the special case where the simplex $\sigma_{a}$ is zero-dimensional
(in which case we must have $m=1$). We can identify $\sigma_{a}$ with
an object $D \in \overline{\calE}$. Since $\sigma_{a}$ is new and closed, the image
of $D$ in $\overline{\calD}$ is equivalent to $(v, \id_{X})$. Let $Z_0 = 
\calE \times_{ \overline{\calD} } \overline{\calD}^{\acti}_{/D}$. It follows from
assumption $(ii)$ that $g_0 | Z_0$ can be can be extended to an operadic
$q$-colimit diagram (compatible with the projection to $\calO^{\otimes})$. 
Since $Z_0$ is a localizatin of $Z$, the inclusion $Z_0 \subseteq Z$ is cofinal; we can therefore extend $g_0$ to a map $g \in \Fun_{\calO^{\otimes}}(Z \star \star \sigma_{a}, \calC^{\otimes})$.
whose restriction to $Z_0^{\triangleright}$ is an operadic $q$-colimit diagram.
Moreover, this construction guarantees that $f^{1}$ will satisfy condition $(\ast')$ for
the object $D$.

Now suppose that $\sigma_{a}$ is a simplex of positive dimension. We again let
$Z_0$ denote the simplicial subset of $Z$ spanned by those vertices 
which correspond to diagrams $\sigma_{a}^{\triangleleft} \rightarrow
\calC^{\otimes}$ which project to a sequence of active morphisms in
$\Nerve(\FinSeg)$. The inclusion $Z_0 \subseteq Z$ admits a left adjoint and is therefore
cofinal; it follows that the induced map $\calC^{\otimes}_{(F_{m-1} | Z)/} \rightarrow
\calC^{\otimes}_{(F_{m-1}|Z_0)/} \times_{ \calO^{\otimes}_{ (q F_{m-1} | Z_0)/} }
\calO^{\otimes}_{(q F_{m-1}|Z)/}$ is a trivial Kan fibration. It therefore suffices to show that
the restriction $g'_0 = g_0 | (Z_0 \star \bd \sigma_{a})$ can be extended to a map
$g' \in \Fun_{\calO^{\otimes}}( Z_0 \star \sigma_{a}, \calC^{\otimes})$ .
Let $D$ denote the initial vertex of $\sigma$. In view of
Proposition \ref{quadly}, it will suffice to show that the restriction
$g'_0 | (Z_0 \star \{D\})$ is an operadic $q$-colimit diagram.
Since the inclusion $\{D\} \subseteq \sigma_{a}$ is left anodyne, the projection
map $Z_0 \rightarrow \calE^{\acti}_{/D} = \calE \times_{ \overline{\calD} } \overline{\calD}^{\acti}_{/D}$ is a trivial Kan fibration. It will therefore suffice to show that $F_{m-1}$ induces an
operadic $q$-colimit diagram $\delta: (\calE^{\acti}_{/D})^{\triangleright} \rightarrow \calC^{\otimes}$.

The object $D \in \overline{\calD}$ determines a semi-inert morphism
$\gamma: X \rightarrow Y$ in $\calO^{\otimes}$. Since $\sigma_{a}$ is new, this
morphism is not null and therefore determines an equivalence $Y \simeq Y_0 \oplus X$.
Let $\calX$ denote the full subcategory of $\calE^{\acti}_{/D}$ spanned by those objects
corresponding to diagrams
$$ \xymatrix{ X' \ar[r] \ar[d] & Y' \ar[d]^{\phi} \\
X \ar[r] & Y_0 \oplus X}$$
where the active morphism $\phi$ exhibits $Y'$ as a sum
$Y'_0 \oplus Y'_1$, where $Y'_0 \simeq Y_0$ and $Y'_1 \rightarrow X$ is an active morphism.
The inclusion $\calX \subseteq \calE^{\acti}_{/D}$ admits a left adjoint, and is therefore cofinal.
Let $D_0 = (v, \id_{X}) \in \overline{\calD}$, so that the operation
$\bigdot \mapsto Y_0 \oplus \bigdot$ determines an equivalence $\calE^{\acti}_{/D_0}
\rightarrow \calX$. It therefore suffices to show that the composite map
$$ (\calE^{\acti}_{/D_0})^{\triangleright} \stackrel{ Y_0 \oplus}{\rightarrow}
\calX^{\triangleright} \subseteq (\calE^{\acti}_{/D})^{\triangleright} \stackrel{\delta}{\rightarrow} \calC^{\otimes}$$ is an operadic $q$-colimit diagram. Using condition $(\star')$, we can identify
this map with the composition
$$ (\calE^{\acti}_{/D_0})^{\triangleright} \rightarrow \calC^{\otimes}
\stackrel{ A(Y_0) \oplus }{\rightarrow} \calC^{\otimes},$$
which is an operadic $q$-colimit diagram by virtue of $(\ast')$.

\item[$(2)$] We now assume that $f^{1}$ has been constructed. 
Since $q$ is a categorical fibration, to produce the desired extension
$f^{2}$ of $f^{1}$ it is sufficient to show that the inclusion
$\overline{K}(1) \subseteq \overline{K}(2)$ is a categorical equivalence.
For each simplex $\sigma$ of 
$\Nerve(\calJ)$ of length $m$ belonging to $G_{(2)}$, let $k(\sigma) < m$ be the integer such that
$\alpha^{\sigma}_{k(\sigma)}$ is inert while
$\alpha^{\sigma}_{j}$ is active for $k(\sigma) < j \leq m$. We will say that
$\sigma$ is {\it good} if $\alpha^{\sigma}_{k(\sigma)}$ induces a map
$\seg{p} \rightarrow \seg{p'}$ in $\FinSeg$ whose restriction to $(\alpha^{\sigma}_{k})^{-1} \nostar{p'}$
is order preserving. Let $J = K^{\triangleright} \times_{ \Fun( \{0\}, \Nerve(\FinSeg)} \Fun( \Delta^1, \Nerve(\FinSeg) )$, so that we can identify $\Nerve(\calJ)$ with the simplicial subset
$J \times_{ K^{\triangleright} } \{v\}$. 
Let $\{ \sigma_{a} \}_{a \in A}$ be the collection of all nondegenerate simplices
of $J$ such that the intersection
$\sigma'_{a} = \sigma_{a} \cap \Nerve(\calJ)$ is nonempty and good. For each $a \in A$, let $k_{a} = k(\sigma'_{a})$ 
and $j_{a} = j^{\sigma'_{a}}_{k}$. Choose a well-ordering of $A$ with the following
properties:
\begin{itemize}
\item The map $a \mapsto k_{a}$ is a (nonstrictly) increasing function of $a \in A$.
\item For each integer $k$, the dimension of the simplex $\sigma_{a}$ is a nonstrictly
increasing function of $a \in A_{k} = \{ a \in A: k_a = k \}$. 
\item Fix integers $k,d \geq 0$, and let $A_{k,d}$
be the collection of elements $a \in A_{k}$ such that $\sigma_{a}$ has dimension $d$.
The map $a \mapsto j_{a}$ is a nonstrictly increasing function on
$A_{k,d}$.
\end{itemize}
Let $J_{< 0}$ be the collection of those simplices of $J$ whose intersection
with $\Nerve(\calJ)$ belongs to $K(1)$. For each
$a \in A$, let $J_{\leq a}$ denote the simplicial subset of
$J$ generated by $J_0$ together with the simplices
the simplices $\{ \sigma_{a'} \}_{a' \leq a}$, and define $J_{<a}$ similarly.
The inclusion $\overline{K}(1) \subseteq \overline{K}(2)$ can be
obtained as a transfinite composition of inclusions
$$ i_{a}: \overline{\calE} \times_{J}
J_{< a} \rightarrow
\overline{\calE} \times_{J} J_{\leq a}.$$
Each $i_{a}$ is a pushout of an inclusion
$i'_{a}: \overline{\calE} \times_{ J} \sigma^{0}_{a}
\overline{\calE} \times_{J} \sigma_{a},$
where $\sigma^{0}_{a} \subseteq \sigma_{a}$ denotes the inner horn
obtained by removing the interior of $\sigma_{a}$ together with the
face opposite the $j_{a}$th vertex of $\sigma'_{a}$. We wish to prove
that $i'_{a}$ is a categorical equivalence.
Using the fact that $\overline{\calE}$ is a fiberwise deformation retract
of $\overline{\calD}$, we are reduced to proving that the map
$\overline{\calD} \times_{ J} \sigma^{0}_{a}
\overline{\calD} \times_{J} \sigma_{a}$ is a categorical equivalence, which follows from
Lemma \ref{skope}. 

\item[$(3)$] To find the desired extension $f^{3}$ of $f^{2}$, it suffices to show that
the inclusion $\overline{K}(2) \subseteq \overline{K}(3)$ is a categorical equivalence. This follows from the argument given in step $(2)$.

\item[$(4)$] Let $\{ \sigma_{a} \}_{a \in A}$ denote the collection of all nondegenerate simplices
$\sigma$ of $\overline{\calE}$ with the property that $\sigma \cap E$ is nonempty and
projects to a simplex of $\Nerve(\calJ)$ of length $m-1$ which belongs to $G'_{4}$.
Choose a well-ordering of $A$ having the property that the dimensions of the simplices $\sigma_{a}$ form a (nonstrictly) increasing function of $a$. For 
each $a \in A$ let $D_{a} \in \overline{\calE}$ denote the final vertex of $\sigma_{a}$
and let $Z_{a}$ denote the full subcategory of $E \times_{ \overline{\calE} } \overline{\calE}_{\sigma_a/}$ spanned by those objects for which the underlying map $D_{a} \rightarrow D$ induces an
inert morphism $\seg{m} \rightarrow \seg{1}$ in $\Nerve(\FinSeg)$.
We have a canonical map
$t_{a}: \sigma_{a} \star Z_{a} \rightarrow \overline{\calE}$. For each
$a \in A$, let $\overline{\calE}_{\leq a} \subseteq \overline{\calE}$ denote the simplicial subset generated by $\overline{K}(3)$ together with the image of $t_{a'}$ for all $a' \leq a$, and define
$\overline{\calE}_{< a}$ similarly. Using our assumption that $\overline{\calE}$ is strictly nital
(and a mild variation on Proposition \ref{appstrict}), we deduce that $\overline{K}(4) = \bigcup_{a \in A} \overline{\calE}_{\leq a}$ and that for each $a \in A$ we have a pushout diagram of simplicial sets
$$ \xymatrix{ \bd \sigma_{a} \star Z_{a} \ar[r] \ar[d] & \overline{\calE}_{< a} \ar[d] \\
\sigma_{a} \star Z_{a} \ar[r] & \overline{\calE}_{\leq a}. }$$
To construct $f^{4}$, we are reduced to the problem of solving a sequence of extension
problems of the form
$$ \xymatrix{ \bd \sigma_{a} \star Z_{a} \ar[r]^{g_0} \ar[d] & \calC^{\otimes} \ar[d]^{q} \\
\sigma_{a} \star Z_{a} \ar[r] \ar@{-->}[ur]^{g} & \calO^{\otimes}. }$$
Note that $Z_a$ decomposes a disjoint union of $\infty$-categories
$\coprod_{1 \leq i \leq m} Z_{a,i}$ (where $\seg{m}$ denotes the image of $D_a$ in $\Nerve(\FinSeg)$).
Each of the $\infty$-categories $Z_{a,i}$ has an initial object $B_{i}$, given by any
map $\sigma_{a} \star \{ D_{i} \} \rightarrow \calC^{\otimes}$ which induces an inert
morphism $D_{a} \rightarrow D_i$ covering $\Colp{i}: \seg{m} \rightarrow \seg{1}$.
Let $h: Z_{a} \rightarrow \calC^{\otimes}$ be the map induced by $g_0$, and let
$h'$ be the restriction of $h$ to the discrete simplicial set $Z'_{a} = \{ B_{i} \}_{1 \leq i \leq m}$.
Since the inclusion $Z'_{a} \rightarrow Z_a$ is left anodyne, we have
a trivial Kan fibration $\calC^{\otimes}_{/h} \rightarrow \calC^{\otimes}_{/h'}
\times_{ \calO^{\otimes}_{/qh'}} \calO^{\otimes}_{/qh}$. Unwinding the definitions, we are reduced to the lifting probelm depicted in the diagram
$$ \xymatrix{ \bd \sigma_{a} \star Z'_{a} \ar[r]^{g'_0} \ar[d] & \calC^{\otimes} \ar[d]^{q} \\
\sigma_{a} \star Z'_{a} \ar[r]^{\overline{g}'} \ar@{-->}[ur]^{g'} & \calO^{\otimes}. }$$
If the dimension of $\sigma_{a}$ is positive, then it suffices to show that
$g'_0$ carries $\{D_a \} \star Z'_{a}$ to a $q$-limit diagram in $\calC^{\otimes}$.
Let $q'$ denote the canonical map $\calO^{\otimes} \rightarrow \Nerve(\FinSeg)$.
In view of Proposition \toposref{basrel}, it suffices to show that
$g'_0$ carries $\{D_a \} \star Z'_{a}$ to a $(q' \circ q)$-limit diagram and that
$q \circ g'_0$ carries $\{D_a \} \star Z'_{a}$ to a $q'$-limit diagram. The first of these
assertions follows from $(\star')$ and from the fact that $\calC^{\otimes}$ is an $\infty$-operad,
and the second follows by the same argument (since $q$ is a map of $\infty$-operads).

It remains to treat the case where $\sigma_{a}$ is zero dimensional (in which case
we have $m=1$). Since $\calC^{\otimes}$ is an $\infty$-operad, we can solve the lifting problem depicted in the diagram
$$ \xymatrix{ \bd \sigma_{a} \star Z'_{a} \ar[r]^{g'_0} \ar[d] & \calC^{\otimes} \ar[d]^{q' \circ q} \\
\sigma_{a} \star Z'_{a} \ar[r] \ar@{-->}[ur]^{g''} & \Nerve(\FinSeg) }$$
in such a way that $g''$ carries edges of $\sigma_{a} \star Z'_{a}$ to inert morphisms
in $\calC^{\otimes}$. Since $q$ is an $\infty$-operad map, it follows that
$q \circ g''$ has the same property. Since $\calO^{\otimes}$ is an $\infty$-operad, we conclude that
$q \circ g''$ and $\overline{g}'$ are both $q'$-limit diagrams in $\calO^{\otimes}$
extending $q \circ g'_0$, and therefore equivalent to one another via
an equivalence which is fixed on $Z'_{a}$ and compatible with the projection to
$\Nerve(\FinSeg)$. Since $q$ is a categorical fibration, we can lift this equivalence to an
equivalence $g'' \simeq g'$, where $g': \sigma_{a} \star Z'_{a} \rightarrow \calC^{\otimes}$
is the desired extension of $g'_0$. We note that this construction ensures that condition
$(\star')$ is satisfied.

\item[$(5)$] To find the desired extension $f^5$ of $f^4$, it suffices to show that
the inclusion $\overline{K}(4) \subseteq \overline{K}(5)$ is a categorical equivalence. This again follows from the argument given in step $(2)$.

\item[$(6)$] The verification that $f^{5}$ can be extended to a map
$f^{6}: \overline{K}(6) \rightarrow \calC^{\otimes}$ proceeds as in step $(4)$, but is slightly easier
(since $G'_{(6)}$ contains no $0$-simplices).
\end{itemize}
\end{proof}

\begin{proof}[Proof of Part $(1)$ of Theorem \ref{colmodd}]
We wish to show that a diagram $\overline{p}: K^{\triangleright} \rightarrow
\Mod^{\calO}_{A}(\calC)^{\otimes}$ satisfies condition $(\ast)$ if and only if it is an
operadic $\psi$-colimit diagram. We will prove the ``only if'' direction; the converse will follow
from assertion $(2)$ together with the uniqueness properties of operadic colimit diagrams.
Using Proposition \toposref{minimod}, we may assume without loss of generality that
$K$ is a minimal $\infty$-category.
Fix an object $\overline{Y} \in \Mod^{\calO}_{A}( \calC)^{\otimes}$ lying over
$Y \in \calO^{\otimes}$ and let
$\overline{p}_{\overline{Y}}: K^{\triangleright} \rightarrow \Mod^{\calO}_{A}(\calC)^{\otimes}$;
by given by the composition
$$ K^{\triangleright} \stackrel{\overline{p}}{\rightarrow} 
\Mod^{\calO}_{A}(\calC)^{\otimes} \stackrel{\oplus \overline{Y}}{\rightarrow}
\Mod^{\calO}_{A}(\calC)^{\otimes};$$
we must show that $\overline{p}_{\overline{Y}}$ is a weak $\psi$-operadic colimit diagram.
Unwinding the definitions, we must show that for $n > 0$, every lifting problem of the form
$$ \xymatrix{ K \star \bd \Delta^n \ar[d] \ar[r]^{p'} & \Mod^{\calO}_{A}(\calC)^{\otimes}_{\acti} \ar[d] \\
K \star \Delta^n \ar[r]^{\overline{p}'_0} \ar@{-->}[ur] & \calO^{\otimes}_{\acti} }$$ 
admits a solution, provided that $p' | K \star \{0\} = \overline{p}_{ \overline{Y} }$ and
$\overline{p}'_0$ carries $\Delta^{ \{1, \ldots, n\} }$ into $\calO \subseteq \calO^{\otimes}_{\acti}$.

Let $\overline{\calT}$ denote the fiber product $(K \star \Delta^n) \times_{ \calO^{\otimes} }
\calK_{\calO}$, let $\overline{\calT}^{0} = (K \star \Delta^n) \times_{ \calO^{\otimes}} \calK_{\calO}$.
Let $\overline{\calT}^{1}$ denote the full subcategory of
$\overline{\calT} = (K \star \{0\}) \times_{ \calO^{\otimes} } \calK_{\calO}$ spanned by
$\calT = K \times_{ \calO^{\otimes} } \calK_{\calO}$ together with those vertices of
$\{0\} \times_{ \calO^{\otimes} } \calK_{\calO}$ which correspond to semi-inert morphisms
$X \oplus Y \rightarrow Y'$ in $\calO^{\otimes}$ such that the underlying map in
$\FinSeg$ is not injective. 
Let $\calT'$ denote the full subcategory of $\overline{\calT}$ spanned by
$\overline{\calT}^{1}$ together with $\overline{\calT}^{0}$. Note that there is a unique map
$z: \calT' \rightarrow \Delta^2$ such that $z^{-1} \Delta^{ \{0,1\} } = \overline{\calT}^{1}$
and $z^{-1} \Delta^{ \{1,2\} } = \overline{\calT}^{0}$. The map $z$ is a coCartesian fibration and therefore flat. The algebra $A$ and the map $\overline{p}$ determine a map $\Lambda^2_{1} \times_{ \Delta^2} \calT' \rightarrow \calC^{\otimes}$. Using the fact that $q$ is a categorical fibration and
that the inclusion $\Lambda^2_1 \times_{ \Delta^2} \calT' \subseteq \calT'$ is a categorical
equivalence (Proposition \bicatref{seda}), we can find a map
$F' \in \Fun_{\calO^{\otimes}}(\calT', \calC^{\otimes})$ compatible with $F$ and
$A$. Amalgamating $F_0$ with the map determined by $p'$, we obtain a map
$F'' \in \Fun_{ \calO^{\otimes}}( \calT'', \calC^{\otimes})$ where
$$\calT'' = (( K \star \bd \Delta^n) \times_{ \calO^{\otimes}} \calK_{\calO})
\coprod_{ \calT } \calT' \subseteq \overline{\calT}$$
To complete the proof, we wish to prove that $F''$ can be extended to a map
$\overline{F} \in \Fun_{ \calO^{\otimes}}( \overline{\calT}, \calC^{\otimes})$. 

Using Proposition \toposref{minimod}, we can choose a
simplicial subset $\overline{\calE} \subseteq \overline{\calT}$ such that the projection
$\overline{\calE} \rightarrow K \star \Delta^n$ is a minimal inner fibration, where
$\widetilde{\calE}$ is a deformation retract
of $\overline{\calT}$; in particular the inclusion
$\overline{\calE} \subseteq \overline{\calT}$ is a categorical equivalence.
It follows as in the previous proof that $\overline{\calE}$ is a strictly unital $\infty$-category.
Let $\calE'' = \overline{\calE} \cap \calT''$ and let $F''_{0} = F'' | \calE'$. We will prove that there exists an extension $\overline{F}'' \in \Fun_{ \calO^{\otimes}}( \overline{\calE}, \calC^{\otimes})$ of $F''_0$. 
It will then follow from Proposition \toposref{princex} that $F''$ admits an extension
$\overline{F} \in \Fun_{ \calO^{\otimes}}( \overline{\calT}, \calC^{\otimes})$ and the proof will be complete.

Let $E$ denote the fiber product $\overline{\calE} \times_{ K \star \Delta^n} \Delta^n$ and let
let $E^0 = E \times_{ \calK_{\calO} } \calK_{\calO}^{0}$. The diagram
$K \star \Delta^n \rightarrow \calO^{\otimes}$ determines a map $\Delta^n \rightarrow
\Nerve(\FinSeg)$, corresponding to a composable sequence of active morphisms
$$ \seg{j+1} \rightarrow \seg{1} \simeq \seg{1} \simeq \cdots \simeq \seg{1}.$$
Here $\seg{j}$ corresponds to the image of $Y$ in $\Nerve(\FinSeg)$. Let
$\calJ$ denote the fiber product category $[n] \times_{ \Fun( \{0\}, \FinSeg)} \Fun( [1], \FinSeg)$.
There is an evident forgetful functor $E \rightarrow \Nerve(\calJ)$, whose second projection
is given by the composition
$$ E \subseteq \Delta^n \times_{ \calO^{\otimes} } \calK_{\calO} 
\rightarrow \Fun( \Delta^1, \calO^{\otimes}) \rightarrow \Nerve( \Fun( [1], \calJ)).$$ 
We will say that a morphism $\alpha$ in $\calJ$ is {\it active} or {\it inert} if its
image in $\FinSeg$ is active or inert, respectively; otherwise, we will say that $\alpha$ is {\it neutral}. Let
$\sigma$ be an $m$-simplex of $\Nerve(\calJ)$. We will say that $\sigma$ is {\it new} if
$\sigma$ is nondegenerate, the map $\sigma \rightarrow \Delta^n$ is surjective, and 
the semi-inert morphism $\seg{j+1} \rightarrow \seg{k_0}$ given by the first vertex of
$\sigma$ is injective. Every such simplex determines a chain of morphisms
$$ \seg{k_0} \stackrel{\alpha(1)}{\rightarrow} \seg{k_1} \stackrel{\alpha(2)}{\rightarrow}
\cdots \stackrel{\alpha(m)}{\rightarrow} \seg{k_m}.$$
in the category $\FinSeg$. We let $J_{\sigma}$ denote the collection of integers $j \in \{1, \ldots, m\}$ for which the map $\alpha(j)$ is not an isomorphism. We will denote the cardinality of
$J_{\sigma}$ by $l(\sigma)$ and refer to it as the {\it length} of $\sigma$ (note that this length
is generally smaller than $m$).
For $1 \leq d \leq l(\sigma)$, we let $j^{\sigma}_{d}$ denote the $d$th element of
$J_{\sigma}$ and set $\alpha^{\sigma}_{d} = \alpha( j^{\sigma}_{d})$.
We will say that $\sigma$ is {\it closed} if $k_{m} = 1$; otherwise we will say that
$\sigma$ is {\it open}. As in the proof of assertion $(2)$ (and the proof of Theorem
\ref{oplk}), we partition the new simplices of $\Nerve(\calJ)$ into eleven groups
$\{ G_{(i)} \}_{2 \leq i \leq 6}$, $\{ G'_{(i)} \}_{1 \leq i \leq 6}$. 

For each integer $m \geq 0$, we let $\Nerve(\calJ)_{m}$ denote the simplicial
subset spanned by those simplices which are either not new, have length
$\leq m$, or have length $m$ and belong to one of the groups 
$G_{(i)}$ for $2 \leq i \leq 6$. Let $E(m)$ denote the inverse image
$E \times_{ \Nerve(\calJ) } \Nerve(\calJ)_{m}$ and let 
$\overline{\calE}(m)$ denote the simplicial subset
of $\overline{\calE}$ spanned by those simplices whose intersection with
$E$ belongs to $E(m)$. Then $\overline{\calE}(0) = \calE''$ is the domain of the map
$F''_0$. We will complete the proof by extending $F'_0$ to a compatible sequence of maps
$F''_{m} \in \Fun_{ \calO^{\otimes}}( \overline{\calE}(m), \calC^{\otimes})$.

Fix $m > 0$ and assume that $F''_{m-1}$ has already been constructed.
We define a filtration
$$ \Nerve(\calJ)_{m-1} = K(0) \subseteq K(1) \subseteq K(2) \subseteq K(3) \subseteq K(4) \subseteq K(5) \subseteq K(6) = \Nerve(\calJ)_{m}$$
as follows:
\begin{itemize}
\item We let $K(1)$ denote the simplicial subset of $\Nerve(\calJ)$ spanned by those simplices which either belong to $K(0)$ or have length $(m-1)$ and belong
to $G'_{(1)}$. 
\item For $2 \leq i \leq 6$, we let $K(i)$ be the simplicial subset of $\Nerve(\calJ)$ spanned by those simplices which either belong to $K(i-1)$, have length $m$ and belong to $G_{(i)}$, or have length $m-1$ and belong to $G'_{(i)}$. 
\end{itemize}
For $0 \leq i \leq 6$, we let $\overline{K}(i)$ denote the simplicial subset of
$\overline{\calE}$ spanned by those simplices whose intersection
with $E$ belongs to the inverse image of $K(i)$. We will define maps
$f^{i}: \overline{K}(i) \rightarrow \calC^{\otimes}$ with $f^{0} = F''_{m-1}$.
We will explain how to construct $f^{1}$ from $f^{0}$; the remaining steps
can be handled as in our proof of part $(2)$.

Assume that $f^{0} = F''_{m-1}$ has been constructed; we wish to define
$f^{1}$. Let $\{ \sigma_{a} \}_{a \in A}$ be the collection of all 
simplices of $E$ whose image in $\Nerve(\calJ)$ have length
$(m-1)$ and belong to $G'_{(1)}$.
Choose a well-ordering of the set $A$ such that the dimensions of the simplices
$\sigma_{a}$ form a (nonstrictly) increasing function of $a$. 
For each $a \in A$, let $\overline{\calE}_{\leq a}$ denote the
simplicial subset of $\overline{\calE}$ spanned by those simplices which
either belong to $\overline{K}(0)$ or whose intersection with
$E$ is contained in $\sigma_{a'}$ for some $a' \leq a$, and define $\overline{\calE}_{< a}$ similarly. We construct a compatible family of maps
$f^{\leq a} \in \Fun_{\calO^{\otimes}}(\overline{\calE}_{\leq a}, \calC^{\otimes})$ extending $f^{0}$, using transfinite induction on $a$. Assume that $f^{\leq a'}$ has been constructed for $a' < a$; these
maps can be amalgamated to obtain a map $f^{<a} \in \Fun_{\calO^{\otimes}}( \overline{\calE}_{< a},\calC^{\otimes})$. Let $\calE = \overline{\calE} \times_{ K \star \Delta^n} K$, and let
$Z = \calE \times_{ \overline{\calE} } \overline{\calE}_{/ \sigma_{a} }$. 
Since $\overline{\calE}$ is strictly unital, Proposition \ref{appstrict} guarantees that we have
a pushout diagram of simplicial sets
$$ \xymatrix{ Z \star \bd \sigma_a \ar[r] \ar[d] & \overline{\calE}_{< a} \ar[d] \\
Z \star \sigma_{a} \ar[r] & \overline{\calE}_{\leq a}. }$$
It will therefore suffice to extend the composition
$$ g_0: Z \star \bd \sigma_{a} \rightarrow \overline{\calE}_{< a} 
\stackrel{f^{<a}}{\rightarrow} \calC^{\otimes}$$
to a map $g \in \Fun_{\calO^{\otimes}}(Z \star \star \sigma_{a}, \calC^{\otimes})$.

Note that $\sigma_{a}$ necessarily has positive dimension (since
the map $\sigma_{a} \rightarrow \Delta^n$ is surjective). Let
$Z_0$ denote the simplicial subset of $Z$ spanned by those vertices 
which correspond to diagrams $\sigma_{a}^{\triangleleft} \rightarrow
\calC^{\otimes}$ which project to a sequence of active morphisms in
$\Nerve(\FinSeg)$. The inclusion $Z_0 \subseteq Z$ admits a left adjoint and is therefore
cofinal; it follows that the induced map $\calC^{\otimes}_{(F_{m-1} | Z)/} \rightarrow
\calC^{\otimes}_{(F_{m-1}|Z_0)/} \times_{ \calO^{\otimes}_{ (q F_{m-1} | Z_0)/} }
\calO^{\otimes}_{(q F_{m-1}|Z)/}$ is a trivial Kan fibration. It therefore suffices to show that
the restriction $g'_0 = g_0 | (Z_0 \star \bd \sigma_{a})$ can be extended to a map
$g' \in \Fun_{\calO^{\otimes}}( Z_0 \star \sigma_{a}, \calC^{\otimes})$ .
Let $D$ denote the initial vertex of $\sigma$. In view of
Proposition \ref{quadly}, it will suffice to show that the restriction
$g'_0 | (Z_0 \star \{D\})$ is an operadic $q$-colimit diagram.
Since the inclusion $\{D\} \subseteq \sigma_{a}$ is left anodyne, the projection
map $Z_0 \rightarrow \calE^{\acti}_{/D} = \calE \times_{ \overline{\calT} } \overline{\calT}^{\acti}_{/D}$ is a trivial Kan fibration. It will therefore suffice to show that $F_{m-1}$ induces an
operadic $q$-colimit diagram $\delta: (\calE^{\acti}_{/D})^{\triangleright} \rightarrow \calC^{\otimes}$.

The object $D \in \overline{\calT}$ determines a semi-inert morphism
$\gamma: X \oplus Y \rightarrow Y'$ in $\calO^{\otimes}$. Since $\sigma_{a}$ is new,
the underlying morphism in $\Nerve(\FinSeg)$ is injective, so that $\gamma$ induces an
equivalence $Y' \simeq X \oplus Y \oplus Y''$, for some $Y'' \in \calO^{\otimes}$.
Let $\calX$ denote the full subcategory of $\calE^{\acti}_{/D}$ spanned by those objects
corresponding to diagrams
$$ \xymatrix{ X_0 \oplus Y \ar[r] \ar[d] & Y_0 \ar[d]^{\phi} \\
X \oplus Y \ar[r] & X \oplus Y \oplus Y''}$$
where the active morphism $\phi$ exhibits $Y_0$ as a sum
$Y'_0 \oplus Y''_0$, where $Y'_0 \rightarrow X$ is an active morphism
and $Y''_0 \simeq Y \oplus Y''$.
The inclusion $\calX \subseteq \calE^{\acti}_{/D}$ admits a left adjoint, and is therefore cofinal.

Let $D_0 \in \overline{\calD} = (K \star \{0\} ) \times_{ \calO^{\otimes} } \calK_{\calO}$ be the object corresponding to $\id_{X}$. The operation $\bigdot \mapsto \bigdot \oplus (Y \oplus Y'')$ determines
an equivalence of $\infty$-categories $\calD^{\acti}_{/D_0} \rightarrow \calX$.
It therefore suffices to show that the composite map
$$ (\calD^{\acti}_{/D_0})^{\triangleright} \rightarrow
\calX^{\triangleright} \subseteq (\calE^{\acti}_{/D})^{\triangleright} \stackrel{\delta}{\rightarrow} \calC^{\otimes}$$ is an operadic $q$-colimit diagram. This composition can be identified
with the map
$$ (\calD^{\acti}_{/D_0})^{\triangleright} \rightarrow \calC^{\otimes}
\stackrel{\oplus A(Y \oplus Y'')}{\rightarrow} \calC^{\otimes},$$
which is an operadic $q$-colimit diagram by virtue of our assumption that $\overline{p}$ satisfies
condition $(\ast)$.
\end{proof}

\subsection{Modules for the Associative $\infty$-Operad}\label{bimid}

The goal of this section is to study the theory of modules in the case where
$\calO^{\otimes}$ is the associative $\infty$-operad $\Ass$ of Example \ref{xe2}.
Our main results are as follows:
\begin{itemize}
\item[$(1)$] The associative $\infty$-operad $\Ass$ is coherent (Proposition \ref{asscoh}), so
there exists a sensible theory of modules over associative algebras.
\item[$(2)$] Let $p: \calC^{\otimes} \rightarrow \Ass$ be a fibration of $\infty$-operads, and let
$A \in \Alg_{\Ass}(\calC)$. Then there is an $\infty$-category $\Biod{A}{A}{\calC}$ of $A$-bimodule objects of $\calC$ (see Definition \ref{bidefn}) which is canonically equivalent
to $\Mod^{\Ass}_{A}(\calC)$ (Theorem \ref{rox}).
\item[$(3)$] Assume that $p$ determines an $\Ass$-monoidal structure on $\calC$ which is
compatible with the formation of geometric realizations. Then $\Mod_{\Ass}(\calC)$ inherits
the structure of an $\Ass$-monoidal $\infty$-category, where the monoidal structure
is given by the relative tensor product $\otimes_{A}$ defined in \S \monoidref{balpair}
(Theorem \ref{qwent}).
\end{itemize}

\begin{proposition}\label{asscoh}
The associative $\infty$-operad $\Ass$ is coherent.
\end{proposition}

\begin{remark}
We will prove Proposition \ref{asscoh} directly from the definition of a coherent $\infty$-operad
using a somewhat complicated combinatorial analysis. Later, we will give a simpler criterion
will enable us to verify the coherence of little cubes operads in general, using geometric arguments. Taking $n=1$, this leads to another proof of Proposition \ref{asscoh} which is somewhat different from the one given below.
\end{remark}

\begin{lemma}\label{combus}
Let $\Lin$ denote the category of finite linearly ordered sets (with morphisms given by
nonstrictly increasing maps). Suppose we are given a diagram $\tau$
$$ \xymatrix{ X_0 \ar[r]^{\theta} \ar[d]^{\alpha} & Y_0 \\
X & }$$
in $\Lin$, where $\alpha$ is injective. Let $\calC$ denote the full subcategory of
$\Lin_{\tau/}$ spanned by those commutative diagrams $\overline{\tau}$
$$ \xymatrix{ X_0 \ar[r]^{\theta} \ar[d]^{\alpha} & Y_0 \ar[d]^{\beta} \\
X \ar[r]^{\theta'} & Y}$$
in $\Lin$ for which $\beta$ is injective. Then the simplicial set $\Nerve(\calC)$ is weakly contractible.
\end{lemma}

\begin{proof}
We define an equivalence relation $\sim$ on $X$ as follows: $x \sim y$ if there exist
elements $a,b \in X_0$ such that $a \leq x \leq b$, $a \leq y \leq b$, and $\theta(a) = \theta(b)$.
Let $X' = X / \sim$. We observe that any commutative diagram $\overline{\tau}$ as above can be extended uniquely to a diagram
$$ \xymatrix{ X_0 \ar[r] \ar[d]^{\alpha} & \theta(X_0) \ar[d] \ar[r] & Y_0 \ar[d]^{\beta} \\
X \ar[r] & X' \ar[r] & Y}$$
We may therefore replace $X_0$ by $\theta(X_0)$ and $X$ by $X'$, thereby reducing to the case where the map $\theta$ is injective.

Let $\calC'$ denote the full subcategory of $\calC$ spanned by those diagrams $\tau$ for which
the map $Y_0 \coprod X \rightarrow Y$ is surjective. The inclusion $\calC' \subseteq \calC$
has a right adjoint, given by replacing $Y$ by the image of the map $Y_0 \coprod X$. 
It will therefore suffice to show that $\Nerve( \calC' )$ is weakly contractible. 

Let $Y_0^{\ast}$ denote the collection of all subsets $A \subseteq Y_0^{\ast}$ which
are closed downward in the following sense: if $y \in A$ and $y' \leq y$, then $y' \in A$.
We regard $Y_0^{\ast}$ as a partially ordered set with respect to inclusion. The
disjoint union $Y_0 \coprod Y_0^{\ast}$ also has the structure of a linear ordered
set, if we say that $y \leq A$ if and only if $y \in A$ for $y \in Y_0$, $A \in Y_0^{\ast}$.
Let $S = X - \alpha(X_0)$. Any commutative diagram $\overline{\tau}$ as above
gives rise to a nondecreasing map $f: S \rightarrow Y_0 \coprod Y_0^{\ast}$, which can
be defined as follows: for $s \in S$, we let $f(s) = y \in Y_0$ if $\theta'(s) = \beta(y)$ for
some $y \in Y_0$ (the element $y$ is then uniquely determined, since $\beta$ is injective),
and $f(s) = \{ y \in Y_0: \beta(y) < \theta'(s) \} \in Y_0^{\ast}$ otherwise. The function $f$
satisfies the following conditions:
\begin{itemize}
\item[$(a)$] Let $s \in S$ and $x \in X_0$ be such that $s < \alpha(x)$. Then
$f(s) \leq \theta(x)$.
\item[$(b)$] Let $s \in S$ and $x \in X_0$ be such that $s > \alpha(x)$. Then
$f(s) \geq \theta(x)$.
\end{itemize}

Conversely, given a nondecreasing map $f: S \rightarrow Y_0 \coprod Y_0^{\ast}$
satisfying conditions $(a)$ and $(b)$, we can uniquely construct a commutative diagram $\overline{\tau}$
$$ \xymatrix{ X_0 \ar[r]^{\theta} \ar[d]^{\alpha} & Y_0 \ar[d]^{\beta} \\
X \ar[r]^{\theta'} & Y}$$
corresponding to an element of $\calC'$. Namely, let $S' = f^{-1} Y_0^{\ast}$, and
set $Y = Y_0 \coprod S'$. We endow $Y$ with the unique linear ordering extending
that of $Y_0$ and $S'$, for which $y \leq s$ if and only if $y \in f(s)$
for $y \in Y_0$ and $s \in S'$. We take $\beta$ to be the evident inclusion, and
define $\theta'$ by the formula
$$ \theta'(s) = \begin{cases} (\beta \circ \alpha)(x) & \text{ if } s = \alpha(x) \\
f(s) & \text{ if } s \in f^{-1} Y_0 \\
s & \text{ if } s \in S'. \end{cases}$$
It is not difficult to see that this construction provides an inverse to the definition of $f$ from
$\overline{\tau}$. Consequently, we may identify objects of $\calC'$ with maps
$f: S \rightarrow Y_0 \coprod Y_0^{\ast}$ satisfying conditions $(a)$ and $(b)$.

Let $X_0^{\ast}$ denote the set of downward closed subsets of $X_0$. Since
the map $\theta$ is injective, it induces a surjection $Y_0^{\ast} \rightarrow X_0^{\ast}$.
Choose a section $\phi: X_0^{\ast} \rightarrow Y_0^{\ast}$ of this surjection.
We now choose a sequence of elements $t_i \in Y_0 \coprod Y_0^{\ast}$
using induction on $i \geq 1$. Assuming that $t_j$ has been chosen for
$j < i$, let $t_i$ denote an element of $Y_0 \coprod Y_0^{\ast} - \{ t_1, \ldots, t_{i-1} \}$ satisfying one of the following conditions:
\begin{itemize}
\item[$(1)$] The element $t_i$ belongs to $\theta(X_0) \subseteq Y_0$.
\item[$(2)$] There exists an element $A \in X_0^{\ast}$ such that
$t_i$ is a minimal element of
$\{ t \in Y_0 \coprod Y_0^{\ast}: (\phi(A) < t) \wedge ( \forall x \in X_0)[ (\theta(x) \leq t) \Rightarrow x \in A] \} $.
\item[$(3)$] There exists an element $A \in X_0^{\ast}$ such that
$t_i$ is a maximal element of
$\{ t \in Y_0 \coprod Y_0^{\ast}: (\phi(A) > t) \wedge ( \forall x \in X_0)[ (\theta(x) \geq t) \Rightarrow x \notin A] \}$.
\end{itemize}
Choosing the elements $t_i$ in this manner, we eventually obtain a sequence
$t_1, t_2, \ldots, t_n \in Y_0 \coprod Y_0^{\ast}$ such that 
$(Y_0 \coprod Y_0^{\ast} ) - \{ t_1, \ldots, t_n \} = \phi( X_0^{\ast})$. 
For $0 \leq i \leq n$, let $\calC'_{i}$ denote the full subcategory
of $\calC'$ spanned by those functions $f: S \rightarrow Y_0 \coprod Y_0^{\ast}$ whose
images are disjoint from $\{ t_1, \ldots, t_i \}$. The category
$\calC'_{0}$ coincides with $\calC'$, while the category
$\calC'_{n}$ consists of a single object: the map $f: S \rightarrow \phi(X_0^{\ast})$ given by $f(s) = \phi( \{ x \in X_0: \alpha(x) < s \})$. Consequently, to prove that $\Nerve(\calC')$ is weakly contractible,
it will suffice to show that each of the inclusions
$\Nerve( \calC'_{i-1} ) \subseteq \calC'_{i}$ is a weak homotopy equivalence.
We have several cases to consider, depending on whether $t_i$ satisfies condition
$(1)$, $(2)$, or $(3)$. We will assume that $t_i$ satisfies either condition $(1)$ or $(2)$;
the case where $t_i$ satisfies condition $(3)$ is treated using the same argument.

Suppose first that $t_i$ satisfies condition $(1)$, so that $t_i = \theta(x)$ for some
$x \in X_0$. Let $A_{-} = \{ y \in Y_0: y < \theta(x) \} \in Y_0^{\ast}$ and let $A_{+} = \{ y \in Y_0: y \leq \theta(x) \} \in Y_0^{\ast}$. For $f \in \calC'_{i-1}$, define $f': S \rightarrow Y_0 \coprod Y_0^{\ast}$
by the formula
$$ f'(s) = \begin{cases} f(s) & \text{ if } f(s) \neq t_i \\
A_{-} & \text{ if } f(s) = t_i \text{ and } s < \alpha(x) \\
A_{+} & \text{ if } f(s) = t_i \text{ and } s > \alpha(x). \end{cases}$$
The construction $f \mapsto f'$ is a right adjoint to the inclusion $\calC_{i} \subseteq \calC_{i-1}$,
so the inclusion $\Nerve( \calC_{i} ) \subseteq \Nerve( \calC_{i-1} )$ is a weak homotopy equivalence.

Suppose now that $t_i$ satisfies condition $(2)$ and that
$t_i \in Y_0$. Let $B = \{ y \in Y_0: y < t_i \} \in Y_0^{\ast}$. Given $f \in \calC'_{i-1}$, we define
$f': S \rightarrow Y_0 \coprod Y_0^{\ast}$ by the formula
$$ f'(s) = \begin{cases} f(s) & \text{ if } f(s) \neq t_i \\
B & \text{ if } f(s) = t_i. \end{cases}$$
We again observe that the construction $f \mapsto f'$ is a right adjoint to the inclusion
$\calC_i \subseteq \calC_{i-1}$, so that the inclusion $\Nerve( \calC_{i} ) \subseteq \Nerve( \calC_{i-1} )$ is a weak homotopy equivalence.

Suppose finally that $t_i$ satisfies condition $(2)$ and $t_i = B \in Y_0^{\ast}$, we proceed a little bit differently.
Note that the inequality $t_i > \phi(A)$ guarantees that $B$ is nonempty; let $y \in Y_0$ be
a maximal element of $B$. Given $f \in \calC'_{i-1}$, we define
$f': S \rightarrow Y_0 \coprod Y_0^{\ast}$ by the formula
$$ f'(s) = \begin{cases} f(s) & \text{ if } f(s) \neq t_i \\
y & \text{ if } f(s) = t_i. \end{cases}$$
The construction $f \mapsto f'$ is a left adjoint to the inclusion $\calC_{i} \subseteq \calC_{i-1}$,
so that $\Nerve( \calC_{i} ) \subseteq \Nerve( \calC_{i-1} )$ is a weak homotopy equivalence.
\end{proof}

\begin{proof}[Proof of Proposition \ref{asscoh}]
It is obvious that $\Ass$ is unital (it is the nerve of a category $\CatAss$ which contains a zero object). 
To prove that $\Ass$ is coherent, it will suffice to show that the evaluation map
$e_0: \calK_{\Ass} \rightarrow \Ass$ is a flat categorical fibration. Fix a $2$-simplex
$\Delta^2 \rightarrow \Ass$, corresponding to a commutative diagram $\tau$
$$ \xymatrix{ & Y_0 \ar[dr] & \\
X_0 \ar[ur] \ar[rr] & & Z_0 }$$
in $\CatAss$. We wish to show that the induced map $\calK_{\Ass} \times_{ \Ass} \Delta^2 \rightarrow \Delta^2$ is a flat categorical fibration. To prove this, choose a morphism
$\widetilde{X}_0 \rightarrow \widetilde{Z_0}$ in $\calK_{\Ass}$ lifting the map $X_0 \rightarrow Z_0$,
corresponding to a commutative diagram $\sigma_0$ 
$$ \xymatrix{ X_0 \ar[r] \ar[d]^{\alpha} & Z_0 \ar[d]^{\gamma} \\
X \ar[r] & Z }$$
in $\CatAss$ where the vertical maps are semi-inert. We wish to show that the $\infty$-category
$$\overline{\calC} = ( \calK_{\Ass} \times_{\Ass} \Delta^2)_{ \widetilde{X}_0 / \, / \widetilde{Z}_0} \times_{ \Delta^2 } \{1\}$$ is weakly contractible. Unwinding the definitions, we can identify
$\overline{\calC}$ with the nerve of the category $\calC$ of objects $Y \in \CatAss$ which
fit into a commutative diagram $\sigma$
$$ \xymatrix{ X_0 \ar[r] \ar[d]^{\alpha} & Y_0 \ar[r] \ar[d]^{\beta} & Z_0 \ar[d]^{\gamma} \\
X \ar[r] & Y \ar[r] & Z }$$
which extends $\sigma_0$ and $\tau$, where the map $\beta$ is semi-inert. 

Let us identify $X_0, Y_0, Z_0, X,$ and $Z$ with finite pointed sets. Let $\ast$ denote
the base point of $Z$. For every element $z \in Z$, let
$X_0^{z}, Y_0^{z}$, and $X^{z}$ denote the inverse images of $\{z\}$ in $X_0$, $Y_0$, and
$X$, respectively. We have a diagram of sets
$$ X^{z} \stackrel{\alpha^{z}}{\leftarrow} X_0^{z} \stackrel{\theta^{z}}{\rightarrow} Y_0^{z}.$$
If $z \neq \ast$, then these sets are linearly ordered and the maps $\alpha^{z}$ and $\theta^{z}$ are nonstrictly increasing; moreover, the condition that $\alpha$ is semi-inert guarantees that
$\alpha^{z}$ is injective. If $z = \ast$, then we can identify $\alpha^{z}$ and $\theta^{z}$ with morphisms in $\CatAss$. Producing the commutative diagram $\sigma$ as depicted above is equivalent
to giving a family of commutative diagrams $\sigma_{z}$
$$ \xymatrix{ X_0^{z} \ar[r]^{ \theta^{z}} \ar[d]^{\alpha^{z}} & Y_0^{z} \ar[d]^{\beta^{z}} \\
X^{z} \ar[r] & Y^{z}}$$
satisfying the following conditions:
\begin{itemize}
\item[$(a)$] If $z \neq \ast$, then $\sigma_{z}$ is a commutative diagram of linearly ordered
sets, and the map $\beta^{z}$ is injective.
\item[$(b)$] If $z = \ast$, then $\sigma_{z}$ is a commutative diagram in $\CatAss$, and
$\beta^{z}$ is inert.
\end{itemize}
For $z \in Z$, let $\calC_{z}$ denote the category whose objects are sets $Y_{z}$ fitting into
a commutative diagram $\sigma_{z}$ as above. Our reasoning shows that $\calC$ is equivalent
to the product $\prod_{z \in Z} \calC_{z}$. It will therefore suffice to show that
each $\calC_{z}$ has a weakly contractible nerve. If $z = \ast$, this follows from the observation
that $\calC_{z}$ has a final object (take $Y = \{ \ast \}$). If $z \neq \ast$, we invoke Lemma \ref{combus}.
\end{proof}

Let $\calC^{\otimes} \rightarrow \Ass$ be a fibration of $\infty$-operads, and let
$A \in \Alg_{\Ass}(\calC)$. Our next goal is to obtain a more explicit understanding of
the $\infty$-category $\Mod^{\Ass}_{A}(\calC)$. It follows from Proposition \ref{algass} that
$A$ is determined (up to equivalence) by the simplicial object given by the composition
$$ \Nerve(\cDelta)^{op} \rightarrow \Ass \stackrel{A}{\rightarrow} \calC^{\otimes}.$$
We will show that the theory of $A$-modules admits an analogous ``simplicial'' presentation.

\begin{notation}\label{urpas2}
Let $\star$ denote the {\em join} functor on linearly ordered sets: if $A$ and $B$
are linearly ordered sets, then $A \star B$ denotes the linearly ordered set
$A \coprod B$, where $a < b$ for all $a \in A$, $b \in B$. This construction
determines a functor $\cDelta \times \cDelta \rightarrow \cDelta$, where
$\cDelta$ denotes the category of finite simplices. We will denote this functor also
by $\star$: concretely, $\star$ is given by the formula
$$ [m] \star [n] = [m+n +1].$$
We observe that there are canonical maps $[m] \rightarrow [m] \star [n] \leftarrow [n]$, which
depend functorially on $m$ and $n$.

Let $\phi: \cDelta^{op} \rightarrow \CatAss$ be the functor described in Construction \ref{urpas}.
We define functors $\phi_{L}, \phi_{R}, \phi_{+}: \cDelta^{op} \times \cDelta^{op} \rightarrow \CatAss$
by the formulas
$$ \phi_{L}( [m], [n] ) = \phi( [m] ) \quad \phi_{+}( [m], [n] ) = \phi( [m] \star [n] )
\quad \phi_{R}( [m], [n] ) = \phi([n]).$$
We have natural transformations of functors 
$$ \phi_{L} \stackrel{ \beta_{L}}{\leftarrow} \phi_{+} \stackrel{\beta_{R}}{\rightarrow} \phi_{R},$$
which determine a diagram
$$\Phi: \Lambda^2_{0} \times \Nerve( \cDelta \times \cDelta)^{op} \rightarrow \Ass.$$
\end{notation}

\begin{definition}\label{bidefn}
Let $\calC^{\otimes} \rightarrow \Ass$ be a fibration of $\infty$-operads.
Let $\overline{M}:  \Lambda^2_0 \times \Nerve(\cDelta \times \cDelta)^{op} \rightarrow \calC^{\otimes}$
be a map of simplicial sets, which we think of as a triple
$(M, u: M \rightarrow A, v: M \rightarrow B)$ where $M$, $A$, and $B$
are functors $\Nerve( \cDelta \times \cDelta)^{op} \rightarrow \calC^{\otimes}$ and
$u$ and $v$ are natural transformations. We will say that $\overline{M}$ is a 
{\it prebimodule object of $\calC$} if the following conditions are satisfied:

\begin{itemize}
\item[$(1)$] The diagram
$$ \xymatrix{ & \calC^{\otimes} \ar[d] \\
\Lambda^2_0 \times ( \cDelta \times \cDelta)^{op} \ar[ur]^{\overline{M}} \ar[r]^-{\Phi} & \Ass }$$
is commutative, where $\Phi$ is as in Notation \ref{urpas2}.

\item[$(2)$] For every object $([m], [n]) \in \cDelta \times \cDelta$, the maps
$u_{[m], [n]}: M( [m], [n]) \rightarrow A( [m], [n])$ and
$v_{[m], [n]}: M( [m], [n]) \rightarrow B( [m], [n])$ are inert.

\item[$(3)$] Let $[0] \rightarrow [m]$ be the morphism in $\cDelta$ which carries
$0 \in [0]$ to $m \in [m]$, and let $[0] \rightarrow [n]$ be the morphism which
carries $0 \in [0]$ to $0 \in [n]$. Then the induced map
$M( [m], [n] ) \rightarrow M( [0], [0])$ is an inert morphism of $\calC^{\otimes}$.

\item[$(4)$] Given a morphism $([m], [n]) \rightarrow ([m'], [n'])$ in $\cDelta \times \cDelta$ such that 
the map $[m] \rightarrow [m']$ is convex, the induced map
$A( [m'], [n']) \rightarrow A([m], [n])$ is inert.

\item[$(5)$] Given a morphism $([m], [n]) \rightarrow ([m'], [n'])$ in $\cDelta \times \cDelta$ such that 
the map $[n] \rightarrow [n']$ is convex, the induced map
$B( [m'], [n']) \rightarrow B([m], [n])$ is inert.
\end{itemize}

We let $\Bimodd(\calC)$ denote the full
subcategory of $\Fun_{ \Ass}( \Lambda^2_0 \times \Nerve( \cDelta \times \cDelta)^{op}, \calC^{\otimes})$ spanned by the prebimodule objects of $\calC$.

Let
$\overline{M} = (M, u: M \rightarrow A, v: M \rightarrow B)$ be a bimodule
object of $\calC^{\otimes}$. We can regard $A$ as a functor from
$a: \Nerve( \cDelta)^{op} \rightarrow \Fun_{ \Ass}( \Nerve(\cDelta)^{op}, \calC^{\otimes})$.
Condition $(4)$ is equivalent to the assertions that $a$ carries every morphism
in $\cDelta^{op}$ to an equivalence, and that the essential image of $a$ is contained in
$\Alg( \calC) \subseteq \Fun_{ \Ass}( \Nerve(\cDelta)^{op}, \calC^{\otimes})$; here
$\Alg(\calC)$ denotes the full subcategory of $\Fun_{ \Ass}( \Nerve(\cDelta)^{op}, \calC^{\otimes})$ spanned by those functors which carry convex morphisms in $\cDelta^{op}$ to inert morphisms in
$\calC^{\otimes}$. Using the fact that $\Nerve(\cDelta)^{op}$ is weakly contractible, we conclude
that projection onto the first factor induces an equivalence $e_{L}$ from
$\Alg(\calD)$ to the full subcategory of $\Fun_{ \Ass}( \Nerve( \cDelta \times \cDelta)^{op}, \calC^{\otimes})$ spanned by those functors which satisfy $(4)$. Similarly,
projection onto the second factor induces an equivalence $e_{R}$ from
$\Alg(\calD)$ to the full subcategory of $\Fun_{ \Ass} ( \Nerve( \cDelta \times \cDelta)^{op}, \calC^{\otimes})$ spanned by those functors which satisfy $(5)$. Form a pullback diagram
$$ \xymatrix{ \Bimod(\calC) \ar[r] \ar[d] & \Bimodd(\calC) \ar[d]^{r} \\
\Alg(\calC) \times \Alg(\calC) \ar[r]^-{ (e_L, e_{R}) } & \Fun_{ \Ass}( \{ 1,2 \} \times \Nerve(\cDelta \times \cDelta)^{op}, \calC^{\otimes} ). }$$
Since $\Alg(\calC) \times \Alg(\calC)$ is an $\infty$-category and the map $r$ is a categorical fibration between $\infty$-categories, this map is a homotopy pullback diagram. It follows from the above arguments that the inclusion $\Bimod(\calC) \rightarrow \Bimodd(\calC)$ is an equivalence of $\infty$-categories. We will refer $\Bimod(\calC)$ as the {\it $\infty$-category of bimodule objects of $\calC$}.
Given a pair of algebra objects $A, B \in \Alg(\calC)$, we let
$\Biod{A}{B}{\calC}$ denote the fiber product $\Bimod(\calC) \times_{ \Alg(\calC) \times \Alg(\calC) }
\{ (A,B) \}$; we will refer to $\Biod{A}{B}{\calC}$ as the {\it $\infty$-category of $A$-$B$ bimodule
objects of $\calC$}.
\end{definition}

\begin{remark}
According to Proposition \ref{algass}, composition with the functor
$\phi: \cDelta^{op} \rightarrow \CatAss$ of Construction \ref{urpas} induces an equivalence
of $\infty$-categories $\theta: \Alg_{\Ass}(\calC) \rightarrow \Alg(\calC)$. Given a pair
of $\Ass$-algebra objects $A,B \in \Alg_{\Ass}(\calC)$, we let
$\Biod{A}{B}{\calC}$ denote the $\infty$-category $\Biod{\theta(A)}{\theta(B)}{\calC}$.
\end{remark}

\begin{construction}\label{kilus}
Let $\calK_{\Ass}^{\seg{1}} = \calK_{\Ass} \times_{ \Fun( \{0\}, \Ass} \{ \seg{1} \}$ denote the full subcategory of $\Ass^{ \seg{1}/ }$ spanned by the semi-inert morphisms $\seg{1} \rightarrow \seg{n}$ in
$\Ass$. Let $\phi: \cDelta^{op} \rightarrow \CatAss$ be the functor described in Construciton \ref{urpas}.
For every pair of objects $[m], [n] \in \cDelta$, there is a canonical map
$\seg{1} \rightarrow \phi( [m] \star [n] ) = \phi( [m + n + 1] ) = \seg{m+n+1}$, 
which carries the element $1 \in \nostar{1}$ to $m+1 \in \seg{m+n+1}$.
These maps fit into a commutative diagram
$$ \xymatrix{ \seg{1} \ar[d] & \seg{1} \ar[d] \ar[l] \ar[r] & \seg{1} \ar[d] \\
\phi_{L}( [m], [n]) & \phi_{+}( [m], [n]) \ar[l] \ar[r] & \phi_{R}( [m], [n] ), }$$
where the vertical maps are semi-inert. We can therefore factor the
map $\Phi$ of Notation \ref{urpas2} as a composition
$$ \Lambda^2_{0} \times \Nerve( \cDelta \times \cDelta)^{op} \stackrel{ \Phi'}{\rightarrow}
\calK_{\Ass}^{\seg{1}} \rightarrow \Ass.$$
\end{construction}

\begin{theorem}\label{rox}
Let $p: \calC^{\otimes} \rightarrow \Ass$ be a fibration of $\infty$-operads. Then
we have a homotopy pullback diagram of $\infty$-categories
$$ \xymatrix{ \Mod^{\Ass}(\calC) \ar[r] \ar[d] & \Bimod(\calC) \ar[d] \\
\Alg_{\Ass}(\calC) \ar[r] & \Alg(\calC) \times \Alg(\calC), }$$
where the upper horizontal map is given by composition with
the map $\Phi'$ of Construction \ref{kilus}. 
In particular, for every $\Ass$-algebra object $A \in \Alg_{\Ass}(\calC)$, we have an equivalence of
$\infty$-categories $\Mod^{\Ass}_{A}(\calC) \rightarrow \Biod{A}{A}{\calC}$.
\end{theorem}

\begin{proof}
We will employ the basic strategy of the proof of Proposition \ref{algass}. We begin by defining
a category $\calJ$ as follows:
\begin{itemize}
\item An object of $\calJ$ consists of either a semi-inert morphism $\alpha: \seg{1} \rightarrow \seg{n}$ in $\CatAss$ or a triple $( s, [m], [n])$ where $s$ is a vertex of $\Lambda^2_0$ and
$[m], [n] \in \cDelta$.
\item Morphisms in $\calJ$ are defined as follows: 
$$ \Hom_{\calJ}( \alpha, \beta) = \Hom_{ \CatAss_{\seg{1}/}}( \alpha, \beta) \quad \quad
\Hom_{\calJ}( (s, [m], [n]), \beta) = \emptyset $$
$$ \Hom_{\calJ}( \alpha, (s, [m], [n]) ) = \Hom_{ \CatAss_{/\seg{1}} }( \alpha,
\Phi'(s, [m], [n] )$$
$$\Hom_{\calJ}( (s, [m], [n]), (s', [m'], [n']))
= \begin{cases} \Hom_{ \cDelta^{op}}( [m], [m']) \times \Hom_{ \cDelta^{op}}( [n], [n']) & 
\text{ if } s=0 \text{ or } s = s' \\
\emptyset & \text{ otherwise.} \end{cases}$$
\end{itemize}
Let $\calJ_{0}$ denote the full subcategory of $\calJ$ spanned by the objects of the
form $(s, [m], [n])$. Let $\calJ_{1}$ denote the full subcategory of $\calJ$ spanned
by the null morphisms $\alpha: \seg{1} \rightarrow \seg{n}$ in $\CatAss$ together
with objects of the form $(s, [m], [n])$ where $s \neq 0$. Let
$\calJ_{01} = \calJ_0 \cap \calJ_1$.

We observe that there is a canonical retraction $r$ of $\calJ$ onto the full subcategory
$\calJ_{2} \subseteq \calJ$ spanned by the objects of the form $\alpha: \seg{1} \rightarrow \seg{n}$, which carries $(s, [m], [n])$ to $\Phi'(s, [m], [n])$. This retraction determines a map
$\Nerve(\calJ) \rightarrow \Ass$. Let $\calX$ denote the full subcategory of $\Fun_{\Ass}( \Nerve(\calJ), \calC^{\otimes})$ spanned by those functors $F$ with the following properties:
\begin{itemize}
\item[$(a)$] For every object $([m], [n]) \in \cDelta \times \cDelta$, the maps
$F(0, [m], [n]) \rightarrow F(1, [m], [n])$ and $F(0, [m], [n]) \rightarrow F(2, [m], [n])$ are inert.

\item[$(b)$] Let $[0] \rightarrow [m]$ be the morphism in $\cDelta$ which carries
$0 \in [0]$ to $m \in [m]$, and let $[0] \rightarrow [n]$ be the morphism in $\cDelta$ which
carries $0 \in [0]$ to $0 \in [n]$. Then the induced map
$F(0, [m], [n]) \rightarrow F(0, [0], [0])$ is inert.

\item[$(c)$] Given a morphism $([m], [n]) \rightarrow ([m'], [n'])$ in $\cDelta \times \cDelta$ such that 
the map $[m] \rightarrow [m']$ is convex, the induced map
$F(1, [m'], [n']) \rightarrow F(1, [m], [n])$ is inert.

\item[$(d)$] Given a morphism $([m], [n]) \rightarrow ([m'], [n'])$ in $\cDelta \times \cDelta$ such that 
the map $[n] \rightarrow [n']$ is convex, the induced map
$F( 2, [m'], [n']) \rightarrow F(2, [m], [n])$ is inert.

\item[$(e_{s})$] Let $(s, [m], [n]) \in \calJ_0$, and let $u: \Phi'( s, [m], [n]) \rightarrow (s, [m], [n])$ be the canonical map in $\calJ$. Then $F(u)$ is an equivalence in $\calC^{\otimes}$. 
\end{itemize}

We let $\calX_0$ denote the full subcategory of $\Fun_{\Ass}( \Nerve(\calJ_0), \calC^{\otimes})$ spanned by those functors which satisfy $(a)$, $(b)$, $(c)$, and $(d)$, 
$\calX_1$ the full subcategory of $\Fun_{\Ass}( \Nerve(\calJ_1), \calC^{\otimes})$ spanned by
those functors which satisfy $(c)$, $(d)$, $(e_1)$ and $(e_2)$, and
$\calX_{01}$ the full subcategory of $\Fun_{\Ass}( \Nerve(\calJ_{01}), \calC^{\otimes})$ spanned
by those functors which satisfy $(c)$ and $(d)$.

As explained in Definition \ref{bidefn}, we have a canonical equivalence
$g_{01}: \Alg(\calC) \times \Alg(\calC) \rightarrow \calX_{01}$. 
The retraction $r$ determines a map $g_1: \Alg_{\Ass}(\calC) \rightarrow \calX_{1}$.
Note that an object $F \in \Fun_{\Ass}( \Nerve(\calJ_1), \calC^{\otimes})$ satisfies
conditions $(e_1)$ and $(e_2)$ if and only if $F$ is a $p$-left Kan extension of
its restriction to $\Nerve(\calJ_1 \cap \calJ_2)$. Using Lemma \ref{preswope}, we see
that $F| \Nerve(\calJ_1 \cap \calJ_2)$ is equivalent to a functor which factors as a composition
$\Nerve(\calJ_1 \cap \calJ_2) \rightarrow \Ass \stackrel{F'}{\rightarrow} \calC^{\otimes}.$
Using Proposition \ref{algass}, we deduce that if $F$ satisfies $(e_1)$ and $(e_2)$, then
$F$ satisfies $(c)$ if and only if $F'$ is an $\Ass$-algebra object of $\calC^{\otimes}$; similarly,
$F$ satisfies $(d)$ if and only if $F' \in \Alg_{\Ass}(\calC)$. Using Proposition \toposref{lklk},
we deduce that the restriction functor $\calX_1 \rightarrow \calX_{12}$ is a trivial Kan fibration,
where $\calX_{12}$ is the essential image of the functor
$\Alg_{\Ass}(\calC) \rightarrow \Fun_{\Ass}( \Nerve(\calJ_1 \cap \calJ_2), \calC^{\otimes})$.
It follows that $g_1$ is an equivalence of $\infty$-categories.
We observe that $\calX_0 = \Bimodd(\calC)$, so we have a canonical equivalence
of $\infty$-categories $g_0: \Bimod(\calC) \rightarrow \calX_0$ (see Definition \ref{bidefn}).

Composition with $r$ yields a map $g': \MMod^{\Ass}(\calC) \rightarrow \calX$;
let $g$ be the composition of $g'$ with the equivalence $\Mod^{\Ass}(\calC)
\rightarrow \MMod^{\Ass}(\calC)$ provided by Remark \ref{user2}.
We claim that $g$ is an equivalence of $\infty$-categories. To prove this, it suffices to show
that $g'$ is a categorical equivalence. Note that a functor
$F \in \Fun_{\Ass}( \Nerve(\calJ), \calC^{\otimes})$ is a $p$-left Kan extension of
$F | \Nerve(\calJ_2)$ if and only if it satisfies conditions $(e_0)$, $(e_1)$, and $(e_2)$.
Moreover, any functor $F_0 \in \Fun_{\Ass}( \Nerve(\calJ_2), \calC^{\otimes})$
admits a $p$-left Kan extension $F \in \Fun_{\Ass}( \Nerve(\calJ), \calC^{\otimes})$:
we can take $F = F_0 \circ r$. It follows from Proposition \toposref{lklk} that
the restriction map $\calX \rightarrow \Fun_{\Ass}( \Nerve(\calJ_2), \calC^{\otimes})$
is a trivial Kan fibration onto the full subcategory of $\Fun_{\Ass}( \Nerve(\calJ_2), \calC^{\otimes})$
spanned by those functors $F_0$ such that $F_0 \circ r$ satisfies $(a)$, $(b)$, $(c)$, and $(d)$.
Since this restriction map is a left inverse to $g'$, it will suffice to show that this
full subcategory coincides with $\MMod^{\Ass}(\calC)$. In other words, we claim the following:
\begin{itemize}
\item[$(\ast)$] Let $F_0 \in \Fun_{\Ass}( \Nerve(\calJ_2), \calC^{\otimes})$. Then
$F_0 \in \MMod^{\Ass}(\calC)$ if and only if $F_0 \circ r$ satisfies $(a)$, $(b)$, $(c)$, and $(d)$.
\end{itemize}
The ``only if'' direction is clear. For the converse, let us suppose that
$F_0 \circ r$ satisfies $(a)$, $(b)$, $(c)$, and $(d)$. Let $u$ be a natural transformation
between semi-inert morphisms $\alpha: \seg{1} \rightarrow \seg{m}$ and
$\alpha': \seg{1} \rightarrow \seg{m'}$ in $\Ass$ such that the induced map
$\seg{m} \rightarrow \seg{m'}$ is inert; we wish to show that $F_0(u)$ is an inert morphism
in $\calC^{\otimes}$. The proof proceeds by induction on $m$. If $m' = m$, then $u$ is an equivalence and there is nothing to prove. We may therefore assume that $m' < m$. For $1 \leq i \leq m'$, let $\beta_i$ denote the composition of $\alpha'$ with the morphism $\Colp{i}: \seg{m'} \rightarrow \seg{1}$, and let $v_i: \alpha' \rightarrow \beta_i$ denote the induced morphism of $\calJ_2$.
By the inductive hypothesis, each $F_0(v_i)$ is an inert morphism in $\calC^{\otimes}$.
Consequently, to prove that $F_0(u)$ is inert, it will suffice to show that $F_0(v_i \circ u)$ is inert
for $1 \leq i \leq m'$. Replacing $u$ by $v_i \circ u$, we may reduce to the case where $m' = 1$.
There are three cases to consider. If $\alpha$ and $\alpha'$
are both null, then we can write $u= r( \overline{u} )$ where $\overline{u}: (1, [m], [0]) \rightarrow
( 1, [1], [0])$ is such that the underlying map $[1] \rightarrow [m]$ is a convex morphism
in $\cDelta$, so that $F_0(u) = (F \circ r)( \overline{u})$ is inert by virtue of $(c)$.
If neither $\alpha$ nor $\alpha'$ is null, then we can write $u = r( \overline{u})$ where
$\overline{u}: (0, [m-1], [0]) \rightarrow (0, [0], [0])$ is such that the underlying map
$[0] \rightarrow [m-1]$ carries $0 \in [0]$ to $m-1 \in [m-1]$, so that
$F_0(u) = (F_0 \circ r)( \overline{u})$ is inert by virtue of $(b)$. If $\alpha$ is not null but
$\alpha'$ is, then we can factor $u$ as a composition 
$$ \alpha \stackrel{u_0}{\rightarrow} (\alpha'': \seg{1} \rightarrow \seg{m-1}) \stackrel{u_1}{\rightarrow} \alpha',$$
where $\alpha''$ is null and the maps $\seg{m} \rightarrow \seg{m-1} \rightarrow \seg{m'}$ are inert.
The inductive hypothesis shows that $F_0( u_1)$ is inert, so it will suffice to show that
$F_0(u_0)$ is inert. For this, we observe that $u_0 = r( \overline{u})$, where
$\overline{u}$ is the map $(0, [m-1], [0]) \rightarrow (1, [m-1], [0])$ in $\calJ_0$.
Since $F_0$ satisfies $(a)$, we conclude that $F_0(u_0) = (F_0 \circ r)( \overline{u})$ is inert.

It follows that the maps $g$, $g_0$, $g_1$, and $g_{01}$ determine a categorical
equivalence of diagrams
$$ \xymatrix{ \Mod^{\Ass}(\calC) \ar[r] \ar[d] & \Bimod(\calC) \ar[d] & \calX \ar[r] \ar[d] & \calX_0 \ar[d] \\
\Alg_{\Ass}(\calC) \ar[r] & \Alg(\calC) \times \Alg(\calC) & \calX_1 \ar[r] & \calX_{01}. }$$
Consequently, it will suffice to show that the diagram on the right is a homotopy pullback square.
Since the diagram on the right consists of $\infty$-categories and categorical fibrations, this is equivalent to the assertion that the map $\theta: \calX \rightarrow \calX_0 \times_{ \calX_{01} } \calX_{1}$ is a categorical equivalence. We will prove that $\theta$ is a trivial Kan fibration.

Let $\calJ'$ denote the full subcategory of $\calJ$ spanned by $\calJ_0$ and $\calJ_1$, and let
$\calX'$ denote the full subcategory of $\Fun_{\Ass}( \Nerve(\calJ'), \calC^{\otimes})$ spanned by those 
functors $F_0$ which satisfy conditions $(a)$, $(b)$, $(c)$, $(d)$, $(e_1)$, and $(e_2)$.
We have an evident isomorphism $\calX' \simeq \calX_0 \times_{ \calX_{01}} \calX_1$.
To prove that $\theta$ is a trivial Kan fibration, it will suffice to prove the following:
\begin{itemize}
\item[$(i)$] A functor $F \in \Fun_{\Ass}( \Nerve(\calJ), \calC^{\otimes})$ such
that $F_0 = F | \Nerve(\calJ') \in \calX'$ is a $p$-right Kan extension of
$F_0$ if and only if it satisfies condition $(e_0)$.
\item[$(ii)$] For every object $F_0 \in \calX'$, there exists an object
$F \in \Fun_{\Ass}( \Nerve(\calJ), \calC^{\otimes})$ satisfying the equivalent conditions of $(i)$.

\end{itemize}

To prove these assertions, fix $F_0 \in \calX'$ and consider an object
$\alpha: \seg{1} \rightarrow \seg{n}$ such that $\alpha(1) = i \in \nostar{n}$ (so that
$\alpha \notin \calJ'$). Let $\alpha' \in \calJ_1$ denote the composition of
$\alpha$ with an inert morphism $\beta: \seg{n} \rightarrow \seg{n-1}$ such that
$\beta(i) = \ast$, so that $\beta$ induces a morphism $u: \alpha \rightarrow \alpha'$ in $\calJ'$.
Let $v: \alpha \rightarrow (0, [0], [0])$ be the unique
morphism in $\calJ$ such that the underlying map $\alpha \rightarrow \Phi'(0, [0], [0])$
induces the inert morphism $\Colp{i}: \seg{n} \rightarrow \seg{1} \simeq \Phi(0, [0], [0])$.
We will prove the following analogue of condition $(i)$:
\begin{itemize}
\item[$(i')$] Let $F \in \Fun_{ \Ass}( \Nerve(\calJ), \calC^{\otimes})$ be a functor such
that $F_0 = F | \Nerve(\calJ') \in \calX'$. Then $F$ is a $p$-right Kan extension of
$F_0$ if and only if, for each object $\alpha: \seg{1} \rightarrow \seg{n}$ of
$\calJ$ which does not belong to $\calJ'$, the functor $F$ carries the morphisms
$u: \alpha \rightarrow \alpha'$ and $v: \alpha \rightarrow (0, [0], [0])$ defined above
to inert morphisms in $\calC^{\otimes}$.
\end{itemize}

Assuming $(i')$ for the moment, let us prove $(i)$. Let $F \in \Fun_{\Ass}( \Nerve(\calJ, \calC^{\otimes})$
be such that $F_0 = F | \Nerve(\calJ') \in \calX'$; we will show that $F$ satisfies $(e_0)$ if and only if
it satisfies the condition described in $(i')$. We first prove the ``if'' direction. Let
$F$ satisfy $(e_0)$. Choose an object $\alpha: \seg{1} \rightarrow \seg{n}$ of
$\calJ$ which does not belong to $\calJ'$, and let $u: \alpha \rightarrow \alpha'$
and $v: \alpha \rightarrow (0, [0], [0])$ be defined as above. We have a commutative
diagram in $\calJ$
$$ \xymatrix{ \alpha \ar[r]^{u} \ar[d] & \alpha' \ar[d] \\
(0, [n-1], [0]) \ar[r]^{u'} & (1, [n-1], [0]). }$$
Since $F$ satisfies conditions $(e_0)$ and $(e_1)$, it carries the vertical morphisms
to equivalences in $\calC^{\otimes}$. Consequently, $F(u)$ is inert if and only if
$F_0(u')$ is inert, which follows from assumption $(a)$. To prove that $F(v)$ is inert, we write
$v$ as a composition
$$ \alpha \stackrel{v'}{\rightarrow} (0, [n-1], [0]) \stackrel{v''}{\rightarrow} (0, [0], [0]),$$
where $v''$ induces the map $[0] \rightarrow [n-1]$ in $\cDelta$ which carries
$0 \in [0]$ to $n-1 \in [n-1]$. The map $F(v')$ is an equivalence by virtue of
$(e_0)$, and the map $F(v'') = F_0(v'')$ is inert by virtue of $(b)$. It follows that
$F(v)$ is inert, as desired.

Conversely, suppose that $F$ satisfies the condition described in $(i')$; we will show that
$F$ satisfies $(e_0)$. Fix $m, n \geq 0$ and let $w: \Phi'(0, [m], [n]) \rightarrow (0, [m], [n])$ be the canonical map; we wish to show that $F(w)$ is an equivalence. Let
$j: (0, [m], [n]) \rightarrow (0, [0], [0])$ be the map in $\calJ$ given by the maps in $\cDelta$
$$[0] \rightarrow [m] \quad \quad [0] \rightarrow [n] $$
$$0 \in [0] \mapsto m \in [m] \quad \quad 0 \in [0] \mapsto 0 \in [n].$$
Let $j_{L}$ denote the map $(0, [m], [n]) \rightarrow (1, [m], [n])$ in $\calJ$, and
$j_{R}$ the map $(0, [m], [n]) \rightarrow (2, [m], [n])$ in $\calJ$. Since
$F_0$ satisfies $(a)$ and $(b)$, the maps $F_0(j)$, $F_0(j_{L})$, and $F_0(j_R)$ are inert.
Since $\calC^{\otimes}$ is an $\infty$-operad, to prove that $F(w)$ is an equivalence it will suffice to show that $F(j \circ w)$, $F(j_{L} \circ w)$, and $F( j_{R} \circ w)$ are inert. The first
follows from our assumption, since $j \circ w$ is the morphism $v$ associated $\alpha = \Phi'( 0, [m], [n])$. We will show that $F(j_{L} \circ w)$ is inert; the same argument will show that $F( j_{R} \circ w)$ is inert as well. The map $j_L \circ w$ can be factored as a composition
$$ \alpha \stackrel{u}{\rightarrow} \alpha' \stackrel{u'}{\rightarrow} 
\Phi'(1, [m], [n]) \stackrel{u''}{\rightarrow} (1, [m], [n]).$$
The map $F(u)$ is inert by assumption and $F(u'')$ is an equivalence by virtue of
assumption $(e_1)$. It therefore suffices to show that $F(u')$ is inert, which
follows from $(e_1)$ and $(c)$ because $u'$ can be identified with
$\Phi'( \overline{u}')$ for a suitable morphism $\overline{u'}: (1, [m+n], [0]) \rightarrow (1, [m], [0])$
which induces a convex map $[m] \rightarrow [m+n]$ in $\cDelta$. This completes the proof
that $(i') \Rightarrow (i)$.

We now prove $(i')$. Fix $F \in \Fun_{\Ass}( \Nerve(\calJ), \calC^{\otimes})$ such that
$F_0 = F | \Nerve(\calJ')$ belongs to $\calX'$. Let $\alpha: \seg{n} \rightarrow \seg{1}$ be an object of $\calJ$ which does not belong to $\calJ'$, and let $u: \alpha \rightarrow \alpha'$ and
$v: \alpha \rightarrow (0, [0], [0])$ be defined as above. We will show that $F$ is a $p$-right Kan extension of $F_0$ at $\alpha$ if and only if the morphisms $F(u)$ and $F(v)$ are inert.
Let $\calI$ denote the category $\calJ_{ \alpha / } \times_{\calJ} \calJ'$. We wish to prove that
the induced map $\Nerve(\calI)^{\triangleleft} \rightarrow \calC^{\otimes}$ is a $p$-limit diagram.
Let $\calI_0$ denote the full subcategory of $\calI$ spanned by those morphisms
$\alpha \rightarrow X$ in $\calJ$ which induce inert morphisms in $\FinSeg$.
It is not difficult to see that the inclusion $\calI_0 \subseteq \calI$ admits
a right adjoint, so that the induced map $\Nerve(\calI_0)^{\op} \rightarrow \Nerve(\calI)^{op}$ is cofinal.
Consequently, it will suffice to show that $F$ induces a $p$-limit diagram
$F': \Nerve(\calI_0)^{\triangleleft} \rightarrow \calC^{\otimes}$.

Let $\calI_1 \subseteq \calI_0$ be the full subcategory spanned by those morphisms 
$\alpha \rightarrow X$ in $\calJ$ satisfying the following additional conditions: 
\begin{itemize}
\item If $X$ has the form $(0, [m], [m'])$, then $m = m' = 0$.
\item If $X$ has the form $(1, [m], [m'])$ then $m > 0$.
\item If $X$ has the form $(2, [m], [m'])$, then $m' > 0$.
\end{itemize}
We claim that $F'| \Nerve(\calI_0)$ is a $p$-right Kan extension of
$F'| \Nerve(\calI_1)$. To prove this, consider an arbitrary object
$\overline{X}: \alpha \rightarrow X$ in $\calI_0$; we will show that $F'| \Nerve(\calI_0)$ is
a $p$-right Kan extension of $F' | \Nerve(\calI_1)$ at $\overline{X}$. If $\overline{X} \in \calI_1$
there is nothing to prove. If $X = (1, [0], [m'])$, then $( \calI_0)_{X/} \times_{ \calI_0} \calI_1$ is empty,
so it suffices to prove that $F'$ carries $\overline{X}$ to a $p$-final object of
$\calC^{\otimes}$. This follows from the assumption that $p$ is a fibration of $\infty$-operads, since
$p( \overline{X} ) \in \calC^{\otimes}_{\seg{0}}$. A similar argument applies if
$X$ has the form $(2, [m], [0])$ for some $m \geq 0$. We may therefore assume that
$X$ has the form $(0, [m], [m'])$, where $m + m' > 0$. Unwinding the definitions,
we deduce that $( \calI_0)_{\overline{X}/ } \times_{ \calI_0} \calI_1$ can be written
as a disjoint union $\calD_0 \coprod \calD_1 \coprod \calD_2$, where
$\calD_i$ is spanned by morphisms in $\calI_0$ of the form $\overline{X} \rightarrow \overline{X}'$
in $\calI_{0}$, where $\overline{X}'$ is a morphism in $\calJ$ of the form $\alpha \rightarrow (i, [m_0], [m'_0])$. We observe that $\calD_0$ contains a unique object (and only the identity morphism),
which we will denote by $f_0: \overline{X} \rightarrow \overline{X}_0$.
The category $\calD_1$ is empty if $m = 0$, and otherwise contains an initial object
$f_1: \overline{X} \rightarrow \overline{X}_1$, where $\overline{X}_1$ is the map
$\alpha \rightarrow (1, [m], [m'])$ determined by $\overline{X}$. Similarly,
$\calD_2$ is empty if $m' = 0$, and otherwise contains an initial object
$f_2: \overline{X} \rightarrow \overline{X}_2$. We will assume for simplicity
that $m, m' > 0$ (the cases where $m = 0$ or $m' = 0$ are handled similarly).
To prove that $F' | \Nerve(\calI_0)$ is a
$p$-right Kan extension of $F' | \Nerve(\calI_1)$ at $\overline{X}$, it suffices to show that
$F'$ exhibits $\overline{X}$ as a $p$-product of the objects $F'( \overline{X}_i )$ for
$i \in \{0,1,2\}$. Since $p$ is a fibration of $\infty$-operads, it will suffice to show that
each of the maps $F'( f_i)$ is an inert morphism in $\calC^{\otimes}$. This
follows from assumption $(b)$ when $i=0$ and from assumption $(a)$ when $i =1$ or $i=2$.

Invoking Lemma \toposref{kan0}, we are reduced to proving that the restriction
$F' | \Nerve( \calI_1)^{\triangleleft}$ is a $p$-limit diagram.
We observe that $\calI_1$ decomposes as a disjoint union of full subcategories
$\calI'_{1} \coprod \calI''_{1}$, where $\calI'_{1}$ is the subcategory consisting of the
single object $v: \alpha \rightarrow (0, [0], [0])$, and $\calI''_{1} = \calI_1 \times_{ \calJ} \calJ_{1}$.
The subcategory $\calI''_{1}$ has an initial object, given by the map $u: \alpha \rightarrow \alpha''$.
It follows that $F' | \Nerve( \calI_1)^{\triangleleft}$ is a limit diagram if and only if
the maps $F(u)$ and $F(v)$ exhibit $F(\alpha)$ as a $p$-product of the objects
$F( \alpha')$ and $F(0, [0], [0])$. Since $p$ is a fibration of $\infty$-operads, 
this is equivalent to the requirement that $F(u)$ and $F(v)$ are inert. This completes the proof
of $(i')$ (and therefore also the proof of $(i)$).

We now prove $(ii)$. Let $F_0 \in \calX'$; we wish to prove that $F_0$ admits a 
$p$-right Kan extension $F \in \Fun_{\Ass}( \Nerve(\calJ), \calC^{\otimes})$. 
In view of Lemma \toposref{kan2}, it will suffice to prove that for each
object $\alpha \in \calJ$ which does not belong to $\calJ'$, if we define
$\calI$ as above, then the diagram
$$ F'_0: \Nerve(\calI) \rightarrow \Nerve(\calJ') \stackrel{F_0}{\rightarrow} \calC^{\otimes}$$
can be extended to a $p$-limit diagram lifting the evident map
$G': \Nerve( \calI)^{\triangleleft} \rightarrow \Nerve(\calJ) \rightarrow \Ass.$
Since the inclusion $\Nerve( \calI_0) \subseteq \Nerve(\calI)$ is cofinal, it suffices
to show that $F'_0 | \Nerve(\calI_0)$ can be extended to a $p$-limit diagram
compatible with $G'$. Using Lemma \toposref{kan0}, we can reduce to showing
that $F'_0 | \Nerve(\calI_1)$ can be extended to a $p$-limit diagram compatible with $G'$.
Another cofinality argument shows that this is equivalent to proving the existence of a
$p$-limit diagram
$$ F_0(0, [0], [0]) \stackrel{ F(v)}{\leftarrow} F( \alpha) \stackrel{ F(u)}{\rightarrow} F_0(\alpha')$$
in $\calC^{\otimes}$ which is compatible with the projection to $\Ass$: the existence of this diagram follows easily from the assumption that $p$ is a fibration of $\infty$-operads.
\end{proof}

\begin{definition}
Let $s_{L}: \cDelta^{op} \rightarrow \cDelta^{op}$ be the functor
described by the formula $[n] \mapsto [n] \star [0] \simeq [n+1]$,
and let $\alpha_{L}: s_{L} \rightarrow \id$ be the natural transformation
induced by the inclusions $[n] \hookrightarrow [n] \star [0]$.
The natural transformation $\alpha$ and the functor
$\phi: \cDelta^{op} \rightarrow \CatAss$ determine a map of simplicial sets
$$ \Nerve( \cDelta)^{op} \times \Delta^1 \rightarrow \Ass.$$
Let $\calC^{\otimes} \rightarrow \Ass$ be a fibration of $\infty$-operads.
We let $$\Mod^{L}(\calC) \subseteq \Fun_{ \Ass}( \Nerve(\cDelta)^{op} \times \Delta^1, \calC^{\otimes})$$
be the full subcategory spanned
by those natural transformations $M \rightarrow A$ of functors
$M, A: \Nerve(\cDelta)^{op} \rightarrow \calC^{\otimes}$ which satisfy the following
conditions:
\begin{itemize}
\item[$(1)$] The functor $A$ carries convex morphisms in $\Nerve(\cDelta)^{op}$ to
inert morphisms in $\calC^{\otimes}$ (so that $A \in \Alg(\calC)$).
\item[$(2)$] For every object $[n] \in \cDelta^{op}$, the map
$M( [n] ) \rightarrow A( [n])$ is inert.
\item[$(3)$] Let $i: [0] \rightarrow [n]$ be the morphism in $\cDelta$ given by $i(0) = n$.
Then the induced map $M( [n]) \rightarrow M([0])$ is inert.
\end{itemize}
We will refer to $\Mod^{L}(\calC)$ as the {\it $\infty$-category of left module objects
of $\calC$}. Similarly, if we let $s_{R}$ denote the shift functor $[n] \mapsto [0] \star [n]$
and $\alpha_{R}: s_{R} \rightarrow \id$ the induced natural transformation, we
define $\Mod^{R}(\calC)$ to be the full subcategory of
$\Fun_{ \Ass}( \Nerve(\cDelta)^{op} \times \Delta^1, \calC^{\otimes})$ spanned by those functors
which satisfy conditions $(1)$ and $(2)$, together with the following analogue of $(3)$:
\begin{itemize}
\item[$(3')$] Let $i: [0] \rightarrow [n]$ be the morphism in $\cDelta$ given by $i(0) = 0$.
Then the induced map $M( [n]) \rightarrow M([0])$ is inert.
\end{itemize}
We will refer to $\Mod^{R}(\calC)$ as the {\it $\infty$-category of right $R$-module objects of $\calC$.}
\end{definition}

\begin{remark}
Let $\calC^{\otimes} \rightarrow \Ass$ be a fibration of $\infty$-operads. There are evident categorical fibrations $\Mod^{L}(\calC) \rightarrow \Alg(\calC) \leftarrow \Mod^{R}(\calC)$.
If $A \in \Alg(\calC)$, we let $\Mod^{L}_{A}(\calC)$ and $\Mod^{R}_{A}(\calC)$ denote the fibers of
these maps over $A$. We will refer to $\Mod^{L}_{A}(\calC)$ and
$\Mod^{R}_{A}(\calC)$ as {\it the $\infty$-category of left $A$-module objects of $\calC$
and the $\infty$-category of right $A$-module objects of $\calC$}, respectively.
\end{remark}

\begin{remark}
Suppose that $\calC^{\otimes} \rightarrow \Ass$ is a coCartesian fibration of $\infty$-operads.
Then $\Mod^{L}(\calC)$ and $\Mod^{R}(\calC)$ can be identified with the $\infty$-categories
of left and right module objects of the associated monoidal $\infty$-category
$\calC^{\otimes} \times_{\Ass} \Nerve( \cDelta)^{op}$, as defined in Remark \monoidref{capsul}.
\end{remark}

\begin{remark}
Let $\calC^{\otimes} \rightarrow \Ass$ be a fibration of $\infty$-operads. The inclusion
$$ \Nerve( \cDelta)^{op} \times \Delta^1
\simeq \Nerve( \cDelta \times \{ [0] \})^{op} \times \Delta^1 \subseteq \Nerve( \cDelta \times \cDelta)^{op} \times \Lambda^2_0$$
induces a forgetful functor $\Bimod(\calC) \rightarrow \Mod^{L}(\calC)$. Similarly,
we have a forgetful functor $\Bimod(\calC) \rightarrow \Mod^{R}(\calC)$. Given a pair of
object $A, B \in \Alg(\calC)$, these forgetful functors induce maps
$$ \Mod_{A}^{L}(\calC) \leftarrow \Biod{A}{B}{\calC} \rightarrow \Mod_{B}^{R}(\calC).$$
\end{remark}

\begin{theorem}\label{qwent}
Let $p: \calC^{\otimes} \rightarrow \Ass$ be a coCartesian fibration of $\infty$-operads which is
compatible with $\Nerve(\cDelta)^{op}$-indexed colimits, and let $A \in \Alg_{\Ass}(\calC)$.
Then:
\begin{itemize}
\item[$(1)$] The induced map $q: \Mod_{A}^{\Ass}(\calC)^{\otimes} \rightarrow \Ass$ is a coCartesian fibration of $\infty$-operads.
\item[$(2)$] The tensor product functor $\otimes: \Mod_{A}^{\Ass}(\calC) \times \Mod_{A}^{\Ass}(\calC) 
\simeq \Mod_{A}^{\Ass}(\calC)^{\otimes}_{\seg{2}} \rightarrow
\Mod_{A}^{\Ass}(\calC)$ fits into a diagram
$$ \xymatrix{ \Mod_{A}^{\Ass}(\calC) \times \Mod_{A}^{\Ass}(\calC) \ar[r]^-{\otimes} \ar[d] & \Mod_{A}^{\Ass}(\calC) \ar[dd] \\
\Biod{A}{A}{\calC} \times \Biod{A}{A}{\calC} \ar[d] & \\
\Mod^{R}_{A}(\calC) \times \Mod^{L}_{A}(\calC) \ar[r]^-{ \otimes_{A} } & \calC }$$
which commutes up to canonical homotopy, where the relative tensor product functor $\otimes_{A}: \Mod^{R}_{A}(\calC) \times \Mod^{L}_{A}(\calC) \rightarrow \calC$ is as defined in \S \monoidref{balpair}.
\end{itemize}
\end{theorem}

The proof of Theorem \ref{qwent} will use the following criterion for detecting the existence of $\Ass$-monoidal structures:

\begin{proposition}\label{breadtime}
Let $q: \calC^{\otimes} \rightarrow \Ass$ be a fibration of $\infty$-operads. Then
$q$ is a coCartesian fibration if and only if the following conditions are satisfied:
\begin{itemize}
\item[$(1)$] The fibration $q$ has units (in the sense of Definition \ref{cuspe}). 
\item[$(2)$] Let $\alpha: \seg{2} \rightarrow \seg{1}$ be an active morphism in $\Ass$,
and let $X \in \calC^{\otimes}_{\seg{2}}$. Then there exists a map
$\overline{\alpha}: X \rightarrow X'$ with $q( \overline{\alpha}) = \alpha$, which determines
an operadic $q$-colimit diagram $\Delta^1 \rightarrow \calC^{\otimes}$.
\end{itemize}
\end{proposition}

\begin{proof}
If $q$ is a coCartesian fibration, then conditions $(1)$ and $(2)$ follow
immediately from Proposition \ref{justica}. Conversely, suppose that
$(1)$ and $(2)$ are satisfied; we wish to prove that $q$ is a coCartesian fibration. In view of Proposition \ref{justica}, it will suffice to prove that the following stronger
version of $(2)$ is satisfied:
\begin{itemize}
\item[$(2')$] Let $\alpha: \seg{n} \rightarrow \seg{1}$ be an active morphism in
$\Ass$, and let $X \in \calC^{\otimes}_{\seg{n}}$. Then there a map $\overline{\alpha}: X \rightarrow X'$ lifting $\alpha$ which is given by an operadic $q$-colimit diagram.
\end{itemize}
We verify $(2')$ using induction on $n$. If $n = 0$, the desired result follows from $(1)$.
If $n=1$, then $\alpha$ is the identity and we can choose $\overline{\alpha}$ to be a degenerate edge. Assume that $n \geq 2$. The map $\alpha$ determines a linear ordering on the set $\nostar{n}$.
Choose any partition of this set into nonempty disjoint subsets $\nostar{n}_{-}$ and
$\nostar{n}_{+}$, where $\nostar{n}_{-}$ is closed downwards and $\nostar{n}_{+}$ is closed
upwards. This decomposition determines a factorization of
$\alpha$ as a composition
$$ \seg{n} \stackrel{ \alpha'}{\rightarrow} \seg{2} \stackrel{ \alpha''}{\rightarrow} \seg{1}.$$
Let $n_{-}$ and $n_{+}$ be the cardinalities of $\nostar{n}_{-}$ and $\nostar{n}_{+}$, respectively, so that our decomposition induces inert morphisms $\seg{n} \rightarrow \seg{n_{+}}$ and
$\seg{n} \rightarrow \seg{n_{-}}$. Since $q$ is a fibration of $\infty$-operads, we can lift these
maps to inert morphisms $X \rightarrow X_{-}$ and $X \rightarrow X_{+}$ in $\calC^{\otimes}$.
Using the inductive hypothesis, we can choose morphisms $X_{-} \rightarrow X'_{-}$ and
$X_{+} \rightarrow X'_{-}$ lying over the maps $\seg{n_{-}} \rightarrow \seg{1} \leftarrow \seg{n_{+}}$
in $\Ass$, which are classified by operadic $q$-colimit diagrams $\Delta^1 \rightarrow \calC^{\otimes}$.
Let $\overline{\alpha}'$ denote the induced morphism
$X \simeq X_{-} \oplus X_{+} \rightarrow X'_{-} \oplus X'_{+}.$
Using $(2)$, we can choose a morphism
$\overline{\alpha}'': X'_{-} \oplus X'_{+} \rightarrow X'$ lifting $\alpha''$, which is classified
by an operadic $q$-colimit diagram. To complete the proof, it will suffice to show that
$\overline{\alpha}'' \circ \overline{\alpha}'$ is classified by an operadic $q$-colimit diagram.
Choose an object $Y \in \calC^{\otimes}$; we must show that the composite map
$$ X \oplus Y \simeq X_{-} \oplus X_{+} \oplus Y
\stackrel{\beta}{\rightarrow} X'_{-} \oplus X'_{+} \oplus Y \stackrel{\gamma}{\rightarrow} X' \oplus Y$$
is classified by a weak operadic $q$-colimit diagram. To prove this, it suffices to show that
$\gamma \circ \beta$ is $q$-coCartesian. In view of Proposition \toposref{protohermes}, it suffices to show that $\gamma$ and $\beta$ are $q$-coCartesian, which follows from
Proposition \ref{copan}.
\end{proof}

\begin{variant}\label{subless2}
Proposition \ref{breadtime} remains valid (with the same proof) if we replace the
associative $\infty$-operad $\Ass$ by the commutative $\infty$-operad $\CommOp$.
\end{variant}

\begin{variant}\label{subless}
The proof of Proposition \ref{breadtime} also establishes the following:
\begin{itemize}
\item[$(\ast)$] Suppose given a commutative diagram of $\infty$-operads
$$ \xymatrix{ \calC^{\otimes} \ar[rr]^{f} \ar[dr]^{p} & & \calD^{\otimes} \ar[dl]^{q} \\
& \Ass & }$$
where $p$ and $q$ are coCartesian fibrations. Then $f$ is an $\Ass$-monoidal functor
(that is, $f$ carries $p$-coCartesian morphisms to $q$-coCartesian morphisms) if and only
if the following conditions are satisfied:
\begin{itemize}
\item[$(i)$] The functor $f$ preserves unit objects (that is, if 
$\eta: C_0 \rightarrow C$ is a morphism in $\calC^{\otimes}$ which exhibits
$C \in \calC$ as a unit object of $\calC^{\otimes}$, then $f(\eta)$ exhibits $f(C)$ as a unit object
of $\calD^{\otimes}$). 
\item[$(ii)$] The functor $f$ preserves tensor products. That is, for every pair of objects
$C, C' \in \calC$, the canonical map $f(C) \otimes f(C') \rightarrow f( C \otimes C')$ is an equivalence.
(Equivalently, if $\eta$ is a $p$-coCartesian morphism in $\calC^{\otimes}$ lying over an
active map $\seg{2} \rightarrow \seg{1}$ in $\Ass$, then $f(\eta)$ is $q$-coCartesian.)
\end{itemize}
\end{itemize}
\end{variant}

\begin{proof}[Proof of Theorem \ref{qwent}]
Example \ref{swagg} shows that $q$ has units. In view of Proposition \ref{breadtime},
assertion $(1)$ is equivalent to the following:
\begin{itemize}
\item[$(1')$] Let $\alpha: \seg{2} \rightarrow \seg{1}$ be an active morphism in $\Ass$,
and let $X \in \Mod^{\Ass}_{A}(\calC)^{\otimes}_{\seg{2}}$. Then there exists a map
$\overline{\alpha}: X \rightarrow X'$ with $q( \overline{\alpha}) = \alpha$, which determines
an operadic $q$-colimit diagram $\Delta^1 \rightarrow \Mod^{\Ass}_{A}(\calC)^{\otimes}$.
\end{itemize}
We will deduce $(1')$ using the criterion of Theorem \ref{colmodd}. 
Let $\overline{\calD} = \calK_{\Ass} \times_{ \Ass} \Delta^1$, where
$\Delta^1 \rightarrow \Ass$ is the map determined by $\alpha$. Let
$\calD = \overline{\calD} \times_{ \Delta^1} \{0\}$, and let $Y \in \overline{\calD}$
be the object corresponding to the identity map from $\seg{1}$ to itself.
The map $X$ determines an diagram
$F: \calD' \rightarrow \calC^{\otimes}$, where $\calD' = \overline{\calD}_{/Y} \times_{ \overline{\calD}} \calD$.
We must prove that $F$ can be extended to an operadic $p$-colimit diagram
$\overline{F}: {\calD'}^{\triangleright} \rightarrow \calC^{\otimes}$ which is compatible with
the evident map
$$g: { \calD'}^{\triangleright} \rightarrow \overline{\calD}_{/Y}^{\triangleright} \rightarrow
\overline{\calD} \rightarrow \Ass.$$

We can identify objects of $\calD'$ with the nerve of the category of diagrams
$$ \xymatrix{ \seg{2} \ar[r]^{v} \ar[d]^{\alpha} & \seg{m} \ar[d]^{u} \\
\seg{1} \ar[r]^{\id} & \seg{1} }$$
in $\Ass$, where $v$ is semi-inert. Let $\calD'_0$ denote the full subcategory
of $\calD'$ spanned by those diagrams where $u$ is active. The inclusion $\calD'_0 \subseteq \calD'$ admits a left adjoint, and is therefore cofinal. Consequently, it will suffice to show that
$F| \calD'_0$ can be extended to an operadic $p$-colimit diagram ${ \calD'}^{\triangleright}_0 \rightarrow \calC^{\otimes}$ compatible with $g$.

For every object of $\calD'_0$ (corresponding to a diagram as above), the map $u$ determines
a linear ordering of $\nostar{m}$, and we can identify $v$ with an order-preserving injection $\nostar{2} \rightarrow \nostar{m}$. Let $\calD'_1$ denote the full subcategory of $\calD'_0$ spanned by those
objects for which $v(1)$ is a minimal element of $\nostar{m}$ and $v(2)$ is a maximal element of
$\nostar{m}$. The inclusion $\calD'_0 \subseteq \calD'_1$ admits a left adjoint, and is therefore cofinal.
We are therefore reduced to proving that $F | \calD'_1$ can be extended to an operadic
$p$-colimit diagram ${\calD'}_{1}^{\triangleright} \rightarrow \calC^{\otimes}$ which is compatible with $g$.

The functor $\phi$ of Construction \ref{urpas} determines an equivalence of $\infty$-categories
$\Nerve( \cDelta)^{op} \rightarrow \calD'_{1}$. Let $F': \Nerve( \cDelta)^{op} \rightarrow \calC^{\otimes}$ be the functor obtained by composing this equivalence with $F$. It will now suffice to show that
$F'$ can be extended to an operadic $p$-colimit diagram $( \Nerve(\cDelta)^{op})^{\triangleright} \rightarrow \calC^{\otimes}$ compatible with $g$. Since $p$ is a coCartesian fibration compatible
with $\Nerve(\cDelta)^{op}$-indexed colimits, the existence of this extension follows from
Proposition \ref{slavewell}. This completes the proof of $(1)$. Moreover, it shows that
the comosition $\Mod^{\Ass}_{A}(\calC) \times \Mod^{\Ass}_{A}(\calC) \rightarrow
\Mod_{A}^{\Ass}(\calC) \rightarrow \calC$ can be computed as the composition
$$ \Mod_{A}^{\Ass}(\calC) \times \Mod_{A}^{\Ass}(\calC)
\stackrel{\theta}{\simeq} \Mod_{A}^{\Ass}(\calC)^{\otimes}_{\seg{2}}
\stackrel{\theta'}{\rightarrow} \Fun( \Nerve(\cDelta)^{op}, \calC^{\otimes}_{\seg{2}})
\stackrel{\theta''}{\rightarrow} \Fun( \Nerve(\cDelta)^{op}, \calC)
\stackrel{| |}{\rightarrow} \calC.$$
To prove $(2)$, it suffices to observe that $\theta'' \circ \theta' \circ \theta$ is equivalent to the composition
$$ \Mod_{A}^{\Ass}(\calC) \times \Mod_{A}^{\Ass}(\calC)
\rightarrow \Mod_{A}^{R}(\calC) \times \Mod_{A}^{R}(\calC)
\stackrel{\psi}{\rightarrow} \Fun( \Nerve(\cDelta)^{op}, \calC),$$
where $\psi$ is the functor $(M,N) \mapsto \Baar_{A}(M,N)_{\bigdot}$ of
Definition \monoidref{barcon}.
\end{proof}

\subsection{Modules for the Commutative $\infty$-Operad}\label{modcom}

In this section, we will study the theory of modules over the commutative $\infty$-operad
$\CommOp = \Nerve(\FinSeg)$ of Example \ref{xe1}. Our main results are as follows:

\begin{itemize}
\item[$(1)$] The $\infty$-operad $\CommOp$ is coherent (Proposition \ref{commcoh}), so that there is a well-behaved theory of module objects in symmetric monoidal $\infty$-categories.
\item[$(2)$] Let $\calC^{\otimes}$ be an $\infty$-operad, and let $A \in \CAlg(\calC)$. Then
the $\infty$-category $\Mod_{A}(\calC) = \Mod^{\CommOp}_{A}(\calC)$ is equivalent to the $\infty$-category
$\Mod^{L}_{A}(\calC)$ of left $A$-module objects of $\calC$ (Proposition \ref{socrime}).
\item[$(3)$] Suppose that $\calC^{\otimes}$ is a symmetric monoidal $\infty$-category, and that
the symmetric monoidal structure on $\calC^{\otimes}$ is compatible with the formation of geometric realizations of simplicial objects. Then, for every commutative algebra object $A \in \CAlg(\calC)$, the $\infty$-category $\Mod_{A}(\calC)$ inherits the structure of a symmetric monoidal $\infty$-category.
The tensor product on $\Mod_{A}(\calC)$ is given by the relative tensor product $\otimes_{A}$
defined in \S \monoidref{balpair} (this follows from Proposition \ref{stancer} and Theorem \ref{qwent}).
\item[$(4)$] Let $\calC^{\otimes}$ be as in $(3)$, and let $f: A \rightarrow B$ be a map of commutative algebra objects of $\calC^{\otimes}$. Then the forgetful functor $\Mod_{B}(\calC) \rightarrow
\Mod_{A}(\calC)$ has a left adjoint $M \mapsto M \otimes_{A} B$, which is a symmetric monoidal functor from $\Mod_{A}(\calC)$ to $\Mod_{B}(\calC)$ (Theorem \ref{podus}).
\end{itemize}

We begin with the coherence of $\CommOp$.

\begin{proposition}\label{commcoh}
The commutative $\infty$-operad $\CommOp$ is coherent.
\end{proposition}

\begin{proof}
Let $\overline{e}_0: \Fun( \Delta^1, \CommOp) \rightarrow \CommOp$ be the map given by evaluation at $0$. Since $\CommOp$ is the nerve of the category of pointed
finite sets, it admits pushouts. Applying Lemma \toposref{charpull}, we deduce that
$\overline{e}_0$ is a coCartesian fibration. Moreover, a morphism in $\Fun( \Delta^1, \CommOp)$ is
$p$-coCartesian if and only if it corresponds to a pushout diagram
$$ \xymatrix{ \seg{m} \ar[r] \ar[d]^{f} & \seg{n} \ar[d]^{g} \\
\seg{m_0} \ar[r] & \seg{n_0} }$$
of pointed finite sets. Note that if $f$ is semi-inert, then $g$ is also semi-inert. In other words,
if $\alpha: f \rightarrow g$ is an $\overline{e}_0$-coCartesian morphism in
$\Fun( \Delta^1, \CommOp)$ such that $f \in \calK_{\CommOp}$, then $g \in \calK_{\CommOp}$.
It follows that $\overline{e}_0$ restricts to a coCartesian fibration $e_0: \calK_{\CommOp} \rightarrow \CommOp$. Using Example \bicatref{cartflat} and Proposition \toposref{funkyfibcatfib}, we conclude that
$e_0$ is a flat categorical fibration. Since $\CommOp$ is pointed (it contains a zero object $\seg{0}$),
it is a unital $\infty$-operad, and is therefore coherent as desired.
\end{proof}

We now analyze the $\infty$-category of modules over a commutative algebra object
$A$ of an $\infty$-operad $\calC^{\otimes}$. It turns out that the $\infty$-category of $A$-modules
depends only on the underlying associative algebra of $A$.

\begin{remark}\label{swilk}
Let $f: {\calO'}^{\otimes} \rightarrow \calO^{\otimes}$ be a map of coherent $\infty$-operads, and let $\calC^{\otimes} \rightarrow \calO^{\otimes}$ be a fibration of $\infty$-operads.
We let $\Mod^{\calO'}(\calC)^{\otimes}$ denote the $\infty$-category
$\Mod^{\calO'}(\calC')^{\otimes}$, where ${\calC'}^{\otimes} = \calC^{\otimes} \times_{ \calO^{\otimes}}
{\calO'}^{\otimes}$.

The map $f$ induces a map $F: \calK_{\calO'} \rightarrow \calK_{ \calO} \times_{ \calO^{\otimes} } {\calO'}^{\otimes}$. Composition with $F$ determines a functor
$$\Mod^{\calO}(\calC)^{\otimes} \times_{ {\calO^{\otimes}}} {\calO'}^{\otimes}
\rightarrow \Mod^{\calO'}(\calC)^{\otimes}.$$
\end{remark}

\begin{notation}\label{crime}
Let $\calC^{\otimes}$ be an $\infty$-operad. We let
$\Alg(\calC)$ and $\Mod^{L}(\calC)$ denote the $\infty$-categories of algebra
objects and left module objects of the associated fibration of $\infty$-operads
$\calC^{\otimes} \times_{ \CommOp} \Ass \rightarrow \Ass$. The construction
of Remark \ref{swilk} determines a functor
$\Mod^{\CommOp}(\calC) \rightarrow \Mod^{\Ass}(\calC)$. We let
$\psi$ denote the composite functor
$$ \Mod^{\CommOp}(\calC) \rightarrow \Mod^{\Ass}(\calC) \rightarrow
\Bimod(\calC) \rightarrow \Mod^{L}(\calC).$$
\end{notation}

\begin{proposition}\label{socrime}
Let $p: \calC^{\otimes} \rightarrow \Nerve(\FinSeg)$ be an $\infty$-operad. Then the functor
$\psi$ of Notation \ref{crime} fits into a homotopy pullback diagram
of $\infty$-categories
$$ \xymatrix{ \Mod^{\CommOp}(\calC) \ar[r]^{\psi} \ar[d] & \Mod^{L}(\calC) \ar[d] \\
\CAlg(\calC) \ar[r] & \Alg(\calC). }$$
\end{proposition}

The proof of Proposition \ref{socrime} is somewhat complicated, and will be given at the end of this section.

\begin{remark}
Let $\calC^{\otimes}$ be an $\infty$-operad, and let $A$ be a commutative
algebra object of $\calC^{\otimes}$. Proposition \ref{socrime} implies in particular that
the induced map $\Mod^{\CommOp}_{A} \rightarrow \Mod^{L}_{A}(\calC)$ is an equivalence
of $\infty$-categories. In other words, if $A \in \Alg(\calC)$ can be promoted to a commutative
algebra object of $\calC$, then the $\infty$-category $\Mod^{L}_{A}(\calC)$ of left $A$-modules
itself admits the structure of an $\infty$-operad (which is often a symmetric monoidal $\infty$-category, as we will see in a moment.
\end{remark}

Our next goal is to study the $\infty$-operad structure on $\Mod_{A}^{\CommOp}(\calC)$, where
$A$ is a commutative algebra object of a symmetric monoidal $\infty$-category $\calC$. 
The following result reduces the analysis of the tensor product to the associative
setting: 

\begin{proposition}\label{stancer}
Let $p: \calC^{\otimes} \rightarrow \Nerve(\FinSeg)$ be a symmetric monoidal category
compatible with $\Nerve(\cDelta)^{op}$-indexed colimits, and let $A \in \CAlg(\calC)$.
Then:
\begin{itemize}
\item[$(1)$] The $\infty$-operad $\Mod^{\CommOp}_{A}(\calC)^{\otimes}$ is a
symmetric monoidal $\infty$-category.
\item[$(2)$] The functor
$$ \psi: \Mod^{\CommOp}_{A}(\calC)^{\otimes} \times_{ \CommOp} \Ass \rightarrow
\Mod^{\Ass}_{A}(\calC)^{\otimes}$$
is an $\Ass$-monoidal functor. Here we identify $A \in \CAlg(\calC)$ with its image
in $\Alg_{\Ass}(\calC)$.
\end{itemize}
\end{proposition}

\begin{proof}
We will prove $(1)$ by applying Variant \ref{subless2} of Proposition \ref{breadtime}.
Note that $\Mod^{\CommOp}_{A}(\calC)^{\otimes} \rightarrow \CommOp$ automatically
has units (Example \ref{swagg}). It will therefore suffice to prove the following:
\begin{itemize}
\item[$(\ast)$] For every object $X \in \Mod^{\CommOp}_{A}(\calC)^{\otimes}_{\seg{2}}$, there
exists a morphism $f: X \rightarrow Y$ in $\Mod^{\CommOp}_{A}(\calC)^{\otimes}$ covering
the active morphism $\seg{2} \rightarrow \seg{1}$ in $\FinSeg$ and classified by an
operadic $q$-colimit diagram, where $q: \Mod^{\CommOp}_{A}(\calC)^{\otimes} \rightarrow \CommOp$ denotes the projection. 
\end{itemize}

Let $\calD$ denote the subcategory
of $\Nerve(\FinSeg)_{ \seg{2}/ }$ whose objects are injective maps $\seg{2} \rightarrow \seg{n}$
and whose morphisms are diagrams
$$ \xymatrix{ & \seg{m} \ar[dr]^{u} &  \\
\seg{2} \ar[ur] \ar[rr] & & \seg{n}  }$$
where $u$ is active, so that $\calD^{\triangleright}$ maps to $\Nerve( \FinSeg)$ by a map
carrying the cone point of $\calD^{\triangleright}$ to $\seg{1}$. The category
$\calD$ has an evident forgetful functor $\calD \rightarrow \calK_{\CommOp} \times_{ \CommOp} \{ \seg{2} \}$, so that $X$ determines a functor $F_0 \in \Fun_{ \CommOp}( \calD, \calC^{\otimes})$. 
In view of Theorem \ref{colmodd}, to prove $(\ast)$ it will suffice to show that
$F_0$ can be extended to an operadic $p$-colimit diagram
$F \in \Fun_{ \CommOp}( \calD^{\triangleright}, \calC^{\otimes})$. 

In view of our assumption on $\calC^{\otimes}$ and Proposition \ref{slavewell}, the existence
of $F$ will follow provided that we can exhibit a cofinal map $\phi: \Nerve(\cDelta)^{op} \rightarrow \calD$. We define $\phi$ by a variation of Construction \ref{urpas}. Let $\calJ$ denote the category
whose objects are finite sets $S$ containing a pair of distinct points $x,y \in S$, so we have
a canonical equivalence $\calD \simeq \Nerve(\calJ)$. We will obtain $\phi$ as
the nerve of a functor $\phi_0: \cDelta^{op} \rightarrow \calJ$, where
$\phi_0( [n] )$ is the set of all downward-closed subsets of $[n]$
(with distinguished points given by $\emptyset, [n] \subseteq [n]$). To prove
that $\phi$ is cofinal, it will suffice to show that for each $\overline{S} = (S,x,y) \in \calJ$, the 
category $\calI = \cDelta^{op} \times_{ \calJ} \calJ_{ \overline{S}/}$ has weakly contractible nerve.
Writing $S = S_0 \coprod \{ x, y \}$, we can identify $\calI^{op}$ with the category
of simplices of the simplicial set $( \Delta^1)^{S_0}$. Since $( \Delta^1)^{S_0}$ is weakly contractible, it follows that $\Nerve(\calI)$ is weakly contractible. This proves $(1)$.

We will prove $(2)$ by applying Variant \ref{subless} of Proposition \ref{breadtime}.
Since the functor $\psi$ clearly preserves units (Corollary \ref{swink}), it suffices to verify the following:
\begin{itemize}
\item[$(\ast')$] Let $X \in \Mod^{\CommOp}_{A}(\calC)^{\otimes}_{\seg{2}}$, let
$\alpha: \seg{2} \rightarrow \seg{1}$ be an active morphism in $\Ass$, let
$\alpha_0$ be its image in $\CommOp$, and let $\overline{\alpha}_0:
X \rightarrow Y$ be a $q$-coCartesian morphism in $\Mod^{\CommOp}_{A}(\calC)^{\otimes}$
lifting $\alpha_0$. Then the induced morphism $\overline{\alpha}$ in $\Mod^{\Ass}_{A}(\calC)^{\otimes}$
is $q'$-coCartesian (where $q': \Mod^{\Ass}_{A}(\calC)^{\otimes} \rightarrow \Ass$ denotes the projection). 
\end{itemize}
Using Proposition \ref{chocolateoperad}, we deduce that $\overline{\alpha}_0$ is
classified by an operadic $q$-colimit diagram $$\Delta^1 \rightarrow \Mod^{\CommOp}_{A}(\calC)^{\otimes}.$$ Applying Theorem \ref{colmodd}, we deduce that the underlying diagram
$\calD^{\triangleright} \rightarrow \calC^{\otimes}$ is an operadic $p$-colimit diagram.
By the cofinality argument above, we conclude that $\overline{\alpha}_0$ induces
an operadic $p$-colimit diagram $\Nerve( \cDelta^{op} )^{\triangleright} \rightarrow \calC^{\otimes}$.
It follows from the proof of Theorem \ref{qwent} that $\overline{\alpha}$ is classified
by an operadic $q'$-colimit diagram $\Delta^1 \rightarrow \Mod^{\Ass}_{A}(\calC)^{\otimes}$, so
that $\overline{\alpha}$ is $q'$-coCartesian by virtue of Proposition \ref{chocolateoperad}.
\end{proof}

Our next result describes an important special feature of the commutative $\infty$-operad,
which makes the theory of commutative algebras (and their modules) very pleasant to work with:

\begin{theorem}\label{podus}
Let $\calC^{\otimes}$ be a symmetric monoidal $\infty$-category whose monoidal structure
is compatible with $\Nerve(\cDelta)^{op}$-indexed colimits. Then the map $p: \Mod^{\CommOp}(\calC)^{\otimes} \rightarrow \CAlg(\calC) \times \Nerve(\FinSeg)$ is a coCartesian fibration.
\end{theorem}

The proof of Theorem \ref{podus} will require some preliminaries.

\begin{lemma}\label{capse}
Let $S$ be an $\infty$-category and let
$p: \calC^{\otimes} \rightarrow S \times \Nerve(\FinSeg)$ be an $S$-family of $\infty$-operads.
Assume that:
\begin{itemize}
\item[$(a)$] The composite map $p': \calC^{\otimes} \rightarrow S$ is a Cartesian fibration;
moreover, the image in $\Nerve(\FinSeg)$ of any $p'$-Cartesian morphism of $\calC^{\otimes}$
is an equivalence.
\item[$(b)$] For each $s \in S$, the induced map $p_{s}: \calC^{\otimes}_{s} \rightarrow \Nerve(\FinSeg)$
is a coCartesian fibration.
\item[$(c)$] The underlying map $p_0: \calC \rightarrow S$ is a coCartesian fibration. 
\end{itemize}
Then:
\begin{itemize}
\item[$(1)$] The map  $p$ is a locally coCartesian fibration.
\item[$(2)$] A morphism $f$ in $\calC^{\otimes}$
is locally $p$-coCartesian if and only if it factors as a composition
$f'' \circ f'$, where $f'$ is a $p_{s}$-coCartesian morphism in
$\calC^{\otimes}_{s}$ for some $s \in S$ (here $p_{s}: \calC^{\otimes}_{s} \rightarrow \Nerve(\FinSeg)$
denotes the restriction of $p$) and $f'': Y \rightarrow Z$ is a morphism in $\calC^{\otimes}_{\seg{n}}$
with the following property: for $1 \leq i \leq n$, there exists a commutative diagram
$$ \xymatrix{ Y \ar[r]^{f''} \ar[d] & Z \ar[d] \\
Y_i \ar[r]^{f''_{i}} & Z_i }$$
where the vertical maps are inert morphisms of $\calC^{\otimes}$ lying over
$\Colp{i}: \seg{n} \rightarrow \seg{1}$ and the map $f''_{i}$ is a locally
$p_0$-coCartesian morphism in $\calC$.
\end{itemize}
\end{lemma}

\begin{proof}
Let $f \simeq f' \circ f''$ be as in $(2)$, and let
$y, z \in S$ denote the images of $Y, Z \in \calC^{\otimes}$. Using the equivalences
$$ \calC^{\otimes} \times_{ S \times \Nerve(\FinSeg)} (y, \seg{n}) \simeq
\calC_{y}^{n}
\quad \quad \calC^{\otimes} \times_{ S \times \Nerve(\FinSeg)} (z, \seg{n} ) \simeq \calC_{z}^{n},$$
we deduce that if $f''$ satisfies the stated condition, then $f''$ is locally $p$-coCartesian.
If $f'$ is a $p_{s}$-coCartesian morphism in $\calC^{\otimes}_{s}$, then condition
$(a)$ and Corollary \toposref{pannaheave} guarantee that $f'$ is $p$-coCartesian.
Replacing $p$ by the induced fibration $\calC^{\otimes} \times_{ S \times \Nerve(\FinSeg)} \Delta^2 \rightarrow \Delta^2$ and applying Proposition \toposref{protohermes}, we deduce that
$f$ is locally $p$-coCartesian. This proves the ``if'' direction of $(2)$.

To prove $(1)$, consider an object $X \in \calC^{\otimes}$ lying over
$(s, \seg{n}) \in S \times \Nerve(\FinSeg)$, and let $f_0: (s, \seg{n}) \rightarrow (s', \seg{n'})$
be a morphism in $S \times \Nerve(\FinSeg)$. Then $f_0$ factors canonically as a composition
$$ (s, \seg{n}) \stackrel{ f'_0}{\rightarrow} (s, \seg{n'}) \stackrel{ f''_0}{\rightarrow} ( s', \seg{n'}).$$
Using assumption $(b)$, we can lift $f'_0$ to a $p_{s}$-coCartesian morphism $f': X \rightarrow Y$ in $\calC^{\otimes}_{s}$, and using $(c)$ we can lift $f''_0$ to a morphism $f'': Y \rightarrow Z$ in
$\calC^{\otimes}_{\seg{n}}$ satisfying the condition given in $(2)$. It follows from the above argument
that $f = f'' \circ f'$ is a locally $p$-coCartesian morphism of $\calC^{\otimes}$ lifting
$f_0$. This proves $(1)$.

The ``only if'' direction of $(2)$ now follows from the above arguments together with the uniqueness properties of locally $p$-coCartesian morphisms.
\end{proof}

For the statement of the next lemma, we need to introduce a bit of terminology. Let
$p: \calC^{\otimes} \rightarrow S \times \Nerve(\FinSeg)$ be as in Lemma \ref{capse}, and
let $f: X \rightarrow Y$ be a locally $p$-coCartesian morphism in 
$\calC^{\otimes}_{\seg{n}}$ for some $n \geq 0$. If $\alpha: \seg{n} \rightarrow \seg{n'}$ is a
morphism in $\FinSeg$, we will say that $f$ is {\it $\alpha$-good} if the following condition
is satsified:
\begin{itemize}
\item[$(\ast)$] Let $x, y \in S$ denote the images of $X, Y \in \calC^{\otimes}$, and choose
morphisms $\overline{\alpha}_{x}: X \rightarrow X'$ in
$\calC^{\otimes}_{x}$ and $\overline{\alpha}_{y}: Y \rightarrow Y'$ in $\calC^{\otimes}_{y}$
which are $p_{x}$ and $p_{y}$-coCartesian lifts of $\alpha$, respectively, so that we have a commutative diagram
$$ \xymatrix{ X \ar[r]^{f} \ar[d]^{\overline{\alpha}_{x}} & Y \ar[d]^{\overline{\alpha}_{y}} \\
X' \ar[r]^{f'} & Y'. }$$
Then $f'$ is locally $p$-coCartesian.
\end{itemize}

\begin{lemma}\label{soapse}
Let $p: \calC^{\otimes} \rightarrow S \times \Nerve(\FinSeg)$ be as in Lemma \ref{capse}.
The following conditions are equivalent:
\begin{itemize}
\item[$(1)$] The map $p$ is a coCartesian fibration.
\item[$(2)$] For every morphism $f: X \rightarrow Y$ in $\calC^{\otimes}_{\seg{n}}$ and
every morphism $\alpha: \seg{n} \rightarrow \seg{n'}$, the morphism $f$ is $\alpha$-good.
\item[$(3)$] For every morphism $f: X \rightarrow Y$ in $\calC^{\otimes}_{\seg{n}}$ where
$n \in \{0,2\}$ and every active morphism $\alpha: \seg{n} \rightarrow \seg{1}$, the morphism $f$ is $\alpha$-good.
\end{itemize}
\end{lemma}

\begin{proof}
The implication $(1) \Rightarrow (2)$ follows from Proposition \toposref{protohermes} and
the implication $(2) \Rightarrow (3)$ is obvious. We next prove that $(2) \Rightarrow (1)$.
Since $p$ is a locally coCartesian fibration (Lemma \ref{capse}), it will suffice to show that if
$f: X \rightarrow Y$ and $g: Y \rightarrow Z$ are locally $p$-coCartesian morphsims in $\calC^{\otimes}$, then $g \circ f$ is locally $p$-coCartesian (Proposition \toposref{gotta}). 
Using Lemma \ref{capse}, we can assume that $f = f'' \circ f'$, where
$f'$ is $p$-coCartesian and $f''$ is a locally $p$-coCartesian morphism in 
$\calC^{\otimes}_{\seg{n}}$ for some $n \geq 0$. To prove that
$g \circ f \simeq (g \circ f'') \circ f'$ is locally $p$-coCartesian, it will suffice to show
that $g \circ f''$ is locally $p$-coCartesian. We may therefore replace $f$ by $f''$
and thereby assume that $f$ has degenerate image in $\Nerve(\FinSeg)$.
Applying Lemma \ref{capse} again, we can write $g = g'' \circ g'$ where
$g': Y \rightarrow Y'$ is a $p$-coCartesian morphism in $\calC^{\otimes}_{y}$
covering some map $\alpha: \seg{n} \rightarrow \seg{n'}$, and 
$g''$ is a locally $p$-coCartesian morphism in $\calC^{\otimes}_{\seg{n''}}$.
Let $x$ denote the image of $X$ in $S$, and let $k: X \rightarrow X'$ be a
$p$-coCartesian morphism in $\calC^{\otimes}_{x}$ lying over $\alpha$. 
We have a commutative diagram
$$ \xymatrix{ X \ar[r]^{k} \ar[d] & Y \ar[d]^{g'} \ar[dr]^{g} \\
X' \ar[r]^{h} & Y' \ar[r]^{g''} & Z. }$$
Since $f$ is $\alpha$-good (by virtue of assumption $(2)$), we deduce that
$h$ is locally $p$-coCartesian. To prove that $g \circ f \simeq (g'' \circ h) \circ k$
is locally $p$-coCartesian, it will suffice (by virtue of Lemma \ref{capse}) to show
that $g'' \circ h$ is locally $p$-coCartesian. We claim more generally that the collection of
locally $p$-coCartesian edges in $\calC^{\otimes}_{\seg{n}}$ is closed under composition.
This follows from the observation that $\calC^{\otimes}_{\seg{n}} \rightarrow S$ is a coCartesian fibration (being equivalent to a fiber power of the coCartesian fibration $p_0: \calC \rightarrow S$).
This completes the proof that $(2) \Rightarrow (3)$.

We now prove that $(3) \Rightarrow (2)$. Let us say that a morphism
$\alpha: \seg{n} \rightarrow \seg{m}$ is {\it perfect} if every locally $p$-coCartesian morphism
$f$ in $\calC^{\otimes}_{\seg{n}}$ is $\alpha$-good. Assumption $(3)$ guarantees
that the active morphisms $\seg{2} \rightarrow \seg{1}$ and $\seg{0} \rightarrow \seg{1}$
are perfect; we wish to prove that every morphism in $\FinSeg$ is perfect. This follows immediately
from the following three claims:
\begin{itemize}
\item[$(i)$] If $\alpha: \seg{n} \rightarrow \seg{m}$ is perfect, then for each $k \geq 0$ the induced map
$\seg{n+k} \rightarrow \seg{m+k}$ is perfect.
\item[$(ii)$] The unique (inert) morphism $\seg{1} \rightarrow \seg{0}$ is perfect.
\item[$(iii)$] The collection of perfect morphisms in $\FinSeg$ is closed under composition.
\end{itemize}
Assertions $(i)$ and $(ii)$ are obvious. To prove $(iii)$, suppose that
we are given perfect morphisms $\alpha: \seg{n} \rightarrow \seg{m}$ and
$\beta: \seg{m} \rightarrow \seg{k}$. Let $f: X \rightarrow Y$ be a locally $p$-coCartesian
morphism in $\calC^{\otimes}_{\seg{n}}$, and form a commutative diagram
$$ \xymatrix{ X \ar[r]^{f} \ar[d] & Y\ar[d] \\
X' \ar[r]^{f'} \ar[d] & Y' \ar[d] \\
X'' \ar[r]^{f''} & Y'' }$$
where the upper vertical maps are $p$-coCartesian lifts of $\alpha$ and the lower vertical maps are $p$-coCartesian lifts of $\beta$. Since $\alpha$ is perfect, the map $f'$ is locally $p$-coCartesian.
Since $\beta$ is perfect, the map $f''$ is locally $p$-coCartesian, from which it follows that
$f$ is $\beta \circ \alpha$-good as desired.
\end{proof}

\begin{lemma}\label{costick}
Let $\calC^{\otimes} \rightarrow \Ass$ be a coCartesian fibration of $\infty$-operads such that the induced $\Ass$-monoidal structure on $\calC$ is compatible with $\Nerve(\cDelta)^{op}$-indexed colimits. Then the forgetful map $\theta: \Mod^{L}(\calC) \rightarrow \Alg(\calC)$ is a coCartesian fibration.
\end{lemma}

\begin{proof}
Corollary \monoidref{thetacart} guarantees that $\theta$ is a Cartesian fibration.
The assertion that $\theta$ is a coCartesian fibration follows from Proposition
\toposref{getcocart} and Lemma \monoidref{peacestick}.
\end{proof}

\begin{proof}[Proof of Theorem \ref{podus}]
We first show that $p$ is a locally coCartesian fibration by verifying the hypotheses of Lemma \ref{capse}:
\begin{itemize}
\item[$(a)$] The map $p': \Mod^{\CommOp}(\calC)^{\otimes} \rightarrow \CAlg(\calC)$ is
a Cartesian fibration by virtue of Corollary \ref{postsabel}; moreover, a morphism
$f$ in $\Mod^{\CommOp}(\calC)^{\otimes}$ is $p'$-Cartesian if and only if its image in
$\calC^{\otimes}$ is an equivalence (which implies that its image in $\Nerve( \FinSeg)$ is an equivalence).
\item[$(b)$] For each $A \in \CAlg(\calC)$, the fiber $\Mod^{\CommOp}_{A}(\calC)^{\otimes}$
is a symmetric monoidal $\infty$-category (Proposition \ref{stancer}).
\item[$(c)$] The map $p_0: \Mod^{\CommOp}(\calC) \rightarrow \CAlg(\calC)$ is a
coCartesian fibration. This follows from Proposition \ref{socrime}, because the forgetful
functor $\Mod^{L}(\calC) \rightarrow \Alg(\calC)$ is a coCartesian fibration (Lemma \ref{costick}).
\end{itemize}
To complete the proof, it will suffice to show that $p$ satisfies condition
$(3)$ of Lemma \ref{soapse}. Fix a map $f_0: A \rightarrow B$ in $\CAlg(\calC)$
and a locally $p$-coCartesian morphism $f: M_{A} \rightarrow M_{B}$ in
$\Mod^{\CommOp}(\calC)^{\otimes}_{\seg{n}}$ covering $f$; we wish to prove
that $f$ is $\alpha$-good, where $\alpha: \seg{n} \rightarrow \seg{1}$ is the unique active morphism
(here $n = 0$ or $n=2$). Let $\calD^{\otimes} = \Mod^{\calO}_{A}(\calC)^{\otimes}$.
Proposition \ref{stancer} implies that $\calD^{\otimes}$ is a symmetric monoidal
$\infty$-category, Propositions \ref{socrime} and Corollary \monoidref{gloop} imply that
$\calD$ admits $\Nerve(\cDelta)^{op}$-indexed colimits, and Proposition \monoidref{colimp}
implies that the symmetric monoidal structure on $\calD$ is compatible with
$\Nerve(\cDelta)^{op}$-indexed colimits. Corollary \ref{costabs} implies that
the forgetful functor $\Mod^{\CommOp}(\calD)^{\otimes} \rightarrow
\Mod^{\CommOp}(\calC)^{\otimes} \times_{ \CAlg(\calC)} \CAlg(\calC)^{A/}$ is
an equivalence of $\infty$-categories. Replacing $\calC^{\otimes}$ by
$\calD^{\otimes}$, we may assume without loss of generality that
$A$ is a unit algebra in $\calC^{\otimes}$.

Suppose first that $n=0$. Unwinding the definitions, we are required to show that
if $M_{A} \rightarrow M'_{A}$ and $M_{B} \rightarrow M'_{B}$ are morphisms 
$\Mod^{\CommOp}(\calC)^{\otimes}$
which exhibit $M'_{A} \in \Mod^{\CommOp}_{A}(\calC)$ as a $p_{A}$-unit
(see Definition \ref{cuspe}; here $p_{A}: \Mod^{\calO}_{A}(\calC)^{\otimes} \rightarrow
\calO^{\otimes}$ denotes the restriction of $p$) and $M'_{B} \in \Mod^{\CommOp}_{B}(\calC)$
as a $p_{B}$-unit (where $p_B$ is defined similarly), then the induced map
$f': M'_{A} \rightarrow M'_{B}$ is locally $p$-coCartesian. Using Corollary
\ref{postsabel}, it suffices to show that $\phi(f')$ is an equivalence, where
$\phi: \Mod^{\calO}(\calC) \rightarrow \calC$ is the forgetful functor.
Using Corollary \ref{swink} and Proposition \ref{socrime}, this translates into the following
assertion: the map $A \rightarrow B$ exhibits $B$ as the free left $B$-module generated
by $A$. This follows from Proposition \monoidref{pretara}.

We now treat the case $n=2$. Since $A$ is a unit algebra, Proposition \ref{stupe} allows us to identify $M_{A}$ with a pair of objects $P, Q \in \calC$. We can identify $M_{B}$ with a pair
of objects $P_{B}, Q_{B} \in \Mod_{B}^{\CommOp}(\calC)$. Using the
equivalence of Proposition \ref{socrime}, we will view $P_{B}$ as a right $B$-module
object of $\calC$, and $Q_{B}$ as a left $B$-module object of $\calC$. Unwinding
the definitions (and using Proposition \ref{stancer}), it suffices to show that the canonical map
$f': P \otimes Q \rightarrow (P_{B} \otimes_{B} Q_{B})$ exhibits the relative tensor product
$P_{B} \otimes_{B} Q_{B}$ as the free $B$-module generated by $P \otimes Q$.
Using the push-pull formula (Proposition \monoidref{ussr}), we obtain a canonical
identification $P_{B} \otimes_{B} Q_{B} \simeq P \otimes_{A} Q_{B} \simeq P \otimes Q_{B}$.
According to Proposition \monoidref{pretara}, it will suffice to show that $f'$ induces an equivalence
$\phi: P \otimes Q \otimes B \rightarrow P \otimes Q_{B}$. It now suffices to observe that
$\phi$ is equivalent to the tensor product of the identity map $\id_{P}$ with the
equivalence $Q \otimes B \simeq Q_{B}$ provided by Proposition \monoidref{pretara}.
\end{proof}

We now return to the proof of Proposition \ref{socrime}.

\begin{proof}[Proof of Proposition \ref{socrime}]
Consider the functor $\Phi: \cDelta^{op} \times [1] \rightarrow (\FinSeg)_{ \seg{1}/ }$ described as follows:
\begin{itemize}
\item For $n \geq 0$, we have $\Phi( [n], 0) = ( \alpha: \seg{1} \rightarrow \psi( [n] \star [0] ) )$, where
$\psi$ denotes the composition of the functor $\phi: \cDelta^{op} \rightarrow \CatAss$ with the forgetful functor $\CatAss \rightarrow \FinSeg$ and $\alpha: \seg{1} \rightarrow \psi( [n] \star [0])$
carries $1 \in \seg{1}$ to the point $[n] \in \psi( [n] \star [0] )$.
\item For $n \geq 0$, we have $\Phi( [n], 1) = ( \alpha: \seg{1} \rightarrow \psi( [n] ) )$, where
$\alpha$ is the null morphism given by the composition $\seg{1} \rightarrow \seg{0} \rightarrow
\psi( [n] ) \simeq \seg{n}$.
\end{itemize}

Let $\calI_0$ denote the full subcategory of $( \FinSeg)_{ \seg{1}/}$ spanned by the semi-inert morphisms $\seg{1} \rightarrow \seg{n}$. We observe that $\Phi$ defines a functor
$\cDelta^{op} \times [1] \rightarrow \calI_0$. 
Let $\calI$ denote the categorical mapping cylinder of the functor $\Phi$. More precisely, the category $\calI$ is defined as follows:

\begin{itemize}
\item An object of $\calI$ is either an object $\alpha: \seg{1} \rightarrow \seg{n}$ of
$\calI_0$ or an object $([n], i)$ of $\cDelta^{op} \times [1]$.
\item Morphisms in $\calI$ are defined as follows:
$$ \Hom_{\calI}( \alpha, \alpha') = \Hom_{\calI_0}( \alpha, \alpha') \quad \quad
\Hom_{\calI}( ([n], i), ([n'],i')) = \Hom_{ \cDelta^{op} \times [1]}( ([n], i), ([n'], i'))$$
$$ \Hom_{ \calI}( \alpha, ([n],i)) = \Hom_{\calI_0}( \alpha, \Phi( [n],i)) \quad \quad
\Hom_{\calI}( ([n],i), \alpha) = \emptyset.$$
\end{itemize}

The full subcategory of $\calK_{\CommOp} \times_{ \CommOp} \seg{1}$ spanned by the null morphisms
$\seg{1} \rightarrow \seg{n}$ is actually isomorphic (rather than merely equivalent) to
$\CommOp$. Consequently, we have an isomorphism $\Mod^{\CommOp}(\calC) \simeq
\MMod^{\CommOp}(\calC)$, where the latter can be identified with a full subcategory of
$\Fun_{ \CommOp}( \Nerve(\calI_0), \calC^{\otimes})$. We regard
$\Nerve(\calI)$ as equipped with a forgetful functor to $\CommOp$, given by composing the retraction
$r: \calI \rightarrow \calI_0$ with the forgetful functor $\Nerve(\calI_0) \rightarrow \CommOp$.
Let $\calD$ denote the full subcategory of
$\Fun_{ \CommOp}( \Nerve(\calI), \calC^{\otimes})$ spanned by those functors $F$ which satisfy the following conditions:
\begin{itemize}
\item[$(i)$] The restriction of $F$ to $\Nerve(\calI_0)$ belongs to $\Mod^{\CommOp}(\calC)$.
\item[$(ii)$] For every object $([n], i) \in \cDelta^{op} \times [1]$, the canonical map
$F( \Phi([n],i) ) \rightarrow F( [n], i)$ is an equivalence in $\calC^{\otimes}$ (equivalently:
$F$ is a $p$-left Kan extension of $F | \Nerve(\calI_0)$).
\end{itemize}
Let $\calI'$ denote the full subcategory of $\calI$ spanned by those objects corresponding
to null morphisms $\seg{1} \rightarrow \seg{n}$ in $\FinSeg$ or 
having the form $\alpha: \seg{1} \rightarrow \seg{n}$ where $\alpha$ is null,
or pairs $([n],i) \in \cDelta^{op} \times [1]$ where $i = 1$. Let $\calD'$ denote the full subcategory of
$\Fun_{ \CommOp}( \Nerve(\calI'), \calC^{\otimes})$ spanned by those functors $F$ satisfying the following
conditions:
\begin{itemize}
\item[$(i')$] The restriction of $F$ to $\Nerve(\calI_0 \cap \calI')$ belongs to $\CAlg(\calC)$.
\item[$(ii')$] For $n \geq 0$, the canonical map $F( \Phi([n],1) ) \rightarrow F( [n], 1)$ is an equivalence in $\calC^{\otimes}$ (equivalently: $F$ is a $p$-left Kan extension of $F | \Nerve(\calI_0 \cap \calI')$).
\end{itemize}

Using Proposition \toposref{lklk}, we deduce that inclusion $\calI_0 \hookrightarrow \calI$
induces a trivial Kan fibrations 
$$\calD \rightarrow \Mod^{\CommOp}(\calC) \quad \quad \calD_0 \rightarrow \CAlg(\calC).$$
We have a commutative diagram
$$ \xymatrix{ \Mod^{\CommOp}(\calC) \ar[r] \ar[d] & \calD \ar[d] \ar[r] & \Mod^{L}(\calC) \ar[d] \\
\CAlg(\calC) \ar[r] & \calD' \ar[r] & \Alg(\calC) }$$
where the horizontal maps on the left are given by composition with the retraction $r$
(since these are sections of the trivial Kan fibrations above, they are categorical equivalences)
and the horizontal maps on the right are given by composition with the inclusion
$\cDelta^{op} \times [1] \hookrightarrow \calI$. Consequently, it will suffice to prove
that the square on the right is a homotopy pullback diagram. Since the vertical maps in this
square are categorical fibrations between $\infty$-categories, it will suffice to show that
the map $\calD \rightarrow \Mod^{L}(\calC) \times_{ \Alg(\calC) } \calD'$ is a trivial Kan fibration.
Let $\calI''$ denote the full subcategory of $\calI$ spanned by the objects of
$\calI'$ and $\cDelta^{op} \times [1]$. We observe that
$\Mod^{L}(\calC) \times_{ \Alg(\calC)} \calD'$ can be identified with the full subcategory of
$\calD'' \subseteq \Fun_{ \CommOp}( \Nerve(\calI''), \calC^{\otimes})$ spanned by those functors
$F$ which satisfy conditions $(i')$, $(ii')$, and the following additional condition:
\begin{itemize}
\item[$(iii')$] The restriction of $F$ to $\Nerve(\cDelta)^{op} \times \Delta^1$ is a
left module object of $\calC$. 
\end{itemize}

We wish to prove that $\calD \rightarrow \calD''$ is a trivial Kan fibration. In view of Proposition \toposref{lklk}, it will suffice to prove the following:
\begin{itemize}
\item[$(a)$] Every functor $F_0 \in \calD''$ admits a $p$-right Kan extension
$F \in \Fun_{\CommOp}( \Nerve(\calI), \calC^{\otimes})$.
\item[$(b)$] Let $F \in \Fun_{\CommOp}( \Nerve(\calI), \calC^{\otimes})$ be
a functor such that $F_0 = F | \Nerve(\calI'') \in \calD''$. Then
$F \in \calD$ if and only if $F$ is a $p$-right Kan extension of $F_0$.
\end{itemize}

We first prove $(a)$. Let $F_0 \in \calD''$, and consider an object
$\alpha: \seg{1} \rightarrow \seg{n}$ of $\calI$ which does not belong to
$\calI''$, and let $i = \alpha(1) \in \nostar{n}$. Let $\calJ = \calI'' \times_{ \calI} \calI_{\alpha/}$. We wish to prove that the diagram $\Nerve(\calJ) \rightarrow \calC^{\otimes}$ determined by $F_0$
can be extended to a $p$-limit diagram covering the map
$$\Nerve(\calJ)^{\triangleleft} \rightarrow \Nerve( \calI_{\alpha/})^{\triangleleft}
\rightarrow \Nerve(\calI) \rightarrow \Nerve(\FinSeg).$$
Let $\calJ_0$ denote the full subcategory of $\calJ$ spanned by the objects
of $\calI' \times_{ \calI} \calI_{\alpha/}$ together with those maps
$\alpha \rightarrow ( [m], 0)$ in $\calI$ for which the underlying map
$u: \seg{n} \rightarrow \psi( [m] \star [0] )$ satisfies the following condition:
if $u(j) = u(i)$ for $j \in \nostar{n}$, then $i = j$. It is not difficult to see that the inclusion
$\calJ_0 \subseteq \calJ$ admits a right adjoint, so that
$\Nerve(\calJ_0)^{op} \subseteq \Nerve(\calJ)^{op}$ is cofinal.
Consequently, it will suffice to show that the induced map
$G: \Nerve(\calJ_0) \rightarrow \calC^{\otimes}$ can be extended
to a $p$-limit diagram (compatible with the underlying map
$\Nerve(\calJ_0)^{\triangleleft} \rightarrow \Nerve(\FinSeg)$).

Let $\calJ_1$ denote the full subcategory of $\calJ_0$ spanned by
the objects of $\calI' \times_{ \calI} \calI_{\alpha/}$ together with the morphism
$s: \alpha \rightarrow ([0],0)$ in $\calI$ given by the map
$\Colp{i}: \seg{n} \rightarrow \psi( [0] \star [0] ) \simeq \seg{1}$.
We claim that $G$ is a $p$-right Kan extension of $G | \Nerve(\calJ_1)$.
To prove this, consider an arbitrary object $J \in \calJ_0$ which does not belong
to $\calJ_1$, which we identify with a map $t: \alpha \rightarrow ( [m], 0)$ in $\calI$.
We have a commutative diagram
$$ \xymatrix{ (t: \alpha \rightarrow ( [m], 0)) \ar[r] \ar[d] & (s: \alpha \rightarrow ([0],0)) \ar[d] \\
(t': \alpha \rightarrow ([m], 1) ) \ar[r] & (s' \alpha \rightarrow ([0],1) }$$
in $\calJ_0$, which we can identify with a diagram $s \rightarrow s' \leftarrow t'$
in $(\calJ_0)_{J/} \times_{ \calJ_0} \calJ_1$. Using Theorem \toposref{hollowtt}, we deduce
that this diagram determines a cofinal map $\Lambda^{2}_{0} \rightarrow
\Nerve( (\calJ_0)_{J/} \times_{ \calJ_0} \calJ_1 )^{op}$. Consequently, to prove
that $G$ is a $p$-right Kan extension of $G | \Nerve(\calJ_1)$, it suffices to verify that
the diagram
$$ \xymatrix{ F_0( [m], 0) \ar[r] \ar[d] & F_0([0],0) \ar[d] \\
F_0( [m], 1) \ar[r] & F_0( [m],1) }$$
is a $p$-limit diagram, which follows from $(iii')$. Using Lemma \toposref{kan0}, we are
are reduced to proving that $G_1 = G | \Nerve(\calJ_1)$ can be extended to
a $p$-limit diagram lifting the evident map $\Nerve(\calJ_1)^{\triangleleft} \rightarrow \Nerve(\FinSeg)$.

Let $\beta: \seg{n} \rightarrow \seg{n-1}$ be an inert morphism such that
$\beta(i) = \ast$, and let $\alpha': \seg{1} \rightarrow \seg{n-1}$ be given by the
composition $\beta \circ \alpha$. Let $\calJ_2$ denote the full subcategory of $\calJ_1$ spanned by the map $s: \alpha \rightarrow ([0],0)$, the induced map $s': \alpha \rightarrow ([0],1)$, and the natural
transformation $\alpha \rightarrow \alpha'$ induced by $\beta$. Using Theorem \toposref{hollowtt}, we
deduce that the inclusion $\Nerve(\calJ_2)^{op} \subseteq \Nerve(\calJ_1)^{op}$ is cofinal.
Consequently, we are reduced to proving that $G_2 = G | \Nerve(\calJ_2)$ can be extended
to a $p$-limit diagram lifting the evident map $\Nerve(\calJ_2)^{\triangleleft} \rightarrow \Nerve(\FinSeg)$. In other words, we must find a $p$-limit diagram
$$ \xymatrix{ F(\alpha) \ar[r] \ar[d] & F_0( [0],0) \ar[d] \\
F_0( \alpha') \ar[r] & F_0( [0], 1) }$$
covering the diagram
$$ \xymatrix{ \seg{n} \ar[d]^{\beta} \ar[r]^{ \Colp{i} } & \seg{1} \ar[d] \\
\seg{n-1} \ar[r] & \seg{0} }$$
in $\FinSeg$. The existence of such a diagram follows immediately from our assumption
that $p: \calC^{\otimes} \rightarrow \Nerve(\FinSeg)$ is an $\infty$-operad (see Proposition \ref{prejule}).
This completes the proof of $(a)$. Moreover, the proof also gives the following analogue of
$(b)$:

\begin{itemize}
\item[$(b')$] Let $F \in \Fun_{\CommOp}( \Nerve(\calI), \calC^{\otimes})$ be
a functor such that $F_0 = F | \Nerve(\calI'') \in \calD''$. Then
$F$ is a $p$-right Kan extension of $F_0$ if and only if it satisfies
the following condition:
\begin{itemize}
\item[$(\ast)$] For every
object $\alpha: \seg{1} \rightarrow \seg{n}$ of $\calI$ which does not belong
to $\calI''$, the induced diagram
$$ \xymatrix{ F( \alpha) \ar[r]^{u} \ar[d]^{v} & F_0( [0], 0) \ar[d] \\
F_0( \alpha') \ar[r] & F_0( [0],1) }$$ is a $p$-limit diagram.
\end{itemize}
\end{itemize}

We observe that $(\ast)$ can be reformulated as follows:
\begin{itemize}
\item[$(\ast')$]  For every
object $\alpha: \seg{1} \rightarrow \seg{n}$ of $\calI$ which does not belong
to $\calI''$, the induced morphisms $u: F( \alpha) \rightarrow F_0([0],0)$ and
$v: F(\alpha) \rightarrow F_0(\alpha')$ are inert.
\end{itemize}

To complete the proof, it will suffice to show that  if $F \in \Fun_{\CommOp}( \Nerve(\calI), \calC^{\otimes})$
satisfies $(i')$, $(ii')$, and $(iii')$, then $F$ satisfies condition $(\ast')$ if and only if it
satisfies conditions $(i)$ and $(ii)$. We first prove the ``if'' direction. The assertion
that $v$ is inert follows immediately from $(i)$. To see that $u$ is inert, choose a map
$u_0: \alpha \rightarrow ( [n-1], 0)$ in $\calI$ corresponding to an isomorphism
$\seg{n} \simeq \psi( [n-1] \star [0])$, so that $u$ factors as a composition
$$ F(\alpha) \stackrel{F(u_0)}{\rightarrow} F( [n-1],0) \stackrel{F(u_1)}{\rightarrow} F([0], 0).$$
The map $F(u_0)$ is inert by virtue of $(ii)$, and the map $F(u_1)$ is inert by virtue of
$(iii')$. 

We now prove the ``only if'' direction. Suppose that $F$ satisfies $(i')$, $(ii')$, $(iii')$,
and $(\ast')$. We first claim that $F$ satisfies $(i)$. Let $f: \alpha \rightarrow \beta$
be a morphism in $\calI_0$ whose image in $\FinSeg$ is inert; we wish to prove
that $F(f)$ is an inert morphism in $\calC^{\otimes}$. If $\alpha \in \calI'$, then
this follows from assumption $(i')$. If $\beta \in \calI'$ and $\alpha \notin \calI'$, then
we can factor $f$ as a composition
$$ \alpha \stackrel{f'}{\rightarrow} \alpha' \stackrel{f''}{\rightarrow} \beta$$
where $\alpha' \in \calI'$ is defined as above. Then $F(f'')$ is inert by
virtue of $(i')$, while $F(f')$ is inert by $(\ast')$, so that $F(f)$ is inert as desired.
Finally, suppose that $\beta: \seg{1} \rightarrow \seg{m}$ does not belong to
$\calI'$, so that $\beta(1) = i \in \nostar{m}$. To prove that $F(f)$ is inert, it will suffice
to show that for each $j \in \nostar{m}$, there exists an inert morphism
$\gamma_{j}: F( \beta) \rightarrow C_{j'}$ in $\calC^{\otimes}$ covering
$\Colp{j}: \seg{m} \rightarrow \seg{1}$ such that
$\gamma_{j} \circ F(f)$ is inert. If $j \neq i$, we can take
$\gamma_{j} = F( g)$, where $\beta_{j} = \Colp{j} \circ \beta$ and
$g: \beta \rightarrow \beta_{j}$ is the induced map; then $\beta_{j} \in \calI'$ so that
$F(g)$ and $F(g \circ f)$ are both inert by the arguments presented above.
If $i = j$, we take $\gamma_{j} = F(g)$ where $g$ is the map
$\beta \rightarrow ([0], 0)$ in $\calI$ determined by $\Colp{i}$. Then
$F(g)$ and $F(g \circ f)$ are both inert by virtue of $(\ast')$. This completes the proof
of $(i)$.

We now prove $(ii)$. Let $([n], i) \in \Nerve(\cDelta)^{op} \times \Delta^1$, let
$\alpha = \Psi( [n], i) \in \calI_0$, and let $f: \alpha \rightarrow ([n],i)$ be the canonical
morphism in $\calI$. We wish to prove that $F(f)$ is an equivalence in $\calC^{\otimes}$.
If $i = 1$, this follows from $(ii')$. Assume therefore that $i = 0$. Observe that $\alpha$ can be identified with the morphism $\seg{1} \rightarrow \seg{n+1}$ carrying $1 \in \seg{1}$ to $n+1 \in \seg{n+1}$. 
To prove that $F(f)$ is an equivalence, it will suffice to show that for each $j \in \nostar{n+1}$, there
exists an inert morphism $\gamma_{j}: F( [n],0) \rightarrow C_i$ in $\calC^{\otimes}$
covering $\Colp{j}: \seg{n+1} \rightarrow \seg{1}$
such that $\gamma_{j} \circ F(f)$ is also inert. If $j \neq n+1$, we take
$\gamma_{j} = F(g)$, where $g: ([n], 0) \rightarrow ([1], 1)$ is a morphism in
$\cDelta^{op} \times [1]$ covering the map $\Colp{i}$. 
Then $\gamma_{j}$ is inert by virtue of $(iii')$, while $\gamma_{j} \circ F(f) \simeq F(g \circ f)$
can be written as a composition
$$ F( \alpha) \stackrel{ F( h')}{\rightarrow} F( \alpha') \stackrel{F(h'')}{\rightarrow}
F( r([1],1)) \stackrel{ F(h''')}{\rightarrow} F( [1],1).$$
Assumption $(\ast)$ guarantees that $F(h')$ is inert, assumption
$(i')$ guarantees that $F(h'')$ is inert, and assumption $(ii')$ guarantees that
$F(h''')$ is inert. It follows that $F( h''' \circ h'' \circ h') \simeq F(g \circ f)$ is inert as desired.
In the case $j = n+1$, we instead take $\gamma_j = F(g)$ where
$g: ([n],0) \rightarrow ([0],0)$ is the morphism in $\calI$ determined by the map $[0] \rightarrow [n]$ in $\cDelta$ carrying $0 \in [0]$ to $n \in [n]$. Then $\gamma_j = F(g)$ is inert by virtue of
$(iii')$, while $\gamma_{j} \circ F(f) \simeq F(g \circ f)$ is inert by virtue of $(\ast')$.
\end{proof}

\section{Commutative Ring Spectra}\label{crung}

Our goal in this section is to apply our general formalism of symmetric monoidal $\infty$-categories, algebras, and modules in the setting of stable homotopy theory. We will begin in \S \ref{comm7} by showing that the $\infty$-category $\Spectra$ of spectra admits a symmetric monoidal structure,
given concretely by the classical {\it smash product} of spectra. This allows us to make sense of
commutative algebra objects of $\Spectra$. In \S \ref{comm8} we will study these
commutative algebras, which we refer to as {\it $\OpE{\infty}$-rings} (or sometimes
as {\it $\OpE{\infty}$-ring spectra}). 

To facilitate the comparison of our theory of 
$\OpE{\infty}$-rings with more classical approaches to the study of structured ring spectra,
we will prove in \S \ref{comm9} a rectification result (Theorem \ref{cbeckify}) which implies,
in particular, that the $\infty$-category of $\OpE{\infty}$-rings can be identified with the underlying
$\infty$-category of a model category whose objects are commutative monoids in the category of
symmetric spectra (Example \ref{bulwump}).

\subsection{Commutativity of the Smash Product}\label{comm7}

In \S \monoidref{monoid6.0}, we showed that the $\infty$-category $\Spectra$ of spectra admits an essentially unique monoidal structure, subject to the condition that the tensor product
$\otimes: \Spectra \times \Spectra \rightarrow \Spectra$ preserve colimits separately in each variable (Corollary \monoidref{surcoi}). The induced monoidal structure on the stable homotopy category
$\h{\Spectra}$ agrees with the classical {\it smash product} operation on spectra. In fact, the classical smash product endows $\h{\Spectra}$ with the structure of a symmetric monoidal category: 
given a pair of spectra $X$ and $Y$, there is a homotopy equivalence of spectra $X \otimes Y \simeq Y \otimes X$, which is well-defined up to homotopy. The goal of this section is to explain this observation by showing that the monoidal structure on $\Spectra$ can be promoted to a symmetric monoidal structure. We will obtain this result by a refinement of the method used in \S \monoidref{monoid6.0}.

We begin by reviewing a bit of classical category theory. Let $\calC$ be a monoidal category. The category $\Fun( \calC^{op}, \Set)$ of presheaves of
sets on $\calC$ then inherits a monoidal structure, which is characterized by the following assertions:
\begin{itemize}
\item[$(1)$] The Yoneda embedding $\calC \rightarrow \Fun( \calC^{op}, \Set)$ is a monoidal functor.
\item[$(2)$] The tensor product functor $\otimes: \Fun( \calC^{op}, \Set) \times
\Fun( \calC^{op}, \Set) \rightarrow \Fun( \calC^{op}, \Set)$ preserves colimits separately in each variable.
\end{itemize}
The bifunctor $\otimes: \Fun( \calC^{op}, \Set) \times
\Fun( \calC^{op}, \Set) \rightarrow \Fun( \calC^{op}, \Set)$ is called the {\it Day convolution product}.
More concretely, this operation is given by the composition
\begin{eqnarray*} \Fun( \calC^{op}, \Set) 
\times \Fun( \calC^{op}, \Set) 
& \rightarrow & \Fun( \calC^{op} \times \calC^{op}, \Set \times \Set) \\
& \stackrel{\phi}{\rightarrow} & \Fun( \calC^{op} \times \calC^{op}, \Set) \\
& \stackrel{\phi'}{\rightarrow} & \Fun( \calC^{op}, \Set)
\end{eqnarray*}
where $\phi$ is induced by the Cartesian product functor $\Set \times \Set \rightarrow \Set$, and
$\phi'$ is given by left Kan extension along the functor $\calC^{op} \times \calC^{op} \rightarrow \calC^{op}$ determined by the monoidal structure on $\calC$.

Our first goal in this section is to introduce an $\infty$-categorical version of the Day convolution product. For this, we need to recall a bit of terminology from \S \toposref{agileco}.


\begin{notation}
Given a pair of collections of diagrams $\calR_1 = \{ \overline{p}_{\alpha}: K^{\triangleright}_{\alpha} \rightarrow \calC_1 \}$ and $\calR_2 = \{ \overline{q}_{\alpha}: K^{\triangleright}_{\beta} \rightarrow \calC_2 \}$, we let $\calR_1 \boxtimes \calR_2$ denote the collection of all diagrams 
$\overline{r}: K^{\triangleright} \rightarrow \calC_1 \times \calC_2$ satisfying one of the following conditions:
\begin{itemize}
\item[$(1)$] There exists an index $\alpha$ and an object $C_2 \in \calC_2$ such that
$K = K_{\alpha}$ and $\overline{r}$ is given by the composition
$$ K^{\triangleright} \simeq K_{\alpha}^{\triangleright} \times \{ C_2 \} 
\stackrel{ \overline{p}_{\alpha}}{\rightarrow} \calC_1 \times \calC_2.$$
\item[$(2)$] There exists an index $\beta$ and an object $C_1 \in \calC_1$ such that
$K = K_{\beta}$ and $\overline{r}$ is given by the composition
$$ K^{\triangleright} \simeq \{ C_1 \} \times K_{\beta}^{\triangleright}
\stackrel{\overline{q}_{\beta}}{\rightarrow} \calC_1 \times \calC_2.$$
\end{itemize}

The operation $\boxtimes$ is coherently associative (and commutative). In other words, if for each $1 \leq i \leq n$ we are given an $\infty$-category $\calC_i$ and a collection $\calR_i$ of diagrams in $\calC_i$, then we have a well-defined collection $\calR_1 \boxtimes \ldots \boxtimes \calR_n$ of diagrams in the product $\calC_1 \times \ldots \times \calC_n$. 
so there is no ambiguity in writing expressions such as $\calR_1 \boxtimes \ldots \boxtimes \calR_{n}$. In the case $n=0$, we agree that this product coincides with the empty set of diagrams in the final $\infty$-category $\Delta^0$.
\end{notation}

We now proceed to construct a very large diagram of large symmetric monoidal $\infty$-categories.
First, we need a bit more notation.

\begin{notation}\label{surpint}
Let $\widehat{\Cat}_{\infty}$ denote the $\infty$-category of (not necessarily small) $\infty$-categories. Then $\widehat{\Cat}_{\infty}$ admits finite products. Consequently, there a Cartesian
symmetric monoidal structure on $\widehat{\Cat}_{\infty}$. For convenience, we will give an explicit construction of this symmetric monoidal structure.

Let $\bfA$ be the category of (not necessarily small) marked simplicial sets, as defined in \S \toposref{twuf}. Then $\bfA$ has the structure of a simplicial symmetric monoidal model category, where the symmetric monoidal structure is given by Cartesian product. We let
$\widehat{\Cat}_{\infty}^{\otimes}$ denote the $\infty$-category $\Nerve( \bfA^{\otimes, \degree})$ appearing in the statement of Proposition \ref{hurgoff2}. More concretely: 
\begin{itemize}
\item[$(i)$] The objects of $\widehat{\Cat}_{\infty}^{\otimes}$ are finite sequences
$[ X_1, \ldots, X_n ]$, where each $X_i$ is an $\infty$-category.
\item[$(ii)$] Given a pair of objects $[X_1, \ldots, X_n], [Y_1, \ldots, Y_m] \in \widehat{\Cat}_{\infty}^{\otimes}$, a morphism from $[X_1, \ldots, X_n]$ to $[Y_1, \ldots, Y_m]$ consists of a map 
$\alpha: \seg{n} \rightarrow \seg{m}$ and a collection of functors 
$\eta_{j}: \prod_{ \alpha(i) = j} X_{i} \rightarrow Y_j$.
\end{itemize}

Let $P$ denote the collection of all classes of simplicial sets, partially ordered by inclusion.
We define a subcategory $\calM$ of $\widehat{\Cat}_{\infty}^{\otimes} \times \Nerve(P)$ as follows:

\begin{itemize}
\item[$(iii)$] An object $( [X_1, \ldots, X_n], \calK)$ of $\widehat{\Cat}_{\infty}^{\otimes} \times \Nerve(P)$ belongs to $\calM$ if and only if each of the $\infty$-categories $X_i$ admits $\calK$-indexed colimits.

\item[$(iv)$] Let $f: ( [X_1, \ldots, X_n], \calK) \rightarrow ( [Y_1, \ldots, Y_m], \calK')$ be a morphism
in $\widehat{\Cat}_{\infty}^{\otimes} \times \Nerve(P)$, covering a map $\alpha: \seg{n} \rightarrow \seg{m}$ in $\FinSeg$. Then $f$ belongs to $\calM$ if and only if each of the associated functors
$$\prod_{ \alpha(i) = j } X_i \rightarrow Y_j$$
preserves $\calK$-indexed colimits separately in each variable.
\end{itemize}
\end{notation}

\begin{proposition}\label{successs}
In the situation of Notation \ref{surpint}, the canonical map
$p: \calM \rightarrow \Nerve(\FinSeg) \times \Nerve(P)$ is a coCartesian fibration.
\end{proposition}

\begin{proof}
We first show that $p$ is a locally coCartesian fibration.
Suppose given an object $( [ \calC_1, \ldots, \calC_n], \calK)$ in $\calM$, and a morphism
$\alpha: ( \seg{n}, \calK) \rightarrow ( \seg{m}, \calK')$ in $\Nerve( \FinSeg) \times \Nerve(P)$.
We wish to show that $\alpha$ can be lifted to a locally $p$-coCartesian morphism in $\calM$.
Supplying a lift of $\alpha$ is tantamount to choosing a collection of functors
$$ f_j: \prod_{ \alpha(i) = j } \calC_i \rightarrow \calD_j$$
for $1 \leq j \leq m$, such that $\calD_j$ admits $\calK'$-indexed colimits, and $f_j$ preserves $\calK$-indexed colimits separately in each variable. For $1 \leq i \leq n$, let $\calR_i$ denote the collection of
all colimit diagrams in $\calC_i$, indexed by simplicial sets belonging to $\calK$.
We now set $\calD_j = \calP_{\calR}^{\calK'}( \prod_{ \alpha(i) = j} \calC_i )$, where $\calR$ denotes the $\boxtimes$-product of $\{ \calR_i \}_{ \alpha(i) = j }$ (we refer the reader to \S \toposref{agileco} for an explanation of this notation). It follows from Proposition \toposref{cupper1} that the functors $\{ f_j \}_{1 \leq j \leq m}$ assemble to a locally $p$-coCartesian morphism in $\calM$. 

To complete the proof that $p$ is a coCartesian fibration, it will suffice to show that the locally $p$-coCartesian morphisms are closed under composition (Proposition \toposref{gotta}). In view of the construction of locally $p$-coCartesian morphisms given above, this follows immediately from Proposition \toposref{transit1}.
\end{proof}

\begin{remark}
Let $\calK$ be a collection of simplicial sets, and define
$\widehat{\Cat}_{\infty}^{\otimes}(\calK)$ denote the fiber product $\calM \times_{ \Nerve(P) } \{ \calK \};$ here $\calM$ and $P$ are defined as
in Notation \ref{surpint}. It follows from Proposition \ref{successs} that the projection
$\widehat{\Cat}_{\infty}^{\otimes}(\calK) \rightarrow \Nerve(\FinSeg)$ determines a symmetric monoidal structure on the $\infty$-category $$\widehat{\Cat}_{\infty}(\calK) = \calM \times_{ \Nerve(\FinSeg \times P) }
\{ ( \seg{1}, \calK) \}.$$ Here we may identify $\widehat{\Cat}_{\infty}(\calK)$ with the $\infty$-category of (not necessarily small) $\infty$-categories which admit $\calK$-indexed colimits, with morphisms given by functors which preserve $\calK$-indexed colimits. In the case $\calK = \emptyset$, this monoidal structure coincides with the Cartesian monoidal structure (see Notation \ref{surpint}). 
\end{remark}

The coCartesian fibration $p$ of Proposition \ref{successs} classifies a functor
$\Nerve(\FinSeg) \times \Nerve(P) \rightarrow \widehat{\Cat}_{\infty}$, which we may identify
with a functor from $\Nerve(P)$ to the $\infty$-category of commutative monoid objects of
$\widehat{\Cat}_{\infty}$. Consequently, we obtain a functor from $P$ to the $\infty$-category
of symmetric monoidal $\infty$-categories. In other words, if $\calK \subseteq \calK'$ are collections of simplicial sets, then we obtain a symmetric monoidal functor
from $\widehat{\Cat}_{\infty}( \calK)$ to $\widehat{\Cat}_{\infty}(\calK')$. It follows from the proof of Proposition \ref{successs} that this functor is given on objects by the formula
$\calC \mapsto \calP^{\calK'}_{\calK}(\calC)$.




\begin{remark}
Let $\calK \subseteq \calK'$ be collections of simplicial sets. We have an evident inclusion
of $\infty$-categories $\widehat{\Cat}_{\infty}^{\otimes}( \calK') \subseteq
\widehat{\Cat}_{\infty}^{\otimes}( \calK)$, which we can view as a lax symmetric monoidal
functor from $\widehat{\Cat}_{\infty}(\calK')$ to $\widehat{\Cat}_{\infty}(\calK)$. The $\infty$-category $\calM$ of Notation \ref{surpint} determines a correspondence which realizes this forgetful functor as the right adjoint to the symmetric monoidal functor $\calP^{\calK'}_{\calK}$. 
\end{remark}

Let $\calK$ be a collection of simplicial sets.
Recall that the $\infty$-category $\widehat{\Cat}_{\infty}^{\sMon}$ of symmetric monoidal $\infty$-categories is canonically equivalent to the $\infty$-category of commutative algebra objects 
$\CAlg( \widehat{\Cat}_{\infty})$ (Example \ref{interplace} and Proposition \ref{ungbatt}). 
Under this equivalence, the subcategory $\CAlg( \widehat{ \Cat}_{\infty}( \calK) )$ corresponds to the
subcategory of $\widehat{\Cat}_{\infty}^{\sMon}$ whose objects are symmetric monoidal $\infty$-categories $\calC$ which admit $\calK$-indexed colimits such that the tensor product
$\calC \times \calC \rightarrow \calC$ preserves $\calK$-indexed colimits separately in each variable, and whose morphisms are symmetric monoidal functors which preserve $\calK$-indexed colimits. Given an inclusion $\calK \subseteq \calK'$ of collections of simplicial sets, the induced inclusion
$$ \CAlg( \widehat{\Cat}_{\infty}( \calK') ) \subseteq \CAlg( \widehat{\Cat}_{\infty}(\calK) )$$ admits a left adjoint, given by composition with the symmetric monoidal functor $\calP^{\calK'}_{\calK}$.
Contemplating the unit of this adjunction, we arrive at the following result:

\begin{proposition}\label{skimmy}
Let $\calK \subseteq \calK'$ be collections of simplicial sets. Let $\calC$ be a
symmetric monoidal $\infty$-category. Assume that $\calC$ admits $\calK$-indexed colimits, and
that the tensor product $\otimes: \calC \times \calC \rightarrow \calC$ preserves $\calK$-indexed colimits separately in each variable. Then there exists a symmetric monoidal functor $\calC^{\otimes} \rightarrow \calD^{\otimes}$ with the following properties:
\begin{itemize}
\item[$(1)$] The $\infty$-category $\calD$ admits $\calK'$-indexed colimits. Moreover, the tensor product $\otimes: \calD \times \calD \rightarrow \calD$ preserves $\calK'$-indexed colimits separately in each variable.
\item[$(2)$] The underlying functor $f: \calC \rightarrow \calD$ preserves $\calK$-indexed colimits.
\item[$(3)$] For every $\infty$-category $\calE$ which admits $\calK'$-indexed colimits, composition with $f$ induces an equivalence of $\infty$-categories $\Fun_{\calK'}( \calD, \calE) \rightarrow \Fun_{\calK}( \calC, \calE)$. In particular, $f$ induces an identification $\calD \simeq \calP^{\calK'}_{\calK}(\calC)$, and is therefore fully faithful (Proposition \toposref{cupper1}). 
\item[$(4)$] If $\calE$ is a symmetric monoidal $\infty$-category which admits $\calK'$-indexed colimits, such that the tensor product $\otimes: \calE \times \calE \rightarrow \calE$ preserves $\calK'$-indexed colimits separately in each variable, then the induced map
$$ \Fun_{\calK'}^{\sMon}( \calD, \calE) \rightarrow \Fun_{\calK}( \calC, \calE)$$
is an equivalence of $\infty$-categories.
\end{itemize}
\end{proposition}

Specializing to the case where $\calK = \emptyset$ and $\calK'$ is the class of all small simplicial sets, we deduce the existence of the following analogue of the Day convolution product:

\begin{corollary}\label{banpan}
Let $\calC$ be a small symmetric monoidal $\infty$-category. Then there exists a symmetric monoidal structure on the $\infty$-category $\calP(\calC)$ of presheaves on $\calC$. It is characterized up to (symmetric monoidal) equivalence by the following properties:
\begin{itemize}
\item[$(1)$] The Yoneda embedding $j: \calC \rightarrow \calP(\calC)$ can be extended to a symmetric monoidal functor.
\item[$(2)$] The tensor product $\otimes: \calP(\calC) \times \calP(\calC) \rightarrow \calP(\calC)$ preserves small colimits separately in each variable.
\end{itemize}
\end{corollary}

\begin{proof}
The existence of the desired symmetric monoidal structures on $\calP(\calC)$ and $j$ follows from
Proposition \ref{skimmy}. The uniqueness follows from the universal property given in assertion $(3)$ of
Proposition \ref{skimmy}). 
\end{proof}

Applying the same argument to the case where $\calK = \emptyset$ and $\calK'$ is the class of all small {\em filtered} simplicial sets, we deduce the following:

\begin{corollary}\label{panban}
Let $\calC$ be a small symmetric monoidal $\infty$-category. Then there exists a symmetric monoidal structure on the $\infty$-category $\Ind(\calC)$ of $\Ind$-objects on $\calC$. It is characterized up to
(symmetric monoidal) equivalence by the following properties:
\begin{itemize}
\item[$(1)$] The Yoneda embedding $j: \calC \rightarrow \Ind(\calC)$ can be extended to a symmetric monoidal functor.
\item[$(2)$] The tensor product $\otimes: \Ind(\calC) \times \Ind(\calC) \rightarrow \Ind(\calC)$ preserves small filtered colimits separately in each variable.
\end{itemize}
Moreover, if $\calC$ admits finite colimits, and the tensor product $\otimes: \calC \times \calC \rightarrow \calC$ preserves finite colimits separately in each variable, then assertion $(2)$ can be strengthened as follows:
\begin{itemize}
\item[$(2')$] The tensor product functor $\otimes: \Ind(\calC) \times \Ind(\calC) \rightarrow \Ind(\calC)$
preserves small colimits separately in each variable.
\end{itemize}
\end{corollary}

\begin{proof}
Only assertion $(2')$ requires proof. For each $D \in \Ind(\calC)$, let $e_{D}$ denote the functor
$C \mapsto C \otimes D$. Let $\calD$ denote the full subcategory spanned by those objects
$D \in \Ind(\calC)$ such that the functor $e_{D}$ preserves small colimits separately in each variable.
We wish to prove that $\calD = \Ind(\calC)$.
Assertion $(2)$ implies that the correspondence $D \mapsto e_{D}$ is given by a
functor $\Ind(\calC) \rightarrow \Fun( \Ind(\calC), \Ind(\calC) )$ which preserves filtered colimits.
Consequently, $\calD$ is stable under filtered colimits in $\Ind(\calC)$. It will therefore suffice to show
that $\calD$ contains the essential image of $j$. 

Let $C \in \calC$. Since the Yoneda embedding
$j$ is a symmetric monoidal functor, we have a homotopy commutative diagram
$$ \xymatrix{ \calC \ar[r]^{e_C} \ar[d]^{j} & \calC \ar[d]^{j} \\
\Ind(\calC) \ar[r]^{e_{j(C)} } & \Ind(\calC). }$$
The desired result now follows from Proposition \toposref{sumatch}, since
$e_{j(C)}$ preserves filtered colimits and $j \circ e_{C}$ preserves finite colimits.
\end{proof}

\begin{remark}
The above analysis can be carried out without essential change in the context of (not necessarily symmetric) monoidal $\infty$-categories, and Corollaries \ref{banpan} and \ref{panban} have evident analogues in that setting.
\end{remark}

The $\infty$-category $\LPress$ of presentable $\infty$-categories can be identified with
a full subcategory of $\widehat{\Cat}_{\infty}(\calK)$, where $\calK$ is the collection of all small simplicial sets. 

\begin{proposition}\label{spurk}
Let $\calK$ denote the collection of all small simplicial sets. The $\infty$-category $\LPress$ of presentable $\infty$-categories is closed under the tensor product
$\otimes: \widehat{\Cat}_{\infty}(\calK) \otimes \widehat{\Cat}_{\infty}(\calK) \rightarrow \widehat{\Cat}_{\infty}(\calK)$, and therefore inherits a symmetric monoidal structure (Proposition \ref{yak2}).
\end{proposition}

\begin{proof}
Corollary \ref{banpan} implies that the unit object of $\widehat{\Cat}_{\infty}(\calK)$ is given
by $\calP( \Delta^0) \simeq \SSet$, which is a presentable $\infty$-category. It will therefore suffice
to show that if $\calC$ and $\calC'$ are two presentable $\infty$-categories, then their
tensor product $\calC \otimes \calC'$ in $\widehat{\Cat}_{\infty}(\calK)$ is likewise presentable.
If $\calC = \calP(\calC_0)$ and $\calC' = \calP( \calC'_0)$ for a pair of small $\infty$-categories
$\calC_0$ and $\calC'_0$, then Corollary \ref{banpan} yields an equivalence
$\calC \otimes \calC' \simeq \calP( \calC_0 \times \calC'_0)$, so that $\calC \otimes \calC'$ is presentable as desired. In view of Theorem \toposref{pretop}, every presentable $\infty$-category
is a localization of an $\infty$-category of presheaves. It will therefore suffice to prove the following:
\begin{itemize}
\item[$(\ast)$] Let $\calC$ and $\calC'$ be presentable $\infty$-categories, and assume that the tensor product $\calC \otimes \calC'$ is presentable. Let $S$ be a (small) set of morphisms in $\calC$. Then
the tensor product $S^{-1} \calC \otimes \calC'$ is presentable.
\end{itemize}

To prove $(\ast)$, we choose a (small) set of objects $M = \{ C'_{\alpha} \}$ which generates
$\calC'$ under colimits. Let $f: \calC \times \calC' \rightarrow \calC \otimes \calC'$ be the canonical
map, and let $T$ be the collection of all morphisms in $\calC \otimes \calC'$ having the form
$f( s \times \id_{C'} )$, where $s \in S$ and $C' \in M$. Consider now the composition
$$ g: S^{-1} \calC \times \calC' \subseteq \calC \times \calC' \stackrel{f}{\rightarrow}
\calC \otimes \calC' \stackrel{L}{\rightarrow} T^{-1} (\calC \otimes \calC'),$$
where $L$ is a left adjoint to the inclusion $T^{-1}(\calC \otimes \calC') \subseteq \calC \otimes \calC'$.
We claim that $g$ exhibits $T^{-1} (\calC \otimes \calC')$ as a tensor product of $S^{-1} \calC$ with
$\calC'$. In other words, we claim that if $\calD$ is an arbitrary $\infty$-category which admits small colimits, then composition with $g$ induces an equivalence
$$ \Fun_{\calK}( T^{-1}( \calC \otimes \calC'), \calD) \rightarrow \Fun_{\calK \boxtimes \calK}( S^{-1} \calC \times \calC', \calD).$$
The left hand side can be identified with the full subcategory of
$\Fun_{\calK}( \calC \otimes \calC', \calD)$ spanned by those functors which carry each morphism
in $T$ to an equivalence. Under the equivalence
$$\Fun_{\calK}( \calC \otimes \calC', \calD) \simeq \Fun_{\calK \boxtimes \calK}( \calC \times \calC', \calD)
\simeq \Fun_{\calK}( \calC, \Fun_{\calK}(\calC', \calD) ),$$
this corresponds to the full subcategory spanned by those functors
$F: \calC \rightarrow \Fun_{\calK}(\calC', \calD)$ which carry each morphism in $S$ to an equivalence.
This $\infty$-category is equivalent to $\Fun_{\calK}( S^{-1} \calC, \Fun_{\calK}( \calC', \calD) )
\simeq \Fun_{\calK \boxtimes \calK}( S^{-1} \calC \times \calC', \calD)$, as desired.
\end{proof}

\begin{remark}
Proposition \ref{spurk} yields a symmetric monoidal structure
$\MonLPres \rightarrow \Nerve(\FinSeg)$ on the $\infty$-category $\LPress$ of presentable $\infty$-categories. By restriction, we obtain also a monoidal structure on $\LPress$. Unwinding the definitions, we see that this monoidal structure coincides {\em up to isomorphism} with the monoidal structure constructed in \S \monoidref{jurmit}.
\end{remark}

\begin{notation}
Let $\MonStab \subseteq \MonLPres$ denote the full subcategory spanned by the objects
$[ X_1, \ldots, X_n]$ where each $X_{i}$ is a {\em stable} presentable $\infty$-category.
\end{notation}

\begin{remark}\label{sutyy}
Let $\LPressStab$ denote the full subcategory of $\LPress$ spanned by the stable $\infty$-categories.
In view of Corollary \stableref{mapprop}, the inclusion $\LPressStab \subseteq \LPress$ admits a left
adjoint $L: \LPress \rightarrow \LPressStab$, given by the formula
$$ \calC \mapsto \Stab(\calC).$$
In view of Example \monoidref{srower}, we can identify $L$ with the functor
$\calC \mapsto \calC \otimes \SSet_{\infty}$. It follows easily that $L$
is compatible with the symmetric monoidal structure on $\LPress$, in the sense of Definition \ref{compatsymon}.
\end{remark}

\begin{proposition}\label{tupik}
\begin{itemize}
\item[$(1)$] The projection $\MonStab \rightarrow \Nerve(\FinSeg)$ determines a symmetric monoidal structure on the $\infty$-category $\LPressStab$.
\item[$(2)$] The inclusion functor $i:\MonStab \subseteq \MonLPres$ is lax symmetric monoidal.
\item[$(3)$] The stabilization functor $\Stab: \LPress \rightarrow \LPressStab$ extends to a symmetric monoidal functor from $\MonLPres$ to $\MonStab$.
\item[$(4)$] The underlying monoidal structure on $\LPressStab$ agrees with the monoidal structure
described by Proposition \monoidref{tupid}.
\item[$(5)$] The unit object of $\LPressStab$ is the $\infty$-category $\Spectra$ of spectra.
\end{itemize}
\end{proposition}

\begin{proof}
Assertions $(1)$ through $(3)$ follow from Proposition \ref{localjerk2} and Remark \ref{sutyy}. Assertion $(4)$ is obvious, and $(5)$ follows from $(4)$ and the corresponding assertion of
Proposition \monoidref{tupid}.
\end{proof}

\begin{remark}\label{otherlandire}
Applying the reasoning which precedes Proposition \ref{skimmy}, we deduce that
$\CAlg(\LPressStab)$ can be identified with a full subcategory of $\CAlg( \LPressStab)$. A
symmetric monoidal $\infty$-category $\calC$ belongs to this full subcategory if and only if
$\calC$ is stable, presentable, and the bifunctor $\otimes: \calC \times \calC \rightarrow \calC$ preserves colimits separately in each variable.
\end{remark}

\begin{corollary}\label{surcoi2}
Let $\widehat{\Cat}_{\infty}^{\sMon}$ denote
the $\infty$-category of $($not necessarily small$)$ symmetric monoidal $\infty$-categories. Let
$\widehat{\Cat}_{\infty}^{\Stabb, \sMon}$ denote the subcategory whose objects are
required to be stable presentable monoidal $\infty$-categories such that the bifunctor
$\otimes$ preserves colimits in separately in each variable, and whose morphisms
are given by colimit-preserving symmetric monoidal functors. Then:

\begin{itemize}
\item[$(1)$] The $\infty$-category $\widehat{\Cat}_{\infty}^{\Stabb, \sMon}$ has an initial object $\calC^{\otimes}$.

\item[$(2)$] The underlying $\infty$-category $\calC^{\otimes}_{\seg{1}}$ is equivalent to the $\infty$-category of spectra.

\item[$(3)$] Let $\calD^{\otimes}$ be an arbitrary symmetric monoidal $\infty$-category. Suppose that:
\begin{itemize}
\item[$(i)$] The underlying $\infty$-category $\calD = \calD^{\otimes}_{\seg{1}}$ is stable and presentable.
\item[$(ii)$] The functor $\Spectra \rightarrow \calD$ determined by the unit object of $\calD$ $($see Corollary \stableref{choccrok}$)$ is an equivalence of $\infty$-categories.
\item[$(iii)$] The bifunctor $\otimes: \calD \times \calD \rightarrow \calD$ preserves small colimits in each variable.
\end{itemize}
Then there exists an equivalence of symmetric monoidal $\infty$-categories $\calC^{\otimes} \rightarrow \calD^{\otimes}$. Moreover, the collection of such equivalences is parametrized by a contractible Kan complex.
\end{itemize}
\end{corollary}

In other words, the $\infty$-category $\Spectra$ admits a symmetric monoidal structure, which is uniquely characterized by the following properties:

\begin{itemize}
\item[$(a)$] The bifunctor $\otimes: \Spectra \times \Spectra \rightarrow \Spectra$ preserves small colimits separately in each variable.
\item[$(b)$] The unit object of $\Spectra$ is the sphere spectrum $\Sphere$.
\end{itemize}

We will refer to this monoidal structure on $\Spectra$ as the {\it smash product symmetric monoidal structure}. Using the analogous uniqueness in the associative case (Corollary \monoidref{surcoi}), we deduce immediately that the monoidal structure on $\Spectra$ induces by the smash product symmetric monoidal structure is the monoidal structure constructed in \cite{monoidal}.

\begin{proof}[Proof of Corollary \ref{surcoi2}]
In view of Remark \ref{otherlandire}, we can identify $\widehat{\Cat}_{\infty}^{\Stabb, \sMon}$ with the $\infty$-category of commutative algebra objects of $\LPressStab$. Proposition \ref{gargle1} implies that $\widehat{\Cat}_{\infty}^{\Stabb, \sMon}$ has an initial object $\calC^{\otimes}$. This proves $(1)$. Moreover, Proposition \ref{gargle1} also asserts that the underlying $\infty$-category
$\calC = \calC^{\otimes}_{\seg{1}}$ is equivalent to the unit object of $\LPressStab$, which
is the $\infty$-category of spectra (Proposition \ref{tupik}); this proves $(2)$. 

Suppose that $\calD^{\otimes}$ is a symmetric monoidal $\infty$-category satisfying $(i)$, $(ii)$, and $(iii)$. Since $\calC^{\otimes}$ is an initial object of $\widehat{\Cat}_{\infty}^{\Stabb, \sMon}$, there exists a symmetric monoidal functor $F: \calC^{\otimes} \rightarrow \calD^{\otimes}$, unique up to a contractible space of choices. We claim that $F$ is an equivalence of symmetric monoidal $\infty$-categories. According to Remark \ref{symmetricegg}, it will suffice to show $F$ induces an equivalence
$f: \calC \rightarrow \calD$ on the level of the underlying $\infty$-categories.
Corollary \stableref{choccrok} implies that $f$ is determined, up to equivalence, by the image of the sphere spectrum $\Sphere \in \calC$. Since $\Sphere$ is the unit object of $\calC$, $f$ carries $\Sphere$ to the unit object of $\calD$. Condition $(ii)$ now implies that $f$ is an equivalence, as desired.
\end{proof}

\begin{warning}
Given a pair of objects $X,Y \in \h{\Spectra}$, the smash product of $X$ and $Y$ is usually denoted by $X \wedge Y$. We will depart from this convention by writing instead $X \otimes Y$ for the smash product. 
\end{warning}

\begin{remark}
According to Proposition \monoidref{discur}, the forgetful functor $\theta: \Mod^{L}_{ \Spectra}( \LPress ) \rightarrow \LPressStab$ is an equivalence (in other words, every presentable stable $\infty$-category is canonically tensored over $\Spectra$). Since $\Spectra$ has the structure of a commutative algebra object of $\LPress$, the results of
\S \ref{modcom} show that $\Mod^{L}_{\Spectra}( \LPress)$ 
 inherits the a symmetric monoidal structure. Under the equivalence $\theta$, this symmetric monoidal structure coincides (up to equivalence) with the symmetric monoidal structure of Proposition \ref{tupik}.
\end{remark}

\subsection{$\EInfty$-Ring Spectra}\label{comm8}

In this section, we finally come to the main object of study of this paper: the theory of commutative ring spectra, or {\it $E_{\infty}$-rings}. In this section, we will define the $\infty$-category of $E_{\infty}$-rings and establish some of its basic properties. In this, we will primarily follow the treatment of associative ring spectra given in \S \monoidref{monoid9}.

\begin{definition}
An {\it $E_{\infty}$-ring} is a commutative algebra object of the $\infty$-category of $\Spectra$ (endowed with its smash product monoidal structure). We let $\EInfty$ denote the $\infty$-category
$\CAlg(\Spectra)$ of $E_{\infty}$-rings.
\end{definition}

By construction, we have forgetful functors
$$ \CAlg( \Spectra) \rightarrow \Alg(\Spectra) \rightarrow \Spectra.$$
We will generally not distinguish in notation between $E_{\infty}$-ring $R$ and its underlying
$A_{\infty}$-ring, or its underlying spectrum. In particular, to any $E_{\infty}$-ring $R$ we can associate
a collection of homotopy groups $\{ \pi_{n} R \}_{n \in \Z }$. As explained in \S \monoidref{monoid9}, an
$A_{\infty}$-structure on $R$ endows the direct sum $\pi_{\ast} R = \bigoplus_{n \in \Z} \pi_{n} R$ with the structure of a graded ring. If $R$ is an $E_{\infty}$-ring, then we get a bit more: the multiplication
on $\pi_{\ast} R$ is {\it graded commutative}. That is, for $x \in \pi_{n} R$ and
$y \in \pi_{m} R$, we have $xy = (-1)^{nm} yx$. Here the sign results from the fact
that the composition
$$ \Sphere[n+m] \simeq \Sphere[n] \otimes \Sphere[m]
\stackrel{\sigma}{\simeq} \Sphere[m] \otimes \Sphere[n]  \simeq \Sphere[n+m]$$
is given by the sign $(-1)^{nm}$. In particular, the homotopy group $\pi_0 R$ is equipped with the structure of a discrete commutative ring, and every other homotopy group $\pi_n R$ has the structure of a module over $\pi_0 R$.

Let $\Omega^{\infty}: \Spectra \rightarrow \SSet$ be the $0$th space functor.
If $R$ is an $E_{\infty}$-ring, we will refer to $X=\Omega^{\infty} R$ as the {\it underlying space} of $R$. The underlying space $X$ is equipped with an addition $X \times X \rightarrow X$ (determined by the fact that it is the $0$th space of a spectrum) and a multiplication $X \times X \rightarrow X$ (determined by the map $R \otimes R \rightarrow R$); these maps endow $X$ with the structure of a commutative ring object in the homotopy category $\calH$ of spaces. However, the structure of an $E_{\infty}$-ring is much richer: not only do the ring axioms on $X$ hold up to homotopy, they hold up to {\em coherent} homotopy. 

The functor $\Omega^{\infty}: \EInfty \rightarrow \SSet$ is not conservative: a map of $E_{\infty}$-rings $f: A \rightarrow B$ which induces a homotopy equivalence of underlying spaces need not be an equivalence in $\EInfty$. We observe that $f$ is an equivalence of $E_{\infty}$-rings if and only if it is an equivalence of spectra; that is, if and only if $\pi_n(f): \pi_n A \rightarrow \pi_n B$ is an isomorphism of abelian groups for all $n \in \Z$. By contrast, $\Omega^{\infty}(f)$ is a homotopy equivalence of spaces provided only that
$\pi_{n}(f)$ is an isomorphism for $n \geq 0$; this is generally a weaker condition.

Recall that a spectrum $X$ is said to be {\it connective} if $\pi_{n} X \simeq 0$ for $n < 0$. We will say that an $E_{\infty}$-ring $R$ is {\it connective} if its underlying spectrum is connective. We let $\EInfty^{\conn}$ denote the full subcategory of $\EInfty$ spanned by the connective objects. Equivalently, we may view $\EInfty^{\conn}$ as the $\infty$-category $\CAlg( \connSpectra )$ of commutative algebra objects in connective spectra (the full subcategory $\connSpectra \subset \Spectra$ inherits a symmetric monoidal structure in view of Proposition \ref{yak2} and Lemma \monoidref{smashstab}).
When restricted to connective $E_{\infty}$-rings, the functor $\Omega^{\infty}$ detects equivalences: if $f: A \rightarrow B$ is a morphism in $\EInfty^{\conn}$ such that $\Omega^{\infty}(f)$ is an equivalence, then $f$ is an equivalence. We observe that
the functor $\Omega^{\infty}: \EInfty^{\conn} \rightarrow \SSet$ is a composition of a pair of functors
$ \CAlg( \connSpectra ) \rightarrow \connSpectra \rightarrow \SSet,$
both of which preserve geometric realizations (Corollary \ref{filtfemme} and \stableref{denkmal}) and admit left adjoints. It follows from Theorem \monoidref{barbeq} that $\EInfty^{\conn}$ can be identified with the $\infty$-category of modules over a suitable monad on $\SSet$. In other words, we can view connective $E_{\infty}$-rings as spaces equipped with some additional structures. Roughly speaking, these additional structures consist of an addition and multiplication which satisfy the axioms for a commutative ring, up to coherent homotopy.

\begin{definition}\label{caoncon}
Let $R$ be an $E_{\infty}$-ring. A {\it connective cover} of $R$ is a morphism
$\phi: R' \rightarrow R$ of $E_{\infty}$-rings with the following properties:
\begin{itemize}
\item[$(1)$] The $E_{\infty}$-ring $R'$ is connective.
\item[$(2)$] For every connective $A_{\infty}$-ring $E''$, composition with $\phi$ induces a homotopy equivalence
$$ \bHom_{ \EInfty}(R'', R') \rightarrow \bHom_{\EInfty}(R'',R). $$
\end{itemize}
\end{definition}

\begin{remark}
In the situation of Definition \ref{caoncon}, we will generally abuse terminology and simply refer to $R'$ as a connective cover of $R$, in the case where the map $\phi$ is implicitly understood.
\end{remark}

\begin{proposition}\label{canconexit}
\begin{itemize}
\item[$(1)$] Every $E_{\infty}$-ring $R$ admits a connective cover.
\item[$(2)$] An arbitrary map $\phi: R' \rightarrow R$ of $E_{\infty}$-rings is a connective cover of $R$ if and only if $R'$ is connective, and the induced map $\pi_{n} R' \rightarrow \pi_{n} R$ is an isomorphism for $n \geq 0$.
\item[$(3)$] The inclusion $\EInfty^{\conn} \subseteq \EInfty$ admits a right adjoint $G$, which carries each $E_{\infty}$-ring $R$ to a connective cover $R'$ of $R$.
\end{itemize}
\end{proposition}

\begin{proof}
Combine Lemma \monoidref{smashstab}, Proposition \ref{yak2} and Remark \ref{yukan}.
\end{proof}

Recall that an object $X$ of an $\infty$-category $\calC$ is said to be {\it $n$-truncated} if the
mapping spaces $\bHom_{\calC}(Y,X)$ are $n$-truncated, for every $Y \in \calC$ (see \S \toposref{truncintro}). Corollary \ref{algprice1} implies that the $\infty$-category  $\EInfty^{\conn}$ is presentable, so that we have a good theory of truncation functors.

\begin{lemma}\label{swerve}
Let $\calC$ be a symmetric monoidal $\infty$-category. Assume that $\calC$ admits countable colimits, and that the tensor product $\otimes: \calC \times \calC \rightarrow \calC$ preserves
countable colimits separately in each variable. An object $A \in \CAlg(\calC)$ is
$n$-truncated if and only if it is $n$-truncated as an object of $\calC$.
\end{lemma}

\begin{proof}
According to Corollary \ref{spaltwell}, the forgetful functor $\CAlg(\calC) \rightarrow \calC$
admits a left adjoint $F$. If $A$ is an $n$-truncated object of $\CAlg(\calC)$, then
for every $C \in \calC$, the space
$\bHom_{\calC}( C, A) \simeq \bHom_{ \CAlg(\calC)}( FC, A)$
is $n$-truncated. Thus $A$ is an $n$-truncated object of $\calC$.
Conversely, assume that $A$ is $n$-truncated in
$\calC$, and let $\CAlg'(\calC)$ denote the full subcategory of
$\CAlg(\calC)$ spanned by those objects $B$ such that
$\bHom_{\CAlg(\calC)}( B, A)$ is $n$-truncated. Since the collection of $n$-truncated
spaces is stable under limits, the full subcategory $\CAlg'(\calC)$ is stable under colimits
in $\CAlg(\calC)$. Moreover, the preceding argument shows that $\CAlg'(\calC)$ contains the essential image of the functor $F$. Using Proposition \monoidref{littlebeck}, we conclude that
$\CAlg'(\calC) = \CAlg(\calC)$, so that $A$ is $n$-truncated in $\CAlg(\calC)$.
\end{proof}

\begin{remark}
Using a slightly more involved argument, one can establish ``if'' direction of Lemma \ref{swerve} without assuming the existence of any colimits in $\calC$.
\end{remark}



\begin{proposition}\label{equitrick}
Let $R$ be a connective $E_{\infty}$-ring and let $n$ be a nonnegative integer. The following conditions are equivalent:
\begin{itemize}
\item[$(1)$] As an object of $\EInfty^{\conn}$, $R$ is $n$-truncated. 
\item[$(2)$] As an object of $\connSpectra$, $R$ is $n$-truncated.
\item[$(3)$] The space $\Omega^{\infty}(R)$ is $n$-truncated.
\item[$(4)$] For all $m > n$, the homotopy group $\pi_{m} R$ is trivial.
\end{itemize}
\end{proposition}

\begin{proof}
The equivalences $(2) \Leftrightarrow (3) \Leftrightarrow (4)$ follow from
Proposition \monoidref{equitrunc}, and the equivalence $(1) \Leftrightarrow (2)$
follows from Lemma \ref{swerve}.
\end{proof}

\begin{remark}
An $E_{\infty}$-ring $R$ is $n$-truncated as an object of $\EInfty$ if and only if it is equivalent to zero (since the $\infty$-category $\Spectra$ has no nontrivial $n$-truncated objects).
\end{remark}

Let $\tau_{\leq n}: \connSpectra \rightarrow \connSpectra$ be the truncation functor on connective spectra, and let $\tau^{\CAlg}_{\leq n}: \EInfty^{\conn} \rightarrow \EInfty^{\conn}$ be
the truncation functor on connective $E_{\infty}$-rings. Since the forgetful functor
$\theta: \EInfty^{\conn} \rightarrow \connSpectra$ preserves $n$-truncated objects, there is a canonical natural transformation
$\alpha: \tau_{\leq n} \circ \theta \rightarrow \theta \circ \tau_{\leq n}^{\CAlg}.$
Our next goal is to show that $\alpha$ is an equivalence.

\begin{proposition}\label{superjumpp}
Let $\calC$ be a symmetric monoidal $\infty$-category. Assume that $\calC$ is presentable
and that the tensor product $\otimes: \calC \times \calC \rightarrow \calC$ preserves colimits separately in each variable.
\begin{itemize}
\item[$(1)$] The localization functor $\tau_{\leq n}: \calC \rightarrow \calC$ is compatible with the symmetric monoidal structure on $\calC$, in the sense of Definition \ref{compatsymon}.
\item[$(2)$] The symmetric monoidal structure on $\calC$ induces a symmetric monoidal structure on $\tau_{\leq n}$ and an identification
$\CAlg( \tau_{\leq n} \calC ) \simeq \tau_{\leq n} \CAlg(\calC)$.
\end{itemize}
\end{proposition}

\begin{proof}
Assertion $(1)$ follows from Example \monoidref{jumper}. Assertion $(2)$ follows from $(1)$ and Lemma \ref{swerve}.
\end{proof}

Consequently, if $R$ is a connective $E_{\infty}$-ring, then the commutative algebra structure on $R$ determines a commutative algebra structure on $\tau_{\leq n} R$ for all $n \geq 0$. 

\begin{definition}
We will say that an $E_{\infty}$-ring is {\it discrete} if it is connective and $0$-truncated.
We let $\EInfty^{\disc}$ denote the full subcategory of $\EInfty$ spanned by the discrete objects.
\end{definition}

Since the mapping spaces in $\EInfty^{\disc}$ are $0$-truncated, it follows that
$\EInfty^{\disc}$ is equivalent to the nerve of an ordinary category. Our next result identifies this underlying category:

\begin{proposition}\label{umberhilt}
The functor $R \mapsto \pi_0 R$ determines an equivalence from
$\EInfty^{\disc}$ to the $($nerve of the$)$ category of commutative rings.
\end{proposition}

\begin{proof}
According to Proposition \ref{superjumpp}, we can identify $\EInfty^{\disc}$ with the
$\infty$-category of commutative algebra objects of the heart $\Spectralheart \subseteq \Spectra$, which inherits a symmetric monoidal structure from $\Spectra$ in view of Example \ref{swoon} and Lemma \monoidref{smashstab}. Proposition \stableref{specster1} allows us to identify $\Spectralheart$ with (the nerve of) the category of abelian groups. Moreover, the induced symmetric monoidal structure on $\Spectralheart$ has $\pi_0 \Sphere \simeq \Z$ as unit object, and the tensor product functor $\otimes$ preserves colimits separately in each variable. It follows that this symmetric monoidal structure coincides (up to canonical equivalence) with the usual symmetric monoidal structure on $\Spectralheart$, given by tensor products of abelian groups. Consequently, we may identify $\CAlg( \Spectralheart )$ with the (nerve of the) category of commutative rings.
\end{proof}

Let $A$ be an $E_{\infty}$-ring. Using the forgetful functor $\CAlg( \Spectra) \rightarrow \Alg(\Spectra)$, we see that $A$ can also be regarded as an $A_{\infty}$-ring, so that we can consider 
$\infty$-categories of left and right $R$-module spectra $\Mod_{A}^{L}(\Spectra)$ and $\Mod_{A}^{R}(\Spectra)$ studied in \S \monoidref{monoid11}. In view of Proposition \ref{socrime}, the commutativity of $A$ yields an equivalence between the theory of left and right $A$-modules; more precisely, we have a diagram of trivial Kan fibrations
$$ \Mod_{A}^{L}(\Spectra) \leftarrow \Mod^{\CommOp}_{A}(\Spectra) \rightarrow \Mod_{A}^{R}(\Spectra).$$
We will denote the $\infty$-category $\Mod^{\CommOp}_{A}(\Spectra)$ simply by $\Mod_{A}$. 

\begin{definition}
Let $R$ be an $E_{\infty}$-ring. We will say that an $R$-module $M$ is {\it flat} 
if it is flat when regarded as a left module over the underlying $A_{\infty}$-ring of $R$, in the sense of Definition \monoidref{defflat}. In other words, $M$ is flat if and only if the following conditions are satisfied:
\begin{itemize}
\item[$(1)$] The abelian group $\pi_0 M$ is flat as a module over the commutative ring
$\pi_0 R$.
\item[$(2)$] For every $n \in Z$, the canonical map
$$ \pi_{i} R \otimes_{ \pi_0 R} \pi_0 M \rightarrow \pi_{i} M$$ is an isomorphism of abelian groups.
\end{itemize}

We will say that a map $f: R \rightarrow R'$ of $E_{\infty}$-rings is {\it flat} if $R'$ is flat when regarded as an $R$-module. We let $\Mod^{\flat}_{R}$ denote the full subcategory of $\Mod_{R}$ spanned by
the flat $R$-module spectra.
\end{definition}

\begin{remark}\label{smacker}
Let $R$ be an $E_{\infty}$-ring. The $\infty$-category
$\Mod_{R}$ of $R$-module spectra admits a symmetric monoidal structure, given by relative tensor products over $R$ (Proposition \ref{stancer} and Theorem \ref{qwent}). If $M$ and $N$ are $R$-modules, then Proposition \monoidref{siiwe} yields a convergent spectral sequence 
$$ E_2^{p,q} = \Tor_{p}^{\pi_{\ast} R}( \pi_{\ast} M, \pi_{\ast} N)_{q}
\Rightarrow \pi_{p+q} ( M \otimes_{R} N ).$$
If $M$ or $N$ is flat over $R$, then $E_{2}^{p,q}$ vanishes for $p \neq 0$. It follows that this spectral sequence degenerates at the $E_{2}$-page, and we get a canonical isomorphism
$$ \pi_{q}( M \otimes_{R} N) \simeq E_{2}^{0,q}
\simeq \pi_{0} M \otimes_{ \pi_0 R} \pi_q R \otimes_{ \pi_{0} R} \pi_{0} N.$$
If $M$ and $N$ are {\em both} flat over $R$, then the relative tensor product $M \otimes_{R} N$ is again flat over $R$. Moreover, the unit object
$R \in \Mod_{R}$ clearly belongs to $\Mod_{R}^{\flat}$. 
Invoking Proposition \ref{yak2}, we conclude that 
the full subcategory $\Mod_{R}^{\flat} \subseteq \Mod_{R}$ inherits the structure of a symmetric monoidal category.
\end{remark}

\begin{remark}\label{smucker}
Combining Corollary \ref{skoke} with Remark \ref{smacker}, we conclude that for every
$E_{\infty}$-ring $R$, we have a fully faithful functor $\CAlg( \Mod_{R}^{\flat} ) \rightarrow
(\EInfty)^{R/}$ whose essential image consists those maps $R \rightarrow R'$ which exhibit
$R'$ as flat over $R$.
\end{remark}

\begin{remark}
Let $f: R \rightarrow R'$ be a flat map of $E_{\infty}$-rings. If $R$ is connective (discrete), then
$R'$ is also connective (discrete).
\end{remark}

\begin{proposition}\label{flatner}
Let $f: R \rightarrow R'$ be a map of $E_{\infty}$-rings. Suppose that $f$ induces an isomorphism
$\pi_{i} R \rightarrow \pi_{i} R'$ for all $i \geq 0$. Let
$\phi: (\EInfty)_{R'/} \rightarrow (\EInfty)_{R/}$ be the functor defined by composition with $f$. Then:
\begin{itemize}
\item[$(1)$] The functor $\phi$ admits a left adjoint $\psi$, given on objects by
$A \mapsto A \otimes_{R} R'$.
\item[$(2)$] The functor $\psi$ induces an equivalence
from the full subcategory of $(\EInfty)_{R/}$ spanned by the flat $R$-algebras to
the full subcategory of $( \EInfty)_{R'/}$ spanned by the flat $R'$-algebras.
\end{itemize}
\end{proposition}

\begin{proof}
Assertion $(1)$ follows immediately from Corollary \ref{skoke} and Proposition \ref{cocarten}.
To prove $(2)$, we combine Remark \ref{smucker} with Proposition \monoidref{urwise}.
\end{proof}

Our next goal is to study the $\infty$-category of commutative algebras over a fixed {\em connective} $E_{\infty}$-ring $R$. 

\begin{definition}
Let $R$ be a connective $E_{\infty}$-ring.
Let $\Mod_{R}^{\conn}$ denote the full subcategory of $\Mod_{R}$ spanned by the connective objects.
We will refer to an object of $\CAlg( \Mod_{R} )$ as an {\it commutative $R$-algebra}, and an object of 
$\CAlg( \Mod_{R}^{\conn} )$ as a {\it connective commutative $R$-algebra}.
We will say that a connective commutative $R$-algebra $A$ is:
\begin{itemize}
\item {\it finitely generated and free} if there exists a finitely generated free module $M$ and an equivalence $A \simeq \Sym^{\ast}(M)$ in $\CAlg(\Mod^{\conn}_{R} )$ (see Definition \monoidref{mofree}). We let $\CAlg( \Mod^{\conn}_{R} )^{\free}$ denote the full subcategory spanned by those modules which are free on finitely many generators. 
\item {\it of finite presentation} if $A$ belongs to the smallest full subcategory
of $\CAlg( \Mod^{\conn}_{R} )$ which contains $\CAlg( \Mod^{\conn}_{R} )^{\free}$
and is stable under finite colimits and retracts.
\item {\it almost of finite presentation} if $A$ is a compact object of $\tau_{\leq n} \CAlg( \Mod^{\conn}_{R} )$, for all $-2 \leq n < \infty$.
\end{itemize}
\end{definition}

The basic properties of these notions are summarized in the following result:

\begin{proposition}\label{baseprops}
Let $R$ be a connective $E_{\infty}$-ring. Then:
\begin{itemize}
\item[$(1)$] The $\infty$-category $\CAlg( \Mod^{\conn}_{R} )$ is projectively generated $($Definition \toposref{defpro}$)$. Moreover, a connective $R$-algebra $A$ is a compact projective object
of $\CAlg( \Mod^{\conn}_{R} )$ if and only if $A$ is a retract of a finitely generated free commutative $R$-algebra.

\item[$(2)$] Let $\CAlg( \Mod^{\conn}_{R})^{\finite} \subseteq \CAlg( \Mod^{\conn}_{R} )$ denote the smallest full subcategory of $\CAlg( \Mod^{\conn}_{R})$ which contains all finitely generated free algebras, and is stable under finite colimits. Then the preceding inclusion induces an equivalence $\Ind( \CAlg( \Mod^{\conn}_{R})^{\finite} ) \simeq \CAlg( \Mod^{\conn}_{R} )$.

\item[$(3)$] The $\infty$-category $\CAlg( \Mod^{\conn}_{R})$ is compactly generated. The compact objects of $\CAlg( \Mod^{\conn}_{R})$ are precisely the finitely presented $R$-algebras.

\item[$(4)$] An $R$-algebra $A$ is almost of finite presentation if and only if, for every $n \geq 0$, there exists an $R$-algebra $A'$ of finite presentation such that $\tau_{ \leq n} A$ is a retract of
$\tau_{\leq n} A'$ in $\CAlg( \Mod_{R}^{\conn} )$. 
\end{itemize}
\end{proposition}

\begin{proof}
Assertion $(1)$ from Corollaries \monoidref{chadwick} and \monoidref{progen}, where the latter is applied to the forgetful functor $\CAlg( \Mod^{\conn}_{R}) \rightarrow \Mod^{\conn}_{R}$.
Assertion $(2)$ then follows from Proposition \toposref{uterr}, $(3)$ from Lemma \toposref{stylus}, and $(4)$ from Corollary \toposref{hunterygreen}.
\end{proof}

\begin{remark}
Suppose given a pushout diagram
$$ \xymatrix{ R \ar[r] \ar[d] & R' \ar[d] \\
A \ar[r] & A' }$$
in the $\infty$-category of connective $E_{\infty}$-rings. If $A$ is
free and finitely generated (finitely presented, almost finitely presented), as an $R$-algebra, then $A'$ is free and finitely generated (finitely presented, almost finitely presented). In the first three cases, this is obvious. The last two
follow from assertions $(3)$ and $(4)$ of Proposition \ref{baseprops} (together with the observation that $\tau_{\leq n} A' \simeq \tau_{\leq n} ( \tau_{\leq n} A \otimes_{R} R' )$ ). 
\end{remark}

\begin{remark}\label{unkid}
Suppose given a commutative diagram
$$ \xymatrix{ & B \ar[dr] & \\
A \ar[rr] \ar[ur] & & C }$$
of $E_{\infty}$-rings, where $B$ is of finite presentation over $A$. Then $C$ is of finite presentation over $B$ if and only if $C$ is of finite presentation over $A$. This follows immediately from Propositions \ref{baseprops} and \toposref{accessforwardslice}. In \cite{deformation}, we will prove the analogue of this statement for morphisms which are almost of finite presentation.
\end{remark}

We now introduce some absolute finiteness properties of $E_{\infty}$-rings.

\begin{definition}
Let $R$ be a connective $E_{\infty}$-ring. We will say that
$R$ is {\it coherent} if it is left coherent when regarded as an $A_{\infty}$-ring (Definition \monoidref{pokus}). In other words, $R$ is coherent if $\pi_0 R$ is a coherent ring (in the sense of ordinary commutative algebra), and each $\pi_i R$ is a finitely presented module over $\pi_0 R$.
We will say that $R$ is {\it Noetherian} if $R$ is coherent and the ordinary commutative ring
$\pi_0 R$ is Noetherian.
\end{definition}

We close this section with the following result, which relates the absolute and relative finiteness conditions introduced above:

\begin{proposition}[Hilbert Basis Theorem]
Let $f: R \rightarrow R'$ be a map of connective $E_{\infty}$-rings. Suppose that
$R$ is Noetherian. Then $R'$ is almost of finite presentation as an $R$-algebra if and only if the following conditions are satisfied:
\begin{itemize}
\item[$(1)$] The ring $\pi_{0} R'$ is a finitely generated commutative $\pi_0 R$-algebra.
\item[$(2)$] The $E_{\infty}$-ring $R'$ is Noetherian.
\end{itemize}
\end{proposition}

\begin{proof}
We first prove the ``only if'' direction. We first note that $\pi_0 R' = \tau_{\leq 0} R'$ is a compact object
in the ordinary category of commutative $\pi_0 R$-algebras. This proves $(1)$. The classical
Hilbert basis theorem implies that $\pi_0 R'$ is Noetherian. It remains only to show that each
$\pi_{n} R'$ is a finitely generated module over $\pi_ 0 R'$. This condition depends only on
the truncation $\tau_{\leq n} R'$. In view of part $(4)$ of Proposition \ref{baseprops}, we may assume that $R'$ is a finitely presented commutative $R$-algebra.

Let $\calC$ be the full subcategory of $\CAlg( \Mod_{R}^{\conn} )$ spanned by those commutative $R$-algebras which satisfy $(1)$ and $(2)$. To show that $\calC$ contains all finitely presented $R$-algebras, it will suffice to show that $\calC$ is stable under finite colimits, and contains the free commutative $R$-algebra on a single generator. We first prove the stability under finite colimits. In view of Corollary \toposref{allfin}, it will suffice to show that $R \in \calC$ and that $\calC$ is stable under pushouts. The inclusion $R \in \calC$ is obvious (since $R$ is Noetherian by assumption). Suppose given a pushout diagram
$$ \xymatrix{ A \ar[r] \ar[d] & A' \ar[d] \\
A'' \ar[r] & A' \otimes_{A} A'' }$$
in $\CAlg( \Mod_{R}^{\conn} )$, where $A, A', A'' \in \calC$. We wish to prove that
$A' \otimes_{A} A'' \in \calC$. According to Corollary \monoidref{huty}, the commutative ring
$B = \pi_0 (A' \otimes_{A} A'')$ is canonically isomorphic to the classical tensor product
$\pi_0 A' \otimes_{ \pi_0 A} \pi_0 A''$, and is therefore a finitely generated commutative $\pi_0 R$-algebra. Thus $A' \otimes_{A} A''$ satisfies $(1)$.
Moreover, Proposition \monoidref{siiwe} yields a convergent spectral sequence
$$ E_2^{p,q} = \Tor_{p}^{\pi_{\ast} A}( \pi_{\ast} A', \pi_{\ast} A'')_{q}
\Rightarrow \pi_{p+q} ( A' \otimes_{A} A'' ).$$
It is not difficult to see that this is a spectral sequence of modules over the commutative ring $B$.
Each of the $B$-modules $E_{2}^{p,q}$ is a finitely generated $B$-module. It follows that
the subquotients $E_{\infty}^{p,q}$ are likewise finitely generated over $B$. Consequently,
each homotopy group $\pi_{n} (A' \otimes_{A} A'')$ has a finite filtration by finitely generated $B$-modules, and is therefore itself finitely generated over $B$. This proves that $A' \otimes_{A} A''$ satisfies $(2)$, and therefore belongs to the $\infty$-category $\calC$.

Let $R[X] = \Sym^{\ast}(R)$ denote the free commutative $R$-algebra on a single generator. We wish to show that $R[X] \in \calC$. We first treat the case where $R$ is discrete. Let $k$ denote the commutative ring $\pi_0 R$. Remark \ref{tappus} implies that we have canonical isomorphisms
$$ \pi_{m} R[X] \simeq \pi_{m} ( \coprod_{n} \Sym^{n}(R) )
\simeq \oplus_{n \geq 0} \HH_{m}( \Sigma_{n}; k). $$
In particular, $\pi_0 R[X]$ is canonically isomorphic to the polynomial ring $k[x]$, which proves $(1)$. To prove $(2)$, we must show that each of the groups $\oplus_{n \geq 0} \HH_{m}( \Sigma_{n}; k)$ is a finitely generated $k[x]$-module. Since the homology of a finite group can be computed using a finite complex, each of the groups $\HH_{m}( \Sigma_{n}; k)$ is individually a finite $k$-module. Moreover, multiplication by $x$ is given by the maps
$$ \eta_{m,n}: \HH_{m}( \Sigma_{n}; k) \rightarrow \HH_{m}( \Sigma_{n+1}; k)$$
induced by the inclusions $\Sigma_{n} \subseteq \Sigma_{n+1}$. To complete the proof, it will suffice
to show that for fixed $m$, $\eta_{m,n}$ is surjective for $n \gg 0$. This follows from Nakaoka's homological stability theorem for the symmetric groups (for a simple proof, we refer to the reader to
\cite{kerz}).

Let us now treat the general case. Let $B$ denote the polynomial ring $\pi_0 R[X] = (\pi_0 R)[X]$. Then $B$ is finitely generated over $\pi_0 R$, so $R[X]$ satisfies $(1)$. We wish to prove that $R[X]$ satisfies $(2)$. Assume otherwise, and choose $m$ minimal such that $\pi_{m} R[X]$ is not a finitely generated $B$-module; note that $m$ is necessarily positive.
Set $R_0 = \tau_{\leq 0}$. Then the argument above proves that the tensor product
$R_0 \otimes_{R} R[X] \simeq R_0[X]$ belongs to $\calC$. Invoking Proposition \monoidref{siiwe} again, we obtain a convergent spectral sequence of $B$-modules 
$$ E_2^{p,q} = \Tor_{p}^{\pi_{\ast} R}( \pi_{0} R, \pi_{\ast} (R[X]))_{q}
\Rightarrow \pi_{p+q} (R_0[X]).$$
Using the minimality of $m$, we deduce that $E_{2}^{p,q}$ is finitely generated over $B$ for all $q < m$. It follows that
$E_{\infty}^{0,m}$ is the quotient of $E_{2}^{0,m}$ by a finitely generated submodule. Since
$E_{\infty}^{0,m}$ is a submodule of the finitely generated $B$-module $\pi_{m} (R_0[X])$, we conclude that $E_{2}^{0,m}$ is itself finitely generated over $B$. But $E_{2}^{0,m}$ contains
$\pi_{m}(R[X])$ as a summand, contradicting our choice of $m$. This completes the proof of the ``only if'' direction.

Now suppose that $R'$ satisfies $(1)$ and $(2)$. We will construct a sequence of
maps $\psi_{n}: A_{n} \rightarrow R'$ in $\CAlg( \Mod_{R}^{\conn})$ with the property that
the maps $\psi_{n}$ induce isomorphisms $\pi_{m} A_{n} \rightarrow \pi_{m} R'$ for
$m < n$ and surjections $\pi_{n} A_{n} \rightarrow \pi_{n} R'$. To construct $A_0$, we invoke
assumption $(1)$: choose a finite set of elements $x_1, \ldots, x_i \in \pi_{0} R'$ which
generate $\pi_0 R'$ as an $\pi_0 R$-algebra. These elements determine a map of $R$-algebras
$\psi_0: \Sym^{\ast}( R^{i} ) \rightarrow R'$ with the desired property.

Let us now suppose that $\psi_{n}: A_{n} \rightarrow R'$ has already been constructed. 
Let $K$ denote the fiber of the map $\psi_{n}$, formed in the $\infty$-category of
$A_{n}$-modules. Our assumption on $\psi_{n}$ implies that
$\pi_{i} K \simeq 0$ for $i < n$. Let $B = \pi_0 A_{n}$. Since $A_{n}$ is almost of finite presentation over $R$, the first part of the proof shows that $B$ is a Noetherian ring and that each
of the homotopy groups $\pi_{i} A_{n}$ is a finitely generated $B$-module. Moreover, since the map
$B \rightarrow \pi_0 R'$ is surjective and $R'$ satisfies $(2)$, we conclude that each
of the homotopy groups $\pi_{i} R'$ is a finitely generated $B$-module. Using the long exact sequence
$$ \ldots \rightarrow \pi_{n+1} R' \rightarrow \pi_{n} K \rightarrow \pi_{n} A_{n} \rightarrow \ldots $$
we deduce that $\pi_{n} K$ is a finitely generated $B$-module. Consequently, there exists
a finitely generated free $A_{n}$-module $M$ and a map $M[n] \rightarrow K$ which is surjective on $\pi_{n}$. Let $f: \Sym^{\ast}( M[n] ) \rightarrow A_{n}$ be the induced map, and form a pushout diagram
$$ \xymatrix{ \Sym^{\ast}( M[n] ) \ar[r]^{f} \ar[d]^{f_0} & A_{n} \ar[d]^{f_0'} \\
A_{n} \ar[r] & A_{n+1}, }$$
where $f_0$ classifies the zero from $M[n]$ to $A_{n}$. By construction, 
we have a canonical homotopy from $\psi_{n} \circ f$ to $\psi_{n} \circ f_0$, which
determines a map $\psi_{n+1}: A_{n+1} \rightarrow R'$. We observe that there is a distinguished triangle of $A_{n}$-modules
$$ \ker(f'_0) \rightarrow  \ker( \psi_{n} ) \rightarrow \ker(\psi_{n+1}) \rightarrow \ker(f'_0)[1].$$
To show that $\psi_{n+1}$ has the desired properties, it suffices to show that
$\pi_{i} \ker( \psi_{n+1} ) \simeq 0$ for $i \leq n$. Using the long exact sequence associated to the distinguished triangle above, we may reduce to proving the following pair of assertions:
\begin{itemize}
\item[$(a)$] The homotopy groups $\pi_{i} \ker(f'_0)$ vanish for $i < n$.
\item[$(b)$] The canonical map $\pi_{n} \ker(f'_0) \rightarrow \pi_n K$ is surjective.
\end{itemize}

We now observe that there is an equivalence $\ker( f'_0) \simeq \ker(f_0) \otimes_{ \Sym^{\ast}( M[n] ) } A_{n}$. In view of Corollary \monoidref{huty}, it will suffice to prove the same assertions after replacing $f'_0$ by $f_0$. Using Remark \ref{tappus}, we obtain a canonical equivalence
$$ \ker( f_0 ) \simeq \coprod_{ m > 0} \Sym^{m}( M[n] ).$$
Because the $\infty$-category $\Mod_{A_{n}}^{\conn}$ is closed under colimitsin $\Mod_{A_{n}}$, it follows that the homotopy groups $\pi_{i} \Sym^{m}( M[n] )$ vanish for 
$i < mn$. This proves $(a)$. To prove $(b)$, it will suffice to show that the composite map
$$ M[n] \simeq \Sym^{1}( M[n] ) \rightarrow
\coprod_{m > 0} \Sym^{m}( M[n] ) \rightarrow K$$
induces a surjection on $\pi_{n}$, which follows immediately from our construction.
\end{proof}

\subsection{Symmetric Monoidal Model Categories}\label{comm9}

If $\bfA$ is a simplicial model category, then the collection of fibrant-cofibrant objects
of $\bfA$ can be organized into an $\infty$-category $\Nerve( \bfA^{\degree})$. In this section,
we will show that that $\bfA$ is equipped with a symmetric monoidal structure (which is compatible with the model structure and the simplicial structure in a suitable sense), then the underlying
$\infty$-category $\Nerve( \bfA^{\degree})$ inherits the structure of a symmetric monoidal
$\infty$-category (Proposition \ref{hurgoff2}). Moreover, under suitable assumptions on $\bfA$ (see Definition \ref{powerdef}), the category $\CAlg(\bfA)$ of commutative algebras in $\bfA$ inherits a simplicial model structure (Proposition \ref{combocheese}) whose underlying $\infty$-category
$\Nerve( \CAlg(\bfA)^{\degree})$ is equivalent to the $\infty$-category of
commutative algebra objects $\CAlg( \Nerve(\bfA^{\degree}) )$ (Theorem \ref{cbeckify}). This can be regarded as a rectification theorem: an object of $\CAlg( \Nerve(\bfA^{\degree}))$ can be regarded
as an object $A \in \bfA$ equipped with a multiplication $A \otimes A \rightarrow A$ which
is commutative and associative (and unital) up to coherent homotopy, and Theorem \ref{cbeckify} guarantees that $A$ is weakly equivalent to an algebra $A'$ which is strictly commutative and associative.

We begin with some very general considerations.

\begin{notation}\label{capin}
If $\calO$ is a simplicial colored operad (see Variation \ref{col} on Definition \ref{colop}), we let $\calO^{\otimes}$ denote the simplicial category given by Construction \ref{coco}: 
\begin{itemize}
\item[$(i)$] The objects of
$\calO^{\otimes}$ are pairs $(\seg{n}, (C_1, \ldots, C_n) )$, where
$\seg{n} \in \FinSeg$ and $C_1, \ldots, C_n$ are colored of $\calO$.
\item[$(ii)$] Given a pair of objects $C = ( \seg{m}, (C_1, \ldots, C_m))$ to
$C' = (\seg{n}, (C'_1, \ldots, C'_n))$ in $\calO^{\otimes}$, the simplicial set
$\bHom_{\calO^{\otimes}}(C, C')$ is defined to be
$$\coprod_{ \alpha: \seg{m} \rightarrow \seg{n}} \prod_{ 1 \leq j \leq n} \Mul_{\calO}( \{ C_i \}_{\alpha(i) = j}, C'_j ).$$
\item[$(iii)$] Composition in $\calO^{\otimes}$ is defined in the obvious way.
\end{itemize}
\end{notation}

\begin{definition}
Let $\calO$ be a simplicial colored operad. We will denote the simplicial nerve of the
category $\calO^{\otimes}$ by $\Nerve^{\otimes}(\calO)$; we will refer to
$\Nerve^{\otimes}(\calO)$ as the {\it operadic nerve of $\calO$}.
\end{definition}

\begin{remark}
Let $\calO$ be a simplicial colored operad. Then there is an evident forgetful functor
from $\calO^{\otimes}$ to the ordinary category $\FinSeg$ (regarded as a simplicial category).
This forgetful functor induces a canonical map $\Nerve^{\otimes}( \calO) \rightarrow \Nerve(\FinSeg)$. 
We may therefore regard $\Nerve^{\otimes}$ as a functor from the category $\SCop$ of simplicial
colored operads to the category $(\sSet)_{/ \Nerve(\FinSeg)}$.
\end{remark}

\begin{remark}
Let $\calO$ be a simplicial colored operad. The fiber product $\Nerve^{\otimes}(\calO)
\times_{ \Nerve(\FinSeg) } \{ \seg{1} \}$ is canonically isomorphic to the nerve of
the simplicial category underlying $\calO$; we will denote this simplicial set
by $\Nerve(\calO)$.
\end{remark}

\begin{definition}
We will say that a simplicial colored operad $\calO$ is {\it fibrant} if each of the simplicial sets
$\Mul_{\calO}( \{ X_i \}_{i \in I}, Y)$ is fibrant.
\end{definition}

\begin{proposition}\label{calp}
Let $\calO$ be a fibrant simplicial colored operad. Then the operadic nerve
$\Nerve^{\otimes}(\calO)$ is an $\infty$-operad.
\end{proposition}

\begin{proof}
If $\calO$ is a fibrant simplicial colored operad, then $\calO^{\otimes}$ is a fibrant simplicial category
so that $\Nerve^{\otimes}(\calO)$ is an $\infty$-category. Let 
$C = (\seg{m}, (C_1, \ldots, C_m))$ be an object of $\Nerve^{\otimes}(\calC)$ and let
$\alpha: \seg{m} \rightarrow \seg{n}$ be an inert morphism in $\FinSeg$. Then
we have a canonical map $C \rightarrow C' = (\seg{n}, (C_{ \alpha^{-1}(1)},
\ldots, C_{ \alpha^{-1}(n)}))$ in $\calO^{\otimes}$, which we can identify with an
edge $\overline{\alpha}$ of $\Nerve^{\otimes}(\calO)$ lying over $\alpha$.
Using Proposition \toposref{trainedg}, we deduce that $\overline{\alpha}$ is
$p$-coCartesian, where $p: \Nerve^{\otimes}(\calO) \rightarrow \Nerve(\FinSeg)$.

As a special case, we observe that there are $p$-coCartesian morphisms
$\overline{\alpha}^{i, \nostar{m}}: C \rightarrow (\seg{1}, C_i)$ covering
$\Colp{i}: C \rightarrow C_i$ for $1 \leq i \leq m$. To prove that $\Nerve^{\otimes}(\calO)$ is an $\infty$-operad, we must show that these maps determine a $p$-limit diagram
$\nostar{m}^{\triangleleft} \rightarrow \Nerve^{\otimes}(\calO)$. Unwinding the definitions, we
must show that for every object $D = ( \seg{n}, (D_1, \ldots, D_n))$ and every
morphism $\beta: \seg{n} \rightarrow \seg{m}$, the canonical map
$$ \bHom_{\calO^{\otimes}}^{\beta}( D, C) \rightarrow
\prod_{1 \leq i \leq m} \bHom_{\calO^{\otimes}}^{\Colp{i} \circ \beta}( D, (\seg{1}, C_i))$$
is a homotopy equivalence; here $\bHom_{\calO^{\otimes}}^{\beta}( D, C)$ denotes the inverse
image of $\{ \beta \}$ in $\bHom_{\calO^{\otimes}}(D,C)$ and 
$\bHom_{\calO^{\otimes}}^{\Colp{i} \circ \beta}( D, (\seg{1}, C_i))$ the inverse
image of $\{ \Colp{i} \circ \beta \}$ in $\bHom_{\calO^{\otimes}}( D, (\seg{1}, C_i))$.
We now observe that this map is an isomorphism of simplicial sets.

To complete the proof that $\Nerve^{\otimes}(\calO)$ is an $\infty$-operad, it suffices to show that
for each $m \geq 0$ the functors $\Colll{i}{!} $ associated to $p$ induces an
essentially surjective map
$$ \Nerve^{\otimes}(\calO) \times_{ \Nerve(\FinSeg)} \{ \seg{m} \} \rightarrow
\prod_{1 \leq i \leq m} \Nerve(\calO).$$
In fact, $\Nerve^{\otimes}(\calO) \times_{ \Nerve(\FinSeg) } \{ \seg{m} \}$ is
canonically isomorphic with $\Nerve(\calO)^{m}$ and this isomorphism identifies the above map with the identity.
\end{proof}

\begin{definition}\label{hirt}
Let $\calC$ be simplicial category. We will say that a symmetric monoidal structure on $\calC$ is {\it weakly compatible} with the simplicial structure on $\calC$ provided that
the operation $\otimes: \calC \times \calC \rightarrow \calC$ is endowed with the structure of a simplicial functor, which is compatible with associativity, commutativity, and unit constraints of $\calC$.

Suppose furthermore that $\calC$ is {\it closed}: that is, for every pair objects $X, Y \in \calC$,
there exists an exponential object $X^{Y} \in \calC$ and an evaluation map
$e: X^{Y} \otimes Y \rightarrow X$ which induces bijections
$$ \Hom_{\calC}( Z, X^{Y} ) \rightarrow \Hom_{\calC}( Z \otimes Y \rightarrow X)$$
for every $Z \in \calC$. We will say that the symmetric monoidal structure on $\calC$ is {\it compatible} with a simplicial structure on $\calC$ if it is weakly compatible, and the map $e$ induces isomorphisms of simplicial sets
$$ \bHom_{\calC}( Z, X^Y ) \rightarrow \bHom_{\calC}( Z \otimes Y, X) $$
for every $Z \in \calC$.
\end{definition}

\begin{remark}
As with Definition \monoidref{hyrt}, the compatibility of Definition \ref{hirt} is not merely a condition, but additional data.
\end{remark}

\begin{remark}\label{kil}
Let $\calC$ be a symmetric monoidal category equipped with a weakly compatible simplicial struct.re
We can then regard $\calC$ as a simplicial colored operad by setting
$$ \Mul_{\calC}( \{ X_i \}, Y) = \bHom_{\calC}( \otimes_{i} X_i, Y).$$
We let $\calC^{\otimes}$ be the simplicial category obtained by applying
the construction of Notation \ref{capin} to this simplicial colored operad.
\end{remark}

\begin{proposition}\label{hurgoff}
Let $\calC$ be a fibrant simplicial category with a weakly compatible symmetric monoidal structure, and let $\calC^{\otimes}$ be the simplicial category of Remark \ref{kil}. The forgetful functor
$p: \Nerve( \calC^{\otimes} ) \rightarrow \Nerve(\FinSeg)$
determines a symmetric monoidal structure on the simplicial nerve $\Nerve(\calC)$. 
\end{proposition}

\begin{proof}
Proposition \ref{calp} guarantees that $p$ exhibits $\Nerve( \calC^{\otimes})$ as an $\infty$-operad.
To complete the proof, it will suffice to show that $p$ is a coCartesian fibration. It is obviously an inner fibration, since
$\Nerve( \calC^{\otimes})$ is an $\infty$-category and $\Nerve( \FinSeg)$ is the nerve of an ordinary category. Choose an object
$(C_1, \ldots, C_{n}) \in \calC^{\otimes}$ and a morphism $\alpha: \seg{n} \rightarrow \seg{m}$ in $\FinSeg$.
Then there exists a morphism $\overline{\alpha}: (C_1, \ldots, C_{n}) \rightarrow (D_{1}, \ldots, D_{m})$ in $\Nerve( \calC^{\otimes})$ which covers $\alpha$, which induces an isomorphism
$D_{j} \simeq \otimes_{\alpha(i) = j} C_{i}$ for each $j \in \nostar{m}$.
It follows immediately from Proposition \toposref{trainedg} that $\overline{\alpha}$ is $p$-coCartesian. \end{proof}

In practice, fibrant simplicial categories often arise as the categories of fibrant-cofibrant objects in simplicial model categories. In these cases, the hypotheses of Proposition \ref{hurgoff} are unnaturally strong. If $\bfA$ is a simplicial model category equipped with a compatible symmetric monoidal structure (Definition \ref{hirs} below), then the full subcategory of fibrant-cofibrant objects is usually not stable under tensor products. Nevertheless, a mild variant of Proposition \ref{hurgoff} can be used to endow $\Nerve( \bfA^{\degree})$ with the structure of a symmetric monoidal $\infty$-category.

\begin{definition}\label{hirs}
A {\it symmetric monoidal model category} is a model category $\bfA$ equipped with a symmetric monoidal structure which satisfies the following conditions:
\begin{itemize}
\item[$(1)$] For every pair of cofibrations $i: A \rightarrow A'$, $j: B \rightarrow B'$ in $\bfA$, the induced map
$$k: (A \otimes B') \coprod_{ A \otimes B} (A' \otimes B) \rightarrow A' \otimes B'$$
is a cofibration. Moreover, if either $i$ or $j$ is a weak equivalence, then $k$ is a weak equivalence.

\item[$(2)$] The unit object $1$ of $\bfA$ is cofibrant.

\item[$(3)$] The symmetric monoidal structure on $\bfA$ is closed (see Definition \ref{hirt}).
\end{itemize}

A {\it simplicial symmetric monoidal model category} is a simplicial model category
$\bfA$ equipped with a symmetric monoidal structure satisfying $(1)$, $(2)$, and $(3)$, which is compatible with the simplicial structure in the sense of Definition \ref{hirt}.
\end{definition}

\begin{remark}
One can reformulate Definition \ref{hirs} as follows. Let $\bfA$ be a symmetric monoidal model category. Then a compatible simplicial structure on $\bfA$ can be identified with
symmetric monoidal left Quillen functor $\psi: \sSet \rightarrow \bfA$. Given such a functor, 
$\bfA$ inherits the structure of a simplicial category, where the mapping spaces are characterized by the existence of a natural bijection
$$ \Hom_{\sSet}(K , \bHom_{\bfA}(A, B) ) \simeq \Hom_{\bfA}( \psi(K) \otimes A, B).$$
\end{remark}

\begin{proposition}\label{hurgoff2}
Let $\bfA$ be a simplicial symmetric monoidal model category and let
$\bfA^{\otimes}$ be defined as in Remark \ref{kil}. Let $\bfA^{\degree}$ be the full subcategory of $\bfA$ spanned by the fibrant-cofibrant objects, and let $\bfA^{\otimes,\degree}$ be the full subcategory of $\bfA^{\otimes}$ spanned by those objects $(C_1, \ldots, C_{n} )$ such that each
$C_i$ belongs to $\bfA^{\degree}$. Then the natural map
$p: \Nerve( \bfA^{\otimes,\degree} ) \rightarrow \Nerve( \FinSeg)$
determines a symmetric monoidal structure on the $\infty$-category $\Nerve( \bfA^{\degree} )$.
\end{proposition}

\begin{proof}
Our first step is to prove that $\Nerve( \bfA^{\otimes, \degree})$ is an $\infty$-category.
To prove this, it will suffice to show that $\bfA^{\otimes, \degree}$ is a fibrant simplicial category. 
Let $(C_1, \ldots, C_{n})$ and $(D_1, \ldots, D_{m})$ be objects of $\bfA^{\otimes, \degree}$. Then
$$\bHom_{ \bfA^{\otimes, \degree} }( (C_1, \ldots, C_{n}), (D_1, \ldots, D_{m}) )$$
is a disjoint union of products of simplicial sets of the form
$\bHom_{\bfA}( C_{i_1} \otimes \ldots \otimes C_{i_k}, D_j)$. Each of these simplicial sets
is a Kan complex, since $C_{i_1} \otimes \ldots \otimes C_{i_k}$ is cofibrant and $D_j$ is fibrant.

Since $p$ is a map from an $\infty$-category to the nerve of an ordinary category, it is automatically an inner fibration. We next claim that $p$ is a coCartesian fibration. Let
$(C_1, \ldots, C_{n})$ be an object of $\bfA^{\otimes, \degree}$, and let 
$\alpha: \seg{n} \rightarrow \seg{m}$ be a map in $\FinSeg$. For each $j \in \nostar{m}$,
choose a trivial cofibration
$$\eta_j: \bigotimes_{ \alpha(i) = j } C_{i} \rightarrow D_j,$$
where $D_{j} \in \bfA$ is fibrant. Together these determine a map $\overline{\alpha}: (C_1, \ldots, C_{n}) \rightarrow (D_1, \ldots, D_{m})$ in $\Nerve( \bfA^{\otimes, \degree})$. We claim that $\overline{\alpha}$ is $p$-coCartesian. In view of Proposition \toposref{trainedg}, it will suffice to show that for every morphism $\beta: \seg{m} \rightarrow \seg{k}$ in $\FinSeg$ and every 
$(E_1, \ldots, E_{k} ) \in \bfA^{\otimes, \degree}$, the induced map
$$\bHom_{ \bfA^{\otimes} }( (D_1, \ldots, D_{m}) , (E_1, \ldots, E_{k}) )
\times_{ \Hom_{\FinSeg}( \seg{m}, \seg{k} ) } \{ \beta \} 
\rightarrow \bHom_{ \bfA^{\otimes} }( (C_1, \ldots, C_{n}), (E_1, \ldots, E_{k}))$$
determines a homotopy equivalence onto the summand of
$\bHom_{ \bfA^{\otimes} }( (C_1, \ldots, C_{n}), (E_{1}, \ldots, E_{k}) )$ spanned
by those morphisms which cover the map $\beta \circ \alpha: \seg{n} \rightarrow \seg{k}$.
Unwinding the definitions, it will suffice to prove that for $l \in \nostar{k}$, the induced map
$$ \bHom_{\bfA}( \bigotimes_{ \beta(j) = l } D_{j}, E_{l} )
\rightarrow \bHom_{\bfA}( \bigotimes_{ (\beta \circ \alpha)(i) = l} C_{i} , E_{l} )$$
is a homotopy equivalence of Kan complexes. Since each $E_{l}$ is a fibrant object of 
$\bfA$, it will suffice to show that the map
$$\eta: \bigotimes_{ (\beta \circ \alpha)(i) = l } C_{i} \rightarrow \bigotimes_{ \beta(j) = l} D_j$$
is a weak equivalence of cofibrant in $\bfA$. This follows from the observation that $\eta$ can be identified with the tensor product of the maps $\{ \eta_{j} \}_{g(i-1) < j \leq g(i) }$, each of which is individually a weak equivalence between cofibrant objects.

We now conclude the proof by observing that the fiber $\Nerve( \bfA^{\otimes, \degree} )_{\seg{n}}$
is isomorphic to an $n$-fold product of $\Nerve( \bfA^{\degree})$ with itself, and that the projection onto the $j$th factor can be identified with the functor associated to the map
$\Colp{j-1}: \seg{n} \rightarrow \seg{1}$.
\end{proof}

\begin{example}\label{exalcut}
Let $\bfA$ be a simplicial model category. Suppose that the Cartesian symmetric monoidal structure on $\bfA$ is compatible with the model structure (in other words, that the final object of $\bfA$ is cofibrant, and that for any pair of cofibrations $i: A \rightarrow A'$, $j: B \rightarrow B'$, the induced map
$i \wedge j: (A \times B') \coprod_{A \times B} (A' \times B) \rightarrow A' \times B'$ is a cofibration, trivial if either $i$ or $j$ is trivial). Then there is a canonical map of simplicial categories
$\theta: \bfA^{\otimes} \rightarrow \bfA$, given on objects by the formula
$ \theta(A_1, \ldots, A_{n}) = A_1 \times \ldots \times A_{n}.$
Since the collection of fibrant-cofibrant objects of $\bfA$ is stable under finite products,
$\theta$ induces a map $\bfA^{\otimes, \degree} \rightarrow \bfA^{\degree}$. 
It is easy to see that $\theta$ induces a Cartesian structure on the symmetric monoidal $\infty$-category $\Nerve( \bfA^{\otimes, \degree} )$. Consequently, 
the induced symmetric monoidal structure on $\Nerve( \bfA^{\degree} )$ is Cartesian.
\end{example}

\begin{remark}
Let $\bfA$ be a simplicial symmetric monoidal model category, so that the underlying $\infty$-category $\Nerve(\bfA^{\degree})$ inherits the a symmetric monoidal structure (Proposition \ref{hurgoff2}). The underlying monoidal structure on $\Nerve( \bfA^{\degree})$ agrees with the monoidal structure described by Proposition \monoidref{hurgoven} (the associative analogue of Proposition \ref{hurgoff2}). 
\end{remark}

Let $\bfA$ be a simplicial model category equipped with a compatible symmetric monoidal structure.
Our goal is to show that, up to equivalence, every commutative algebra object $A$ of $\Nerve( \bfA^{\degree})$ arises from a {\em strict} commutative algebra in $\bfA$. To prove this, we will need some rather strong assumptions about the behavior of tensor products in $\bfA$.

\begin{notation}\label{puffstar}
Let $\bfA$ be a symmetric monoidal category which admits colimits.
Given a pair of morphisms $f: A \rightarrow A'$, $g: B \rightarrow B'$, we let
$f \wedge g$ denote the induced map
$$ (A \otimes B') \coprod_{ A \otimes B} (A' \otimes B) \rightarrow A' \otimes B'.$$
We observe that the operation $\wedge$ determines a symmetric monoidal structure on the category of morphisms in $\bfA$. In particular, for every morphism $f: X \rightarrow Y$, we iterate the above construction to obtain a map
$$ \wedge^{n}(f): \bbox^{n}(f) \rightarrow Y^{\otimes n}.$$
Here the source and target of $\wedge^{n}(f)$ carry actions of the symmetric group $\Sigma_{n}$, and
$\wedge^{n}(f)$ is a $\Sigma_{n}$-equivariant map. Passing to $\Sigma_{n}$-coinvariants, we obtain
a new map, which we will denote by $\sigma^{n}(f): \Sym^{n}(Y; X) \rightarrow \Sym^{n}(Y).$
\end{notation}

Before giving the next definition, we need to review a bit of terminology. 
Recall that a collection $S$ of morphisms in a presentable category $\bfA$ is {\it weakly saturated} if it is stable under pushouts, retracts, and transfinite composition (see Definition \toposref{saturated}).
For every collection $S$ of morphisms in $\bfA$, there is a smallest weakly saturated collection of morphisms $\overline{S}$ containing $S$. In this case, we will say that $\overline{S}$ is generated by $S$.

\begin{definition}\label{powerdef}
Let $\bfA$ be a combinatorial symmetric monoidal model category. We will say that a morphism
$f: X \rightarrow Y$ is a {\it power cofibration} if the following condition is satisfied:
\begin{itemize}
\item[$(\star)$] For every $n \geq 0$, the induced map $\wedge^{n}(f): \bbox^{n}(f) \rightarrow Y^{\otimes n}$ is a cofibration in $\bfA^{\Sigma_{n}}$. Here $\bfA^{\Sigma_{n}}$ denotes the category of
objects of $\bfA$ equipped with an action of the symmetric group $\Sigma_{n}$, endowed with the {\em projective} model structure (see \S \toposref{quasilimit3}).  
\end{itemize}
We will say that an object $X \in \bfA$ is {\it power cofibrant} if the map $\emptyset \rightarrow X$ is a power cofibration, where $\emptyset$ is an initial object of $\bfA$. 
We will say that $\bfA$ is {\it freely powered} if the following conditions are satisfied:
\begin{itemize}
\item[$(F1)$] The symmetric monoidal model category $\bfA$ satisfies the {\it monoid axiom} of \cite{monmod}. That is, if $\overline{U}$ denotes the weakly saturated class of morphisms generated by morphisms of the form
$$ \id_{X} \otimes f: X \otimes Y \rightarrow X \otimes Y',$$
where $X$ is arbitrary and $f$ is a trivial cofibration, then every morphism in $\overline{U}$ is
a weak equivalence.
\item[$(F2)$] The model category $\bfA$ is left proper, and the collection of cofibrations in $\bfA$ is generated (as a weakly saturated class) by cofibrations between cofibrant objects.
\item[$(F3)$] Every cofibration in $\bfA$ is a power cofibration.
\end{itemize}
\end{definition}

\begin{remark}
In the situation of Definition \ref{powerdef}, if every object of $\bfA$ is cofibrant, then
$\bfA$ automatically satisfies $(F1)$ and $(F2)$.
\end{remark}

\begin{remark}
Let $\bfA$ be a combinatorial symmetric monoidal model category. Then every power cofibration in $\bfA$ is a cofibration (take $n=1$). It is not difficult to show that the collection of power cofibrations in $\bfA$ is weakly saturated. Consequently, to show that $\bfA$ satisfies $(F3)$, it suffices to show that there
is a set of generating cofibrations for $\bfA$ which consists of power cofibrations.
\end{remark}

\begin{remark}\label{spuksy}
If $f: X \rightarrow Y$ is a cofibration in $\bfA$, then the definition of a monoidal model category guarantees that $\wedge^{n}(f)$ is a cofibration {\em in $\bfA$}, which is trivial if $f$
is trivial. Condition $(\star)$ is much stronger: roughly speaking, it guarantees that the symmetric group $\Sigma_{n}$ acts {\em freely} on the object $Y^{\otimes n}$ (see Lemma \ref{bullow} below).
\end{remark}

We now turn to the study of commutative algebras in a symmetric monoidal simplicial model category $\bfA$. Let $\CAlg(\bfA)$ denote the category whose objects are commutative algebras in $\bfA$.
The category $\CAlg(\bfA)$ inherits the structure of a simplicial category from $\bfA$: namely, we regard $\CAlg(\bfA)$ as cotensored over simplicial sets using the observation that any commutative
algebra structure on an object $A \in \bfA$ induces a commutative algebra structure on
$A^{K}$ for every simplicial set $K$.

The main results of this section can be stated as follows:

\begin{proposition}\label{combocheese}
Let $\bfA$ be a symmetric monoidal model category. Assume that $\bfA$ is combinatorial and freely powered. Then:
\begin{itemize}
\item[$(1)$] The category $\CAlg(\bfA)$ admits a combinatorial model structure, where:
\begin{itemize}
\item[$(W)$] A morphism $f: A \rightarrow B$ of commutative algebra objects of $\bfA$ is a weak equivalence if it is a weak equivalence when regarded as a morphism in $\bfA$.
\item[$(F)$] A morphism $f: A \rightarrow B$ of commutative algebra objects of $\bfA$ is a fibration if it is a fibration when regarded as a morphism in $\bfA$.
\end{itemize}
\item[$(2)$] The forgetful functor $\theta: \CAlg(\bfA) \rightarrow \bfA$ is a right Quillen functor.
\item[$(3)$] If $\bfA$ is a simplicial symmetric monoidal model category, then $\CAlg(\bfA)$ inherits the structure of a simplicial model category.
\end{itemize}
\end{proposition}

\begin{theorem}\label{cbeckify}
Let $\bfA$ be a simplicial symmetric monoidal model category.
Assume that $\bfA$ is combinatorial and freely powered. Then the canonical map
$$ \Nerve(\CAlg(\bfA)^{\degree}) \rightarrow \CAlg( \Nerve(\bfA^{\degree}) )$$
is an equivalence of $\infty$-categories.
\end{theorem}

We will prove Proposition \ref{combocheese} and Theorem \ref{cbeckify} at the end of this section.

\begin{example}\label{fielding}
Let $k$ be a field and let $\bfA$ denote the category of complexes of $k$-vector spaces
$$ \ldots \rightarrow M_{n+1} \rightarrow M_{n} \rightarrow M_{n-1} \rightarrow \ldots. $$
The category $\bfA$ has a symmetric monoidal structure, given by tensor products of complexes (here the symmetry isomorphism $M_{\bigdot} \otimes N_{\bigdot} \simeq N_{\bigdot} \otimes M_{\bigdot}$ involves the insertion of a sign).
The category $\bfA$ also admits the structure of a combinatorial simplicial model category, as explained in Example \monoidref{fielder} and S \stableref{stable10}:
\begin{itemize}
\item[$(C)$] A map of complexes $f: M_{\bigdot} \rightarrow N_{\bigdot}$ is a {\it cofibration} if it
induces an injection $M_{n} \rightarrow N_{n}$ of $k$-vector spaces, for each
$n \in \Z$.
\item[$(F)$] A map of complexes $f: M_{\bigdot} \rightarrow N_{\bigdot}$ is a {\it fibration} if it
induces an surjection $M_{n} \rightarrow N_{n}$ of $k$-vector spaces, for each
$n \in \Z$.
\item[$(W)$] A map of complexes $f: M_{\bigdot} \rightarrow N_{\bigdot}$ is a {\it weak equivalence} if it is a quasi-isomorphism; that is, if $f$ induces an isomorphism on homology groups
$H_{n}( M_{\bigdot}) \rightarrow H_{n}( N_{\bigdot})$ for each $n \in \Z$.
\end{itemize}
It follows from Proposition \ref{hurgoff2} that the underlying $\infty$-category $\Nerve(\bfA^{\degree})$ (which we can identify with the (unbounded) {\it derived $\infty$-category of $k$-vector spaces} $\calD(k)$) inherits a symmetric monoidal structure.

Every object of $\bfA$ is cofibrant. Moreover, {\em if $k$ is of characteristic zero}, then $\bfA$ satisfies axiom $(F3)$ of Definition \ref{powerdef}. To see this, let us consider an arbitrary cofibration $f: X \rightarrow Y$ in $\bfA$.
We observe that the induced map $\wedge^{n}(f): \bbox^{n}(f) \rightarrow Y^{\otimes n}$ is automatically a cofibration in $\bfA$ (Remark \ref{spuksy}). Moreover, the category
$\bfA^{\Sigma_{n}}$ can be identified with the category of complexes of modules over the group algebra $k[ \Sigma_{n} ]$, equipped with the {\it projective} model structure (so that fibrations and weak equivalences are as defined above). If $k$ is of characteristic zero, then the group algebra
$k[ \Sigma_{n} ]$ is semisimple. It follows that the projective and injective model structure on $\bfA^{\Sigma_{n} }$ coincide, so that every monomorphism of complexes of $k[ \Sigma_{n} ]$-modules is automatically a cofibration in $\bfA^{\Sigma_n}$.
Theorem \ref{cbeckify} implies that, when $k$ is of characteristic zero, the $\infty$-category
$\CAlg( \calD(k) )$ is equivalent to the $\infty$-category underlying the model category
$\CAlg( \bfA)$ of {\it commutative differential graded $k$-algebras}.

With slightly more effort, one can show that the {\em projective} model structure on
$\bfA$ is freely powered whenever $k$ is a commutative algebra over the field $\Q$ of rational numbers. However, the results of this section fail drastically outside of characteristic zero. If $k$ is a field of positive characteristic, then it is not even possible to construct the model structure on $\CAlg(\bfA)$ described in Proposition \ref{combocheese}. If such a model structure did exist, then the free commutative algebra functor $F: \bfA \rightarrow \CAlg(\bfA)$ would be a left Quillen functor. However, the functor $F$ does not preserve weak equivalences between cofibrant objects: in characteristic $p$, the $p$th symmetric power of an acyclic complex need not be acyclic.
\end{example}

\begin{example}\label{bulwump}
Let $\bfA$ be the category of {\it symmetric spectra}, as defined in \cite{symmetricspectra}. 
We will regard $\bfA$ as endowed with the {\it positive $S$-model structure} described in
\cite{shipley}. Then $\bfA$ is a simplicial symmetric monoidal model category, which is combinatorial and freely powered. Using Corollary \ref{surcoi2}, we deduce that $\Nerve(\bfA^{\degree})$ is equivalent, as a symmetric monoidal $\infty$-category, to
the $\infty$-category $\Spectra$ of spectra (endowed with the smash product monoidal structure; see \S \ref{comm7}). Using Theorem \ref{cbeckify}, we deduce that our $\infty$-category
$\EInfty = \CAlg( \Spectra)$ of $E_{\infty}$-rings is equivalent to the $\infty$-category underlying the category $\CAlg(\bfA)$ of {\it commutative symmetric ring spectra}.
\end{example}

\begin{remark}
Let $\calC$ be a symmetric monoidal $\infty$-category. In order for $\calC$ to arise from the situation described in Proposition \ref{hurgoff2}, it is necessary for $\calC$ to be presentable and for the
tensor product $\otimes: \calC \times \calC \rightarrow \calC$ to preserve colimits separately in each variable. It seems likely that these conditions are also sufficient.
\end{remark}

We now turn to the proof of Proposition \ref{combocheese}. First, we need a few preliminaries.

\begin{lemma}\label{fulltime}
Let $\bfA$ be a combinatorial symmetric monoidal model category. For every category $\calC$, let $\bfA^{\calC}$ denote the associated diagram category, endowed with the projective model structure (see \S \toposref{quasilimit3}). 
Let $f$ be a cofibration in $\bfA^{\calC}$, and $g$ a cofibration in $\bfA^{\calC'}$. Then
the smash product
$f \wedge g$ (see Notation \ref{puffstar}) is a cofibration in $\bfA^{ \calC \times \calC'}$.
In particular, if $X \in \bfA^{\calC}$ and $Y \in \bfA^{\calC'}$ are cofibrant, then
$X \otimes Y$ is a cofibrant object of $\bfA^{\calC \times \calC'}$.
\end{lemma}

\begin{proof}
Let $S$ denote the collection of all morphisms $f$ in $\bfA^{\calC}$ for which the conclusion of the Lemma holds. It is not difficult to see that $S$ is weakly saturated, in the sense of Definition \toposref{saturated}. Consequently, it will suffice to prove that $S$ contains a set of generating cofibrations for 
$\bfA^{\calC}$. Let $i^{C}: \{ \ast \} \rightarrow \calC$ be the inclusion of an object $C \in \calC$, and let
$i^{C}_{!}: \bfA \rightarrow \bfA^{\calC}$ be the corresponding left Kan extension functor (a left adjoint to the evaluation at $C$). Then the collection of cofibrations in $\bfA^{\calC}$ is generated by morphisms of the form $i^{C}_{!}(f_0)$, where $f_0$ is a cofibration in $\bfA$. We may therefore assume that $f = i^{C}_{!}(f_0)$. Using the same argument, we may assume that
$g = i^{C'}_{!}(g_0)$, where $C' \in \calC$ and $g_0$ is a cofibration in $\bfA$. We now observe that
$f \wedge g$ is isomorphic to $i^{(C,C')}_{!} ( f_0 \wedge g_0 )$. Since $i^{(C,C')}_{!}: \bfA \rightarrow \bfA^{\calC \times \calC' }$ is a left Quillen functor, it will suffice to show that $f_0 \wedge g_0$ is a cofibration in $\bfA$, which follows from our assumption that $\bfA$ is a (symmetric) monoidal model category.
\end{proof}

\begin{lemma}\label{bullow}
Let $\bfA$ be a combinatorial symmetric monoidal model category. Then:
\begin{itemize}
\item[$(1)$] Let $f: X \rightarrow Y$ be a power cofibration in $\bfA$. Then the induced map
$$ \sigma^{n}(f): \Sym^{n}(Y;X) \rightarrow \Sym^{n}(X)$$
is a cofibration, which is weak equivalence if $f$ is a weak equivalence $($see Notation \ref{puffstar}$)$ and $n > 0$.

\item[$(2)$] Let $Y$ be a power cofibrant object of $\bfA$. Then $\Sym^{n}(Y)$ coincides with
the homotopy colimit of the action of $\Sigma_{n}$ on $Y^{\otimes n}$.

\item[$(3)$] Let $f: X \rightarrow Y$ be a power cofibration between power-cofibrant objects
of $\bfA$. Then the object $\bbox^{n}(f) \in \bfA^{ \Sigma_{n} }$ is cofibrant $($with respect to the projective model structure$)$.
\end{itemize}
\end{lemma}

\begin{proof}
Let $F: \bfA^{\Sigma_{n} } \rightarrow \bfA$ be a left adjoint to the diagonal functor, so that
$F$ carries an object $X \in \bfA^{\Sigma_{n}}$ to the object of coinvariants $X_{\Sigma_{n}}$. 
We observe that $F$ is a left Quillen functor, and that $\sigma^{n}(f) = F( \wedge^{n}(f) )$ for every morphism in $\bfA$. Assertion $(1)$ now follows immediately from the definitions (and the observation that $\wedge^{n}(f)$ is a weak equivalence if $n > 0$ and $f$ is a trivial cofibration). Assertion $(2)$ follows immediately from the definition of a homotopy colimit.

We now prove $(3)$. Let $f: X \rightarrow Y$ be a power cofibration between power cofibrant objects of $\bfA$. We observe that $Y^{\otimes n}$ admits a $\Sigma_{n}$-equivariant filtration
$$ X^{\otimes n} = Z_0 \stackrel{\gamma_1}{\rightarrow} Z_1 \stackrel{\gamma_2}{\rightarrow}
\ldots \stackrel{\gamma_{n-1}}{\rightarrow} Z_{n-1} = \bbox^{n}(f) \rightarrow
Z_{n} = Y^{\otimes n}.$$
It will therefore suffice to prove that each of the maps $\{ \gamma_{i} \}_{ 1 \leq i \leq n-1 }$ is a cofibration in $\bfA^{ \Sigma_{n} }$. For this, we observe that there is a pushout diagram
$$ \xymatrix{ \pi_{!} ( \bbox^{i}(f) \times X^{ \otimes (n-i) } )
\ar[rr]^{\pi_{!} ( \wedge^{i}(f) \otimes \id)} \ar[d] & &  \pi_{!} (Y^{\otimes i} \otimes X^{\otimes n-i}) \ar[d] \\
Z_{i-1} \ar[rr]^{\gamma_i} & & Z_{i}, }$$
where $\pi_{!}$ denotes the left adjoint to the forgetful functor
$\bfA^{\Sigma_{n}} \rightarrow \bfA^{ \Sigma_{i} \times \Sigma_{n-i} }$. Since $\pi_{!}$ is a left Quillen functor, it suffices to show that
each $\wedge^{i}(f) \otimes \id$ is a cofibration in $\bfA^{ \Sigma_{i} \times \Sigma_{n-i} }$. This
follows from our assumption that $X$ is power cofibrant, our assumption that $f$ is a power cofibration, and Lemma \ref{fulltime}.
\end{proof}

\begin{proof}[Proof of Proposition \ref{combocheese}]
We first observe that the category $\CAlg(\bfA)$ is presentable (this is a special case of Corollary \ref{algprice1}). Since $\bfA$ is combinatorial, there exists a (small) collection of morphisms $I = \{ i_{\alpha}: C \rightarrow C' \}$ which generates the class of cofibrations in $\calC$, and a (small) collection of morphisms $J = \{ j_{\alpha}: D \rightarrow D' \}$ which generates the class of trivial cofibrations in $\calC$.

Let $F: \bfA \rightarrow \CAlg(\bfA)$ be a left adjoint to the forgetful functor. Let $\overline{F(I)}$ be the weakly saturated class of morphisms in $\CAlg(\bfA)$ generated by $\{ F(i): i \in I \}$, and let $\overline{F(J)}$ be defined similarly. Unwinding the definitions, we see that a morphism in $\CAlg(\bfA)$ is a trivial fibration if and only if it has the right lifting property with respect to $F(i)$, for every $i \in I$. Invoking the small object argument, we deduce that every morphism $f: A \rightarrow C$ in $\CAlg(\bfA)$ admits a factorization
$ A \stackrel{f'}{\rightarrow} B \stackrel{f''}{\rightarrow} C$ where $f' \in \overline{F(I)}$ and
$f''$ is a trivial fibration. Similarly, we can find an analogous factorization where
$f' \in \overline{F(J)}$ and $f''$ is a fibration.

Using a standard argument, we may reduce the proof of $(1)$ to the problem of showing that every morphism belonging to $\overline{ F(J) }$ is a weak equivalence in $\CAlg(\bfA)$. 
Let $\overline{U}$ be as in Definition \ref{powerdef}, and let $S$ be the collection of morphisms
in $\CAlg(\bfA)$ such that the underlying morphism in $\bfA$ belongs to $\overline{U}$.
Since $\bfA$ satisfies the monoid axiom, $S$ consists of weak equivalences in $\bfA$.
It will therefore suffice to show that $\overline{F(J)} \subseteq S$. Because $S$ is 
weakly saturated, it will suffice to show that $F(J) \subseteq S$. Unwinding the definitions, we are reduced to proving the following:
\begin{itemize}
\item[$(\ast)$] Let 
$$ \xymatrix{ F(C) \ar[r]^{F(i)} \ar[d] & F(C') \ar[d] \\
A \ar[r]^{f} & A' }$$ be a pushout diagram in $\CAlg(\bfA)$. If $i$ is a trivial cofibration in
$\bfA$, then $f \in S$.
\end{itemize}
To prove $(\ast)$, we observe that $F(C')$ admits a filtration by $F(C)$-modules
$$ F(C) \simeq B_0 \rightarrow B_1 \rightarrow B_2 \rightarrow \ldots, $$
where $F(C') \simeq \colim \{ B_i \}$ and for each $n > 0$ there is a pushout diagram
$$ \xymatrix{ F(C) \otimes \Sym^{n}(C';C) \ar[rr]^{ \id_{F(C)} \otimes \sigma^{n}(i) } \ar[d] & &  F(C) \otimes \Sym^{n}(C') \ar[d] \\
B_{n-1} \ar[rr] & & B_{n}. }$$
(See Notation \ref{puffstar}.)
It follows that $A'$ admits a filtration by $A$-modules
$$ A \simeq B'_0 \rightarrow B'_1 \rightarrow B'_2 \rightarrow \ldots,$$
where $A' \simeq \colim \{B'_i \}$ and for each $n > 0$ there is a pushout diagram
$$ \xymatrix{ A \otimes \Sym^{n}(C';C) \ar[rr]^{ \id_{A} \otimes \sigma^{n}(i) } \ar[d] & & A \otimes \Sym^{n}(C') \ar[d] \\ B'_{n-1} \ar[rr] & &  B'_{n}. }$$
Since $i$ is a trivial power cofibration, Lemma \ref{bullow} implies that $\sigma^{n}(i)$ is a trivial cofibration. It follows that $\id_{A} \otimes \sigma^{n}(i)$ belongs to $\overline{U}$.
Since $\overline{U}$ is stable under transfinite composition, we conclude that $f$ belongs to $S$.
This completes the proof of $(1)$.

Assertion $(2)$ is obvious. To prove $(3)$, we observe both $\bfA$ and $\CAlg(\bfA)$ are cotensored  over simplicial sets, and that we have canonical isomorphisms $\theta( A^K ) \simeq \theta(A)^K$
for $A \in \CAlg(\bfA)$, $K \in \sSet$. To prove that $\CAlg(\bfA)$ is a simplicial model category, it will suffice to show that $\CAlg(\bfA)$ is tensored over simplicial sets, and that given a fibration
$i: A \rightarrow A'$ in $\CAlg(\bfA)$ and a cofibration $j: K \rightarrow K'$ in $\sSet$, the induced map $A^{K'} \rightarrow A^{K} \times_{ {A'}^{K}} {A'}^{K'}$ is a fibration, trivial if either $i$ or
$j$ is a fibration. The second claim follows from the fact that $\theta$ detects fibrations and trivial fibrations. For the first, it suffices to prove that for $K \in \sSet$, the functor $A \mapsto A^{K}$
has a left adjoint; this follows from the adjoint functor theorem.
\end{proof}

The proof of Theorem \ref{cbeckify} rests on the following analogue of
Lemma \monoidref{GGG}:

\begin{lemma}\label{GGGG}
Let $\bfA$ be a simplicial symmetric monoidal model category, and let $\calC$ be a small category. 
Assume that $\bfA$ is combinatorial and freely powered, and that the simplicial set
$\Nerve(\calC)$ is sifted $($Definition \toposref{siftdef}$)$.
Then the forgetful functor $\Nerve( \CAlg( \bfA)^{\degree} ) \rightarrow \Nerve( \bfA^{\degree} )$
preserves $\Nerve(\calC)$-indexed colimits.
\end{lemma}

\begin{proof}
In view of Theorem \toposref{colimcomparee} and Proposition \toposref{gumby444}, it will suffice to prove that the forgetful functor $\theta: \CAlg(\bfA) \rightarrow \bfA$ preserves homotopy colimits indexed by $\calC$. Let us regard $\CAlg(\bfA)^{\calC}$ and $\bfA^{\calC}$ as endowed with
the projective model structure (see \S \toposref{quasilimit3}). Let $F: \bfA^{\calC} \rightarrow \bfA$
and $F_{\CAlg}: \CAlg(\bfA)^{\calC} \rightarrow \CAlg(\bfA)$ be colimit functors, and let
$\theta^{\calC}: \CAlg(\bfA)^{\calC} \rightarrow \bfA^{\calC}$ be given by composition with
$\theta$. Since $\Nerve(\calC)$ is sifted, there is a canonical isomorphism of functors
$\alpha: F \circ \theta^{\calC} \simeq \theta \circ F_{\CAlg}$. We wish to prove that this isomorphism persists after deriving all of the relevant functors. Since $\theta$ and $\theta^{\calC}$ preserve weak equivalences, they can be identified with their right derived functors. Let $LF$ and $LF_{\CAlg}$
be the left derived functors of $F$ and $F_{\CAlg}$, respectively. Then $\alpha$ induces a natural transformation $\overline{\alpha}: LF \circ \theta^{\calC} \rightarrow \theta \circ LF_{\CAlg}$; we wish to show that $\overline{\alpha}$ is an isomorphism. Let $A: \calC \rightarrow \CAlg(\bfA)$ be a projectively cofibrant object of $\CAlg(\bfA)^{\calC}$; we must show that the natural map
$$LF( \theta^{\calC}(A)) \rightarrow \theta( LF_{\CAlg}(A)) \simeq  \theta( F_{\CAlg}(A))
\simeq F( \theta^{\calC}(A))$$ is a weak equivalence in $\bfA$. 

Let us say that an object $X \in \bfA^{\calC}$ is {\it good} if each $X(C) \in \bfA$ is cofibrant, the colimit $F(X) \in \bfA$ is cofibrant, and the canonical map
the natural map $LF(X) \rightarrow F(X)$
is an isomorphism in the homotopy category $\h{\bfA}$ (in other words, the colimit of $X$ is also a homotopy colimit of $X$).
To complete the proof, it will suffice to show that $\theta^{\calC}(A)$ is good, whenever $A$
is a projectively cofibrant object of $\CAlg(\bfA)^{\calC}$. This is not obvious, since $\theta^{\calC}$ is a right Quillen functor and does not preserve projectively cofibrant objects in general (note that we have not yet used the full strength of our assumption that $\Nerve(\calC)$ is sifted). To continue the proof, we will need a relative version of the preceding condition. We will say that a morphism
$f: X \rightarrow Y$ in $\bfA^{\calC}$ is {\it good} if the following conditions are satisfied:
\begin{itemize}
\item[$(i)$] The objects $X,Y \in \bfA^{\calC}$ are good.
\item[$(ii)$] For each $C \in \calC$, the induced map $X(C) \rightarrow Y(C)$ is a cofibration in $\bfA$.
\item[$(iii)$] The map $F(X) \rightarrow F(Y)$ is a cofibration in $\bfA$.
\end{itemize}

As in the proof of Lemma \monoidref{GGG}, we have the following:

\begin{itemize}
\item[$(1)$] The collection of good morphisms is stable under transfinite composition.

\item[$(2)$] Suppose given a pushout diagram
$$ \xymatrix{ X \ar[r]^{f} \ar[d] & Y \ar[d] \\
X' \ar[r]^{f'} & Y' }$$
in $\bfA^{\calC}$. If $f$ is good and $X'$ is good, then $f'$ is good. 

\item[$(3)$] Let $F: \calC \rightarrow \bfA$ be a constant functor whose value is a cofibrant object of $\bfA$. Then $F$ is good. 

\item[$(4)$] Every projectively cofibrant object of $\bfA^{\calC}$ is good. Every strong
cofibration between projectively cofibrant objects of $\bfA^{\calC}$ is good.

\item[$(5)$] If $X$ and $Y$ are good objects of $\bfA^{\calC}$, then $X \otimes Y$ is good. 

\item[$(6)$] Let $f: X \rightarrow X'$ be a good morphism in $\bfA^{\calC}$, and let
$Y$ be a good object of $\bfA^{\calC}$. Then the morphism
$f \otimes \id_{Y}$ is good. 

\item[$(7)$] Let $f: X \rightarrow X'$ and $g: Y \rightarrow Y'$ be good morphisms in $\bfA^{\calC}$. Then $$f \wedge g: (X \otimes Y') \coprod_{ X \otimes Y} (X' \otimes Y) \rightarrow X' \otimes Y'$$
is good. 
\end{itemize}

Moreover, our assumption that $\bfA$ is freely powered ensures that the class of good morphisms has the following additional property:

\begin{itemize}
\item[$(8)$] Let $f: X \rightarrow Y$ be a good morphism in $\bfA^{\calC}$. Then the induced map
$\sigma^{n}(f): \Sym^{n}(Y;X) \rightarrow \Sym^{n}(Y)$ (see Notation \ref{puffstar}) is good. Condition $(ii)$ follows immediately from Lemma \ref{bullow}, and condition $(iii)$ follows from Lemma \ref{bullow} and the observation that $F( \sigma^{n}(f) ) = \sigma^{n}( F(f) )$ (since the functor $F$ commutes with colimits and tensor products). It will therefore suffice to show that the objects $\Sym^{n}(Y;X)$ and $\Sym^{n}(Y)$ are good. Let $D: \bfA^{\Sigma_{n} } \rightarrow \bfA$ be the coinvariants functor, and consider the following diagram of left Quillen functors (which commutes up to canonical isomorphism):
$$ \xymatrix{ \bfA^{ \calC \times \Sigma_{n} } \ar[r]^{D^{\calC}} \ar[d]^{F} & \bfA^{\calC} \ar[d]^{F} \\
\bfA^{\Sigma_{n}} \ar[r]^{D} & \bfA. }$$
Let us regard $\wedge^{n}(f)$ as an object in the category of arrows of $\h{ \bfA^{\calC \times \Sigma_{n} }}$. We wish to show that the canonical map
$LF( \sigma^{n}(f) ) \rightarrow F( \sigma^{n}(f) )$ is an isomorphism (in the category of morphisms in $\h{\bfA}$). We now observe that $\sigma^{n}(f) = D^{\calC} ( \wedge^{n}(f) )$. Lemma \ref{bullow} implies that the canonical map $LD^{\calC} \wedge^{n}(f) \rightarrow D^{\calC} \wedge^{n}(f)$
is a weak equivalence. It will therefore suffice to show that the transformation
$$\alpha: L( F \circ D^{\calC} )(\wedge^{n}(f)) \rightarrow (F \circ D^{\calC})(\wedge^{n}(f))$$
is an isomorphism (in the category of morphisms of $\h{\bfA}$). Using the commutativity of the above diagram, we can identify $\alpha$ with the map
$$ (LD \circ LF) \wedge^{n}(f) \rightarrow (D \circ F)(\wedge^{n}(f) ) = D \wedge^{n}( F(f) ).$$
Using Lemma \ref{bullow} again, we can identify the right hand side with $LD \wedge^{n}( F(f) )$.
It will therefore suffice to show that the map
$LF \wedge^{n}(f) \rightarrow F \wedge^{n}(f)$ is an isomorphism in the category of morphisms
of $\h{ \bfA^{\Sigma_{n} }}$. Since the forgetful functor $\bfA^{\Sigma_{n}} \rightarrow \bfA$ preserves homotopy colimits (it is also a left Quillen functor) and detects equivalences, we are reduced to proving that the morphism $\wedge^{n}(f)$ is good. This follows from $(7)$ using induction on $n$.
\end{itemize}

We observe that axiom $(F2)$ of Definition \ref{powerdef} has the following consequence:
\begin{itemize}
\item[$(F2')$] The collection of all projective cofibrations in $\bfA^{\calC}$ is generated by
projective cofibrations between projectively cofibrant objects.
\end{itemize}

Let $T: \bfA^{\calC} \rightarrow \CAlg(\bfA)^{\calC}$ be a left adjoint to $\theta^{\calC}$. 
Using the small object argument and $(B')$, we conclude that for every projectively cofibrant object $A \in \CAlg(\bfA)^{\calC}$ there exists a transfinite sequence $\{ A^{\beta} \}_{\beta \leq \alpha}$ in $\CAlg(\bfA)^{\calC}$ with the following properties:
\begin{itemize}
\item[$(a)$] The object $A^{0}$ is initial in $\CAlg(\bfA)^{\calC}$.
\item[$(b)$] The object $A$ is a retract of $A^{\alpha}$.
\item[$(c)$] If $\lambda \leq \alpha$ is a limit ordinal, then
$A^{\lambda} \simeq \colim \{ A^{\beta} \}_{ \beta < \lambda}$.  
\item[$(d)$] For each $\beta < \alpha$, there is a pushout diagram
$$ \xymatrix{ T(X') \ar[r]^{T(f)} \ar[d] & T(X) \ar[d] \\
A^{\beta} \ar[r] & A^{\beta+1} }$$
where $f$ is a projective cofibration between projectively cofibrant objects of $\bfA^{\calC}$. 
\end{itemize}

We wish to prove that $\theta^{\calC}(A)$ is good. In view of $(b)$, it will suffice to show
that $\theta^{\calC}(A^{\alpha})$ is good. We will prove a more general assertion: for
every $\gamma \leq \beta \leq \alpha$, the induced morphism $u_{\gamma, \beta}: 
\theta^{\calC}(A^{\gamma}) \rightarrow \theta^{\calC}( A^{\beta})$ is good. The proof is by induction on $\beta$. If $\beta = 0$, then we are reduced to proving that $\theta^{\calC}(A^0)$ is good. This follows from $(a)$ and $(3)$. If $\beta$ is a nonzero limit ordinal, then the desired result follows from $(c)$ and $(1)$. It therefore suffices to treat the case where $\beta = \beta' + 1$ is a successor ordinal. Moreover, we may suppose that
$\gamma = \beta'$: if $\gamma < \beta'$, then we observe that
$u_{\gamma, \beta} = u_{\beta', \beta} \circ u_{ \gamma, \beta'}$ and invoke $(1)$, while
if $\gamma > \beta'$, then $\gamma = \beta$ and we are reduced to proving that
$\theta^{\calC}(A^{\beta})$ is good, which follows from the assertion that
$u_{\beta', \beta}$ is good. We are now reduced to proving the following:

\begin{itemize}
\item[$(\ast)$] Let $$ \xymatrix{ T(X') \ar[r]^{T(f)} \ar[d] & T(X) \ar[d] \\
B' \ar[r]^{v} & B }$$
be a pushout diagram in $\CAlg(\bfA)^{\calC}$, where $f: X' \rightarrow X$ is a projective cofibration between projectively cofibrant objects of $\bfA^{\calC}$. If $\theta^{\calC}(B')$ is good, then $\theta^{\calC}(v)$ is good.
\end{itemize}

To prove $(\ast)$, we set $Y = \theta^{\calC}(B) \in \bfA^{\calC}$, $Y' = \theta^{\calC}(B') \in \bfA^{\calC}$. Let $g: \emptyset \rightarrow Y'$ the unique morphism, where $\emptyset$
denotes an initial object of $\bfA^{\calC}$. As in the proof of Proposition \ref{combocheese}, 
$Y$ can be identified with the colimit of a sequence
$$ Y' = Y^{(0)} \stackrel{w_1}{\rightarrow} Y^{(1)} \stackrel{w_2}{\rightarrow} \ldots $$
where $Y^{(0)} = Y'$, and $w_{k}$ is a pushout of the morphism
$f^{(k)} = B' \otimes \sigma^{k}(f)$. The desired result now follows immediately from $(4)$, $(6)$ and $(8)$.
\end{proof}

We are now ready to prove our main result:

\begin{proof}[Proof of Theorem \ref{cbeckify}]
Consider the diagram
$$ \xymatrix{ \Nerve( \CAlg(\bfA)^{\degree}) \ar[rr] \ar[dr]^{G} & & \CAlg( \Nerve( \bfA^{\degree}) )\ar[dl]^{G'} \\
& \Nerve(\bfA^{\degree}). & }$$
It will suffice to show that this diagram satisfies the hypotheses of Corollary \monoidref{littlerbeck}:
\begin{itemize}
\item[$(a)$] The $\infty$-categories $\Nerve( \CAlg(\bfA)^{\degree} )$ and
$\CAlg( \Nerve(\bfA^{\degree}))$ admit geometric realizations of simplicial objects. In fact, both of these $\infty$-categories are presentable. For $\Nerve( \CAlg(\bfA)^{\degree})$, this follows from
Propositions \toposref{notthereyet} and \ref{combocheese}. For $\CAlg( \Nerve(\bfA)^{\degree})$, we first observe that $\Nerve(\bfA)^{\degree}$ is presentable (Proposition \toposref{notthereyet}) and that
the tensor product preserves colimits separately in each variable, and then apply Corollary \ref{algprice1}.

\item[$(b)$] The functors $G$ and $G'$ admit left adjoints $F$ and $F'$. The existence of a left adjoint to $G$ follows from the fact that $G$ is given by a right Quillen functor. The existence of a left adjoint to $G'$ follows from Corollary \ref{spaltwell}.

\item[$(c)$] The functor $G'$ is conservative and preserves geometric realizations of simplicial objects. This follows from Corollaries \ref{filtfemme} and \ref{jumunj22}.

\item[$(d)$] The functor $G$ is conservative and preserves geometric realizations of simplicial objects. The first assertion is immediate from the definition of the weak equivalences in
$\CAlg(\bfA)$, and the second follows from Lemma \ref{GGGG}.

\item[$(e)$] The canonical map $G' \circ F' \rightarrow G \circ F$ is an equivalence of functors. 
In other words, we must show that for every object $C \in \calC$, a fibrant replacement for the free strictly commutative algebra $\coprod_{ n} \Sym^{n}(C) \in \CAlg(\bfA)$ is a free algebra generated by $C$, in the sense of Definition \ref{dfree}. In view of Remark \ref{tappus}, it suffices to show that
the colimit defining the total symmetric power $\coprod_{n} \Sym^{n}(C)$ in $\bfA$ is also a homotopy colimit. This follows immediately from part $(3)$ of Lemma \ref{bullow}.
\end{itemize}
\end{proof}


\end{document}